\documentclass[12pt]{amsbook}
\usepackage{amssymb}
\usepackage[all]{xy}

\usepackage[colorlinks=true,linkcolor=magenta,citecolor=magenta]{hyperref}
\makeindex

\newtheorem{theorem}{Theorem}[chapter]
\newtheorem{answer}[theorem]{Answer}
\newtheorem{comment}[theorem]{Comment}
\newtheorem{conclusion}[theorem]{Conclusion}
\newtheorem{conjecture}[theorem]{Conjecture}

\newtheorem{definition}[theorem]{Definition}
\newtheorem{example}[theorem]{Example}
\newtheorem{exercise}[theorem]{Exercise}
\newtheorem{fact}[theorem]{Fact}
\newtheorem{principle}[theorem]{Principle}
\newtheorem{proposition}[theorem]{Proposition}
\newtheorem{question}[theorem]{Question}
\newtheorem{thought}[theorem]{Thought}
\newtheorem{warning}[theorem]{Warning}

\begin{document}

\title{Principles of operator algebras}

\author{Teo Banica}
\address{Department of Mathematics, University of Cergy-Pontoise, F-95000 Cergy-Pontoise, France. {\tt teo.banica@gmail.com}}


\begin{abstract}
This is an introduction to the algebras $A\subset B(H)$ that the linear operators $T:H\to H$ can form, once a complex Hilbert space $H$ is given. Motivated by quantum mechanics, we are mainly interested in the von Neumann algebras, which are stable under taking adjoints, $T\to T^*$, and are weakly closed. When the algebra has a trace $tr:A\to\mathbb C$, we can think of it as being of the form $A=L^\infty(X)$, with $X$ being a quantum measured space. Of particular interest is the free case, where the center of the algebra reduces to the scalars, $Z(A)=\mathbb C$. Following von Neumann, Connes, Jones, Voiculescu and others, we discuss the basic properties of such algebras $A$, and how to do algebra, geometry, analysis and probability on the underlying quantum spaces $X$.
\end{abstract}

\maketitle

\chapter*{Preface}

Quantum mechanics as we know it is the source of many puzzling questions. The simplest quantum mechanical system is the hydrogen atom, consisting of a negative charge, an electron, moving around a positive charge, a proton. This reminds electrodynamics, and accepting the fact that the electron is a bit of a slippery particle, whose position and speed are described by probability, rather than by exact formulae, the hydrogen atom can indeed be solved, by starting with electrodynamics, and making a long series of corrections, for the most coming from experiments, but sometimes coming as well from intuition, with the idea in mind that beautiful mathematics should correspond to true physics. The solution, as we presently know it, is something quite complicated.

\bigskip

Mathematically, the commonly accepted belief is that the good framework for the study of quantum mechanics is an infinite dimensional complex Hilbert space $H$, whose vectors can be thought of as being states of the system, and with the linear operators $T:H\to H$ corresponding to the observables. This is however to be taken with care, because in order to do ``true physics'', things must be far sharper than that. Always remember indeed that the simplest object of quantum mechanics is the hydrogen atom, whose simplest states and observables are something quite complicated. Thus when talking about ``states and observables'', we have a whole continuum of possible considerations and theories, ranging from true physics to very abstract mathematics.

\bigskip

For making things worse, even the existence and exact relevance of the Hilbert space $H$ is subject to debate. This is something more philosophical, related to the 2-body hydrogen problem evoked above, which has twisted the minds of many scientists, starting with Einstein and others. Can we get someday to a better quantum mechanics, by adding more variables to those available inside $H$? No one really knows the answer here.

\bigskip

The present book is an introduction to the algebras $A\subset B(H)$ that the bounded linear operators $T:H\to H$ can form, once a Hilbert space $H$ is given. There has been an enormous amount of work on such algebras, starting with von Neumann in the 1930s, and we will insist here on the aspects which are beautiful. With the idea, or rather hope in mind, that beautiful mathematics should correspond to true physics.

\bigskip

So, what is beauty, in the operator algebra framework? In our opinion, the source of all possible beauty is an old result of von Neumann, related to the Spectral Theorem for normal operators, which states that any commutative von Neumann algebra $A\subset B(H)$ must be of the form $A=L^\infty(X)$, with $X$ being a measured space.

\bigskip

This is something subtle and interesting, which suggests doing several things with the von Neumann algebras $A\subset B(H)$. Given such an algebra we can write the center as $Z(A)=L^\infty(X)$, we have then a decomposition of type $A=\int_XA_xdx$, and the problem is that of understanding the structure of the fibers, called ``factors''. This is what von Neumann himself, and then Connes and others, did. Another idea, more speculative, following later work of Connes, and in parallel work of Voiculescu, is that of writing $A=L^\infty(X)$, with $X$ being an abstract ``quantum measured space'', and then trying to understand the geometry and  probabilistic theory of $X$. Finally, yet another beautiful idea, due this time to Jones, is that of looking at the inclusions $A_0\subset A_1$ of von Neumann algebras, instead at the von Neumann algebras themselves, the point being that the ``symmetries'' of such an inclusion lead to interesting combinatorics.

\bigskip

All in all, many things that can be done with a von Neumann algebra $A\subset B(H)$, and explaining the basics, plus having a look at the above 4 directions of research, is already what a medium sized book can cover. And this book is written exactly with this idea in mind. We will talk about all the above, keeping things as simple as possible, and with everything being accessible with a minimal knowledge of undergraduate mathematics.

\bigskip

The book is organized in 4 parts, with Part I explaining the basics of operator theory, Part II explaining the basics of operator algebras, with a look into geometry and probability too, then Part III going into the structure of the von Neumann factors, and finally Part IV being an introduction to the subfactor theory of Jones.

\bigskip

This book contains, besides the basics of the operator algebra theory, some modern material as well, namely quantum group illustrations for pretty much everything, and I am grateful to Julien Bichon, Beno\^ it Collins, Steve Curran and the others, for our joint work. Many thanks go as well to my cats. Their views and opinions on mathematics, and knowledge of advanced functional analysis, have always been of great help.

\bigskip

\

{\em Cergy, August 2024}

\smallskip

{\em Teo Banica}

\baselineskip=15.95pt
\tableofcontents
\baselineskip=14pt

\part{Linear operators}

\ \vskip50mm

\begin{center}
{\em Does anybody here remember Vera Lynn

Remember how she said that

We would meet again

Some sunny day}
\end{center}

\chapter{Linear algebra}

\section*{1a. Linear maps}

According to various findings in physics, starting with those of Heisenberg from the early 1920s, basic quantum mechanics involves linear operators $T:H\to H$ from a complex Hilbert space $H$ to itself. The space $H$ is typically infinite dimensional, a basic example being the Schr\"odinger space $H=L^2(\mathbb R^3)$ of the wave functions $\psi:\mathbb R^3\to\mathbb C$ of the electron. In fact, in what regards the electron, this space $H=L^2(\mathbb R^3)$ is basically the correct one, with the only adjustment needed, due to Pauli and others, being that of tensoring with a copy of $K=\mathbb C^2$, in order to account for the electron spin.

\bigskip

But more on this later. Let us start this book more modestly, as follows:

\begin{fact}
We are interested in quantum mechanics, taking place in infinite dimensions, but as a main source of inspiration we will have $H=\mathbb C^N$, with scalar product
$$<x,y>=\sum_ix_i\bar{y}_i$$
with the linearity at left being the standard mathematical convention. More specifically, we will be interested in the mathematics of the linear operators $T:H\to H$.
\end{fact}

The point now, that you surely know about, is that the above operators $T:H\to H$ correspond to the square matrices $A\in M_N(\mathbb C)$. Thus, as a preliminary to what we want to do in this book, we need a good knowledge of linear algebra over $\mathbb C$. 

\bigskip

You probably know well linear algebra, but always good to recall this, and this will be the purpose of the present chapter. Let us start with the very basics:

\begin{theorem}
The linear maps $T:\mathbb C^N\to\mathbb C^N$ are in correspondence with the square matrices $A\in M_N(\mathbb C)$, with the linear map associated to such a matrix being
$$Tx=Ax$$
and with the matrix associated to a linear map being $A_{ij}=<Te_j,e_i>$.
\end{theorem}

\begin{proof}
The first assertion is clear, because a linear map $T:\mathbb C^N\to\mathbb C^N$ must send a vector $x\in\mathbb C^N$ to a certain vector $Tx\in\mathbb C^N$, all whose components are linear combinations of the components of $x$. Thus, we can write, for certain complex numbers $A_{ij}\in\mathbb C$:
$$T\begin{pmatrix}
x_1\\
\vdots\\
\vdots\\
x_N
\end{pmatrix}
=\begin{pmatrix}
A_{11}x_1+\ldots+A_{1N}x_N\\
\vdots\\
\vdots\\
A_{N1}x_1+\ldots+A_{NN}x_N
\end{pmatrix}$$

Now the parameters $A_{ij}\in\mathbb C$ can be regarded as being the entries of a square matrix $A\in M_N(\mathbb C)$, and with the usual convention for matrix multiplication, we have:
$$Tx=Ax$$

Regarding the second assertion, with $Tx=Ax$ as above, if we denote by $e_1,\ldots,e_N$ the standard basis of $\mathbb C^N$, then we have the following formula:
$$Te_j
=\begin{pmatrix}
A_{1j}\\
\vdots\\
\vdots\\
A_{Nj}
\end{pmatrix}$$

But this gives the second formula, $<Te_j,e_i>=A_{ij}$, as desired.
\end{proof}

Our claim now is that, no matter what we want to do with $T$ or $A$, of advanced type, we will run at some point into their adjoints $T^*$ and $A^*$, constructed as follows:

\index{adjoint operator}

\begin{theorem}
The adjoint operator $T^*:\mathbb C^N\to\mathbb C^N$, which is given by
$$<Tx,y>=<x,T^*y>$$
corresponds to the adjoint matrix $A^*\in M_N(\mathbb C)$, given by
$$(A^*)_{ij}=\bar{A}_{ji}$$
via the correspondence between linear maps and matrices constructed above.
\end{theorem}

\begin{proof}
Given a linear map $T:\mathbb C^N\to\mathbb C^N$, fix $y\in\mathbb C^N$, and consider the linear form $\varphi(x)=<Tx,y>$. This form must be as follows, for a certain vector $T^*y\in\mathbb C^N$:
$$\varphi(x)=<x,T^*y>$$

Thus, we have constructed a map $y\to T^*y$ as in the statement, which is obviously linear, and that we can call $T^*$. Now by taking the vectors $x,y\in\mathbb C^N$ to be elements of the standard basis of $\mathbb C^N$, our defining formula for $T^*$ reads:
$$<Te_i,e_j>=<e_i,T^*e_j>$$

By reversing the scalar product on the right, this formula can be written as:
$$<T^*e_j,e_i>=\overline{<Te_i,e_j>}$$

But this means that the matrix of $T^*$ is given by $(A^*)_{ij}=\bar{A}_{ji}$, as desired.
\end{proof}

Getting back to our claim, the adjoints $*$ are indeed ubiquitous, as shown by:

\index{scalar product}
\index{unitary}
\index{projection}

\begin{theorem}
The following happen:
\begin{enumerate}
\item $T(x)=Ux$ with $U\in M_N(\mathbb C)$ is an isometry precisely when $U^*=U^{-1}$.

\item $T(x)=Px$ with $P\in M_N(\mathbb C)$ is a projection precisely when $P^2=P^*=P$.
\end{enumerate}
\end{theorem}

\begin{proof}
Let us first recall that the lengths, or norms, of the vectors $x\in\mathbb C^N$ can be recovered from the knowledge of the scalar products, as follows:
$$||x||=\sqrt{<x,x>}$$

Conversely, we can recover the scalar products out of norms, by using the following difficult to remember formula, called complex polarization identity:
$$4<x,y>
=||x+y||^2-||x-y||^2+i||x+iy||^2-i||x-iy||^2$$

The proof of this latter formula is indeed elementary, as follows:
\begin{eqnarray*}
&&||x+y||^2-||x-y||^2+i||x+iy||^2-i||x-iy||^2\\
&=&||x||^2+||y||^2-||x||^2-||y||^2+i||x||^2+i||y||^2-i||x||^2-i||y||^2\\
&&+2Re(<x,y>)+2Re(<x,y>)+2iIm(<x,y>)+2iIm(<x,y>)\\
&=&4<x,y>
\end{eqnarray*}

Finally, we will use Theorem 1.3, and more specifically the following formula coming from there, valid for any matrix $A\in M_N(\mathbb C)$ and any two vectors $x,y\in\mathbb C^N$:
$$<Ax,y>=<x,A^*y>$$

(1) Given a matrix $U\in M_N(\mathbb C)$, we have indeed the following equivalences, with the first one coming from the polarization identity, and the other ones being clear:
\begin{eqnarray*}
||Ux||=||x||
&\iff&<Ux,Uy>=<x,y>\\
&\iff&<x,U^*Uy>=<x,y>\\
&\iff&U^*Uy=y\\
&\iff&U^*U=1\\
&\iff&U^*=U^{-1}
\end{eqnarray*}

(2) Given a matrix $P\in M_N(\mathbb C)$, in order for $x\to Px$ to be an oblique projection, we must have $P^2=P$. Now observe that this projection is orthogonal when:
\begin{eqnarray*}
<Px-x,Py>=0
&\iff&<P^*Px-P^*x,y>=0\\
&\iff&P^*Px-P^*x=0\\
&\iff&P^*P-P^*=0\\
&\iff&P^*P=P^*
\end{eqnarray*}

The point now is that by conjugating the last formula, we obtain $P^*P=P$. Thus we must have $P=P^*$, and this gives the result. 
\end{proof}

Summarizing, the linear operators come in pairs $T,T^*$, and the associated matrices come as well in pairs $A,A^*$. This is something quite interesting, philosophically speaking, and will keep this in mind, and come back to it later, on numerous occasions.

\section*{1b. Diagonalization}

Let us discuss now the diagonalization question for the linear maps and matrices. Again, we will be quite brief here, and for more, we refer to any standard linear algebra book. By the way, there will be some complex analysis involved too, and here we refer to Rudin \cite{rud}. Which book of Rudin will be in fact the one and only true prerequisite for reading the present book, but more on references and reading later. 

\bigskip

The basic diagonalization theory, formulated in terms of matrices, is as follows:

\index{eigenvalue}
\index{eigenvector}
\index{diagonalization}

\begin{proposition}
A vector $v\in\mathbb C^N$ is called eigenvector of $A\in M_N(\mathbb C)$, with corresponding eigenvalue $\lambda$, when $A$ multiplies by $\lambda$ in the direction of $v$:
$$Av=\lambda v$$
In the case where $\mathbb C^N$ has a basis $v_1,\ldots,v_N$ formed by eigenvectors of $A$, with corresponding eigenvalues $\lambda_1,\ldots,\lambda_N$, in this new basis $A$ becomes diagonal, as follows:
$$A\sim\begin{pmatrix}\lambda_1\\&\ddots\\&&\lambda_N\end{pmatrix}$$
Equivalently, if we denote by $D=diag(\lambda_1,\ldots,\lambda_N)$ the above diagonal matrix, and by $P=[v_1\ldots v_N]$ the square matrix formed by the eigenvectors of $A$, we have:
$$A=PDP^{-1}$$
In this case we say that the matrix $A$ is diagonalizable.
\end{proposition}

\begin{proof}
This is something which is clear, the idea being as follows:

\medskip

(1) The first assertion is clear, because the matrix which multiplies each basis element $v_i$ by a number $\lambda_i$ is precisely the diagonal matrix $D=diag(\lambda_1,\ldots,\lambda_N)$.

\medskip

(2) The second assertion follows from the first one, by changing the basis. We can prove this by a direct computation as well, because we have $Pe_i=v_i$, and so:
\begin{eqnarray*}
PDP^{-1}v_i
&=&PDe_i\\
&=&P\lambda_ie_i\\
&=&\lambda_iPe_i\\
&=&\lambda_iv_i
\end{eqnarray*}

Thus, the matrices $A$ and $PDP^{-1}$ coincide, as stated.
\end{proof}

Let us recall as well that the basic example of a non diagonalizable matrix, over the complex numbers as above, is the following matrix:
$$J=\begin{pmatrix}0&1\\0&0\end{pmatrix}$$

Indeed, we have $J\binom{x}{y}=\binom{y}{0}$, so the eigenvectors are the vectors of type $\binom{x}{0}$, all with eigenvalue $0$. Thus, we have not enough eigenvectors for constructing a basis of $\mathbb C^2$.

\bigskip

In general, in order to study the diagonalization problem, the idea is that the eigenvectors can be grouped into linear spaces, called eigenspaces, as follows:

\begin{theorem}
Let $A\in M_N(\mathbb C)$, and for any eigenvalue $\lambda\in\mathbb C$ define the corresponding eigenspace as being the vector space formed by the corresponding eigenvectors:
$$E_\lambda=\left\{v\in\mathbb C^N\Big|Av=\lambda v\right\}$$
These eigenspaces $E_\lambda$ are then in a direct sum position, in the sense that given vectors $v_1\in E_{\lambda_1},\ldots,v_k\in E_{\lambda_k}$ corresponding to different eigenvalues $\lambda_1,\ldots,\lambda_k$, we have:
$$\sum_ic_iv_i=0\implies c_i=0$$
In particular we have the following estimate, with sum over all the eigenvalues, 
$$\sum_\lambda\dim(E_\lambda)\leq N$$
and our matrix is diagonalizable precisely when we have equality.
\end{theorem}

\begin{proof}
We prove the first assertion by recurrence on $k\in\mathbb N$. Assume by contradiction that we have a formula as follows, with the scalars $c_1,\ldots,c_k$ being not all zero:
$$c_1v_1+\ldots+c_kv_k=0$$

By dividing by one of these scalars, we can assume that our formula is:
$$v_k=c_1v_1+\ldots+c_{k-1}v_{k-1}$$

Now let us apply $A$ to this vector. On the left we obtain:
$$Av_k
=\lambda_kv_k
=\lambda_kc_1v_1+\ldots+\lambda_kc_{k-1}v_{k-1}$$

On the right we obtain something different, as follows:
\begin{eqnarray*}
A(c_1v_1+\ldots+c_{k-1}v_{k-1})
&=&c_1Av_1+\ldots+c_{k-1}Av_{k-1}\\
&=&c_1\lambda_1v_1+\ldots+c_{k-1}\lambda_{k-1}v_{k-1}
\end{eqnarray*}

We conclude from this that the following equality must hold:
$$\lambda_kc_1v_1+\ldots+\lambda_kc_{k-1}v_{k-1}=c_1\lambda_1v_1+\ldots+c_{k-1}\lambda_{k-1}v_{k-1}$$

On the other hand, we know by recurrence that the vectors $v_1,\ldots,v_{k-1}$ must be linearly independent. Thus, the coefficients must be equal, at right and at left:
$$\lambda_kc_1=c_1\lambda_1$$
$$\vdots$$
$$\lambda_kc_{k-1}=c_{k-1}\lambda_{k-1}$$

Now since at least one of the numbers $c_i$ must be nonzero, from $\lambda_kc_i=c_i\lambda_i$ we obtain $\lambda_k=\lambda_i$, which is a contradiction. Thus our proof by recurrence of the first assertion is complete. As for the second assertion, this follows from the first one.
\end{proof}

In order to reach now to more advanced results, we can use the characteristic polynomial, which appears via the following fundamental result:

\index{characteristic polynomial}

\begin{theorem}
Given a matrix $A\in M_N(\mathbb C)$, consider its characteristic polynomial:
$$P(x)=\det(A-x1_N)$$
The eigenvalues of $A$ are then the roots of $P$. Also, we have the inequality
$$\dim(E_\lambda)\leq m_\lambda$$
where $m_\lambda$ is the multiplicity of $\lambda$, as root of $P$.
\end{theorem}

\begin{proof}
The first assertion follows from the following computation, using the fact that a linear map is bijective when the determinant of the associated matrix is nonzero:
\begin{eqnarray*}
\exists v,Av=\lambda v
&\iff&\exists v,(A-\lambda 1_N)v=0\\
&\iff&\det(A-\lambda 1_N)=0
\end{eqnarray*}

Regarding now the second assertion, given an eigenvalue $\lambda$ of our matrix $A$, consider the dimension $d_\lambda=\dim(E_\lambda)$ of the corresponding eigenspace. By changing the basis of $\mathbb C^N$, as for the eigenspace $E_\lambda$ to be spanned by the first $d_\lambda$ basis elements, our matrix becomes as follows, with $B$ being a certain smaller matrix:
$$A\sim\begin{pmatrix}\lambda 1_{d_\lambda}&0\\0&B\end{pmatrix}$$

We conclude that the characteristic polynomial of $A$ is of the following form:
$$P_A
=P_{\lambda 1_{d_\lambda}}P_B
=(\lambda-x)^{d_\lambda}P_B$$

Thus the multiplicity $m_\lambda$ of our eigenvalue $\lambda$, as a root of $P$, satisfies $m_\lambda\geq d_\lambda$, and this leads to the conclusion in the statement.
\end{proof}

Now recall that we are over $\mathbb C$, which is something that we have not used yet, in our last two statements. And the point here is that we have the following key result:

\begin{theorem}
Any polynomial $P\in\mathbb C[X]$ decomposes as
$$P=c(X-a_1)\ldots (X-a_N)$$
with $c\in\mathbb C$ and with $a_1,\ldots,a_N\in\mathbb C$.
\end{theorem}

\begin{proof}
It is enough to prove that $P$ has one root, and we do this by contradiction. Assume that $P$ has no roots, and pick a number $z\in\mathbb C$ where $|P|$ attains its minimum:
$$|P(z)|=\min_{x\in\mathbb C}|P(x)|>0$$ 

Since $Q(t)=P(z+t)-P(z)$ is a polynomial which vanishes at $t=0$, this polynomial must be of the form $ct^k$ + higher terms, with $c\neq0$, and with $k\geq1$ being an integer. We obtain from this that, with $t\in\mathbb C$ small, we have the following estimate:
$$P(z+t)\simeq P(z)+ct^k$$

Now let us write $t=rw$, with $r>0$ small, and with $|w|=1$. Our estimate becomes:
$$P(z+rw)\simeq P(z)+cr^kw^k$$

Now recall that we have assumed $P(z)\neq0$. We can therefore choose $w\in\mathbb T$ such that $cw^k$ points in the opposite direction to that of $P(z)$, and we obtain in this way:
$$|P(z+rw)|
\simeq|P(z)+cr^kw^k|
=|P(z)|(1-|c|r^k)$$

Now by choosing $r>0$ small enough, as for the error in the first estimate to be small, and overcame by the negative quantity $-|c|r^k$, we obtain from this:
$$|P(z+rw)|<|P(z)|$$

But this contradicts our definition of $z\in\mathbb C$, as a point where $|P|$ attains its minimum. Thus $P$ has a root, and by recurrence it has $N$ roots, as stated.
\end{proof}

Now by putting everything together, we obtain the following result:

\index{eigenvalue}
\index{eigenvector}
\index{characteristic polynomial}
\index{diagonalization}

\begin{theorem}
Given a matrix $A\in M_N(\mathbb C)$, consider its characteristic polynomial
$$P(X)=\det(A-X1_N)$$ 
then factorize this polynomial, by computing the complex roots, with multiplicities,
$$P(X)=(-1)^N(X-\lambda_1)^{n_1}\ldots(X-\lambda_k)^{n_k}$$
and finally compute the corresponding eigenspaces, for each eigenvalue found:
$$E_i=\left\{v\in\mathbb C^N\Big|Av=\lambda_iv\right\}$$
The dimensions of these eigenspaces satisfy then the following inequalities,
$$\dim(E_i)\leq n_i$$
and $A$ is diagonalizable precisely when we have equality for any $i$.
\end{theorem}

\begin{proof}
This follows by combining Theorem 1.6, Theorem 1.7 and Theorem 1.8. Indeed, the statement is well formulated, thanks to Theorem 1.8. By summing the inequalities $\dim(E_\lambda)\leq m_\lambda$ from Theorem 1.7, we obtain an inequality as follows:
$$\sum_\lambda\dim(E_\lambda)\leq\sum_\lambda m_\lambda\leq N$$

On the other hand, we know from Theorem 1.6 that our matrix is diagonalizable when we have global equality. Thus, we are led to the conclusion in the statement.
\end{proof}

This was for the main result of linear algebra. There are countless applications of this, and generally speaking, advanced linear algebra consists in building on Theorem 1.9.

\bigskip

In practice, diagonalizing a matrix remains something quite complicated. Let us record a useful algorithmic version of the above result, as follows:

\begin{theorem}
The square matrices $A\in M_N(\mathbb C)$ can be diagonalized as follows:
\begin{enumerate}
\item Compute the characteristic polynomial.

\item Factorize the characteristic polynomial.

\item Compute the eigenvectors, for each eigenvalue found.

\item If there are no $N$ eigenvectors, $A$ is not diagonalizable.

\item Otherwise, $A$ is diagonalizable, $A=PDP^{-1}$.
\end{enumerate}
\end{theorem}

\begin{proof}
This is an informal reformulation of Theorem 1.9, with (4) referring to the total number of linearly independent eigenvectors found in (3), and with $A=PDP^{-1}$ in (5) being the usual diagonalization formula, with $P,D$ being as before.
\end{proof}

As an illustration for all this, which is a must-know computation, we have:

\begin{proposition}
The rotation of angle $t\in\mathbb R$ in the plane diagonalizes as:
$$\begin{pmatrix}\cos t&-\sin t\\ \sin t&\cos t\end{pmatrix}
=\frac{1}{2}\begin{pmatrix}1&1\\i&-i\end{pmatrix}
\begin{pmatrix}e^{-it}&0\\0&e^{it}\end{pmatrix}
\begin{pmatrix}1&-i\\1&i\end{pmatrix}$$
Over the reals this is impossible, unless $t=0,\pi$, where the rotation is diagonal.
\end{proposition}

\begin{proof}
Observe first that, as indicated, unlike we are in the case $t=0,\pi$, where our rotation is $\pm1_2$, our rotation is a ``true'' rotation, having no eigenvectors in the plane. Fortunately the complex numbers come to the rescue, via the following computation:
$$\begin{pmatrix}\cos t&-\sin t\\ \sin t&\cos t\end{pmatrix}\binom{1}{i}
=\binom{\cos t-i\sin t}{i\cos t+\sin t}
=e^{-it}\binom{1}{i}$$

We have as well a second complex eigenvector, coming from:
$$\begin{pmatrix}\cos t&-\sin t\\ \sin t&\cos t\end{pmatrix}\binom{1}{-i}
=\binom{\cos t+i\sin t}{-i\cos t+\sin t}
=e^{it}\binom{1}{-i}$$

Thus, we are led to the conclusion in the statement.
\end{proof}

\section*{1c. Matrix tricks}

At the level of basic examples of diagonalizable matrices, we first have the following result, which provides us with the ``generic'' examples:

\begin{theorem}
For a matrix $A\in M_N(\mathbb C)$ the following conditions are equivalent,
\begin{enumerate}
\item The eigenvalues are different, $\lambda_i\neq\lambda_j$,

\item The characteristic polynomial $P$ has simple roots,

\item The characteristic polynomial satisfies $(P,P')=1$,

\item The resultant of $P,P'$ is nonzero, $R(P,P')\neq0$,

\item The discriminant of $P$ is nonzero, $\Delta(P)\neq0$,
\end{enumerate}
and in this case, the matrix is diagonalizable.
\end{theorem}

\begin{proof}
The last assertion holds indeed, due to Theorem 1.9. As for the equivalences in the statement, these are all standard, the idea for their proofs, along with some more theory, needed for using in practice the present result, being as follows:

\medskip

$(1)\iff(2)$ This follows from Theorem 1.9.

\medskip

$(2)\iff(3)$ This is standard, the double roots of $P$ being roots of $P'$.

\medskip

$(3)\iff(4)$ The idea here is that associated to any two polynomials $P,Q$ is their resultant $R(P,Q)$, which checks whether $P,Q$ have a common root. Let us write:
$$P=c(X-a_1)\ldots(X-a_k)$$
$$Q=d(X-b_1)\ldots(X-b_l)$$

We can define then the resultant as being the following quantity:
$$R(P,Q)=c^ld^k\prod_{ij}(a_i-b_j)$$

The point now, that we will explain as well, is that this is a polynomial in the coefficients of $P,Q$, with integer coefficients. Indeed, this can be checked as follows:

\medskip

-- We can expand the formula of $R(P,Q)$, and in what regards $a_1,\ldots,a_k$, which are the roots of $P$, we obtain in this way certain symmetric functions in these variables, which will be therefore polynomials in the coefficients of $P$, with integer coefficients.

\medskip

-- We can then look what happens with respect to the remaining variables $b_1,\ldots,b_l$, which are the roots of $Q$. Once again what we have here are certain symmetric functions, and so polynomials in the coefficients of $Q$, with integer coefficients.

\medskip

-- Thus, we are led to the above conclusion, that $R(P,Q)$ is a polynomial in the coefficients of $P,Q$, with integer coefficients, and with the remark that the $c^ld^k$ factor is there for these latter coefficients to be indeed integers, instead of rationals.

\medskip

Alternatively, let us write our two polynomials in usual form, as follows:
$$P=p_kX^k+\ldots+p_1X+p_0$$
$$Q=q_lX^l+\ldots+q_1X+q_0$$

The corresponding resultant appears then as the determinant of an associated matrix, having size $k+l$, and having $0$ coefficients at the blank spaces, as follows:
$$R(P,Q)=
\begin{vmatrix}
p_k&&&q_l\\
\vdots&\ddots&&\vdots&\ddots\\
p_0&&p_k&q_0&&q_l\\
&\ddots&\vdots&&\ddots&\vdots\\
&&p_0&&&q_0
\end{vmatrix}
$$

$(4)\iff(5)$ Once again this is something standard, the idea here being that the discriminant $\Delta(P)$ of a polynomial $P\in\mathbb C[X]$ is, modulo scalars, the resultant $R(P,P')$. To be more precise, let us write our polynomial as follows:
$$P(X)=cX^N+dX^{N-1}+\ldots$$

Its discriminant is then defined as being the following quantity:
$$\Delta(P)=\frac{(-1)^{\binom{N}{2}}}{c}R(P,P')$$

This is a polynomial in the coefficients of $P$, with integer coefficients, with the division by $c$ being indeed possible, under $\mathbb Z$, and with the sign being there for various reasons, including the compatibility with some well-known formulae, at small values of $N$.
\end{proof}

All the above might seem a bit complicated, so as an illustration, let us work out an example. Consider the case of a polynomial of degree 2, and a polynomial of degree 1:
$$P=ax^2+bx+c\quad,\quad 
Q=dx+e$$

In order to compute the resultant, let us factorize our polynomials:
$$P=a(x-p)(x-q)\quad,\quad 
Q=d(x-r)$$

The resultant can be then computed as follows, by using the two-step method:
\begin{eqnarray*}
R(P,Q)
&=&ad^2(p-r)(q-r)\\
&=&ad^2(pq-(p+q)r+r^2)\\
&=&cd^2+bd^2r+ad^2r^2\\
&=&cd^2-bde+ae^2
\end{eqnarray*}

Observe that $R(P,Q)=0$ corresponds indeed to the fact that $P,Q$ have a common root. Indeed, the root of $Q$ is $r=-e/d$, and we have:
$$P(r)
=\frac{ae^2}{d^2}-\frac{be}{d}+c
=\frac{R(P,Q)}{d^2}$$

We can recover as well the resultant as a determinant, as follows:
$$R(P,Q)
=\begin{vmatrix}
a&d&0\\
b&e&d\\
c&0&e
\end{vmatrix}
=ae^2-bde+cd^2$$

Finally, in what regards the discriminant, let us see what happens in degree 2. Here we must compute the resultant of the following two polynomials:
$$P=aX^2+bX+c\quad,\quad 
P'=2aX+b$$

The resultant is then given by the following formula:
\begin{eqnarray*}
R(P,P')
&=&ab^2-b(2a)b+c(2a)^2\\
&=&4a^2c-ab^2\\
&=&-a(b^2-4ac)
\end{eqnarray*}

Now by doing the discriminant normalizations, we obtain, as we should:
$$\Delta(P)=b^2-4ac$$

As already mentioned, one can prove that the matrices having distinct eigenvalues are ``generic'', and so the above result basically captures the whole situation. We have in fact the following collection of density results, which are quite advanced:

\begin{theorem}
The following happen, inside $M_N(\mathbb C)$:
\begin{enumerate}
\item The invertible matrices are dense.

\item The matrices having distinct eigenvalues are dense.

\item The diagonalizable matrices are dense.
\end{enumerate}
\end{theorem}

\begin{proof}
These are quite advanced results, which can be proved as follows:

\medskip

(1) This is clear, intuitively speaking, because the invertible matrices are given by the condition $\det A\neq 0$. Thus, the set formed by these matrices appears as the complement of the hypersurface $\det A=0$, and so must be dense inside $M_N(\mathbb C)$, as claimed. 

\medskip

(2) Here we can use a similar argument, this time by saying that the set formed by the matrices having distinct eigenvalues appears as the complement of the hypersurface given by $\Delta(P_A)=0$, and so must be dense inside $M_N(\mathbb C)$, as claimed. 

\medskip

(3) This follows from (2), via the fact that the matrices having distinct eigenvalues are diagonalizable, that we know from Theorem 1.12. There are of course some other proofs as well, for instance by putting the matrix in Jordan form.
\end{proof}

As an application of the above results, and of our methods in general, we have:

\index{functional calculus}

\begin{theorem}
The following happen:
\begin{enumerate}
\item We have $P_{AB}=P_{BA}$, for any two matrices $A,B\in M_N(\mathbb C)$.

\item $AB,BA$ have the same eigenvalues, with the same multiplicities.

\item If $A$ has eigenvalues $\lambda_1,\ldots,\lambda_N$, then $f(A)$ has eigenvalues $f(\lambda_1),\ldots,f(\lambda_N)$.
\end{enumerate}
\end{theorem}

\begin{proof}
These results can be deduced by using Theorem 1.13, as follows:

\medskip

(1) It follows from definitions that the characteristic polynomial of a matrix is invariant under conjugation, in the sense that we have the following formula:
$$P_C=P_{ACA^{-1}}$$

Now observe that, when assuming that $A$ is invertible, we have:
$$AB=A(BA)A^{-1}$$

Thus, we have the result when $A$ is invertible. By using now Theorem 1.13 (1), we conclude that this formula holds for any matrix $A$, by continuity.  

\medskip

(2) This is a reformulation of (1), via the fact that $P$ encodes the eigenvalues, with multiplicities, which is hard to prove with bare hands.

\medskip

(3) This is something quite informal, clear for the diagonal matrices $D$, then for the diagonalizable matrices $PDP^{-1}$, and finally for all matrices, by using Theorem 1.13 (3), provided that $f$ has suitable regularity properties. We will be back to this.
\end{proof}

Let us go back to the main problem raised by the diagonalization procedure, namely the computation of the roots of characteristic polynomials. We have here:

\index{eigenvalue}

\begin{theorem}
The complex eigenvalues of a matrix $A\in M_N(\mathbb C)$, counted with multiplicities, have the following properties:
\begin{enumerate}
\item Their sum is the trace.

\item Their product is the determinant.
\end{enumerate}
\end{theorem}

\begin{proof}
Consider indeed the characteristic polynomial $P$ of the matrix:
\begin{eqnarray*}
P(X)
&=&\det(A-X1_N)\\
&=&(-1)^NX^N+(-1)^{N-1}Tr(A)X^{N-1}+\ldots+\det(A)
\end{eqnarray*}

We can factorize this polynomial, by using its $N$ complex roots, and we obtain:
\begin{eqnarray*}
P(X)
&=&(-1)^N(X-\lambda_1)\ldots(X-\lambda_N)\\
&=&(-1)^NX^N+(-1)^{N-1}\left(\sum_i\lambda_i\right)X^{N-1}+\ldots+\prod_i\lambda_i
\end{eqnarray*}

Thus, we are led to the conclusion in the statement.
\end{proof}

Regarding now the intermediate terms, we have here:

\index{characteristic polynomial}

\begin{theorem}
Assume that $A\in M_N(\mathbb C)$ has eigenvalues $\lambda_1,\ldots,\lambda_N\in\mathbb C$, counted with multiplicities. The basic symmetric functions of these eigenvalues, namely
$$c_k=\sum_{i_1<\ldots<i_k}\lambda_{i_1}\ldots\lambda_{i_k}$$
are then given by the fact that the characteristic polynomial of the matrix is:
$$P(X)=(-1)^N\sum_{k=0}^N(-1)^kc_kX^k$$
Moreover, all symmetric functions of the eigenvalues, such as the sums of powers
$$d_s=\lambda_1^s+\ldots+\lambda_N^s$$
appear as polynomials in these characteristic polynomial coefficients $c_k$.
\end{theorem}

\begin{proof}
These results can be proved by doing some algebra, as follows:

\medskip

(1) Consider indeed the characteristic polynomial $P$ of the matrix, factorized by using its $N$ complex roots, taken with multiplicities. By expanding, we obtain:
\begin{eqnarray*}
P(X)
&=&(-1)^N(X-\lambda_1)\ldots(X-\lambda_N)\\
&=&(-1)^NX^N+(-1)^{N-1}\left(\sum_i\lambda_i\right)X^{N-1}+\ldots+\prod_i\lambda_i\\
&=&(-1)^NX^N+(-1)^{N-1}c_1X^{N-1}+\ldots+(-1)^0c_N\\
&=&(-1)^N\left(X^N-c_1X^{N-1}+\ldots+(-1)^Nc_N\right)
\end{eqnarray*}

With the convention $c_0=1$, we are led to the conclusion in the statement.

\medskip

(2) This is something standard, coming by doing some abstract algebra. Working out the formulae for the sums of powers $d_s=\sum_i\lambda_i^s$, at small values of the exponent $s\in\mathbb N$, is an excellent exercise, which shows how to proceed in general, by recurrence.
\end{proof}

\section*{1d. Spectral theorems}

Let us go back now to the diagonalization question. Here is a key result:

\index{self-adjoint operator}

\begin{theorem}
Any matrix $A\in M_N(\mathbb C)$ which is self-adjoint, $A=A^*$, is diagonalizable, with the diagonalization being of the following type,
$$A=UDU^*$$
with $U\in U_N$, and with $D\in M_N(\mathbb R)$ diagonal. The converse holds too.
\end{theorem}

\begin{proof}
As a first remark, the converse trivially holds, because if we take a matrix of the form $A=UDU^*$, with $U$ unitary and $D$ diagonal and real, then we have:
\begin{eqnarray*}
A^*
&=&(UDU^*)^*\\
&=&UD^*U^*\\
&=&UDU^*\\
&=&A
\end{eqnarray*}

In the other sense now, assume that $A$ is self-adjoint, $A=A^*$.  Our first claim is that the eigenvalues are real. Indeed, assuming $Av=\lambda v$, we have:
\begin{eqnarray*}
\lambda<v,v>
&=&<\lambda v,v>\\
&=&<Av,v>\\
&=&<v,Av>\\
&=&<v,\lambda v>\\
&=&\bar{\lambda}<v,v>
\end{eqnarray*}

Thus we obtain $\lambda\in\mathbb R$, as claimed. Our next claim now is that the eigenspaces corresponding to different eigenvalues are pairwise orthogonal. Assume indeed that:
$$Av=\lambda v\quad,\quad 
Aw=\mu w$$

We have then the following computation, using $\lambda,\mu\in\mathbb R$:
\begin{eqnarray*}
\lambda<v,w>
&=&<\lambda v,w>\\
&=&<Av,w>\\
&=&<v,Aw>\\
&=&<v,\mu w>\\
&=&\mu<v,w>
\end{eqnarray*}

Thus $\lambda\neq\mu$ implies $v\perp w$, as claimed. In order now to finish the proof, it remains to prove that the eigenspaces of $A$ span the whole space $\mathbb C^N$. For this purpose, we will use a recurrence method. Let us pick an eigenvector of our matrix:
$$Av=\lambda v$$

Assuming now that we have a vector $w$ orthogonal to it, $v\perp w$, we have:
\begin{eqnarray*}
<Aw,v>
&=&<w,Av>\\
&=&<w,\lambda v>\\
&=&\lambda<w,v>\\
&=&0
\end{eqnarray*}

Thus, if $v$ is an eigenvector, then the vector space $v^\perp$ is invariant under $A$. Moreover, since a matrix $A$ is self-adjoint precisely when $<Av,v>\in\mathbb R$ for any vector $v\in\mathbb C^N$, as one can see by expanding the scalar product, the restriction of $A$ to the subspace $v^\perp$ is self-adjoint. Thus, we can proceed by recurrence, and we obtain the result.
\end{proof}

As basic examples of self-adjoint matrices, we have the orthogonal projections. The diagonalization result regarding them is as follows:

\index{projection}

\begin{proposition}
The matrices $P\in M_N(\mathbb C)$ which are projections,
$$P^2=P^*=P$$
are precisely those which diagonalize as follows,
$$P=UDU^*$$
with $U\in U_N$, and with $D\in M_N(0,1)$ being diagonal.
\end{proposition}

\begin{proof}
The equation for the projections being $P^2=P^*=P$, the eigenvalues $\lambda$ are real, and we have as well the following condition, coming from $P^2=P$:
\begin{eqnarray*}
\lambda<v,v>
&=&<\lambda v,v>\\
&=&<Pv,v>\\
&=&<P^2v,v>\\
&=&<Pv,Pv>\\
&=&<\lambda v,\lambda v>\\
&=&\lambda^2<v,v>
\end{eqnarray*}

Thus we obtain $\lambda\in\{0,1\}$, as claimed, and as a final conclusion here, the diagonalization of the self-adjoint matrices is as follows, with $e_i\in \{0,1\}$:
$$P\sim\begin{pmatrix}
e_1\\
&\ddots\\
&&e_N
\end{pmatrix}$$

To be more precise, the number of 1 values is the dimension of the image of $P$, and the number of 0 values is the dimension of space of vectors sent to 0 by $P$.
\end{proof}

An important class of self-adjoint matrices, which includes for instance all the projections, are the positive matrices. The theory here is as follows:

\index{positive operator}

\begin{theorem}
For a matrix $A\in M_N(\mathbb C)$ the following conditions are equivalent, and if they are satisfied, we say that $A$ is positive:
\begin{enumerate}
\item $A=B^2$, with $B=B^*$.

\item $A=CC^*$, for some $C\in M_N(\mathbb C)$.

\item $<Ax,x>\geq0$, for any vector $x\in\mathbb C^N$.

\item $A=A^*$, and the eigenvalues are positive, $\lambda_i\geq0$.

\item $A=UDU^*$, with $U\in U_N$ and with $D\in M_N(\mathbb R_+)$ diagonal.
\end{enumerate}
\end{theorem}

\begin{proof}
The idea is that the equivalences in the statement basically follow from some elementary computations, with only Theorem 1.17 needed, at some point:

\medskip

$(1)\implies(2)$ This is clear, because we can take $C=B$.

\medskip

$(2)\implies(3)$ This follows from the following computation:
\begin{eqnarray*}
<Ax,x>
&=&<CC^*x,x>\\
&=&<C^*x,C^*x>\\
&\geq&0
\end{eqnarray*}

$(3)\implies(4)$ By using the fact that $<Ax,x>$ is real, we have:
\begin{eqnarray*}
<Ax,x>
&=&<x,A^*x>\\
&=&<A^*x,x>
\end{eqnarray*}

Thus we have $A=A^*$, and the remaining assertion, regarding the eigenvalues, follows from the following computation, assuming $Ax=\lambda x$:
\begin{eqnarray*}
<Ax,x>
&=&<\lambda x,x>\\
&=&\lambda<x,x>\\
&\geq&0
\end{eqnarray*}

$(4)\implies(5)$ This follows indeed by using Theorem 1.17.

\medskip

$(5)\implies(1)$ Assuming $A=UDU^*$, with $U\in U_N$, and with $D\in M_N(\mathbb R_+)$ being diagonal, we can set $B=U\sqrt{D}U^*$. Then $B$ is self-adjoint, and its square is given by:
\begin{eqnarray*}
B^2
&=&U\sqrt{D}U^*\cdot U\sqrt{D}U^*\\
&=&UDU^*\\
&=&A
\end{eqnarray*}

Thus, we are led to the conclusion in the statement.
\end{proof}

Let us record as well the following technical version of the above result:

\index{strictly positive operator}

\begin{theorem}
For a matrix $A\in M_N(\mathbb C)$ the following conditions are equivalent, and if they are satisfied, we say that $A$ is strictly positive:
\begin{enumerate}
\item $A=B^2$, with $B=B^*$, invertible.

\item $A=CC^*$, for some $C\in M_N(\mathbb C)$ invertible.

\item $<Ax,x>>0$, for any nonzero vector $x\in\mathbb C^N$.

\item $A=A^*$, and the eigenvalues are strictly positive, $\lambda_i>0$.

\item $A=UDU^*$, with $U\in U_N$ and with $D\in M_N(\mathbb R_+^*)$ diagonal.
\end{enumerate}
\end{theorem}

\begin{proof}
This follows either from Theorem 1.19, by adding the various extra assumptions in the statement, or from the proof of Theorem 1.19, by modifying where needed.
\end{proof}

Let us discuss now the case of the unitary matrices. We have here:

\index{unitary}

\begin{theorem}
Any matrix $U\in M_N(\mathbb C)$ which is unitary, $U^*=U^{-1}$, is diagonalizable, with the eigenvalues on $\mathbb T$. More precisely we have
$$U=VDV^*$$
with $V\in U_N$, and with $D\in M_N(\mathbb T)$ diagonal. The converse holds too.
\end{theorem}

\begin{proof}
As a first remark, the converse trivially holds, because given a matrix of type $U=VDV^*$, with $V\in U_N$, and with $D\in M_N(\mathbb T)$ being diagonal, we have:
\begin{eqnarray*}
U^*
&=&(VDV^*)^*\\
&=&VD^*V^*\\
&=&VD^{-1}V^{-1}\\
&=&(V^*)^{-1}D^{-1}V^{-1}\\
&=&(VDV^*)^{-1}\\
&=&U^{-1}
\end{eqnarray*}

Let us prove now the first assertion, stating that the eigenvalues of a unitary matrix $U\in U_N$ belong to $\mathbb T$. Indeed, assuming $Uv=\lambda v$, we have:
\begin{eqnarray*}
<v,v>
&=&<U^*Uv,v>\\
&=&<Uv,Uv>\\
&=&<\lambda v,\lambda v>\\
&=&|\lambda|^2<v,v>
\end{eqnarray*}

Thus we obtain $\lambda\in\mathbb T$, as claimed. Our next claim now is that the eigenspaces corresponding to different eigenvalues are pairwise orthogonal. Assume indeed that:
$$Uv=\lambda v\quad,\quad 
Uw=\mu w$$

We have then the following computation, using $U^*=U^{-1}$ and $\lambda,\mu\in\mathbb T$:
\begin{eqnarray*}
\lambda<v,w>
&=&<\lambda v,w>\\
&=&<Uv,w>\\
&=&<v,U^*w>\\
&=&<v,U^{-1}w>\\
&=&<v,\mu^{-1}w>\\
&=&\mu<v,w>
\end{eqnarray*}

Thus $\lambda\neq\mu$ implies $v\perp w$, as claimed. In order now to finish the proof, it remains to prove that the eigenspaces of $U$ span the whole space $\mathbb C^N$. For this purpose, we will use a recurrence method. Let us pick an eigenvector of our matrix:
$$Uv=\lambda v$$

Assuming that we have a vector $w$ orthogonal to it, $v\perp w$, we have:
\begin{eqnarray*}
<Uw,v>
&=&<w,U^*v>\\
&=&<w,U^{-1}v>\\
&=&<w,\lambda^{-1}v>\\
&=&\lambda<w,v>\\
&=&0
\end{eqnarray*}

Thus, if $v$ is an eigenvector, then the vector space $v^\perp$ is invariant under $U$. Now since $U$ is an isometry, so is its restriction to this space $v^\perp$. Thus this restriction is a unitary, and so we can proceed by recurrence, and we obtain the result.
\end{proof}

The self-adjoint matrices and the unitary matrices are particular cases of the general notion of a ``normal matrix'', and we have here:

\index{normal operator}

\begin{theorem}
Any matrix $A\in M_N(\mathbb C)$ which is normal, $AA^*=A^*A$, is diagonalizable, with the diagonalization being of the following type,
$$A=UDU^*$$
with $U\in U_N$, and with $D\in M_N(\mathbb C)$ diagonal. The converse holds too.
\end{theorem}

\begin{proof}
As a first remark, the converse trivially holds, because if we take a matrix of the form $A=UDU^*$, with $U$ unitary and $D$ diagonal, then we have:
\begin{eqnarray*}
AA^*
&=&UDU^*\cdot UD^*U^*\\
&=&UDD^*U^*\\
&=&UD^*DU^*\\
&=&UD^*U^*\cdot UDU^*\\
&=&A^*A
\end{eqnarray*}

In the other sense now, this is something more technical. Our first claim is that a matrix $A$ is normal precisely when the following happens, for any vector $v$:
$$||Av||=||A^*v||$$

Indeed, the above equality can be written as follows:
$$<AA^*v,v>=<A^*Av,v>$$

But this is equivalent to $AA^*=A^*A$, by expanding the scalar products. Our next claim is that $A,A^*$ have the same eigenvectors, with conjugate eigenvalues:
$$Av=\lambda v\implies A^*v=\bar{\lambda}v$$

Indeed, this follows from the following computation, and from the trivial fact that if $A$ is normal, then so is any matrix of type $A-\lambda 1_N$:
\begin{eqnarray*}
||(A^*-\bar{\lambda}1_N)v||
&=&||(A-\lambda 1_N)^*v||\\
&=&||(A-\lambda 1_N)v||\\
&=&0
\end{eqnarray*}

Let us prove now, by using this, that the eigenspaces of $A$ are pairwise orthogonal. Assume that we have two eigenvectors, corresponding to different eigenvalues, $\lambda\neq\mu$:
$$Av=\lambda v\quad,\quad 
Aw=\mu w$$

We have the following computation, which shows that $\lambda\neq\mu$ implies $v\perp w$:
\begin{eqnarray*}
\lambda<v,w>
&=&<\lambda v,w>\\
&=&<Av,w>\\
&=&<v,A^*w>\\
&=&<v,\bar{\mu}w>\\
&=&\mu<v,w>
\end{eqnarray*}

In order to finish, it remains to prove that the eigenspaces of $A$ span the whole $\mathbb C^N$. This is something that we have already seen for the self-adjoint matrices, and for unitaries, and we will use here these results, in order to deal with the general normal case. As a first observation, given an arbitrary matrix $A$, the matrix $AA^*$ is self-adjoint:
$$(AA^*)^*=AA^*$$

Thus, we can diagonalize this matrix $AA^*$, as follows, with the passage matrix being a unitary, $V\in U_N$, and with the diagonal form being real, $E\in M_N(\mathbb R)$:
$$AA^*=VEV^*$$

Now observe that, for matrices of type $A=UDU^*$, which are those that we supposed to deal with, we have the following formulae:
$$V=U\quad,\quad 
E=D\bar{D}$$

In particular, the matrices $A$ and $AA^*$ have the same eigenspaces. So, this will be our idea, proving that the eigenspaces of $AA^*$ are eigenspaces of $A$. In order to do so, let us pick two eigenvectors $v,w$ of the matrix $AA^*$, corresponding to different eigenvalues, $\lambda\neq\mu$. The eigenvalue equations are then as follows:
$$AA^*v=\lambda v\quad,\quad 
AA^*w=\mu w$$

We have the following computation, using the normality condition $AA^*=A^*A$, and the fact that the eigenvalues of $AA^*$, and in particular $\mu$, are real:
\begin{eqnarray*}
\lambda<Av,w>
&=&<\lambda Av,w>\\
&=&<A\lambda v,w>\\
&=&<AAA^*v,w>\\
&=&<AA^*Av,w>\\
&=&<Av,AA^*w>\\
&=&<Av,\mu w>\\
&=&\mu<Av,w>
\end{eqnarray*}

We conclude that we have $<Av,w>=0$. But this reformulates as follows:
$$\lambda\neq\mu\implies A(E_\lambda)\perp E_\mu$$

Now since the eigenspaces of $AA^*$ are pairwise orthogonal, and span the whole $\mathbb C^N$, we deduce from this that these eigenspaces are invariant under $A$:
$$A(E_\lambda)\subset E_\lambda$$

But with this result in hand, we can finish. Indeed, we can decompose the problem, and the matrix $A$ itself, following these eigenspaces of $AA^*$, which in practice amounts in saying that we can assume that we only have 1 eigenspace. Now by rescaling, this is the same as assuming that we have $AA^*=1$. But with this, we are now into the unitary case, that we know how to solve, as explained in Theorem 1.21, and so done.
\end{proof}

As a first application, we have the following result:

\index{absolute value}
\index{modulus of operator}

\begin{theorem}
Given a matrix $A\in M_N(\mathbb C)$, we can construct a matrix $|A|$ as follows, by using the fact that $A^*A$ is diagonalizable, with positive eigenvalues:
$$|A|=\sqrt{A^*A}$$
This matrix $|A|$ is then positive, and its square is $|A|^2=A^*A$. In the case $N=1$, we obtain in this way the usual absolute value of the complex numbers.
\end{theorem}

\begin{proof}
Consider indeed the matrix $A^*A$, which is normal. According to Theorem 1.22, we can diagonalize this matrix as follows, with $U\in U_N$, and with $D$ diagonal:
$$A=UDU^*$$

From $A^*A\geq0$ we obtain $D\geq0$. But this means that the entries of $D$ are real, and positive. Thus we can extract the square root $\sqrt{D}$, and then set:
$$\sqrt{A^*A}=U\sqrt{D}U^*$$

Thus, we are basically done. Indeed, if we call this latter matrix $|A|$, then we are led to the conclusions in the statement. Finally, the last assertion is clear from definitions.
\end{proof}

We can now formulate a first polar decomposition result, as follows:

\index{polar decomposition}
\index{partial isometry}

\begin{theorem}
Any invertible matrix $A\in M_N(\mathbb C)$ decomposes as
$$A=U|A|$$
with $U\in U_N$, and with $|A|=\sqrt{A^*A}$ as above.
\end{theorem}

\begin{proof}
This is routine, and follows by comparing the actions of $A,|A|$ on the vectors $v\in\mathbb C^N$, and deducing from this the existence of a unitary $U\in U_N$ as above. We will be back to this, later on, directly in the case of the linear operators on Hilbert spaces.
\end{proof}

Observe that at $N=1$ we obtain in this way the usual polar decomposition of the nonzero complex numbers. More generally now, we have the following result:

\begin{theorem}
Any square matrix $A\in M_N(\mathbb C)$ decomposes as
$$A=U|A|$$
with $U$ being a partial isometry, and with $|A|=\sqrt{A^*A}$ as above.
\end{theorem}

\begin{proof}
Again, this follows by comparing the actions of $A,|A|$ on the vectors $v\in\mathbb C^N$, and deducing from this the existence of a partial isometry $U$ as above. Alternatively, we can get this from Theorem 1.24, applied on the complement of the 0-eigenvectors. 
\end{proof}

This was for our basic presentation of linear algebra. There are of course many other things that can be said, but we will come back to some of them in what follows, directly in the case of the linear operators on the arbitrary Hilbert spaces.

\section*{1e. Exercises}

Linear algebra is a wide topic, and there are countless interesting matrices, and exercises about them. As a continuation of our discussion about rotations, we have:

\begin{exercise}
Prove that the symmetry and projection with respect to the $Ox$ axis rotated by an angle $t/2\in\mathbb R$ are given by the matrices
$$S_t=\begin{pmatrix}\cos t&\sin t\\ \sin t&-\cos t\end{pmatrix}$$
$$P_t=\frac{1}{2}\begin{pmatrix}1+\cos t&\sin t\\ \sin t&1-\cos t\end{pmatrix}$$
and then diagonalize these matrices, and if possible without computations.
\end{exercise}

Here the first part can only be clear on pictures, and by the way, prior to this, do not forget to verify as well that our formula of $R_t$ is the good one. As for the second part, just don't go head-first into computations, there might be some geometry over there.

\begin{exercise}
Prove that the isometries of $\mathbb R^2$ are rotations or symmetries,
$$R_t=\begin{pmatrix}\cos t&-\sin t\\ \sin t&\cos t\end{pmatrix}\quad,\quad 
S_t=\begin{pmatrix}\cos t&\sin t\\ \sin t&-\cos t\end{pmatrix}$$
and then try as well to find a formula for the isometries of $\mathbb R^3$.
\end{exercise}

Here for the first question you should look first at the determinant of such an isometry. As for the second question, this is something quite difficult. If you're good at computers, you can look into the code of 3D games, the rotation formula is probably there.

\begin{exercise}
Prove that the isometries of $\mathbb C^2$ of determinant $1$ are
$$U=\begin{pmatrix}a&b\\ -\bar{b}&\bar{a}\end{pmatrix}\quad,\quad |a|^2+|b|^2=1$$
then work out as well the general case, of arbitrary determinant. 
\end{exercise}

As a comment here, if done with this exercise about $\mathbb C^2$, but not yet with the previous one about $\mathbb R^3$, you can go back to that exercise, by using a $\mathbb C^2\simeq\mathbb R^4$ trick. And in case this trick leads to tough computations and big headache, look it up.

\begin{exercise}
Prove that the flat matrix, which is the all-one $N\times N$ matrix, diagonalizes over the complex numbers as follows,
$$\begin{pmatrix}
1&\ldots&\ldots&1\\
\vdots&&&\vdots\\
\vdots&&&\vdots\\
1&\ldots&\ldots&1\end{pmatrix}=\frac{1}{N}\,F_N
\begin{pmatrix}
N\\
&0\\
&&\ddots\\
&&&0\end{pmatrix}F_N^*$$
where $F_N=(w^{ij})_{ij}$ with $w=e^{2\pi i/N}$ is the Fourier matrix, with the convention that the indices are taken to be $i,j=0,1,\ldots,N-1$.
\end{exercise}

This is something very instructive. Normally you have to look for eigenvectors for the flat matrix, and you are led in this way to the equation $x_0+\ldots+x_{N-1}=0$. The problem however is that this equation, while looking very gentle, has no ``canonical'' solutions over the real numbers. Thus you are led to the complex numbers, and more specifically to the roots of unity, and their magic, leading to the above result. Enjoy.

\chapter{Linear operators}

\section*{2a. Hilbert spaces}

We discuss in what follows an extension of the linear algebra results from the previous chapter, obtained by looking at the linear operators $T:H\to H$, with the space $H$ being no longer assumed to be finite dimensional. Our motivations come from quantum mechanics, and in order to get motivated, here is some suggested reading:

\bigskip

(1) Generally speaking, physics is best learned from Feynman \cite{fey}. If you already know some, and want to learn quantum mechanics, go with Griffiths \cite{gri}. And if you're already a bit familiar with quantum mechanics, a good book is Weinberg \cite{wei}.

\bigskip

(2) A look at classics like Dirac \cite{dir}, von Neumann \cite{vn4} or Weyl \cite{wey} can be instructive too. On the opposite, you have as well modern, fancy books on quantum information, such as Bengtsson-\.Zyczkowski \cite{bzy}, Nielsen-Chuang \cite{nch} or Watrous \cite{wat}.

\bigskip

(3) In short, many ways of getting familiar with this big mess which is quantum mechanics, and as long as you stay away from books advertised as ``rigorous'', ``axiomatic'', ``mathematical'', things fine. By the way, you can try as well my book \cite{ba4}.

\bigskip

Getting to work now, physics tells us to look at infinite dimensional complex spaces, such as the space of wave functions $\psi:\mathbb R^3\to\mathbb C$ of the electron. In order to do some mathematics on these spaces, we will need scalar products. So, let us start with:

\index{scalar product}

\begin{definition}
A scalar product on a complex vector space $H$ is a binary operation $H\times H\to\mathbb C$, denoted $(x,y)\to <x,y>$, satisfying the following conditions:
\begin{enumerate}
\item $<x,y>$ is linear in $x$, and antilinear in $y$.

\item $\overline{<x,y>}=<y,x>$, for any $x,y$.

\item $<x,x>>0$, for any $x\neq0$.
\end{enumerate}
\end{definition}

As before in chapter 1, we use here mathematicians' convention for scalar products, that is, $<\,,>$ linear at left, as opposed to physicists' convention, $<\,,>$ linear at right. The reasons for this are quite subtle, coming from the fact that, while basic quantum mechanics looks better with $<\,,>$ linear at right, advanced quantum mechanics looks better with $<\,,>$ linear at left. Or at least that's what my cats say.

\bigskip

As a basic example for Definition 2.1, we have the finite dimensional vector space $H=\mathbb C^N$, with its usual scalar product, namely:
$$<x,y>=\sum_ix_i\bar{y}_i$$ 

There are many other examples, and notably various spaces of $L^2$ functions, which naturally appear in problems coming from physics. We will discuss them later on. In order to study now the scalar products, let us formulate the following definition:

\index{norm of vector}

\begin{definition}
The norm of a vector $x\in H$ is the following quantity:
$$||x||=\sqrt{<x,x>}$$
We also call this number length of $x$, or distance from $x$ to the origin.
\end{definition}

The terminology comes from what happens in $\mathbb C^N$, where the length of the vector, as defined above, coincides with the usual length, given by:
$$||x||=\sqrt{\sum_i|x_i|^2}$$

In analogy with what happens in finite dimensions, we have two important results regarding the norms. First we have the Cauchy-Schwarz inequality, as follows:

\index{Cauchy-Schwarz}

\begin{theorem}
We have the Cauchy-Schwarz inequality
$$|<x,y>|\leq||x||\cdot||y||$$
and the equality case holds precisely when $x,y$ are proportional.
\end{theorem}

\begin{proof}
This is something very standard. Consider indeed the following quantity, depending on a real variable $t\in\mathbb R$, and on a variable on the unit circle, $w\in\mathbb T$:
$$f(t)=||twx+y||^2$$

By developing $f$, we see that this is a degree 2 polynomial in $t$:
\begin{eqnarray*}
f(t)
&=&<twx+y,twx+y>\\
&=&t^2<x,x>+tw<x,y>+t\bar{w}<y,x>+<y,y>\\
&=&t^2||x||^2+2tRe(w<x,y>)+||y||^2
\end{eqnarray*}

Since $f$ is obviously positive, its discriminant must be negative:
$$4Re(w<x,y>)^2-4||x||^2\cdot||y||^2\leq0$$

But this is equivalent to the following condition:
$$|Re(w<x,y>)|\leq||x||\cdot||y||$$

Now the point is that we can arrange for the number $w\in\mathbb T$ to be such that the quantity $w<x,y>$ is real. Thus, we obtain the following inequality:
$$|<x,y>|\leq||x||\cdot||y||$$

Finally, the study of the equality case is straightforward, by using the fact that the discriminant of $f$ vanishes precisely when we have a root. But this leads to the conclusion in the statement, namely that the vectors $x,y$ must be proportional. 
\end{proof}
 
As a second main result now, we have the Minkowski inequality:

\index{Minkowski inequality}

\begin{theorem}
We have the Minkowski inequality
$$||x+y||\leq||x||+||y||$$
and the equality case holds precisely when $x,y$ are proportional.
\end{theorem}

\begin{proof}
This follows indeed from the Cauchy-Schwarz inequality, as follows:
\begin{eqnarray*}
&&||x+y||\leq||x||+||y||\\
&\iff&||x+y||^2\leq(||x||+||y||)^2\\
&\iff&||x||^2+||y||^2+2Re<x,y>\leq||x||^2+||y||^2+2||x||\cdot||y||\\
&\iff&Re<x,y>\leq||x||\cdot||y||
\end{eqnarray*}

As for the equality case, this is clear from Cauchy-Schwarz as well.
\end{proof}
  
As a consequence of this, we have the following result:

\begin{theorem}
The following function is a distance on $H$,
$$d(x,y)=||x-y||$$
in the usual sense, that of the abstract metric spaces.
\end{theorem}

\begin{proof}
This follows indeed from the Minkowski inequality, which corresponds to the triangle inequality, the other two axioms for a distance being trivially satisfied.
\end{proof}

The above result is quite important, because it shows that we can do geometry and analysis in our present setting, with distances and angles, a bit as in the finite dimensional case. In order to do such abstract geometry, we will often need the following key result, which shows that everything can be recovered in terms of distances:

\index{polarization identity}

\begin{proposition}
The scalar products can be recovered from distances, via the formula
$$4<x,y>
=||x+y||^2-||x-y||^2
+i||x+iy||^2-i||x-iy||^2$$
called complex polarization identity.
\end{proposition}

\begin{proof}
This is something that we have already met in finite dimensions. In arbitrary dimensions the proof is similar, as follows:
\begin{eqnarray*}
&&||x+y||^2-||x-y||^2+i||x+iy||^2-i||x-iy||^2\\
&=&||x||^2+||y||^2-||x||^2-||y||^2+i||x||^2+i||y||^2-i||x||^2-i||y||^2\\
&&+2Re(<x,y>)+2Re(<x,y>)+2iIm(<x,y>)+2iIm(<x,y>)\\
&=&4<x,y>
\end{eqnarray*}

Thus, we are led to the conclusion in the statement.
\end{proof}

In order to do analysis on our spaces, we need the Cauchy sequences that we construct to converge. This is something which is automatic in finite dimensions, but in arbitrary dimensions, this can fail. It is convenient here to formulate a detailed new definition, as follows, which will be the starting point for our various considerations to follow:

\index{Hilbert space}

\begin{definition}
A Hilbert space is a complex vector space $H$ given with a scalar product $<x,y>$, satisfying the following conditions:
\begin{enumerate}
\item $<x,y>$ is linear in $x$, and antilinear in $y$.

\item $\overline{<x,y>}=<y,x>$, for any $x,y$.

\item $<x,x>>0$, for any $x\neq0$.

\item $H$ is complete with respect to the norm $||x||=\sqrt{<x,x>}$.
\end{enumerate}
\end{definition}

In other words, we have taken here Definition 2.1 above, and added the condition that $H$ must be complete with respect to the norm $||x||=\sqrt{<x,x>}$, that we know indeed to be a norm, according to the Minkowski inequality proved above. As a basic example, as before, we have the space $H=\mathbb C^N$, with its usual scalar product, namely:
$$<x,y>=\sum_ix_i\bar{y}_i$$

More generally now, we have the following construction of Hilbert spaces:

\begin{proposition}
The sequences of complex numbers $(x_i)$ which are square-summable, 
$$\sum_i|x_i|^2<\infty$$
form a Hilbert space $l^2(\mathbb N)$, with the following scalar product:
$$<x,y>=\sum_ix_i\bar{y}_i$$
In fact, given any index set $I$, we can construct a Hilbert space $l^2(I)$, in this way.
\end{proposition}

\begin{proof}
There are several things to be proved, as follows:

\medskip

(1) Our first claim is that $l^2(\mathbb N)$ is a vector space. For this purpose, we must prove that $x,y\in l^2(\mathbb N)$ implies $x+y\in l^2(\mathbb N)$. But this leads us into proving $||x+y||\leq||x||+||y||$, where $||x||=\sqrt{<x,x>}$. Now since we know this inequality to hold on each subspace $\mathbb C^N\subset l^2(\mathbb N)$ obtained by truncating, this inequality holds everywhere, as desired.

\medskip

(2) Our second claim is that $<\,,>$ is well-defined on $l^2(\mathbb N)$. But this follows from the Cauchy-Schwarz inequality, $|<x,y>|\leq||x||\cdot||y||$, which can be established by truncating, a bit like we established the Minkowski inequality in (1) above.

\medskip

(3) It is also clear that $<\,,>$ is a scalar product on $l^2(\mathbb N)$, so it remains to prove that $l^2(\mathbb N)$ is complete with respect to $||x||=\sqrt{<x,x>}$. But this is clear, because if we pick a Cauchy sequence $\{x^n\}_{n\in\mathbb N}\subset l^2(\mathbb N)$, then each numeric sequence $\{x^n_i\}_{i\in\mathbb N}\subset\mathbb C$ is Cauchy, and by setting $x_i=\lim_{n\to\infty}x^n_i$, we have $x^n\to x$ inside $l^2(\mathbb N)$, as desired.

\medskip

(4) Finally, the same arguments extend to the case of an arbitrary index set $I$, leading to a Hilbert space $l^2(I)$, and with the remark here that there is absolutely no problem of taking about quantities of type $||x||^2=\sum_{i\in I}|x_i|^2\in[0,\infty]$, even if the index set $I$ is uncountable, because we are summing positive numbers.
\end{proof}

Even more generally, we have the following construction of Hilbert spaces:

\begin{theorem}
Given a measured space $X$, the functions $f:X\to\mathbb C$, taken up to equality almost everywhere, which are square-summable, 
$$\int_X|f(x)|^2dx<\infty$$
form a Hilbert space $L^2(X)$, with the following scalar product:
$$<f,g>=\int_Xf(x)\overline{g(x)}dx$$
In the case $X=I$, with the counting measure, we obtain in this way the space $l^2(I)$.
\end{theorem}

\begin{proof}
This is a straightforward generalization of Proposition 2.8, with the arguments from the proof of Proposition 2.8 carrying over in our case, as follows:

\medskip

(1) The first part, regarding Cauchy-Schwarz and Minkowski, extends without problems, by using this time approximation by step functions.

\medskip

(2) Regarding the fact that $<\,,>$ is indeed a scalar product on $L^2(X)$, there is a subtlety here, because if we want $<f,f>>0$ for $f\neq 0$, we must declare that $f=0$ when $f=0$ almost everywhere, and so that $f=g$ when $f=g$ almost everywhere.

\medskip

(3) Regarding the fact that $L^2(X)$ is complete with respect to $||f||=\sqrt{<f,f>}$, this is again basic measure theory, by picking a Cauchy sequence $\{f_n\}_{n\in\mathbb N}\subset L^2(X)$, and then constructing a pointwise, and hence $L^2$ limit, $f_n\to f$, almost everywhere. 

\medskip

(4) Finally, the last assertion is clear, because the integration with respect to the counting measure is by definition a sum, and so $L^2(I)=l^2(I)$ in this case.
\end{proof}

Quite remarkably, any Hilbert space must be of the form $L^2(X)$, and even of the particular form $l^2(I)$. This follows indeed from the following key result:

\index{orthogonal basis}
\index{Gram-Schmidt}

\begin{theorem}
Let $H$ be a Hilbert space.
\begin{enumerate}
\item Any algebraic basis of this space $\{f_i\}_{i\in I}$ can be turned into an orthonormal basis $\{e_i\}_{i\in I}$, by using the Gram-Schmidt procedure.

\item Thus, $H$ has an orthonormal basis, and so we have $H\simeq l^2(I)$, with $I$ being the indexing set for this orthonormal basis.
\end{enumerate}
\end{theorem}

\begin{proof}
All this is standard by Gram-Schmidt, the idea being as follows:

\medskip

(1) First of all, in finite dimensions an orthonormal basis $\{e_i\}_{i\in I}$ is by definition a usual algebraic basis, satisfying $<e_i,e_j>=\delta_{ij}$. But the existence of such a basis follows by applying the Gram-Schmidt procedure to any algebraic basis $\{f_i\}_{i\in I}$, as claimed.

\medskip

(2) In infinite dimensions, a first issue comes from the fact that the standard basis $\{\delta_i\}_{i\in\mathbb N}$ of the space $l^2(\mathbb N)$ is not an algebraic basis in the usual sense, with the finite linear combinations of the functions $\delta_i$ producing only a dense subspace of $l^2(\mathbb N)$, that of the functions having finite support. Thus, we must fine-tune our definition of ``basis''.

\medskip

(3) But this can be done in two ways, by saying that $\{f_i\}_{i\in I}$ is a basis of $H$ when  the functions $f_i$ are linearly independent, and when either the finite linear combinations of these functions $f_i$ form a dense subspace of $H$, or the linear combinations with $l^2(I)$ coefficients of these functions $f_i$ form the whole $H$. For orthogonal bases $\{e_i\}_{i\in I}$ these definitions are equivalent, and in any case, our statement makes now sense.

\medskip

(4) Regarding now the proof, in infinite dimensions, this follows again from Gram-Schmidt, exactly as in the finite dimensional case, but by using this time a tool from logic, called Zorn lemma, in order to correctly do the recurrence.
\end{proof}

The above result, and its relation with Theorem 2.9, is something quite subtle, so let us further get into this. First, we have the following definition, based on the above:

\begin{definition}
A Hilbert space $H$ is called separable when the following equivalent conditions are satisfied:
\begin{enumerate}
\item $H$ has a countable algebraic basis $\{f_i\}_{i\in\mathbb N}$.

\item $H$ has a countable orthonormal basis $\{e_i\}_{i\in\mathbb N}$.

\item We have $H\simeq l^2(\mathbb N)$, isomorphism of Hilbert spaces.
\end{enumerate}
\end{definition}

In what follows we will be mainly interested in the separable Hilbert spaces, where most of the questions coming from quantum physics take place. In view of the above, the following philosophical question appears: why not simply talking about $l^2(\mathbb N)$?

\bigskip

In answer to this, we cannot really do so, because many of the separable spaces that we are interested in appear as spaces of functions, and such spaces do not necessarily have a very simple or explicit orthonormal basis, as shown by the following result: 

\index{orthogonal polynomials}

\begin{proposition}
The Hilbert space $H=L^2[0,1]$ is separable, having as orthonormal basis the orthonormalized version of the algebraic basis $f_n=x^n$ with $n\in\mathbb N$.
\end{proposition}
 
\begin{proof}
This follows from the Weierstrass theorem, which provides us with the basis $f_n=x^n$, which can be orthogonalized by using the Gram-Schmidt procedure, as explained in Theorem 2.10. Working out the details here is actually an excellent exercise.
\end{proof}

As a conclusion to all this, we are interested in 1 space, namely the unique separable Hilbert space $H$, but due to various technical reasons, it is often better to forget that we have $H=l^2(\mathbb N)$, and say instead that we have $H=L^2(X)$, with $X$ being a separable measured space, or simply say that $H$ is an abstract separable Hilbert space.

\section*{2b. Linear operators}

Let us get now into the study of linear operators $T:H\to H$. Before anything, we should mention that things are quite tricky with respect to quantum mechanics, and physics in general. Indeed, if there is a central operator in physics, this is the Laplace operator on the smooth functions $f:\mathbb R^N\to\mathbb C$, given by:
$$\Delta f(x)=\sum_i\frac{d^2f}{dx_i^2}$$

And the problem is that what we have here is an operator $\Delta:C^\infty(\mathbb R^N)\to C^\infty(\mathbb R^N)$, which does not extend into an operator $\Delta:L^2(\mathbb R^N)\to L^2(\mathbb R^N)$. Thus, we should perhaps look at operators $T:H\to H$ which are densely defined, instead of looking at operators $T:H\to H$ which are everywhere defined. We will not do so, for two reasons:

\bigskip

(1) Tactical retreat. When physics looks too complicated, as it is the case now, you can always declare that mathematics comes first. So, let us be pure mathematicians, simply looking in generalizing linear algebra to infinite dimensions. And from this viewpoint, it is a no-brainer to look at everywhere defined operators $T:H\to H$.

\bigskip

(2) Modern physics. We will see later, towards the middle of the present book, when talking about various mathematical physics findings of Connes, Jones, Voiculescu and others, that a lot of interesting mathematics, which is definitely related to modern physics, can be developed by using the everywhere defined operators $T:H\to H$.

\bigskip

In short, you'll have to trust me here. And hang on, we are not done yet, because with this choice made, there is one more problem, mathematical this time. The problem comes from the fact that in infinite dimensions the everywhere defined operators $T:H\to H$ can be bounded or not, and for reasons which are mathematically intuitive and obvious, and physically acceptable too, we want to deal with the bounded case only. 

\bigskip

Long story short, let us avoid too much thinking, and start in a simple way, with:

\index{operator norm}
\index{bounded operator}

\begin{proposition}
For a linear operator $T:H\to H$, the following are equivalent:
\begin{enumerate}
\item $T$ is continuous.

\item $T$ is continuous at $0$.

\item $T(B)\subset cB$ for some $c<\infty$, where $B\subset H$ is the unit ball.

\item $T$ is bounded, in the sense that $||T||=\sup_{||x||\leq1}||Tx||$ satisfies $||T||<\infty$.
\end{enumerate}
\end{proposition}

\begin{proof}
This is elementary, with $(1)\iff(2)$ coming from the linearity of $T$, then $(2)\iff(3)$ coming from definitions, and finally $(3)\iff(4)$ coming from the fact that the number $||T||$ from (4) is the infimum of the numbers $c$ making (3) work.
\end{proof}

Regarding such operators, we have the following result:

\index{operator algebra}

\begin{theorem}
The linear operators $T:H\to H$ which are bounded,
$$||T||=\sup_{||x||\leq1}||Tx||<\infty$$
form a complex algebra with unit $B(H)$, having the property 
$$||ST||\leq||S||\cdot||T||$$
and which is complete with respect to the norm.
\end{theorem}

\begin{proof}
The fact that we have indeed an algebra, satisfying the product condition in the statement, follows from the following estimates, which are all elementary:
$$||S+T||\leq||S||+||T||$$
$$||\lambda T||=|\lambda|\cdot||T||$$
$$||ST||\leq||S||\cdot||T||$$

Regarding now the last assertion, if $\{T_n\}\subset B(H)$ is Cauchy then $\{T_nx\}$ is Cauchy for any $x\in H$, so we can define the limit $T=\lim_{n\to\infty}T_n$ by setting:
$$Tx=\lim_{n\to\infty}T_nx$$

Let us first check that the application $x\to Tx$ is linear. We have:
\begin{eqnarray*}
T(x+y)
&=&\lim_{n\to\infty}T_n(x+y)\\
&=&\lim_{n\to\infty}T_n(x)+T_n(y)\\
&=&\lim_{n\to\infty}T_n(x)+\lim_{n\to\infty}T_n(y)\\
&=&T(x)+T(y)
\end{eqnarray*}

Similarly, we have as well the following computation:
\begin{eqnarray*}
T(\lambda x)
&=&\lim_{n\to\infty}T_n(\lambda x)\\
&=&\lambda\lim_{n\to\infty}T_n(x)\\
&=&\lambda T(x)
\end{eqnarray*}

Thus we have a linear map $T:A\to A$. It remains to prove that we have $T\in B(H)$, and that we have $T_n\to T$ in norm. For this purpose, observe that we have:
\begin{eqnarray*}
&&||T_n-T_m||\leq\varepsilon\ ,\ \forall n,m\geq N\\
&\implies&||T_nx-T_mx||\leq\varepsilon\ ,\ \forall||x||=1\ ,\ \forall n,m\geq N\\
&\implies&||T_nx-Tx||\leq\varepsilon\ ,\ \forall||x||=1\ ,\ \forall n\geq N\\
&\implies&||T_Nx-Tx||\leq\varepsilon\ ,\ \forall||x||=1\\
&\implies&||T_N-T||\leq\varepsilon
\end{eqnarray*}

As a first consequence, we obtain $T\in B(H)$, because we have:
\begin{eqnarray*}
||T||
&=&||T_N+(T-T_N)||\\
&\leq&||T_N||+||T-T_N||\\
&\leq&||T_N||+\varepsilon\\
&<&\infty
\end{eqnarray*}

As a second consequence, we obtain $T_N\to T$ in norm, and we are done.
\end{proof}

In the case where $H$ comes with a basis $\{e_i\}_{i\in I}$, we can talk about the infinite matrices $M\in M_I(\mathbb C)$, with the remark that the multiplication of such matrices is not always defined, in the case $|I|=\infty$. In this context, we have the following result:

\index{linear operator}

\begin{theorem}
Let $H$ be a Hilbert space, with orthonormal basis $\{e_i\}_{i\in I}$. The bounded operators $T\in B(H)$ can be then identified with matrices $M\in M_I(\mathbb C)$ via
$$Tx=Mx\quad,\quad M_{ij}=<Te_j,e_i>$$
and we obtain in this way an embedding as follows, which is multiplicative:
$$B(H)\subset M_I(\mathbb C)$$
In the case $H=\mathbb C^N$ we obtain in this way the usual isomorphism $B(H)\simeq M_N(\mathbb C)$. In the separable case we obtain in this way a proper embedding $B(H)\subset M_\infty(\mathbb C)$.
\end{theorem}

\begin{proof}
We have several assertions to be proved, the idea being as follows:

\medskip

(1) Regarding the first assertion, given a bounded operator $T:H\to H$, let us associate to it a matrix $M\in M_I(\mathbb C)$ as in the statement, by the following formula:
$$M_{ij}=<Te_j,e_i>$$

It is clear that this correspondence $T\to M$ is linear, and also that its kernel is $\{0\}$. Thus, we have an embedding of linear spaces $B(H)\subset M_I(\mathbb C)$.

(2) Our claim now is that this embedding is multiplicative. But this is clear too, because if we denote by $T\to M_T$ our correspondence, we have:
\begin{eqnarray*}
(M_{ST})_{ij}
&=&<STe_j,e_i>\\
&=&\left<S\sum_k<Te_j,e_k>e_k,e_i\right>\\
&=&\sum_k<Se_k,e_i><Te_j,e_k>\\
&=&\sum_k(M_S)_{ik}(M_T)_{kj}\\
&=&(M_SM_T)_{ij}
\end{eqnarray*}

(3) Finally, we must prove that the original operator $T:H\to H$ can be recovered from its matrix $M\in M_I(\mathbb C)$ via the formula in the statement, namely $Tx=Mx$. But this latter formula holds for the vectors of the basis, $x=e_j$, because we have:
\begin{eqnarray*}
(Te_j)_i
&=&<Te_j,e_i>\\
&=&M_{ij}\\
&=&(Me_j)_i
\end{eqnarray*}

Now by linearity we obtain from this that the formula $Tx=Mx$ holds everywhere, on any vector $x\in H$, and this finishes the proof of the first assertion.

\medskip

(4) In finite dimensions we obtain an isomorphism, because any matrix $M\in M_N(\mathbb C)$ determines an operator $T:\mathbb C^N\to\mathbb C^N$, according to the formula $<Te_j,e_i>=M_{ij}$. In infinite dimensions, however, we do not have an isomorphism. For instance on $H=l^2(\mathbb N)$ the following matrix does not define an operator:
$$M=\begin{pmatrix}1&1&\ldots\\
1&1&\ldots\\
\vdots&\vdots
\end{pmatrix}$$

Indeed, $T(e_1)$ should be the all-one vector, which is not square-summable.
\end{proof}

In connection with our previous comments on bases, the above result is something quite theoretical, because for basic Hilbert spaces like $L^2[0,1]$, which do not have a simple orthonormal basis, the embedding $B(H)\subset M_\infty(\mathbb C)$ that we obtain is not something very useful. In short, while the bounded operators $T:H\to H$ are basically some infinite matrices, it is better to think of these operators as being objects on their own.

\bigskip

As another comment, the construction $T\to M$ makes sense for any linear operator $T:H\to H$, but when $\dim H=\infty$, we do not obtain an embedding $\mathcal L(H)\subset M_I(\mathbb C)$ in this way. Indeed, set $H=l^2(\mathbb N)$, let $E=span(e_i)$ be the linear space spanned by the standard basis, and pick an algebraic complement $F$ of this space $E$, so that we have $H=E\oplus F$, as an algebraic direct sum. Then any linear operator $S:F\to F$ gives rise to a linear operator $T:H\to H$, given by $T(e,f)=(0,S(f))$, whose associated matrix is $0$. And, restrospectively speaking, it is in order to avoid such pathologies that we decided some time ago to restrict the attention to the bounded case, $T\in B(H)$.

\bigskip

As in the finite dimensional case, we can talk about adjoint operators, in this setting, the definition and main properties of the construction $T\to T^*$ being as follows:

\index{adjoint operator}

\begin{theorem}
Given a bounded operator $T\in B(H)$, the following formula defines a bounded operator $T^*\in B(H)$, called adjoint of $H$:
$$<Tx,y>=<x,T^*y>$$
The correspondence $T\to T^*$ is antilinear, antimultiplicative, and is an involution, and an isometry. In finite dimensions, we recover the usual adjoint operator.
\end{theorem}

\begin{proof}
There are several things to be done here, the idea being as follows:

\medskip

(1) We will need a standard functional analysis result, stating that the continuous linear forms $\varphi:H\to\mathbb C$ appear as scalar products, as follows, with $z\in H$:
$$\varphi(x)=<x,z>$$

Indeed, in one sense this is clear, because given $z\in H$, the application $\varphi(x)=<x,z>$ is linear, and continuous as well, because by Cauchy-Schwarz we have:
$$|\varphi(x)|\leq||x||\cdot||z||$$

Conversely now, by using a basis we can assume $H=l^2(\mathbb N)$, and our linear form $\varphi:H\to\mathbb C$ must be then, by linearity, given by a formula of the following type:
$$\varphi(x)=\sum_ix_i\bar{z}_i$$

But, again by Cauchy-Schwarz, in order for such a formula to define indeed a continuous linear form $\varphi:H\to\mathbb C$ we must have $z\in l^2(\mathbb N)$, and so $z\in H$, as desired.

\medskip

(2) With this in hand, we can now construct the adjoint $T^*$, by the formula in the statement. Indeed, given $y\in H$, the formula $\varphi(x)=<Tx,y>$ defines a linear map $H\to\mathbb C$. Thus, we must have a formula as follows, for a certain vector $T^*y\in H$:
$$\varphi(x)=<x,T^*y>$$

Moreover, this vector $T^*y\in H$ is unique with this property, and we conclude from this that the formula $y\to T^*y$ defines a certain map $T^*:H\to H$, which is unique with the property in the statement, namely $<Tx,y>=<x,T^*y>$ for any $x,y$.

\medskip

(3) Let us prove that we have $T^*\in B(H)$. By using once again the uniqueness of $T^*$, we conclude that we have the following formulae, which show that $T^*$ is linear:
$$T^*(x+y)=T^*x+T^*y\quad,\quad 
T^*(\lambda x)=\lambda T^*x$$

Observe also that $T^*$ is bounded as well, because we have:
\begin{eqnarray*}
||T||
&=&\sup_{||x||=1}\sup_{||y||=1}<Tx,y>\\
&=&\sup_{||y||=1}\sup_{||x||=1}<x,T^*y>\\
&=&||T^*||
\end{eqnarray*}

(4) The fact that the correspondence $T\to T^*$ is antilinear, antimultiplicative, and is an involution comes from the following formulae, coming from uniqueness:
$$(S+T)^*=S^*+T^*\quad,\quad 
(\lambda T)^*=\bar{\lambda}T^*$$
$$(ST)^*=T^*S^*\quad,\quad 
(T^*)^*=T$$

As for the isometry property with respect to the operator norm, $||T||=||T^*||$, this is something that we already know, from the proof of (3) above.

\medskip

(5) Regarding finite dimensions, let us first examine the general case where our Hilbert space comes with a basis, $H=l^2(I)$. We can compute the matrix $M^*\in M_I(\mathbb C)$ associated to the operator $T^*\in B(H)$, by using $<Tx,y>=<x,T^*y>$, in the following way:
\begin{eqnarray*}
(M^*)_{ij}
&=&<T^*e_j,e_i>\\
&=&\overline{<e_i,T^*e_j>}\\
&=&\overline{<Te_i,e_j>}\\
&=&\overline{M}_{ji}
\end{eqnarray*}

Thus, we have reached to the usual formula for the adjoints of matrices, and in the particular case $H=\mathbb C^N$, we conclude that $T^*$ comes indeed from the usual $M^*$.
\end{proof}

As in finite dimensions, the operators $T,T^*$ can be thought of as being ``twin brothers'', and there is a lot of interesting mathematics connecting them. We first have:

\index{adjoint operator}

\begin{proposition}
Given a bounded operator $T\in B(H)$, the following happen:
\begin{enumerate}
\item $\ker T^*=(Im T)^\perp$.

\item $\overline{Im T^*}=(\ker T)^\perp$.
\end{enumerate}
\end{proposition}

\begin{proof}
Both these assertions are elementary, as follows:

\medskip

(1) Let us first prove ``$\subset$''. Assuming $T^*x=0$, we have indeed $x\perp ImT$, because:
$$<x,Ty>
=<T^*x,y>
=0$$

As for ``$\supset$'', assuming $<x,Ty>=0$ for any $y$, we have $T^*x=0$, because:
$$<T^*x,y>
=<x,Ty>
=0$$

(2) This can be deduced from (1), applied to the operator $T^*$, as follows:
$$(\ker T)^\perp
=(Im T^*)^{\perp\perp}
=\overline{Im T^*}$$

Here we have used the formula $K^{\perp\perp}=\bar{K}$, valid for any linear subspace $K\subset H$ of a Hilbert space, which for $K$ closed reads $K^{\perp\perp}=K$, and comes from $H=K\oplus K^\perp$, and which in general follows from $K^{\perp\perp}\subset\bar{K}^{\perp\perp}=\bar{K}$, the reverse inclusion being clear.
\end{proof}

Let us record as well the following useful formula, relating $T$ and $T^*$:

\index{norm of operators}

\begin{theorem}
We have the following formula,
$$||TT^*||=||T||^2$$
valid for any operator $T\in B(H)$.
\end{theorem}

\begin{proof}
We recall from Theorem 2.16 that the correspondence $T\to T^*$ is an isometry with respect to the operator norm, in the sense that we have:
$$||T||=||T^*||$$

In order to prove now the formula in the statement, observe first that we have:
$$||TT^*||
\leq||T||\cdot||T^*||
=||T||^2$$

On the other hand, we have as well the following estimate:
\begin{eqnarray*}
||T||^2
&=&\sup_{||x||=1}|<Tx,Tx>|\\
&=&\sup_{||x||=1}|<x,T^*Tx>|\\
&\leq&||T^*T||
\end{eqnarray*}

By replacing $T\to T^*$ we obtain from this that we have:
$$||T||^2\leq||TT^*||$$

Thus, we have obtained the needed inequality, and we are done.
\end{proof}

\section*{2c. Unitaries, projections}

Let us discuss now some explicit examples of operators, in analogy with what happens in finite dimensions. The most basic examples of linear transformations are the rotations, symmetries and projections. Then, we have certain remarkable classes of linear transformations, such as the positive, self-adjoint and normal ones. In what follows we will develop the basic theory of such transformations, in the present Hilbert space setting.

\bigskip

Let us begin with the rotations. The situation here is quite tricky in arbitrary dimensions, and we have several notions instead of one. We first have the following result:

\index{isometry}

\begin{theorem}
For a linear operator $U\in B(H)$ the following conditions are equivalent, and if they are satisfied, we say that $U$ is an isometry:
\begin{enumerate}
\item $U$ is a metric space isometry, $d(Ux,Uy)=d(x,y)$.

\item $U$ is a normed space isometry, $||Ux||=||x||$.

\item $U$ preserves the scalar product, $<Ux,Uy>=<x,y>$.

\item $U$ satisfies the isometry condition $U^*U=1$.
\end{enumerate}
In finite dimensions, we recover in this way the usual unitary transformations.
\end{theorem}

\begin{proof}
The proofs are similar to those in finite dimensions, as follows:

\medskip

$(1)\iff(2)$ This follows indeed from the formula of the distances, namely:
$$d(x,y)=||x-y||$$

$(2)\iff(3)$ This is again standard, because we can pass from scalar products to distances, and vice versa, by using $||x||=\sqrt{<x,x>}$, and the polarization formula.

\medskip

$(3)\iff(4)$ We have indeed the following equivalences, by using the standard formula $<Tx,y>=<x,T^*y>$, which defines the adjoint operator:
\begin{eqnarray*}
<Ux,Uy>=<x,y>
&\iff&<x,U^*Uy>=<x,y>\\
&\iff&U^*Uy=y\\
&\iff&U^*U=1
\end{eqnarray*}

Thus, we are led to the conclusions in the statement.
\end{proof}

The point now is that the condition $U^*U=1$ does not imply in general $UU^*=1$, the simplest counterexample here being the shift operator on $l^2(\mathbb N)$:

\index{shift}

\begin{proposition}
The shift operator on the space $l^2(\mathbb N)$, given by
$$S(e_i)=e_{i+1}$$
is an isometry, $S^*S=1$. However, we have $SS^*\neq1$.
\end{proposition}

\begin{proof}
The adjoint of the shift is given by the following formula:
$$S^*(e_i)=\begin{cases}
e_{i-1}&{\rm if}\ i>0\\
0&{\rm if}\ i=0
\end{cases}$$

When composing $S,S^*$, in one sense we obtain the following formula:
$$S^*S(e_i)=e_i$$

In other other sense now, we obtain the following formula:
$$SS^*(e_i)=\begin{cases}
e_i&{\rm if}\ i>0\\
0&{\rm if}\ i=0
\end{cases}$$

Summarizing, the compositions are given by the following formulae:
$$S^*S=1\quad,\quad 
SS^*=Proj(e_0^\perp)$$

Thus, we are led to the conclusions in the statement.
\end{proof}

As a conclusion, the notion of isometry is not the correct infinite dimensional analogue of the notion of unitary, and the unitary operators must be introduced as follows:

\index{unitary}
\index{isometry}

\begin{theorem}
For a linear operator $U\in B(H)$ the following conditions are equivalent, and if they are satisfied, we say that $U$ is a unitary:
\begin{enumerate}
\item $U$ is an isometry, which is invertible.

\item $U$, $U^{-1}$ are both isometries.

\item $U$, $U^*$ are both isometries.

\item $UU^*=U^*U=1$.

\item $U^*=U^{-1}$.
\end{enumerate}
Moreover, the unitary operators from a group $U(H)\subset B(H)$.
\end{theorem}

\begin{proof}
There are several statements here, the idea being as follows:

\medskip

(1) The various equivalences in the statement are all clear from definitions, and from Theorem 2.19 in what regards the various possible notions of isometries which can be used, by using the formula $(ST)^*=T^*S^*$ for the adjoints of the products of operators.

\medskip

(2) The fact that the products and inverses of unitaries are unitaries is also clear, and we conclude that the unitary operators from a group $U(H)\subset B(H)$, as stated.
\end{proof}

Let us discuss now the projections. Modulo the fact that all the subspaces $K\subset H$ where these projections project must be assumed to be closed, in the present setting, here the result is perfectly similar to the one in finite dimensions, as follows:

\index{projection}

\begin{theorem}
For a linear operator $P\in B(H)$ the following conditions are equivalent, and if they are satisfied, we say that $P$ is a projection:
\begin{enumerate}
\item $P$ is the orthogonal projection on a closed subspace $K\subset H$.

\item $P$ satisfies the projection equations $P^2=P^*=P$.
\end{enumerate}
\end{theorem}

\begin{proof}
As in finite dimensions, $P$ is an abstract projection, not necessarily orthogonal, when it is an idempotent, algebrically speaking, in the sense that we have:
$$P^2=P$$

The point now is that this projection is orthogonal when:
\begin{eqnarray*}
<Px-x,Py>=0
&\iff&<P^*Px-P^*x,y>=0\\
&\iff&P^*Px-P^*x=0\\
&\iff&P^*P-P^*=0\\
&\iff&P^*P=P^*
\end{eqnarray*}

Now observe that by conjugating, we obtain $P^*P=P$. Thus, we must have $P=P^*$, and so we have shown that any orthogonal projection must satisfy, as claimed:
$$P^2=P^*=P$$

Conversely, if this condition is satisfied, $P^2=P$ shows that $P$ is a projection, and $P=P^*$ shows via the above computation that $P$ is indeed orthogonal.
\end{proof}

There is a relation between the projections and the general isometries, such as the shift $S$ that we met before, and we have the following result:

\begin{proposition}
Given an isometry $U\in B(H)$, the operator 
$$P=UU^*$$
is a projection, namely the orthogonal projection on $Im(U)$.
\end{proposition}

\begin{proof}
Assume indeed that we have an isometry, $U^*U=1$. The fact that $P=UU^*$ is indeed a projection can be checked abstractly, as follows:
$$(UU^*)^*=UU^*$$
$$UU^*UU^*=UU^*$$

As for the last assertion, this is something that we already met, for the shift, and the situation in general is similar, with the result itself being clear.
\end{proof}

More generally now, along the same lines, and clarifying the whole situation with the unitaries and isometries, we have the following result:

\index{partial isometry}

\begin{theorem}
An operator $U\in B(H)$ is a partial isometry, in the usual geometric sense, when the following two operators are projections:
$$P=UU^*\quad,\quad 
Q=U^*U$$
Moreover, the isometries, adjoints of isometries and unitaries are respectively characterized by the conditions $Q=1$, $P=1$, $P=Q=1$.
\end{theorem}

\begin{proof}
The first assertion is a straightforward extension of Proposition 2.23, and the second assertion follows from various results regarding isometries established above.
\end{proof}

It is possible to talk as well about symmetries, in the following way:

\index{symmetry}

\begin{definition}
An operator $S\in B(H)$ is called a symmetry when $S^2=1$, and a unitary symmetry when one of the following equivalent conditions is satisfied:
\begin{enumerate}
\item $S$ is a unitary, $S^*=S^{-1}$, and a symmetry as well, $S^2=1$.

\item $S$ satisfies the equations $S=S^*=S^{-1}$.
\end{enumerate}
\end{definition}

Here the terminology is a bit non-standard, because even in finite dimensions, $S^2=1$ is not exactly what you would require for a ``true'' symmetry, as shown by the following transformation, which is a symmetry in our sense, but not a unitary symmetry:
$$\begin{pmatrix}0&2\\ 1/2&0\end{pmatrix}\binom{x}{y}=\binom{2y}{x/2}$$ 

Let us study now some larger classes of operators, which are of particular importance, namely the self-adjoint, positive and normal ones. We first have:

\index{self-adjoint operator}

\begin{theorem}
For an operator $T\in B(H)$, the following conditions are equivalent, and if they are satisfied, we call $T$ self-adjoint:
\begin{enumerate}
\item $T=T^*$.

\item $<Tx,x>\in\mathbb R$.
\end{enumerate}
In finite dimensions, we recover in this way the usual self-adjointness notion.
\end{theorem}

\begin{proof}
There are several assertions here, the idea being as follows:

\medskip

$(1)\implies(2)$ This is clear, because we have:
\begin{eqnarray*}
\overline{<Tx,x>}
&=&<x,Tx>\\
&=&<T^*x,x>\\
&=&<Tx,x>
\end{eqnarray*}

$(2)\implies(1)$ In order to prove this, observe that the beginning of the above computation shows that, when assuming $<Tx,x>\in\mathbb R$, the following happens:
$$<Tx,x>=<T^*x,x>$$

Thus, in terms of the operator $S=T-T^*$, we have:
$$<Sx,x>=0$$

In order to finish, we use a polarization trick. We have the following formula:
$$<S(x+y),x+y>=<Sx,x>+<Sy,y>+<Sx,y>+<Sy,x>$$

Since the first 3 terms vanish, the sum of the 2 last terms vanishes too. But, by using $S^*=-S$, coming from $S=T-T^*$, we can process this latter vanishing as follows:
\begin{eqnarray*}
<Sx,y>
&=&-<Sy,x>\\
&=&<y,Sx>\\
&=&\overline{<Sx,y>}
\end{eqnarray*}

Thus we must have $<Sx,y>\in\mathbb R$, and with $y\to iy$ we obtain $<Sx,y>\in i\mathbb R$ too, and so $<Sx,y>=0$. Thus $S=0$, which gives $T=T^*$, as desired.

\medskip

(3) Finally, in what regards the finite dimensions, or more generally the case where our Hilbert space comes with a basis, $H=l^2(I)$, here the condition $T=T^*$ corresponds to the usual self-adjointness condition $M=M^*$ at the level of the associated matrices.
\end{proof}

At the level of the basic examples, the situation is as follows:

\index{unitary symmetry}

\begin{proposition}
The folowing operators are self-adjoint:
\begin{enumerate}
\item The projections, $P^2=P^*=P$. In fact, an abstract, algebraic projection is an orthogonal projection precisely when it is self-adjoint.

\item The unitary symmetries, $S=S^*=S^{-1}$. In fact, a unitary is a unitary symmetry precisely when it is self-adjoint.
\end{enumerate}
\end{proposition}

\begin{proof}
These assertions are indeed all clear from definitions.
\end{proof}

Next in line, we have the notion of positive operator. We have here:

\index{positive operator}

\begin{theorem}
The positive operators, which are the operators $T\in B(H)$ satisfying $<Tx,x>\geq0$, have the following properties: 
\begin{enumerate}
\item They are self-adjoint, $T=T^*$.

\item As examples, we have the projections, $P^2=P^*=P$.

\item More generally, $T=S^*S$ is positive, for any $S\in B(H)$.

\item In finite dimensions, we recover the usual positive operators.
\end{enumerate}
\end{theorem}

\begin{proof}
All these assertions are elementary, the idea being as follows:

\medskip

(1) This follows from Theorem 2.26, because $<Tx,x>\geq0$ implies $<Tx,x>\in\mathbb R$.

\medskip

(2) This is clear from $P^2=P=P^*$, because we have:
\begin{eqnarray*}
<Px,x>
&=&<P^2x,x>\\
&=&<Px,Px>\\
&=&||Px||^2
\end{eqnarray*}

(3) This follows from a similar computation, namely:
$$<S^*Sx,x>
=<Sx,Sx>
=||Sx||^2$$

(4) This is well-known, the idea being that the condition $<Tx,x>\geq0$ corresponds to the usual positivity condition $A\geq0$, at the level of the associated matrix.
\end{proof}

It is possible to talk as well about strictly positive operators, and we have here:

\begin{theorem}
The strictly positive operators, which are the operators $T\in B(H)$ satisfying $<Tx,x>>0$, for any $x\neq0$, have the following properties: 
\begin{enumerate}
\item They are self-adjoint, $T=T^*$.

\item As examples, $T=S^*S$ is positive, for any $S\in B(H)$ injective.

\item In finite dimensions, we recover the usual strictly positive operators.
\end{enumerate}
\end{theorem}

\begin{proof}
As before, all these assertions are elementary, the idea being as follows:

\medskip

(1) This is something that we know, from Theorem 2.28.

\medskip

(2) This follows from the injectivity of $S$, because for any $x\neq0$ we have:
\begin{eqnarray*}
<S^*Sx,x>
&=&<Sx,Sx>\\
&=&||Sx||^2\\
&>&0
\end{eqnarray*}

(3) This is well-known, the idea being that the condition $<Tx,x>>0$ corresponds to the usual strict positivity condition $A>0$, at the level of the associated matrix.
\end{proof}

As a comment, while any strictly positive matrix $A>0$ is well-known to be invertible, the analogue of this fact does not hold in infinite dimensions, a counterexample here being the following operator on $l^2(\mathbb N)$, which is strictly positive but not invertible:
$$T=\begin{pmatrix}
1\\
&\frac{1}{2}\\
&&\frac{1}{3}\\
&&&\ddots
\end{pmatrix}$$

As a last remarkable class of operators, we have the normal ones. We have here:

\index{normal operator}

\begin{theorem}
For an operator $T\in B(H)$, the following conditions are equivalent, and if they are satisfied, we call $T$ normal:
\begin{enumerate}
\item $TT^*=T^*T$.

\item $||Tx||=||T^*x||$.
\end{enumerate}
In finite dimensions, we recover in this way the usual normality notion.
\end{theorem}

\begin{proof}
There are several assertions here, the idea being as follows:

\medskip

$(1)\implies(2)$ This is clear, due to the following computation:
\begin{eqnarray*}
||Tx||^2
&=&<Tx,Tx>\\
&=&<T^*Tx,x>\\
&=&<TT^*x,x>\\
&=&<T^*x,T^*x>\\
&=&||T^*x||^2
\end{eqnarray*}

$(2)\implies(1)$ This is clear as well, because the above computation shows that, when assuming $||Tx||=||T^*x||$, the following happens:
$$<TT^*x,x>=<T^*Tx,x>$$

Thus, in terms of the operator $S=TT^*-T^*T$, we have:
$$<Sx,x>=0$$

In order to finish, we use a polarization trick. We have the following formula:
$$<S(x+y),x+y>=<Sx,x>+<Sy,y>+<Sx,y>+<Sy,x>$$

Since the first 3 terms vanish, the sum of the 2 last terms vanishes too. But, by using $S=S^*$, coming from $S=TT^*-T^*T$, we can process this latter vanishing as follows:
\begin{eqnarray*}
<Sx,y>
&=&-<Sy,x>\\
&=&-<y,Sx>\\
&=&-\overline{<Sx,y>}
\end{eqnarray*}

Thus we must have $<Sx,y>\in i\mathbb R$, and with $y\to iy$ we obtain $<Sx,y>\in \mathbb R$ too, and so $<Sx,y>=0$. Thus $S=0$, which gives $TT^*=T^*T$, as desired.

\medskip

(3) Finally, in what regards finite dimensions, or more generally the case where our Hilbert space comes with a basis, $H=l^2(I)$, here the condition $TT^*=T^*T$ corresponds to the usual normality condition $MM^*=M^*M$ at the level of the associated matrices.
\end{proof}

Observe that the normal operators generalize both the self-adjoint operators, and the unitaries. We will be back to such operators, on many occassions, in what follows.

\section*{2d. Diagonal operators}

Let us work out now what happens in the case that we are mostly interested in, namely $H=L^2(X)$, with $X$ being a measured space. We first have:

\index{multiplication operator}

\begin{theorem}
Given a measured space $X$, consider the Hilbert space $H=L^2(X)$. Associated to any function $f\in L^\infty(X)$ is then the multiplication operator
$$T_f:H\to H\quad,\quad 
T_f(g)=fg$$
which is well-defined, linear and bounded, having norm as follows:
$$||T_f||=||f||_\infty$$
Moreover, the correspondence $f\to T_f$ is linear, multiplicative and involutive.
\end{theorem}

\begin{proof}
There are several assertions here, the idea being as follows:

\medskip

(1) We must first prove that the formula in the statement, $T_f(g)=fg$, defines indeed an operator $H\to H$, which amounts in saying that we have:
$$f\in L^\infty(X),\ g\in L^2(X)\implies fg\in L^2(X)$$

But this follows from the following explicit estimate:
\begin{eqnarray*}
||fg||_2
&=&\sqrt{\int_X|f(x)|^2|g(x)|^2d\mu(x)}\\
&\leq&\sup_{x\in X}|f(x)|^2\sqrt{\int_X|g(x)|^2d\mu(x)}\\
&=&||f||_\infty||g||_2\\
&<&\infty
\end{eqnarray*}

(2) Next in line, we must prove that $T$ is linear and bounded. We have:
$$T_f(g+h)=T_f(g)+T_f(h)\quad,\quad 
T_f(\lambda g)=\lambda T_f(g)$$

As for the boundedness condition, this follows from the estimate from the proof of (1), which gives, in terms of the operator norm of $B(H)$:
$$||T_f||\leq||f||_\infty$$

(3) Let us prove now that we have equality, $||T_f||=||f||_\infty$, in the above estimate. For this purpose, we use the well-known fact that the $L^\infty$ functions can be approximated by $L^2$ functions. Indeed, with such an approximation $g_n\to f$ we obtain:
\begin{eqnarray*}
||fg_n||_2
&=&\sqrt{\int_X|f(x)|^2|g_n(x)|^2d\mu(x)}\\
&\simeq&\sup_{x\in X}|f(x)|^2\sqrt{\int_X|g_n(x)|^2d\mu(x)}\\
&=&||f||_\infty||g_n||_2
\end{eqnarray*}

Thus, with $n\to\infty$ we obtain $||T_f||\geq||f||_\infty$, which is reverse to the inequality obtained in the proof of (2), and this leads to the conclusion in the statement.

\medskip

(4) Regarding now the fact that the correspondence $f\to T_f$ is indeed linear and multiplicative, the corresponding formulae are as follows, both clear:
$$T_{f+h}(g)=T_f(g)+T_h(g)\quad,\quad 
T_{\lambda f}(g)=\lambda T_f(g)$$

(5) Finally, let us prove that the correspondence $f\to T_f$ is involutive, in the sense that it transforms the standard involution $f\to\bar{f}$ of the algebra $L^\infty(X)$ into the standard involution $T\to T^*$ of the algebra $B(H)$. We must prove that we have:
$$T_f^*=T_{\bar{f}}$$

But this follows from the following computation:
\begin{eqnarray*}
<T_fg,h>
&=&<fg,h>\\
&=&\int_Xf(x)g(x)\bar{h}(x)d\mu(x)\\
&=&\int_Xg(x)f(x)\bar{h}(x)d\mu(x)\\
&=&<g,\bar{f}h>\\
&=&<g,T_{\bar{f}}h>
\end{eqnarray*}

Indeed, since the adjoint is unique, we obtain from this $T_f^*=T_{\bar{f}}$. Thus the correspondence $f\to T_f$ is indeed involutive, as claimed.
\end{proof}

In what regards now the basic classes of operators, the above construction provides us with many new examples, which are very explicit, and are complementary to the usual finite dimensional examples that we usually have in mind, as follows:

\index{multiplication operator}
\index{unitary}
\index{symmetry}
\index{projection}
\index{isometry}
\index{positive operator}
\index{self-adjoint operator}
\index{normal operator}

\begin{theorem}
The multiplication operators $T_f(g)=fg$ on the Hilbert space $H=L^2(X)$ associated to the functions $f\in L^\infty(X)$ are as follows:
\begin{enumerate}
\item $T_f$ is unitary when $f:X\to\mathbb T$.

\item $T_f$ is a symmetry when $f:X\to\{-1,1\}$.

\item $T_f$ is a projection when $f=\chi_Y$ with $Y\in X$.

\item There are no non-unitary isometries.

\item There are no non-unitary symmetries.

\item $T_f$ is positive when $f:X\to\mathbb R_+$.

\item $T_f$ is self-adjoint when $f:X\to\mathbb R$.

\item $T_f$ is always normal, for any $f:X\to\mathbb C$.
\end{enumerate}
\end{theorem}

\begin{proof}
All these assertions are clear from definitions, and from the various properties of the correspondence $f\to T_f$, established above, as follows:

\medskip

(1) The unitarity condition $U^*=U^{-1}$ for the operator $T_f$ reads $\bar{f}=f^{-1}$, which means that we must have $f:X\to\mathbb T$, as claimed.

\medskip

(2) The symmetry condition $S^2=1$ for the operator $T_f$ reads $f^2=1$, which means that we must have $f:X\to\{-1,1\}$, as claimed.

\medskip

(3) The projection condition $P^2=P^*=P$ for the operator $T_f$ reads $f^2=f=\bar{f}$, which means that we must have $f:X\to\{0,1\}$, or equivalently, $f=\chi_Y$ with $Y\subset X$.

\medskip

(4) A non-unitary isometry must satisfy by definition $U^*U=1,UU^*\neq1$, and for the operator $T_f$ this means that we must have $|f|^2=1,|f|^2\neq1$, which is impossible.

\medskip

(5) This follows from (1) and (2), because the solutions found in (2) for the symmetry problem are included in the solutions found in (1) for the unitarity problem.

\medskip

(6) The fact that $T_f$ is positive amounts in saying that we must have $<fg,g>\geq0$ for any $g\in L^2(X)$, and this is equivalent to the fact that we must have $f\geq0$, as desired.

\medskip

(7) The self-adjointness condition $T=T^*$ for the operator $T_f$ reads $f=\bar{f}$, which means that we must have $f:X\to\mathbb R$, as claimed.

\medskip

(8) The normality condition $TT^*=T^*T$ for the operator $T_f$ reads $f\bar{f}=\bar{f}f$, which is automatic for any function $f:X\to\mathbb C$, as claimed.
\end{proof}

The above result might look quite puzzling, at a first glance, messing up our intuition with various classes of operators, coming from usual linear algebra. However, a bit of further thinking tells us that there is no contradiction, and that Theorem 2.32 in fact is very similar to what we know about the diagonal matrices. To be more precise, the diagonal matrices are unitaries precisely when their entries are in $\mathbb T$, there are no non-unitary isometries, all such matrices are normal, and so on. In order to understand all this, let us work out what happens with the correspondence $f\to T_f$, in finite dimensions. The situation here is in fact extremely simple, and illuminating, as follows:

\index{multiplication operator}
\index{diagonal operator}

\begin{theorem}
Assuming $X=\{1,\ldots,N\}$ with the counting measure, the embedding
$$L^\infty(X)\subset B(L^2(X))$$
constructed via multiplication operators, $T_f(g)=fg$, corresponds to the embedding
$$\mathbb C^N\subset M_N(\mathbb C)$$
given by the diagonal matrices, constructed as follows:
$$f\to diag(f_1,\ldots,f_N)$$
Thus, Theorem 2.32 generalizes what we know about the diagonal matrices.
\end{theorem} 

\begin{proof}
The idea is that all this is trivial, with not a single new computation needed, modulo some algebraic thinking, of quite soft type. Let us go back indeed to Theorem 2.31 above and its proof, with the abstract measured space $X$ appearing there being now the following finite space, with its counting mesure:
$$X=\{1,\ldots,N\}$$

Regarding the functions $f\in L^\infty(X)$, these are now functions as follows:
$$f:\{1,\ldots,N\}\to\mathbb C$$

We can identify such a function with the corresponding vector $(f(i))_i\in\mathbb C^N$, and so we conclude that our input algebra $L^\infty(X)$ is the algebra $\mathbb C^N$:
$$L^\infty(X)=\mathbb C^N$$

Regarding now the Hilbert space $H=L^2(X)$, this is equal as well to $\mathbb C^N$, and for the same reasons, namely that $g\in L^2(X)$ can be identified with the vector $(g(i))_i\in\mathbb C^N$:
$$L^2(X)=\mathbb C^N$$

Observe that, due to our assumption that $X$ comes with its counting measure, the scalar product that we obtain on $\mathbb C^N$ is the usual one, without weights. Now, let us identify the operators on $L^2(X)=\mathbb C^N$ with the square matrices, in the usual way:
$$B(L^2(X))=M_N(\mathbb C)$$

This was our final identification, in order to get started. Now by getting back to Theorem 2.31, the embedding $L^\infty(X)\subset B(L^2(X))$ constructed there reads:
$$\mathbb C^N\subset M_N(\mathbb C)$$

But this can only be the embedding given by the diagonal matrices, so are basically done. In order to finish, however, let us understand what the operator associated to an arbitrary vector $f\in\mathbb C^N$ is. We can regard this vector as a function, $f(i)=f_i$, and so the action $T_f(g)=fg$ on the vectors of $L^2(X)=\mathbb C^N$ is by componentwise multiplication by the numbers $f_1,\ldots,f_N$. But this is exactly the action of the diagonal matrix $diag(f_1,\ldots,f_N)$, and so we are led to the conclusion in the statement.
\end{proof}

There are other things that can be said about the embedding $L^\infty(X)\subset B(L^2(X))$, a key observation here, which is elementary to prove, being the fact that the image of $L^\infty(X)$ is closed with respect to the weak topology, the one where $T_n\to T$ when $T_nx\to Tx$ for any $x\in H$. And with this meaning that $L^\infty(X)$ is a so-called von Neumann algebra on $L^2(X)$. We will be back to this, on numerous occasions, in what follows.

\section*{2e. Exercises} 

As before with linear algebra, operator theory is a wide area of mathematics, and there are many interesting operators, and exercises about them. We first have:

\begin{exercise}
Find an explicit orthonormal basis for the Hilbert space 
$$H=L^2[0,1]$$
by starting with the algebraic basic $f_n=x^n$ with $n\in\mathbb N$, and applying Gram-Schmidt.
\end{exercise}

This is actually quite non-trivial, and in case you're stuck with complicated computations, better look it up, preferably in the physics literature, physicists being well-known to adore such things, and then write a brief account of what you found.

\begin{exercise}
Find all the $2\times2$ complex matrices 
$$S=\begin{pmatrix}a&b\\ c&d\end{pmatrix}$$
which are symmetries, $S^2=1$, and interpret them geometrically.
\end{exercise}

Here you can of course start with the real case first, $S\in M_2(\mathbb R)$. Also, you can have a look at 3 dimensions too, real or complex, and beware of the computations here.

\begin{exercise}
Prove that any positive operator $T\geq0$ appears as
$$T=S^2$$
with $S$ self-adjoint, first in finite dimensions, then in general.
\end{exercise}

Here the discussion in finite dimensions involves positive eigenvalues and their square roots, which is something quite standard. In infinite dimensions things are a bit more complicated, because we don't have yet such eigenvalue technology, and with this being actually to come in the next chapter, but you can try of course some other tricks. 

\chapter{Spectral theorems}

\section*{3a. Basic theory}

We discuss in this chapter the diagonalization problem for the operators $T\in B(H)$, in analogy with the diagonalization problem for the usual matrices $A\in M_N(\mathbb C)$. As a first observation, we can talk about eigenvalues and eigenvectors, as follows:

\index{eigenvector}
\index{eigenvalue}

\begin{definition}
Given an operator $T\in B(H)$, assuming that we have
$$Tx=\lambda x$$
we say that $x\in H$ is an eigenvector of $T$, with eigenvalue $\lambda\in\mathbb C$.
\end{definition}

We know many things about eigenvalues and eigenvectors, in the finite dimensional case. However, most of these will not extend to the infinite dimensional case, or at least not extend in a straightforward way, due to a number of reasons:

\medskip

\begin{enumerate}
\item Most of basic linear algebra is based on the fact that $Tx=\lambda x$ is equivalent to $(T-\lambda)x=0$, so that $\lambda$ is an eigenvalue when $T-\lambda$ is not invertible. In the infinite dimensional setting $T-\lambda$ might be injective and not surjective, or vice versa, or invertible with $(T-\lambda)^{-1}$ not bounded, and so on.

\medskip

\item Also, in linear algebra $T-\lambda$ is not invertible when $\det(T-\lambda)=0$, and with this leading to most of the advanced results about eigenvalues and eigenvectors. In infinite dimensions, however, it is impossible to construct a determinant function $\det:B(H)\to\mathbb C$, and this even for the diagonal operators on $l^2(\mathbb N)$.
\end{enumerate}

\medskip

Summarizing, we are in trouble with our extension program, and this right from the beginning. In order to have some theory started, however, let us forget about (2), which obviously leads nowhere, and focus on the difficulties in (1). 

\bigskip

In order to cut short the discussion there, regarding the various properties of $T-\lambda$, we can just say that $T-\lambda$ is either invertible with bounded inverse, the ``good case'', or not. We are led in this way to the following definition:

\index{spectrum}
\index{invertible operator}

\begin{definition}
The spectrum of an operator $T\in B(H)$ is the set
$$\sigma(T)=\left\{\lambda\in\mathbb C\Big|T-\lambda\not\in B(H)^{-1}\right\}$$
where $B(H)^{-1}\subset B(H)$ is the set of invertible operators.
\end{definition}

As a basic example, in the finite dimensional case, $H=\mathbb C^N$, the spectrum of a usual matrix $A\in M_N(\mathbb C)$ is the collection of its eigenvalues, taken without multiplicities. We will see many other examples. In general, the spectrum has the following properties:

\index{eigenvalue}

\begin{proposition}
The spectrum of $T\in B(H)$ contains the eigenvalue set
$$\varepsilon(T)=\left\{\lambda\in\mathbb C\Big|\ker(T-\lambda)\neq\{0\}\right\}$$
and $\varepsilon(T)\subset\sigma(T)$ is an equality in finite dimensions, but not in infinite dimensions.
\end{proposition}

\begin{proof}
We have several assertions here, the idea being as follows:

\medskip

(1) First of all, the eigenvalue set is indeed the one in the statement, because $Tx=\lambda x$ tells us precisely that $T-\lambda$ must be not injective. The fact that we have $\varepsilon(T)\subset\sigma(T)$ is clear as well, because if $T-\lambda$ is not injective, it is not bijective.

\medskip

(2) In finite dimensions we have $\varepsilon(T)=\sigma(T)$, because $T-\lambda$ is injective if and only if it is bijective, with the boundedness of the inverse being automatic. 

\medskip

(3) In infinite dimensions we can assume $H=l^2(\mathbb N)$, and the shift operator $S(e_i)=e_{i+1}$ is injective but not surjective. Thus $0\in\sigma(T)-\varepsilon(T)$.
\end{proof}

We will see more examples and counterexamples, and some general theory, in a moment. Philosophically speaking, the best way of thinking at all this is as follows:

\medskip

-- The numbers $\lambda\notin\sigma(T)$ are good, because we can invert $T-\lambda$.

\medskip

-- The numbers $\lambda\in\sigma(T)-\varepsilon(T)$ are bad.

\medskip

-- The eigenvalues $\lambda\in\varepsilon(T)$ are evil.

\medskip

Note that this is somewhat contrary to what happens in linear algebra, where the eigenvalues are highly valued, and cherished, and regarded as being the source of all good things on Earth. Welcome to operator theory, where some things are upside down.

\bigskip

Let us develop now some general theory for the spectrum, or perhaps for its complement, with the promise to come back to eigenvalues later. As a first result, we would like to prove that the spectra are non-empty. This is something tricky, and we will need:

\index{invertible operator}

\begin{proposition}
The following happen:
\begin{enumerate}
\item $||T||<1\implies(1-T)^{-1}=1+T+T^2+\ldots$

\item The set $B(H)^{-1}$ is open.

\item The map $T\to T^{-1}$ is differentiable.
\end{enumerate}
\end{proposition}

\begin{proof}
All these assertions are elementary, as follows:

\medskip

(1) This follows as in the scalar case, the computation being as follows, provided that everything converges under the norm, which amounts in saying that $||T||<1$:
\begin{eqnarray*}
(1-T)(1+T+T^2+\ldots)
&=&1-T+T-T^2+T^2-T^3+\ldots\\
&=&1
\end{eqnarray*}

(2) Assuming $T\in B(H)^{-1}$, let us pick $S\in B(H)$ such that:
$$||T-S||<\frac{1}{||T^{-1}||}$$

We have then the following estimate:
\begin{eqnarray*}
||1-T^{-1}S||
&=&||T^{-1}(T-S)||\\
&\leq&||T^{-1}||\cdot||T-S||\\
&<&1
\end{eqnarray*}

Thus we have $T^{-1}S\in B(H)^{-1}$, and so $S\in B(H)^{-1}$, as desired.

\medskip

(3) In the scalar case, the derivative of $f(t)=t^{-1}$ is $f'(t)=-t^{-2}$. In the present normed space setting the derivative is no longer a number, but rather a linear transformation, which can be found by developing $f(T)=T^{-1}$ at order 1, as follows:
\begin{eqnarray*}
(T+S)^{-1}
&=&((1+ST^{-1})T)^{-1}\\
&=&T^{-1}(1+ST^{-1})^{-1}\\
&=&T^{-1}(1-ST^{-1}+(ST^{-1})^2-\ldots)\\
&\simeq&T^{-1}(1-ST^{-1})\\
&=&T^{-1}-T^{-1}ST^{-1}
\end{eqnarray*}

Thus $f(T)=T^{-1}$ is indeed differentiable, with derivative $f'(T)S=-T^{-1}ST^{-1}$.
\end{proof}

We can now formulate our first theorem about spectra, as follows:

\index{spectrum}

\begin{theorem}
The spectrum of a bounded operator $T\in B(H)$ is:
\begin{enumerate}
\item Compact.

\item Contained in the disc $D_0(||T||)$.

\item Non-empty.
\end{enumerate}
\end{theorem}

\begin{proof}
This can be proved by using Proposition 3.4, along with a bit of complex and functional analysis, for which we refer to Rudin \cite{rud} and Lax \cite{lax}, as follows:

\medskip

(1) In view of (2) below, it is enough to prove that $\sigma(T)$ is closed. But this follows from the following computation, with $|\varepsilon|$ being small:
\begin{eqnarray*}
\lambda\notin\sigma(T)
&\implies&T-\lambda\in B(H)^{-1}\\
&\implies&T-\lambda-\varepsilon\in B(H)^{-1}\\
&\implies&\lambda+\varepsilon\notin\sigma(T)
\end{eqnarray*}

(2) This follows from the following computation:
\begin{eqnarray*}
\lambda>||T||
&\implies&\Big|\Big|\frac{T}{\lambda}\Big|\Big|<1\\
&\implies&1-\frac{T}{\lambda}\in B(H)^{-1}\\
&\implies&\lambda-T\in B(H)^{-1}\\
&\implies&\lambda\notin\sigma(T)
\end{eqnarray*}

(3) Assume by contradiction $\sigma(T)=\emptyset$. Given a linear form $f\in B(H)^*$, consider the following map, which is well-defined, due to our assumption $\sigma(T)=\emptyset$:
$$\varphi:\mathbb C\to\mathbb C\quad,\quad 
\lambda\to f((T-\lambda)^{-1})$$

By using the fact that $T\to T^{-1}$ is differentiable, that we know from Proposition 3.4, we conclude that this map is differentiable, and so holomorphic. Also, we have:
\begin{eqnarray*}
\lambda\to\infty
&\implies&T-\lambda\to\infty\\
&\implies&(T-\lambda)^{-1}\to0\\
&\implies&f((T-\lambda))^{-1}\to0
\end{eqnarray*}

Thus by the Liouville theorem we obtain $\varphi=0$. But, in view of the definition of $\varphi$, this gives $(T-\lambda)^{-1}=0$, which is a contradiction, as desired.
\end{proof}

Here is now a second basic result regarding the spectra, inspired from what happens in finite dimensions, for the usual complex matrices, and which shows that things do not necessarily extend without troubles to the infinite dimensional setting:

\index{spectrum of products}

\begin{theorem}
We have the following formula, valid for any operators $S,T$:
$$\sigma(ST)\cup\{0\}=\sigma(TS)\cup\{0\}$$
In finite dimensions we have $\sigma(ST)=\sigma(TS)$, but this fails in infinite dimensions.
\end{theorem}

\begin{proof}
There are several assertions here, the idea being as follows:

\medskip

(1) This is something that we know in finite dimensions, coming from the fact that the characteristic polynomials of the associated matrices $A,B$ coincide:
$$P_{AB}=P_{BA}$$

Thus we obtain $\sigma(ST)=\sigma(TS)$ in this case, as claimed. Observe that this improves twice the general formula in the statement, first because we have no issues at 0, and second because what we obtain is actually an equality of sets with mutiplicities.

\medskip

(2) In general now, let us first prove the main assertion, stating that $\sigma(ST),\sigma(TS)$ coincide outside 0. We first prove that we have the following implication:
$$1\notin\sigma(ST)\implies1\notin\sigma(TS)$$

Assume indeed that $1-ST$ is invertible, with inverse denoted $R$:
$$R=(1-ST)^{-1}$$

We have then the following formulae, relating our variables $R,S,T$:
$$RST=STR=R-1$$

By using $RST=R-1$, we have the following computation:
\begin{eqnarray*}
(1+TRS)(1-TS)
&=&1+TRS-TS-TRSTS\\
&=&1+TRS-TS-TRS+TS\\
&=&1
\end{eqnarray*}

A similar computation, using $STR=R-1$, shows that we have:
$$(1-TS)(1+TRS)=1$$

Thus $1-TS$ is invertible, with inverse $1+TRS$, which proves our claim. Now by multiplying by scalars, we deduce from this that for any $\lambda\in\mathbb C-\{0\}$ we have:
$$\lambda\notin\sigma(ST)\implies\lambda\notin\sigma(TS)$$ 

But this leads to the conclusion in the statement.

\medskip

(3) Regarding now the counterexample to the formula $\sigma(ST)=\sigma(TS)$, in general, let us take $S$ to be the shift on $H=L^2(\mathbb N)$, given by the following formula:
$$S(e_i)=e_{i+1}$$

As for $T$, we can take it to be the adjoint of $S$, which is the following operator:
$$S^*(e_i)=\begin{cases}
e_{i-1}&{\rm if}\ i>0\\
0&{\rm if}\ i=0
\end{cases}$$

Let us compose now these two operators. In one sense, we have:
$$S^*S=1\implies 0\notin\sigma(S^*S)$$

In the other sense, however, the situation is different, as follows:
$$SS^*=Proj(e_0^\perp)\implies 0\in\sigma(SS^*)$$

Thus, the spectra do not match on $0$, and we have our counterexample, as desired.
\end{proof}

\section*{3b. Spectral radius}

Let us develop now some systematic theory for the computation of the spectra, based on what we know about the eigenvalues of the usual complex matrices. As a first result, which is well-known for the usual matrices, and extends well, we have:

\index{polynomial calculus}

\begin{theorem}
We have the ``polynomial functional calculus'' formula
$$\sigma(P(T))=P(\sigma(T))$$
valid for any polynomial $P\in\mathbb C[X]$, and any operator $T\in B(H)$.
\end{theorem}

\begin{proof}
We pick a scalar $\lambda\in\mathbb C$, and we decompose the polynomial $P-\lambda$:
$$P(X)-\lambda=c(X-r_1)\ldots(X-r_n)$$

We have then the following equivalences:
\begin{eqnarray*}
\lambda\notin\sigma(P(T))
&\iff&P(T)-\lambda\in B(H)^{-1}\\
&\iff&c(T-r_1)\ldots(T-r_n)\in B(H)^{-1}\\
&\iff&T-r_1,\ldots,T-r_n\in B(H)^{-1}\\
&\iff&r_1,\ldots,r_n\notin\sigma(T)\\
&\iff&\lambda\notin P(\sigma(T))
\end{eqnarray*}

Thus, we are led to the formula in the statement.
\end{proof}

The above result is something very useful, and generalizing it will be our next task. As a first ingredient here, assuming that $A\in M_N(\mathbb C)$ is invertible, we have:
$$\sigma(A^{-1})=\sigma(A)^{-1}$$

It is possible to extend this formula to the arbitrary operators, and we will do this in a moment. Before starting, however, we have to think in advance on how to unify this potential result, that we have in mind, with Theorem 3.7 itself. 

\bigskip

What we have to do here is to find a class of functions generalizing both the polynomials $P\in\mathbb C[X]$ and the inverse function $x\to x^{-1}$, and the answer to this question is provided by the rational functions, which are as follows:

\begin{definition}
A rational function $f\in\mathbb C(X)$ is a quotient of polynomials:
$$f=\frac{P}{Q}$$
Assuming that $P,Q$ are prime to each other, we can regard $f$ as a usual function,
$$f:\mathbb C-X\to\mathbb C$$
with $X$ being the set of zeros of $Q$, also called poles of $f$.
\end{definition}

Here the term ``poles'' comes from the fact that, if you want to imagine the graph of such a rational function $f$, in two complex dimensions, what you get is some sort of tent, supported by poles of infinite height, situated at the zeros of $Q$. For more on all this, and on complex analysis in general, we refer as usual to Rudin \cite{rud}. Although a look at an abstract algebra book can be interesting as well.

\bigskip

Now that we have our class of functions, the next step consists in applying them to operators. Here we cannot expect $f(T)$ to make sense for any $f$ and any $T$, for instance because $T^{-1}$ is defined only when $T$ is invertible. We are led in this way to:

\begin{definition}
Given an operator $T\in B(H)$, and a rational function $f=P/Q$ having poles outside $\sigma(T)$, we can construct the following operator,
$$f(T)=P(T)Q(T)^{-1}$$
that we can denote as a usual fraction, as follows,
$$f(T)=\frac{P(T)}{Q(T)}$$
due to the fact that $P(T),Q(T)$ commute, so that the order is irrelevant.
\end{definition}

To be more precise, $f(T)$ is indeed well-defined, and the fraction notation is justified too. In more formal terms, we can say that we have a morphism of complex algebras as follows, with $\mathbb C(X)^T$ standing for the rational functions having poles outside $\sigma(T)$:
$$\mathbb C(X)^T\to B(H)\quad,\quad f\to f(T)$$

Summarizing, we have now a good class of functions, generalizing both the polynomials and the inverse map $x\to x^{-1}$. We can now extend Theorem 3.7, as follows:

\index{rational calculus}

\begin{theorem}
We have the ``rational functional calculus'' formula
$$\sigma(f(T))=f(\sigma(T))$$
valid for any rational function $f\in\mathbb C(X)$ having poles outside $\sigma(T)$.
\end{theorem}

\begin{proof}
We pick a scalar $\lambda\in\mathbb C$, we write $f=P/Q$, and we set:
$$F=P-\lambda Q$$

By using now Theorem 3.7, for this polynomial, we obtain:
\begin{eqnarray*}
\lambda\in\sigma(f(T))
&\iff&F(T)\notin B(H)^{-1}\\
&\iff&0\in\sigma(F(T))\\
&\iff&0\in F(\sigma(T))\\
&\iff&\exists\mu\in\sigma(T),F(\mu)=0\\
&\iff&\lambda\in f(\sigma(T))
\end{eqnarray*}

Thus, we are led to the formula in the statement.
\end{proof}

As an application of the above methods, we can investigate certain special classes of operators, such as the self-adjoint ones, and the unitary ones. Let us start with:

\index{adjoint operator}

\begin{proposition}
The following happen:
\begin{enumerate}
\item We have $\sigma(T^*)=\overline{\sigma(T)}$, for any $T\in B(H)$.

\item If $T=T^*$ then $X=\sigma(T)$ satisfies $X=\overline{X}$.

\item If $U^*=U^{-1}$ then $X=\sigma(U)$ satisfies $X^{-1}=\overline{X}$.
\end{enumerate}
\end{proposition}

\begin{proof}
We have several assertions here, the idea being as follows:

\medskip

(1) The spectrum of the adjoint operator $T^*$ can be computed as follows:
\begin{eqnarray*}
\sigma(T^*)
&=&\left\{\lambda\in\mathbb C\Big|T^*-\lambda\notin B(H)^{-1}\right\}\\
&=&\left\{\lambda\in\mathbb C\Big|T-\bar{\lambda}\notin B(H)^{-1}\right\}\\
&=&\overline{\sigma(T)}
\end{eqnarray*}

(2) This is clear indeed from (1).

\medskip

(3) For a unitary operator, $U^*=U^{-1}$, Theorem 3.10 and (1) give:
$$\sigma(U)^{-1}=\sigma(U^{-1})=\sigma(U^*)=\overline{\sigma(U)}$$

Thus, we are led to the conclusion in the statement.
\end{proof}

In analogy with what happens for the usual matrices, we would like to improve now (2,3) above, with results stating that the spectrum $X=\sigma(T)$ satisfies $X\subset\mathbb R$ for self-adjoints, and $X\subset\mathbb T$ for unitaries. This will be tricky. Let us start with:

\index{unitary}

\begin{theorem}
The spectrum of a unitary operator 
$$U^*=U^{-1}$$
is on the unit circle, $\sigma(U)\subset\mathbb T$. 
\end{theorem}

\begin{proof}
Assuming $U^*=U^{-1}$, we have the following norm computation:
$$||U||
=\sqrt{||UU^*||}
=\sqrt{1}
=1$$

Now if we denote by $D$ the unit disk, we obtain from this:
$$\sigma(U)\subset D$$

On the other hand, once again by using $U^*=U^{-1}$, we have as well:
$$||U^{-1}||
=||U^*||
=||U||
=1$$

Thus, as before with $D$ being the unit disk in the complex plane, we have:
$$\sigma(U^{-1})\subset D$$

Now by using Theorem 3.10, we obtain $\sigma(U)
\subset D\cap D^{-1}
=\mathbb T$, as desired.
\end{proof}

We have as well a similar result for self-adjoints, as follows:

\index{self-adjoint operator}

\begin{theorem}
The spectrum of a self-adjoint operator
$$T=T^*$$
consists of real numbers, $\sigma(T)\subset\mathbb R$.
\end{theorem}

\begin{proof}
The idea is that we can deduce the result from Theorem 3.12, by using the following remarkable rational function, depending on a parameter $r\in\mathbb R$:
$$f(z)=\frac{z+ir}{z-ir}$$

Indeed, for $r>>0$ the operator $f(T)$ is well-defined, and we have:
$$\left(\frac{T+ir}{T-ir}\right)^*
=\frac{T-ir}{T+ir}
=\left(\frac{T+ir}{T-ir}\right)^{-1}$$

Thus $f(T)$ is unitary, and by using Theorem 3.12 we obtain:
\begin{eqnarray*}
\sigma(T)
&\subset&f^{-1}(f(\sigma(T)))\\
&=&f^{-1}(\sigma(f(T)))\\
&\subset&f^{-1}(\mathbb T)\\
&=&\mathbb R
\end{eqnarray*}

Thus, we are led to the conclusion in the statement.
\end{proof}

As a theoretical remark, it is possible to deduce as well Theorem 3.12 from Theorem 3.13, by performing the above computation in the other sense. Indeed, by assuming that Theorem 3.13 holds indeed, and starting with a unitary $U\in B(H)$, we obtain:
\begin{eqnarray*}
\sigma(U)
&\subset&f(f^{-1}(\sigma(U)))\\
&=&f(\sigma(f^{-1}(U)))\\
&\subset&f(\mathbb R)\\
&=&\mathbb T
\end{eqnarray*} 

As a conclusion now, we have so far a beginning of spectral theory, with results allowing us to investigate the unitaries and the self-adjoints, and with the remark that these two classes of operators are related by a certain wizarding rational function, namely:
$$f(z)=\frac{z+ir}{z-ir}$$

Let us keep building on this, with more complex analysis involved. One key thing that we know about matrices, and which follows for instance by using the fact that the diagonalizable matrices are dense, is the following formula:
$$\sigma(e^A)=e^{\sigma(A)}$$

We would like to have such formulae for the general operators $T\in B(H)$, but this is something quite technical. Consider the rational calculus morphism from Definition 3.9, which is as follows, with the exponent standing for ``having poles outside $\sigma(T)$'':
$$\mathbb C(X)^T\to B(H)\quad,\quad 
f\to f(T)$$

As mentioned before, the rational functions are holomorphic outside their poles, and this raises the question of extending this morphism, as follows:
$$Hol(\sigma(T))\to B(H)\quad,\quad
f\to f(T)$$

Normally this can be done in several steps. Let us start with:

\begin{proposition}
We can exponentiate any operator $T\in B(H)$, by setting:
$$e^T=\sum_{k=0}^\infty\frac{T^k}{k!}$$ 
Similarly, we can define $f(T)$, for any holomorphic function $f:\mathbb C\to\mathbb C$.
\end{proposition}

\begin{proof}
We must prove that the series defining $e^T$ converges, and this follows from:
$$||e^T||
\leq\sum_{k=0}^\infty\frac{||T||^k}{k!}
=e^{||T||}$$

The case of the arbitrary holomorphic functions $f:\mathbb C\to\mathbb C$ is similar.
\end{proof}

In general, the holomorphic functions are not entire, and the above method won't cover the rational functions $f\in\mathbb C(X)^T$ that we want to generalize. Thus, we must use something else. And the answer here comes from the Cauchy formula:
$$f(t)=\frac{1}{2\pi i}\int_\gamma\frac{f(z)}{z-t}\,dz$$

Indeed, given a rational function $f\in\mathbb C(X)^T$, the operator $f(T)\in B(H)$, constructed in Definition 3.9, can be recaptured in an analytic way, as follows:
$$f(T)=\frac{1}{2\pi i}\int_\gamma\frac{f(z)}{z-T}\,dz$$

Now given an arbitrary function $f\in Hol(\sigma(T))$, we can define $f(T)\in B(H)$ by the exactly same formula, and we obtain in this way the desired correspondence:
$$Hol(\sigma(T))\to B(H)\quad,\quad 
f\to f(T)$$

This was for the plan. In practice now, all this needs a bit of care, with many verifications needed, and with the technical remark that a winding number must be added to the above Cauchy formulae, for things to be correct. The result is as follows: 

\index{holomorphic calculus}
\index{Cauchy formula}

\begin{theorem}
We have the ``holomorphic functional calculus'' formula
$$\sigma(f(T))=f(\sigma(T))$$
valid for any holomorphic function $f\in Hol(\sigma(T))$.
\end{theorem}

\begin{proof}
This is something that we will not really need, for the purposes of the present book, which is more algebraic than analytic, but here is the general idea:

\medskip

(1) As explained above, given a rational function $f\in\mathbb C(X)^T$, the corresponding operator $f(T)\in B(H)$ can be recaptured in an analytic way, as follows:
$$f(T)=\frac{1}{2\pi i}\int_\gamma\frac{f(z)}{z-T}\,dz$$

(2) Now given an arbitrary function $f\in Hol(\sigma(T))$, we can define $f(T)\in B(H)$ by the exactly same formula, and we obtain in this way the desired correspondence:
$$Hol(\sigma(T))\to B(H)\quad,\quad 
f\to f(T)$$

(3) In practice now, all this needs a bit of care, notably with the verification of the fact that the operator $f(T)\in B(H)$ does not depend on $\gamma$, and with the technical remark that a winding number must be added to the above Cauchy formulae, for things to be correct. But this can be done via a standard study, keeping in mind the fact that in the case $H=\mathbb C$, where our operators are usual numbers, $B(H)=\mathbb C$, what we want to do is simply proving that the usual Cauchy formula holds indeed.

\medskip

(4) Now with this correspondence $f\to f(T)$ constructed, and so with the formula in the statement, namely $\sigma(f(T))=f(\sigma(T))$, making now sense, it remains to prove that this formula holds indeed. But this follows as well via a careful use of the Cauchy formula, or by using approximation by polynomials, or rational functions.
\end{proof}

As already said, the above result is important for advanced operator theory and applications, and we will not get further into this subject. We will be back, however, to all this in the special case of the normal operators, which is of particular interest for us. 

\bigskip

In order to formulate now our next result, we will need the following notion:

\index{spectral radius}

\begin{definition}
Given an operator $T\in B(H)$, its spectral radius 
$$\rho(T)\in\big[0,||T||\big]$$
is the radius of the smallest disk centered at $0$ containing $\sigma(T)$. 
\end{definition}

Here we have included for convenience a number of basic results from Theorem 3.5, namely the fact that the spectrum is non-empty, and is contained in the disk $D_0(||T||)$, which provide us respectively with the inequalities $\rho(T)\geq0$, with the usual convention $\sup\emptyset=-\infty$, and $\rho(T)\leq||T||$. Now with this notion in hand, we have the following key result, improving our key result so far, namely $\sigma(T)\neq\emptyset$, from Theorem 3.5:

\begin{theorem}
The spectral radius of an operator $T\in B(H)$ is given by
$$\rho(T)=\lim_{n\to\infty}||T^n||^{1/n}$$
and in this formula, we can replace the limit by an inf.
\end{theorem}

\begin{proof}
We have several things to be proved, the idea being as follows:

\medskip

(1) Our first claim is that the numbers $u_n=||T^n||^{1/n}$ satisfy:
$$(n+m)u_{n+m}\leq nu_n+mu_m$$

Indeed, we have the following estimate, using the Young inequality $ab\leq a^p/p+b^q/q$, with exponents $p=(n+m)/n$ and $q=(n+m)/m$:
\begin{eqnarray*}
u_{n+m}
&=&||T^{n+m}||^{1/(n+m)}\\
&\leq&||T^n||^{1/(n+m)}||T^m||^{1/(n+m)}\\
&\leq&||T^n||^{1/n}\cdot\frac{n}{n+m}+||T^m||^{1/m}\cdot\frac{m}{n+m}\\
&=&\frac{nu_n+mu_m}{n+m}
\end{eqnarray*}

(2) Our second claim is that the second assertion holds, namely:
$$\lim_{n\to\infty}||T^n||^{1/n}=\inf_n||T^n||^{1/n}$$

For this purpose, we just need the inequality found in (1). Indeed, fix $m\geq1$, let $n\geq1$, and write $n=lm+r$ with $0\leq r\leq m-1$. By using twice $u_{ab}\leq u_b$, we get:
\begin{eqnarray*}
u_n
&\leq&\frac{1}{n}( lmu_{lm}+ru_r)\\
&\leq&\frac{1}{n}( lmu_{m}+ru_1)\\
&\leq&u_{m}+\frac{r}{n}\,u_1
\end{eqnarray*}

It follows that we have $\lim\sup_nu_n\leq u_m$, which proves our claim.

\medskip

(3) Summarizing, we are left with proving the main formula, which is as follows, and with the remark that we already know that the sequence on the right converges:
$$\rho(T)=\lim_{n\to\infty}||T^n||^{1/n}$$

In one sense, we can use the polynomial calculus formula $\sigma(T^n)=\sigma(T)^n$. Indeed, this gives the following estimate, valid for any $n$, as desired:
\begin{eqnarray*}
\rho(T)
&=&\sup_{\lambda\in\sigma(T)}|\lambda|\\
&=&\sup_{\rho\in\sigma(T)^n}|\rho|^{1/n}\\
&=&\sup_{\rho\in\sigma(T^n)}|\rho|^{1/n}\\
&=&\rho(T^n)^{1/n}\\
&\leq&||T^n||^{1/n}
\end{eqnarray*}

(4) For the reverse inequality, we fix a number $\rho>\rho(T)$, and we want to prove that we have $\rho\geq\lim_{n\to\infty}||T^n||^{1/n}$. By using the Cauchy formula, we have:
\begin{eqnarray*}
\frac{1}{2\pi i}\int_{|z|=\rho}\frac{z^n}{z-T}\,dz
&=&\frac{1}{2\pi i}\int_{|z|=\rho}\sum_{k=0}^\infty z^{n-k-1}T^k\,dz\\
&=&\sum_{k=0}^\infty\frac{1}{2\pi i}\left(\int_{|z|=\rho}z^{n-k-1}dz\right)T^k\\
&=&\sum_{k=0}^\infty\delta_{n,k+1}T^k\\
&=&T^{n-1}
\end{eqnarray*}

By applying the norm we obtain from this formula:
\begin{eqnarray*}
||T^{n-1}||
&\leq&\frac{1}{2\pi}\int_{|z|=\rho}\left|\left|\frac{z^n}{z-T}\right|\right|\,dz\\
&\leq&\rho^n\cdot\sup_{|z|=\rho}\left|\left|\frac{1}{z-T}\right|\right|
\end{eqnarray*}

Since the sup does not depend on $n$, by taking $n$-th roots, we obtain in the limit:
$$\rho\geq\lim_{n\to\infty}||T^n||^{1/n}$$

Now recall that $\rho$ was by definition an arbitrary number satisfying $\rho>\rho(T)$. Thus, we have obtained the following estimate, valid for any $T\in B(H)$:
$$\rho(T)\geq\lim_{n\to\infty}||T^n||^{1/n}$$

Thus, we are led to the conclusion in the statement.
\end{proof}

In the case of the normal elements, we have the following finer result:

\index{normal operator}

\begin{theorem}
The spectral radius of a normal element,
$$TT^*=T^*T$$
is equal to its norm.
\end{theorem}

\begin{proof}
We can proceed in two steps, as follows:

\medskip

\underline{Step 1}. In the case $T=T^*$ we have $||T^n||=||T||^n$ for any exponent of the form $n=2^k$, by using the formula $||TT^*||=||T||^2$, and by taking $n$-th roots we get:
$$\rho(T)\geq||T||$$

Thus, we are done with the self-adjoint case, with the result $\rho(T)=||T||$.

\medskip

\underline{Step 2}. In the general normal case $TT^*=T^*T$ we have $T^n(T^n)^*=(TT^*)^n$, and by using this, along with the result from Step 1, applied to $TT^*$, we obtain:
\begin{eqnarray*}
\rho(T)
&=&\lim_{n\to\infty}||T^n||^{1/n}\\
&=&\sqrt{\lim_{n\to\infty}||T^n(T^n)^*||^{1/n}}\\
&=&\sqrt{\lim_{n\to\infty}||(TT^*)^n||^{1/n}}\\
&=&\sqrt{\rho(TT^*)}\\
&=&\sqrt{||T||^2}\\
&=&||T||
\end{eqnarray*}

Thus, we are led to the conclusion in the statement.
\end{proof}

As a first comment, the spectral radius formula $\rho(T)=||T||$ does not hold in general, the simplest counterexample being the following non-normal matrix:
$$J=\begin{pmatrix}0&1\\0&0\end{pmatrix}$$

As another comment, we can combine the formula $\rho(T)=||T||$ for normal operators with the formula $||TT^*||=||T||^2$, and we are led to the following statement:

\index{norm of operators}

\begin{theorem}
The norm of $B(H)$ is given by
$$||T||=\sqrt{\sup\left\{\lambda\in\mathbb C\Big|TT^*-\lambda\notin B(H)^{-1}\right\}}$$
and so is a purely algebraic quantity.
\end{theorem}

\begin{proof}
We have the following computation, using the formula $||TT^*||=||T||^2$, then the spectral radius formula for $TT^*$, and finally the definition of the spectral radius:
\begin{eqnarray*}
||T||
&=&\sqrt{||TT^*||}\\
&=&\sqrt{\rho(TT^*)}\\
&=&\sqrt{\sup\left\{\lambda\in\mathbb C\Big| \lambda\in\sigma(TT^*)\right\}}\\
&=&\sqrt{\sup\left\{\lambda\in\mathbb C\Big| TT^*-\lambda\notin B(H)^{-1}\right\}}
\end{eqnarray*}

Thus, we are led to the conclusion in the statement.
\end{proof}

The above result is quite interesting, philosophically speaking. We will be back to this, with further results and comments on $B(H)$, and other algebras of the same type.

\section*{3c. Normal operators}

By using Theorem 3.18 we can say a number of non-trivial things concerning the normal operators, commonly known as ``spectral theorem for normal operators''. As a first result here, we can improve the polynomial functional calculus formula:

\index{normal operator}
\index{polynomial calculus}

\begin{theorem}
Given $T\in B(H)$ normal, we have a morphism of algebras
$$\mathbb C[X]\to B(H)\quad,\quad 
P\to P(T)$$
having the properties $||P(T)||=||P_{|\sigma(T)}||$, and $\sigma(P(T))=P(\sigma(T))$.
\end{theorem}

\begin{proof}
This is an improvement of Theorem 3.7 in the normal case, with the extra assertion being the norm estimate. But the element $P(T)$ being normal, we can apply to it the spectral radius formula for normal elements, and we obtain:
\begin{eqnarray*}
||P(T)||
&=&\rho(P(T))\\
&=&\sup_{\lambda\in\sigma(P(T))}|\lambda|\\
&=&\sup_{\lambda\in P(\sigma(T))}|\lambda|\\
&=&||P_{|\sigma(T)}||
\end{eqnarray*}

Thus, we are led to the conclusions in the statement.
\end{proof}

We can improve as well the rational calculus formula, and the holomorphic calculus formula, in the same way. Importantly now, at a more advanced level, we have:

\index{normal operator}
\index{continuous calculus}

\begin{theorem}
Given $T\in B(H)$ normal, we have a morphism of algebras
$$C(\sigma(T))\to B(H)\quad,\quad 
f\to f(T)$$
which is isometric, $||f(T)||=||f||$, and has the property $\sigma(f(T))=f(\sigma(T))$.
\end{theorem}

\begin{proof}
The idea here is to ``complete'' the morphism in Theorem 3.20, namely:
$$\mathbb C[X]\to B(H)\quad,\quad 
P\to P(T)$$

Indeed, we know from Theorem 3.20 that this morphism is continuous, and is in fact isometric, when regarding the polynomials $P\in\mathbb C[X]$ as functions on $\sigma(T)$:
$$||P(T)||=||P_{|\sigma(T)}||$$

We conclude from this that we have a unique isometric extension, as follows:
$$C(\sigma(T))\to B(H)\quad,\quad  
f\to f(T)$$

It remains to prove $\sigma(f(T))=f(\sigma(T))$, and we can do this by double inclusion:

\medskip

``$\subset$'' Given a continuous function $f\in C(\sigma(T))$, we must prove that we have:
$$\lambda\notin f(\sigma(T))\implies\lambda\notin\sigma(f(T))$$

For this purpose, consider the following function, which is well-defined:
$$\frac{1}{f-\lambda}\in C(\sigma(T))$$

We can therefore apply this function to $T$, and we obtain:
$$\left(\frac{1}{f-\lambda}\right)T=\frac{1}{f(T)-\lambda}$$

In particular $f(T)-\lambda$ is invertible, so  $\lambda\notin\sigma(f(T))$, as desired.

\medskip

``$\supset$'' Given a continuous function $f\in C(\sigma(T))$, we must prove that we have: 
$$\lambda\in f(\sigma(T))\implies\lambda\in\sigma(f(T))$$

But this is the same as proving that we have:
$$\mu\in\sigma(T)\implies f(\mu)\in\sigma(f(T))$$

For this purpose, we approximate our function by polynomials, $P_n\to f$, and we examine the following convergence, which follows from $P_n\to f$:
$$P_n(T)-P_n(\mu)\to f(T)-f(\mu)$$

We know from polynomial functional calculus that we have:
$$P_n(\mu)
\in P_n(\sigma(T))
=\sigma(P_n(T))$$

Thus, the operators $P_n(T)-P_n(\mu)$ are not invertible. On the other hand, we know that the set formed by the invertible operators is open, so its complement is closed. Thus the limit $f(T)-f(\mu)$ is not invertible either, and so $f(\mu)\in\sigma(f(T))$, as desired.
\end{proof}

As an important comment, Theorem 3.21 is not exactly in final form, because it misses an important point, namely that our correspondence maps:
$$\bar{z}\to T^*$$

However, this is something non-trivial, and we will be back to this later. Observe however that Theorem 3.21 is fully powerful for the self-adjoint operators, $T=T^*$, where the spectrum is real, and so where $z=\bar{z}$ on the spectrum. We will be back to this.

\bigskip

As a second result now, along the same lines, we can further extend Theorem 3.21 into a measurable functional calculus theorem, as follows:

\index{normal operator}
\index{measurable calculus}

\begin{theorem}
Given $T\in B(H)$ normal, we have a morphism of algebras as follows, with $L^\infty$ standing for abstract measurable functions, or Borel functions,
$$L^\infty(\sigma(T))\to B(H)\quad,\quad 
f\to f(T)$$
which is isometric, $||f(T)||=||f||$, and has the property $\sigma(f(T))=f(\sigma(T))$.
\end{theorem}

\begin{proof}
As before, the idea will be that of ``completing'' what we have. To be more precise, we can use the Riesz theorem and a polarization trick, as follows:

\medskip

(1) Given a vector $x\in H$, consider the following functional:
$$C(\sigma(T))\to\mathbb C\quad,\quad 
g\to<g(T)x,x>$$

By the Riesz theorem, this functional must be the integration with respect to a certain measure $\mu$ on the space $\sigma(T)$. Thus, we have a formula as follows:
$$<g(T)x,x>=\int_{\sigma(T)}g(z)d\mu(z)$$

Now given an arbitrary Borel function $f\in L^\infty(\sigma(T))$, as in the statement, we can define a number $<f(T)x,x>\in\mathbb C$, by using exactly the same formula, namely:
$$<f(T)x,x>=\int_{\sigma(T)}f(z)d\mu(z)$$

Thus, we have managed to define numbers $<f(T)x,x>\in\mathbb C$, for all vectors $x\in H$, and in addition we can recover these numbers as follows, with $g_n\in C(\sigma(T))$:
$$<f(T)x,x>=\lim_{g_n\to f}<g_n(T)x,x>$$ 

(2) In order to define now numbers $<f(T)x,y>\in\mathbb C$, for all vectors $x,y\in H$, we can use a polarization trick. Indeed, for any operator $S\in B(H)$ we have:
$$<S(x+y),x+y>=<Sx,x>+<Sy,y>+<Sx,y>+<Sy,x>$$

By replacing $y\to iy$, we have as well the following formula:
$$<S(x+iy),x+iy>=<Sx,x>+<Sy,y>-i<Sx,y>+i<Sy,x>$$

By multiplying this latter formula by $i$, we obtain the following formula:
$$i<S(x+iy),x+iy>=i<Sx,x>+i<Sy,y>+<Sx,y>-<Sy,x>$$

Now by summing this latter formula with the first one, we obtain:
\begin{eqnarray*}
<S(x+y),x+y>+i<S(x+iy),x+iy>
&=&(1+i)[<Sx,x>+<Sy,y>]\\
&+&2<Sx,y>
\end{eqnarray*}

(3) But with this, we can now finish. Indeed, by combining (1,2), given a Borel function $f\in L^\infty(\sigma(T))$, we can define numbers $<f(T)x,y>\in\mathbb C$ for any $x,y\in H$, and it is routine to check, by using approximation by continuous functions $g_n\to f$ as in (1), that we obtain in this way an operator $f(T)\in B(H)$, having all the desired properties.
\end{proof}

The same comments as before apply. Theorem 3.22 is not exactly in final form, because it misses an important point, namely that our correspondence maps:
$$\bar{z}\to T^*$$

However, this is something non-trivial, and we will be back to this later. Observe however that Theorem 3.22 is fully powerful for the self-adjoint operators, $T=T^*$, where the spectrum is real, and so where $z=\bar{z}$ on the spectrum. We will be back to this.

\bigskip

As another comment, the above result and its proof provide us with more than a Borel functional calculus, because what we got is a certain measure on the spectrum $\sigma(T)$, along with a functional calculus for the $L^\infty$ functions with respect to this measure. We will be back to this later, and for the moment we will only need Theorem 3.22 as formulated, with $L^\infty(\sigma(T))$ standing, a bit abusively, for the Borel functions on $\sigma(T)$.

\section*{3d. Diagonalization}

We can now diagonalize the normal operators. We will do this in 3 steps, first for the self-adjoint operators, then for the families of commuting self-adjoint operators, and finally for the general normal operators, by using a trick of the following type:
$$T=Re(T)+iIm(T)$$

The diagonalization in infinite dimensions is more tricky than in finite dimensions, and instead of writing a formula of type $T=UDU^*$, with $U,D\in B(H)$ being respectively unitary and diagonal, we will express our operator as $T=U^*MU$, with $U:H\to K$ being a certain unitary, and with $M\in B(K)$ being a certain diagonal operator. 

\bigskip

This is indeed how the spectral theorem is best formulated, in view of applications. In practice, the explicit construction of $U,M$, which will be actually rather part of the proof, is also needed. For the self-adjoint operators, the statement and proof are as follows:

\index{self-adjoint operator}
\index{diagonalization}

\begin{theorem}
Any self-adjoint operator $T\in B(H)$ can be diagonalized,
$$T=U^*M_fU$$
with $U:H\to L^2(X)$ being a unitary operator from $H$ to a certain $L^2$ space associated to $T$, with $f:X\to\mathbb R$ being a certain function, once again associated to $T$, and with
$$M_f(g)=fg$$
being the usual multiplication operator by $f$, on the Hilbert space $L^2(X)$.
\end{theorem}

\begin{proof}
The construction of $U,f$ can be done in several steps, as follows:

\medskip

(1) We first prove the result in the special case where our operator $T$ has a cyclic vector $x\in H$, with this meaning that the following holds:
$$\overline{span\left(T^kx\Big|n\in\mathbb N\right)}=H$$

For this purpose, let us go back to the proof of Theorem 3.22. We will use the following formula from there, with $\mu$ being the measure on $X=\sigma(T)$ associated to $x$:
$$<g(T)x,x>=\int_{\sigma(T)}g(z)d\mu(z)$$

Our claim is that we can define a unitary $U:H\to L^2(X)$, first on the dense part spanned by the vectors $T^kx$, by the following formula, and then by continuity:
$$U[g(T)x]=g$$

Indeed, the following computation shows that $U$ is well-defined, and isometric:
\begin{eqnarray*}
||g(T)x||^2
&=&<g(T)x,g(T)x>\\
&=&<g(T)^*g(T)x,x>\\
&=&<|g|^2(T)x,x>\\
&=&\int_{\sigma(T)}|g(z)|^2d\mu(z)\\
&=&||g||_2^2
\end{eqnarray*}

We can then extend $U$ by continuity into a unitary $U:H\to L^2(X)$, as claimed. Now observe that we have the following formula:
\begin{eqnarray*}
UTU^*g
&=&U[Tg(T)x]\\
&=&U[(zg)(T)x]\\
&=&zg
\end{eqnarray*} 

Thus our result is proved in the present case, with $U$ as above, and with $f(z)=z$.

\medskip

(2) We discuss now the general case. Our first claim is that $H$ has a decomposition as follows, with each $H_i$ being invariant under $T$, and admitting a cyclic vector $x_i$:
$$H=\bigoplus_iH_i$$

Indeed, this is something elementary, the construction being by recurrence in finite dimensions, in the obvious way, and by using the Zorn lemma in general. Now with this decomposition in hand, we can make a direct sum of the diagonalizations obtained in (1), for each of the restrictions $T_{|H_i}$, and we obtain the formula in the statement.
\end{proof}

We have the following technical generalization of the above result:

\index{commuting self-adjoint operators}
\index{diagonalization}

\begin{theorem}
Any family of commuting self-adjoint operators $T_i\in B(H)$ can be jointly diagonalized,
$$T_i=U^*M_{f_i}U$$
with $U:H\to L^2(X)$ being a unitary operator from $H$ to a certain $L^2$ space associated to $\{T_i\}$, with $f_i:X\to\mathbb R$ being certain functions, once again associated to $T_i$, and with
$$M_{f_i}(g)=f_ig$$
being the usual multiplication operator by $f_i$, on the Hilbert space $L^2(X)$.
\end{theorem}

\begin{proof}
This is similar to the proof of Theorem 3.23, by suitably modifying the measurable calculus formula, and the measure $\mu$ itself, as to have this formula working for all the operators $T_i$. With this modification done, everything extends.
\end{proof}

In order to discuss now the case of the arbitrary normal operators, we will need:

\begin{proposition}
Any operator $T\in B(H)$ can be written as
$$T=Re(T)+iIm(T)$$
with $Re(T),Im(T)\in B(H)$ being self-adjoint, and this decomposition is unique.
\end{proposition}

\begin{proof}
This is something elementary, the idea being as follows:

\medskip

(1) As a first observation, in the case $H=\mathbb C$ our operators are usual complex numbers, and the formula in the statement corresponds to the following basic fact:
$$z=Re(z)+iIm(z)$$

(2) In general now, we can use the same formulae for the real and imaginary part as in the complex number case, the decomposition formula being as follows:
$$T=\frac{T+T^*}{2}+i\cdot\frac{T-T^*}{2i}$$

To be more precise, both the operators on the right are self-adjoint, and the summing formula holds indeed, and so we have our decomposition result, as desired.

\medskip

(3) Regarding now the uniqueness, by linearity it is enough to show that $R+iS=0$ with $R,S$ both self-adjoint implies $R=S=0$. But this follows by applying the adjoint to $R+iS=0$, which gives $R-iS=0$, and so $R=S=0$, as desired.
\end{proof}

As a comment here, the above result is just the ``tip of the iceberg'', in what regards decomposition results for the operators $T\in B(H)$, in analogy with decomposition results for the complex numbers $z\in\mathbb C$. As a sample result here, improving Proposition 3.25, we can write any operator $T\in B(H)$ as a linear combination of 4 positive operators, by writing both $Re(T),Im(T)$ as differences of positive operators. More on this later.

\bigskip

Good news, after all these preliminaries, that you enjoyed I hope, as much as I did, we can eventually discuss the case of arbitrary normal operators. We have here the following result, generalizing what we know from chapter 1 about the normal matrices:

\index{normal operator}
\index{diagonalization}

\begin{theorem}
Any normal operator $T\in B(H)$ can be diagonalized,
$$T=U^*M_fU$$
with $U:H\to L^2(X)$ being a unitary operator from $H$ to a certain $L^2$ space associated to $T$, with $f:X\to\mathbb C$ being a certain function, once again associated to $T$, and with
$$M_f(g)=fg$$
being the usual multiplication operator by $f$, on the Hilbert space $L^2(X)$.
\end{theorem}

\begin{proof}
This is our main diagonalization theorem, the idea being as follows:

\medskip

(1) Consider the decomposition of $T$ into its real and imaginary parts, as constructed in the proof of Proposition 3.25, namely:
$$T=\frac{T+T^*}{2}+i\cdot\frac{T-T^*}{2i}$$

We know that the real and imaginary parts are self-adjoint operators. Now since $T$ was assumed to be normal, $TT^*=T^*T$, these real and imaginary parts commute:
$$\left[\frac{T+T^*}{2}\,,\,\frac{T-T^*}{2i}\right]=0$$

Thus Theorem 3.24 applies to these real and imaginary parts, and gives the result.

\medskip

(2) Alternatively, we can use methods similar to those that we used in chapter 1, in order to deal with the usual normal matrices, involving the special relation between $T$ and the operator $TT^*$, which is self-adjoint. We will leave this as an instructive exercise.
\end{proof}

This was for our series of diagonalization theorems. There is of course one more result here, regarding the families of commuting normal operators, as follows:

\index{commuting normal operators}
\index{diagonalization}

\begin{theorem}
Any family of commuting normal operators $T_i\in B(H)$ can be jointly diagonalized,
$$T_i=U^*M_{f_i}U$$
with $U:H\to L^2(X)$ being a unitary operator from $H$ to a certain $L^2$ space associated to $\{T_i\}$, with $f_i:X\to\mathbb C$ being certain functions, once again associated to $T_i$, and with
$$M_{f_i}(g)=f_ig$$
being the usual multiplication operator by $f_i$, on the Hilbert space $L^2(X)$.
\end{theorem}

\begin{proof}
This is similar to the proof of Theorem 3.24 and Theorem 3.26, by combining the arguments there. To be more precise, this follows as Theorem 3.24, by using the decomposition trick from the proof of Theorem 3.26.
\end{proof}

With the above diagonalization results in hand, we can now ``fix'' the continuous and measurable functional calculus theorems, with a key complement, as follows:

\index{polynomial calculus}
\index{rational calculus}
\index{holomorphic calculus}
\index{continuous calculus}
\index{measurable calculus}
\index{adjoint operator}

\begin{theorem}
Given a normal operator $T\in B(H)$, the following hold, for both the functional calculus and the measurable calculus morphisms:
\begin{enumerate}
\item These morphisms are $*$-morphisms.

\item The function $\bar{z}$ gets mapped to $T^*$.

\item The functions $Re(z),Im(z)$ get mapped to $Re(T),Im(T)$.

\item The function $|z|^2$ gets mapped to $TT^*=T^*T$.

\item If $f$ is real, then $f(T)$ is self-adjoint. 
\end{enumerate}
\end{theorem}

\begin{proof}
These assertions are more or less equivalent, with (1) being the main one, which obviously implies everything else. But this assertion (1) follows from the diagonalization result for normal operators, from Theorem 3.26.
\end{proof}

This was for the spectral theory of arbitrary and normal operators, or at least for the basics of this theory. As a conclusion here, our main results are as follows:

\medskip

\begin{enumerate}
\item Regarding the arbitrary operators, the main results here, or rather the most advanced results, are the holomorphic calculus formula from Theorem 3.15, and the spectral radius estimate from Theorem 3.17.

\medskip

\item For the self-adjoint operators, the main results are the spectral radius formula from Theorem 3.18, the measurable calculus formula from Theorem 3.22, and the diagonalization result from Theorem 3.23.

\medskip

\item For general normal operators, the main results are the spectral radius formula from Theorem 3.18, the measurable calculus formula from Theorem 3.22, complemented by Theorem 3.28, and the diagonalization result in Theorem 3.26.
\end{enumerate}

\medskip

There are of course many other things that can be said about the spectral theory of the bounded operators $T\in B(H)$, and on that of the unbounded operators too. As a complement, we recommend any good operator theory book, with the comment however that there is a bewildering choice here, depending on taste, and on what exactly you want to do with your operators $T\in B(H)$. In what concerns us, who are rather into general quantum mechanics, but with our operators being bounded, good choices are the functional analysis book of Lax \cite{lax}, or the operator algebra book of Blackadar \cite{bla}.

\section*{3e. Exercises} 

The main theoretical notion introduced in this chapter was that of the spectrum of an operator, and as a first exercise here, we have:

\begin{exercise}
Prove that for the usual matrices $A,B\in M_N(\mathbb C)$ we have
$$\sigma^+(AB)=\sigma^+(BA)$$
where $\sigma^+$ denotes the set of eigenvalues, taken with multiplicities.
\end{exercise}

As a remark, we have seen in the above that $\sigma(AB)=\sigma(BA)$ holds outside $\{0\}$, and the equality on $\{0\}$ holds as well, because $AB$ is invertible if and only if $BA$ is invertible. However, in what regards the eigenvalues taken with multiplicities, things are more tricky here, and the answer should be somewhere inside your linear algebra knowledge.

\begin{exercise}
Clarify, with examples and counterexamples, the relation between the eigenvalues of an operator $T\in B(H)$, and its spectrum $\sigma(T)\subset\mathbb C$. 
\end{exercise}

Here, as usual, the counterexamples could only come from the shift operator $S$, on the space $H=l^2(\mathbb N)$. As a bonus exercise here, try computing the spectrum of $S$.

\begin{exercise}
Draw the picture of the following function, and of its inverse,
$$f(z)=\frac{z+ir}{z-ir}$$
with $r\in\mathbb R$, and prove that for $r>>0$ and $T=T^*$, the element $f(T)$ is well-defined.
\end{exercise}

This is something that we used in the above, when computing spectra of self-adjoints and unitaries, and the problem is that of working out all the details.

\begin{exercise}
Comment on the spectral radius theorem, stating that for a normal operator, $TT^*=T^*T$, the spectral radius is equal to the norm,
$$\rho(T)=||T||$$
with examples and counterexamples, and simpler proofs of well, in various particular cases of interest, such as the finite dimensional one.
\end{exercise}

This is of course something a bit philosophical, but the spectral radius theorem being our key technical result so far, some further thinking on it is definitely a good thing.

\begin{exercise}
Develop a theory of $*$-algebras $A$ for which the quantity
$$||a||=\sqrt{\sup\left\{\lambda\in\mathbb C\Big|aa^*-\lambda\notin A^{-1}\right\}}$$
defines a norm, for the elements $a\in A$.
\end{exercise}

As pointed out in the above, the spectral radius formula shows that for $A=B(H)$ the norm is given by the above formula, and so there should be such a theory of ``good'' $*$-algebras, with $A=B(H)$ as a main example. However, this is tricky.

\begin{exercise}
Find and write down a proof for the spectral theorem for normal operators in the spirit of the proof for normal matrices from chapter 1, and vice versa.
\end{exercise}

To be more precise, the problem is that the proof of the spectral theorem for the usual matrices, from chapter 1, was using a certain kind of trick, while the proof of the spectral theorem for the arbitrary operators, given in this chapter, was using some other kind of trick. Thus, for fully understanding all this, working out more proofs, both for the usual matrices and for the arbitrary operators, is a useful thing.

\begin{exercise}
Find and write down an enhancement of the proof given above for the spectral theorem, as for $\bar{z}\to T^*$ to appear way before the end of the proof.
\end{exercise} 

This is something a bit philosophical, and check here first the various comments made above, and maybe work out this as well in parallel with the previous exercise.

\chapter{Compact operators}

\section*{4a. Polar decomposition}

We have seen so far the basic theory of bounded operators, in the arbitrary, normal and self-adjoint cases, and in a few other cases of interest. In this chapter we discuss a number of more specialized questions, for the most dealing with the compact operators, which are particularly close, conceptually speaking, to the usual complex matrices.

\bigskip

We have in fact considerably many interesting things that we can talk about, in this final chapter on operator theory, and our choices will be as follows:

\bigskip

(1) Before anything, at the general level, we would like to understand the matrix and operator theory analogues of the various things that we know about the complex numbers $z\in M_1(\mathbb C)$, such as $z\bar{z}=|z|^2$, or $z=re^{it}$ and so on. We will discuss this first.

\bigskip

(2) Then, motivated by advanced linear algebra, we will go on a lengthy discussion on the algebra of compact operators $K(H)\subset B(H)$, which for many advanced operator theory purposes is the correct generalization of the matrix algebra $M_N(\mathbb C)$.

\bigskip

(3) Our discussion on the compact operators will feature as well some more specialized types of operators, $F(H)\subset B_1(H)\subset B_2(H)\subset K(H)$, with $F(H)$ being the finite rank ones, $B_1(H)$ being the trace class ones, and $B_2(H)$ being the Hilbert-Schmidt ones.

\bigskip

And that is pretty much it, all basic things, that must be known. Of course this will be just the tip of the iceberg, and more of an introduction to modern operator theory.

\bigskip

Getting started now, we would first like to systematically develop the theory of positive operators, and then establish polar decomposition results for the operators $T\in B(H)$. We first have the following result, improving our knowledge from chapter 2:

\index{positive operator}
\index{square root}

\begin{theorem}
For an operator $T\in B(H)$, the following are equivalent:
\begin{enumerate}
\item $<Tx,x>\geq0$, for any $x\in H$.

\item $T$ is normal, and $\sigma(T)\subset[0,\infty)$.

\item $T=S^2$, for some $S\in B(H)$ satisfying $S=S^*$.

\item $T=R^*R$, for some $R\in B(H)$.
\end{enumerate}
If these conditions are satisfied, we call $T$ positive, and write $T\geq0$.
\end{theorem}

\begin{proof}
We have already seen some implications in chapter 2, but the best is to forget the few partial results that we know, and prove everything, as follows:

\medskip

$(1)\implies(2)$ Assuming $<Tx,x>\geq0$, with $S=T-T^*$ we have:
\begin{eqnarray*}
<Sx,x>
&=&<Tx,x>-<T^*x,x>\\
&=&<Tx,x>-<x,Tx>\\
&=&<Tx,x>-\overline{<Tx,x>}\\
&=&0
\end{eqnarray*}

The next step is to use a polarization trick, as follows:
\begin{eqnarray*}
<Sx,y>
&=&<S(x+y),x+y>-<Sx,x>-<Sy,y>-<Sy,x>\\
&=&-<Sy,x>\\
&=&<y,Sx>\\
&=&\overline{<Sx,y>}
\end{eqnarray*}

Thus we must have $<Sx,y>\in\mathbb R$, and with $y\to iy$ we obtain $<Sx,y>\in i\mathbb R$ too, and so $<Sx,y>=0$. Thus $S=0$, which gives $T=T^*$. Now since $T$ is self-adjoint, it is normal as claimed. Moreover, by self-adjointness, we have:
$$\sigma(T)\subset\mathbb R$$

In order to prove now that we have indeed $\sigma(T)\subset[0,\infty)$, as claimed, we must invert $T+\lambda$, for any $\lambda>0$. For this purpose, observe that we have:
\begin{eqnarray*}
<(T+\lambda)x,x>
&=&<Tx,x>+<\lambda x,x>\\
&\geq&<\lambda x,x>\\
&=&\lambda||x||^2
\end{eqnarray*}

But this shows that $T+\lambda$ is injective. In order to prove now the surjectivity, and the boundedness of the inverse, observe first that we have:
\begin{eqnarray*}
Im(T+\lambda)^\perp
&=&\ker(T+\lambda)^*\\
&=&\ker(T+\lambda)\\
&=&\{0\}
\end{eqnarray*}

Thus $Im(T+\lambda)$ is dense. On the other hand, observe that we have:
\begin{eqnarray*}
||(T+\lambda)x||^2
&=&<Tx+\lambda x,Tx+\lambda x>\\
&=&||Tx||^2+2\lambda<Tx,x>+\lambda^2||x||^2\\
&\geq&\lambda^2||x||^2
\end{eqnarray*}

Thus for any vector in the image $y\in Im(T+\lambda)$ we have:
$$||y||\geq\lambda\big|\big|(T+\lambda)^{-1}y\big|\big|$$

As a conclusion to what we have so far, $T+\lambda$ is bijective and invertible as a bounded operator from $H$ onto its image, with the following norm bound:
$$||(T+\lambda)^{-1}||\leq\lambda^{-1}$$

But this shows that $Im(T+\lambda)$ is complete, hence closed, and since we already knew that $Im(T+\lambda)$ is dense, our operator $T+\lambda$ is surjective, and we are done.

\medskip

$(2)\implies(3)$ Since $T$ is normal, and with spectrum contained in $[0,\infty)$, we can use the continuous functional calculus formula for the normal operators from chapter 3, with the function $f(x)=\sqrt{x}$, as to construct a square root $S=\sqrt{T}$. 

\medskip

$(3)\implies(4)$ This is trivial, because we can set $R=S$. 

\medskip

$(4)\implies(1)$ This is clear, because we have the following computation:
$$<R^*Rx,x>
=<Rx,Rx>
=||Rx||^2$$

Thus, we have the equivalences in the statement.
\end{proof}

In analogy with what happens in finite dimensions, where among the positive matrices $A\geq0$ we have the strictly positive ones, $A>0$, given by the fact that the eigenvalues are strictly positive, we have as well a ``strict'' version of the above result, as follows:

\index{strictly positive operator}
\index{square root}

\begin{theorem}
For an operator $T\in B(H)$, the following are equivalent:
\begin{enumerate}
\item $T$ is positive and invertible.

\item $T$ is normal, and $\sigma(T)\subset(0,\infty)$.

\item $T=S^2$, for some $S\in B(H)$ invertible, satisfying $S=S^*$.

\item $T=R^*R$, for some $R\in B(H)$ invertible.
\end{enumerate}
If these conditions are satisfied, we call $T$ strictly positive, and write $T>0$.
\end{theorem}

\begin{proof}
Our claim is that the above conditions (1-4) are precisely the conditions (1-4) in Theorem 4.1, with the assumption ``$T$ is invertible'' added. Indeed:

\medskip

(1) This is clear by definition.

\medskip

(2) In the context of Theorem 4.1 (2), namely when $T$ is normal, and $\sigma(T)\subset[0,\infty)$, the invertibility of $T$, which means $0\notin\sigma(T)$, gives $\sigma(T)\subset(0,\infty)$, as desired.

\medskip

(3) In the context of Theorem 4.1 (3), namely when $T=S^2$, with $S=S^*$, by using the basic properties of the functional calculus for normal operators, the invertibility of $T$ is equivalent to the invertibility of its square root $S=\sqrt{T}$, as desired.

\medskip

(4) In the context of Theorem 4.1 (4), namely when $T=RR^*$, the invertibility of $T$ is equivalent to the invertibility of $R$. This can be either checked directly, or deduced via the equivalence $(3)\iff(4)$ from Theorem 4.1, by using the above argument (3).
\end{proof}

As a subtlety now, we have the following complement to the above result:

\begin{proposition}
For a strictly positive operator, $T>0$, we have
$$<Tx,x>>0\quad,\quad\forall x\neq0$$
but the converse of this fact is not true, unless we are in finite dimensions.
\end{proposition}

\begin{proof}
We have several things to be proved, the idea being as follows:

\medskip

(1) Regarding the main assertion, the inequality can be deduced as follows, by using the fact that the operator $S=\sqrt{T}$ is invertible, and in particular injective:
\begin{eqnarray*}
<Tx,x>
&=&<S^2x,x>\\
&=&<Sx,S^*x>\\
&=&<Sx,Sx>\\
&=&||Sx||^2\\
&>&0
\end{eqnarray*}

(2) In finite dimensions, assuming $<Tx,x>>0$ for any $x\neq0$, we know from Theorem 4.1 that we have $T\geq0$. Thus we have $\sigma(T)\subset[0,\infty)$, and assuming by contradiction $0\in\sigma(T)$, we obtain that $T$ has $\lambda=0$ as eigenvalue, and the corresponding eigenvector $x\neq0$ has the property $<Tx,x>=0$, contradiction. Thus $T>0$, as claimed.

\medskip

(3) Regarding now the counterexample, consider the following operator on $l^2(\mathbb N)$:
$$T=\begin{pmatrix}
1\\
&\frac{1}{2}\\
&&\frac{1}{3}\\
&&&\ddots
\end{pmatrix}$$

This operator $T$ is well-defined and bounded, and we have $<Tx,x>>0$ for any $x\neq0$. However $T$ is not invertible, and so the converse does not hold, as stated.
\end{proof}

With this done, let us discuss now some decomposition results for the bounded operators $T\in B(H)$. We know that any $z\in\mathbb C$ can be written as follows, with $a,b\in\mathbb R$:
$$z=a+ib$$

Also, we know that both the real and imaginary parts $a,b\in\mathbb R$, and more generally any real number $c\in\mathbb R$, can be written as follows, with $r,s\geq0$: 
$$c=r-s$$

Here are the operator theoretic generalizations of these results:

\begin{proposition}
Given an operator $T\in B(H)$, the following happen:
\begin{enumerate}
\item We can write $T=A+iB$, with $A,B\in B(H)$ being self-adjoint.

\item When $T=T^*$, we can write $T=R-S$, with $R,S\in B(H)$ being positive.

\item Thus, we can write any $T$ as a linear combination of $4$ positive elements.
\end{enumerate}
\end{proposition}

\begin{proof}
All this follows from basic spectral theory, as follows:

\medskip

(1) This is something that we have already met in chapter 3, when proving the spectral theorem in its general form, the decomposition formula being as follows:
$$T=\frac{T+T^*}{2}+i\cdot\frac{T-T^*}{2i}$$

(2) This follows from the measurable functional calculus. Indeed, assuming $T=T^*$ we have $\sigma(T)\subset\mathbb R$, so we can use the following decomposition formula on $\mathbb R$:
$$1=\chi_{[0,\infty)}+\chi_{(-\infty,0)}$$

To be more precise, let us multiply by $z$, and rewrite this formula as follows:
$$z=\chi_{[0,\infty)}z-\chi_{(-\infty,0)}(-z)$$

Now by applying these measurable functions to $T$, we obtain as formula as follows, with both the operators $T_+,T_-\in B(H)$ being positive, as desired:
$$T=T_+-T_-$$

(3) This follows indeed by combining the results in (1) and (2) above.
\end{proof}

Going ahead with our decomposition results, another basic thing that we know about complex numbers is that any $z\in\mathbb C$ appears as a real multiple of a unitary:
$$z=re^{it}$$

Finding the correct operator theoretic analogue of this is quite tricky, and this even for the usual matrices $A\in M_N(\mathbb C)$. As a basic result here, we have:

\begin{proposition}
Given an operator $T\in B(H)$, the following happen:
\begin{enumerate}
\item When $T=T^*$ and $||T||\leq1$, we can write $T$ as an average of $2$ unitaries:
$$T=\frac{U+V}{2}$$

\item In the general $T=T^*$ case, we can write $T$ as a rescaled sum of unitaries:
$$T=\lambda(U+V)$$

\item Thus, in general, we can write $T$ as a rescaled sum of $4$ unitaries.
\end{enumerate}
\end{proposition}

\begin{proof}
This follows from the results that we have, as follows:

\medskip

(1) Assuming $T=T^*$ and $||T||\leq1$ we have $1-T^2\geq0$, and the decomposition that we are looking for is as follows, with both the components being unitaries:
$$T=\frac{T+i\sqrt{1-T^2}}{2}+\frac{T-i\sqrt{1-T^2}}{2}$$

To be more precise, the square root can be extracted as in Theorem 4.1 (3), and the check of the unitarity of the components goes as follows:
\begin{eqnarray*}
(T+i\sqrt{1-T^2})(T-i\sqrt{1-T^2})
&=&T^2+(1-T^2)\\
&=&1
\end{eqnarray*}

(2) This simply follows by applying (1) to the operator $T/||T||$.

\medskip

(3) Assuming first $||T||\leq1$, we know from Proposition 4.4 (1) that we can write $T=A+iB$, with $A,B$ being self-adjoint, and satisfying $||A||,||B||\leq1$. Now by applying (1) to both $A$ and $B$, we obtain a decomposition of $T$ as follows:
$$T=\frac{U+V+W+X}{2}$$

In general, we can apply this to the operator $T/||T||$, and we obtain the result.
\end{proof}

All this gets us into the multiplicative theory of the complex numbers, that we will attempt to generalize now. As a first construction, that we would like to generalize to the bounded operator setting, we have the construction of the modulus, as follows:
$$|z|=\sqrt{z\bar{z}}$$

The point now is that we can indeed generalize this construction, as follows:

\index{modulus of operator}
\index{absolute value}
\index{square root}

\begin{proposition}
Given an operator $T\in B(H)$, we can construct a positive operator $|T|\in B(H)$ as follows, by using the fact that $T^*T$ is positive:
$$|T|=\sqrt{T^*T}$$
The square of this operator is then $|T|^2=T^*T$. In the case $H=\mathbb C$, we obtain in this way the usual absolute value of the complex numbers:
$$|z|=\sqrt{z\bar{z}}$$
More generally, in the case where $H=\mathbb C^N$ is finite dimensional, we obtain in this way the usual moduli of the complex matrices $A\in M_N(\mathbb C)$.
\end{proposition}

\begin{proof}
We have several things to be proved, the idea being as follows:

\medskip

(1) The first assertion follows from Theorem 4.1. Indeed, according to (4) there the operator $T^*T$ is indeed positive, and then according to (2) there we can extract the square root of this latter positive operator, by applying to it the function $\sqrt{z}$. 

\medskip

(2) By functional calculus we have then $|T|^2=T^*T$, as desired. 

\medskip

(3) In the case $H=\mathbb C$, we obtain indeed the absolute value of complex numbers.

\medskip

(4) In the case where the space $H$ is finite dimensional, $H=\mathbb C^N$, we obtain indeed the usual moduli of the complex matrices $A\in M_N(\mathbb C)$.
\end{proof}

As a comment here, it is possible to talk as well about $\sqrt{TT^*}$, which is in general different from $\sqrt{T^*T}$. Note that when $T$ is normal, no issue, because we have:
$$TT^*=T^*T\implies\sqrt{TT^*}=\sqrt{T^*T}$$

Regarding now the polar decomposition formula, let us start with a weak version of this statement, regarding the invertible operators, as follows:

\index{polar decomposition}

\begin{theorem}
We have the polar decomposition formula
$$T=U\sqrt{T^*T}$$
with $U$ being a unitary, for any $T\in B(H)$ invertible.
\end{theorem}

\begin{proof}
According to our definition of the modulus, $|T|=\sqrt{T^*T}$, we have:
\begin{eqnarray*}
<|T|x,|T|y>
&=&<x,|T|^2y>\\
&=&<x,T^*Ty>\\
&=&<Tx,Ty>
\end{eqnarray*}

Thus we can define a unitary operator $U\in B(H)$ by the following formula:
$$U(|T|x)=Tx$$

But this formula shows that we have $T=U|T|$, as desired.
\end{proof}

Observe that we have uniqueness in the above result, in what regards the choice of the unitary $U\in B(H)$, due to the fact that we can write this unitary as follows:
$$U=T(\sqrt{T^*T})^{-1}$$

More generally now, we have the following result:

\index{polar decomposition}
\index{partial isometry}

\begin{theorem}
We have the polar decomposition formula
$$T=U\sqrt{T^*T}$$
with $U$ being a partial isometry, for any $T\in B(H)$.
\end{theorem}

\begin{proof}
As before, we have the following equality, for any two vectors $x,y\in H$:
$$<|T|x,|T|y>=<Tx,Ty>$$

We conclude that the following linear application is well-defined, and isometric:
$$U:Im|T|\to Im(T)\quad,\quad 
|T|x\to Tx$$

Now by continuity we can extend this isometry $U$ into an isometry between certain Hilbert subspaces of $H$, as follows:
$$U:\overline{Im|T|}\to\overline{Im(T)}\quad,\quad 
|T|x\to Tx$$

Moreover, we can further extend $U$ into a partial isometry $U:H\to H$, by setting $Ux=0$, for any $x\in\overline{Im|T|}^\perp$, and with this convention, the result follows. 
\end{proof}

\section*{4b. Compact operators}

We have seen so far the basic theory of the bounded operators, in the arbitrary, normal and self-adjoint cases, and in a few other cases of interest. We will keep building on this, with a number of more specialized results, regarding the finite rank operators and compact operators, and other special classes of related operators, namely the trace class operators, and the Hilbert-Schmidt operators. Let us start with a basic definition, as follows:

\index{finite rank operator}

\begin{definition}
An operator $T\in B(H)$ is said to be of finite rank if its image 
$$Im(T)\subset H$$ 
is finite dimensional. The set of such operators is denoted $F(H)$.
\end{definition}

There are many interesting examples of finite rank operators, the most basic ones being the finite rank projections, on the finite dimensional subspaces $K\subset H$. Observe also that in the case where $H$ is finite dimensional, any operator $T\in B(H)$ is automatically of finite rank. In general, this is of course wrong, but we have the following result:

\begin{proposition}
The set of finite rank operators
$$F(H)\subset B(H)$$
is a two-sided $*$-ideal.
\end{proposition}

\begin{proof}
We have several assertions to be proved, the idea being as follows:

\medskip

(1) It is clear from definitions that $F(H)$ is indeed a vector space, with this due to the following formulae, valid for any $S,T\in B(H)$, which are both clear:
$$\dim(Im(S+T))\leq\dim(Im(S))+\dim(Im(T))$$
$$\dim(Im(\lambda T))=\dim(Im(T))$$

(2) Let us prove now that $F(H)$ is stable under $*$. Given $T\in F(H)$, we can regard it as an invertible operator between finite dimensional Hilbert spaces, as follows:
$$T:(\ker T)^\perp\to Im(T)$$

We conclude from this that we have the following dimension equality:
$$\dim((\ker T)^\perp)=\dim(Im(T))$$

Our claim now, in relation with our problem, is that we have equalities as follows:
\begin{eqnarray*}
\dim(Im(T^*))
&=&\dim(\overline{Im(T^*)})\\
&=&\dim((\ker T)^\perp)\\
&=&\dim(Im(T))
\end{eqnarray*}

Indeed, the third equality is the one above, and the second equality is something that we know too, from chapter 2. Now by combining these two equalities we deduce that $Im(T^*)$ is finite dimensional, and so the first equality holds as well. Thus, our equalities are proved, and this shows that we have $T^*\in F(H)$, as desired.

\medskip

(3) Finally, regarding the ideal property, this follows from the following two formulae, valid for any $S,T\in B(H)$, which are once again clear from definitions:
$$\dim(Im(ST))\leq\dim(Im(T))$$
$$\dim(Im(TS))\leq\dim(Im(T))$$

Thus, we are led to the conclusion in the statement.
\end{proof}

Let us discuss now the compact operators, which will be the main topic of discussion, for the present chapter. These are best introduced as follows:

\index{compact operator}

\begin{definition}
An operator $T\in B(H)$ is said to be compact if the closed set
$$\overline{T(B_1)}\subset H$$
is compact, where $B_1\subset H$ is the unit ball. The set of such operators is denoted $K(H)$.
\end{definition}

Equivalently, an operator $T\in B(H)$ is compact when for any sequence $\{x_n\}\subset B_1$, or more generally for any bounded sequence $\{x_n\}\subset H$, the sequence $\{T(x_n)\}$ has a convergence subsequence. We will see later some further criteria of compactness.

\bigskip

In finite dimensions any operator is compact. In general, as a first observation, any finite rank operator is compact. We have in fact the following result:

\begin{proposition}
Any finite rank operator is compact,
$$F(H)\subset K(H)$$
and the finite rank operators are dense inside the compact operators.
\end{proposition}

\begin{proof}
The first assertion is clear, because if $Im(T)$ is finite dimensional, then the following subset is closed and bounded, and so it is compact:
$$\overline{T(B_1)}\subset Im(T)$$

Regarding the second assertion, let us pick a compact operator $T\in K(H)$, and a number $\varepsilon>0$. By compactness of $T$ we can find a finite set $S\subset B_1$ such that:
$$T(B_1)\subset\bigcup_{x\in S}B_\varepsilon(Tx)$$

Consider now the orthogonal projection $P$ onto the following finite dimensional space:
$$E=span\left(Tx\Big|x\in S\right)$$

Since the set $S$ is finite, this space $E$ is finite dimensional, and so $P$ is of finite rank, $P\in F(H)$. Now observe that for any norm one $y\in H$ and any $x\in S$ we have:
\begin{eqnarray*}
||Ty-Tx||^2
&=&||Ty-PTx||^2\\
&=&||Ty-PTy+PTy-PTx||^2\\
&=&||Ty-PTy||^2+||PTx-PTy||^2
\end{eqnarray*}

Now by picking $x\in S$ such that the ball $B_\varepsilon(Tx)$ covers the point $Ty$, we conclude from this that we have the following estimate:
$$||Ty-PTy||\leq||Ty-Tx||\leq\varepsilon$$ 

Thus we have $||T-PT||\leq\varepsilon$, which gives the density result.
\end{proof}

Quite remarkably, the set of compact operators is closed, and we have:

\begin{theorem}
The set of compact operators
$$K(H)\subset B(H)$$
is a closed two-sided $*$-ideal.
\end{theorem}

\begin{proof}
We have several assertions here, the idea being as follows:

\medskip

(1) It is clear from definitions that $K(H)$ is indeed a vector space, with this due to the following formulae, valid for any $S,T\in B(H)$, which are both clear:
$$(S+T)(B_1)\subset S(B_1)+T(B_1)$$
$$(\lambda T)(B_1)=|\lambda|\cdot T(B_1)$$

(2) In order to prove now that $K(H)$ is closed, assume that a sequence $T_n\in K(H)$ converges to $T\in B(H)$. Given $\varepsilon>0$, let us pick $N\in\mathbb N$ such that:
$$||T-T_N||\leq\varepsilon$$

By compactness of $T_N$ we can find a finite set $S\subset B_1$ such that:
$$T_N(B_1)\subset\bigcup_{x\in S}B_\varepsilon(T_Nx)$$

We conclude that for any $y\in B_1$ there exists $x\in S$ such that:
\begin{eqnarray*}
||Ty-Tx||
&\leq&||Ty-T_Ny||+||T_Ny-T_Nx||+||T_Nx-Tx||\\
&\leq&\varepsilon+\varepsilon+\varepsilon\\
&=&3\varepsilon
\end{eqnarray*}

Thus, we have an inclusion as follows, with $S\subset B_1$ being finite:
$$T(B_1)\subset\bigcup_{x\in S}B_{3\varepsilon}(Tx)$$

But this shows that our limiting operator $T$ is compact, as desired.

\medskip

(3) Regarding the fact that $K(H)$ is stable under involution, this follows from Proposition 4.10, Proposition 4.12 and (2). Indeed, by using Proposition 4.12, given $T\in K(H)$ we can write it as a limit of finite rank operators, as follows:
$$T=\lim_{n\to\infty}T_n$$

Now by applying the adjoint, we obtain that we have as well:
$$T^*=\lim_{n\to\infty}T_n^*$$

We know from Proposition 4.10 that the operators $T_n^*$ are of finite rank, and so compact by Proposition 4.12, and by using (2) we obtain that $T^*$ is compact too, as desired.

\medskip

(4) Finally, regarding the ideal property, this follows from the following two formulae, valid for any $S,T\in B(H)$, which are once again clear from definitions:
$$(ST)(B_1)=S(T(B_1))$$
$$(TS)(B_1)\subset ||S||\cdot T(B_1)$$

Thus, we are led to the conclusion in the statement.
\end{proof}

Here is now a second key result regarding the compact operators:

\begin{theorem}
A bounded operator $T\in B(H)$ is compact precisely when
$$Te_n\to0$$
for any orthonormal system $\{e_n\}\subset H$.
\end{theorem}

\begin{proof}
We have two implications to be proved, the idea being as follows:

\medskip

``$\implies$'' Assume that $T$ is compact. By contradiction, assume $Te_n\not\to0$. This means that there exists $\varepsilon>0$ and a subsequence satisfying $||Te_{n_k}||>\varepsilon$, and by replacing $\{e_n\}$ with this subsequence, we can assume that the following holds, with $\varepsilon>0$:
$$||Te_n||>\varepsilon$$

Since $T$ was assumed to be compact, and the sequence $\{e_n\}$ is bounded, a certain subsequence $\{Te_{n_k}\}$ must converge. Thus, by replacing once again $\{e_n\}$ with a subsequence, we can assume that the following holds, with $x\neq0$:
$$Te_n\to x$$

But this is a contradiction, because we obtain in this way:
\begin{eqnarray*}
<x,x>
&=&\lim_{n\to\infty}<Te_n,x>\\
&=&\lim_{n\to\infty}<e_n,T^*x>\\
&=&0
\end{eqnarray*}

Thus our assumption $Te_n\not\to0$ was wrong, and we obtain the result.

\medskip

``$\Longleftarrow$'' Assume $Te_n\to0$, for any orthonormal system $\{e_n\}\subset H$. In order to prove that $T$ is compact, we use the various results established above, which show that this is the same as proving that $T$ is in the closure of the space of finite rank operators:
$$T\in\overline{F(H)}$$

We do this by contradiction. So, assume that the above is wrong, and so that there exists $\varepsilon>0$ such that the following holds:
$$S\in F(H)\implies||T-S||>\varepsilon$$

As a first observation, by using $S=0$ we obtain $||T||>\varepsilon$. Thus, we can find a norm one vector $e_1\in H$ such that the following holds:
$$||Te_1||>\varepsilon$$

Our claim, which will bring the desired contradiction, is that we can construct by recurrence vectors $e_1,\ldots,e_n$ such that the following holds, for any $i$: 
$$||Te_i||>\varepsilon$$

Indeed, assume that we have constructed such vectors $e_1,\ldots,e_n$. Let $E\subset H$ be the linear space spanned by these vectors, and let us set:
$$P=Proj(E)$$

Since the operator $TP$ has finite rank, our assumption above shows that we have:
$$||T-TP||>\varepsilon$$

Thus, we can find a vector $x\in H$ such that the following holds:
$$||(T-TP)x||>\varepsilon$$

We have then $x\not\in E$, and so we can consider the following nonzero vector:
$$y=(1-P)x$$

With this nonzero vector $y$ constructed, in this way, now let us set:
$$e_{n+1}=\frac{y}{||y||}$$

This vector $e_{n+1}$ is then orthogonal to $E$, has norm one, and satisfies:
$$||Te_{n+1}||\geq||y||^{-1}\varepsilon\geq\varepsilon$$

Thus we are done with our construction by recurrence, and this contradicts our assumption that $Te_n\to0$, for any orthonormal system $\{e_n\}\subset H$, as desired.
\end{proof}

Summarizing, we have so far a number of results regarding the compact operators, in analogy with what we know about the usual complex matrices. Let us discuss now the spectral theory of the compact operators. We first have the following result:

\index{eigenvalue}

\begin{proposition}
Assuming that $T\in B(H)$, with $\dim H=\infty$, is compact and self-adjoint, the following happen:
\begin{enumerate}
\item The eigenvalues of $T$ form a sequence $\lambda_n\to0$.

\item All eigenvalues $\lambda_n\neq0$ have finite multiplicity.
\end{enumerate}
\end{proposition}

\begin{proof}
We prove both the assertions at the same time. For this purpose, we fix a number $\varepsilon>0$, we consider all the eigenvalues satisfying $|\lambda|\geq\varepsilon$, and for each such eigenvalue we consider the corresponding eigenspace $E_\lambda\subset H$. Let us set:
$$E=span\left(E_\lambda\,\Big|\,|\lambda|\geq\varepsilon\right)$$ 

Our claim, which will prove both (1) and (2), is that this space $E$ is finite dimensional. In now to prove now this claim, we can proceed as follows:

\medskip

(1) We know that we have $E\subset Im(T)$. Our claim is that we have:
$$\bar{E}\subset Im(T)$$

Indeed, assume that we have a sequence $g_n\in E$ which converges, $g_n\to g\in\bar{E}$. Let us write $g_n=Tf_n$, with $f_n\in H$. By definition of $E$, the following condition is satisfied:
$$h\in E\implies||Th||\geq\varepsilon||h||$$

Now since the sequence $\{g_n\}$ is Cauchy we obtain from this that the sequence $\{f_n\}$ is Cauchy as well, and with $f_n\to f$ we have $Tf_n\to Tf$, as desired.

\medskip

(2) Consider now the projection $P\in B(H)$ onto the closure $\bar{E}$ of the above vector space $E$. The composition $PT$ is then as follows, surjective on its target:
$$PT:H\to\bar{E}$$

On the other hand since $T$ is compact so must be $PT$, and if follows from this that the space $\bar{E}$ is finite dimensional. Thus $E$ itself must be finite dimensional too, and as explained in the beginning of the proof, this gives (1) and (2), as desired.
\end{proof}

In order to construct now eigenvalues, we will need:

\begin{proposition}
If $T$ is compact and self-adjoint, one of the numbers 
$$||T||\ ,\ -||T||$$
must be an eigenvalue of $T$.
\end{proposition}

\begin{proof}
We know from the spectral theory of the self-adjoint operators that the spectral radius $||T||$ of our operator $T$ is attained, and so one of the numbers $||T||,-||T||$ must be in the spectrum. In order to prove now that one of these numbers must actually appear as an eigenvalue, we must use the compactness of $T$, as follows:

\medskip

(1) First, we can assume $||T||=1$. By functional calculus this implies $||T^3||=1$ too, and so we can find a sequence of norm one vectors $x_n\in H$ such that:
$$|<T^3x_n,x_n>|\to1$$

By using our assumption $T=T^*$, we can rewrite this formula as follows:
$$|<T^2x_n,Tx_n>|\to1$$

Now since $T$ is compact, and $\{x_n\}$ is bounded, we can assume, up to changing the sequence $\{x_n\}$ to one of its subsequences, that the sequence $Tx_n$ converges:
$$Tx_n\to y$$

Thus, the convergence formula found above reformulates as follows, with $y\neq0$:
$$|<Ty,y>|=1$$

(2) Our claim now, which will finish the proof, is that this latter formula implies $Ty=\pm y$. Indeed, by using Cauchy-Schwarz and $||T||=1$, we have:
$$|<Ty,y>|\leq||Ty||\cdot||y||\leq1$$

We know that this must be an equality, so $Ty,y$ must be proportional. But since $T$ is self-adjoint the proportionality factor must be $\pm1$, and so we obtain, as claimed:
$$Ty=\pm y$$

Thus, we have constructed an eigenvector for $\lambda=\pm1$, as desired.
\end{proof}

We can further build on the above results in the following way:

\begin{proposition}
If $T$ is compact and self-adjoint, there is an orthogonal basis of $H$ made of eigenvectors of $T$.
\end{proposition}

\begin{proof}
We use Proposition 4.15. According to the results there, we can arrange the nonzero eigenvalues of $T$, taken with multiplicities, into a sequence $\lambda_n\to0$. Let $y_n\in H$ be the corresponding eigenvectors, and consider the following space:
$$E=\overline{span(y_n)}$$

The result follows then from the following observations:

\medskip

(1) Since we have $T=T^*$, both $E$ and its orthogonal $E^\perp$ are invariant under $T$. 

\medskip

(2) On the space $E$, our operator $T$ is by definition diagonal.

\medskip

(3) On the space $E^\perp$, our claim is that we have $T=0$. Indeed, assuming that the restriction $S=T_{E^\perp}$ is nonzero, we can apply Proposition 4.16 to this restriction, and we obtain an eigenvalue for $S$, and so for $T$, contradicting the maximality of $E$.
\end{proof}

With the above results in hand, we can now formulate a first spectral theory result for compact operators, which closes the discussion in the self-adjoint case:

\begin{theorem}
Assuming that $T\in B(H)$, with $\dim H=\infty$, is compact and self-adjoint, the following happen:
\begin{enumerate}
\item The spectrum $\sigma(T)\subset\mathbb R$ consists of a sequence $\lambda_n\to0$.

\item All spectral values $\lambda\in\sigma(T)-\{0\}$ are eigenvalues.

\item All eigenvalues $\lambda\in\sigma(T)-\{0\}$ have finite multiplicity.

\item There is an orthogonal basis of $H$ made of eigenvectors of $T$.
\end{enumerate}
\end{theorem}

\begin{proof}
This follows from the various results established above:

\medskip

(1) In view of Proposition 4.15 (1), this will follow from (2) below.

\medskip

(2) Assume that $\lambda\neq0$ belongs to the spectrum $\sigma(T)$, but is not an eigenvalue. By using Proposition 4.17, let us pick an orthonormal basis $\{e_n\}$ of $H$ consisting of eigenvectors of $T$, and then consider the following operator:
$$Sx=\sum_n\frac{<x,e_n>}{\lambda_n-\lambda}\,e_n$$

Then $S$ is an inverse for $T-\lambda$, and so we have $\lambda\notin\sigma(T)$, as desired.

\medskip

(3) This is something that we know, from Proposition 4.15 (2).

\medskip

(4) This is something that we know too, from Proposition 4.17.
\end{proof}

Finally, we have the following result, regarding the general case:

\index{compact operator}
\index{singular values}
\index{diagonalization}
\index{eigenvalue}
\index{polar decomposition}

\begin{theorem}
The compact operators $T\in B(H)$, with $\dim H=\infty$, are the operators of the following form, with $\{e_n\}$, $\{f_n\}$ being orthonormal families, and with $\lambda_n\searrow0$:
$$T(x)=\sum_n\lambda_n<x,e_n>f_n$$
The numbers $\lambda_n$, called singular values of $T$, are the eigenvalues of $|T|$. In fact, the polar decomposition of $T$ is given by $T=U|T|$, with
$$|T|(x)=\sum_n\lambda_n<x,e_n>e_n$$
and with $U$ being given by $Ue_n=f_n$, and $U=0$ on the complement of $span(e_i)$.
\end{theorem}

\begin{proof}
This basically follows from Theorem 4.8 and Theorem 4.18, as follows:

\medskip

(1) Given two orthonormal families $\{e_n\}$, $\{f_n\}$, and a sequence of real numbers $\lambda_n\searrow0$, consider the linear operator given by the formula in the statement, namely:
$$T(x)=\sum_n\lambda_n<x,e_n>f_n$$

Our first claim is that $T$ is bounded. Indeed, when assuming $|\lambda_n|\leq\varepsilon$ for any $n$, which is something that we can do if we want to prove that $T$ is bounded, we have:
\begin{eqnarray*}
||T(x)||^2
&=&\left|\sum_n\lambda_n<x,e_n>f_n\right|^2\\
&=&\sum_n|\lambda_n|^2|<x,e_n>|^2\\
&\leq&\varepsilon^2\sum_n|<x,e_n>|^2\\
&\leq&\varepsilon^2||x||^2
\end{eqnarray*}

(2) The next observation is that this operator is indeed compact, because it appears as the norm limit, $T_N\to T$, of the following sequence of finite rank operators:
$$T_N=\sum_{n\leq N}\lambda_n<x,e_n>f_n$$

(3) Regarding now the polar decomposition assertion, for the above operator, this follows once again from definitions. Indeed, the adjoint is given by:
$$T^*(x)=\sum_n\lambda_n<x,f_n>e_n$$

Thus, when composing $T^*$ with $T$, we obtain the following operator:
$$T^*T(x)=\sum_n\lambda_n^2<x,e_n>e_n$$

Now by extracting the square root, we obtain the formula in the statement, namely:
$$|T|(x)=\sum_n\lambda_n<x,e_n>e_n$$

(4) Conversely now, assume that $T\in B(H)$ is compact. Then $T^*T$, which is self-adjoint, must be compact as well, and so by Theorem 4.18 we have a formula as follows, with $\{e_n\}$ being a certain orthonormal family, and with $\lambda_n\searrow0$:
$$T^*T(x)=\sum_n\lambda_n^2<x,e_n>e_n$$

By extracting the square root we obtain the formula of $|T|$ in the statement, and then by setting $U(e_n)=f_n$ we obtain a second orthonormal family, $\{f_n\}$, such that:
$$T(x)
=U|T|
=\sum_n\lambda_n<x,e_n>f_n$$

Thus, our compact operator $T\in B(H)$ appears indeed as in the statement.
\end{proof}

As a technical remark here, it is possible to slightly improve a part of the above statement. Consider indeed an operator of the following form, with $\{e_n\}$, $\{f_n\}$ being orthonormal families as before, and with $\lambda_n\to0$ being now complex numbers:
$$T(x)=\sum_n\lambda_n<x,e_n>f_n$$

Then the same proof as before shows that $T$ is compact, and that the polar decomposition of $T$ is given by $T=U|T|$, with the modulus $|T|$ being as follows:
$$|T|(x)=\sum_n|\lambda_n|<x,e_n>e_n$$

As for the partial isometry $U$, this is given by $Ue_n=w_nf_n$, and $U=0$ on the complement of $span(e_i)$, where $w_n\in\mathbb T$ are such that $\lambda_n=|\lambda_n|w_n$.

\section*{4c. Trace class operators}

We have not talked so far about the trace of operators $T\in B(H)$, in analogy with the trace of the usual matrices $M\in M_N(\mathbb C)$. This is because the trace can be finite or infinite, or even not well-defined, and we will discuss this now. Let us start with:

\index{trace of operators}

\begin{proposition}
Given a positive operator $T\in B(H)$, the quantity
$$Tr(T)=\sum_n<Te_n,e_n>\in[0,\infty]$$
is indpendent on the choice of an orthonormal basis $\{e_n\}$.
\end{proposition}

\begin{proof}
If $\{f_n\}$ is another orthonormal basis, we have:
\begin{eqnarray*}
\sum_n<Tf_n,f_n>
&=&\sum_n<\sqrt{T}f_n,\sqrt{T}f_n>\\
&=&\sum_n||\sqrt{T}f_n||^2\\
&=&\sum_{mn}|<\sqrt{T}f_n,e_m>|^2\\
&=&\sum_{mn}|<T^{1/4}f_n,T^{1/4}e_m>|^2
\end{eqnarray*}

Since this quantity is symmetric in $e,f$, this gives the result.
\end{proof}

We can now introduce the trace class operators, as follows:

\index{trace class operator}

\begin{definition}
An operator $T\in B(H)$ is said to be of trace class if:
$$Tr|T|<\infty$$
The set of such operators, also called integrable, is denoted $B_1(H)$.
\end{definition}

In finite dimensions, any operator is of course of trace class. In arbitrary dimension, finite or not, we first have the following result, regarding such operators:

\index{singular values}
\index{diagonalization}

\begin{proposition}
Any finite rank operator is of trace class, and any trace class operator is compact, so that we have embeddings as follows:
$$F(H)\subset B_1(H)\subset K(H)$$
Moreover, for any compact operator $T\in K(H)$ we have the formula
$$Tr|T|=\sum_n\lambda_n$$
where $\lambda_n\geq0$ are the singular values, and so $T\in B_1(H)$ precisely when $\sum_n\lambda_n<\infty$.
\end{proposition}

\begin{proof}
We have several assertions here, the idea being as follows:

\medskip

(1) If $T$ is of finite rank, it is clearly of trace class.

\medskip

(2) In order to prove now the second assertion, assume first that $T>0$ is of trace class. For any orthonormal basis $\{e_n\}$ we have:
\begin{eqnarray*}
\sum_n||\sqrt{T}e_n||^2
&=&\sum_n<Te_n,e_n>\\
&\leq&Tr(T)\\
&<&\infty
\end{eqnarray*}

But this shows that we have a convergence as follows:
$$\sqrt{T}e_n\to0$$

Thus the operator $\sqrt{T}$ is compact. Now since the compact operators form an ideal, it follows that $T=\sqrt{T}\cdot\sqrt{T}$ is compact as well, as desired.

\medskip

(3) In order to prove now the second assertion in general, assume that $T\in B(H)$ is of trace class. Then $|T|$ is also of trace class, and so compact by (2), and since we have $T=U|T|$ by polar decomposition, it follows that $T$ is compact too.

\medskip

(4) Finally, in order to prove the last assertion, assume that $T$ is compact. The singular value decomposition of $|T|$, from Theorem 4.19, is then as follows:
$$|T|(x)=\sum_n\lambda_n<x,e_n>e_n$$

But this gives the formula for $Tr|T|$ in the statement, and proves the last assertion.
\end{proof}

Here is a useful reformulation of the above result, or rather of the above result coupled with Theorem 4.19, without reference to compact operators:

\begin{theorem}
The trace class operators are precisely the operators of the form
$$|T|(x)=\sum_n\lambda_n<x,e_n>f_n$$
with $\{e_n\}$, $\{f_n\}$ being orthonormal systems, and with $\lambda\searrow0$ being a sequence satisfying:
$$\sum_n\lambda_n<\infty$$
Moreover, for such an operator we have the following estimate:
$$|Tr(T)|\leq Tr|T|=\sum_n\lambda_n$$ 
\end{theorem}

\begin{proof}
This follows indeed from Proposition 4.22, or rather for step (4) in the proof of Proposition 4.22, coupled with Theorem 4.19.
\end{proof}

Next, we have the following result, which comes as a continuation of Proposition 4.22, and is our central result here, regarding the trace class operators:

\index{trace class operator}

\begin{theorem}
The space of trace class operators, which appears as an intermediate space between the finite rank operators and the compact operators,
$$F(H)\subset B_1(H)\subset K(H)$$
is a two-sided $*$-ideal of $K(H)$. The following is a Banach space norm on $B_1(H)$, 
$$||T||_1=Tr|T|$$
satisfying $||T||\leq||T||_1$, and for $T\in B_1(H)$ and $S\in B(H)$ we have:
$$||ST||_1\leq||S||\cdot||T||_1$$
Also, the subspace $F(H)$ is dense inside $B_1(H)$, with respect to this norm.
\end{theorem}

\begin{proof}
There are several assertions here, the idea being as follows:

\medskip

(1) In order to prove that $B_1(H)$ is a linear space, and that $||T||_1=Tr|T|$ is a norm on it, the only non-trivial point is that of proving the following inequality:
$$Tr|S+T|\leq Tr|S|+Tr|T|$$

For this purpose, consider the polar decompositions of these operators:
$$S=U|S|\quad,\quad
T=V|T|\quad,\quad 
S+T=W|S+T|$$

Given an orthonormal basis $\{e_n\}$, we have the following formula:
\begin{eqnarray*}
Tr|S+T|
&=&\sum_n<|S+T|e_n,e_n>\\
&=&\sum_n<W^*(S+T)e_n,e_n>\\
&=&\sum_n<W^*U|S|e_n,e_n>+\sum_n<W^*V|T|e_n,e_n>
\end{eqnarray*}

The point now is that the first sum can be estimated as follows:
\begin{eqnarray*}
&&\sum_n<W^*U|S|e_n,e_n>\\
&=&\sum_n<\sqrt{|S|}e_n,\sqrt{|S|}U^*We_n>\\
&\leq&\sum_n\Big|\Big|\sqrt{|S|}e_n\Big|\Big|\cdot\Big|\Big|\sqrt{|S|}U^*We_n\Big|\Big|\\
&\leq&\sqrt{\sum_n\Big|\Big|\sqrt{|S|}e_n\Big|\Big|^2}\cdot\sqrt{\sum_n\Big|\Big|\sqrt{|S|}U^*We_n\Big|\Big|^2}
\end{eqnarray*}

In order to estimate the terms on the right, we can proceed as follows:
\begin{eqnarray*}
\sum_n\Big|\Big|\sqrt{|S|}U^*We_n\Big|\Big|^2
&=&\sum_n<W^*U|S|U^*We_n,e_n>\\
&=&Tr(W^*U|S|U^*W)\\
&\leq&Tr(U|S|U^*)\\
&\leq&Tr(|S|)
\end{eqnarray*}

The second sum in the above formula of $Tr|S+T|$ can be estimated in the same way, and in the end we obtain, as desired:
$$Tr|S+T|\leq Tr|S|+Tr|T|$$

(2) The estimate $||T||\leq||T||_1$ can be established as follows:
\begin{eqnarray*}
||T||
&=&\big|\big||T|\big|\big|\\
&=&\sup_{||x||=1}<|T|x,x>\\
&\leq&Tr|T|
\end{eqnarray*}

(3) The fact that $B_1(H)$ is indeed a Banach space follows by constructing a limit for any Cauchy sequence, by using the singular value decomposition.

\medskip

(4) The fact that $B_1(H)$ is indeed closed under the involution follows from:
\begin{eqnarray*}
Tr(T^*)
&=&\sum_n<T^*e_n,e_n>\\
&=&\sum_n<e_n,Te_N>\\
&=&\overline{Tr(T)}
\end{eqnarray*}

(5) In order to prove now the ideal property of $B_1(H)$, we use the standard fact, that we know from Proposition 4.5, that any bounded operator $T\in B(H)$ can be written as a linear combination of 4 unitary operators, as follows:
$$T=\lambda_1U_1+\lambda_2U_2+\lambda_3U_3+\lambda_4U_4$$

Indeed, by taking the real and imaginary part we can first write $T$ as a linear combination of 2 self-adjoint operators, and then by functional calculus each of these 2 self-adjoint operators can be written as a linear linear combination of 2 unitary operators.

\medskip

(6) With this trick in hand, we can now prove the ideal property of $B_1(H)$. Indeed, it is enough to prove that we have:
$$T\in B_1(H),U\in U(H)\implies UT,TU\in B_1(H)$$

But this latter result follows by using the polar decomposition theorem.

\medskip

(7) With a bit more care, we obtain from this the estimate $||ST||_1\leq||S||\cdot||T||_1$ from the statement. As for the last assertion, this is clear as well.
\end{proof}

This was for the basic theory of the trace class operators. Much more can be said, and we refer here to the literature, such as Lax \cite{lax}. In what concerns us, we will be back to these operators later in this book, in Part III, when discussing operator algebras.

\section*{4d. Hilbert-Schmidt operators}

As a last topic of this chapter, let us discuss yet another important class of operators, namely the Hilbert-Schmidt ones. These operators, that we will need on several key occasions in what follows, when talking operator algebras, are introduced as follows:

\index{Hilbert-Schmidt operator}

\begin{definition}
An operator $T\in B(H)$ is said to be Hilbert-Schmidt if:
$$Tr(T^*T)<\infty$$
The set of such operators is denoted $B_2(H)$.
\end{definition}

As before with other sets of operators, in finite dimensions we obtain in this way all the operators. In general, we have the following result, regarding such operators:

\index{singular values}

\begin{theorem}
The space $B_2(H)$ of Hilbert-Schmidt operators, which appears as an intermediate space between the trace class operators and the compact operators,
$$F(H)\subset B_1(H)\subset B_2(H)\subset K(H)$$
is a two-sided $*$-ideal of $K(H)$. This ideal has the property
$$S,T\in B_2(H)\implies ST\in B_1(H)$$
and conversely, each $T\in B_1(H)$ appears as product of two operators in $B_2(H)$. In terms of the singular values $(\lambda_n)$, the Hilbert-Schmidt operators are characterized by:
$$\sum_n\lambda_n^2<\infty$$
Also, the following formula, whose output is finite by Cauchy-Schwarz,
$$<S,T>=Tr(ST^*)$$
defines a scalar product of $B_2(H)$, making it a Hilbert space.
\end{theorem}

\begin{proof}
All this is quite standard, from the results that we have already, and more specifically from the singular value decomposition theorem, and its applications. To be more precise, the proof of the various assertions goes as follows:

\medskip

(1) First of all, the fact that the space of Hilbert-Schmidt operators $B_2(H)$ is stable under taking sums, and so is a vector space, follows from:
\begin{eqnarray*}
(S+T)^*(S+T)
&\leq&(S+T)^*(S+T)+(S-T)^*(S-T)\\
&=&(S^*+T^*)(S+T)+(S^*-T^*)(S-T)\\
&=&2(S^*S+T^*T)
\end{eqnarray*}

Regarding now multiplicative properties, we can use here the following inequality:
$$(ST)^*(ST)
=T^*S^*ST
\leq||S||^2T^*T$$

Thus, the space $B_2(H)$ is a two-sided $*$-ideal of $K(H)$, as claimed.

\medskip

(2) In order to prove now that the product of any two Hilbert-Schmidt operators is a trace class operator, we can use the following formula, which is elementary:
$$S^*T=\sum_{k=1}^4i^k(S-iT)^*(S-iT)$$

Conversely, given an arbitrary trace class operator $T\in B_1(H)$, we have:
$$T\in B_1(H)
\implies|T|\in B_1(H)
\implies\sqrt{|T|}\in B_2(H)$$

Thus, by using the polar decomposition $T=U|T|$, we obtain the following decomposition for $T$, with both components being Hilbert-Schmidt operators:
$$T
=U|T|
=U\sqrt{|T|}\cdot\sqrt{|T|}$$

(3) The condition for the singular values is clear.

\medskip

(4) The fact that we have a scalar product is clear as well.

\medskip

(5) The proof of the completness property is routine as well.
\end{proof}

We have as well the following key result, regarding the Hilbert-Schmidt operators:

\begin{theorem}
We have the following formula,
$$Tr(ST)=Tr(TS)$$
valied for any Hilbert-Schmidt operators $S,T\in B_2(H)$.
\end{theorem}

\begin{proof}
We can prove this in two steps, as follows:

\medskip

(1) Assume first that $|S|$ is trace class. Consider the polar decomposition $S=U|S|$, and choose an orthonormal basis $\{x_i\}$ for the image of $U$, suitably extended to an orthonormal basis of $H$. We have then the following computation, as desired:
\begin{eqnarray*}
Tr(ST)
&=&\sum_i<U|S|Tx_i,x_i>\\
&=&\sum_i<|S|TUU^*x_i,U^*x_i>\\
&=&Tr(|S|TU)\\
&=&Tr(TU|S|)\\
&=&Tr(TS)
\end{eqnarray*}

(2) Assume now that we are in the general case, where $S$ is only assumed to be Hilbert-Schmidt. For any finite rank operator $S'$ we have then:
\begin{eqnarray*}
|Tr(ST)-Tr(TS)|
&=&|Tr((S-S')T)-Tr(T(S-S'))|\\
&\leq&2||S-S'||_2\cdot||T||_2
\end{eqnarray*}

Thus by choosing $S'$ with $||S-S'||_2\to0$, we obtain the result.
\end{proof}

This was for the basic theory of bounded operators on a Hilbert space, $T\in B(H)$. In the remainder of this book we will be rather interested in the operator algebras $A\subset B(H)$ that these operators can form. This is of course related to operator theory, because we can, at least in theory, take $A=<T>$, and then study $T$ via the properties of $A$. Actually, this is something that we already did a few times, when doing spectral theory, and notably when talking about functional calculus for normal operators.

\bigskip

For further operator theory, however, nothing beats a good operator theory book, and various ad-hoc methods, depending on the type of operators involved, and especially, on what you want to do with them. As before, in relation with topics to be later discussed in this book, we recommend here the books of Lax \cite{lax} and Blackadar \cite{bla}.

\bigskip

Let us mention as well that there is a lot of interesting theory regarding the unbounded operators $T\in\mathcal L(H)$ too, which is something quite technical, and here once again, we warmly recommend a good operator theory book. In addition, we recommend as well a good PDE book, because most of the questions making appear unbounded operators usually have PDE formulations as well, which are extremely efficient. 

\section*{4e. Exercises} 

There has been a lot of theory in this chapter, with some of the things not really explained in great detail, and we have several exercises about all this. First comes:

\begin{exercise}
Try to find the best operator theoretic analogue of the formula
$$z=re^{it}$$
for the complex numbers, telling us that any number is a real multiple of a unitary.
\end{exercise}

As explained in the above, a weak analogue of this holds, stating that any operator is a linear combination of 4 unitaries. The problem is that of improving this.

\begin{exercise}
Work out a few explicit examples of the polar decomposition formula
$$T=U\sqrt{T^*T}$$
with, if possible, a non-trivial computation for the square root.
\end{exercise}

This is actually something quite tricky, even for the usual matrices. So, as a preliminary exercise here, have some fun with the $2\times2$ matrices.

\begin{exercise}
Look up the various extra general properties of the sets of finite rank, trace class, Hilbert-Schmidt and compact operators,
$$F(H)\subset B_1(H)\subset B_2(H)\subset K(H)$$
coming in addition to what has been said above, about such operators. 
\end{exercise}

This is of course quite vague, and, as good news, it is not indicated either if you should just come with a list of such properties, or with a list of such properties coming with complete proofs. Up to you here, and the more the better.

\part{Operator algebras}

\ \vskip50mm

\begin{center}
{\em There was something in the air that night

The stars were bright, Fernando

They were shining there for you and me

For liberty, Fernando}
\end{center}

\chapter{Operator algebras}

\section*{5a. Normed algebras}

We have seen that the study of the bounded operators $T\in B(H)$ often leads to the consideration of the algebras $<T>\subset B(H)$ generated by such operators, the idea being that the study of $A=<T>$ can lead to results about $T$ itself. In the remainder of this book we focus on the study of such algebras $A\subset B(H)$. Before anything, we should mention that there are countless ways of getting introduced to operator algebras, depending on motivations and taste, with the available books including:

\bigskip

(1) The old book of von Neumann \cite{vn4}, which started everything. This is a very classical book, with mathematical physics content, written at times when mathematics and physics were starting to part ways. A great book, still enjoyable nowadays.

\bigskip

(2) Various post-war treatises, such as Dixmier \cite{dix}, Kadison-Ringrose \cite{kri}, Str\u atil\u a-Zsid\'o \cite{szs} and Takesaki \cite{tak}. As a warning, however, these books are purely mathematical. Also, they sometimes avoid deep results of von Neumann and Connes.

\bigskip

(3) More recent books, including Arveson \cite{arv}, Blackadar \cite{bla}, Brown-Ozawa \cite{boz}, Connes \cite{co3}, Davidson \cite{dav}, Jones \cite{jo6}, Murphy \cite{mur}, Pedersen \cite{ped} and Sakai \cite{sak}. These are well-concieved one-volume books, written with various purposes in mind.

\bigskip

Our presentation below is inspired by Blackadar \cite{bla}, Connes \cite{co3}, Jones \cite{jo6}, but is yet another type of beast, often insisting on probabilistic aspects. But probably enough talking, more on this later, and let us get to work. We are interested in the study of the algebras of bounded operators $A\subset B(H)$. Let us start our discussion with the following broad definition, obtained by imposing the ``minimal'' set of reasonable axioms:

\index{operator algebra}

\begin{definition}
An operator algebra is an algebra of bounded operators $A\subset B(H)$ which contains the unit, is closed under taking adjoints, 
$$T\in A\implies T^*\in A$$
and is closed as well under the norm.
\end{definition}

Here, as in the previous chapters, $H$ is an arbitrary Hilbert space, with the case that we are mostly interested in being the separable one. By separable we mean having a countable orthonormal basis, $\{e_i\}_{i\in I}$ with $I$ countable, and such a space is of course unique. The simplest model is the space $l^2(\mathbb N)$, but in practice, we are particularly interested in the spaces of the form $H=L^2(X)$, which are separable too, but with the basis $\{e_i\}_{i\in\mathbb N}$ and the subsequent identification $H\simeq l^2(\mathbb N)$ being not necessarily very explicit.

\bigskip

Also as in the previous chapters, $B(H)$ is the algebra of linear operators $T:H\to H$ which are bounded, in the sense that the norm $||T||=\sup_{||x||=1}||Tx||$ is finite. This algebra has an involution $T\to T^*$, with the adjoint operator $T^*\in B(H)$ being defined by the formula $<Tx,y>=<x,T^*y>$, and in the above definition, the assumption $T\in A\implies T^*\in A$ refers to this involution. Thus, $A$ must be a $*$-algebra.

\bigskip

As a first result now regarding the operator algebras, in relation with the normal operators, where most of the non-trivial results that we have so far are, we have:

\index{normal operator}
\index{continuous calculus}

\begin{theorem}
The operator algebra $<T>\subset B(H)$ generated by a normal operator $T\in B(H)$ appears as an algebra of continuous functions,
$$<T>=C(\sigma(T))$$
where $\sigma(T)\subset\mathbb C$ denotes as usual the spectrum of $T$.
\end{theorem}

\begin{proof}
This is an abstract reformulation of the continuous functional calculus theorem for the normal operators, that we know from chapter 3. Indeed, that theorem tells us that we have a continuous morphism of $*$-algebras, as follows:
$$C(\sigma(T))\to B(H)\quad,\quad 
f\to f(T)$$

Moreover, by the general properties of the continuous calculus, also established in chapter 3, this morphism is injective, and its image is the norm closed algebra $<T>$ generated by $T,T^*$. Thus, we obtain the isomorphism in the statement.
\end{proof}

The above result is very nice, and it is possible to further build on it, by using this time the spectral theorem for families of normal operators, as follows:

\begin{theorem}
The operator algebra $<T_i>\subset B(H)$ generated by a family of normal operators $T_i\in B(H)$ appears as an algebra of continuous functions,
$$<T>=C(X)$$
where $X\subset\mathbb C$ is a certain compact space associated to the family $\{T_i\}$. Equivalently, any commutative operator algebra $A\subset B(H)$ is of the form $A=C(X)$.
\end{theorem}

\begin{proof}
We have two assertions here, the idea being as follows:

\medskip

(1) Regarding the first assertion, this follows exactly as in the proof of Theorem 5.2, by using this time the spectral theorem for families of normal operators.

\medskip

(2) As for the second assertion, this is clear from the first one, because any commutative algebra $A\subset B(H)$ is generated by its elements $T\in A$, which are all normal.
\end{proof}

All this is good to know, but Theorem 5.2 and Theorem 5.3 remain something quite heavy, based on the spectral theorem. We would like to present now an alternative proof for these results, which is rather elementary, and has the advantage of reconstructing the compact space $X$ directly from the knowledge of the algebra $A$. We will need:

\index{spectrum}
\index{polynomial calculus}
\index{rational calculus}
\index{holomorphic calculus}
\index{unitary}
\index{self-adjoint element}
\index{spectral radius}
\index{continuous calculus}

\begin{theorem}
Given an operator $T\in A\subset B(H)$, define its spectrum as:
$$\sigma(T)=\left\{\lambda\in\mathbb C\Big|T-\lambda\notin A^{-1}\right\}$$
The following spectral theory results hold, exactly as in the $A=B(H)$ case:
\begin{enumerate}
\item We have $\sigma(ST)\cup\{0\}=\sigma(TS)\cup\{0\}$.

\item We have polynomial, rational and holomorphic calculus.

\item As a consequence, the spectra are compact and non-empty.

\item The spectra of unitaries $(U^*=U^{-1})$ and self-adjoints $(T=T^*)$ are on $\mathbb T,\mathbb R$.

\item The spectral radius of normal elements $(TT^*=T^*T)$ is given by $\rho(T)=||T||$.
\end{enumerate}
In addition, assuming $T\in A\subset B$, the spectra of $T$ with respect to $A$ and to $B$ coincide.
\end{theorem}

\begin{proof}
This is something that we know from the beginning of chapter 3, in the case $A=B(H)$. In general the proof is similar, the idea being as follows:

\medskip

(1) Regarding the assertions (1-5), which are of course formulated a bit informally, the proofs here are perfectly similar to those for the full operator algebra $A=B(H)$. All this is standard material, and in fact, things in chapter 3 were written in such a way as for their extension now, to the general operator algebra setting, to be obvious.

\medskip

(2) Regarding the last assertion, the inclusion $\sigma_B(T)\subset\sigma_A(T)$ is clear. For the converse, assume $T-\lambda\in B^{-1}$, and consider the following self-adjoint element:
$$S=(T-\lambda )^*(T-\lambda )$$

The difference between the two spectra of $S\in A\subset B$ is then given by:
$$\sigma_A(S)-\sigma_B(S)=\left\{\mu\in\mathbb C-\sigma_B(S)\Big|(S-\mu)^{-1}\in B-A\right\}$$

Thus this difference in an open subset of $\mathbb C$. On the other hand $S$ being self-adjoint, its two spectra are both real, and so is their difference. Thus the two spectra of $S$ are equal, and in particular $S$ is invertible in $A$, and so $T-\lambda\in A^{-1}$, as desired.

\medskip

(3) As an observation, the last assertion applied with $B=B(H)$ shows that the spectrum $\sigma(T)$ as constructed in the statement coincides with the spectrum $\sigma(T)$ as constructed and studied in chapter 3, so the fact that (1-5) hold indeed is no surprise.

\medskip

(4) Finally, I can hear you screaming that I should have concieved this book differently, matter of not proving the same things twice. Good point, with my distinguished colleague Bourbaki saying the same, and in answer, wait for chapter 7 below, where we will prove exactly the same things a third time. We can discuss pedagogy at that time.
\end{proof}

We can now get back to the commutative algebras, and we have the following result, due to Gelfand, which provides an alternative to Theorem 5.2 and Theorem 5.3:

\index{commutative algebra}
\index{character}
\index{Banach algebra}
\index{spectrum of algebra}
\index{Gelfand theorem}

\begin{theorem}
Any commutative operator algebra $A\subset B(H)$ is of the form 
$$A=C(X)$$
with the ``spectrum'' $X$ of such an algebra being the space of characters $\chi:A\to\mathbb C$, with topology making continuous the evaluation maps $ev_T:\chi\to\chi(T)$.
\end{theorem}

\begin{proof}
Given a commutative operator algebra $A$, we can define $X$ as in the statement. Then $X$ is compact, and $T\to ev_T$ is a morphism of algebras, as follows:
$$ev:A\to C(X)$$

(1) We first prove that $ev$ is involutive. We use the following formula, which is similar to the $z=Re(z)+iIm(z)$ formula for the usual complex numbers:
$$T=\frac{T+T^*}{2}+i\cdot\frac{T-T^*}{2i}$$

Thus it is enough to prove the equality $ev_{T^*}=ev_T^*$ for self-adjoint elements $T$. But this is the same as proving that $T=T^*$ implies that $ev_T$ is a real function, which is in turn true, because $ev_T(\chi)=\chi(T)$ is an element of $\sigma(T)$, contained in $\mathbb R$.

\medskip

(2) Since $A$ is commutative, each element is normal, so $ev$ is isometric:
$$||ev_T||
=\rho(T)
=||T||$$

(3) It remains to prove that $ev$ is surjective. But this follows from the Stone-Weierstrass theorem, because $ev(A)$ is a closed subalgebra of $C(X)$, which separates the points.
\end{proof}

The above theorem of Gelfand is something very beautiful, and far-reaching. It is possible to further build on it, indefinitely high. We will be back to this.

\section*{5b. Von Neumann algebras}

Instead of further building on the above results, which are already quite non-trivial, let us return to our modest status of apprentice operator algebraists, and declare ourselves rather unsatisfied with Definition 5.1, on the following intuitive grounds:

\begin{thought}
Our assumption that $A\subset B(H)$ is norm closed is not satisfying, because we would like $A$ to be stable under polar decomposition, under taking spectral projections, and more generally, under measurable functional calculus.
\end{thought}

Here all these ``defects'' are best visible in the context of Theorem 5.3, with the algebra $A=C(X)$ found there, with $X=\sigma(T)$, being obviously too small. In fact, Theorem 5.3 teaches us that, when looking for a fix, we should look for a weaker topology on $B(H)$, as for the algebra $A=<T>$ generated by a normal operator to be $A=L^\infty(X)$.

\bigskip

So, let us get now into this, topologies on $B(H)$, and fine-tunings of Definition 5.1, based on them. The result that we will need, which is elementary, is as follows: 

\index{weak operator topology}
\index{strong operator topology}
\index{norm topology}

\begin{proposition}
For a subalgebra $A\subset B(H)$, the following are equivalent:
\begin{enumerate}
\item $A$ is closed under the weak operator topology, making each of the linear maps $T\to<Tx,y>$ continuous.

\item $A$ is closed under the strong operator topology, making each of the linear maps $T\to Tx$ continuous.
\end{enumerate}
In the case where these conditions are satisfied, $A$ is closed under the norm topology.
\end{proposition}

\begin{proof}
There are several statements here, the proof being as follows:

\medskip

(1) It is clear that the norm topology is stronger than the strong operator topology, which is in turn stronger than the weak operator topology. At the level of the subsets $S\subset B(H)$ which are closed things get reversed, in the sense that weakly closed implies strongly closed, which in turn implies norm closed. Thus, we are left with proving that for any algebra $A\subset B(H)$, strongly closed implies weakly closed.

\medskip

(2) Consider the Hilbert space obtained by summing $n$ times $H$ with itself:
$$K=H\oplus\ldots\oplus H$$

The operators over $K$ can be regarded as being square matrices with entries in $B(H)$, and in particular, we have a representation $\pi:B(H)\to B(K)$, as follows:
$$\pi(T)=\begin{pmatrix}
T\\
&\ddots\\
&&T
\end{pmatrix}$$

Assume now that we are given an operator $T\in\bar{A}$, with the bar denoting the weak closure. We have then, by using the Hahn-Banach theorem, for any $x\in K$:
\begin{eqnarray*}
T\in\bar{A}
&\implies&\pi(T)\in\overline{\pi(A)}\\
&\implies&\pi(T)x\in\overline{\pi(A)x}\\
&\implies&\pi(T)x\in\overline{\pi(A)x}^{\,||.||}
\end{eqnarray*}

Now observe that the last formula tells us that for any $x=(x_1,\ldots,x_n)$, and any $\varepsilon>0$, we can find $S\in A$ such that the following holds, for any $i$:
$$||Sx_i-Tx_i||<\varepsilon$$

Thus $T$ belongs to the strong operator closure of $A$, as desired.
\end{proof}

Observe that in the above the terminology is a bit confusing, because the norm topology is stronger than the strong operator topology. As a solution, we agree to call the norm topology ``strong'', and the weak and strong operator topologies ``weak'', whenever these two topologies coincide. With this convention made, the algebras $A\subset B(H)$ in Proposition 5.7 are those which are weakly closed. Thus, we can now formulate:

\index{weak topology}
\index{von Neumann algebra}

\begin{definition}
A von Neumann algebra is an operator algebra
$$A\subset B(H)$$
which is closed under the weak topology.
\end{definition}

These algebras will be our main objects of study, in what follows. As basic examples, we have the algebra $B(H)$ itself, then the singly generated algebras, $A=<T>$ with $T\in B(H)$, and then the multiply generated algebras, $A=<T_i>$ with $T_i\in B(H)$. But for the moment, let us keep things simple, and build directly on Definition 5.8, by using basic functional analysis methods. We will need the following key result:

\index{bicommutant}

\begin{theorem}
For an operator algebra $A\subset B(H)$, we have
$$A''=\bar{A}$$
with $A''$ being the bicommutant inside $B(H)$, and $\bar{A}$ being the weak closure.
\end{theorem}

\begin{proof}
We can prove this by double inclusion, as follows:

\medskip

``$\supset$'' Since any operator commutes with the operators that it commutes with, we have a trivial inclusion $S\subset S''$, valid for any set $S\subset B(H)$. In particular, we have:
$$A\subset A''$$

Our claim now is that the algebra $A''$ is closed, with respect to the strong operator topology. Indeed, assuming that we have $T_i\to T$ in this topology, we have:
\begin{eqnarray*}
T_i\in A''
&\implies&ST_i=T_iS,\ \forall S\in A'\\
&\implies&ST=TS,\ \forall S\in A'\\
&\implies&T\in A
\end{eqnarray*}

Thus our claim is proved, and together with Proposition 5.7, which allows us to pass from the strong to the weak operator topology, this gives $\bar{A}\subset A''$, as desired.

\medskip

``$\subset$'' Here we must prove that we have the following implication, valid for any $T\in B(H)$, with the bar denoting as usual the weak operator closure:
$$T\in A''\implies T\in\bar{A}$$

For this purpose, we use the same amplification trick as in the proof of Proposition 5.7. Consider the Hilbert space obtained by summing $n$ times $H$ with itself:
$$K=H\oplus\ldots\oplus H$$

The operators over $K$ can be regarded as being square matrices with entries in $B(H)$, and in particular, we have a representation $\pi:B(H)\to B(K)$, as follows:
$$\pi(T)=\begin{pmatrix}
T\\
&\ddots\\
&&T
\end{pmatrix}$$

The idea will be that of doing the computations in this representation. First, in this representation, the image of our algebra $A\subset B(H)$ is given by:
$$\pi(A)=\left\{\begin{pmatrix}
T\\
&\ddots\\
&&T
\end{pmatrix}\Big|T\in A\right\}$$

We can compute the commutant of this image, exactly as in the usual scalar matrix case, and we obtain the following formula:
$$\pi(A)'=\left\{\begin{pmatrix}
S_{11}&\ldots&S_{1n}\\
\vdots&&\vdots\\
S_{n1}&\ldots&S_{nn}
\end{pmatrix}\Big|S_{ij}\in A'\right\}$$

We conclude from this that, given an operator $T\in A''$ as above, we have:
$$\begin{pmatrix}
T\\
&\ddots\\
&&T
\end{pmatrix}\in\pi(A)''$$

In other words, the conclusion of all this is that we have:
$$T\in A''\implies \pi(T)\in\pi(A)''$$

Now given a vector $x\in K$, consider the orthogonal projection $P\in B(K)$ on the norm closure of the vector space $\pi(A)x\subset K$. Since the subspace $\pi(A)x\subset K$ is invariant under the action of $\pi(A)$, so is its norm closure inside $K$, and we obtain from this:
$$P\in\pi(A)'$$

By combining this with what we found above, we conclude that we have:
$$T\in A''\implies \pi(T)P=P\pi(T)$$

Since this holds for any $x\in K$, we conclude that any operator $T\in A''$ belongs to the strong operator closure of $A$. By using now Proposition 5.7, which allows us to pass from the strong to the weak operator closure, we conclude that we have:
$$A''\subset\bar{A}$$

Thus, we have the desired reverse inclusion, and this finishes the proof.
\end{proof}

Now by getting back to the von Neumann algebras, from Definition 5.8, we have the following result, which is a reformulation of Theorem 5.9, by using this notion:

\index{bicommutant}
\index{von Neumann algebra}

\begin{theorem}
For an operator algebra $A\subset B(H)$, the following are equivalent:
\begin{enumerate}
\item $A$ is weakly closed, so it is a von Neumann algebra.

\item $A$ equals its algebraic bicommutant $A''$, taken inside $B(H)$.
\end{enumerate}
\end{theorem}

\begin{proof}
This follows from the formula $A''=\bar{A}$ from Theorem 5.9, along with the trivial fact that the commutants are automatically weakly closed.
\end{proof}

The above statement, called bicommutant theorem, and due to von Neumann \cite{vn1}, is quite interesting, philosophically speaking. Among others, it shows that the von Neumann algebras are exactly the commutants of the self-adjoint sets of operators:

\index{commutant algebra}
\index{tricommutant formula}

\begin{proposition}
Given a subset $S\subset B(H)$ which is closed under $*$, the commutant
$$A=S'$$
is a von Neumann algebra. Any von Neumann algebra appears in this way.
\end{proposition}

\begin{proof}
We have two assertions here, the idea being as follows:

\medskip

(1) Given $S\subset B(H)$ satisfying $S=S^*$, the commutant $A=S'$ satisfies $A=A^*$, and is also weakly closed. Thus, $A$ is a von Neumann algebra. Note that this follows as well from the following ``tricommutant formula'', which follows from Theorem 5.10:
$$S'''=S'$$

(2) Given a von Neumann algebra $A\subset B(H)$, we can take $S=A'$. Then $S$ is closed under the involution, and we have $S'=A$, as desired.
\end{proof}

Observe that Proposition 5.11 can be regarded as yet another alternative definition for the von Neumann algebras, and with this definition being probably the best one when talking about quantum mechanics, where the self-adjoint operators $T:H\to H$ can be though of as being ``observables'' of the system, and with the commutants $A=S'$ of the sets of such observables $S=\{T_i\}$ being the algebras $A\subset B(H)$ that we are interested in. And with all this actually needing some discussion about self-adjointness, and about  boundedness too, but let us not get into this here, and stay mathematical, as before.

\bigskip

As another interesting consequence of Theorem 5.10, we have:

\index{center of algebra}
\index{commutative algebra}

\begin{proposition}
Given a von Neumann algebra $A\subset B(H)$, its center
$$Z(A)=A\cap A'$$
regarded as an algebra $Z(A)\subset B(H)$, is a von Neumann algebra too.
\end{proposition}

\begin{proof}
This follows from the fact that the commutants are weakly closed, that we know from the above, which shows that $A'\subset B(H)$ is a von Neumann algebra. Thus, the intersection $Z(A)=A\cap A'$ must be a von Neumann algebra too, as claimed.
\end{proof}

In order to develop some general theory, let us start by investigating the finite dimensional case. Here the ambient algebra is $B(H)=M_N(\mathbb C)$, any linear subspace $A\subset B(H)$ is automatically closed, for all 3 topologies in Proposition 5.7, and we have:

\index{finite dimensional algebra}
\index{sum of matrix algebras}
\index{multimatrix algebra}

\begin{theorem}
The $*$-algebras $A\subset M_N(\mathbb C)$ are exactly the algebras of the form
$$A=M_{n_1}(\mathbb C)\oplus\ldots\oplus M_{n_k}(\mathbb C)$$
depending on parameters $k\in\mathbb N$ and $n_1,\ldots,n_k\in\mathbb N$ satisfying
$$n_1+\ldots+n_k=N$$
embedded into $M_N(\mathbb C)$ via the obvious block embedding, twisted by a unitary $U\in U_N$.
\end{theorem}

\begin{proof}
We have two assertions to be proved, the idea being as follows:

\medskip

(1) Given numbers $n_1,\ldots,n_k\in\mathbb N$ satisfying $n_1+\ldots+n_k=N$, we have indeed an obvious embedding of $*$-algebras, via matrix blocks, as follows:
$$M_{n_1}(\mathbb C)\oplus\ldots\oplus M_{n_k}(\mathbb C)\subset M_N(\mathbb C)$$

In addition, we can twist this embedding by a unitary $U\in U_N$, as follows:
$$M\to UMU^*$$

(2) In the other sense now, consider a $*$-algebra $A\subset M_N(\mathbb C)$. It is elementary to prove that the center $Z(A)=A\cap A'$, as an algebra, is of the following form:
$$Z(A)\simeq\mathbb C^k$$

Consider now the standard basis $e_1,\ldots,e_k\in\mathbb C^k$, and let  $p_1,\ldots,p_k\in Z(A)$ be the images of these vectors via the above identification. In other words, these elements $p_1,\ldots,p_k\in A$ are central minimal projections, summing up to 1:
$$p_1+\ldots+p_k=1$$

The idea is then that this partition of the unity will eventually lead to the block decomposition of $A$, as in the statement. We prove this in 4 steps, as follows:

\medskip

\underline{Step 1}. We first construct the matrix blocks, our claim here being that each of the following linear subspaces of $A$ are non-unital $*$-subalgebras of $A$:
$$A_i=p_iAp_i$$

But this is clear, with the fact that each $A_i$ is closed under the various non-unital $*$-subalgebra operations coming from the projection equations $p_i^2=p_i^*=p_i$.

\medskip

\underline{Step 2}. We prove now that the above algebras $A_i\subset A$ are in a direct sum position, in the sense that we have a non-unital $*$-algebra sum decomposition, as follows:
$$A=A_1\oplus\ldots\oplus A_k$$

As with any direct sum question, we have two things to be proved here. First, by using the formula $p_1+\ldots+p_k=1$ and the projection equations $p_i^2=p_i^*=p_i$, we conclude that we have the needed generation property, namely:
$$A_1+\ldots+ A_k=A$$

As for the fact that the sum is indeed direct, this follows as well from the formula $p_1+\ldots+p_k=1$, and from the projection equations $p_i^2=p_i^*=p_i$.

\medskip

\underline{Step 3}. Our claim now, which will finish the proof, is that each of the $*$-subalgebras $A_i=p_iAp_i$ constructed above is a full matrix algebra. To be more precise here, with $n_i=rank(p_i)$, our claim is that we have isomorphisms, as follows:
$$A_i\simeq M_{n_i}(\mathbb C)$$

In order to prove this claim, recall that the projections $p_i\in A$ were chosen central and minimal. Thus, the center of each of the algebras $A_i$ reduces to the scalars:
$$Z(A_i)=\mathbb C$$

But this shows, either via a direct computation, or via the bicommutant theorem, that the each of the algebras $A_i$ is a full matrix algebra, as claimed.

\medskip

\underline{Step 4}. We can now obtain the result, by putting together what we have. Indeed, by using the results from Step 2 and Step 3, we obtain an isomorphism as follows:
$$A\simeq M_{n_1}(\mathbb C)\oplus\ldots\oplus M_{n_k}(\mathbb C)$$

Moreover, a more careful look at the isomorphisms established in Step 3 shows that at the global level, that of the algebra $A$ itself, the above isomorphism simply comes by twisting the following standard multimatrix embedding, discussed in the beginning of the proof, (1) above, by a certain unitary matrix $U\in U_N$:
$$M_{n_1}(\mathbb C)\oplus\ldots\oplus M_{n_k}(\mathbb C)\subset M_N(\mathbb C)$$

Now by putting everything together, we obtain the result.
\end{proof}

In relation with the bicommutant theorem, we have the following result, which fully clarifies the situation, with a very explicit proof, in finite dimensions:

\index{commutant}
\index{bicommutant}

\begin{proposition}
Consider a $*$-algebra $A\subset M_N(\mathbb C)$, written as above:
$$A= M_{n_1}(\mathbb C)\oplus\ldots\oplus M_{n_k}(\mathbb C)$$
The commutant of this algebra is then, with respect with the block decomposition used,
$$A'=\mathbb C\oplus\ldots\oplus\mathbb C$$
and by taking one more time the commutant we obtain $A$ itself, $A=A''$.
\end{proposition}

\begin{proof}
Let us decompose indeed our algebra $A$ as in Theorem 5.13:
$$A=M_{n_1}(\mathbb C)\oplus\ldots\oplus M_{n_k}(\mathbb C)$$

The center of each matrix algebra being reduced to the scalars, the commutant of this algebra is then as follows, with each copy of $\mathbb C$ corresponding to a matrix block:
$$A'=\mathbb C\oplus\ldots\oplus\mathbb C$$

By taking once again the commutant we obtain $A$ itself, and we are done.
\end{proof}

As another interesting application of Theorem 5.13, clarifying this time the relation with operator theory, in finite dimensions, we have the following result:

\index{finite dimensional algebra}

\begin{theorem}
Given an operator $T\in B(H)$ in finite dimensions, $H=\mathbb C^N$, the von Neumann algebra $A=<T>$ that it generates inside $B(H)=M_N(\mathbb C)$ is
$$A=M_{n_1}(\mathbb C)\oplus\ldots\oplus M_{n_k}(\mathbb C)$$
with the sizes of the blocks $n_1,\ldots,n_k\in\mathbb N$ coming from the spectral theory of the associated matrix $M\in M_N(\mathbb C)$. In the normal case $TT^*=T^*T$, this decomposition comes from
$$T=UDU^*$$
with $D\in M_N(\mathbb C)$ diagonal, and with $U\in U_N$ unitary.
\end{theorem}

\begin{proof}
This is something which is routine, by using the linear algebra and spectral theory developed in chapter 1, for the matrices $M\in M_N(\mathbb C)$. To be more precise:

\medskip

(1) The fact that $A=<T>$ decomposes into a direct sum of matrix algebras is something that we already know, coming from Theorem 5.13.

\medskip

(2) By using standard linear algebra, we can compute the block sizes $n_1,\ldots,n_k\in\mathbb N$, from the knowledge of the spectral theory of the associated matrix $M\in M_N(\mathbb C)$.

\medskip

(3) In the normal case, $TT^*=T^*T$, we can simply invoke the spectral theorem, and by suitably changing the basis, we are led to the conclusion in the statement.
\end{proof}

Let us get now to infinite dimensions, with Theorem 5.15 as our main source of inspiration. The same argument applies, provided that we are in the normal case, and we have the following result, summarizing our basic knowledge here:

\index{singly generated algebra}
\index{measurable calculus}

\begin{theorem}
Given a bounded operator $T\in B(H)$ which is normal, $TT^*=T^*T$, the von Neumann algebra $A=<T>$ that it generates inside $B(H)$ is
$$<T>=L^\infty(\sigma(T))$$
with $\sigma(T)\subset\mathbb C$ being as usual its spectrum.
\end{theorem}

\begin{proof}
The measurable functional calculus theorem for the normal operators tells us that we have a weakly continuous morphism of $*$-algebras, as follows:
$$L^\infty(\sigma(T))\to B(H)\quad,\quad 
f\to f(T)$$

Moreover, by the general properties of the measurable calculus, also established in chapter 3, this morphism is injective, and its image is the weakly closed algebra $<T>$ generated by $T,T^*$. Thus, we obtain the isomorphism in the statement.
\end{proof}

More generally now, along the same lines, we have the following result:

\index{commuting normal operators}

\begin{theorem}
Given operators $T_i\in B(H)$ which are normal, and which commute, the von Neumann algebra $A=<T_i>$ that these operators generates inside $B(H)$ is
$$<T_i>=L^\infty(X)$$
with $X$ being a certain measured space, associated to the family $\{T_i\}$.
\end{theorem}

\begin{proof}
This is once again routine, by using the spectral theory for the families of commuting normal operators $T_i\in B(H)$ developed in chapter 3. 
\end{proof}

As a fundamental consequence now of the above results, we have:

\index{commutative algebra}
\index{measured space}

\begin{theorem}
The commutative von Neumann algebras are the algebras
$$A=L^\infty(X)$$
with $X$ being a measured space.
\end{theorem}

\begin{proof}
We have two assertions to be proved, the idea being as follows:

\medskip

(1) In one sense, we must prove that given a measured space $X$, we can realize the $A=L^\infty(X)$ as a von Neumann algebra, on a certain Hilbert space $H$. But this is something that we know since chapter 2, the representation being as follows: 
$$L^\infty(X)\subset B(L^2(X))\quad,\quad f\to(g\to fg)$$

(2) In the other sense, given a commutative von Neumann algebra $A\subset B(H)$, we must construct a certain measured space $X$, and an identification $A=L^\infty(X)$. But this follows from Theorem 5.17, because we can write our algebra as follows:
$$A=<T_i>$$

To be more precise, $A$ being commutative, any element $T\in A$ is normal, so we can pick a basis $\{T_i\}\subset A$, and then we have $A=<T_i>$ as above, with $T_i\in B(H)$ being commuting normal operators. Thus Theorem 5.17 applies, and gives the result.

\medskip

(3) Alternatively, and more explicitly, we can deduce this from Theorem 5.16, applied with $T=T^*$. Indeed, by using $T=Re(T)+iIm(T)$, we conclude that any von Neumann algebra $A\subset B(H)$ is generated by its self-adjoint elements $T\in A$. Moreover, by using measurable functional calculus, we conclude that $A$ is linearly generated by its projections. But then, assuming $A=\overline{span}\{p_i\}$, with $p_i$ being projections, we can set:
$$T=\sum_{i=0}^\infty\frac{p_i}{3^i}$$

Then $T=T^*$, and by functional calculus we have $p_0\in<T>$, then $p_1\in <T>$, and so on. Thus $A=<T>$, and $A=L^\infty(X)$ comes now via Theorem 5.16, as claimed.
\end{proof}

The above result is the foundation for all the advanced von Neumann algebra theory, that we will discuss in the remainder of this book, and there are many things that can be said about it. To start with, in relation with the general theory of the normed closed algebras, that we developed in the beginning of this chapter, we have:

\begin{warning}
Although the von Neumann algebras are norm closed, the theory of norm closed algebras does not always apply well to them. For instance for $A=L^\infty(X)$ Gelfand gives $A=C(\widehat{X})$, with $\widehat{X}$ being a certain technical compactification of $X$.
\end{warning}

In short, this would be my advice, do not mess up the two theories that we will be developing in this book, try finding different rooms for them, in your brain. At least at this stage of things, because later, do not worry, we will be playing with both.

\bigskip

Now forgetting about Gelfand, and taking Theorem 5.18 as such, tentative foundation for the theory that we want to develop, as a first consequence of this, we have:

\index{center of algebra}
\index{commutative algebra}

\begin{theorem}
Given a von Neumann algebra $A\subset B(H)$, we have
$$Z(A)=L^\infty(X)$$
with $X$ being a certain measured space.
\end{theorem}

\begin{proof}
We know from Proposition 5.12 that the center $Z(A)\subset B(H)$ is a von Neumann algebra. Thus Theorem 5.18 applies, and gives the result.
\end{proof}

It is possible to further build on this, with a powerful decomposition result as follows, over the measured space $X$ constructed in Theorem 5.20:
$$A=\int_XA_x\,dx$$

But more on this later, after developing the appropriate tools for this program, which is something non-trivial. Among others, before getting into such things, we will have to study the von Neumann algebras $A$ having trivial center, $Z(A)=\mathbb C$, called factors, which include the fibers $A_x$ in the above decomposition result. More on this later.

\section*{5c. Random matrices}

Our main results so far on the von Neumann algebras concern the finite dimensional case, where the algebra is of the form $A=\oplus_iM_{n_i}(\mathbb C)$, and the commutative case, where the algebra is of the form $A=L^\infty(X)$. In order to advance, we must solve:

\begin{question}
What are the next simplest von Neumann algebras, generalizing at the same time the finite dimensional ones, $A=\oplus_iM_{n_i}(\mathbb C)$, and the commutative ones, $A=L^\infty(X)$, that we can use as input for our study?
\end{question}

In this formulation, our question is a no-brainer, the answer to it being that of looking at the direct integrals of matrix algebras, over an arbitrary measured space $X$:
$$A=\int_XM_{n_x}(\mathbb C)dx$$

However, when thinking a bit, all this looks quite tricky, with most likely lots of technical functional analysis and measure theory involved. So, we will leave the investigation of such algebras, which are indeed quite basic, and called of type I, for later.

\bigskip

Nevermind. Let us replace Question 5.21 with something more modest, as follows:

\begin{question}[update]
What are the next simplest von Neumann algebras, generalizing at the same time the usual matrix algebras, $A=M_N(\mathbb C)$, and the commutative ones, $A=L^\infty(X)$, that we can use as input for our study?
\end{question}

But here, what we have is again a no-brainer, because in relation to what has been said above, we just have to restrict the attention to the ``isotypic'' case, where all fibers are isomorphic. And in this case our algebra is a random matrix algebra:
$$A=\int_XM_N(\mathbb C)dx$$

Which looks quite nice, and so good news, we have our algebras. In practice now, although there is some functional analysis to be done with these algebras, the main questions regard the individual operators $T\in A$, called random matrices. Thus, we are basically back to good old operator theory. Let us begin our discussion with:

\index{random matrix}
\index{random matrix algebra}
\index{type I algebra}

\begin{definition}
A random matrix algebra is a von Neumann algebra of the following type, with $X$ being a probability space, and with $N\in\mathbb N$ being an integer:
$$A=M_N(L^\infty(X))$$
In other words, $A$ appears as a tensor product, as follows, 
$$A=M_N(\mathbb C)\otimes L^\infty(X)$$
of a matrix algebra and a commutative von Neumann algebra.
\end{definition}

As a first observation, our algebra can be written as well as follows, with this latter convention being quite standard in the probability literature:
$$A=L^\infty(X,M_N(\mathbb C))$$

In connection with the tensor product notation, which is often the most useful one for computations, we have as well the following possible writing, also used in probability:
$$A=L^\infty(X)\otimes M_N(\mathbb C)$$

Importantly now, each random matrix algebra $A$ is naturally endowed with a canonical von Neumann algebra trace $tr:A\to\mathbb C$, which appears as follows:

\begin{proposition}
Given a random matrix algebra $A=M_N(L^\infty(X))$, consider the linear form $tr:A\to\mathbb C$ given by:
$$tr(T)=\frac{1}{N}\sum_{i=1}^N\int_X T_{ii}^xdx$$
In tensor product notation, $A=M_N(\mathbb C)\otimes L^\infty(X)$, we have then the formula
$$tr=\frac{1}{N}\,Tr\otimes\int_X$$
and this functional $tr:A\to\mathbb C$ is a faithful positive unital trace.
\end{proposition}

\begin{proof}
The first assertion, regarding the tensor product writing of $tr$, is clear from definitions. As for the second assertion, regarding the various properties of $tr$, this follows from this, because these properties are stable under taking tensor products.
\end{proof}

As before, there is a discussion here in connection with the other possible writings of $A$. With the probabilistic notation $A=L^\infty(X,M_N(\mathbb C))$, the trace appears as:
$$tr(T)=\int_X\frac{1}{N}\,Tr(T^x)\,dx$$

Also, with the probabilistic tensor notation $A=L^\infty(X)\otimes M_N(\mathbb C)$, the trace appears exactly as in the second part of Proposition 5.24, with the order inverted:
$$tr=\int_X\otimes\,\,\frac{1}{N}\,Tr$$

To summarize, the random matrix algebras appear to be very basic objects, and the only difficulty, in the beginning, lies in getting familiar with the 4 possible notations for them. Or perhaps 5 possible notations, because we have $A=\int_XM_N(\mathbb C)dx$ as well.

\bigskip

Getting to work now, as already said, the main questions about random matrix algebras regard the individual operators $T\in A$, called random matrices. To be more precise, we are interested in computing the laws of such matrices, constructed according to:

\index{law}
\index{distribution}

\begin{theorem}
Given an operator algebra $A\subset B(H)$ with a faithful trace $tr:A\to\mathbb C$, any normal element $T\in A$ has a law, namely a probability measure $\mu$ satisfying
$$tr(T^k)=\int_\mathbb Cz^kd\mu(z)$$
with the powers being with respect to colored exponents $k=\circ\bullet\bullet\circ\ldots\,$, defined via
$$a^\emptyset=1\quad,\quad a^\circ=a\quad,\quad a^\bullet=a^*$$
and multiplicativity. This law is unique, and is supported by the spectrum $\sigma(T)\subset\mathbb C$. In the non-normal case, $TT^*\neq T^*T$, such a law does not exist.
\end{theorem}

\begin{proof}
We have two assertions here, the idea being as follows:

\medskip

(1) In the normal case, $TT^*=T^*T$, we know from Theorem 5.2, based on the continuous functional calculus theorem, that we have: 
$$<T>=C(\sigma(T))$$

Thus the functional $f(T)\to tr(f(T))$ can be regarded as an integration functional on the algebra $C(\sigma(T))$, and by the Riesz theorem, this latter functional must come from a probability measure $\mu$ on the spectrum $\sigma(T)$, in the sense that we must have:
$$tr(f(T))=\int_{\sigma(T)}f(z)d\mu(z)$$

We are therefore led to the conclusions in the statement, with the uniqueness assertion coming from the fact that the operators $T^k$, taken as usual with respect to colored integer exponents, $k=\circ\bullet\bullet\circ\ldots$\,, generate the whole operator algebra $C(\sigma(T))$.

\medskip

(2) In the non-normal case now, $TT^*\neq T^*T$, we must show that such a law does not exist. For this purpose, we can use a positivity trick, as follows:
\begin{eqnarray*}
TT^*-T^*T\neq0
&\implies&(TT^*-T^*T)^2>0\\
&\implies&TT^*TT^*-TT^*T^*T-T^*TTT^*+T^*TT^*T>0\\
&\implies&tr(TT^*TT^*-TT^*T^*T-T^*TTT^*+T^*TT^*T)>0\\
&\implies&tr(TT^*TT^*+T^*TT^*T)>tr(TT^*T^*T+T^*TTT^*)\\
&\implies&tr(TT^*TT^*)>tr(TTT^*T^*)
\end{eqnarray*}

Now assuming that $T$ has a law $\mu\in\mathcal P(\mathbb C)$, in the sense that the moment formula in the statement holds, the above two different numbers would have to both appear by integrating $|z|^2$ with respect to this law $\mu$, which is contradictory, as desired.
\end{proof}

Back now to the random matrices, as a basic example, assume $X=\{.\}$, so that we are dealing with a usual scalar matrix, $T\in M_N(\mathbb C)$. By changing the basis of $\mathbb C^N$, which won't affect our trace computations, we can assume that $T$ is diagonal:
$$T\sim
\begin{pmatrix}
\lambda_1\\
&\ddots\\
&&\lambda_N
\end{pmatrix}$$

But for such a diagonal matrix, we have the following formula:
$$tr(T^k)=\frac{1}{N}(\lambda_1^k+\ldots+\lambda_N^k)$$

Thus, the law of $T$ is the average of the Dirac masses at the eigenvalues:
$$\mu=\frac{1}{N}\left(\delta_{\lambda_1}+\ldots+\delta_{\lambda_N}\right)$$

As a second example now, assume $N=1$, and so $T\in L^\infty(X)$. In this case we obtain the usual law of $T$, because the equation to be satisfied by $\mu$ is:
$$\int_X\varphi(T)=\int_\mathbb C\varphi(x)d\mu(x)$$

At a more advanced level, the main problem regarding the random matrices is that of computing the law of various classes of such matrices, coming in series:

\begin{question}
What is the law of random matrices coming in series
$$T_N\in M_N(L^\infty(X))$$
in the $N>>0$ regime?
\end{question}

The general strategy here, coming from physicists, is that of computing first the asymptotic law $\mu^0$, in the $N\to\infty$ limit, and then looking for the higher order terms as well, as to finally reach to a series in $N^{-1}$ giving the law of $T_N$, as follows: 
$$\mu_N=\mu^0+N^{-1}\mu^1+N^{-2}\mu^2+\ldots$$

As a basic example here, of particular interest are the random matrices having i.i.d. complex normal entries, under the constraint $T=T^*$. Here the asymptotic law $\mu^0$ is the Wigner semicircle law on $[-2,2]$. We will discuss this in chapter 6 below, and in the meantime we can only recommend some reading, from the original papers of Marchenko-Pastur \cite{mpa}, Voiculescu \cite{vo2}, Wigner \cite{wig}, and from the books of Anderson-Guionnet-Zeitouni \cite{agz}, Mehta \cite{meh}, Nica-Speicher \cite{nsp}, Voiculescu-Dykema-Nica \cite{vdn}.

\section*{5d. Quantum spaces}

Let us end this preliminary chapter on operator algebras with some philosophy, a bit a la Heisenberg. In relation with general ``quantum space'' goals, Theorem 5.18 is something very interesting, philosophically speaking, suggesting us to formulate:

\index{quantum space}
\index{quantum measured space}

\begin{definition}
Given a von Neumann algebra $A\subset B(H)$, we write
$$A=L^\infty(X)$$
and call $X$ a quantum measured space.
\end{definition}

As an example here, for the simplest noncommutative von Neumann algebra that we know, namely the usual matrix algebra $A=M_N(\mathbb C)$, the formula that we want to write is as follows, with $M_N$ being a certain mysterious quantum space:
$$M_N(\mathbb C)=L^\infty(M_N)$$

So, what can we say about this space $M_N$? As a first observation, this is a finite space, with its cardinality being defined and computed as follows:
$$|M_N|
=\dim_\mathbb CM_N(\mathbb C)
=N^2$$

Now since this is the same as the cardinality of the set $\{1,\ldots,N^2\}$, we are led to the conclusion that we should have a twisting result as follows, with the twisting operation $X\to X^\sigma$ being something that destroys the points, but keeps the cardinality:
$$M_N=\{1,\ldots,N^2\}^\sigma$$

From an analytic viewpoint now, we would like to understand what is the integration over $M_N$, giving rise to the corresponding $L^\infty$ functions. And here, we can set:
$$\int_{M_N}A=tr(A)$$

To be more precise, on the left we have the integral of an arbitrary function on $M_N$, which according to our conventions, should be a usual matrix:
$$A\in L^\infty(M_N)=M_N(\mathbb C)$$

As for the quantity on the right, the outcome of the computation, this can only be the trace of $A$. In addition, it is better to choose this trace to be normalized, by $tr(1)=1$, and this in order for our measure on $M_N$ to have mass 1, as it is ideal:
$$tr(A)=\frac{1}{N}\,Tr(A)$$

We can say even more about this. Indeed, since the traces of positive matrices are positive, we are led to the following formula, to be taken with the above conventions, which shows that the measure on $M_N$ that we constructed is a probability measure:
$$A>0\implies \int_{M_N}A>0$$

Before going further, let us record what we found, for future reference:

\index{quantum space}
\index{twist}

\begin{theorem}
The quantum measured space $M_N$ formally given by
$$M_N(\mathbb C)=L^\infty(M_N)$$
has cardinality $N^2$, appears as a twist, in a purely algebraic sense,
$$M_N=\{1,\ldots,N^2\}^\sigma$$
and is a probability space, its uniform integration being given by
$$\int_{M_N}A=tr(A)$$
where at right we have the normalized trace of matrices, $tr=Tr/N$.
\end{theorem}

\begin{proof}
This is something half-informal, mostly for fun, which basically follows from the above discussion, the details and missing details being as follows:

\medskip

(1) In what regards the formula $|M_N|=N^2$, coming by computing the complex vector space dimension, as explained above, this is obviously something rock-solid.

\medskip

(2) Regarding twisting, we would like to have a formula as follows, with the operation $A\to A^\sigma$ being something that destroys the commutativity of the multiplication:
$$L^\infty(M_N)=L^\infty(1,\ldots,N^2)^\sigma$$

In more familiar terms, with usual complex matrices on the left, and with a better-looking product of sets being used on the right, this formula reads:
$$M_N(\mathbb C)=L^\infty\Big(\{1,\ldots,N\}\times\{1,\ldots,N\}\Big)^\sigma$$

In order to establish this formula, consider the algebra on the right. As a complex vector space, this algebra has the standard basis $\{f_{ij}\}$ formed by the Dirac masses at the points $(i,j)$, and the multiplicative structure of this algebra is given by:
$$f_{ij}f_{kl}=\delta_{ij,kl}$$

Now let us twist this multiplication, according to the formula $e_{ij}e_{kl}=\delta_{jk}e_{il}$. We obtain in this way the usual combination formulae for the standard matrix units $e_{ij}:e_j\to e_i$ of the algebra $M_N(\mathbb C)$, and so we have our twisting result, as claimed. 

\medskip

(3) In what regards the integration formula in the statement, with the conclusion that the underlying measure on $M_N$ is a probability one, this is something that we fully explained before, and as for the result (1) above, it is something rock-solid.

\medskip

(4) As a last technical comment, observe that the twisting operation performed in (2) destroys both the involution, and the trace of the algebra. This is something quite interesting, which cannot be fixed, and we will back to it, later on.
\end{proof}

In order to advance now, based on the above result, the key point there is the construction and interpretation of the trace $tr:M_N(\mathbb C)\to\mathbb C$, as an integration functional. But this leads us into the following natural, and quite puzzling question:

\begin{question}
In the general context of Definition 5.27, where we formally wrote $A=L^\infty(X)$, what is the underlying integration functional $tr:A\to\mathbb C$?
\end{question}

This is a quite subtle question, and there are several possible answers here. For instance, we would like the integration functional to have the following property:
$$tr(ab)=tr(ba)$$

And the problem is that certain von Neumann algebras do not possess such traces. This is actually something quite advanced, that we do not know yet, but by anticipating a bit, we are in trouble, and we must modify Definition 5.27, as follows:

\index{quantum space}
\index{quantum measured space}
\index{quantum probability space}

\begin{definition}[update]
Given a von Neumann algebra $A\subset B(H)$, coming with a faithful positive unital trace $tr:A\to\mathbb C$, we write
$$A=L^\infty(X)$$
and call $X$ a quantum probability space. We also write the trace as $tr=\int_X$, and call it integration with respect to the uniform measure on $X$.
\end{definition}

At the level of examples, passed the classical probability spaces $X$, we know from Theorem 5.28 that the quantum space $M_N$ is a finite quantum probability space. But this raises the question of understanding what the finite quantum probability spaces are, in general. For this purpose, we need to examine the finite dimensional von Neumann algebras. And the result here, extending Theorem 5.13, is as follows:

\begin{theorem}
The finite dimensional von Neumann algebras $A\subset B(H)$ over an arbitrary Hilbert space $H$ are exactly the direct sums of matrix algebras,
$$A=M_{n_1}(\mathbb C)\oplus\ldots\oplus M_{n_k}(\mathbb C)$$
embedded into $B(H)$ by using a partition of unity of $B(H)$ with rank $1$ projections
$$1=P_1+\ldots+P_k$$
with the ``factors'' $M_{n_i}(\mathbb C)$ being each embedded into the algebra $P_iB(H)P_i$.
\end{theorem}

\begin{proof}
This is standard, as in the case $A\subset M_N(\mathbb C)$. Consider the center of $A$, which is a finite dimensional commutative von Neumann algebra, of the following form:
$$Z(A)=\mathbb C^k$$

Now let $P_i$ be the Dirac mass at $i\in\{1,\ldots,k\}$. Then $P_i\in B(H)$ is an orthogonal projection, and these projections form a partition of unity, as follows:
$$1=P_1+\ldots+P_k$$

With $A_i=P_iAP_i$, we have then a non-unital $*$-algebra decomposition, as follows:
$$A=A_1\oplus\ldots\oplus A_k$$

On the other hand, it follows from the minimality of each of the projections $P_i\in Z(A)$ that we have unital $*$-algebra isomorphisms $A_i\simeq M_{n_i}(\mathbb C)$, and this gives the result.
\end{proof}

We can now deduce what the finite quantum measured spaces are, in the sense of the old Definition 5.27. Indeed, we must solve here the following equation:
$$L^\infty(X)=M_{n_1}(\mathbb C)\oplus\ldots\oplus M_{n_k}(\mathbb C)$$

Now since the direct unions of sets correspond to direct sums at the level of the associated algebras of functions, in the classical case, we can take the following formula as a definition for a direct union of sets, in the general, noncommutative case:
$$L^\infty(X_1\sqcup\ldots\sqcup X_k)=L^\infty(X_1)\oplus\ldots\oplus L^\infty(X_k)$$

With this, and by remembering the definition of $M_N$, we are led to the conclusion that the solution to our quantum measured space equation above is as follows:
$$X=M_{n_1}\sqcup\ldots\sqcup M_{n_k}$$

For fully solving our problem, in the spirit of the new Definition 5.30, we still have to discuss the traces on $L^\infty(X)$. We are led in this way to the following statement: 

\index{quantum space}
\index{finite quantum space}
\index{finite dimensional algebra}
\index{multimatrix algebra}
\index{sum of matrix algebras}
\index{canonical trace}

\begin{theorem}
The finite quantum measured spaces are the spaces
$$X=M_{n_1}\sqcup\ldots\sqcup M_{n_k}$$
according to the following formula, for the associated algebras of functions:
$$L^\infty(X)=M_{n_1}(\mathbb C)\oplus\ldots\oplus M_{n_k}(\mathbb C)$$
The cardinality $|X|$ of such a space is the following number,
$$N=n_1^2+\ldots+n_k^2$$
and the possible traces are as follows, with $\lambda_i>0$ summing up to $1$:
$$tr=\lambda_1tr_1\oplus\ldots\oplus\lambda_ktr_k$$
Among these traces, we have the canonical trace, appearing as
$$tr:L^\infty(X)\subset\mathcal L(L^\infty(X))\to\mathbb C$$
via the left regular representation, having weights $\lambda_i=n_i^2/N$.
\end{theorem}

\begin{proof}
We have many assertions here, basically coming from the above discussion, with only the last one needing some explanations. Consider the left regular representation of our algebra $A=L^\infty(X)$, which is given by the following formula:
$$\pi:A\subset\mathcal L(A)\quad,\quad 
\pi(a):b\to ab$$

We know that the algebra $\mathcal L(A)$ of linear operators $T:A\to A$ is isomorphic to a matrix algebra, and more specifically to $M_N(\mathbb C)$, with $N=|X|$ being as before:
$$\mathcal L(A)\simeq M_N(\mathbb C)$$

Thus, this algebra has a trace $tr:\mathcal L(A)\to\mathbb C$, and by composing this trace with the representation $\pi$, we obtain a certain trace $tr:A\to\mathbb C$, that we can call ``canonical'':
$$tr:A\subset\mathcal L(A)\to\mathbb C$$

We can compute the weights of this trace by using a multimatrix basis of $A$, formed by matrix units $e_{ab}^i$, with $i\in\{1,\ldots,k\}$ and with $a,b\in\{1,\ldots,n_i\}$, and we obtain:
$$\lambda_i=\frac{n_i^2}{N}$$

Thus, we are led to the conclusion in the statement.
\end{proof}

We will be back to quantum spaces on several occasions, in what follows. In fact, the present book is as much on operator algebras as it is on quantum spaces, and this because these two points of view are both useful, and complementary to each other.

\section*{5e. Exercises} 

The theory in this chapter has been quite exciting, and we have already run into a number of difficult questions. As a basic exercise on all this, we have:

\begin{exercise}
Find a simple proof for the von Neumann bicommutant theorem, in finite dimensions.
\end{exercise}

This is something quite subjective, and try not to cheat. That is, not to convert the amplification proof that we have in general, by using matrix algebras everywhere, nor by using the structure result for the finite dimensional algebras either.

\begin{exercise}
Again in finite dimensions, $H=\mathbb C^N$, compute explicitly the von Neumann algebra $<T>\subset B(H)$ generated by a single operator.
\end{exercise}

As mentioned above, in the normal case the answer is clear, by diagonalizing $T$. The problem is that of understanding what happens when $T$ is not normal.

\begin{exercise}
Try understanding what the law of the simplest non-normal operator,
$$J=\begin{pmatrix}0&1\\0&0\end{pmatrix}$$
acting on $H=\mathbb C^2$ should be. Look also at more general Jordan blocks.
\end{exercise}

There are many non-trivial computations here. We will be back to this.

\begin{exercise}
Develop a full theory of finite quantum spaces, by enlarging what has been said above, with various geometric topics, of your choice.
\end{exercise}

This is of course a bit vague, but some further thinking at all this is certainly useful, at this point, and this is what the exercise is about.

\chapter{Random matrices}

\section*{6a. Random matrices}

We have seen so far the basics of von Neumann algebras $A\subset B(H)$, with a look into some interesting ramifications too, concerning random matrices and quantum spaces. In what regards these ramifications, the situation is as follows:

\bigskip

(1) The random matrix algebras, $A=M_N(L^\infty(X))$ acting on $H=\mathbb C^N\otimes L^2(X)$, are the simplest von Neumann algebras, from a variety of viewpoints. The main problem regarding them is of operator theoretic nature, regarding the computation of the law of individual elements $T\in A$ with respect to the random matrix trace $tr:A\to\mathbb C$.

\bigskip

(2) The quantum spaces are exciting abstract objects, obtained by looking at an arbitrary von Neumann algebra $A\subset B(H)$ coming with a trace $tr:A\to\mathbb C$, and formally writing the algebra as $A=L^\infty(X)$, and its trace as $tr=\int_X$. In this picture, $X$ is our quantum probability space, and $\int_X$ is the integration over it, or expectation.

\bigskip

All this is quite interesting, and we will further explore these two topics, random matrices and quantum spaces, with some basic theory for them, in this chapter and in the next one. As a first observation, these two topics are closely related, due to:

\begin{fact}
A random matrix algebra can be written in the following way,
\begin{eqnarray*}
M_N(L^\infty(X))
&=&M_N(\mathbb C)\otimes L^\infty(X)\\
&=&L^\infty(M_N)\otimes L^\infty(X)\\
&=&L^\infty(M_N\times X)
\end{eqnarray*}
so the underlying quantum space is something very simple, $Y=M_N\times X$.
\end{fact}

With this understood, the philosophical problem is now, what to do with our quantum spaces, be them of random matrix type $Y=M_N\times X$, or more general. Good question, and do not expect a simple answer to it. Indeed, quantum spaces are more or less the same thing as operator algebras, and from this perspective, our question becomes ``what are the operator algebras, and what is to be done with them'', obviously difficult. 

\bigskip

And there is even worse, because when remembering that operator algebras are more or less the same thing as quantum mechanics, our question becomes something of type ``what is quantum mechanics, and what is to be done with it''. So, modesty. 

\bigskip

Getting back to Earth, now that we have our questions and philosophy, for the whole remainder of this book, let us get into random matrices. Quite remarkably, these provide us with an epsilon of answer to our philosophical questions, as follows:

\begin{answer}
The simplest quantum spaces are those coming from random matrix algebras, which are as follows, with $X$ being a usual probability space,
$$Y=M_N\times X$$
and what is to be done with them is the computation of the law of individual elements, the random matrices $T\in L^\infty(Y)=M_N(L^\infty(X))$, in the $N>>0$ regime.
\end{answer}

Which looks very nice, we eventually reached to some concrete questions, and time now for mathematics and computations. Getting started, we must first further build on the material from chapter 5. We recall from there that given a von Neumann algebra $A\subset B(H)$ coming with a trace $tr:A\to\mathbb C$, any normal element $T\in A$ has a law, which is the complex probability measure $\mu\in\mathcal P(\mathbb C)$ given by the following formula:
$$tr(T^k)=\int_\mathbb Cz^kd\mu(z)$$

In the non-normal case, $TT^*\neq T^*T$, the law does not exist as a complex probability measure $\mu\in\mathcal P(\mathbb C)$, as also explained in chapter 5. However, we can trick a bit, and talk about the law of non-normal elements as well, in the following abstract way:

\index{random variable}
\index{moments}
\index{law}
\index{distribution}

\begin{definition}
Let $A$ be a von Neumann algebra, given with a trace $tr:A\to\mathbb C$.
\begin{enumerate}
\item The elements $T\in A$ are called random variables.

\item The moments of such a variable are the numbers $M_k(T)=tr(T^k)$.

\item The law of such a variable is the functional $\mu:P\to tr(P(T))$.
\end{enumerate}
\end{definition}

Here $k=\circ\bullet\bullet\circ\ldots$ is by definition a colored integer, and the powers $T^k$ are defined by multiplicativity and the usual formulae, namely:
$$T^\emptyset=1\quad,\quad
T^\circ=T\quad,\quad
T^\bullet=T^*$$

As for the polynomial $P$, this is a noncommuting $*$-polynomial in one variable: 
$$P\in\mathbb C<X,X^*>$$

Observe that the law is uniquely determined by the moments, because:
$$P(X)=\sum_k\lambda_kX^k\implies\mu(P)=\sum_k\lambda_kM_k(T)$$

Generally speaking, the above definition, due to Voiculescu \cite{vdn}, is something quite abstract, but there is no other way of doing things, at least at this level of generality. However, in the special case where our variable $T\in A$ is self-adjoint, or more generally normal, the theory simplifies, and we recover more familiar objects, as follows:

\index{spectral measure}

\begin{theorem}
The law of a normal variable $T\in A$ can be identified with the corresponding spectral measure $\mu\in\mathcal P(\mathbb C)$, according to the following formula,
$$tr(f(T))=\int_{\sigma(T)}f(x)d\mu(x)$$
valid for any $f\in L^\infty(\sigma(T))$, coming from the measurable functional calculus. In the self-adjoint case the spectral measure is real, $\mu\in\mathcal P(\mathbb R)$.
\end{theorem}

\begin{proof}
This is something that we know well, from chapter 5, coming from the spectral theorem for the normal operators, as developed in chapter 3.
\end{proof}

Getting back now to the random matrices, we have all we need, as general formalism, and we are ready for doing some computations. As a first observation, we have:

\begin{theorem}
The laws of basic random matrices $T\in M_N(L^\infty(X))$ are as follows:
\begin{enumerate}
\item In the case $N=1$ the random matrix is a usual random variable, $T\in L^\infty(X)$, automatically normal, and its law as defined above is the usual law.

\item In the case $X=\{.\}$ the random matrix is a usual scalar matrix, $T\in M_N(\mathbb C)$, and in the diagonalizable case, the law is $\mu=\frac{1}{N}\left(\delta_{\lambda_1}+\ldots+\delta_{\lambda_N}\right)$.
\end{enumerate}
\end{theorem}

\begin{proof}
This is something that we know, once again, from chapter 5, and which is elementary. Indeed, the first assertion follows from definitions, and the above discussion. As for the second assertion, this follows by diagonalizing the matrix.
\end{proof}

In general, what we have can only be a mixture of (1) and (2) above. Our plan will be that of discussing more in detail (1), and then getting into the general case, or rather into the case of the most interesting random matrices, with inspiration from (2).

\section*{6b. Probability theory}

So, let us set $N=1$. Here our algebra is $A=L^\infty(X)$, an arbitrary commutative von Neumann algebra. The most interesting linear operators $T\in A$, that we will rather denote as complex functions $f:X\to\mathbb C$, and call random variables, as it is customary, are the normal, or Gaussian variables, which are defined as follows:

\index{normal variable}
\index{normal law}
\index{Gaussian law}

\begin{definition}
A variable $f:X\to\mathbb R$ is called standard normal when its law is:
$$g_1=\frac{1}{\sqrt{2\pi}}e^{-x^2/2}dx$$
More generally, the normal law of parameter $t>0$ is the following measure:
$$g_t=\frac{1}{\sqrt{2\pi t}}e^{-x^2/2t}dx$$
These are also called Gaussian distributions, with ``g'' standing for Gauss.
\end{definition}

Observe that these normal laws have indeed mass 1, as they should, as shown by a quick change of variable, and the Gauss formula, namely:
\begin{eqnarray*}
\left(\int_\mathbb Re^{-x^2}dx\right)^2
&=&\int_\mathbb R\int_\mathbb Re^{-x^2-y^2}dxdy\\
&=&\int_0^{2\pi}\int_0^\infty e^{-r^2}rdrdt\\
&=&2\pi\times\frac{1}{2}\\
&=&\pi
\end{eqnarray*}

Let us start with some basic results regarding the normal laws. We first have:

\index{moments}
\index{Fourier transform}
\index{convolution semigroup}
\index{double factorials}

\begin{proposition}
The normal law $g_t$ with $t>0$ has the following properties:
\begin{enumerate}
\item The variance is $V=t$.

\item The density is even, so the odd moments vanish.

\item The even moments are $M_k=t^{k/2}\times k!!$, with $k!!=(k-1)(k-3)(k-5)\ldots\,$.

\item Equivalently, the moments are $M_k=\sum_{\pi\in P_2(k)}t^{|\pi|}$, for any $k\in\mathbb N$.

\item The Fourier transform $F_f(x)=\mathbb E(e^{ixf})$ is given by $F(x)=e^{-tx^2/2}$.

\item We have the convolution semigroup formula $g_s*g_t=g_{s+t}$, for any $s,t>0$.
\end{enumerate}
\end{proposition}

\begin{proof}
All this is very standard, with the various notations used in the statement being explained below, the idea being as follows:

\medskip

(1) The normal law $g_t$ being centered, its variance is the second moment, $V=M_2$. Thus the result follows from (3), proved below, which gives in particular: 
$$M_2=t^{2/2}\times 2!!=t$$

(2) This is indeed something self-explanatory.

\medskip

(3) We have indeed the following computation, by partial integration:
\begin{eqnarray*}
M_k
&=&\frac{1}{\sqrt{2\pi t}}\int_\mathbb Rx^ke^{-x^2/2t}dx\\
&=&\frac{1}{\sqrt{2\pi t}}\int_\mathbb R(tx^{k-1})\left(-e^{-x^2/2t}\right)'dx\\
&=&\frac{1}{\sqrt{2\pi t}}\int_\mathbb Rt(k-1)x^{k-2}e^{-x^2/2t}dx\\
&=&t(k-1)\times\frac{1}{\sqrt{2\pi t}}\int_\mathbb Rx^{k-2}e^{-x^2/2t}dx\\
&=&t(k-1)M_{k-2}
\end{eqnarray*}

The initial value being $M_0=1$, we obtain the result.

\medskip

(4) We know from (2,3) that the moments of the normal law $g_t$ satisfy the following recurrence formula, with the initial data $M_0=1,M_1=0$:
$$M_k=t(k-1)M_{k-2}$$

Now let us look at $P_2(k)$, the set of pairings of $\{1,\ldots,k\}$. In order to have such a pairing, we must pair 1 with a number chosen among $2,\ldots,k$, and then come up with a pairing of the remaining $k-2$ numbers. Thus, the number $N_k=|P_2(k)|$ of such pairings is subject to the following recurrence formula, with initial data $N_0=1,N_1=0$:
$$N_k=(k-1)N_{k-2}$$

But this solves our problem at $t=1$, because in this case we obtain the following formula, with $|.|$ standing as usual for the number of blocks of a partition:
$$M_k=N_k=|P_2(k)|=\sum_{\pi\in P_2(k)}1=\sum_{\pi\in P_2(k)}1^{|\pi|}$$

Now back to the general case, $t>0$, our problem here is solved in fact too, because the number of blocks of a pairing $\pi\in P_2(k)$ being constant, $|\pi|=k/2$, we obtain:
$$M_k=t^{k/2}N_k=\sum_{\pi\in P_2(k)}t^{k/2}=\sum_{\pi\in P_2(k)}t^{|\pi|}$$

(5) The Fourier transform formula can be established as follows:
\begin{eqnarray*}
F(x)
&=&\frac{1}{\sqrt{2\pi t}}\int_\mathbb Re^{-y^2/2t+ixy}dy\\
&=&\frac{1}{\sqrt{2\pi t}}\int_\mathbb Re^{-(y/\sqrt{2t}-\sqrt{t/2}ix)^2-tx^2/2}dy\\
&=&\frac{1}{\sqrt{2\pi t}}\int_\mathbb Re^{-z^2-tx^2/2}\sqrt{2t}dz\\
&=&\frac{1}{\sqrt{\pi}}e^{-tx^2/2}\int_\mathbb Re^{-z^2}dz\\
&=&e^{-tx^2/2}
\end{eqnarray*}

(6) This follows indeed from (5), because $\log F_{g_t}$ is linear in $t$.
\end{proof}

We are now ready to establish the Central Limit Theorem (CLT), which is a key result, telling us why the normal laws appear a bit everywhere, in the real life:

\index{CLT}
\index{Central Limit Theorem}

\begin{theorem}
Given a sequence of real random variables $f_1,f_2,f_3,\ldots\in L^\infty(X)$, which are i.i.d., centered, and with variance $t>0$, we have
$$\frac{1}{\sqrt{n}}\sum_{i=1}^nf_i\sim g_t$$
with $n\to\infty$, in moments.
\end{theorem}

\begin{proof}
In terms of moments, the Fourier transform $F_f(x)=\mathbb E(e^{ixf})$ is given by:
$$F_f(x)
=\mathbb E\left(\sum_{k=0}^\infty\frac{(ixf)^k}{k!}\right)
=\sum_{k=0}^\infty\frac{i^kM_k(f)}{k!}\,x^k$$

Thus, the Fourier transform of the variable in the statement is:
\begin{eqnarray*}
F(x)
&=&\left[F_f\left(\frac{x}{\sqrt{n}}\right)\right]^n\\
&=&\left[1-\frac{tx^2}{2n}+O(n^{-2})\right]^n\\
&\simeq&\left[1-\frac{tx^2}{2n}\right]^n\\
&\simeq&e^{-tx^2/2}
\end{eqnarray*}

But this latter function being the Fourier transform of $g_t$, we obtain the result.
\end{proof}

Let us discuss as well the ``discrete'' counterpart of the above results, that we will need too a bit later, in relation with the random matrices. We have:

\index{Poisson law}

\begin{definition}
The Poisson law of parameter $1$ is the following measure,
$$p_1=\frac{1}{e}\sum_k\frac{\delta_k}{k!}$$
and the Poisson law of parameter $t>0$ is the following measure,
$$p_t=e^{-t}\sum_k\frac{t^k}{k!}\,\delta_k$$
with the letter ``p'' standing for Poisson.
\end{definition}

We will see in a moment why these laws appear everywhere, in discrete probability, the reasons behind this coming from the Poisson Limit Theorem (PLT). Getting started now, in analogy with the normal laws, the Poisson laws have the following properties:

\index{Fourier transform}
\index{convolution semigroup}
\index{Bell numbers}
\index{Stirling numbers}

\begin{proposition}
The Poisson law $p_t$ with $t>0$ has the following properties:
\begin{enumerate}
\item The variance is $V=t$.

\item The moments are $M_k=\sum_{\pi\in P(k)}t^{|\pi|}$.

\item The Fourier transform is $F(x)=\exp\left((e^{ix}-1)t\right)$.

\item We have the semigroup formula $p_s*p_t=p_{s+t}$, for any $s,t>0$.
\end{enumerate}
\end{proposition}

\begin{proof}
We have four formulae to be proved, the idea being as follows:

\medskip

(1) The variance is $V=M_2-M_1^2$, and by using the formulae $M_1=t$ and $M_2=t+t^2$, coming from (2), proved below, we obtain as desired, $V=t$.

\medskip

(2) This is something more tricky. Consider indeed the set $P(k)$ of all partitions of $\{1,\ldots,k\}$. At $t=1$, to start with, the formula that we want to prove is:
$$M_k=|P(k)|$$

We have the following recurrence formula for the moments of $p_1$:
\begin{eqnarray*}
M_{k+1}
&=&\frac{1}{e}\sum_s\frac{(s+1)^{k+1}}{(s+1)!}\\
&=&\frac{1}{e}\sum_s\frac{s^k}{s!}\left(1+\frac{1}{s}\right)^k\\
&=&\frac{1}{e}\sum_s\frac{s^k}{s!}\sum_r\binom{k}{r}s^{-r}\\
&=&\sum_r\binom{k}{r}\cdot\frac{1}{e}\sum_s\frac{s^{k-r}}{s!}\\
&=&\sum_r\binom{k}{r}M_{k-r}
\end{eqnarray*}

Our claim is that the numbers $B_k=|P(k)|$ satisfy the same recurrence formula. Indeed, since a partition of $\{1,\ldots,k+1\}$ appears by choosing $r$ neighbors for $1$, among the $k$ numbers available, and then partitioning the $k-r$ elements left, we have:
$$B_{k+1}=\sum_r\binom{k}{r}B_{k-r}$$

Thus we obtain by recurrence $M_k=B_k$, as desired. Regarding now the general case, $t>0$, we can use here a similar method. We have the following recurrence formula for the moments of $p_t$, obtained by using the binomial formula:
\begin{eqnarray*}
M_{k+1}
&=&e^{-t}\sum_s\frac{t^{s+1}(s+1)^{k+1}}{(s+1)!}\\
&=&e^{-t}\sum_s\frac{t^{s+1}s^k}{s!}\left(1+\frac{1}{s}\right)^k\\
&=&e^{-t}\sum_s\frac{t^{s+1}s^k}{s!}\sum_r\binom{k}{r}s^{-r}\\
&=&\sum_r\binom{k}{r}\cdot e^{-t}\sum_s\frac{t^{s+1}s^{k-r}}{s!}\\
&=&t\sum_r\binom{k}{r}M_{k-r}
\end{eqnarray*}

On the other hand, consider the numbers in the statement, $S_k=\sum_{\pi\in P(k)}t^{|\pi|}$. As before, since a partition of $\{1,\ldots,k+1\}$ appears by choosing $r$ neighbors for $1$, among the $k$ numbers available, and then partitioning the $k-r$ elements left, we have:
$$S_{k+1}=t\sum_r\binom{k}{r}S_{k-r}$$

Thus we obtain by recurrence $M_k=B_k$, as desired. 

\medskip

(3) The Fourier transform formula can be established as follows:
\begin{eqnarray*}
F_{p_t}(x)
&=&e^{-t}\sum_k\frac{t^k}{k!}F_{\delta_k}(x)\\
&=&e^{-t}\sum_k\frac{t^k}{k!}\,e^{ikx}\\
&=&e^{-t}\sum_k\frac{(e^{ix}t)^k}{k!}\\
&=&\exp(-t)\exp(e^{ix}t)\\
&=&\exp\left((e^{ix}-1)t\right)
\end{eqnarray*}

(4) This follows from (3), because $\log F_{p_t}$ is linear in $t$.
\end{proof}

We are now ready to establish the Poisson Limit Theorem (PLT), as follows:

\index{PLT}
\index{Poisson Limit Theorem}
\index{Bernoulli law}

\begin{theorem}
We have the following convergence, in moments,
$$\left(\left(1-\frac{t}{n}\right)\delta_0+\frac{t}{n}\delta_1\right)^{*n}\to p_t$$
for any $t>0$.
\end{theorem}

\begin{proof}
Let us denote by $\mu_n$ the Bernoulli measure appearing under the convolution sign. We have then the following computation: 
\begin{eqnarray*}
F_{\delta_r}(x)=e^{irx}
&\implies&F_{\mu_n}(x)=\left(1-\frac{t}{n}\right)+\frac{t}{n}e^{ix}\\
&\implies&F_{\mu_n^{*n}}(x)=\left(\left(1-\frac{t}{n}\right)+\frac{t}{n}e^{ix}\right)^n\\
&\implies&F_{\mu_n^{*n}}(x)=\left(1+\frac{(e^{ix}-1)t}{n}\right)^n\\
&\implies&F(x)=\exp\left((e^{ix}-1)t\right)
\end{eqnarray*}

Thus, we obtain the Fourier transform of $p_t$, as desired.
\end{proof}

As a third and last topic from classical probability, let us discuss now the complex normal laws, that we will need too. To start with, we have the following definition:

\index{complex Gaussian law}
\index{complex normal law}

\begin{definition}
The complex Gaussian law of parameter $t>0$ is
$$G_t=law\left(\frac{1}{\sqrt{2}}(a+ib)\right)$$
where $a,b$ are independent, each following the law $g_t$.
\end{definition}

As in the real case, these measures form convolution semigroups:

\index{convolution semigroup}

\begin{proposition}
The complex Gaussian laws have the property
$$G_s*G_t=G_{s+t}$$
for any $s,t>0$, and so they form a convolution semigroup.
\end{proposition}

\begin{proof}
This follows indeed from the real result, namely $g_s*g_t=g_{s+t}$, established above, simply by taking real and imaginary parts.
\end{proof}

We have the following complex analogue of the CLT:

\index{CCLT}
\index{Complex CLT}

\begin{theorem}[CCLT]
Given complex random variables $f_1,f_2,f_3,\ldots\in L^\infty(X)$ which are i.i.d., centered, and with variance $t>0$, we have, with $n\to\infty$, in moments,
$$\frac{1}{\sqrt{n}}\sum_{i=1}^nf_i\sim G_t$$
where $G_t$ is the complex Gaussian law of parameter $t$.
\end{theorem}

\begin{proof}
This follows indeed from the real CLT, established above, simply by taking the real and imaginary parts of all the variables involved.
\end{proof}

Regarding now the moments, we use the general formalism from Definition 6.3, involving colored integer exponents $k=\circ\bullet\bullet\circ\ldots\,$ We say that a pairing $\pi\in P_2(k)$ is matching when it pairs $\circ-\bullet$ symbols. With this convention, we have the following result:

\index{matching pairings}

\begin{theorem}
The moments of the complex normal law are the numbers
$$M_k(G_t)=\sum_{\pi\in\mathcal P_2(k)}t^{|\pi|}$$
where $\mathcal P_2(k)$ are the matching pairings of $\{1,\ldots,k\}$, and $|.|$ is the number of blocks.
\end{theorem}

\begin{proof}
This is something well-known, which can be established as follows:

\medskip

(1) As a first observation, by using a standard dilation argument, it is enough to do this at $t=1$. So, let us first recall from the above that the moments of the real Gaussian law $g_1$, with respect to integer exponents $k\in\mathbb N$, are the following numbers:
$$m_k=|P_2(k)|$$

Numerically, we have the following formula, explained as well in the above:
$$m_k=\begin{cases}
k!!&(k\ {\rm even})\\
0&(k\ {\rm odd})
\end{cases}$$

(2) We will show here that in what concerns the complex Gaussian law $G_1$, similar results hold. Numerically, we will prove that we have the following formula, where a colored integer $k=\circ\bullet\bullet\circ\ldots$ is called uniform when it contains the same number of $\circ$ and $\bullet$\,, and where $|k|\in\mathbb N$ is the length of such a colored integer:
$$M_k=\begin{cases}
(|k|/2)!&(k\ {\rm uniform})\\
0&(k\ {\rm not\ uniform})
\end{cases}$$

Now since the matching partitions $\pi\in\mathcal P_2(k)$ are counted by exactly the same numbers, and this for trivial reasons, we will obtain the formula in the statement, namely:
$$M_k=|\mathcal P_2(k)|$$

(3) This was for the plan. In practice now, we must compute the moments, with respect to colored integer exponents $k=\circ\bullet\bullet\circ\ldots$\,, of the variable in the statement:
$$c=\frac{1}{\sqrt{2}}(a+ib)$$

As a first observation, in the case where such an exponent $k=\circ\bullet\bullet\circ\ldots$ is not uniform in $\circ,\bullet$\,, a rotation argument shows that the corresponding moment of $c$ vanishes. To be more precise, the variable $c'=wc$ can be shown to be complex Gaussian too, for any $w\in\mathbb C$, and from $M_k(c)=M_k(c')$ we obtain $M_k(c)=0$, in this case.

\medskip

(4) In the uniform case now, where $k=\circ\bullet\bullet\circ\ldots$ consists of $p$ copies of $\circ$ and $p$ copies of $\bullet$\,, the corresponding moment can be computed as follows:
\begin{eqnarray*}
M_k
&=&\frac{1}{2^p}\int(a^2+b^2)^p\\
&=&\frac{1}{2^p}\sum_s\binom{p}{s}\int a^{2s}\int b^{2p-2s}\\
&=&\frac{1}{2^p}\sum_s\binom{p}{s}(2s)!!(2p-2s)!!\\
&=&\frac{1}{2^p}\sum_s\frac{p!}{s!(p-s)!}\cdot\frac{(2s)!}{2^ss!}\cdot\frac{(2p-2s)!}{2^{p-s}(p-s)!}\\
&=&\frac{p!}{4^p}\sum_s\binom{2s}{s}\binom{2p-2s}{p-s}
\end{eqnarray*}

(5) In order to finish now the computation, let us recall that we have the following formula, coming from the generalized binomial formula, or from the Taylor formula:
$$\frac{1}{\sqrt{1+t}}=\sum_{k=0}^\infty\binom{2k}{k}\left(\frac{-t}{4}\right)^k$$

By taking the square of this series, we obtain the following formula:
\begin{eqnarray*}
\frac{1}{1+t}
&=&\sum_{ks}\binom{2k}{k}\binom{2s}{s}\left(\frac{-t}{4}\right)^{k+s}\\
&=&\sum_p\left(\frac{-t}{4}\right)^p\sum_s\binom{2s}{s}\binom{2p-2s}{p-s}
\end{eqnarray*}

Now by looking at the coefficient of $t^p$ on both sides, we conclude that the sum on the right equals $4^p$. Thus, we can finish the moment computation in (4), as follows:
$$M_p=\frac{p!}{4^p}\times 4^p=p!$$

(6) As a conclusion, if we denote by $|k|$ the length of a colored integer $k=\circ\bullet\bullet\circ\ldots$\,, the moments of the variable $c$ in the statement are given by:
$$M_k=\begin{cases}
(|k|/2)!&(k\ {\rm uniform})\\
0&(k\ {\rm not\ uniform})
\end{cases}$$

On the other hand, the numbers $|\mathcal P_2(k)|$ are given by exactly the same formula. Indeed, in order to have matching pairings of $k$, our exponent $k=\circ\bullet\bullet\circ\ldots$ must be uniform, consisting of $p$ copies of $\circ$ and $p$ copies of $\bullet$, with $p=|k|/2$. But then the matching pairings of $k$ correspond to the permutations of the $\bullet$ symbols, as to be matched with $\circ$ symbols, and so we have $p!$ such matching pairings. Thus, we have the same formula as for the moments of $c$, and we are led to the conclusion in the statement.
\end{proof}

This was for the basic probability theory, which is in a certain sense advanced operator theory, inside the commutative von Neumann algebras, $A=L^\infty(X)$. We will be back to this, with some further limiting theorems, in chapter 8 below.

\section*{6c. Wigner matrices}

Let us exit now the classical world, that of the commutative von Neumann algebras $A=L^\infty(X)$, and do as promised some random matrix theory. We recall that a random matrix algebra is a von Neumann algebra of type $A=M_N(L^\infty(X))$, and that we are interested in the computation of the laws of the operators $T\in A$, called random matrices. Regarding the precise classes of random matrices that we are interested in, first we have the complex Gaussian matrices, which are constructed as follows:

\index{Gaussian matrix}

\begin{definition}
A complex Gaussian matrix is a random matrix of type
$$Z\in M_N(L^\infty(X))$$
which has i.i.d. complex normal entries.
\end{definition}

We will see that the above matrices have an interesting, and ``central'' combinatorics, among all kinds of random matrices, with the study of the other random matrices being usually obtained as a modification of the study of the Gaussian matrices.

\bigskip

As a somewhat surprising remark, using real normal variables in Definition 6.16, instead of the complex ones appearing there, leads nowhere. The correct real versions of the Gaussian matrices are the Wigner random matrices, constructed as follows: 

\index{Wigner matrix}

\begin{definition}
A Wigner matrix is a random matrix of type
$$Z\in M_N(L^\infty(X))$$
which has i.i.d. complex normal entries, up to the constraint $Z=Z^*$.
\end{definition}

In other words, a Wigner matrix must be as follows, with the diagonal entries being real normal variables, $a_i\sim g_t$, for some $t>0$, the upper diagonal entries being complex normal variables, $b_{ij}\sim G_t$, the lower diagonal entries being the conjugates of the upper diagonal entries, as indicated, and with all the variables $a_i,b_{ij}$ being independent: 
$$Z=\begin{pmatrix}
a_1&b_{12}&\ldots&\ldots&b_{1N}\\
\bar{b}_{12}&a_2&\ddots&&\vdots\\
\vdots&\ddots&\ddots&\ddots&\vdots\\
\vdots&&\ddots&a_{N-1}&b_{N-1,N}\\
\bar{b}_{1N}&\ldots&\ldots&\bar{b}_{N-1,N}&a_N
\end{pmatrix}$$

As a comment here, for many concrete applications the Wigner matrices are in fact the central objects in random matrix theory, and in particular, they are often more important than the Gaussian matrices. In fact, these are the random matrices which were first considered and investigated, a long time ago, by Wigner himself \cite{wig}.

\bigskip

Finally, we will be interested as well in the complex Wishart matrices, which are the positive versions of the above random matrices, constructed as follows: 

\index{Wishart matrix}

\begin{definition}
A complex Wishart matrix is a random matrix of type
$$Z=YY^*\in M_N(L^\infty(X))$$
with $Y$ being a complex Gaussian matrix.
\end{definition}

As before with the Gaussian and Wigner matrices, there are many possible comments that can be made here, of technical or historical nature. First, using real Gaussian variables instead of complex ones leads to a less interesting combinatorics. Also, these matrices were introduced and studied by Marchenko-Pastur not long after Wigner, in \cite{mpa}, and so historically came second. Finally, in what regards their combinatorics and applications, these matrices quite often come first, before both the Gaussian and the Wigner ones, with all this being of course a matter of knowledge and taste.

\bigskip

Summarizing, we have three main types of random matrices, which can be somehow designated as ``complex'', ``real'' and ``positive'', and that we will study in what follows. Let us also mention that there are many other interesting classes of random matrices, usually appearing as modifications of the above. More on these later.

\bigskip

In order to compute the asymptotic laws of the above matrices, we will use the moment method. We have the following result, which will be our main tool here:

\index{complex normal law}
\index{colored moments}
\index{matching pairings}
\index{Wick formula}

\begin{theorem}
Given independent variables $X_i$, each following the complex normal law $G_t$, with $t>0$ being a fixed parameter, we have the Wick formula
$$\mathbb E\left(X_{i_1}^{k_1}\ldots X_{i_s}^{k_s}\right)=t^{s/2}\#\left\{\pi\in\mathcal P_2(k)\Big|\pi\leq\ker i\right\}$$
where $k=k_1\ldots k_s$ and $i=i_1\ldots i_s$, for the joint moments of these variables.
\end{theorem}

\begin{proof}
This is something well-known, and the basis for all possible computations with complex normal variables, which can be proved in two steps, as follows:

\medskip

(1) Let us first discuss the case where we have a single complex normal variable $X$, which amounts in taking $X_i=X$ for any $i$ in the formula in the statement. What we have to compute here are the moments of $X$, with respect to colored integer exponents $k=\circ\bullet\bullet\circ\ldots\,$, and the formula in the statement tells us that these moments must be:
$$\mathbb E(X^k)=t^{|k|/2}|\mathcal P_2(k)|$$

But this is something that we know well from the above, the idea being that at $t=1$ this follows by doing some combinatorics and calculus, in analogy with the combinatorics and calculus from the real case, where the moment formula is identical, save for the matching pairings $\mathcal P_2$ being replaced by the usual pairings $P_2$, and then that the general case $t>0$ follows from this, by rescaling. Thus, we are done with this case.

\medskip

(2) In general now, the point is that we obtain the formula in the statement. Indeed, when expanding the product $X_{i_1}^{k_1}\ldots X_{i_s}^{k_s}$ and rearranging the terms, we are left with doing a number of computations as in (1), and then making the product of the expectations that we found. But this amounts precisely in counting the partitions in the statement, with the condition $\pi\leq\ker i$ there standing precisely for the fact that we are doing the various type (1) computations independently, and then making the product.
\end{proof}

Now by getting back to the Gaussian matrices, we have the following result, with $\mathcal{NC}_2(k)=\mathcal P_2(k)\cap NC(k)$ standing for the noncrossing pairings of a colored integer $k$:

\index{Gaussian matrix}

\begin{theorem}
Given a sequence of Gaussian random matrices
$$Z_N\in M_N(L^\infty(X))$$
having independent $G_t$ variables as entries, for some fixed $t>0$, we have
$$M_k\left(\frac{Z_N}{\sqrt{N}}\right)\simeq t^{|k|/2}|\mathcal{NC}_2(k)|$$
for any colored integer $k=\circ\bullet\bullet\circ\ldots\,$, in the $N\to\infty$ limit.
\end{theorem}

\begin{proof}
This is something standard, which can be done as follows:

\medskip

(1) We fix $N\in\mathbb N$, and we let $Z=Z_N$. Let us first compute the trace of $Z^k$. With $k=k_1\ldots k_s$, and with the convention $(ij)^\circ=ij,(ij)^\bullet=ji$, we have:
\begin{eqnarray*}
Tr(Z^k)
&=&Tr(Z^{k_1}\ldots Z^{k_s})\\
&=&\sum_{i_1=1}^N\ldots\sum_{i_s=1}^N(Z^{k_1})_{i_1i_2}(Z^{k_2})_{i_2i_3}\ldots(Z^{k_s})_{i_si_1}\\
&=&\sum_{i_1=1}^N\ldots\sum_{i_s=1}^N(Z_{(i_1i_2)^{k_1}})^{k_1}(Z_{(i_2i_3)^{k_2}})^{k_2}\ldots(Z_{(i_si_1)^{k_s}})^{k_s}
\end{eqnarray*}

(2) Next, we rescale our variable $Z$ by a $\sqrt{N}$ factor, as in the statement, and we also replace the usual trace by its normalized version, $tr=Tr/N$. Our formula becomes:
$$tr\left(\left(\frac{Z}{\sqrt{N}}\right)^k\right)=\frac{1}{N^{s/2+1}}\sum_{i_1=1}^N\ldots\sum_{i_s=1}^N(Z_{(i_1i_2)^{k_1}})^{k_1}(Z_{(i_2i_3)^{k_2}})^{k_2}\ldots(Z_{(i_si_1)^{k_s}})^{k_s}$$

Thus, the moment that we are interested in is given by:
$$M_k\left(\frac{Z}{\sqrt{N}}\right)=\frac{1}{N^{s/2+1}}\sum_{i_1=1}^N\ldots\sum_{i_s=1}^N\int_X(Z_{(i_1i_2)^{k_1}})^{k_1}(Z_{(i_2i_3)^{k_2}})^{k_2}\ldots(Z_{(i_si_1)^{k_s}})^{k_s}$$

(3) Let us apply now the Wick formula, from Theorem 6.19. We conclude that the moment that we are interested in is given by the following formula:
\begin{eqnarray*}
&&M_k\left(\frac{Z}{\sqrt{N}}\right)\\
&=&\frac{t^{s/2}}{N^{s/2+1}}\sum_{i_1=1}^N\ldots\sum_{i_s=1}^N\#\left\{\pi\in\mathcal P_2(k)\Big|\pi\leq\ker\left((i_1i_2)^{k_1},(i_2i_3)^{k_2},\ldots,(i_si_1)^{k_s}\right)\right\}\\
&=&t^{s/2}\sum_{\pi\in\mathcal P_2(k)}\frac{1}{N^{s/2+1}}\#\left\{i\in\{1,\ldots,N\}^s\Big|\pi\leq\ker\left((i_1i_2)^{k_1},(i_2i_3)^{k_2},\ldots,(i_si_1)^{k_s}\right)\right\}
\end{eqnarray*}

(4) Our claim now is that in the $N\to\infty$ limit the combinatorics of the above sum simplifies, with only the noncrossing partitions contributing to the sum, and with each of them contributing precisely with a 1 factor, so that we will have, as desired:
\begin{eqnarray*}
M_k\left(\frac{Z}{\sqrt{N}}\right)
&=&t^{s/2}\sum_{\pi\in\mathcal P_2(k)}\Big(\delta_{\pi\in NC_2(k)}+O(N^{-1})\Big)\\
&\simeq&t^{s/2}\sum_{\pi\in\mathcal P_2(k)}\delta_{\pi\in NC_2(k)}\\
&=&t^{s/2}|\mathcal{NC}_2(k)|
\end{eqnarray*}

(5) In order to prove this, the first observation is that when $k$ is not uniform, in the sense that it contains a different number of $\circ$, $\bullet$ symbols, we have $\mathcal P_2(k)=\emptyset$, and so:
$$M_k\left(\frac{Z}{\sqrt{N}}\right)=t^{s/2}|\mathcal{NC}_2(k)|=0$$

(6) Thus, we are left with the case where $k$ is uniform. Let us examine first the case where $k$ consists of an alternating sequence of $\circ$ and $\bullet$ symbols, as follows:
$$k=\underbrace{\circ\bullet\circ\bullet\ldots\ldots\circ\bullet}_{2p}$$

In this case it is convenient to relabel our multi-index $i=(i_1,\ldots,i_s)$, with $s=2p$, in the form $(j_1,l_1,j_2,l_2,\ldots,j_p,l_p)$. With this done, our moment formula becomes:
$$M_k\left(\frac{Z}{\sqrt{N}}\right)
=t^p\sum_{\pi\in\mathcal P_2(k)}\frac{1}{N^{p+1}}\#\left\{j,l\in\{1,\ldots,N\}^p\Big|\pi\leq\ker\left(j_1l_1,j_2l_1,j_2l_2,\ldots,j_1l_p\right)\right\}$$

Now observe that, with $k$ being as above, we have an identification $\mathcal P_2(k)\simeq S_p$, obtained in the obvious way. With this done too, our moment formula becomes:
$$M_k\left(\frac{Z}{\sqrt{N}}\right)
=t^p\sum_{\pi\in S_p}\frac{1}{N^{p+1}}\#\left\{j,l\in\{1,\ldots,N\}^p\Big|j_r=j_{\pi(r)+1},l_r=l_{\pi(r)},\forall r\right\}$$

(7) We are now ready to do our asymptotic study, and prove the claim in (4). Let indeed $\gamma\in S_p$ be the full cycle, which is by definition the following permutation:
$$\gamma=(1 \, 2 \, \ldots \, p)$$

In terms of $\gamma$, the conditions $j_r=j_{\pi(r)+1}$ and $l_r=l_{\pi(r)}$ found above read:
$$\gamma\pi\leq\ker j\quad,\quad 
\pi\leq\ker l$$

Counting the number of free parameters in our moment formula, we obtain:
$$M_k\left(\frac{Z}{\sqrt{N}}\right)
=\frac{t^p}{N^{p+1}}\sum_{\pi\in S_p}N^{|\pi|+|\gamma\pi|}
=t^p\sum_{\pi\in S_p}N^{|\pi|+|\gamma\pi|-p-1}$$

(8) The point now is that the last exponent is well-known to be $\leq 0$, with equality precisely when the permutation $\pi\in S_p$ is geodesic, which in practice means that $\pi$ must come from a noncrossing partition. Thus we obtain, in the $N\to\infty$ limit, as desired:
$$M_k\left(\frac{Z}{\sqrt{N}}\right)\simeq t^p|\mathcal{NC}_2(k)|$$

This finishes the proof in the case of the exponents $k$ which are alternating, and the case where $k$ is an arbitrary uniform exponent is similar, by permuting everything.
\end{proof}

As a conclusion to this, we have obtained as asymptotic law for the Gaussian matrices a certain mysterious distribution, having as moments some  numbers which are similar to the moments of the usual normal laws, but with the ``underlying matching pairings being now replaced by underlying matching noncrossing pairings''. More on this later.

\bigskip
 
Regarding now the Wigner matrices, we have here the following result, coming as a consequence of Theorem 6.20, via some simple algebraic manipulations:

\begin{theorem}
Given a sequence of Wigner random matrices
$$Z_N\in M_N(L^\infty(X))$$
having independent $G_t$ variables as entries, with $t>0$, up to $Z_N=Z_N^*$, we have
$$M_k\left(\frac{Z_N}{\sqrt{N}}\right)\simeq t^{k/2}|NC_2(k)|$$
for any integer $k\in\mathbb N$, in the $N\to\infty$ limit.
\end{theorem}

\begin{proof}
This can be deduced from a direct computation based on the Wick formula, similar to that from the proof of Theorem 6.20, but the best is to deduce this result from Theorem 6.20 itself. Indeed, we know from there that for Gaussian matrices $Y_N\in M_N(L^\infty(X))$ we have the following formula, valid for any colored integer $K=\circ\bullet\bullet\circ\ldots\,$, in the $N\to\infty$ limit, with $\mathcal{NC}_2$ standing for noncrossing matching pairings:
$$M_K\left(\frac{Y_N}{\sqrt{N}}\right)\simeq t^{|K|/2}|\mathcal{NC}_2(K)|$$

By doing some combinatorics, we deduce from this that we have the following formula for the moments of the matrices $Re(Y_N)$, with respect to usual exponents, $k\in\mathbb N$:
\begin{eqnarray*}
M_k\left(\frac{Re(Y_N)}{\sqrt{N}}\right)
&=&2^{-k}\cdot M_k\left(\frac{Y_N}{\sqrt{N}}+\frac{Y_N^*}{\sqrt{N}}\right)\\
&=&2^{-k}\sum_{|K|=k}M_K\left(\frac{Y_N}{\sqrt{N}}\right)\\
&\simeq&2^{-k}\sum_{|K|=k}t^{k/2}|\mathcal{NC}_2(K)|\\
&=&2^{-k}\cdot t^{k/2}\cdot 2^{k/2}|\mathcal{NC}_2(k)|\\
&=&2^{-k/2}\cdot t^{k/2}|NC_2(k)|
\end{eqnarray*}

Now since the matrices $Z_N=\sqrt{2}Re(Y_N)$ are of Wigner type, this gives the result.
\end{proof}

Summarizing, all this brings us into counting noncrossing pairings. So, let us start with some preliminaries here. We first have the following well-known result:

\index{Catalan numbers}

\begin{theorem}
The Catalan numbers, which are by definition given by
$$C_k=|NC_2(2k)|$$
satisfy the following recurrence formula, with initial data $C_0=C_1=1$,
$$C_{k+1}=\sum_{a+b=k}C_aC_b$$ 
their generating series $f(z)=\sum_{k\geq0}C_kz^k$ satisfies the equation
$$zf^2-f+1=0$$
and is given by the following explicit formula,
$$f(z)=\frac{1-\sqrt{1-4z}}{2z}$$ 
and we have the following explicit formula for these numbers:
$$C_k=\frac{1}{k+1}\binom{2k}{k}$$
Numerically, these numbers are $1,1,2,5,14,42,132,429,1430,4862,16796,\ldots$
\end{theorem}

\begin{proof}
We must count the noncrossing pairings of $\{1,\ldots,2k\}$. Now observe that such a pairing appears by pairing 1 to an odd number, $2a+1$, and then inserting a noncrossing pairing of $\{2,\ldots,2a\}$, and a noncrossing pairing of $\{2a+2,\ldots,2l\}$. We conclude that we have the following recurrence formula for the Catalan numbers:
$$C_k=\sum_{a+b=k-1}C_aC_b$$ 

In terms of the generating series $f(z)=\sum_{k\geq0}C_kz^k$, this recurrence formula reads:
\begin{eqnarray*}
zf^2
&=&\sum_{a,b\geq0}C_aC_bz^{a+b+1}\\
&=&\sum_{k\geq1}\sum_{a+b=k-1}C_aC_bz^k\\
&=&\sum_{k\geq1}C_kz^k\\
&=&f-1
\end{eqnarray*}

Thus $f$ satisfies $zf^2-f+1=0$, and by solving this equation, and choosing the solution which is bounded at $z=0$, we obtain the following formula:
$$f(z)=\frac{1-\sqrt{1-4z}}{2z}$$ 

In order to finish, we use the generalized binomial formula, which gives:
$$\sqrt{1+t}=1-2\sum_{k=1}^\infty\frac{1}{k}\binom{2k-2}{k-1}\left(\frac{-t}{4}\right)^k$$

Now back to our series $f$, we obtain the following formula for it:
\begin{eqnarray*}
f(z)
&=&\frac{1-\sqrt{1-4z}}{2z}\\
&=&\sum_{k=1}^\infty\frac{1}{k}\binom{2k-2}{k-1}z^{k-1}\\
&=&\sum_{k=0}^\infty\frac{1}{k+1}\binom{2k}{k}z^k
\end{eqnarray*}

It follows that the Catalan numbers are given by:
$$C_k=\frac{1}{k+1}\binom{2k}{k}$$

Thus, we are led to the conclusion in the statement.
\end{proof}

In order to recapture now the Wigner measure from its moments, we can use:

\index{semicircle law}
\index{Catalan numbers}

\begin{proposition}
The Catalan numbers are the even moments of 
$$\gamma_1=\frac{1}{2\pi}\sqrt{4-x^2}dx$$
called standard semicircle law. As for the odd moments of $\gamma_1$, these all vanish. 
\end{proposition}

\begin{proof}
The even moments of the semicircle law in the statement can be computed with the change of variable $x=2\cos t$, and we are led to the following formula:
\begin{eqnarray*}
M_{2k}
&=&\frac{1}{\pi}\int_0^2\sqrt{4-x^2}x^{2k}dx\\
&=&\frac{1}{\pi}\int_0^{\pi/2}\sqrt{4-4\cos^2t}\,(2\cos t)^{2k}2\sin t\,dt\\
&=&\frac{4^{k+1}}{\pi}\int_0^{\pi/2}\cos^{2k}t\sin^2t\,dt\\
&=&\frac{4^{k+1}}{\pi}\cdot\frac{\pi}{2}\cdot\frac{(2k)!!2!!}{(2k+3)!!}\\
&=&2\cdot 4^k\cdot\frac{(2k)!/2^kk!}{2^{k+1}(k+1)!}\\
&=&C_k
\end{eqnarray*}

As for the odd moments, these all vanish, because the density of $\gamma_1$ is an even function. Thus, we are led to the conclusion in the statement.
\end{proof}

More generally, we have the following result, involving a parameter $t>0$:

\begin{proposition}
Given $t>0$, the real measure having as even moments the numbers $M_{2k}=t^kC_k$ and having all odd moments $0$ is the measure
$$\gamma_t=\frac{1}{2\pi t}\sqrt{4t-x^2}dx$$
called Wigner semicircle law on $[-2\sqrt{t},2\sqrt{t}]$.
\end{proposition}

\begin{proof}
This follows indeed from Proposition 6.23, via a change of variables.
\end{proof}

Now by putting everything together, we obtain the Wigner theorem, as follows:

\begin{theorem}
Given a sequence of Wigner random matrices
$$Z_N\in M_N(L^\infty(X))$$
which by definition have i.i.d. complex normal entries, up to $Z_N=Z_N^*$, we have
$$Z_N\sim\gamma_t$$
in the $N\to\infty$ limit, where $\gamma_t=\frac{1}{2\pi t}\sqrt{4t-x^2}dx$ is the Wigner semicircle law. 
\end{theorem}

\begin{proof}
This follows indeed from all the above, and more specifically, by combining Theorem 6.21, Theorem 6.22 and Proposition 6.24.
\end{proof}

Regarding now the complex Gaussian matrices, in view of this result, it is natural to think at the law found in Theorem 6.20 as being ``circular''. But this is just a thought, and more on this later, in chapter 8 below, when doing free probability.

\section*{6d. Wishart matrices}

Let us discuss now the Wishart matrices, which are the positive analogues of the Wigner matrices. Quite surprisingly, the computation here leads to the Catalan numbers, but not in the same way as for the Wigner matrices, the result being as follows:

\index{Wishart matrix}
\index{Catalan numbers}

\begin{theorem}
Given a sequence of complex Wishart matrices
$$W_N=Y_NY_N^*\in M_N(L^\infty(X))$$
with $Y_N$ being $N\times N$ complex Gaussian of parameter $t>0$, we have
$$M_k\left(\frac{W_N}{N}\right)\simeq t^kC_k$$
for any exponent $k\in\mathbb N$, in the $N\to\infty$ limit.
\end{theorem}

\begin{proof}
There are several possible proofs for this result, as follows:

\medskip

(1) A first method is by using the formula that we have in Theorem 6.20, for the Gaussian matrices $Y_N$. Indeed, we know from there that we have the following formula, valid for any colored integer $K=\circ\bullet\bullet\circ\ldots\,$, in the $N\to\infty$ limit:
$$M_K\left(\frac{Y_N}{\sqrt{N}}\right)\simeq t^{|K|/2}|\mathcal{NC}_2(K)|$$

With $K=\circ\bullet\circ\bullet\ldots\,$, alternating word of length $2k$, with $k\in\mathbb N$, this gives:
$$M_k\left(\frac{Y_NY_N^*}{N}\right)\simeq t^k|\mathcal{NC}_2(K)|$$

Thus, in terms of the Wishart matrix $W_N=Y_NY_N^*$ we have, for any $k\in\mathbb N$:
$$M_k\left(\frac{W_N}{N}\right)\simeq t^k|\mathcal{NC}_2(K)|$$

The point now is that, by doing some combinatorics, we have:
$$|\mathcal{NC}_2(K)|=|NC_2(2k)|=C_k$$

Thus, we are led to the formula in the statement.

\medskip

(2) A second method, that we will explain now as well, is by proving the result directly, starting from definitions. The matrix entries of our matrix $W=W_N$ are given by:
$$W_{ij}=\sum_{r=1}^NY_{ir}\bar{Y}_{jr}$$

Thus, the normalized traces of powers of $W$ are given by the following formula:
\begin{eqnarray*}
tr(W^k)
&=&\frac{1}{N}\sum_{i_1=1}^N\ldots\sum_{i_k=1}^NW_{i_1i_2}W_{i_2i_3}\ldots W_{i_ki_1}\\
&=&\frac{1}{N}\sum_{i_1=1}^N\ldots\sum_{i_k=1}^N\sum_{r_1=1}^N\ldots\sum_{r_k=1}^NY_{i_1r_1}\bar{Y}_{i_2r_1}Y_{i_2r_2}\bar{Y}_{i_3r_2}\ldots Y_{i_kr_k}\bar{Y}_{i_1r_k}
\end{eqnarray*}

By rescaling now $W$ by a $1/N$ factor, as in the statement, we obtain:
$$tr\left(\left(\frac{W}{N}\right)^k\right)=\frac{1}{N^{k+1}}\sum_{i_1=1}^N\ldots\sum_{i_k=1}^N\sum_{r_1=1}^N\ldots\sum_{r_k=1}^NY_{i_1r_1}\bar{Y}_{i_2r_1}Y_{i_2r_2}\bar{Y}_{i_3r_2}\ldots Y_{i_kr_k}\bar{Y}_{i_1r_k}$$

By using now the Wick rule, we obtain the following formula for the moments, with $K=\circ\bullet\circ\bullet\ldots\,$, alternating word of lenght $2k$, and with $I=(i_1r_1,i_2r_1,\ldots,i_kr_k,i_1r_k)$:
\begin{eqnarray*}
M_k\left(\frac{W}{N}\right)
&=&\frac{t^k}{N^{k+1}}\sum_{i_1=1}^N\ldots\sum_{i_k=1}^N\sum_{r_1=1}^N\ldots\sum_{r_k=1}^N\#\left\{\pi\in\mathcal P_2(K)\Big|\pi\leq\ker(I)\right\}\\
&=&\frac{t^k}{N^{k+1}}\sum_{\pi\in\mathcal P_2(K)}\#\left\{i,r\in\{1,\ldots,N\}^k\Big|\pi\leq\ker(I)\right\}
\end{eqnarray*}

In order to compute this quantity, we use the standard bijection $\mathcal P_2(K)\simeq S_k$. By identifying the pairings $\pi\in\mathcal P_2(K)$ with their counterparts $\pi\in S_k$, we obtain:
\begin{eqnarray*}
M_k\left(\frac{W}{N}\right)
&=&\frac{t^k}{N^{k+1}}\sum_{\pi\in S_k}\#\left\{i,r\in\{1,\ldots,N\}^k\Big|i_s=i_{\pi(s)+1},r_s=r_{\pi(s)},\forall s\right\}
\end{eqnarray*}

Now let $\gamma\in S_k$ be the full cycle, which is by definition the following permutation:
$$\gamma=(1 \, 2 \, \ldots \, k)$$

The general factor in the product computed above is then 1 precisely when following two conditions are simultaneously satisfied:
$$\gamma\pi\leq\ker i\quad,\quad 
\pi\leq\ker r$$

Counting the number of free parameters in our moment formula, we obtain:
$$M_k\left(\frac{W}{N}\right)
=t^k\sum_{\pi\in S_k}N^{|\pi|+|\gamma\pi|-k-1}$$

The point now is that the last exponent is well-known to be $\leq 0$, with equality precisely when the permutation $\pi\in S_k$ is geodesic, which in practice means that $\pi$ must come from a noncrossing partition. Thus we obtain, in the $N\to\infty$ limit:
$$M_k\left(\frac{W}{N}\right)\simeq t^kC_k$$

Thus, we are led to the conclusion in the statement.
\end{proof}

As a consequence of the above result, we have a new look on the Catalan numbers, which is more adapted to our present Wishart matrix considerations, as follows:

\begin{proposition}
The Catalan numbers $C_k=|NC_2(2k)|$ appear as well as
$$C_k=|NC(k)|$$
where $NC(k)$ is the set of all noncrossing partitions of $\{1,\ldots,k\}$.
\end{proposition}

\begin{proof}
This follows indeed from the proof of Theorem 6.26. Observe that we obtain as well a formula in terms of matching pairings of alternating colored integers.
\end{proof}

The direct explanation for the above formula, relating noncrossing partitions and pairings, comes form the following result, which is very useful, and good to know:

\index{fattening of partitions}
\index{shrinking partitions}
\index{noncrossing partitions}
\index{noncrossing pairings}

\begin{proposition}
We have a bijection between noncrossing partitions and pairings
$$NC(k)\simeq NC_2(2k)$$
which is constructed as follows:
\begin{enumerate}
\item The application $NC(k)\to NC_2(2k)$ is the ``fattening'' one, obtained by doubling all the legs, and doubling all the strings as well.

\item Its inverse $NC_2(2k)\to NC(k)$ is the ``shrinking'' application, obtained by collapsing pairs of consecutive neighbors.
\end{enumerate}
\end{proposition}

\begin{proof}
The fact that the two operations in the statement are indeed inverse to each other is clear, by computing the corresponding two compositions, with the remark that the construction of the fattening operation requires the partitions to be noncrossing.
\end{proof}

Getting back now to probability, we are led to the question of finding the law having the Catalan numbers as moments, in the above way. The result here is as follows:

\begin{proposition}
The real measure having the Catalan numbers as moments is
$$\pi_1=\frac{1}{2\pi}\sqrt{4x^{-1}-1}\,dx$$
called Marchenko-Pastur law of parameter $1$.
\end{proposition}

\begin{proof}
The moments of the law $\pi_1$ in the statement can be computed with the change of variable $x=4\cos^2t$, as follows:
\begin{eqnarray*}
M_k
&=&\frac{1}{2\pi}\int_0^4\sqrt{4x^{-1}-1}\,x^kdx\\
&=&\frac{1}{2\pi}\int_0^{\pi/2}\frac{\sin t}{\cos t}\cdot(4\cos^2t)^k\cdot 2\cos t\sin t\,dt\\
&=&\frac{4^{k+1}}{\pi}\int_0^{\pi/2}\cos^{2k}t\sin^2t\,dt\\
&=&\frac{4^{k+1}}{\pi}\cdot\frac{\pi}{2}\cdot\frac{(2k)!!2!!}{(2k+3)!!}\\
&=&2\cdot 4^k\cdot\frac{(2k)!/2^kk!}{2^{k+1}(k+1)!}\\
&=&C_k
\end{eqnarray*}

Thus, we are led to the conclusion in the statement.
\end{proof}

Now back to the Wishart matrices, we are led to the following result:

\index{Wishart matrix}

\begin{theorem}
Given a sequence of complex Wishart matrices
$$W_N=Y_NY_N^*\in M_N(L^\infty(X))$$
with $Y_N$ being $N\times N$ complex Gaussian of parameter $t>0$, we have
$$\frac{W_N}{tN}\sim\frac{1}{2\pi}\sqrt{4x^{-1}-1}\,dx$$
with $N\to\infty$, with the limiting measure being the Marchenko-Pastur law $\pi_1$.
\end{theorem}

\begin{proof}
This follows indeed from Theorem 6.26 and Proposition 6.29.
\end{proof}

As a comment now, while the above result is definitely something interesting at $t=1$, at general $t>0$ this looks more like a ``fake'' generalization of the $t=1$ result, because the law $\pi_1$ stays the same, modulo a trivial rescaling. The reasons behind this phenomenon are quite subtle, and skipping some discussion, the point is that Theorem 6.30 is indeed something ``fake'' at general $t>0$, and the correct generalization of the $t=1$ computation, involving more general classes of complex Wishart matrices, is as follows:

\index{Wishart matrix}
\index{Marchenko-Pastur law}

\begin{theorem}
Given a sequence of general complex Wishart matrices
$$W_N=Y_NY_N^*\in M_N(L^\infty(X))$$
with $Y_N$ being $N\times M$ complex Gaussian of parameter $1$, we have
$$\frac{W_N}{N}\sim\max(1-t,0)\delta_0+\frac{\sqrt{4t-(x-1-t)^2}}{2\pi x}\,dx$$
with $M=tN\to\infty$, with the limiting measure being the Marchenko-Pastur law $\pi_t$.
\end{theorem}

\begin{proof}
This follows once again by using the moment method, the limiting moments in the $M=tN\to\infty$ regime being as follows, after doing the combinatorics:
$$M_k\left(\frac{W_N}{N}\right)\simeq\sum_{\pi\in NC(k)}t^{|\pi|}$$

But these numbers are the moments of the Marchenko-Pastur law $\pi_t$, which in addition has the density given by the formula in the statement, and this gives the result.
\end{proof}

As a philosophical conclusion now, we have 4 main laws in what we have been doing so far, namely the Gaussian laws $g_t$, the Poisson laws $p_t$, the Wigner laws $\gamma_t$ and the Marchenko-Pastur laws $\pi_t$. These laws naturally form a diagram, as follows:
$$
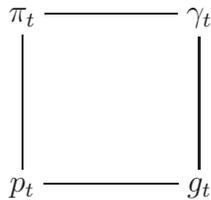
\xymatrix@R=50pt@C=50pt{
\pi_t\ar@{-}[r]\ar@{-}[d]&\gamma_t\ar@{-}[d]\\
p_t\ar@{-}[r]&g_t}$$

We will see in chapter 8 that $\pi_t,\gamma_t$ appear as ``free analogues'' of $p_t,g_t$, and that a full theory can be developed, with central limiting theorems for all 4 laws, convolution semigroup results for all 4 laws too, and Lie group type results for all 4 laws too. And also, we will be back to the random matrices as well, with further results about them.

\section*{6e. Exercises} 

There has been a lot of non-trivial combinatorics and calculus in this chapter, sometimes only briefly explained, and as an exercise on all this, we have:

\begin{exercise}
Clarify all the details in connection with the Wigner and Marchenko-Pastur computations, first at $t=1$, and then for general $t>0$.
\end{exercise}

As before, these are things discussed in the above, but only briefly, this whole chapter having been just a modest introduction to this exciting subject which are the random matrices. In the hope that you will find some time, and do the exercise.

\chapter{Quantum spaces}

\section*{7a. Gelfand theorem}

We have seen that the von Neumann algebras $A\subset B(H)$ are interesting objects, and it is tempting to go ahead with a systematic study of such algebras. This is what Murray and von Neumann did, when first coming across such algebras, back in the 1930s, in their series of papers \cite{mv1}, \cite{mv2}, \cite{mv3}, \cite{vn1}, \cite{vn2}, \cite{vn3}. In what concerns us, we will rather keep this material for later, and talk instead, in this chapter and in the next one, of things which are perhaps more basic, motivated by the following definition:

\index{quantum space}
\index{quantum measured space}
\index{quantum probability space}

\begin{definition}
Given a von Neumann algebra $A\subset B(H)$, coming with a faithful positive unital trace $tr:A\to\mathbb C$, we write
$$A=L^\infty(X)$$
and call $X$ a quantum probability space. We also write the trace as $tr=\int_X$, and call it integration with respect to the uniform measure on $X$.
\end{definition}

Obviously, this is something exciting, and we have seen how some interesting theory can be developed along these lines in the simplest case, that of the random matrix algebras. Thus, all this needs a better understanding, before going ahead with the above-mentioned Murray-von Neumann theory. In order to get started, here are a few comments:

\bigskip

(1) Generally speaking, all this comes from the fact that the commutative von Neumann algebras are those of the form $A=L^\infty(X)$, with $X$ being a measured space. Since in the finite measure case, $\mu(X)<\infty$, the integration can be regarded as being a faithful positive unital trace $tr:L^\infty(X)\to\mathbb C$, we are basically led to Definition 7.1.

\bigskip

(2) Regarding our assumption $\mu(X)<\infty$, making the integration $tr:A\to\mathbb C$ bounded, this is something advanced, coming from deep classification results of von Neumann and Connes, which roughly state that ``modulo classical measure theory, the study of the quantum measured spaces $X$ basically reduces to the case $\mu(X)<\infty$''.

\bigskip

(3) Finally, the traciality of $tr:A\to\mathbb C$ is something advanced too, again coming from that classification results of von Neumann and Connes, which in their more precise formulation state that ``modulo classical measure theory, the study of the quantum measured spaces $X$ basically reduces to the case where $\mu(X)<\infty$, and $\int_X$ is a trace''.

\bigskip

In short, complicated all this, and you will have to trust me here. Moving ahead now, there is one more thing to be discussed in connection with Definition 7.1, and this is physics. Let me formulate here the question that you surely have in mind:

\begin{question}
As physicists we already agreed, without clear evidence, that our operators $T:H\to H$ should be bounded. But what about quantum spaces, is it a good idea to assume that these are as above, of finite mass, and with tracial integration? 
\end{question}

Well, this is certainly an interesting question. In favor of my choice, I would argue that the mathematical physics of Jones \cite{jo1}, \cite{jo2}, \cite{jo3}, \cite{jo5}, \cite{jo6} and Voiculescu \cite{vo1}, \cite{vo2}, \cite{vdn} needs a trace $tr:A\to\mathbb C$, as above. And the same goes for certain theoretical physics continuations of the main work of Connes \cite{co3}, as for instance the basic theory of the Standard Model spectral triple of Chamseddine-Connes, whose free gauge group has tracial Haar integration. Needless to say, all this is quite subjective. But hey, question of theoretical physics you asked, answer of theoretical physics is what you get.

\bigskip

Hang on, we are not done yet. Now that we are convinced that Definition 7.1 is the correct one, be that on mathematical or physical grounds, let us look for examples. And here the situation is quite grim, because even in the classical case, we have:

\begin{fact}
The measure on a classical measured space $X$ cannot come out of nowhere, and is usually a Haar measure, appearing by theorem. Thus, in our picture
$$A\subset B(H)$$
both the Hilbert space $H=L^2(X)$ and the von Neumann algebra $A=L^\infty(X)$ should appear by theorem, not by definition, contrary to what Definition 7.1 says.
\end{fact}

To be more precise, in what regards the first assertion, this is certainly the case with simple objects like Lie groups, or spheres and other homogeneous spaces. Of course you might say that $[0,1]$ with the uniform measure is a measured space, but isn't $[0,1]$ obtained by cutting the Lie group $\mathbb R$, with its Haar measure. And the same goes with $[0,1]$ with an arbitrary measure $f(x)dx$, or with $[0,1]$ being deformed into a curve, and so on, because that $dx$, or what is left from it, will always refer to the Haar measure of $\mathbb R$.

\bigskip

As for the second assertion, nothing much to comment here, mathematics has spoken. So, getting back now to Definition 7.1 as it is, looks like we have two dead bodies there, the Hilbert space $H$ and the operator algebra $A$. So let us try to get rid of at least one of them. But which? In the lack of any obvious idea, let us turn to physics:

\begin{question}
In quantum mechanics, which came first, the Hilbert space $H$, or the operator algebra $A$?
\end{question}

Unfortunately this question is as difficult as the one regarding the chicken and the egg. A look at what various physicists said on this matter, in a direct or indirect way, does not help much, and by the end of the day we are left with guidelines like ``no one understands quantum mechanics'' (Feynman), ``shut up and compute'' (Dirac) and so on. And all this, coming on top on what has been already said on Definition 7.1, of rather unclear nature, is probably too much. That is, the last drop, time to conclude:

\begin{conclusion}
The theory of von Neumann algebras has the same peculiarity as quantum mechanics: it tends to self-destruct, when approached axiomatically.
\end{conclusion}

And we will take this as good news, providing us with warm evidence that the theory of von Neumann algebras is indeed related to quantum mechanics. This is what matters, being on the right track, and difficulties and all the rest, we won't be scared by them.

\bigskip

Back to business now, in practice, we must go back to chapter 5, and examine what we were saying right before introducing the von Neumann algebras. And at that time, we were talking about general operator algebras $A\subset B(H)$, closed with respect to the norm, but not necessarily with respect to the weak topology. But this suggests formulating the following definition, somewhat as a purely mathematical answer to Question 7.4:

\index{Banach algebra}
\index{operator algebra}

\begin{definition}
A $C^*$-algebra is an complex algebra $A$, given with:
\begin{enumerate}
\item A norm $a\to||a||$, making it into a Banach algebra.

\item An involution $a\to a^*$, related to the norm by the formula $||aa^*||=||a||^2$. 
\end{enumerate}
\end{definition}

Here by Banach algebra we mean a complex algebra with a norm satisfying all the conditions for a vector space norm, along with $||ab||\leq||a||\cdot||b||$ and $||1||=1$, and which is such that our algebra is complete, in the sense that the Cauchy sequences converge. As for the involution, this must be antilinear, antimultiplicative, and satisfying $a^{**}=a$.

\bigskip

As basic examples, we have the operator algebra $B(H)$, for any Hilbert space $H$, and more generally, the norm closed $*$-subalgebras $A\subset B(H)$. It is possible to prove that any $C^*$-algebra appears in this way, but this is a non-trivial result, called GNS theorem, and more on this later. Note in passing that this result tells us that there is no need to memorize the above axioms for the $C^*$-algebras, because these are simply the obvious things that can be said about $B(H)$, and its norm closed $*$-subalgebras $A\subset B(H)$.

\bigskip

As a second class of basic examples, which are of great interest for us, we have:

\index{compact space}

\begin{proposition}
If $X$ is a compact space,  the algebra $C(X)$ of continuous functions $f:X\to\mathbb C$ is a $C^*$-algebra, with the usual norm and involution, namely:
$$||f||=\sup_{x\in X}|f(x)|\quad,\quad 
f^*(x)=\overline{f(x)}$$
This algebra is commutative, in the sense that $fg=gf$, for any $f,g\in C(X)$.
\end{proposition}

\begin{proof}
All this is clear from definitions. Observe that we have indeed:
$$||ff^*||
=\sup_{x\in X}|f(x)|^2
=||f||^2$$

Thus, the axioms are satisfied, and finally $fg=gf$ is clear.
\end{proof}

In general, the $C^*$-algebras can be thought of as being algebras of operators, over some Hilbert space which is not present. By using this philosophy, one can emulate spectral theory in this setting, with extensions of the various results from chapters 3,5:

\index{spectrum}
\index{polynomial calculus}
\index{rational calculus}
\index{holomorphic calculus}
\index{unitary}
\index{self-adjoint element}
\index{spectral radius}
\index{continuous calculus}

\begin{theorem}
Given element $a\in A$ of a $C^*$-algebra, define its spectrum as:
$$\sigma(a)=\left\{\lambda\in\mathbb C\Big|a-\lambda\notin A^{-1}\right\}$$
The following spectral theory results hold, exactly as in the $A=B(H)$ case:
\begin{enumerate}
\item We have $\sigma(ab)\cup\{0\}=\sigma(ba)\cup\{0\}$.

\item We have polynomial, rational and holomorphic calculus.

\item As a consequence, the spectra are compact and non-empty.

\item The spectra of unitaries $(u^*=u^{-1})$ and self-adjoints $(a=a^*)$ are on $\mathbb T,\mathbb R$.

\item The spectral radius of normal elements $(aa^*=a^*a)$ is given by $\rho(a)=||a||$.
\end{enumerate}
In addition, assuming $a\in A\subset B$, the spectra of $a$ with respect to $A$ and to $B$ coincide.
\end{theorem}

\begin{proof}
This is something that we know from chapter 3, in the case $A=B(H)$, and then from chapter 5, in the case $A\subset B(H)$. In general, the proof is similar:

\medskip

(1) Regarding the assertions (1-5), which are of course formulated a bit informally, the proofs here are perfectly similar to those for the full operator algebra $A=B(H)$. All this is standard material, and in fact, things in chapters 3 were written in such a way as for their extension now, to the general $C^*$-algebra setting, to be obvious.

\medskip

(2) Regarding the last assertion, we know this from chapter 5 for $A\subset B\subset B(H)$, and the proof in general is similar. Indeed, the inclusion $\sigma_B(a)\subset\sigma_A(a)$ is clear. For the converse, assume $a-\lambda\in B^{-1}$, and consider the following self-adjoint element:
$$b=(a-\lambda )^*(a-\lambda )$$

The difference between the two spectra of $b\in A\subset B$ is then given by:
$$\sigma_A(b)-\sigma_B(b)=\left\{\mu\in\mathbb C-\sigma_B(b)\Big|(b-\mu)^{-1}\in B-A\right\}$$

Thus this difference in an open subset of $\mathbb C$. On the other hand $b$ being self-adjoint, its two spectra are both real, and so is their difference. Thus the two spectra of $b$ are equal, and in particular $b$ is invertible in $A$, and so $a-\lambda\in A^{-1}$, as desired.
\end{proof}

We can now get back to the commutative $C^*$-algebras, and we have the following result, due to Gelfand, which will be of crucial importance for us:

\index{commutative algebra}
\index{character}
\index{Banach algebra}
\index{spectrum of algebra}
\index{Gelfand theorem}

\begin{theorem}
The commutative $C^*$-algebras are exactly the algebras of the form 
$$A=C(X)$$
with the ``spectrum'' $X$ of such an algebra being the space of characters $\chi:A\to\mathbb C$, with topology making continuous the evaluation maps $ev_a:\chi\to\chi(a)$.
\end{theorem}

\begin{proof}
This is something that we basically know from chapter 5, but always good to talk about it again. Given a commutative $C^*$-algebra $A$, we can define $X$ as in the statement. Then $X$ is compact, and $a\to ev_a$ is a morphism of algebras, as follows:
$$ev:A\to C(X)$$

(1) We first prove that $ev$ is involutive. We use the following formula, which is similar to the $z=Re(z)+iIm(z)$ formula for the usual complex numbers:
$$a=\frac{a+a^*}{2}+i\cdot\frac{a-a^*}{2i}$$

Thus it is enough to prove the equality $ev_{a^*}=ev_a^*$ for self-adjoint elements $a$. But this is the same as proving that $a=a^*$ implies that $ev_a$ is a real function, which is in turn true, because $ev_a(\chi)=\chi(a)$ is an element of $\sigma(a)$, contained in $\mathbb R$.

\medskip

(2) Since $A$ is commutative, each element is normal, so $ev$ is isometric:
$$||ev_a||
=\rho(a)
=||a||$$

(3) It remains to prove that $ev$ is surjective. But this follows from the Stone-Weierstrass theorem, because $ev(A)$ is a closed subalgebra of $C(X)$, which separates the points.
\end{proof}

In view of the Gelfand theorem, we can formulate the following key definition:

\index{compact quantum space}
\index{quantum space}

\begin{definition}
Given an arbitrary $C^*$-algebra $A$, we write 
$$A=C(X)$$
and call $X$ a compact quantum space.
\end{definition}

This might look like something informal, but it is not. Indeed, we can define the category of compact quantum spaces to be the category of the $C^*$-algebras, with the arrows reversed. When $A$ is commutative, the above space $X$ exists indeed, as a Gelfand spectrum, $X=Spec(A)$. In general, $X$ is something rather abstract, and our philosophy here will be that of studying of course $A$, but formulating our results in terms of $X$. For instance whenever we have a morphism $\Phi:A\to B$, we will write $A=C(X),B=C(Y)$, and rather speak of the corresponding morphism $\phi:Y\to X$. And so on.

\bigskip

Technically speaking, we will see later that the above formalism has its limitations, and needs a fix. To be more precise, when looking at compact quantum spaces having a probability measure, there are more of them in the sense of Definition 7.10, than in the von Neumann algebra sense. Thus, all this needs a fix. But more on this later.

\bigskip

As a first concrete consequence of the Gelfand theorem, we have:

\index{normal element}
\index{continuous functional calculus}

\begin{proposition}
Assume that $a\in A$ is normal, and let $f\in C(\sigma(a))$.
\begin{enumerate}
\item We can define $f(a)\in A$, with $f\to f(a)$ being a morphism of $C^*$-algebras.

\item We have the ``continuous functional calculus'' formula $\sigma(f(a))=f(\sigma(a))$.
\end{enumerate}
\end{proposition}

\begin{proof}
Since $a$ is normal, the $C^*$-algebra $<a>$ that is generates is commutative, so if we denote by $X$ the space formed by the characters $\chi:<a>\to\mathbb C$, we have:
$$<a>=C(X)$$

Now since the map $X\to\sigma(a)$ given by evaluation at $a$ is bijective, we obtain:
$$<a>=C(\sigma(a))$$

Thus, we are dealing with usual functions, and this gives all the assertions.
\end{proof}

As another consequence of the Gelfand theorem, we have:

\index{positive element}

\begin{proposition}
For a normal element $a\in A$, the following are equivalent:
\begin{enumerate}
\item $a$ is positive, in the sense that $\sigma(a)\subset[0,\infty)$.

\item $a=b^2$, for some $b\in A$ satisfying $b=b^*$.

\item $a=cc^*$, for some $c\in A$.
\end{enumerate}
\end{proposition}

\begin{proof}
This is very standard, exactly as in $A=B(H)$ case, as follows:

\medskip

$(1)\implies(2)$ Since $f(z)=\sqrt{z}$ is well-defined on $\sigma(a)\subset[0,\infty)$, we can set $b=\sqrt{a}$. 

\medskip

$(2)\implies(3)$ This is trivial, because we can set $c=b$. 

\medskip

$(3)\implies(1)$ We proceed by contradiction. By multiplying $c$ by a suitable element of $<cc^*>$, we are led to the existence of an element $d\neq0$ satisfying $-dd^*\geq0$. By writing now $d=x+iy$ with $x=x^*,y=y^*$ we have:
$$dd^*+d^*d=2(x^2+y^2)\geq0$$

Thus $d^*d\geq0$, contradicting the fact that $\sigma(dd^*),\sigma(d^*d)$ must coincide outside $\{0\}$.
\end{proof}

Let us clarify now the relation between $C^*$-algebras and von Neumann algebras. In order to do so, we need a prove a key result, called GNS representation theorem, stating that any $C^*$-algebra appears as an operator algebra. As a first result, we have:

\begin{proposition}
Let $A$ be a commutative $C^*$-algebra, write $A=C(X)$, with $X$ being a compact space, and let $\mu$ be a positive measure on $X$. We have then
$$A\subset B(H)$$
where $H=L^2(X)$, with $f\in A$ corresponding to the operator $g\to fg$.
\end{proposition}

\begin{proof}
Given a continuous function $f\in C(X)$, consider the operator $T_f(g)=fg$, on $H=L^2(X)$. Observe that $T_f$ is indeed well-defined, and bounded as well, because:
$$||fg||_2
=\sqrt{\int_X|f(x)|^2|g(x)|^2d\mu(x)}
\leq||f||_\infty||g||_2$$

The application $f\to T_f$ being linear, involutive, continuous, and injective as well, we obtain in this way a $C^*$-algebra embedding $A\subset B(H)$, as claimed.
\end{proof}

In order to prove the GNS representation theorem, we must extend the above construction, to the case where $A$ is not necessarily commutative. Let us start with:

\index{positive element}
\index{positive linear form}
\index{faithful form}

\begin{definition}
Consider a $C^*$-algebra $A$.
\begin{enumerate}
\item $\varphi:A\to\mathbb C$ is called positive when $a\geq0\implies\varphi(a)\geq0$.

\item $\varphi:A\to\mathbb C$ is called faithful and positive when $a\geq0,a\neq0\implies\varphi(a)>0$.
\end{enumerate}
\end{definition}

In the commutative case, $A=C(X)$, the positive elements are the positive functions, $f:X\to[0,\infty)$. As for the positive linear forms $\varphi:A\to\mathbb C$, these appear as follows, with $\mu$ being positive, and strictly positive if we want $\varphi$ to be faithful and positive:
$$\varphi(f)=\int_Xf(x)d\mu(x)$$

In general, the positive linear forms can be thought of as being integration functionals with respect to some underlying ``positive measures''. We can use them as follows:

\begin{proposition}
Let $\varphi:A\to\mathbb C$ be a positive linear form.
\begin{enumerate}
\item $<a,b>=\varphi(ab^*)$ defines a generalized scalar product on $A$.

\item By separating and completing we obtain a Hilbert space $H$.

\item $\pi(a):b\to ab$ defines a representation $\pi:A\to B(H)$.

\item If $\varphi$ is faithful in the above sense, then $\pi$ is faithful.
\end{enumerate}
\end{proposition}

\begin{proof}
Almost everything here is straightforward, as follows:

\medskip

(1) This is clear from definitions, and from the basic properties of the positive elements $a\geq0$, which can be established exactly as in the $A=B(H)$ case.

\medskip

(2) This is a standard procedure, which works for any scalar product, the idea being that of dividing by the vectors satisfying $<x,x>=0$, then completing.

\medskip

(3) All the verifications here are standard algebraic computations, in analogy with what we have seen many times, for multiplication operators, or group algebras.

\medskip

(4) Assuming that we have $a\neq0$, we have then $\pi(aa^*)\neq0$, which in turn implies by faithfulness that we have $\pi(a)\neq0$, which gives the result.
\end{proof}

In order to establish the embedding theorem, it remains to prove that any $C^*$-algebra has a faithful positive linear form $\varphi:A\to\mathbb C$. This is something more technical:

\begin{proposition}
Let $A$ be a $C^*$-algebra.
\begin{enumerate}
\item Any positive linear form $\varphi:A\to\mathbb C$ is continuous.

\item A linear form $\varphi$ is positive iff there is a norm one $h\in A_+$ such that $||\varphi||=\varphi(h)$.

\item For any $a\in A$ there exists a positive norm one form $\varphi$ such that $\varphi(aa^*)=||a||^2$.

\item If $A$ is separable there is a faithful positive form $\varphi:A\to\mathbb C$.
\end{enumerate}
\end{proposition}

\begin{proof}
The proof here is quite technical, inspired from the existence proof of the probability measures on abstract compact spaces, the idea being as follows:

\medskip

(1) This follows from Proposition 7.15, via the following estimate:
$$|\varphi(a)|
\leq||\pi(a)||\varphi(1)
\leq||a||\varphi(1)$$

(2) In one sense we can take $h=1$. Conversely, let $a\in A_+$, $||a||\leq1$. We have:
$$|\varphi(h)-\varphi(a)|
\leq||\varphi||\cdot||h-a||
\leq\varphi(h)$$

Thus we have $Re(\varphi(a))\geq0$, and with $a=1-h$ we obtain:
$$Re(\varphi(1-h))\geq0$$

Thus $Re(\varphi(1))\geq||\varphi||$, and so $\varphi(1)=||\varphi||$, so we can assume $h=1$. Now observe that for any self-adjoint element $a$, and any $t\in\mathbb R$ we have, with $\varphi(a)=x+iy$:
\begin{eqnarray*}
\varphi(1)^2(1+t^2||a||^2)
&\geq&\varphi(1)^2||1+t^2a^2||\\
&=&||\varphi||^2\cdot||1+ita||^2\\
&\geq&|\varphi(1+ita)|^2\\
&=&|\varphi(1)-ty+itx|\\
&\geq&(\varphi(1)-ty)^2
\end{eqnarray*}

Thus we have $y=0$, and this finishes the proof of our remaining claim.

\medskip

(3) We can set $\varphi(\lambda aa^*)=\lambda||a||^2$ on the linear space spanned by $aa^*$, then extend this functional by Hahn-Banach, to the whole $A$. The positivity follows from (2).

\medskip

(4) This is standard, by starting with a dense sequence $(a_n)$, and taking the Ces\`aro limit of the functionals constructed in (3). We have $\varphi(aa^*)>0$, and we are done.
\end{proof}

With these ingredients in hand, we can now state and prove:

\index{GNS theorem}
\index{Gelfand-Naimark-Segal}

\begin{theorem}
Any $C^*$-algebra appears as a norm closed $*$-algebra of operators
$$A\subset B(H)$$
over a certain Hilbert space $H$. When $A$ is separable, $H$ can be taken to be separable.
\end{theorem}

\begin{proof}
This result, called called GNS representation theorem after Gelfand, Naimark and Segal, follows indeed by combining Proposition 7.15 with Proposition 7.16.
\end{proof}

All this might seem quite surprising, and your first reaction would be to say what have we been we doing here, with our $C^*$-algebra theory, because we are now back to operator algebras $A\subset B(H)$, and everything that we did with $C^*$-algebras, extending things that we knew about operator algebras $A\subset B(H)$, looks more like a waste of time.

\bigskip

Error. The axioms in Definition 7.6, coupled with the writing $A=C(X)$ in Definition 7.10, are something powerful, because they do not involve any kind of $L^2$ or $L^\infty$ functions on our quantum spaces $X$. Thus, we can start hunting for such spaces, just by defining $C^*$-algebras with generators and relations, then look for Haar measures on such spaces, and use the GNS construction in order to reach to von Neumann algebras. Before getting into this, however, let us summarize the above discussion as follows:

\index{compact quantum space}
\index{compact quantum measured space}
\index{quantum measured space}

\begin{theorem}
We can talk about compact quantum measured spaces, as follows:
\begin{enumerate}
\item The category of compact quantum measured spaces $(X,\mu)$ is the category of the $C^*$-algebras with faithful traces $(A,\varphi)$, with the arrows reversed. 

\item In the case where we have a non-faithful trace $\varphi$, we can still talk about the corresponding space $(X,\mu)$, by performing the GNS construction.

\item By taking the weak closure in the GNS representation, we obtain the von Neumann algebra $A''=L^\infty(X)$, in the previous general measured space sense.
\end{enumerate}
\end{theorem}

\begin{proof}
All this follows from Theorem 7.17, and from the other things that we already know, with the whole result itself being something rather philosophical.
\end{proof}

\section*{7b. Tori, amenability}

In the remainder of this chapter we explore the whole new world opened by the $C^*$-algebra theory, with the study of several key examples. We will first discuss the group duals, also called noncommutative tori. Let us start with a well-known result:

\index{compact abelian group}
\index{discrete abelian group}
\index{finite abelian group}
\index{abelian group}
\index{dual group}
\index{Pontrjagin duality}

\begin{theorem}
The compact abelian groups $G$ are in correspondence with the discrete abelian groups $\Gamma$, via Pontrjagin duality, 
$$G=\widehat{\Gamma}\quad,\quad 
\Gamma=\widehat{G}$$
with the dual of a locally compact group $L$ being the locally compact group $\widehat{L}$ consisting of the continuous group characters $\chi:L\to\mathbb T$.
\end{theorem}

\begin{proof}
This is something very standard, the idea being that, given a group $L$ as above, its continuous characters $\chi:L\to\mathbb T$ form indeed a group, that we can call $\widehat{L}$. The correspondence $L\to\widehat{L}$ constructed in this way has then the following properties:

\medskip

(1) We have $\widehat{\mathbb Z}_N=\mathbb Z_N$. This is the basic computation to be performed, before anything else, and which is something algebraic, with roots of unity.

\medskip

(2) More generally, the dual of a finite abelian group $G=\mathbb Z_{N_1}\times\ldots\times\mathbb Z_{N_k}$ is the group $G$ itself. This comes indeed from (1) and from $\widehat{G\times H}=\widehat{G}\times\widehat{H}$.

\medskip

(3) At the opposite end now, that of the locally compact groups which are not compact, nor discrete, the main example, which is standard, is $\widehat{\mathbb R}=\mathbb R$.

\medskip

(4) Getting now to what we are interested in, it follows from the definition of the correspondence $L\to\widehat{L}$ that when $L$ is compact $\widehat{L}$ is discrete, and vice versa.

\medskip

(5) Finally, in order to best understand this latter phenomenon, the best is to work out the main pair of examples, which are $\widehat{\mathbb T}=\mathbb Z$ and $\widehat{\mathbb Z}=\mathbb T$.
\end{proof}

Our claim now is that, by using operator algebra theory, we can talk about the dual $G=\widehat{\Gamma}$ of any discrete group $\Gamma$. Let us start our discussion in the von Neumann algebra setting, where things are particularly simple. We have here:

\index{group algebra}
\index{group von Neumann algebra}
\index{left regular representation}
\index{Fourier transform}

\begin{theorem}
Given a discrete group $\Gamma$, we can construct its von Neumann algebra
$$L(\Gamma)\subset B(l^2(\Gamma))$$
by using the left regular representation. This algebra has a faithful positive trace, $tr(g)=\delta_{g,1}$, and when $\Gamma$ is abelian we have an isomorphism of tracial von Neumann algebras
$$L(\Gamma)\simeq L^\infty(G)$$
given by a Fourier type transform, where $G=\widehat{\Gamma}$ is the compact dual of $\Gamma$.
\end{theorem}

\begin{proof}
There are many assertions here, the idea being as follows:

\medskip

(1) The first part is standard, with the left regular representation of $\Gamma$ working as expected, and being a unitary representation, as follows:
$$\Gamma\subset B(l^2(\Gamma))\quad,\quad
\pi(g):h\to gh$$

(2) The positivity of the trace comes from the following alternative formula for it, with the equivalence with the definition in the statement being clear:
$$tr(T)=<T1,1>$$

(3) The third part is standard as well, because when $\Gamma$ is abelian the algebra $L(\Gamma)$ is commutative, and its spectral decomposition leads by delinearization to the group characters $\chi:\Gamma\to\mathbb T$, and so the dual group $G=\widehat{\Gamma}$, as indicated.

\medskip

(4) Finally, the fact that our isomorphism transforms the trace of $L(\Gamma)$ into the Haar integration functional of $L^\infty(G)$ is clear. Moreover, the study of various examples show that what we constructed is in fact the Fourier transform, in its various incarnations.
\end{proof}

Getting back now to our quantum space questions, we have a beginning of answer, because based on the above, we can formulate the following definition:

\begin{definition}
Given a discrete group $\Gamma$, not necessarily abelian, we can construct its abstract dual $G=\widehat{\Gamma}$ as a quantum measured space, via the following formula:
$$L^\infty(G)=L(\Gamma)$$
In the case where $\Gamma$ happens to be abelian, this quantum space $G=\widehat{\Gamma}$ is a classical space, namely the usual Pontrjagin dual of $\Gamma$, endowed with its Haar measure.
\end{definition}

Let us discuss now the same questions, in the $C^*$-algebra setting. The situation here is more complicated than in the von Neumann algebra setting, as follows:

\index{group algebra}
\index{full group algebra}
\index{reduced group algebra}
\index{left regular representation}
\index{regular representation}

\begin{proposition}
Associated to any discrete group $\Gamma$ are several group $C^*$-algebras,
$$C^*(\Gamma)\to C^*_\pi(\Gamma)\to C^*_{red}(\Gamma)$$
which are constructed as follows:
\begin{enumerate}
\item $C^*(\Gamma)$ is the closure of the group algebra $\mathbb C[\Gamma]$, with involution $g^*=g^{-1}$, with respect to the maximal $C^*$-seminorm on this $*$-algebra, which is a $C^*$-norm. 

\item $C^*_{red}(\Gamma)$ is the norm closure of the group algebra $\mathbb C[\Gamma]$ in the left regular representation, on the Hilbert space $l^2(\Gamma)$, given by $\lambda(g)(h)=gh$ and linearity.

\item $C^*_\pi(\Gamma)$ can be any intermediate $C^*$-algebra, but for best results, the indexing object $\pi$ must be a unitary group representation, satisfying $\pi\otimes\pi\subset\pi$.
\end{enumerate}
\end{proposition}

\begin{proof}
This is something quite technical, with (2) being very similar to the von Neumann algebra construction from Theorem 7.20, with (1) being something new, with the norm property there coming from (2), and finally with (3) being an informal statement, that we will comment on later, once we will know about compact quantum groups.
\end{proof}

When $\Gamma$ is finite, or abelian, or more generally amenable, all the above group algebras coincide. In the abelian case, that we are particularly interested in here, the precise result is as follows, complementing the $L^\infty$ analysis from Theorem 7.20:

\index{abelian group}
\index{discrete abelian group}
\index{compacy abelian group}
\index{Pontrjagin duality}
\index{Fourier transform}

\begin{theorem}
When $\Gamma$ is abelian all its group $C^*$-algebras coincide, and we have an isomorphism as follows, given by a Fourier type transform, 
$$C^*(\Gamma)\simeq C(G)$$
where $G=\widehat{\Gamma}$ is the compact dual of $\Gamma$. Moreover, this isomorphism transforms the standard group algebra trace $tr(g)=\delta_{g,1}$ into the Haar integration of $G$.
\end{theorem}

\begin{proof}
Since $\Gamma$ is abelian, any of its group $C^*$-algebras $A=C^*_\pi(\Gamma)$ is commutative. Thus, we can apply the Gelfand theorem, and we obtain $A=C(X)$, with $X=Spec(A)$. But the spectrum $X=Spec(A)$, consisting of the characters $\chi:A\to\mathbb C$, can be identified by delinearizing with the Pontrjagin dual $G=\widehat{\Gamma}$, and this gives the results.
\end{proof}

At a more advanced level now, we have the following result:

\index{full algebra}
\index{reduced algebra}
\index{amenable group}
\index{Kesten amenability}

\begin{theorem}
For a discrete group $\Gamma=<g_1,\ldots,g_N>$, the following conditions are equivalent, and if they are satisfied, we say that $\Gamma$ is amenable:
\begin{enumerate}
\item The projection map $C^*(\Gamma)\to C^*_{red}(\Gamma)$ is an isomorphism.

\item The morphism $\varepsilon:C^*(\Gamma)\to\mathbb C$ given by $g\to1$ factorizes through $C^*_{red}(\Gamma)$.

\item We have $N\in\sigma(Re(g_1+\ldots+g_N))$, the spectrum being taken inside $C^*_{red}(\Gamma)$.
\end{enumerate}
The amenable groups include all finite groups, and all abelian groups. As a basic example of a non-amenable group, we have the free group $F_N$, with $N\geq2$.
\end{theorem}

\begin{proof}
There are several things to be proved, the idea being as follows:

\medskip

(1) The implication $(1)\implies(2)$ is trivial, and $(2)\implies(3)$ comes from the following computation, which shows that $N-Re(g_1+\ldots+g_N)$ is not invertible inside $C^*_{red}(\Gamma)$:
\begin{eqnarray*}
\varepsilon[N-Re(g_1+\ldots+g_N)]
&=&N-Re[\varepsilon(g_1)+\ldots+\varepsilon(g_n)]\\
&=&N-N\\
&=&0
\end{eqnarray*}

As for $(3)\implies(1)$, this is something more advanced, that we will not need for the moment. We will be back to this later, directly in a more general setting.

\medskip

(2) The fact that any finite group $G$ is amenable is clear, because all the group $C^*$-algebras are equal to the usual group $*$-algebra $\mathbb C[G]$, in this case. As for the case of the abelian groups, these are all amenable as well, as shown by Theorem 7.23. 

\medskip

(3) It remains to prove that $F_N$ with $N\geq2$ is not amenable. By using $F_2\subset F_N$, it is enough to do this at $N=2$. So, consider the free group $F_2=<g,h>$. In order to prove that $F_2$ is not amenable, we use $(1)\implies(3)$. To be more precise, it is enough to show that 4 is not in the spectrum of the following operator:
$$T=\lambda(g)+\lambda(g^{-1})+\lambda(h)+\lambda(h^{-1})$$

This is a sum of four terms, each of them acting via $\delta_w\to\delta_{ew}$, with $e$ being a certain length one word. Thus if $w\neq 1$ has length $n$ then $T(\delta_w)$ is a sum of four Dirac masses, three of them at words of length $n+1$ and the remaining one at a length $n-1$ word. We can therefore decompose $T$ as a sum $T_++T_-$, where $T_+$ adds and $T_-$ cuts:
$$T=T_++T_-$$

That is, if $w\neq 1$ is a word, say beginning with $h$, then $T_\pm$ act on $\delta_w$ as follows:
$$T_+(\delta_w)=\delta_{gw}+\delta_{g^{-1}w}+\delta_{hw}\quad,\quad 
T_-(\delta_w)=\delta_{h^{-1}w}$$ 

It follows from definitions that we have $T_+^*=T_-$. We can use the following trick:
$$(T_++T_-)^2+\left( i(T_+-T_-)\right)^2=2(T_+T_-+T_-T_+)$$

Indeed, this gives $(T_++T_-)^2\leq 2(T_+T_-+T_-T_+)$, and we obtain in this way:
$$||T||^2=||T_++T_-||^2\leq 2||T_+T_-+T_-T_+||$$

Let $w\neq 1$ be a word, say beginning with $h$. We have then:
$$T_-T_+(\delta_w)=T_-(\delta_{gw}+\delta_{g^{-1}w}+\delta_{hw})=3\delta_w$$

The action of $T_-T_+$ on the remaining vector $\delta_1$ is computed as follows:
$$T_-T_+(\delta_1)=T_-(\delta_g+\delta_{g^{-1}}+\delta_{h}+\delta_{h^{-1}})=4\delta_1$$

Summing up, with $P:\delta_w\to\delta_1$ being the projection onto $\mathbb C\delta_1$, we have:
$$T_-T_+=3+P$$

On the other hand we have $T_+T_-(\delta_1)=T_+(0)=0$, so the subspace $\mathbb C\delta_1$ is invariant under the operator $T_+T_-+T_-T_+$. We have the following norm estimate:
$$||T||^2
\leq2||T_+T_-+T_-T_+||
\leq2\cdot\max\left\{||3+P||,\,\,\, ||(T_+T_-+T_-T_+)(1-P)||\right\}$$

The norm of $3+P$ is equal to $4$, and the other norm is estimated as follows:
\begin{eqnarray*}
||(T_+T_-+T_-T_+)(1-P)||
&\leq&||T_+T_-||+||(3+P)(1-P)||\\
&=&||T_-T_+||+3\\
&=&7
\end{eqnarray*}

Thus we have $||T||\leq\sqrt{14}<4$, and this finishes the proof.
\end{proof}

\section*{7c. Quantum groups}

The duals of discrete groups have several similarities with the compact groups, and our goal now will be that of unifying these two classes of compact quantum spaces. Let us start with the following definition, due to Woronowicz \cite{wo1}:

\index{Hopf algebra}
\index{Woronowicz algebra}

\begin{definition}
A Woronowicz algebra is a $C^*$-algebra $A$, given with a unitary matrix $u\in M_N(A)$ whose coefficients generate $A$, such that the formulae
$$\Delta(u_{ij})=\sum_ku_{ik}\otimes u_{kj}\quad,\quad
\varepsilon(u_{ij})=\delta_{ij}\quad,\quad
S(u_{ij})=u_{ji}^*$$
define morphisms of $C^*$-algebras $\Delta:A\to A\otimes A$, $\varepsilon:A\to\mathbb C$, $S:A\to A^{opp}$.
\end{definition}

We say that $A$ is cocommutative when $\Sigma\Delta=\Delta$, where $\Sigma(a\otimes b)=b\otimes a$ is the flip.  We have the following result, which justifies the terminology and axioms:

\index{cocommutative algebra}

\begin{proposition}
The following are Woronowicz algebras:
\begin{enumerate}
\item $C(G)$, with $G\subset U_N$ compact Lie group. Here the structural maps are:
$$\Delta(\varphi)=(g,h)\to \varphi(gh)\quad,\quad 
\varepsilon(\varphi)=\varphi(1)\quad,\quad
S(\varphi)=g\to\varphi(g^{-1})$$

\item $C^*(\Gamma)$, with $F_N\to\Gamma$ finitely generated group. Here the structural maps are:
$$\Delta(g)=g\otimes g\quad,\quad
\varepsilon(g)=1\quad,\quad
S(g)=g^{-1}$$
\end{enumerate}
Moreover, we obtain in this way all the commutative/cocommutative algebras.
\end{proposition}

\begin{proof}
In both cases, we have to exhibit a certain matrix $u$. For the first assertion, we can use the matrix $u=(u_{ij})$ formed by matrix coordinates of $G$, given by:
$$g=\begin{pmatrix}
u_{11}(g)&\ldots&u_{1N}(g)\\
\vdots&&\vdots\\
u_{N1}(g)&\ldots&u_{NN}(g)
\end{pmatrix}$$

As for the second assertion, here we can use the diagonal matrix formed by generators, $u=diag(g_1,\ldots,g_N)$. Finally, the last assertion follows from the Gelfand theorem, in the commutative case, and in the cocommutative case, we will be back to this later.
\end{proof}

In general now, the structural maps $\Delta,\varepsilon,S$ have the following properties:

\index{comultiplication}
\index{counit}
\index{antipode}
\index{square of antipode}

\begin{proposition}
Let $(A,u)$ be a Woronowicz algebra.
\begin{enumerate} 
\item $\Delta,\varepsilon$ satisfy the usual axioms for a comultiplication and a counit, namely:
$$(\Delta\otimes id)\Delta=(id\otimes \Delta)\Delta$$
$$(\varepsilon\otimes id)\Delta=(id\otimes\varepsilon)\Delta=id$$

\item $S$ satisfies the antipode axiom, on the $*$-subalgebra generated by entries of $u$: 
$$m(S\otimes id)\Delta=m(id\otimes S)\Delta=\varepsilon(.)1$$

\item In addition, the square of the antipode is the identity, $S^2=id$.
\end{enumerate}
\end{proposition}

\begin{proof}
When $A$ is commutative, by using Proposition 7.26 we can write:
$$\Delta=m^t\quad,\quad
\varepsilon=u^t\quad,\quad
S=i^t$$

The above 3 conditions come then by transposition from the basic 3 group theory conditions satisfied by $m,u,i$, which are as follows, with $\delta(g)=(g,g)$:
$$m(m\times id)=m(id\times m)$$
$$m(id\times u)=m(u\times id)=id$$
$$m(id\times i)\delta=m(i\times id)\delta=1$$

Observe that $S^2=id$ is satisfied as well, coming from $i^2=id$. In general now, all the formulae in the statement are satisfied on the generators $u_{ij}$, and so by linearity, multiplicativity and continuity they are satisfied everywhere, as desired.
\end{proof}

In view of Proposition 7.26, we can formulate the following definition:

\index{compact quantum group}
\index{discrete quantum group}
\index{Pontrjagin duality}

\begin{definition}
Given a Woronowicz algebra $A$, we formally write
$$A=C(G)=C^*(\Gamma)$$
and call $G$ compact quantum group, and $\Gamma$ discrete quantum group.
\end{definition}

When $A$ is both commutative and cocommutative, $G$ is a compact abelian group, $\Gamma$ is a discrete abelian group, and these groups are dual to each other:
$$G=\widehat{\Gamma}\quad,\quad 
\Gamma=\widehat{G}$$

In general, we still agree to write $G=\widehat{\Gamma},\Gamma=\widehat{G}$, in a formal sense. Finally, in relation with functoriality issues, let us complement Definitions 7.25 and 7.28 with:

\begin{definition}
Given two Woronowicz algebras $(A,u)$ and $(B,v)$, we write 
$$A\simeq B$$
and we identify as well the corresponding compact and discrete quantum groups, when we have an isomorphism of $*$-algebras $<u_{ij}>\simeq <v_{ij}>$, mapping $u_{ij}\to v_{ij}$.
\end{definition}

In order to develop now some theory, let us call corepresentation of $A$ any unitary matrix $v\in M_n(\mathcal A)$, with $\mathcal A=<u_{ij}>$, satisfying the same conditions as $u$, namely:
$$\Delta(v_{ij})=\sum_kv_{ik}\otimes v_{kj}\quad,\quad 
\varepsilon(v_{ij})=\delta_{ij}\quad,\quad 
S(v_{ij})=v_{ji}^*$$

These can be thought of as corresponding to the unitary representations of the underlying compact quantum group $G$. Following Woronowicz \cite{wo1}, we have:

\index{Haar integration}
\index{corepresentation}
\index{Ces\`aro limit}

\begin{theorem}
Any Woronowicz algebra has a unique Haar integration functional, 
$$\left(\int_G\otimes id\right)\Delta=\left(id\otimes\int_G\right)\Delta=\int_G(.)1$$
which can be constructed by starting with any faithful positive form $\varphi\in A^*$, and setting
$$\int_G=\lim_{n\to\infty}\frac{1}{n}\sum_{k=1}^n\varphi^{*k}$$
where $\phi*\psi=(\phi\otimes\psi)\Delta$. Moreover, for any corepresentation $v\in M_n(\mathbb C)\otimes A$ we have
$$\left(id\otimes\int_G\right)v=P$$
where $P$ is the orthogonal projection onto $Fix(v)=\{\xi\in\mathbb C^n|v\xi=\xi\}$.
\end{theorem}

\begin{proof}
Following \cite{wo1}, this can be done in 3 steps, as follows:

\medskip

(1) Given $\varphi\in A^*$, our claim is that the following limit converges, for any $a\in A$:
$$\int_\varphi a=\lim_{n\to\infty}\frac{1}{n}\sum_{k=1}^n\varphi^{*k}(a)$$

Indeed, by linearity we can assume that $a$ is the coefficient of corepresentation, $a=(\tau\otimes id)v$. But in this case, an elementary computation shows that we have the following formula, where $P_\varphi$ is the orthogonal projection onto the $1$-eigenspace of $(id\otimes\varphi)v$:
$$\left(id\otimes\int_\varphi\right)v=P_\varphi$$

(2) Since $v\xi=\xi$ implies $[(id\otimes\varphi)v]\xi=\xi$, we have $P_\varphi\geq P$, where $P$ is the orthogonal projection onto the space $Fix(v)=\{\xi\in\mathbb C^n|v\xi=\xi\}$. The point now is that when $\varphi\in A^*$ is faithful, by using a positivity trick, one can prove that we have $P_\varphi=P$. Thus our linear form $\int_\varphi$ is independent of $\varphi$, and is given on coefficients $a=(\tau\otimes id)v$ by:
$$\left(id\otimes\int_\varphi\right)v=P$$

(3) With the above formula in hand, the left and right invariance of $\int_G=\int_\varphi$ is clear on coefficients, and so in general, and this gives all the assertions. See \cite{wo1}. 
\end{proof}

As a main application, we can develop a Peter-Weyl type theory for the corepresentations of $A$. Consider the dense $*$-subalgebra $\mathcal A\subset A$ generated by the coefficients of the fundamental corepresentation $u$, and endow it with the following scalar product: 
$$<a,b>=\int_Gab^*$$

With this convention, we have the following result, also from Woronowicz \cite{wo1}:

\index{representation}
\index{corepresentation}
\index{character}
\index{Peter-Weyl}

\begin{theorem}
We have the following Peter-Weyl type results:
\begin{enumerate}
\item Any corepresentation decomposes as a sum of irreducible corepresentations.

\item Each irreducible corepresentation appears inside a certain $u^{\otimes k}$.

\item $\mathcal A=\bigoplus_{v\in Irr(A)}M_{\dim(v)}(\mathbb C)$, the summands being pairwise orthogonal.

\item The characters of irreducible corepresentations form an orthonormal system.
\end{enumerate}
\end{theorem}

\begin{proof}
All these results are from \cite{wo1}, the idea being as follows:

\medskip

(1) Given $v\in M_n(A)$, its intertwiner algebra $End(v)=\{T\in M_n(\mathbb C)|Tv=vT\}$ is a finite dimensional $C^*$-algebra, and so decomposes as $End(v)=M_{n_1}(\mathbb C)\oplus\ldots\oplus M_{n_r}(\mathbb C)$. But this gives a decomposition of type $v=v_1+\ldots+v_r$, as desired.

\medskip

(2) Consider indeed the Peter-Weyl corepresentations, $u^{\otimes k}$ with $k$ colored integer, defined by $u^{\otimes\emptyset}=1$, $u^{\otimes\circ}=u$, $u^{\otimes\bullet}=\bar{u}$ and multiplicativity. The coefficients of these corepresentations span the dense algebra $\mathcal A$, and by using (1), this gives the result.

\medskip

(3) Here the direct sum decomposition, which is technically a $*$-coalgebra isomorphism, follows from (2). As for the second assertion, this follows from the fact that $(id\otimes\int_G)v$ is the orthogonal projection $P_v$ onto the space $Fix(v)$, for any corepresentation $v$.

\medskip

(4) Let us define indeed the character of $v\in M_n(A)$ to be the matrix trace, $\chi_v=Tr(v)$. Since this character is a coefficient of $v$, the orthogonality assertion follows from (3). As for the norm 1 claim, this follows once again from $(id\otimes\int_G)v=P_v$. 
\end{proof}

We can now solve a problem that we left open before, namely:

\index{cocommutative algebra}

\begin{proposition}
The cocommutative Woronowicz algebras appear as the quotients
$$C^*(\Gamma)\to A\to C^*_{red}(\Gamma)$$
given by $A=C^*_\pi(\Gamma)$ with $\pi\otimes\pi\subset\pi$, with $\Gamma$ being a discrete group.
\end{proposition}

\begin{proof}
This follows from the Peter-Weyl theory, and clarifies a number of things said before, notably in Proposition 7.26. Indeed, for a cocommutative Woronowicz algebra the irreducible corepresentations are all 1-dimensional, and this gives the results.
\end{proof}

As another consequence of the above results, once again by following Woronowicz \cite{wo1}, we have the following statement, dealing with functional analysis aspects, and extending what we already knew about the $C^*$-algebras of the usual discrete groups:

\index{full algebra}
\index{reduced algebra}
\index{amenable quantum group}
\index{coamenable quantum group}
\index{Kesten amenability}

\begin{theorem}
Let $A_{full}$ be the enveloping $C^*$-algebra of $\mathcal A$, and $A_{red}$ be the quotient of $A$ by the null ideal of the Haar integration. The following are then equivalent:
\begin{enumerate}
\item The Haar functional of $A_{full}$ is faithful.

\item The projection map $A_{full}\to A_{red}$ is an isomorphism.

\item The counit map $\varepsilon:A_{full}\to\mathbb C$ factorizes through $A_{red}$.

\item We have $N\in\sigma(Re(\chi_u))$, the spectrum being taken inside $A_{red}$.
\end{enumerate}
If this is the case, we say that the underlying discrete quantum group $\Gamma$ is amenable.
\end{theorem}

\begin{proof}
This is well-known in the group dual case, $A=C^*(\Gamma)$, with $\Gamma$ being a usual discrete group. In general, the result follows by adapting the group dual case proof:

\medskip

$(1)\iff(2)$ This simply follows from the fact that the GNS construction for the algebra $A_{full}$ with respect to the Haar functional produces the algebra $A_{red}$.

\medskip

$(2)\iff(3)$ Here $\implies$ is trivial, and conversely, a counit map $\varepsilon:A_{red}\to\mathbb C$ produces an isomorphism $A_{red}\to A_{full}$, via a formula of type $(\varepsilon\otimes id)\Phi$. See \cite{wo1}.

\medskip

$(3)\iff(4)$ Here $\implies$ is clear, coming from $\varepsilon(N-Re(\chi (u)))=0$, and the converse can be proved by doing some functional analysis. Once again, we refer here to \cite{wo1}.
\end{proof}

Let us discuss now some interesting examples. Following Wang \cite{wan}, we have:

\index{free rotation}
\index{free orthogonal group}
\index{free unitary group}

\begin{proposition}
The following universal algebras are Woronowicz algebras,
$$C(O_N^+)=C^*\left((u_{ij})_{i,j=1,\ldots,N}\Big|u=\bar{u},u^t=u^{-1}\right)$$
$$C(U_N^+)=C^*\left((u_{ij})_{i,j=1,\ldots,N}\Big|u^*=u^{-1},u^t=\bar{u}^{-1}\right)$$
so the underlying spaces $O_N^+,U_N^+$ are compact quantum groups.
\end{proposition}

\begin{proof}
This follows from the elementary fact that if a matrix $u=(u_{ij})$ is orthogonal or biunitary, then so must be the following matrices:
$$u^\Delta_{ij}=\sum_ku_{ik}\otimes u_{kj}\quad,\quad 
u^\varepsilon_{ij}=\delta_{ij}\quad,\quad
u^S_{ij}=u_{ji}^*$$

Thus, we can indeed define morphisms $\Delta,\varepsilon,S$ as in Definition 7.25, by using the universal properties of $C(O_N^+)$, $C(U_N^+)$, and this gives the result.
\end{proof}

There is a connection here with group duals, coming from:

\begin{proposition}
Given a closed subgroup $G\subset U_N^+$, consider its ``diagonal torus'', which is the closed subgroup $T\subset G$ constructed as follows:
$$C(T)=C(G)\Big/\left<u_{ij}=0\Big|\forall i\neq j\right>$$
This torus is then a group dual, $T=\widehat{\Lambda}$, where $\Lambda=<g_1,\ldots,g_N>$ is the discrete group generated by the elements $g_i=u_{ii}$, which are unitaries inside $C(T)$.
\end{proposition}

\begin{proof}
Since $u$ is unitary, its diagonal entries $g_i=u_{ii}$ are unitaries inside $C(T)$. Moreover, from $\Delta(u_{ij})=\sum_ku_{ik}\otimes u_{kj}$ we obtain, when passing inside the quotient:
$$\Delta(g_i)=g_i\otimes g_i$$

It follows that we have $C(T)=C^*(\Lambda)$, modulo identifying as usual the $C^*$-completions of the various group algebras, and so that we have $T=\widehat{\Lambda}$, as claimed.
\end{proof}

With this notion in hand, we have the following result:

\begin{theorem}
The diagonal tori of the basic rotation groups are as follows,
$$\xymatrix@R=15mm@C=15mm{
U_N\ar[r]&U_N^+\\
O_N\ar[r]\ar[u]&O_N^+\ar[u]
}\qquad
\xymatrix@R=7mm@C=15mm{\\ :\\}
\qquad
\xymatrix@R=14mm@C=15mm{
\mathbb T^N\ar[r]&\widehat{F_N}\\
\mathbb Z_2^N\ar[r]\ar[u]&\widehat{\mathbb Z_2^{*N}}\ar[u]
}$$
where $F_N$ is the free group on $N$ generators, and $*$ is a group-theoretical free product.
\end{theorem}

\begin{proof}
This is clear indeed from $U_N^+$, and the other results can be obtained by imposing to the generators of $F_N$ the relations defining the corresponding quantum groups.
\end{proof}

As a conclusion to all this, the $C^*$-algebra theory suggests developing a theory of ``noncommutative geometry'', covering both the classical and the free geometry, by using compact quantum groups. We will be back to this in chapter 8.

\section*{7d. Cuntz algebras}

We would like to end this chapter with an interesting class of $C^*$-algebras, discovered by Cuntz in \cite{cun}, and heavily used since then, for various technical purposes. These algebras are not obviously related to the quantum space program that we have been developing so far, and might even look like some sort of Devil's invention, orthogonal to what is beautiful in operator algebras, but believe me, if planning to do some serious operator algebra work, you will certainly run into them. Their definition is as follows:

\index{Cuntz algebra}
\index{isometries}

\begin{definition}
The Cuntz algebra $O_n$ is the $C^*$-algebra generated by isometries $S_1,\ldots,S_n$ satisfying the following condition:
$$S_1S_1^*+\ldots+S_nS_n^*=1$$
That is, $O_n\subset B(H)$ is generated by $n$ isometries whose ranges sum up to $H$.
\end{definition}

Observe that $H$ must be infinite dimensional, in order to have isometries as above. In what follows we will prove that $O_n$ is independent on the choice of such isometries, and also that this algebra is simple. We will restrict the attention to the case $n=2$, the proof in general being similar. Let us start with some simple computations, as follows:

\begin{proposition}
Given a word $i=i_1\ldots i_k$ with $i_l\in\{1,2\}$, we associate to it the element $S_i=S_{i_1}\ldots S_{i_k}$ of the algebra $O_2$. Then $S_i$ are isometries, and we have
$$S_i^*S_j=\delta_{ij}1$$ 
for any two words $i,j$ having the same lenght.
\end{proposition}

\begin{proof}
We use the relations defining the algebra $O_2$, namely:
$$S_1^*S_1=S_2^*S_2=1\quad,\quad S_1S_1^*+S_2S_2^*=1$$

The fact that $S_i$ are isometries is clear, here being the check for $i=12$:
$$S_{12}^*S_{12}=(S_1S_2)^*(S_1S_2)=S_2^{*}S_1^*S_1S_2=S_2^{*}S_2=1$$

Regarding the last assertion, by recurrence we just have to establish the formula there for the words of length 1. That is, we want to prove the following formulae:
$$S_1^*S_2=S_2^*S_1=0$$

But these two formulae follow from the fact that the projections $P_i=S_iS_i^*$ satisfy by definition $P_1+P_2=1$. Indeed, we have the following computation:
\begin{eqnarray*}
P_1+P_2=1
&\implies&P_1P_2=0\\
&\implies&S_1S_1^*S_2S_2^*=0\\
&\implies&S_1^*S_2=S_1^*S_1S_1^*S_2S_2^*S_2=0
\end{eqnarray*}

Thus, we have the first formula, and the proof of the second one is similar.
\end{proof}

We can use the formulae in Proposition 7.38 as follows:

\begin{proposition}
Consider words in $O_2$, meaning products of $S_1,S_1^*,S_2,S_2^*$.
\begin{enumerate}
\item Each word in $O_2$ is of form $0$ or $S_iS_j^*$ for some words $i,j$.

\item Words of type $S_iS_j^*$ with $l(i)=l(j)=k$ form a system of $2^k\times 2^k$ matrix units.

\item The algebra $A_k$ generated by matrix units in (2) is a subalgebra of $A_{k+1}$.
\end{enumerate}
\end{proposition}

\begin{proof}
Here the first two assertions follow from the formulae in Proposition 7.38, and for the last assertion, we can use the following formula:
$$S_iS_j^*=S_i1S_j^*=S_i(S_1S_1^*+S_2S_2^*)S_j^*$$

Thus, we obtain an embedding of algebras $A_k$, as in the statement.
\end{proof}

Observe now that the embedding constructed in (3) above is compatible with the matrix unit systems in (2). Consider indeed the following diagram:
$$\begin{matrix}
A_{k+1}&\simeq &M_{2^{k+1}}(\mathbb C)\\
\ &\ &\ \\
\cup&\ &\cup\\
\ &\ &\ \\
A_k&\simeq&M_{2^k}(\mathbb C)
\end{matrix}$$

With the notation $e_{ix,yj}=e_{ij}\otimes e_{xy}$, the inclusion on the right is given by:
\begin{eqnarray*}
e_{ij}
&\to&e_{i1,1h}+e_{i2,2j}\\
&=&e_{ij}\otimes e_{11}+e_{ij}\otimes e_{22}\\
&=&e_{ij}\otimes 1
\end{eqnarray*}

Thus, with standard tensor product notations, the inclusion on the right is the canonical inclusion $m\to m\otimes1$, and so the above diagram becomes:
$$\begin{matrix}
A_{k+1}&\simeq &M_2(\mathbb C)^{\otimes k+1}\\
\ &\ &\ \\
\cup&\ &\cup\\
\ &\ &\ \\
A_k&\simeq &M_2(\mathbb C)^{\otimes k}
\end{matrix}$$

The passage from the algebra $A=\cup_kA_k\simeq M_2(\mathbb C)^{\otimes\infty}$ coming from this observation to the full the algebra $O_2$ that we are interested in can be done by using:

\begin{proposition}
Each element $X\in<S_1,S_2>\subset O_2$ decomposes as a finite sum
$$X=\sum_{i>0}S_1^{*i}X_{-i}+X_0+\sum_{i>0}X_iS_1^i$$
where each $X_i$ is in the union $A$ of algebras $A_k$. 
\end{proposition}

\begin{proof}
By linearity and by using Proposition 7.39 we may assume that $X$ is a nonzero word, say $X=S_iS_j^*$. In the case $l(i)=l(j)$ we can set $X_0=X$ and we are done. Otherwise, we just have to add at left or at right terms of the form $1=S_1^*S_1$. For instance $X=S_2$ is equal to $S_2S_1^*S_1$, and we can take $X_1=S_2S_1^*\in A_1$.
\end{proof}

We must show now that the decomposition $X\to (X_i)$ found above is unique, and then prove that each application $X\to X_i$ has good continuity properties. The following formulae show that in both problems we may restrict attention to the case $i=0$:
$$X_{i+1}=(XS_1^*)_i\hskip 2cm X_{-i-1}=(S_1X)_i$$

In order to solve these questions, we use the following fact:

\begin{proposition}
If $P$ is a nonzero projection in $\mathcal O_2=<S_1,S_2>\subset O_2$, its $k$-th average, given by the formula
$$Q=\sum_{l(i)=k}S_iPS_i^*$$
is a nonzero projection in $\mathcal O_2$ having the property that the linear subspace $QA_kQ$ is isomorphic to a matrix algebra, and $Y\to QYQ$ is an isomorphism of $A_k$ onto it.
\end{proposition}

\begin{proof}
We know that the words of form $S_iS_j^*$ with $l(i)=l(j)=k$ are a system of matrix units in $A_k$. We apply to them the map $Y\to QYQ$, and we obtain:
\begin{eqnarray*}
QS_iS_j^*Q
&=&\sum_{pq} S_pPS_p^*S_iS_j^*S_qPS_q^*\\
&=&\sum_{pq}\delta_{ip}\delta_{jq}S_pP^2S_q^*\\
&=&S_iPS_j^*
\end{eqnarray*}

The output being a system of matrix units, $Y\to QYQ$ is an isomorphism from the algebra of matrices $A_k$ to another algebra of matrices $QA_kQ$, and this gives the result.
\end{proof}

Thus any map $Y\to QYQ$ behaves well on the $i=0$ part of the decomposition on $X$. It remains to find $P$ such that $Y\to QYQ$ destroys all $i\neq 0$ terms, and we have here:

\begin{proposition}
Assuming $X_0\in A_k$, there is a nonzero projection $P\in A$ such that $QXQ=QX_0Q$, where $Q$ is the $k$-th average of $P$.
\end{proposition}

\begin{proof}
We want $Y\to QYQ$ to map to zero all terms in the decomposition of $X$, except for $X_0$. Let us call $M_1,\ldots,M_t\in\mathcal O_2-A$ the terms to be destroyed. We want the following equalities to hold, with the sum over all pairs of length $k$ indices:
$$\sum_{ij}S_iPS_i^*M_qS_jPS_j^*=0$$

The simplest way is to look for $P$ such that all terms of all sums are $0$:
$$S_iPS_i^*M_qS_jPS_j^*=0$$

By multiplying to the left by $S_i^*$ and to the right by $S_j$, we want to have:
$$PS_i^*M_qS_jP=0$$

With $N_z=S_i^*M_qS_j$, where $z$ belongs to some new index set, we want to have:
$$PN_zP=0$$

Since $N_z\in\mathcal O_2-A$, we can write $N_z=S_{m_z}S_{n_z}^*$ with $l(m_z)\neq l(n_z)$, and we want:
$$PS_{m_z}S_{n_z}^*P=0$$

In order to do this, we can the projections of form $P=S_rS_r^*$. We want:
$$S_rS_r^*S_{m_z}S_{n_z}^*S_rS_r^*=0$$

Let $K$ be the biggest length of all $m_z,n_z$. Assume that we have fixed $r$, of length bigger than $K$. If the above product is nonzero then both $S_r^*S_{m_z}$ and $S_{n_z}^*S_r$ must be nonzero, which gives the following equalities of words:
$$r_1\ldots r_{l(m_z)}=m_z\quad,\quad r_1\ldots r_{l(n_z)}=n_z$$

Assuming that these equalities hold indeed, the above product reduces as follows:
$$S_rS_{r_{l(r)}}^*\ldots S_{r_{l(m_z)+1}}^*S_{r_{l(n_z)+1}}\ldots S_{r_{l(r)}}S_r^*$$

Now if this product is nonzero, the middle term must be nonzero:
$$S_{r_{l(r)}}^*\ldots S_{r_{l(m_z)+1}}^*S_{r_{l(n_z)+1}}\ldots S_{r_{l(r)}}\neq 0$$

In order for this for hold, the indices starting from the middle to the right must be equal to the indices starting from the middle to the left. Thus $r$ must be periodic, of period $|l(m_z)-l(n_z)|>0$. But this is certainly possible, because we can take any aperiodic infinite word, and let $r$ be the sequence of first $M$ letters, with $M$ big enough.
\end{proof}

We can now start solving our problems. We first have:

\begin{proposition}
The decomposition of $X$ is unique, and we have
$$||X_i||\leq||X||$$
for any $i$.
\end{proposition}

\begin{proof}
It is enough to do this for $i=0$. But this follows from the previous result, via the following sequence of equalities and inequalities:
$$||X_0||=||QX_0Q||=||QXQ||\leq||X||$$

Thus we got the inequality in the statement. As for the uniqueness part, this follows from the fact that $X_0\to QX_0Q=QXQ$ is an isomorphism.
\end{proof}

Remember now we want to prove that the Cuntz algebra $O_2$ does not depend on the choice of the isometries $S_1,S_2$. In order to do so, let $\overline{\mathcal O}_2$ be the completion of the $*$-algebra $\mathcal O_2=<S_1,S_2>\subset O_2$ with respect to the biggest $C^*$-norm. We have:

\begin{proposition}
We have the equivalence
$$X=0\iff X_i=0,\forall i$$
valid for any element $X\in\overline{\mathcal O}_2$.
\end{proposition}

\begin{proof}
Assume $X_i=0$ for any $i$, and choose a sequence $X^k\to X$ with $X^k\in\mathcal O_2$. For $\lambda\in\mathbb T$ we define a representation $\rho_\lambda$ in the following way:
$$\rho_\lambda :S_i\to\lambda S_i$$

We have then $\rho_\lambda(Y)=Y$ for any element $Y\in A$. We fix norm one vectors $\xi,\eta$ and we consider the following continuous functions $f:\mathbb T\to\mathbb C$:
$$f^k(\lambda)=<\rho_\lambda (X^k)\xi,\eta>$$

From $X^k\to X$ we get, with respect to the usual sup norm of $C(\mathbb T)$:
$$f^k\to f$$

Each $X^k\in\mathcal O_2$ can be decomposed, and $f^k$ is given by the following formula:
$$f^k(\lambda )=\sum_{i>0}\lambda^{-i}<S_1^{*i}X^k_{-i}\xi,\eta>+<X_0\xi,\eta>+
\sum_{i>0}\lambda^i<X_i^kS_1^i\xi,\eta>$$

This is a Fourier type expansion of $f^k$, that can we write in the following way:
$$f^k(\lambda)=\sum_{j=-\infty}^\infty a_j^k\lambda^j$$

By using Proposition 7.43 we obtain that with $k\to\infty$, we have:
$$|a_j^k|\leq||X_j^k||\to||X_j^\infty||=0$$

On the other hand we have $a_j^k\to a_j$ with $k\to\infty$. Thus all Fourier coefficients $a_j$ of $f$ are zero, so $f=0$. With $\lambda=1$ this gives the following equality:
$$<X\xi,\eta>=0$$

This is true for arbitrary norm one vectors $\xi,\eta$, so $X=0$ and we are done.
\end{proof}

We can now formulate the Cuntz theorem, from \cite{cun}, as follows:

\index{Cuntz algebra}
\index{isometries}
\index{simple algebra}

\begin{theorem}[Cuntz]
Let $S_1,S_2$ be isometries satisfying $S_1S_1^*+S_2S_2^*=1$.
\begin{enumerate}

\item The $C^*$-algebra $O_2$ generated by $S_1,S_2$ does not depend on the choice of $S_1,S_2$.

\item For any nonzero $X\in O_2$ there are $A,B\in O_2$ with $AXB=1$.

\item In particular $O_2$ is simple.
\end{enumerate}
\end{theorem}

\begin{proof}
This basically follows from the various results established above:

\medskip

(1) Consider the canonical projection map $\pi:\overline{O}_2\to O_2$. We know that $\pi$ is surjective, and we will prove now that $\pi$ is injective. Indeed, if $\pi(X)=0$ then $\pi(X)_i=0$ for any $i$. But $\pi(X)_i$ is in the dense $*$-algebra $A$, so it can be regarded as an element of $\overline{O}_2$, and with this identification, we have $\pi(X)_i=X_i$ in $\overline{O}_2$. Thus $X_i=0$ for any $i$, so $X=0$. Thus $\pi$ is an isomorphism. On the other hand $\overline{O}_2$ depends only on $\mathcal O_2$, and the above formulae in $\mathcal O_2$, for algebraic calculus and for decomposition of an arbitrary $X\in\mathcal O_2$, show that $\mathcal O_2$ does not depend on the choice of $S_1,S_2$. Thus, we obtain the result.

\medskip

(2) Choose a sequence $X^k\to X$ with $X^k\in\mathcal O_2$. We have the following formula:
$$(X^*X)_0=\lim_{k\to\infty}\left( \sum_{i>0}X^{k*}_{-i}X^k_{-i}+X^{k*}_0X^k_0+
\sum_{i>0}S_1^{*i}X^{k*}_iX^k_iS_1^i\right)$$

Thus $X\neq 0$ implies $(X^*X)_0\neq 0$. By linearity we can assume that we have:
$$||(X^*X)_0||=1$$

Now choose a positive element $Y\in\mathcal O_2$ which is close enough to $X^*X$:
$$||X^*X-Y||<\varepsilon$$

Since $Z\to Z_0$ is norm decreasing, we have the following estimate:
$$||Y_0||>1-\varepsilon$$

We apply Proposition 7.42 to our positive element $Y\in\mathcal O_2$. We obtain in this way a certain projection $Q$ such that $QY_0Q=QYQ$ belongs to a certain matrix algebra. We have $QYQ>0$, so we can diagonalize this latter element, as follows:
$$QYQ=\sum\lambda_iR_i$$

Here $\lambda_i$ are positive numbers and $R_i$ are minimal projections in the matrix algebra. Now since $||QYQ||=||Y_0||$, there must be an eigenvalue greater that $1-\varepsilon$:
$$\lambda_0>1-\varepsilon$$

By linear algebra, we can pass from a minimal projection to another:
$$U^*U=R_i\quad,\quad UU^*=S_1^kS_1^{*k}$$

The element $B=QU^*S_1^k$ has norm $\leq 1$, and we get the following inequality:
\begin{eqnarray*}
||1-B^*X^*XB||
&\leq&||1-B^*YB||+||B^*YB-B^*X^*XB||\\
&<&||1-B^*YB||+\varepsilon
\end{eqnarray*}

The last term can be computed by using the diagonalization of $QYQ$, as follows:
\begin{eqnarray*}
B^*YB
&=&S_1^{*k}UQYQU^*S_1^k\\
&=&S_1^{*k}\left(\sum \lambda_i UR_iU^*\right) S_1^k\\
&=&\lambda_0S_1^{*k}S_1^kS_1^{*k}S_1^k\\
&=&\lambda_0
\end{eqnarray*}

From $\lambda_0> 1-\varepsilon$ we get $||1-B^*YB||<\varepsilon$, and we obtain the following estimate:
$$||1-B^*X^*XB||<2\varepsilon$$

Thus $B^*X^*XB$ is invertible, say with inverse $C$, and we have $(B^*X^*)X(BC)=1$.

\medskip

(3) This is clear from the formula $AXB=1$ established in (2).
\end{proof}

\section*{7e. Exercises}

We have seen many things in this chapter, and there are many potential exercises, on all this. We will be however short, and as unique, key exercise, we have:

\begin{exercise}
Work out the proof of the existence result for the Haar measure on a compact group $G$, as a particular case of the result proved for quantum groups.
\end{exercise}

This is of course something very standard, the problem being that of eliminating algebras, linear forms and other functional analysis notions from the proof for the quantum groups, as to have in the end something talking about spaces, and measures on them.

\chapter{Geometric aspects}

\section*{8a. Topology, K-theory}

This chapter is a continuation of the previous one, meant to be a grand finale to the $C^*$-algebra theory that we started to develop there, before getting back to more traditional von Neumann algebra material, following Murray, von Neumann and others. There are countless things to be said, and possible paths to be taken. En hommage to Connes, and his book \cite{co3}, which is probably the finest ever on $C^*$-algebras, we will adopt a geometric viewpoint. To be more precise, we know that a $C^*$-algebra is a beast of type $A=C(X)$, with $X$ being a compact quantum space. So, it is about the ``geometry'' of $X$ that we would like to talk about, everything else being rather of administrative nature.

\bigskip

Let us first look at the classical case, where $X$ is a usual compact space. You might say right away that wrong way, what we need for doing geometry is a manifold. But my answer here is modesty, and no hurry. It is right that you cannot do much geometry with a compact space $X$, but you can do some, and we have here, for instance:

\index{K-theory}
\index{vector bundle}

\begin{definition}
Given a compact space $X$, its first $K$-theory group $K_0(X)$ is the group of formal differences of complex vector bundles over $X$.
\end{definition}

This notion is quite interesting, and we can talk in fact about higher $K$-theory groups $K_n(X)$ as well, and all this is related to the homotopy groups $\pi_n(X)$ too. There are many non-trivial results on the subject, the end of the game being of course that of understanding the ``shape'' of $X$, that you need to know a bit about, before getting into serious geometry, in the case where $X$ happens to be a manifold.

\bigskip

As a question for us now, operator algebra theorists, we have:

\begin{question}
Can we talk about the first $K$-theory group $K_0(X)$ of a compact quantum space $X$?
\end{question}

We will see that this is a quite subtle question. To be more precise, we will see that we can talk, in a quite straightforward way, of the group $K_0(A)$ of an arbitrary $C^*$-algebra $A$, which is constructed as to have $K_0(A)=K_0(X)$ in the commutative case, where $A=C(X)$, with $X$ being a usual compact space. In the noncommutative case, however, $K_0(A)$ will sometimes depend on the choice of $A$ satisfying $A=C(X)$, and so all this will eventually lead to a sort of dead end, and to a rather ``no'' answer to Question 8.2.

\bigskip

Getting started now, in order to talk about the first $K$-theory group $K_0(A)$ of an arbitrary $C^*$-algebra $A$, we will need the following simple fact:

\index{projective module}

\begin{proposition}
Given a $C^*$-algebra $A$, the finitely generated projective $A$-modules $E$ appear via quotient maps $f:A^n\to E$, so are of the form
$$E=pA^n$$
with $p\in M_n(A)$ being an idempotent. In the commutative case, $A=C(X)$ with $X$ classical, these $A$-modules consist of sections of the complex vector bundles over $X$.
\end{proposition}

\begin{proof}
Here the first assertion is clear from definitions, via some standard algebra, and the second assertion is clear from definitions too, again via some algebra.
\end{proof}

With this in hand, let us go back to Definition 8.1. Given a compact space $X$, it is now clear that its $K$-theory group $K_0(X)$ can be recaptured from the knowledge of the associated $C^*$-algebra $A=C(X)$, and to be more precise we have $K_0(X)=K_0(A)$, when the first $K$-theory group of an arbitrary $C^*$-algebra is constructed as follows:

\begin{definition}
The first $K$-theory group of a $C^*$-algebra $A$ is the group of formal differences
$$K_0(A)=\big\{p-q\big\}$$
of equivalence classes of projections $p\in M_n(A)$, with the equivalence being given by
$$p\sim q\iff\exists u, uu^*=p,u^*u=q$$
and with the additive structure being the obvious one, by diagonal concatenation.
\end{definition}

This is very nice, and as a first example, we have $K_0(\mathbb C)=\mathbb Z$. More generally, as already mentioned above, it follows from Proposition 8.3 that in the commutative case, where $A=C(X)$ with $X$ being a compact space, we have $K_0(A)=K_0(X)$. Observe also that we have, by definition, the following formula, valid for any $n\in\mathbb N$:
$$K_0(A)=K_0(M_n(A))$$

Some further elementary observations include the fact that $K_0$ behaves well with respect to direct sums and with inductive limits, and also that $K_0$ is a homotopy invariant, and for details here, we refer to any introductory book on the subject, such as \cite{bla}.

\bigskip

In what concerns us, back to our Question 8.2, what has been said above is certainly not enough for investigating our question, and we need more examples. However, these examples are not easy to find, and for getting them, we need more theory. We have:

\begin{definition}
The second $K$-theory group of a $C^*$-algebra $A$ is the group of connected components of the unitary group of $GL_\infty(A)$, with 
$$GL_n(A)\subset GL_{n+1}(A)\quad,\quad a\to\begin{pmatrix}a&0\\0&1\end{pmatrix}$$
being the embeddings producing the inductive limit $GL_\infty(A)$.
\end{definition}

Again, for a basic example we can take $A=\mathbb C$, and we have here $K_1(\mathbb C)=\{1\}$, trivially. In fact, in the commutative case, where $A=C(X)$, with $X$ being a usual compact space, it is possible to establish a formula of type $K_1(A)=K_1(X)$. Further elementary observations include the fact that $K_1$ behaves well with respect to direct sums and with inductive limits, and also that $K_1$ is a homotopy invariant.

\bigskip

Importantly, the first and second $K$-theory groups are related, as follows:

\begin{theorem}
Given a $C^*$-algebra $A$, we have isomorphisms as follows, with
$$SA=\left\{f\in C([0,1],A)\Big|f(0)=0\right\}$$
standing for the suspension operation for the $C^*$-algebras:
\begin{enumerate}
\item $K_1(A)=K_0(SA)$.

\item $K_0(A)=K_1(SA)$.
\end{enumerate}
\end{theorem}

\begin{proof}
Here the isomorphism in (1) is something rather elementary, and the isomorphism in (2) is something more complicated. In both cases, the idea is to start first with the commutative case, where $A=C(X)$ with $X$ being a compact space, and understand there the isomorphisms (1,2), called Bott periodicity isomorphisms. Then, with this understood, the extension to the general $C^*$-algebra case is quite straightforward.
\end{proof}

The above result is quite interesting, making it clear that the groups $K_0,K_1$ are of the same nature. In fact, it is possible to be a bit more abstract here, and talk in various clever ways about the higher $K$-theory groups, $K_n(A)$ with $n\in\mathbb N$, of an arbitrary $C^*$-algebra, with the result that these higher $K$-theory groups are subject to Bott periodicity:
$$K_n(A)=K_{n+2}(A)$$

However, in practice, this leads us back to Definition 8.4, Definition 8.5 and Theorem 8.6, with these statements containing in fact all we need to know, at $n=0,1$.

\bigskip

Going ahead with examples, following Cuntz \cite{cun} and related papers, we have:

\begin{theorem}
The $K$-theory groups of the Cuntz algebra $O_n$ are given by
$$K_0(O_n)=\mathbb Z_{n-1}\quad,\quad K_1(O_n)=\{1\}$$
with the equivalent projections $P_i=S_iS_i^*$ standing for the standard generator of $\mathbb Z_{n-1}$.
\end{theorem}

\begin{proof}
We recall that the Cuntz algebra $O_n$ is generated by isometries $S_1,\ldots,S_n$ satisfying $S_1S_1^*+\ldots+S_nS_n^*=1$. Since we have $S_i^*S_i=1$, with $P_i=S_iS_i^*$, we have:
$$P_1\sim\ldots\sim P_n\sim 1$$

On the other hand, we also know that we have $P_1+\ldots+P_n=1$, and the conclusion is that, in the first $K$-theory group $K_1(O_n)$, the following equality happens:
$$n[1]=[1]$$

Thus $(n-1)[1]=0$, and it is quite elementary to prove that $k[1]=0$ happens in fact precisely when $k$ is a multiple of $n-1$. Thus, we have a group embedding, as follows:
$$\mathbb Z_{n-1}\subset K_0(O_n)$$

The whole point now is that of proving that this group embedding is an isomorphism, which in practice amounts in proving that any projection in $O_n$ is equivalent to a sum of the form $P_1+\ldots+P_k$, with $P_i=S_iS_i^*$ as above. Which is something non-trivial, requiring the use of Bott periodicity, and the consideration of the second $K$-theory group $K_1(O_n)$ as well, and for details here, we refer to Cuntz \cite{cun} and related papers.
\end{proof}

The above result is very interesting, for various reasons. First, it shows that the structure of the first $K$-theory groups $K_0(A)$ of the arbitrary $C^*$-algebras can be more complicated than that of the first $K$-theory groups $K_0(X)$ of the usual compact spaces $X$, with the group $K_0(A)$ being for instance not ordered, in the case $A=O_n$, and with this being the first in a series of no-go observations that can be formulated.

\bigskip

Second, and on a positive note now, what we have in Theorem 8.7 is a true noncommutative computation, dealing with an algebra which is rather of ``free'' type. The outcome of the computation is something nice and clear, suggesting that, modulo the small technical issues mentioned above, we are on our way of developing a nice theory, and that the answer to Question 8.2 might be ``yes''. However, as bad news, we have:

\begin{theorem}
There are discrete groups $\Gamma$ having the property that the projection
$$\pi:C^*(\Gamma)\to C^*_{red}(\Gamma)$$
is not an isomorphism, at the level of $K$-theory groups.
\end{theorem}

\begin{proof}
For constructing such a counterexample, the group $\Gamma$ must be definitely non-amenable, and the first thought goes to the free group $F_2$. But it is possible to prove that $F_2$ is $K$-amenable, in the sense that $\pi$ is an isomorphism at the $K$-theory level. However, counterexamples do exist, such as the infinite groups $\Gamma$ having Kazhdan's property $(T)$. Indeed, for such a group the asssociated Kazhdan projection $p\in K_0(C^*(\Gamma))$ is nonzero, while mapping to the zero element $0\in K_0(C^*_{red}(\Gamma))$, so we have our counterexample.  
\end{proof}

As a conclusion to all this, which might seem a bit dissapointing, we have:

\begin{conclusion}
The answer to Question 8.2 is no.
\end{conclusion}

Of course, the answer to Question 8.2 remains ``yes'' in many cases, the general idea being that, as long as we don't get too far away from the classical case, the answer remains ``yes'', so we can talk about the $K$-theory groups of our compact quantum spaces $X$, and also, about countless other invariants inspired from the classical theory. For a survey of what can be done here, including applications too, we refer to Connes' book \cite{co3}.

\bigskip

In what concerns us, however, we will not take this path. For various reasons, coming from certain quantum physics beliefs, which can be informally summarized as ``at sufficiently tiny scales, freeness rules'', we will be rather interested in this book in compact quantum spaces $X$ which are of ``free'' type, and we will only accept geometric invariants for them which are well-defined. And $K$-theory, unfortunately, does not qualify.

\section*{8b. Free probability}

As a solution to the difficulties met in the previous section, let us turn to probability. This is surely not geometry, in a standard sense, but at a more advanced level, geometry that is. For instance if you have a quantum manifold $X$, and you want to talk about its Laplacian, or its Dirac operator, you will certainly need to know a bit about $L^2(X)$. And isn't advanced measure theory the same as probability theory, hope we agree on this.

\bigskip

Let us start our discussion with something that we know since chapter 5:

\index{random variable}
\index{moments}
\index{law}
\index{distribution}

\begin{definition}
Let $A$ be a $C^*$-algebra, given with a trace $tr:A\to\mathbb C$.
\begin{enumerate}
\item The elements $a\in A$ are called random variables.

\item The moments of such a variable are the numbers $M_k(a)=tr(a^k)$.

\item The law of such a variable is the functional $\mu:P\to tr(P(a))$.
\end{enumerate}
\end{definition}

Here the exponent $k=\circ\bullet\bullet\circ\ldots$ is as before a colored integer, with the powers $a^k$ being defined by multiplicativity and the usual formulae, namely:
$$a^\emptyset=1\quad,\quad
a^\circ=a\quad,\quad
a^\bullet=a^*$$

As for the polynomial $P$, this is a noncommuting $*$-polynomial in one variable: 
$$P\in\mathbb C<X,X^*>$$

Generally speaking, the above definition is something quite abstract, but there is no other way of doing things, at least at this level of generality. However, in the special case where our variable $a\in A$ is self-adjoint, or more generally normal, we have:

\index{spectral measure}

\begin{proposition}
The law of a normal variable $a\in A$ can be identified with the corresponding spectral measure $\mu\in\mathcal P(\mathbb C)$, according to the following formula,
$$tr(f(a))=\int_{\sigma(a)}f(x)d\mu(x)$$
valid for any $f\in L^\infty(\sigma(a))$, coming from the measurable functional calculus. In the self-adjoint case the spectral measure is real, $\mu\in\mathcal P(\mathbb R)$.
\end{proposition}

\begin{proof}
This is something that we again know well, either from chapter 5, or simply from chapter 3, coming from the spectral theorem for normal operators.
\end{proof}

Let us discuss now independence, and its noncommutative versions. As a starting point, we have the following update of the classical notion of independence:

\index{independence}

\begin{definition}
We call two subalgebras $B,C\subset A$ independent when the following condition is satisfied, for any $x\in B$ and $y\in C$: 
$$tr(xy)=tr(x)tr(y)$$
Equivalently, the following condition must be satisfied, for any $x\in B$ and $y\in C$: 
$$tr(x)=tr(y)=0\implies tr(xy)=0$$
Also, $b,c\in A$ are called independent when $B=<b>$ and $C=<c>$ are independent.
\end{definition}

It is possible to develop some theory here, but this leads to the usual CLT. As a much more interesting notion now, we have Voiculescu's freeness \cite{vo1}:

\index{freeness}

\begin{definition}
Given a pair $(A,tr)$, we call two subalgebras $B,C\subset A$ free when the following condition is satisfied, for any $x_i\in B$ and $y_i\in C$:
$$tr(x_i)=tr(y_i)=0\implies tr(x_1y_1x_2y_2\ldots)=0$$
Also, $b,c\in A$ are called free when $B=<b>$ and $C=<c>$ are free.
\end{definition}

As a first observation, there is a certain lack of symmetry between Definition 8.12 and Definition 8.13, because the latter does not include an explicit formula for quantities of type $tr(x_1y_1x_2y_2\ldots)$. But this can be done, the precise result being as follows:

\begin{proposition}
If $B,C\subset A$ are free, the restriction of $tr$ to $<B,C>$ can be computed in terms of the restrictions of $tr$ to $B,C$. To be more precise, we have
$$tr(x_1y_1x_2y_2\ldots)=P\Big(\{tr(x_{i_1}x_{i_2}\ldots)\}_i,\{tr(y_{j_1}y_{j_2}\ldots)\}_j\Big)$$
where $P$ is certain polynomial, depending on the length of $x_1y_1x_2y_2\ldots\,$, having as variables the traces of products $x_{i_1}x_{i_2}\ldots$ and $y_{j_1}y_{j_2}\ldots\,$, with $i_1<i_2<\ldots$ and $j_1<j_2<\ldots$
\end{proposition}

\begin{proof}
With $x'=x-tr(x)$, we can start our computation as follows:
\begin{eqnarray*}
tr(x_1y_1x_2y_2\ldots)
&=&tr\big[(x_1'+tr(x_1))(y_1'+tr(y_1))(x_2'+tr(x_2))\ldots\big]\\
&=&tr(x_1'y_1'x_2'y_2'\ldots)+{\rm other\ terms}\\
&=&{\rm other\ terms}
\end{eqnarray*}

Thus, we are led to a kind of recurrence, and this gives the result.
\end{proof}

Let us discuss now some examples of independence and freeness. We first have the following result, from \cite{vo1}, which is something elementary:

\index{tensor product}
\index{free product}

\begin{proposition}
Given two algebras $(A,tr)$ and $(B,tr)$, the following hold:
\begin{enumerate}
\item $A,B$ are independent inside their tensor product $A\otimes B$, endowed with its canonical tensor product trace, given on basic tensors by $tr(a\otimes b)=tr(a)tr(b)$.

\item $A,B$ are free inside their free product $A*B$, endowed with its canonical free product trace, given by the formulae in Proposition 8.14.
\end{enumerate}
\end{proposition}

\begin{proof}
Both the assertions are indeed clear from definitions, with just some standard discussion needed for (2), in connection with the free product trace. See \cite{vo1}.
\end{proof}

More concretely now, we have the following result, also from Voiculescu \cite{vo1}:

\index{group algebra}
\index{independence}
\index{freeness}
\index{tensor product}
\index{free product}

\begin{proposition}
We have the following results, valid for group algebras:
\begin{enumerate}
\item $L(\Gamma),L(\Lambda)$ are independent inside $L(\Gamma\times\Lambda)$.

\item $L(\Gamma),L(\Lambda)$ are free inside $L(\Gamma*\Lambda)$.
\end{enumerate}
\end{proposition}

\begin{proof}
In order to prove these results, we can use the general results in Proposition 8.15, along with the following two isomorphisms, which are both standard:
$$L(\Gamma\times\Lambda)=L(\Lambda)\otimes L(\Gamma)\quad,\quad
L(\Gamma*\Lambda)=L(\Lambda)*L(\Gamma)$$

Alternatively, we can check the independence and freeness formulae on group elements, which is something trivial, and then conclude by linearity. See \cite{vo1}.
\end{proof}

We have already seen limiting theorems in classical probability, in chapter 6. In order to deal now with freeness, let us develop some tools. First, we have:

\index{free convolution}

\begin{proposition}
We have a well-defined operation $\boxplus$, given by
$$\mu_a\boxplus\mu_b=\mu_{a+b}$$
with $a,b$ being free, called free convolution.
\end{proposition}

\begin{proof}
We need to check here that if $a,b$ are free, then the distribution $\mu_{a+b}$ depends only on the distributions $\mu_a,\mu_b$. But for this purpose, we can use the formula in Proposition 8.14. Indeed, by plugging in arbitrary powers of $a,b$ as variables $x_i,y_j$, we obtain a family of formulae of the following type, with $Q$ being certain polyomials:
$$tr(a^{k_1}b^{l_1}a^{k_2}b^{l_2}\ldots)=P\Big(\{tr(a^k)\}_k,\{tr(b^l)\}_l\Big)$$

Thus the moments of $a+b$ depend only on the moments of $a,b$, and the same argument shows that the same holds for $*$-moments, and this gives the result.
\end{proof}

In order to advance now, we would need an analogue of the Fourier transform, or rather of the log of the Fourier transform. Quite remarkably, such a transform exists indeed, the precise result here, due to Voiculescu \cite{vo1}, being as follows:

\index{Cauchy transform}
\index{R-transform}

\begin{theorem}
Given a probability measure $\mu$, define its $R$-transform as follows:
$$G_\mu(\xi)=\int_\mathbb R\frac{d\mu(t)}{\xi-t}\implies G_\mu\left(\ R_\mu(\xi)+\frac{1}{\xi}\right)=\xi$$
The free convolution operation is then linearized by the $R$-transform.
\end{theorem}

\begin{proof}
This is something quite tricky, the idea being as follows:

\medskip

(1) In order to model the free convolution, the best is to use creation operators on free Fock spaces, corresponding to the semigroup von Neumann algebras $L(\mathbb N^{*k})$. Indeed, we have some freeness here, a bit in the same way as in the free group algebras $L(F_k)$.

\medskip

(2) The point now, motivating this choice, is that the variables of type $S^*+f(S)$, with $S\in L(\mathbb N)$ being the shift, and with $f\in\mathbb C[X]$ being an arbitrary  polynomial, are easily seen to model in moments all the possible distributions $\mu:\mathbb C[X]\to\mathbb C$.

\medskip

(3) Now let $f,g\in\mathbb C[X]$ and consider the variables $S^*+f(S)$ and $T^*+g(T)$, where $S,T\in L(\mathbb N*\mathbb N)$ are the shifts corresponding to the generators of $\mathbb N*\mathbb N$. These variables are free, and by using a $45^\circ$ argument, their sum has the same law as $S^*+(f+g)(S)$.

\medskip

(4) Thus the operation $\mu\to f$ linearizes the free convolution. We are therefore left with a computation inside $L(\mathbb N)$, which is elementary, and whose conclusion is that $R_\mu=f$ can be recaptured from $\mu$ via the Cauchy transform $G_\mu$, as in the statement.
\end{proof}

With the above linearization technology in hand, we can now establish the following remarkable free analogue of the CLT, also due to Voiculescu \cite{vo1}: 

\index{free CLT}
\index{semicircle law}
\index{Wigner law}
\index{free Gaussian law}

\begin{theorem}[Free CLT]
Given self-adjoint variables $x_1,x_2,x_3,\ldots,$ which are f.i.d., centered, with variance $t>0$, we have, with $n\to\infty$, in moments,
$$\frac{1}{\sqrt{n}}\sum_{i=1}^nx_i\sim\gamma_t$$\
where $\gamma_t=\frac{1}{2\pi t}\sqrt{4t-x^2}dx$ is the Wigner semicircle law of parameter $t$.
\end{theorem}

\begin{proof}
We follow the same idea as in the proof of the CLT:

\medskip 

(1) At $t=1$, the $R$-transform of the variable in the statement can be computed by using the linearization property from Theorem 8.18, and is given by:
$$R(\xi)
=nR_x\left(\frac{\xi}{\sqrt{n}}\right)
\simeq\xi$$

(2) On the other hand, some standard computations show that the Cauchy transform of the Wigner law $\gamma_1$ satisfies the following equation:
$$G_{\gamma_1}\left(\xi+\frac{1}{\xi}\right)=\xi$$

Thus, by using Theorem 8.18, we have the following formula:
$$R_{\gamma_1}(\xi)=\xi$$

(3) We conclude that the laws in the statement have the same $R$-transforms, and so they are equal. The passage to the general case, $t>0$, is routine, by dilation.
\end{proof}

In the complex case now, we have a similar result, also from \cite{vo1}, as follows:

\index{free CCLT}

\begin{theorem}[Free CCLT]
Given random variables $x_1,x_2,x_3,\ldots$ which are f.i.d., centered, with variance $t>0$, we have, with $n\to\infty$, in moments,
$$\frac{1}{\sqrt{n}}\sum_{i=1}^nx_i\sim\Gamma_t$$
where $\Gamma_t=law\big((a+ib)/\sqrt{2}\big)$, with $a,b$ being free, each following the Wigner semicircle law $\gamma_t$, is the Voiculescu circular law of parameter $t$.
\end{theorem}

\begin{proof}
This follows indeed from the free CLT, established before, simply by taking real and imaginary parts of all the variables involved.
\end{proof}

Now that we are done with the basic results in continuous case, let us discuss as well the discrete case. We can establish a free version of the PLT, as follows:

\index{FPLT}
\index{free PLT}
\index{free Poisson law}
\index{Marchenko-Pastur law}

\begin{theorem}[Free PLT]
The following limit converges, for any $t>0$,
$$\lim_{n\to\infty}\left(\left(1-\frac{t}{n}\right)\delta_0+\frac{t}{n}\delta_1\right)^{\boxplus n}$$
and we obtain the Marchenko-Pastur law of parameter $t$, 
$$\pi_t=\max(1-t,0)\delta_0+\frac{\sqrt{4t-(x-1-t)^2}}{2\pi x}\,dx$$
also called free Poisson law of parameter $t$.
\end{theorem}

\begin{proof}
Let $\mu$ be the measure in the statement, appearing under the convolution sign. The Cauchy transform of this measure is elementary to compute, given by:
$$G_{\mu}(\xi)=\left(1-\frac{t}{n}\right)\frac{1}{\xi}+\frac{t}{n}\cdot\frac{1}{\xi-1}$$

By using Theorem 8.18, we want to compute the following $R$-transform:
$$R
=R_{\mu^{\boxplus n}}(y)
=nR_\mu(y)$$

We know that the equation for this function $R$ is as follows:
$$\left(1-\frac{t}{n}\right)\frac{1}{y^{-1}+R/n}+\frac{t}{n}\cdot\frac{1}{y^{-1}+R/n-1}=y$$

With $n\to\infty$ we obtain from this the following formula:
$$R=\frac{t}{1-y}$$

But this being the $R$-transform of $\pi_t$, via some calculus, we are done.
\end{proof}

As a first application now of all this, following Voiculescu \cite{vo2}, we have:

\begin{theorem}
Given a sequence of complex Gaussian matrices $Z_N\in M_N(L^\infty(X))$, having independent $G_t$ variables as entries, with $t>0$, we have
$$\frac{Z_N}{\sqrt{N}}\sim\Gamma_t$$
in the $N\to\infty$ limit, with the limiting measure being Voiculescu's circular law.
\end{theorem}

\begin{proof}
We know from chapter 6 that the asymptotic moments are:
$$M_k\left(\frac{Z_N}{\sqrt{N}}\right)\simeq t^{|k|/2}|\mathcal{NC}_2(k)|$$

On the other hand, the free Fock space analysis done in the proof of Theorem 8.18 shows that we have, with the notations there, the following formulae:
$$S+S^*\sim\gamma_1\quad,\quad S+T^*\sim\Gamma_1$$

By doing some combinatorics, this shows that an abstract noncommutative variable $a\in A$ is circular, following the law $\Gamma_t$, precisely when its moments are:
$$M_k(a)=t^{|k|/2}|\mathcal{NC}_2(k)|$$

Thus, we are led to the conclusion in the statement. See \cite{vo2}.
\end{proof}

Next in line, comes the main result of Voiculescu in \cite{vo2}, as follows:

\index{Wigner matrix}
\index{asymptotic freeness}

\begin{theorem}
Given a family of sequences of Wigner matrices, 
$$Z^i_N\in M_N(L^\infty(X))\quad,\quad i\in I$$
with pairwise independent entries, each following the complex normal law $G_t$, with $t>0$, up to the constraint $Z_N^i=(Z_N^i)^*$, the rescaled sequences of matrices
$$\frac{Z^i_N}{\sqrt{N}}\in M_N(L^\infty(X))\quad,\quad i\in I$$
become with $N\to\infty$ semicircular, each following the Wigner law $\gamma_t$, and free.
\end{theorem}

\begin{proof}
We can assume that we are dealing with 2 sequences of matrices, $Z_N,Z_N'$. In order to prove the asymptotic freeness, consider the following matrix:
$$Y_N=\frac{1}{\sqrt{2}}(Z_N+iZ_N')$$

This is then a complex Gaussian matrix, so by using Theorem 8.22, we have:
$$\frac{Y_N}{\sqrt{N}}\sim\Gamma_t$$

We are therefore in the situation where $(Z_N+iZ_N')/\sqrt{N}$, which has asymptotically semicircular real and imaginary parts, converges to the distribution of a free combination of such variables. Thus $Z_N,Z_N'$ become asymptotically free, as desired.
\end{proof}

Getting now to the complex case, we have a similar result here, as follows:

\index{Gaussian matrix}
\index{asymptotic freeness}

\begin{theorem}
Given a family of sequences of complex Gaussian matrices, 
$$Z^i_N\in M_N(L^\infty(X))\quad,\quad i\in I$$
with pairwise independent entries, each following the law $G_t$, with $t>0$, the matrices
$$\frac{Z^i_N}{\sqrt{N}}\in M_N(L^\infty(X))\quad,\quad i\in I$$
become with $N\to\infty$ circular, each following the Voiculescu law $\Gamma_t$, and free.
\end{theorem}

\begin{proof}
This follows indeed from Theorem 8.23, which applies to the real and imaginary parts of our complex Gaussian matrices, and gives the result.
\end{proof}

Finally, we have as well a similar result for the Wishart matrices, as follows:

\begin{theorem}
Given a family of sequences of complex Wishart matrices, 
$$Z^i_N=Y^i_N(Y^i_N)^*\in M_N(L^\infty(X))\quad,\quad i\in I$$
with each $Y^i_N$ being a $N\times M$ matrix, with entries following the normal law $G_1$, and with all these entries being pairwise independent, the rescaled sequences of matrices
$$\frac{Z^i_N}{N}\in M_N(L^\infty(X))\quad,\quad i\in I$$
become with $M=tN\to\infty$ Marchenko-Pastur, each following the law $\pi_t$, and free.
\end{theorem}

\begin{proof}
Here the first assertion is the Marchenko-Pastur theorem, from chapter 6, and the second assertion follows from Theorem 8.23, or from Theorem 8.24.
\end{proof}

Let us develop now some further limiting theorems, classical and free. We have the following definition, extending the Poisson limit theory developed before:

\index{compound Poisson law}
\index{compound free Poisson law}

\begin{definition}
Associated to any compactly supported positive measure $\rho$ on $\mathbb C$ are the probability measures
$$p_\rho=\lim_{n\to\infty}\left(\left(1-\frac{c}{n}\right)\delta_0+\frac{1}{n}\rho\right)^{*n}\quad,\quad 
\pi_\rho=\lim_{n\to\infty}\left(\left(1-\frac{c}{n}\right)\delta_0+\frac{1}{n}\rho\right)^{\boxplus n}$$
where $c=mass(\rho)$, called compound Poisson and compound free Poisson laws.
\end{definition}

In what follows we will be interested in the case where $\rho$ is discrete, as is for instance the case for $\rho=t\delta_1$ with $t>0$, which produces the Poisson and free Poisson laws. The following result allows one to detect compound Poisson/free Poisson laws:

\index{Fourier transform}
\index{R-transform}

\begin{proposition}
For $\rho=\sum_{i=1}^sc_i\delta_{z_i}$ with $c_i>0$ and $z_i\in\mathbb C$, we have
$$F_{p_\rho}(y)=\exp\left(\sum_{i=1}^sc_i(e^{iyz_i}-1)\right)\quad,\quad 
R_{\pi_\rho}(y)=\sum_{i=1}^s\frac{c_iz_i}{1-yz_i}$$
where $F,R$ denote respectively the Fourier transform, and Voiculescu's $R$-transform.
\end{proposition}

\begin{proof}
Let $\mu_n$ be the measure appearing in Definition 8.26. We have:
\begin{eqnarray*}
F_{\mu_n}(y)=\left(1-\frac{c}{n}\right)+\frac{1}{n}\sum_{i=1}^sc_ie^{iyz_i}
&\implies&F_{\mu_n^{*n}}(y)=\left(\left(1-\frac{c}{n}\right)+\frac{1}{n}\sum_{i=1}^sc_ie^{iyz_i}\right)^n\\
&\implies&F_{p_\rho}(y)=\exp\left(\sum_{i=1}^sc_i(e^{iyz_i}-1)\right)
\end{eqnarray*}

In the free case we can use a similar method, and we obtain the above formula.
\end{proof}

We have the following result, providing an alternative to Definition 8.26, which will be our formulation here of the Compond Poisson Limit Theorem, classical and free:

\index{compound PLT}
\index{free compound PLT}

\begin{theorem}[CPLT]
For $\rho=\sum_{i=1}^sc_i\delta_{z_i}$ with $c_i>0$ and $z_i\in\mathbb C$, we have
$$p_\rho/\pi_\rho={\rm law}\left(\sum_{i=1}^sz_i\alpha_i\right)$$
where the variables $\alpha_i$ are Poisson/free Poisson$(c_i)$, independent/free.
\end{theorem}

\begin{proof}
This follows indeed from the fact that the the Fourier/$R$-transform of the variable in the statement is given by the formulae in Proposition 8.27.
\end{proof}

Following \cite{bb+}, \cite{bbc}, we will be interested here in the main examples of classical and free compound Poisson laws, which are constructed as follows:

\index{Bessel law}
\index{free Bessel law}

\begin{definition}
The Bessel and free Bessel laws are the compound Poisson laws
$$b^s_t=p_{t\varepsilon_s}\quad,\quad \beta^s_t=\pi_{t\varepsilon_s}$$
where $\varepsilon_s$ is the uniform measure on the $s$-th roots unity. In particular:
\begin{enumerate}
\item At $s=1$ we obtain the usual Poisson and free Poisson laws, $p_t,\pi_t$.

\item At $s=2$ we obtain the ``real'' Bessel and free Bessel laws, denoted $b_t,\beta_t$.

\item At $s=\infty$ we obtain the ``complex'' Bessel and free Bessel laws, denoted $B_t,\mathfrak B_t$.
\end{enumerate}
\end{definition}

There is a lot of theory regarding these laws, and we refer here to \cite{bb+}, \cite{bbc}, where these laws were introduced. We will be back to these laws, in a moment.

\section*{8c. Algebraic manifolds}

We are now ready, or almost, to develop some basic noncommutative geometry. The idea will be that of further building on the material from chapter 7, by enlarging the class of compact quantum groups studied there, with the consideration of quantum homogeneous spaces, $X=G/H$, and with classical and free probability as our main tools. 

\bigskip

But let us start with something intuitive, namely basic algebraic geometry, in a basic sense. The simplest compact manifolds that we know are the spheres, and if we want to have free analogues of these spheres, there are not many choices here, and we have:

\index{free sphere}

\begin{definition}
We have compact quantum spaces, constructed as follows,
$$C(S^{N-1}_{\mathbb R,+})=C^*\left(x_1,\ldots,x_N\Big|x_i=x_i^*,\sum_ix_i^2=1\right)$$
$$C(S^{N-1}_{\mathbb C,+})=C^*\left(x_1,\ldots,x_N\Big|\sum_ix_ix_i^*=\sum_ix_i^*x_i=1\right)$$
called respectively free real sphere, and free complex sphere.
\end{definition}

Observe that our spheres are indeed well-defined, due to the following estimate:
$$||x_i||^2=||x_ix_i^*||\leq\left|\left|\sum_ix_ix_i^*\right|\right|=1$$

Given a compact quantum space $X$, meaning as usual the abstract spectrum of a $C^*$-algebra, we define its classical version to be the classical space $X_{class}$ obtained by dividing $C(X)$ by its commutator ideal, then applying the Gelfand theorem:
$$C(X_{class})=C(X)/I\quad,\quad 
I=<[a,b]>$$

Observe that we have an embedding of compact quantum spaces $X_{class}\subset X$. In this situation, we also say that $X$ appears as a ``liberation'' of $X$. We have:

\index{liberation}
\index{classical version}
\index{commutator ideal}

\begin{proposition}
We have embeddings of compact quantum spaces
$$\xymatrix@R=15mm@C=15mm{
S^{N-1}_\mathbb C\ar[r]&S^{N-1}_{\mathbb C,+}\\
S^{N-1}_\mathbb R\ar[r]\ar[u]&S^{N-1}_{\mathbb R,+}\ar[u]
}$$
and the spaces on the right appear as liberations of the spaces of the left.
\end{proposition}

\begin{proof}
In order to prove this, we must establish the following isomorphisms:
$$C(S^{N-1}_\mathbb R)=C^*_{comm}\left(x_1,\ldots,x_N\Big|x_i=x_i^*,\sum_ix_i^2=1\right)$$
$$C(S^{N-1}_\mathbb C)=C^*_{comm}\left(x_1,\ldots,x_N\Big|\sum_ix_ix_i^*=\sum_ix_i^*x_i=1\right)$$

But these isomorphisms are both clear, by using the Gelfand theorem.
\end{proof}

We can now introduce a broad class of compact quantum manifolds, as follows:

\begin{definition}
A real algebraic submanifold $X\subset S^{N-1}_{\mathbb C,+}$ is a closed quantum space defined, at the level of the corresponding $C^*$-algebra, by a formula of type
$$C(X)=C(S^{N-1}_{\mathbb C,+})\Big/\Big<f_i(x_1,\ldots,x_N)=0\Big>$$
for certain noncommutative polynomials $f_i\in\mathbb C<X_1,\ldots,X_N>$. We identify two such manifolds, $X\simeq Y$, when we have an isomorphism of $*$-algebras of coordinates
$$\mathcal C(X)\simeq\mathcal C(Y)$$
mapping standard coordinates to standard coordinates.
\end{definition}

In practice, while our assumption $X\subset S^{N-1}_{\mathbb C,+}$ is definitely something technical, we are not losing much when imposing it, and we have the following list of examples:

\begin{proposition}
The following are algebraic submanifolds $X\subset S^{N-1}_{\mathbb C,+}$:
\begin{enumerate}
\item The spheres $S^{N-1}_\mathbb R\subset S^{N-1}_\mathbb C,S^{N-1}_{\mathbb R,+}\subset S^{N-1}_{\mathbb C,+}$.

\item Any compact Lie group, $G\subset U_n$, with $N=n^2$.

\item The duals $\widehat{\Gamma}$ of finitely generated groups, $\Gamma=<g_1,\ldots,g_N>$.

\item More generally, the closed subgroups $G\subset U_n^+$, with $N=n^2$.
\end{enumerate}
\end{proposition}

\begin{proof}
These facts are all well-known, the proofs being as follows:

\medskip

(1) This is indeed true by definition of our various spheres.

\medskip

(2) Given a closed subgroup $G\subset U_n$, we have an embedding $G\subset S^{N-1}_\mathbb C$, with $N=n^2$, given in double indices by $x_{ij}=u_{ij}/\sqrt{n}$, that we can further compose with the standard embedding $S^{N-1}_\mathbb C\subset S^{N-1}_{\mathbb C,+}$. As for the fact that we obtain indeed a real algebraic manifold, this is standard too, coming either from Lie theory or from Tannakian duality.

\medskip

(3) Given a group $\Gamma=<g_1,\ldots,g_N>$, consider the variables $x_i=g_i/\sqrt{N}$. These variables satisfy then the quadratic relations $\sum_ix_ix_i^*=\sum_ix_i^*x_i=1$ defining $S^{N-1}_{\mathbb C,+}$, and the algebricity claim for the manifold $\widehat{\Gamma}\subset S^{N-1}_{\mathbb C,+}$ is clear.

\medskip

(4) Given a closed subgroup $G\subset U_n^+$, we have indeed an embedding $G\subset S^{N-1}_{\mathbb C,+}$, with $N=n^2$, given by $x_{ij}=u_{ij}/\sqrt{n}$. As for the fact that we obtain indeed a real algebraic manifold, this comes from the Tannakian duality results in \cite{mal}, \cite{wo2}.
\end{proof}

Summarizing, what we have in Definition 8.32 is something quite fruitful, covering many interesting examples. In addition, all this is nice too at the axiomatic level, because the equivalence relation for our algebraic manifolds, as formulated in Definition 8.32, fixes in a quite clever way the functoriality issues of the Gelfand correspondence.

\bigskip

At the level of the general theory now, as a first tool that we can use, for the study of our manifolds, we have the following version of the Gelfand theorem:

\index{classical version}
\index{liberation}

\begin{theorem}
Assuming that $X\subset S^{N-1}_{\mathbb C,+}$ is an algebraic manifold, given by
$$C(X)=C(S^{N-1}_{\mathbb C,+})\Big/\Big<f_i(x_1,\ldots,x_N)=0\Big>$$
for certain noncommutative polynomials $f_i\in\mathbb C<X_1,\ldots,X_N>$, we have
$$X_{class}=\left\{x\in S^{N-1}_\mathbb C\Big|f_i(x_1,\ldots,x_N)=0\right\}$$
and $X$ itself appears as a liberation of $X_{class}$.
\end{theorem}

\begin{proof}
This is something that we know well for the spheres, from Proposition 8.31. In general, the proof is similar, coming from the Gelfand theorem.
\end{proof}

There are of course many other things that can be said about our manifolds, at the purely algebraic level. But in what follows we will be rather going towards analysis. 

\section*{8d. Free geometry}

We have now all the needed tools in our bag for developing ``free geometry''. The idea will be that of going back to the free quantum groups from chapter 7, and further building on that material, with a beginning of free geometry. Let us start with:

\index{quantum permutation}
\index{quantum reflection}
\index{reflection group}
\index{free quantum group}

\begin{theorem}
The classical and free, real and complex quantum rotation groups can be complemented with quantum reflection groups, as follows,
$$
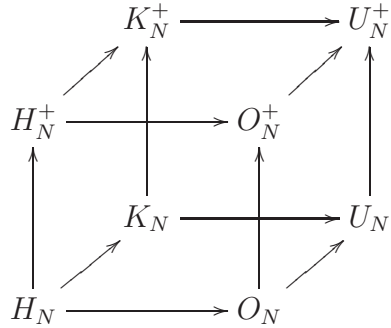
\xymatrix@R=18pt@C=18pt{
&K_N^+\ar[rr]&&U_N^+\\
H_N^+\ar[rr]\ar[ur]&&O_N^+\ar[ur]\\
&K_N\ar[rr]\ar[uu]&&U_N\ar[uu]\\
H_N\ar[uu]\ar[ur]\ar[rr]&&O_N\ar[uu]\ar[ur]
}$$
with $H_N=\mathbb Z_2\wr S_N$ and $K_N=\mathbb T\wr S_N$ being the hyperoctahedral group and the full complex reflection group, and $H_N^+=\mathbb Z_2\wr_*S_N^+$ and $K_N^+=\mathbb T\wr_*S_N^+$ being their free versions.
\end{theorem}

\begin{proof}
This is something quite tricky, the idea being as follows:

\medskip

(1) The first observation is that $S_N$, regarded as group of permutations of the $N$ coordinate axes of $\mathbb R^N$, is a group of orthogonal matrices, $S_N\subset O_N$. The corresponding coordinate functions $u_{ij}:S_N\to\{0,1\}$ form a matrix $u=(u_{ij})$ which is ``magic'', in the sense that its entries are projections, summing up to 1 on each row and each column. In fact, by using the Gelfand theorem, we have the following presentation result:
$$C(S_N)=C^*_{comm}\left((u_{ij})_{i,j=1,\ldots,N}\Big|u={\rm magic}\right)$$

(2) Based on the above, and following Wang's paper \cite{wan}, we can construct the free analogue $S_N^+$ of the symmetric group $S_N$ via the following formula:
$$C(S_N^+)=C^*\left((u_{ij})_{i,j=1,\ldots,N}\Big|u={\rm magic}\right)$$

Here the fact that we have indeed a Woronowicz algebra is standard, exactly as for the free rotation groups in chapter 7, because if a matrix $u=(u_{ij})$ is magic, then so are the matrices $u^\Delta,u^\varepsilon,u^S$ constructed there, and this gives the existence of $\Delta,u,S$.

\medskip

(3) Consider now the group $H_N^s\subset U_N$ consisting of permutation-like matrices having as entries the $s$-th roots of unity. This group decomposes as follows:
$$H_N^s=\mathbb Z_s\wr S_N$$

It is straightforward then to construct a free analogue $H_N^{s+}\subset U_N^+$ of this group, for instance by formulating a definition as follows, with $\wr_*$ being a free wreath product:
$$H_N^{s+}=\mathbb Z_s\wr_*S_N^+$$

(4) In order to finish, besides the case $s=1$, of particular interest are the cases $s=2,\infty$. Here the corresponding reflection groups are as follows:
$$H_N=\mathbb Z_2\wr S_N\quad,\quad K_N=\mathbb T\wr S_N$$

As for the corresponding quantum groups, these are denoted as follows:
$$H_N^+=\mathbb Z_2\wr_*S_N^+\quad,\quad K_N^+=\mathbb T\wr_*S_N^+$$

Thus, we are led to the conclusions in the statement. See \cite{bb+}, \cite{bbc}. 
\end{proof}

The point now is that we can add to the picture spheres and tori, as follows:

\index{free sphere}
\index{free torus}
\index{free manifold}

\begin{fact}
The basic quantum groups can be complemented with spheres and tori,
$$
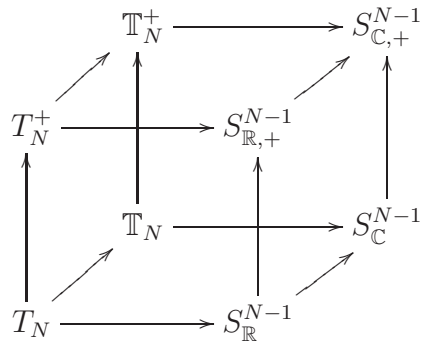
\xymatrix@R=16pt@C=15pt{
&\ \mathbb T_N^+\ar[rr]&&S^{N-1}_{\mathbb C,+}\\
\ T_N^+\ar[rr]\ar[ur]&&S^{N-1}_{\mathbb R,+}\ar[ur]\\
&\ \mathbb T_N\ar[rr]\ar[uu]&&S^{N-1}_\mathbb C\ar[uu]\\
\ T_N\ar[uu]\ar[ur]\ar[rr]&&S^{N-1}_\mathbb R\ar[uu]\ar[ur]
}$$
with $T_N=\mathbb Z_2^N,\mathbb T_N=\mathbb T^N$, and with $T_N^+,\mathbb T_N^+$ standing for the duals of $\mathbb Z_2^{*N},F_N$.
\end{fact}

Again, this is something quite tricky, and there is a long story with all this. We already know from chapter 7 that the diagonal subgroups of the rotation groups are the tori in the statement, but this is just an epsilon of what can be said, and this type of result can be extended as well to the reflection groups, and then we can make the spheres come into play too, with various results connecting them to the quantum groups, and to the tori.

\bigskip

Instead of getting into details here, let us formulate, again a bit informally:

\begin{fact}
The various quantum manifolds that we have, namely spheres $S$, tori $T$, unitary groups $U$, and reflection groups $K$, arrange into $4$ diagrams, as follows,
$$
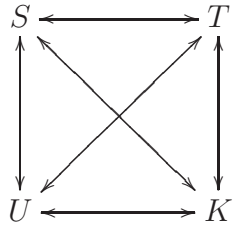
\xymatrix@R=58pt@C=58pt{
S\ar[r]\ar[d]\ar[dr]&T\ar[l]\ar[d]\ar[dl]\\
U\ar[u]\ar[ur]\ar[r]&K\ar[l]\ar[ul]\ar[u]
}$$
with the arrows standing for various correspondences between $(S,T,U,K)$. These diagrams correspond to $4$ main noncommutative geometries, real and complex, classical and free,
$$
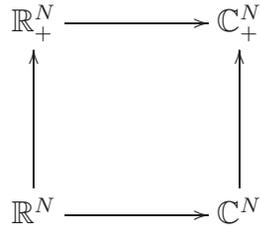
\xymatrix@R=54pt@C=54pt{
\mathbb R^N_+\ar[r]&\mathbb C^N_+\\
\mathbb R^N\ar[u]\ar[r]&\mathbb C^N\ar[u]
}$$
with the remark that, technically speaking, $\mathbb R^N_+$, $\mathbb C^N_+$ do not exist, as quantum spaces.
\end{fact}

As before, things here are quite long and tricky, but we already have some good evidence for all this, so I guess you can just trust me. And if truly interested in all this, later after finishing this book, you can check \cite{bgo} and subsequent papers for details. 

\bigskip

Summarizing, we have some beginning of theory. Now with this understood, let us try to integrate on our manifolds. In order to deal with quantum groups, we will need:

\index{Tannakian category}
\index{Peter-Weyl representations}
\index{tensor category}

\begin{definition}
The Tannakian category associated to a Woronowicz algebra $(A,u)$ is the collection $C_A=(C_A(k,l))$ of vector spaces
$$C_A(k,l)=Hom(u^{\otimes k},u^{\otimes l})$$
where the corepresentations $u^{\otimes k}$ with $k=\circ\bullet\bullet\circ\ldots$ colored integer, defined by
$$u^{\otimes\emptyset}=1\quad,\quad
u^{\otimes\circ}=u\quad,\quad 
u^{\otimes\bullet}=\bar{u}$$
and multiplicativity, $u^{\otimes kl}=u^{\otimes k}\otimes u^{\otimes l}$, are the Peter-Weyl corepresentations.
\end{definition}

As a key remark, the fact that $u\in M_N(A)$ is biunitary translates into the following conditions, where $R:\mathbb C\to\mathbb C^N\otimes\mathbb C^N$ is the linear map given by $R(1)=\sum_ie_i\otimes e_i$:
$$R\in Hom(1,u\otimes\bar{u})\quad,\quad 
R\in Hom(1,\bar{u}\otimes u)$$
$$R^*\in Hom(u\otimes\bar{u},1)\quad,\quad 
R^*\in Hom(\bar{u}\otimes u,1)$$

We are therefore led to the following abstract definition, summarizing the main properties of the categories appearing from Woronowicz algebras:

\begin{definition}
Let $H$ be a finite dimensional Hilbert space. A tensor category over $H$ is a collection $C=(C(k,l))$ of subspaces 
$$C(k,l)\subset\mathcal L(H^{\otimes k},H^{\otimes l})$$
satisfying the following conditions:
\begin{enumerate}
\item $S,T\in C$ implies $S\otimes T\in C$.

\item If $S,T\in C$ are composable, then $ST\in C$.

\item $T\in C$ implies $T^*\in C$.

\item Each $C(k,k)$ contains the identity operator.

\item $C(\emptyset,\circ\bullet)$ and $C(\emptyset,\bullet\circ)$ contain the operator $R:1\to\sum_ie_i\otimes e_i$.
\end{enumerate}
\end{definition}

The point now is that conversely, we can associate a Woronowicz algebra to any tensor category in the sense of Definition 8.39, in the following way:

\begin{proposition}
Given a tensor category $C=(C(k,l))$ over $\mathbb C^N$, as above,
$$A_C=C^*\left((u_{ij})_{i,j=1,\ldots,N}\Big|T\in Hom(u^{\otimes k},u^{\otimes l}),\forall k,l,\forall T\in C(k,l)\right)$$
is a Woronowicz algebra. 
\end{proposition}

\begin{proof}
This is something standard, because the relations $T\in Hom(u^{\otimes k},u^{\otimes l})$ determine a Hopf ideal, so they allow the construction of $\Delta,\varepsilon,S$ as in chapter 7.
\end{proof}

With the above constructions in hand, we have the following result:

\index{Tannakian duality}
\index{bicommutant}

\begin{theorem}
The Tannakian duality constructions 
$$C\to A_C\quad,\quad 
A\to C_A$$
are inverse to each other, modulo identifying full and reduced versions.
\end{theorem}

\begin{proof}
The idea is that we have $C\subset C_{A_C}$, for any algebra $A$, and so we are left with proving that we have $C_{A_C}\subset C$, for any category $C$. But this follows from a long series of algebraic manipulations, and for details we refer to Malacarne \cite{mal}, and also to Woronowicz \cite{wo2}, where this result was first proved, by using other methods.
\end{proof}

In practice now, all this is quite abstract, and we will rather need Brauer type results, for the specific quantum groups that we are interested in. Let us start with:

\index{category of partitions}

\begin{definition}
Let $P(k,l)$ be the set of partitions between an upper colored integer $k$, and a lower colored integer $l$. A collection of subsets 
$$D=\bigsqcup_{k,l}D(k,l)$$
with $D(k,l)\subset P(k,l)$ is called a category of partitions when it has the following properties:
\begin{enumerate}
\item Stability under the horizontal concatenation, $(\pi,\sigma)\to[\pi\sigma]$.

\item Stability under vertical concatenation $(\pi,\sigma)\to[^\sigma_\pi]$, with matching middle symbols.

\item Stability under the upside-down turning $*$, with switching of colors, $\circ\leftrightarrow\bullet$.

\item Each set $P(k,k)$ contains the identity partition $||\ldots||$.

\item The sets $P(\emptyset,\circ\bullet)$ and $P(\emptyset,\bullet\circ)$ both contain the semicircle $\cap$.
\end{enumerate}
\end{definition} 

Observe the similarity with Definition 8.39. In fact Definition 8.42 is a delinearized version of Definition 8.39, the relation with the Tannakian categories coming from:

\index{Kronecker symbol}

\begin{proposition}
Given a partition $\pi\in P(k,l)$, consider the linear map 
$$T_\pi:(\mathbb C^N)^{\otimes k}\to(\mathbb C^N)^{\otimes l}$$
given by the following formula, where $e_1,\ldots,e_N$ is the standard basis of $\mathbb C^N$, 
$$T_\pi(e_{i_1}\otimes\ldots\otimes e_{i_k})=\sum_{j_1\ldots j_l}\delta_\pi\begin{pmatrix}i_1&\ldots&i_k\\ j_1&\ldots&j_l\end{pmatrix}e_{j_1}\otimes\ldots\otimes e_{j_l}$$
and with the Kronecker type symbols $\delta_\pi\in\{0,1\}$ depending on whether the indices fit or not. The assignement $\pi\to T_\pi$ is then categorical, in the sense that we have
$$T_\pi\otimes T_\sigma=T_{[\pi\sigma]}\quad,\quad 
T_\pi T_\sigma=N^{c(\pi,\sigma)}T_{[^\sigma_\pi]}\quad,\quad 
T_\pi^*=T_{\pi^*}$$
where $c(\pi,\sigma)$ are certain integers, coming from the erased components in the middle.
\end{proposition}

\begin{proof}
The concatenation property follows from the following computation:
\begin{eqnarray*}
&&(T_\pi\otimes T_\sigma)(e_{i_1}\otimes\ldots\otimes e_{i_p}\otimes e_{k_1}\otimes\ldots\otimes e_{k_r})\\
&=&\sum_{j_1\ldots j_q}\sum_{l_1\ldots l_s}\delta_\pi\begin{pmatrix}i_1&\ldots&i_p\\j_1&\ldots&j_q\end{pmatrix}\delta_\sigma\begin{pmatrix}k_1&\ldots&k_r\\l_1&\ldots&l_s\end{pmatrix}e_{j_1}\otimes\ldots\otimes e_{j_q}\otimes e_{l_1}\otimes\ldots\otimes e_{l_s}\\
&=&\sum_{j_1\ldots j_q}\sum_{l_1\ldots l_s}\delta_{[\pi\sigma]}\begin{pmatrix}i_1&\ldots&i_p&k_1&\ldots&k_r\\j_1&\ldots&j_q&l_1&\ldots&l_s\end{pmatrix}e_{j_1}\otimes\ldots\otimes e_{j_q}\otimes e_{l_1}\otimes\ldots\otimes e_{l_s}\\
&=&T_{[\pi\sigma]}(e_{i_1}\otimes\ldots\otimes e_{i_p}\otimes e_{k_1}\otimes\ldots\otimes e_{k_r})
\end{eqnarray*}

As for the other two formulae in the statement, their proofs are similar.
\end{proof}

In relation with quantum groups, we have the following result, from \cite{bsp}:

\index{Tannakian duality}

\begin{theorem}
Each category of partitions $D=(D(k,l))$ produces a family of compact quantum groups $G=(G_N)$, one for each $N\in\mathbb N$, via the following formula:
$$Hom(u^{\otimes k},u^{\otimes l})=span\left(T_\pi\Big|\pi\in D(k,l)\right)$$
To be more precise, the spaces on the right form a Tannakian category, and so produce a certain closed subgroup $G_N\subset U_N^+$, via the Tannakian duality correspondence.
\end{theorem}

\begin{proof}
This follows indeed from Woronowicz's Tannakian duality, in its ``soft'' form from Malacarne \cite{mal}, as explained in Theorem 8.41. Indeed, let us set:
$$C(k,l)=span\left(T_\pi\Big|\pi\in D(k,l)\right)$$

By using the axioms in Definition 8.42, and the categorical properties of the operation $\pi\to T_\pi$, from Proposition 8.43, we deduce that $C=(C(k,l))$ is a Tannakian category. Thus the Tannakian duality applies, and gives the result.
\end{proof}

Philosophically speaking, the quantum groups appearing as in Theorem 8.44 are the simplest, from the perspective of Tannakian duality, so let us formulate:

\index{easy quantum group}

\begin{definition}
A closed subgroup $G\subset U_N^+$ is called easy when we have
$$Hom(u^{\otimes k},u^{\otimes l})=span\left(T_\pi\Big|\pi\in D(k,l)\right)$$
for any colored integers $k,l$, for a certain category of partitions $D\subset P$.
\end{definition}

All this might seem a bit complicated, but we will see examples in a moment. Getting back now to integration questions, we have the following key result:

\index{Weingarten formula}
\index{Haar integration}

\begin{theorem}
For an easy quantum group $G\subset U_N^+$, coming from a category of partitions $D=(D(k,l))$, we have the Weingarten integration formula
$$\int_Gu_{i_1j_1}^{e_1}\ldots u_{i_kj_k}^{e_k}=\sum_{\pi,\sigma\in D(k)}\delta_\pi(i)\delta_\sigma(j)W_{kN}(\pi,\sigma)$$
for any $k=e_1\ldots e_k$ and any $i,j$, where $D(k)=D(\emptyset,k)$, $\delta$ are usual Kronecker symbols, and $W_{kN}=G_{kN}^{-1}$, with $G_{kN}(\pi,\sigma)=N^{|\pi\vee\sigma|}$, where $|.|$ is the number of blocks. 
\end{theorem}

\begin{proof}
We know from chapter 7 that the integrals in the statement form altogether the orthogonal projection $P$ onto the space $Fix(u^{\otimes k})=span(D(k))$. Let us set:
$$E(x)=\sum_{\pi\in D(k)}<x,T_\pi>T_\pi$$

By standard linear algebra, it follows that we have $P=WE$, where $W$ is the inverse on $span(T_\pi|\pi\in D(k))$ of the restriction of $E$. But this restriction is the linear map given by $G_{kN}$, and so $W$ is the linear map given by $W_{kN}$, and this gives the result.
\end{proof}

All this is very nice. However, before enjoying the Weingarten formula, we still have to prove that our main quantum groups are easy. The result here is as follows:

\index{orthogonal group}
\index{unitary group}
\index{free orthogonal group}
\index{free unitary group}
\index{easy quantum group}
\index{Brauer theorem}

\begin{theorem}
The basic quantum unitary and reflection groups
$$
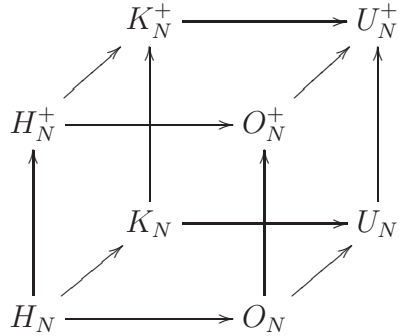
\xymatrix@R=19pt@C=19pt{
&K_N^+\ar[rr]&&U_N^+\\
H_N^+\ar[rr]\ar[ur]&&O_N^+\ar[ur]\\
&K_N\ar[rr]\ar[uu]&&U_N\ar[uu]\\
H_N\ar[uu]\ar[ur]\ar[rr]&&O_N\ar[uu]\ar[ur]
}$$
are all easy, the corresponding categories of partitions being as follows,
$$
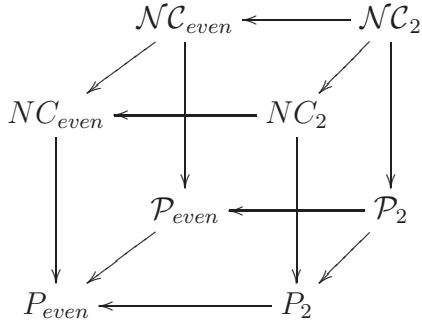
\xymatrix@R=19pt@C5pt{
&\mathcal{NC}_{even}\ar[dl]\ar[dd]&&\mathcal {NC}_2\ar[dl]\ar[ll]\ar[dd]\\
NC_{even}\ar[dd]&&NC_2\ar[dd]\ar[ll]\\
&\mathcal P_{even}\ar[dl]&&\mathcal P_2\ar[dl]\ar[ll]\\
P_{even}&&P_2\ar[ll]
}$$
with $P,NC$ standing for partitions and noncrosssing partitions, $2,even$ standing for pairings, and partitions with even blocks, and with calligraphic standing for matching.
\end{theorem}

\begin{proof}
The quantum group $U_N^+$ is defined via the following relations:
$$u^*=u^{-1}\quad,\quad 
u^t=\bar{u}^{-1}$$ 

Thus, the following operators must be in the associated Tannakian category:
$$T_\pi\quad,\quad \pi={\ }^{\,\cap}_{\circ\bullet}\ ,{\ }^{\,\cap}_{\bullet\circ}$$

We conclude that the associated Tannakian category is $span(T_\pi|\pi\in D)$, with:
$$D
=<{\ }^{\,\cap}_{\circ\bullet}\,\,,{\ }^{\,\cap}_{\bullet\circ}>
={\mathcal NC}_2$$

Thus, we have one result, and the other ones are similar. See \cite{bb+}, \cite{bbc}.
\end{proof}

We are not ready yet for applications, because we still have to understand which assumptions on $N\in\mathbb N$ make the vectors $T_\pi$ linearly independent. We will need:

\index{M\"obius function}

\begin{definition}
The M\"obius function of any lattice, and so of $P$, is given by
$$\mu(\pi,\sigma)=\begin{cases}
1&{\rm if}\ \pi=\sigma\\
-\sum_{\pi\leq\tau<\sigma}\mu(\pi,\tau)&{\rm if}\ \pi<\sigma\\
0&{\rm if}\ \pi\not\leq\sigma
\end{cases}$$
with the construction being performed by recurrence.
\end{definition}

The main interest in this function comes from the M\"obius inversion formula:
$$f(\sigma)=\sum_{\pi\leq\sigma}g(\pi)
\implies g(\sigma)=\sum_{\pi\leq\sigma}\mu(\pi,\sigma)f(\pi)$$

In linear algebra terms, the statement and proof of this formula are as follows:

\begin{proposition}
The inverse of the adjacency matrix of $P$, given by
$$A_{\pi\sigma}=\begin{cases}
1&{\rm if}\ \pi\leq\sigma\\
0&{\rm if}\ \pi\not\leq\sigma
\end{cases}$$
is the M\"obius matrix of $P$, given by $M_{\pi\sigma}=\mu(\pi,\sigma)$.
\end{proposition}

\begin{proof}
This is well-known, coming for instance from the fact that $A$ is upper triangular. Indeed, when inverting, we are led into the recurrence from Definition 8.48.
\end{proof}

Now back to our Gram and Weingarten matrix considerations, with $W_{kN}=G_{kN}^{-1}$, as in the statement of Theorem 8.46, we have the following result:

\begin{proposition}
The Gram matrix is given by $G_{kN}=AL$, where
$$L(\pi,\sigma)=
\begin{cases}
N(N-1)\ldots(N-|\pi|+1)&{\rm if}\ \sigma\leq\pi\\
0&{\rm otherwise}
\end{cases}$$
and where $A=M^{-1}$ is the adjacency matrix of $P(k)$.
\end{proposition}

\begin{proof}
We have indeed the following computation:
\begin{eqnarray*}
N^{|\pi\vee\sigma|}
&=&\#\left\{i_1,\ldots,i_k\in\{1,\ldots,N\}\Big|\ker i\geq\pi\vee\sigma\right\}\\
&=&\sum_{\tau\geq\pi\vee\sigma}\#\left\{i_1,\ldots,i_k\in\{1,\ldots,N\}\Big|\ker i=\tau\right\}\\
&=&\sum_{\tau\geq\pi\vee\sigma}N(N-1)\ldots(N-|\tau|+1)
\end{eqnarray*}

According to the definition of $G_{kN}$ and of $A,L$, this formula reads:
$$(G_{kN})_{\pi\sigma}
=\sum_{\tau\geq\pi}L_{\tau\sigma}
=\sum_\tau A_{\pi\tau}L_{\tau\sigma}
=(AL)_{\pi\sigma}$$

Thus, we obtain the formula in the statement.
\end{proof}

With the above result in hand, we can now formulate:

\index{Gram determinant}

\begin{theorem}
The determinant of the Gram matrix $G_{kN}$ is given by:
$$\det(G_{kN})=\prod_{\pi\in P(k)}\frac{N!}{(N-|\pi|)!}$$
In particular, the vectors $\left\{\xi_\pi|\pi\in P(k)\right\}$
are linearly independent  for $N\geq k$.
\end{theorem}

\begin{proof}
This is an old formula from the 60s, due to Lindst\"om and others, having many things behind it. By using the formula in Proposition 8.50, we have:
$$\det(G_{kN})=\det(A)\det(L)$$

Now if we order $P(k)$ with respect to the number of blocks, then lexicographically, $A$ is upper triangular, and $L$ is lower triangular, and we obtain the above formula.
\end{proof}

Now back to our quantum groups, let us start with:

\index{moments}
\index{main character}
\index{asymptotic moments}

\begin{theorem}
For an easy quantum group $G=(G_N)$, coming from a category of partitions $D=(D(k,l))$, the asymptotic moments of the character $\chi=\sum_iu_{ii}$ are
$$\lim_{N\to\infty}\int_{G_N}\chi^k=|D(k)|$$
where $D(k)=D(\emptyset,k)$, with the limiting sequence on the left consisting of certain integers, and being stationary at least starting from the $k$-th term.
\end{theorem}

\begin{proof}
This is something elementary, which follows straight from Peter-Weyl theory, by using the linear independence result from Theorem 8.51.
\end{proof}

In practice now, for the basic rotation and reflection groups, we obtain:

\begin{theorem}
The character laws for basic rotation and reflection groups are
$$
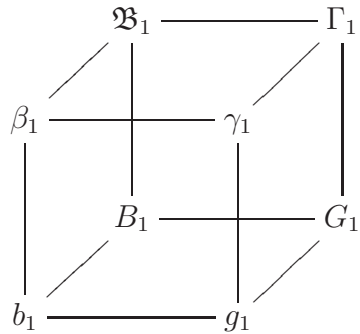
\xymatrix@R=20pt@C=20pt{
&\mathfrak B_1\ar@{-}[rr]\ar@{-}[dd]&&\Gamma_1\ar@{-}[dd]\\
\beta_1\ar@{-}[rr]\ar@{-}[dd]\ar@{-}[ur]&&\gamma_1\ar@{-}[dd]\ar@{-}[ur]\\
&B_1\ar@{-}[rr]\ar@{-}[uu]&&G_1\ar@{-}[uu]\\
b_1\ar@{-}[uu]\ar@{-}[ur]\ar@{-}[rr]&&g_1\ar@{-}[uu]\ar@{-}[ur]
}$$
in the $N\to\infty$ limit, corresponding to the basic probabilistic limiting theorems, at $t=1$.
\end{theorem}

\begin{proof}
This follows indeed from Theorem 8.47 and Theorem 8.52, by using the known moment formulae for the laws in the statement, at $t=1$.
\end{proof}

In the free case, the convergence can be shown to be stationary starting from $N=4$.  The ``fix'' comes by looking at truncated characters, constructed as follows:
$$\chi_t=\sum_{i=1}^{[tN]}u_{ii}$$

With this convention, we have the following final result on the subject, with the convergence being non-stationary at $t<1$, in both the classical and free cases:

\index{truncated character}
\index{normal law}
\index{Poisson law}
\index{semicircle law}
\index{free Poisson law}
\index{reflection group}

\begin{theorem}
The truncated character laws for the basic quantum groups are
$$
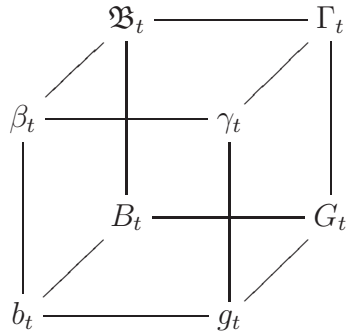
\xymatrix@R=20pt@C=20pt{
&\mathfrak B_t\ar@{-}[rr]\ar@{-}[dd]&&\Gamma_t\ar@{-}[dd]\\
\beta_t\ar@{-}[rr]\ar@{-}[dd]\ar@{-}[ur]&&\gamma_t\ar@{-}[dd]\ar@{-}[ur]\\
&B_t\ar@{-}[rr]\ar@{-}[uu]&&G_t\ar@{-}[uu]\\
b_t\ar@{-}[uu]\ar@{-}[ur]\ar@{-}[rr]&&g_t\ar@{-}[uu]\ar@{-}[ur]
}$$
in the $N\to\infty$ limit, corresponding to the basic probabilistic limiting theorems. 
\end{theorem}

\begin{proof}
We already know that the result holds at $t=1$, and the proof at arbitrary $t>0$ is once again based on easiness, but this time by using the Weingarten formula for the computation of the moments. We refer here to \cite{bb+}, \cite{bbc}, \cite{bco}, \cite{bsp}.
\end{proof}

All this is very nice, as a beginning. Of course, still left for this chapter would be the extension of all this to the case of more general homogeneous spaces $X=G/H$, and other free manifolds, in the sense of the free real and complex geometry axiomatized before. 

\bigskip

But hey, we learned enough math in this chapter, time for a beer. We refer here to the 2010 paper \cite{bgo}, which started things with the computation for $S^{N-1}_{\mathbb R,+}$, and to the book \cite{ba3}, which explains what was found on the subject, in the 10s. And if interested in this, the hot topic, waiting for input from you, are the applications to quantum physics.

\section*{8e. Exercises}

There has been a lot of exciting theory in this chapter, and as exercise, we have:

\begin{exercise}
Prove that $S_N^+$ is easy, coming from the category of all noncrossing partitions $NC$, and compute the asymptotic law of the main character.
\end{exercise}

As bonus exercise, try as well the truncated characters. Also, don't forget about $S_N$.

\part{Theory of factors}

\ \vskip50mm

\begin{center}
{\em And the story tellers say

That the score brave souls inside

For many a lonely day sailed across the milky seas

Never looked back, never feared, never cried}
\end{center}

\chapter{Functional analysis}

\section*{9a. Kaplansky density}

Welcome to this second half of the present book. We will get back here to a more normal pace, at least for most of the text to follow, our goal being to discuss the basics of the von Neumann algebra theory, due to Murray, von Neumann and Connes \cite{co1}, \cite{co2}, \cite{mv1}, \cite{mv2}, \cite{mv3}, \cite{vn1}, \cite{vn2}, \cite{vn3}, or at least the ``basics of the basics'', the whole theory being quite complex, and then the most beautiful advanced  theory which can be built on this, which is the subfactor theory of Jones \cite{jo1}, \cite{jo2}, \cite{jo3}, \cite{jo4}, \cite{jo5}, \cite{jms}, \cite{jsu}.

\bigskip

The material here will be in direct continuation of what we learned in chapter 5, namely bicommutant theorem, commutative case, finite dimensions, and a handful of other things. The idea will be that of building directly on that material, and using the same basic techniques, namely functional analysis and operator theory.

\bigskip 

As an important point, all this is related, but in a subtle way, to what we learned in chapters 6-8 too. To be more precise, what we will be doing in chapters 9-12 here will be more or less orthogonal to what we did in chapters 6-8. However, and here comes our point, the continuation of all this, chapters 13-16 below following Jones, will stand as a direct continuation of what we did in chapters 6-8, with Jones' subfactors being something more general than the random matrices and quantum groups from there.

\bigskip

Getting started, as a first objective we would like to have a better understanding of the precise difference between the norm closed $*$-algebras, or $C^*$-algebras, $A\subset B(H)$, and the weakly closed such algebras, which are the von Neumann algebras, from a functional analytic viewpoint. Let us begin with some generalities. We first have:

\begin{proposition}
The weak operator topology on $B(H)$ is the topology having the following equivalent properties:
\begin{enumerate}
\item It makes $T\to<Tx,y>$ continuous, for any $x,y\in H$.

\item It makes $T_n\to T$ when $<T_nx,y>\to<Tx,y>$, for any $x,y\in H$.

\item Has as subbase the sets $U_T(x,y,\varepsilon)=\{S:|<(S-T)x,y>|<\varepsilon\}$.

\item Has as base $U_T(x_1,\ldots,x_n,y_1,\ldots,y_n,\varepsilon)=\{S:|<(S-T)x_i,y_i>|<\varepsilon,\forall i\}$.
\end{enumerate}
\end{proposition}

\begin{proof}
The equivalences $(1)\iff(2)\iff(3)\iff(4)$ all follow from definitions, with of course (1,2) referring to the coarsest topology making that things happen.
\end{proof}

Similarly, in what regards the strong operator topology, we have:

\begin{proposition}
The strong operator topology on $B(H)$ is the topology having the following equivalent properties:
\begin{enumerate}
\item It makes $T\to Tx$ continuous, for any $x\in H$.

\item It makes $T_n\to T$ when $T_nx\to Tx$, for any $x\in H$.

\item Has as subbase the sets $V_T(x,\varepsilon)=\{S:||(S-T)x||<\varepsilon\}$.

\item Has as base the sets $V_T(x_1,\ldots,x_n,\varepsilon)=\{S:||(S-T)x_i||<\varepsilon,\forall i\}$.
\end{enumerate}
\end{proposition}

\begin{proof}
Again, the equivalences $(1)\iff(2)\iff(3)\iff(4)$ are all clear, and with (1,2) referring to the coarsest topology making that things happen.
\end{proof}

We know from chapter 5 that an operator algebra $A\subset B(H)$ is weakly closed if and only if it is strongly closed. Here is a useful generalization of this fact:

\begin{theorem}
Given a convex set $C\subset B(H)$, its weak operator closure and strong operator closure coincide.
\end{theorem}

\begin{proof}
Since the weak operator topology on $B(H)$ is weaker by definition than the strong operator topology on $B(H)$, we have, for any subset $C\subset B(H)$:
$$\overline{C}^{\,strong}\subset\overline{C}^{\,weak}$$

Now by assuming that $C\subset B(H)$ is convex, we must prove that:
$$T\in\overline{C}^{\,weak}\implies T\in\overline{C}^{\,strong}$$

In order to do so, let us pick vectors $x_1,\ldots,x_n\in H$ and $\varepsilon>0$. We let $K=H^{\oplus n}$, and we consider the standard embedding $i:B(H)\subset B(K)$, given by:
$$i T(y_1,\ldots,y_n)=(Ty_1,\ldots,Ty_n)$$

We have then the following implications, which are all trivial:
$$T\in\overline{C}^{\,weak}\implies iT\in\overline{iC}^{\,weak}\implies iT(x)\in\overline{iC(x)}^{\,weak}$$

Now since the set $C\subset B(H)$ was assumed to be convex, the set $iC(x)\subset K$ is convex too, and by the Hahn-Banach theorem, for compact sets, it follows that we have:
$$iT(x)\in\overline{iC(x)}^{\,||.||}$$

Thus, there exists an operator $S\in C$ such that we have, for any $i$:
$$||Sx_i-Tx_i||<\varepsilon$$

But this shows that we have $S\in V_T(x_1,\ldots,x_n,\varepsilon)$, and since $x_1,\ldots,x_n\in H$ and $\varepsilon>0$ were arbitrary, by Proposition 9.2 it follows that we have $T\in \overline{C}^{\,strong}$, as desired.
\end{proof}

We will need as well the following standard result:

\index{linear form}
\index{weak continuity}

\begin{proposition}
Given a vector space $E\subset B(H)$, and a linear form $f:E\to\mathbb C$, the following conditions are equivalent:
\begin{enumerate}
\item $f$ is weakly continuous.

\item $f$ is strongly continuous.

\item $f(T)=\sum_{i=1}^n<Tx_i,y_i>$, for certain vectors $x_i,y_i\in H$.
\end{enumerate}
\end{proposition}

\begin{proof}
This is something standard, using the same tools at those already used in chapter 5, namely basic functional analysis, and amplification tricks:

\medskip

$(1)\implies(2)$ Since the weak operator topology on $B(H)$ is weaker than the strong operator topology on $B(H)$, weakly continuous implies strongly continuous. To be more precise, assume $T_n\to T$ strongly. Then $T_n\to T$ weakly, and since $f$ was assumed to be weakly continuous, we have $f(T_n)\to f(T)$. Thus $f$ is strongly continuous, as desired.

\medskip

$(2)\implies(3)$ Assume indeed that our linear form $f:E\to\mathbb C$ is strongly continuous. In particular $f$ is strongly continuous at 0, and Proposition 9.2 provides us with vectors $x_1,\ldots,x_n\in H$ and a number $\varepsilon>0$ such that, with the notations there:
$$f(V_0(x_1,\ldots,x_n,\varepsilon))\subset D_0(1)$$

That is, we can find vectors $x_1,\ldots,x_n\in H$ and a number $\varepsilon>0$ such that:
$$||Tx_i||<\varepsilon,\forall i\implies|f(T)|<1$$

But this shows that we have the following estimate:
$$\sum_{i=1}^n||Tx_i||^2<\varepsilon^2\implies |f(T)|<1$$

By linearity, it follows from this that we have the following estimate:
$$|f(T)|<\frac{1}{\varepsilon}\sqrt{\sum_{i=1}^n||Tx_i||^2}$$

Consider now the direct sum $H^{\oplus n}$, and inside it, the following vector:
$$x=(x_1,\ldots,x_n)\in H^{\oplus n}$$

Consider also the following linear space, written in tensor product notation:
$$K=\overline{(E\otimes 1)x}\subset H^{\oplus n}$$

We can define a linear form $f':K\to\mathbb C$ by the following formula, and continuity:
$$f'(Tx_1,\ldots,Tx_n)=f(T)$$

We conclude that there exists a vector $y\in K$ such that the following happens:
$$f'\big((T\otimes1)y\big)=<(T\otimes 1)x,y>$$

But in terms of the original linear form $f:E\to\mathbb C$, this means that we have:
$$f(T)=\sum_{i=1}^n<Tx_i,y_i>$$

$(3)\implies(1)$ This is clear, because we have, with respect to the weak topology:
\begin{eqnarray*}
T_n\to T
&\implies&<T_nx_i,y_i>\to<Tx_i,y_i>,\forall i\\
&\implies&\sum_{i=1}^n<T_nx_i,y_i>\to\sum_{i=1}^n<Tx_i,y_i>\\
&\implies&f(T_n)\to f(T)
\end{eqnarray*}

Thus, our linear form $f$ is weakly continuous, as desired.
\end{proof}

Here is one more well-known result, that we will need as well:

\begin{theorem}
The unit ball of $B(H)$ is weakly compact.
\end{theorem}

\begin{proof}
If we denote by $B_1\subset B(H)$ the unit ball, and by $D_1\subset\mathbb C$ the unit disk, we have a morphism as follows, which is continuous with respect to the weak topology on $B_1$, and with respect to the product topology on the set on the right:
$$B_1\subset\prod_{||x||,||y||\leq1}D_1\quad,\quad T\to(<Tx,y>)_{x,y}$$

Since the set on the right is compact, by Tychonoff, it is enough to show that the image of $B_1$ is closed. So, let $(c_{xy})\in\bar{B}_1$. We can then find $T_i\in B_1$ such that:
$$<T_ix,y>\to c_{xy}\quad,\quad\forall x,y$$

But this shows that the following map is a bounded sesquilinear form:
$$H\times H\to\mathbb C\quad,\quad (x,y)\to c_{xy}$$

Thus, we can find an operator $T\in B(H)$, and so $T\in B_1$, such that $<Tx,y>=c_{xy}$ for any $x,y\in H$, and this shows that we have $(c_{xy})\in B_1$, as desired.
\end{proof}

Getting back to operator algebras, we have the following result, due to Kaplansky, which is something very useful, and of independent interest as well:

\index{unit ball}
\index{Kaplansky density}

\begin{theorem}
Given an operator algebra $A\subset B(H)$, the following happen:
\begin{enumerate}
\item The unit ball of $A$ is strongly dense in the unit ball of $A''$.

\item The same happens for the self-adjoint parts of the above unit balls.
\end{enumerate}
\end{theorem}

\begin{proof}
This is something quite tricky, the idea being as follows:

\medskip

(1) Consider the self-adjoint part $A_{sa}\subset A$. By taking real parts of operators, and using the fact that $T\to T^*$ is weakly continuous, we have then:
$$\overline{A_{sa}}^{\,w}\subset\left(\overline{A}^{\,w}\right)_{sa}$$

Now since the set $A_{sa}$ is convex, and by Theorem 9.3 all weak operator topologies coincide on the convex sets, we conclude that we have in fact equality:
$$\overline{A_{sa}}^{\,w}=\left(\overline{A}^{\,w}\right)_{sa}$$

(2) With this result in hand, let us prove now the second assertion of the theorem. For this purpose, consider an element $T\in\overline{A}^{\,w}$, satisfying $T=T^*$ and $||T||\leq1$. Consider as well the following function, going from the interval $[-1,1]$ to itself:
$$f(t)=\frac{2t}{1+t^2}$$

By functional calculus we can find an element $S\in\left(\overline{A}^{\,w}\right)_{sa}$ such that:
$$f(S)=T$$

In other words, we can find an element $S\in\left(\overline{A}^{\,w}\right)_{sa}$ such that:
$$T=\frac{2S}{1+S^2}$$

Now given arbitrary vectors $x_1,\ldots,x_n\in H$ and an arbitrary number $\varepsilon>0$, let us pick an element $R\in A_{sa}$, subject to the following two inequalities:
$$||RTx_i-STx_i||\leq\varepsilon\quad,\quad
\left|\left|\frac{R}{1+S^2}x_i-\frac{S}{1+S^2}x_i\right|\right|\leq\varepsilon$$

Finally, consider the following element, which has norm $\leq1$:
$$L=\frac{2R}{1+R^2}$$

We have then the following computation, using the above formulae:
\begin{eqnarray*}
L-T
&=&\frac{2R}{1+R^2}-\frac{2S}{1+S^2}\\
&=&2\left(\frac{1}{1+R^2}\big(R(1+S^2)-(1+S^2)R\big)\frac{1}{1+S^2}\right)\\
&=&2\left(\frac{1}{1+R^2}(R-S)\frac{1}{1+S^2}+\frac{R}{1+R^2}(S-R)\frac{S}{1+S^2}\right)\\
&=&\frac{2}{1+R^2}(R-S)\frac{1}{1+S^2}+\frac{L}{2}(S-R)T
\end{eqnarray*}

Thus, we have the following estimate, for any $i\in\{1,\ldots,n\}$:
$$||(L-T)x_i||\leq\varepsilon$$

But this gives the second assertion of the theorem, as desired.

\medskip

(3) Let us prove now the first assertion of the theorem. Given an arbitrary element $T\in\overline{A}^{\,w}$, satisfying $||T||\leq1$, let us look at the following element:
$$T'=\begin{pmatrix}0&T\\ T^*&0\end{pmatrix}\in M_2(\overline{A}^{\,w})$$

This element is then self-adjoint, and we can use what we proved in the above, and we are led in this way to the first assertion in the statement, as desired.
\end{proof}

We can go back now to our original question, from the beginning of the present chapter, namely that of abstractly characterizing the von Neumann algebras, and we have:

\begin{theorem}
An operator algebra $A\subset B(H)$ is a von Neumann algebra precisely when its unit ball is weakly compact.
\end{theorem}

\begin{proof}
This is something which is now clear, coming from the Kaplansky density results established in Theorem 9.6. To be more precise:

\medskip

(1) In one sense, assuming that $A\subset B(H)$ is a von Neumann algebra, this algebra is weakly closed. But since the unit ball of $B(H)$ is weakly compact, we are led to the conclusion that the unit ball of $A$ is weakly compact too.

\medskip

(2) Conversely, assume that an operator algebra $A\subset B(H)$ is such that  its unit ball is weakly compact. In particular, the unit ball of $A$ is weakly closed. Now if $T$ satisfying $||T||\leq1$ belongs to the weak closure of $A$, by Kaplansky density we conclude that we have $T\in A$. Thus our algebra $A$ must be a von Neumann algebra, as claimed.
\end{proof}

There are several other abstract characterizations of the von Neumann algebras, inside the class of $C^*$-algebras, and we will be back to this, on several occasions, and notably at the end of the present chapter, with such a characterization involving the predual.

\section*{9b. Projections, order}

In order to further investigate the von Neumann algebras, the key idea, coming from the analysis of the finite dimensional algebras from chapter 5, will be that of looking at the projections. Let us start with some generalities. In analogy with what happens in finite dimensions, we have the following notions, over an arbitrary Hilbert space $H$:

\begin{definition}
Associated to any two projections $P,Q\in B(H)$ are:
\begin{enumerate}
\item The projection $P\wedge Q$, projecting on the common range.

\item The projection $P\vee Q$, projecting on the span of the ranges.
\end{enumerate}
\end{definition} 

Abstractly speaking, these two operations can be thought of as being inf and sup type operations, and all the known algebraic formulae for inf and sup hold in this setting. For the moment we will not need all this, and we will be back to it later. Let us record however the following basic formula, which is something very useful:

\begin{proposition}
We have the following formula,
$$P+Q=P\wedge Q+P\vee Q$$
valid for any two projections $P,Q\in B(H)$.
\end{proposition}

\begin{proof}
This is clear from definitions, because when computing $P+Q$ we obtain the projection $P\vee Q$ on the span on the ranges, modulo the fact that the vectors in the common range are obtained twice, which amounts in saying that we must add $P\wedge Q$.
\end{proof}

With the above notions in hand, we have the following result:

\begin{theorem}
Consider two projections $P,Q\in B(H)$.
\begin{enumerate}
\item In finite dimensions, over $H=\mathbb C^N$, we have, in norm: 
$$(PQ)^n\to P\wedge Q$$

\item In infinite dimensions, we have the following convergence, for any $x\in H$,
$$(PQ)^nx\to (P\wedge Q)x$$
but the operators $(PQ)^n$ do not necessarily converge in norm.
\end{enumerate}
\end{theorem}

\begin{proof}
We have several assertions here, the proof being as follows:

\medskip

(1) Assume that we are in the case $P,Q\in M_N(\mathbb C)$. By substracting $P\wedge Q$ from both $P,Q$, we can assume $P\wedge Q=0$, and we must prove that we have:
$$P\wedge Q=0\implies(PQ)^n\to 0$$

Our claim is that we have $||PQ||<1$. Indeed, we know that we have:
$$||PQ||\leq||P||\cdot||Q||=1$$

Assuming now by contradiction that we have $||PQ||=1$, since we are in finite dimensions, we must have, for a certain norm one vector, $||x||=1$:
$$||PQx||=1$$

Thus, we must have equalities in the following estimate:
$$||PQx||\leq||Qx||\leq||x||$$

But the second equality tells us that we must have $x\in Im(Q)$, and with this in hand, the first equality tells us that we must have $x\in Im(P)$. But this contradicts $P\wedge Q=0$, so we have proved our claim, and the convergence $(PQ)^n\to 0$ follows.

\medskip

(2) In infinite dimensions now, as before by substracting $P\wedge Q$ from both $P,Q$, we can assume $P\wedge Q=0$, and we must prove that we have, for any $x\in H$:
$$P\wedge Q=0\implies(PQ)^nx\to 0$$

For this purpose, we use a trick. Consider the following operator:
$$R=PQP$$

This operator is positive, because we have $R=(PQ)(PQ)^*$, and we have:
$$||R||\leq||P||\cdot||Q||\cdot||P||=1$$

Our claim, which will finish the proof, is that for any $x\in H$ we have:
$$R^nx\to0$$

In order to prove this claim, let us diagonalize $R$, by using the spectral theorem for self-adjoint operators, from chapter 3. If all the eigenvalues are $<1$ then we are done. If not, this means that we can find a nonzero vector $x\in H$ such that:
$$||Rx||=||x||$$

But this condition means that we must have equalities in the following estimate:
$$||PQPx||\leq||QPx||\leq||Px||\leq||x||$$

The point now is that this is impossible, due to our assumption $P\wedge Q=0$. Indeed, the last equality tells us that we must have $x\in Im(P)$, and with this in hand, the middle equality tells us that we must have $x\in Im(Q)$. But this contradicts $P\wedge Q=0$, so we have proved our claim, and the convergence $(PQ)^nx\to 0$ follows.

\medskip

(3) Finally, for a counterexample to $(PQ)^n\to0$, in infinite dimensions, we can take $H=l^2(\mathbb N)$, and then find projections $P,Q$ such that $(PQ)^ne_k\to 0$ for any $k$, but with the convergence arbitrarily slowing down with $k\to\infty$. Thus, $(PQ)^n\not\to0$.
\end{proof}

As a consequence, in connection with the von Neumann algebras, we have:

\begin{theorem}
Given two projections $P,Q\in B(H)$, the projections
$$P\wedge Q\quad,\quad P\vee Q$$
both belong to the von Neumann algebra generated by $P,Q$.
\end{theorem}

\begin{proof}
This comes from the above. Indeed, in what regards $P\wedge Q$, this is something that follows from Theorem 9.10. As for $P\vee Q$, here the result follows from the result for $P\wedge Q$, and from the formula $P+Q=P\wedge Q+P\vee Q$, from Proposition 9.9.
\end{proof}

The idea now will be that of studying the von Neumann algebras $A\subset B(H)$ by using their projections, $p\in A$. Let us start with the following result:

\begin{theorem}
Any von Neumann algebra is generated by its projections.
\end{theorem}

\begin{proof}
This is something that we know from chapter 5, coming from the measurable functional calculus, which can cut any normal operator into projections.
\end{proof}

There are many other things that can be said about projections, in the general setting. In what follows we will just discuss the most important and useful such results. A first such result, providing us with some geometric intuition on projections, is as follows:

\index{corner of algebra}

\begin{theorem}
Given a von Neumann algebra $A\subset B(H)$, and a projection $p\in A$, we have the following equalities, between von Neumann algebras on $pH$:
\begin{enumerate}
\item $pAp=(A'p)'$.

\item $(pAp)'=A'p$.
\end{enumerate}
\end{theorem}

\begin{proof}
This is not exactly obvious, but can be proved as follows:

\medskip

(1) As a first observation, the von Neumann algebras $pAp$ and $A'p$ commute on $pH$. Thus, we must prove that we have the following implication:
$$x\in(A'p)'\implies x\in pAp$$

For this purpose, consider the element $y=xp$. Then for any $z\in A'$ we have:
\begin{eqnarray*}
zy
&=&zxp\\
&=&zpxp\\
&=&xpzp\\
&=&xpz\\
&=&yz
\end{eqnarray*}

But this shows that we have $y\in A$, and so we obtain, as desired:
$$x=pyp\in pAp$$

(2) As before, one of the inclusions being clear, we must prove that we have:
$$x\in(pAp)'\implies x\in A'p$$

By using the standard fact that any bounded operator appears as a linear combination of 4 unitaries, that we know from the beginning of chapter 4, it is enough to prove this for a unitary element, $x=u$. So, assume that we have a unitary as follows:
$$u\in(pAp)'$$

In order to prove our claim, consider the following vector space:
$$K=\overline{ApH}$$

This space being invariant under both the algebras $A,A'$, we conclude that the projection $q=Proj(K)$ onto it belongs to the center of our von Neumann algebra:
$$q\in Z(A)$$

Our claim now, which will quickly lead to the result that we want to prove, is that we can extend the above unitary $u\in(pAp)'$ to the space $K=\overline{ApH}$ via the following formula, valid for any elements $x_i\in A$, and any vectors $\xi_i\in pH$:
$$v\left(\sum_ix_i\xi_i\right)=\sum_ix_iu\xi_i$$

In order to prove this latter claim, we can use the following computation:
\begin{eqnarray*}
\left|\left|v\left(\sum_ix_i\xi_i\right)\right|\right|^2
&=&\sum_{ij}<x_iu\xi_i,x_ju\xi_j>\\
&=&\sum_{ij}<x_j^*x_iu\xi_i,u\xi_j>\\
&=&\sum_{ij}<px_j^*x_ipu\xi_i,u\xi_j>\\
&=&\sum_{ij}<upx_j^*x_ip\xi_i,u\xi_j>\\
&=&\sum_{ij}<px_j^*x_ip\xi_i,\xi_j>\\
&=&\sum_{ij}<x_j^*x_i\xi_i,\xi_j>\\
&=&\sum_{ij}<x_i\xi_i,x_j\xi_j>\\
&=&\left|\left|\sum_ix_i\xi_i\right|\right|^2
\end{eqnarray*}

Thus $v$ is well-defined by the above formula, and is an isometry of $K$. Now observe that this element $v$ commutes with the algebra $A$ on the space $ApH$, and so on $K$. Thus $vq\in A'$, and so $u=vqp$, which proves that we have $u\in A'p$, as desired.
\end{proof}

As a second result now, once again in the general setting, we have:

\index{equivalence of projections}

\begin{proposition}
Given a von Neumann algebra $A\subset B(H)$, the formula
$$p\simeq q\iff\exists u,\ \begin{cases}
uu^*=p\\
u^*u=q
\end{cases}$$
defines an equivalence relation for the projections $p\in A$.
\end{proposition}

\begin{proof}
This is something elementary, which follows from definitions, with the transitivity coming by composing the corresponding partial isometries.
\end{proof}

As a third result, once again in the general setting, which once again provides us with some intuition, but this time of somewhat abstract type, we have:

\index{lattice of projections}
\index{partial isometry}

\begin{theorem}
Given a von Neumann algebra $A\subset B(H)$, we have a partial order on the projections $p\in A$, constructed as follows, with $u$ being a partial isometry,
$$p\preceq q\iff\exists u,\ \begin{cases}
uu^*=p\\
u^*u\leq q
\end{cases}$$
which is related to the equivalence relation $\simeq$ constructed above by:
$$p\simeq q\iff p\preceq q,\ q\preceq p$$
Thus, $\preceq$ is a partial order on the equivalence classes of projections $p\in A$.
\end{theorem}

\begin{proof}
We have several assertions here, the idea being as follows:

\medskip

(1) The fact that we have indeed a partial order is clear, with the transitivity coming, as before, by composing the corresponding partial isometries.

\medskip

(2) Regarding now the relation with $\simeq$, via the equivalence in the statement, the implication $\implies$ is clear. Thus, we are left with proving $\Longleftarrow$, which reads:
$$p\preceq q,\ q\preceq p\implies p\simeq q$$

Our assumption is that we have partial isometries $u,v$ such that:
$$uu^*=p\quad,\quad u^*u\leq q$$
$$v^*v\leq p\quad,\quad vv^*=q$$

We can construct then two sequences of decreasing projections, as follows:
$$p_0=p\quad,\quad p_{n+1}=v^*q_nv$$
$$q_0=q\quad,\quad q_{n+1}=u^*p_nu$$

Consider now the limits of these two sequences of projections, namely:
$$p_\infty=\bigwedge_ip_i\quad,\quad q_\infty=\bigwedge_iq_i$$

In terms of all these projections that we constructed, we have the following decomposition formulae for the original projections $p,q$:
$$p=(p-p_1)+(p_1-p_2)+\ldots+p_\infty$$
$$q=(q-q_1)+(q_1-q_2)+\ldots+q_\infty$$

Now observe that the summands are equivalent, with this being clear from the  definition of $p_n,q_n$ at the finite indices $n<\infty$, and with $p_\infty\simeq q_\infty$ coming from:
$$v^*q_\infty v=p_\infty\quad,\quad q_\infty vv^*q_\infty=q_\infty$$

Thus we obtain that we have $p\simeq q$, as desired, by summing.

\medskip

(3) Finally, the fact that the order relation $\preceq$ factorizes indeed to the equivalence classes under $\simeq$ follows from the equivalence established in (2).
\end{proof}

Summarizing, in view of Theorem 9.12, and of Theorem 9.15, we can formulate:

\begin{conclusion}
We can think of a von Neumann algebra $A\subset B(H)$ as being a kind of object belonging to ``mathematical logic'', consisting of equivalence classes of projections $p\in A$, ordered via the relation $\preceq$, and producing $A$ itself via transport by partial isometries, and then linear combinations, and weak limits.
\end{conclusion}

Which is something quite remarkable, who on Earth could have guessed, say when we were back in chapter 5, struggling with the basics of the von Neumann algebra theory, or even at the beginning of the present chapter 9, again struggling with some sort of basics, of the more advanced theory, that we will end up with something that luminous.

\bigskip

Well, that person on Earth who found this was von Neumann himself, back in the 1930s. And his Conclusion 9.16, called ``von Neumann vision'' of the operator algebras, has been extremely useful ever since, and is still largely used nowadays. 

\bigskip

Very nice all this, first class mathematics, but in what concerns us, however, we will rather stick to our $A=L^\infty(X)$ viewpoint, with $X$ being a quantum measured space, and the most often being a ``quantum manifold''. This is more of a ``continuous'' philosophy, and in order to keep it intact, and powerful, we will have to take sometimes distances with the von Neumann philosophy, especially in what concerns the terminology. 

\bigskip

In short, we will be definitely users of the von Neumann projection technology, which is extremely powerful, and is quite often the only available tool, but keeping in mind however that we are dealing with continuous objects $X$, and choosing the terminology and notations accordingly, inspired from continuous geometry.

\section*{9c. States, isomorphism}

Getting back now to general questions concerning the von Neumann algebras, one question that we met on several occasions, and that we would like to clarify now, is the relation between abstract isomorphism and spatial isomorphism. 

\bigskip

To be more precise, we would like to understand when two von Neumann algebras $A\subset B(H)$ and $B\subset B(K)$ are isomorphic, in an algebraic and topological sense, but without reference to the ambient Hilbert spaces $H,K$. With the idea in mind that, once this understood, we will be able to talk about the von Neumann algebras $A$ as being abstract objects, a bit as were the $C^*$-algebras, discussed in chapter 7.

\bigskip

In order to discuss this, let us start with some technical preliminaries. Here is a definition that I have been postponing for long, but which is now unavoidable:

\begin{definition}
We call ultraweak and ultrastrong topologies on $B(H)$ the topologies defined exactly as the weak and strong operator topologies, but by using infinite families of vectors $(x_i)_{i\in\mathbb N}\subset H$ instead of finite families $(x_i)_{i=1,\ldots,N}\subset H$. 
\end{definition}

And up to you to tell me if you love such things or not, and I will be here listening, like Sigmund Freud. Anyway. With this convention, we have the following result:

\index{linear form}
\index{weak continuity}

\begin{proposition}
Given a vector space $E\subset B(H)$, and a linear form $f:E\to\mathbb C$, the following conditions are equivalent:
\begin{enumerate}
\item $f$ is ultraweakly continuous.

\item $f$ is ultrastrongly continuous.

\item $f(T)=\sum_{i=1}^\infty<Tx_i,y_i>$, for certain vectors $x_i,y_i\in H$.
\end{enumerate}
\end{proposition}

\begin{proof}
This is similar to the proof of Proposition 9.4, as follows:

\medskip

$(1)\implies(2)$ Since the ultraweak operator topology is weaker than the ultrastrong operator topology, ultraweakly continuous implies ultrastrongly continuous.

\medskip

$(2)\implies(3)$ Assume that $f:E\to\mathbb C$ is ultrastrongly continuous. By continuity we can find vectors $x_i\in H$ and a number $\varepsilon>0$ such that:
$$\sum_{i=1}^\infty||Tx_i||^2<\varepsilon^2\implies |f(T)|<1$$

It follows from this that we have the following estimate:
$$|f(T)|<\frac{1}{\varepsilon}\sqrt{\sum_{i=1}^\infty||Tx_i||^2}$$

Consider now the direct sum $H^{\oplus\infty}$, and inside it, the following vector:
$$x=(x_i)\in H^{\oplus\infty}$$

Consider also the following linear space, written in tensor product notation:
$$K=\overline{(E\otimes 1)x}\subset H^{\oplus\infty}$$

We can define a linear form $f':K\to\mathbb C$ by the following formula, and continuity:
$$f'((Tx_i)_i)=f(T)$$

We conclude that there exists a vector $y\in K$ such that:
$$f'\big((T\otimes1)y\big)=<(T\otimes 1)x,y>$$

But in terms of the original linear form $f:E\to\mathbb C$, this means that we have:
$$f(T)=\sum_{i=1}^\infty<Tx_i,y_i>$$

$(3)\implies(1)$ This is indeed clear from definitions.
\end{proof}

As a consequence of the above result, we have:

\index{complete additivity}
\index{vector state}

\begin{theorem}
Given a von Neumann algebra $A\subset B(H)$, and a positive linear form $f:A\to\mathbb C$, the following are equivalent:
\begin{enumerate}
\item $f$ is normal, in the sense that $f\left(\sup_ix_i\right)=\sup_i f(x_i)$, for any increasing sequence of positive elements $x_i\in A$.

\item $f$ is completely additive, in the sense that $f\left(\bigvee_ip_i\right)=\sum_if(p_i)$, for any family of pairwise orthogonal projections $p_i\in A$.

\item $f$ is ultraweakly continuous, or equivalently, $f$ is a vector state, $f=<Tx,x>$, when suitably extending it to the space $H\otimes l^2(\mathbb N)$.
\end{enumerate}
\end{theorem}

\begin{proof}
This is something very standard, as follows:

\medskip

$(1)\implies(2)$ Given a family of pairwise orthogonal projections $\{p_i\}$, we can consider the following increasing sequence of positive elements:
$$x_n=\sum_{i=1}^np_i$$

By using now the formula in (1) for these elements we obtain, as desired:
\begin{eqnarray*}
f\left(\bigvee_ip_i\right)
&=&f\left(\sup_nx_n\right)\\
&=&\sup_nf(x_n)\\
&=&\sup_n\sum_{i=1}^nf(p_i)\\
&=&\sum_if(p_i)
\end{eqnarray*}

$(2)\implies(3)$ This is something more technical, that we will prove in several steps. Let us fix a projection $q\in A$, and consider a vector $\xi\in Im(q)$ such that:
$$<q\xi,\xi>>1$$

Our claim is that there exists a projection $p\leq q$ such that, for any $x\in A$:
$$f(pxp)\leq<pxp\xi,\xi>$$

Indeed, in order to prove this, let us pick, by using the Zorn lemma, a maximal family of pairwise orthogonal projections $\{p_i\}\subset A$ such that, for any $i$, we have:
$$f(p_i)\geq<p_i\xi,\xi>$$  

By using our complete additivity assumption, we have then:
\begin{eqnarray*}
f\left(\bigvee_ip_i\right)
&=&\sum_if(p_i)\\
&\geq&\sum_i<p_i\xi,\xi>\\
&=&\left<\left(\bigvee_ip_i\right)\xi,\xi\right>
\end{eqnarray*}

Now consider the following projection, which is nonzero:
$$p=q-\bigvee_ip_i$$

By maximality of the family $\{p_i\}$, for any nonzero projection $r\leq p$, we have:
$$f(r)<<r\xi,\xi>$$

We therefore obtain the following estimate, valid for any $x\in A_+$, as desired:
$$f(pxp)\leq<pxp\xi,\xi>$$

Now by Cauchy-Schwarz we obtain that for any $x\in A$, $||x||\leq1$, we have:
\begin{eqnarray*}
|f(xp)|^2
&\leq&f(px^*xp)f(1)\\
&\leq&<px^*xp\xi,\xi>\\
&=&||xp\xi||^2
\end{eqnarray*}

Thus the following linear form is strongly continuous on the unit ball of $A$:
$$x\to f(px)$$

In order to finish now, once again by using the Zorn lemma, let us pick a maximal family of pairwise orthogonal projections $\{p_i\}\subset A$ such that $x\to f(p_ix)$ is strongly continuous on the unit ball of $A$, for any $i$. By maximality we have then:
$$\sum_if(p_i)
=f\left(\sum_ip_i\right)
=f(1)
=1$$

Now given $\varepsilon>0$, let us choose a finite subset of our index set, $F\subset I$, such that for all the finite subsets $F\subset J\subset I$, we have an inequality as follows:
$$1-f\left(\sum_{j\in J}p_j\right)\leq\varepsilon$$

By Cauchy-Schwarz we have then, for any $x\in A$, $||x||=1$, the following estimate:
\begin{eqnarray*}
\left|f\left(x\left(1-\sum_{j\in J}p_j\right)\right)\right|^2
&\leq&f\left(1-\sum_{j\in J}p_j\right)f(xx^*)\\
&\leq&\varepsilon
\end{eqnarray*}

We conclude from this that we have the following estimate:
$$\left|\left|f-f\left(.\left(1-\sum_{j\in J}p_j\right)\right)\right|\right|\leq\sqrt{\varepsilon}$$

Thus we obtain $f\in A_*$, as desired.

\medskip

$(3)\implies(1)$ This is something trivial, coming from definitions.
\end{proof}

We can now go back to our original question, and we have:

\index{algebraic isomorphism}
\index{spatial isomorphism}

\begin{theorem}
Given two von Neumann algebras $A\subset B(H)$ and $B\subset B(K)$, acting on possibly different Hilbert spaces $H,K$, any algebraic isomorphism 
$$\Phi:A\simeq B$$
is spatial up to amplification, in the sense that we have a formula as follows,
$$\Phi(T)\otimes1=U(T\otimes1)U^*$$
for a certain Hilbert space $L$, and a certain unitary $U:H\otimes L\to K\otimes L$. 
\end{theorem}

\begin{proof}
This is something standard, coming from Theorem 9.19, as follows:

\medskip

(1) As a first observation, assuming that a positive unital linear form $f:A\to\mathbb C$ is a vector state, given by a certain vector $x\in H$, then by Theorem 9.19 the linear form $f\Phi^{-1}$ is also a vector state, say given by a vector $y\in K$. 

\medskip

(2) We conclude from this that we have a unitary as follows, intertwining the corresponding actions of the von Neumann algebras $A$ and $B$:
$$U_x:\overline{Ax}\to\overline{By}$$

Now by making the above vector $x\in H$ vary, and performing a direct sum, we obtain with $L=l^2(\mathbb N)$ an isometry as in the statement, namely:
$$U:H\otimes L\to K\otimes L$$

Our construction shows that $U$ intertwines indeed the actions of the von Neumann algebras $A$ and $B$, and what is left to do is to study the unitarity of $U$.

\medskip

(3) We will prove now that, up to a suitable replacement, the above operator $U$ can be taken to be unitary, still intertwining the actions of the von Neumann algebras $A$ and $B$. For this purpose, consider the action of von Neumann algebra $A$ on the direct sum Hilbert space $(H\otimes L)\oplus(K\otimes L)$ given by the following matrices:
$$x'=\begin{pmatrix}
x\otimes1&0\\
0&\Phi(x)\otimes1
\end{pmatrix}$$

Since $U$ intertwines the actions of the von Neumann algebras $A$ and $B$, in terms of $2\times 2$ matrices, we are led to the following conclusion:
$$\begin{pmatrix}
0&0\\
U&0
\end{pmatrix}\in A'$$

Thus, the following happens inside the von Neumann algebra $A'$:
$$\begin{pmatrix}
1&0\\
0&0
\end{pmatrix}\preceq 
\begin{pmatrix}
0&0\\
0&1
\end{pmatrix}$$

On the other hand, the same reasoning applied to the isomorphism $\Phi^{-1}$ shows that we have as well, once again inside the von Neumann algebra $A'$:
$$\begin{pmatrix}
0&0\\
0&1
\end{pmatrix}\preceq 
\begin{pmatrix}
1&0\\
0&0
\end{pmatrix}$$

(4) We are now in position to finish. By combining the above two conclusions, we obtain an equivalence of projections inside $A'$, as follows:
$$\begin{pmatrix}
1&0\\
0&0
\end{pmatrix}\simeq 
\begin{pmatrix}
0&0\\
0&1
\end{pmatrix}$$

Now pick a partial isometry implementing this equivalence. This partial isometry must be of the following form, with $U'$ being now a unitary:
$$V=\begin{pmatrix}
0&0\\
U'&0
\end{pmatrix}$$

Thus, we have a unitary as follows, which intertwines the actions of $A$ and $B$:
$$U':H\otimes L\to K\otimes L$$

But this is the unitary we were looking for, and we are done.
\end{proof}

The above result is something quite fundamental, allowing us to talk about von Neumann algebras $A$ as abstract objects, without reference to the exact Hilbert space $H$ where the elements $a\in A$ live as operators $a\in B(H)$, and with this being of course possible modulo some functional analysis knowledge. We will heavily use this point of view in chapter 10 below, and then in chapters 13-16, when talking about ${\rm II}_1$ factors.

\section*{9d. Predual theory}

We have seen so far, as a consequence of the Kaplansky density theorem, that an operator algebra $A\subset B(H)$ is a von Neumann algebra precisely when its unit ball is weakly compact. This is certainly useful, but there are many other possible characterizations of the von Neumann algebras, as operator algebras, which are useful as well.

\bigskip

To be more precise, going ahead now with more abstract functional analysis, that we will be using in what follows, on several occasions, let us formulate:

\begin{definition}
Given a von Neumann algebra $A\subset B(H)$, we set
$$A_*=\Big\{f:A\to\mathbb C,\ {\rm weakly\ continuous}\Big\}$$
regarded as a linear subspace, $A_*\subset A^*$, of the usual dual, given by
$$A^*=\Big\{f:A\to\mathbb C,\ {\rm norm\ continuous}\Big\}$$
and we call this space $A_*$ predual of our von Neumann algebra $A$.
\end{definition}

Our first goal will be that of proving that we have the following duality formula, between the linear space $A_*$ constructed above, and the algebra $A$ itself:
$$A=(A_*)^*$$

In order to do so, let us first discuss the case of the full operator algebra $A=B(H)$ itself. This is actually the key case, with the extension to the arbitrary von Neumann algebras $A\subset B(H)$ being something coming afterwards, quite straightforward. 

\bigskip

We will need some standard operator theory, developed in chapter 4. First, we have the following result, regarding the trace class operators, established there:

\index{trace class operator}

\begin{theorem}
The space of trace class operators, which appears as an intermediate space between the finite rank operators and the compact operators,
$$F(H)\subset B_1(H)\subset K(H)$$
is a two-sided $*$-ideal of $K(H)$. The following is a Banach space norm on $B_1(H)$, 
$$||T||_1=Tr|T|$$
satisfying $||T||\leq||T||_1$, and for $T\in B_1(H)$ and $S\in B(H)$ we have:
$$||ST||_1\leq||S||\cdot||T||_1$$
Also, the subspace $F(H)$ is dense inside $B_1(H)$, with respect to this norm.
\end{theorem}

\begin{proof}
This is indeed something standard, explained in chapter 4.
\end{proof}

We will need as well the following result, regarding this time the Hilbert-Schmidt operators, which is also from chapter 4:

\index{Hilbert-Schmidt operator}

\begin{theorem}
The space of Hilbert-Schmidt operators, which appears as an intermediate space between the trace class operators and the compact operators,
$$F(H)\subset B_1(H)\subset B_2(H)\subset K(H)$$
is a two-sided $*$-ideal of $K(H)$. In terms of the singular values $(\lambda_n)$, the Hilbert-Schmidt operators are characterized by the following formula:
$$\sum_n\lambda_n^2<\infty$$
Also, the following formula, taking as input two Hilbert-Schmidt operators,
$$<S,T>=Tr(ST^*)$$
defines a scalar product of $B_2(H)$, making it a Hilbert space.
\end{theorem}

\begin{proof}
As before, this is something standard, explained in chapter 4.
\end{proof}

We will need as well the following technical result, also from chapter 4:

\begin{theorem}
We have the following commutation formula,
$$Tr(ST)=Tr(TS)$$
valied for any two Hilbert-Schmidt operators $S,T\in B_2(H)$.
\end{theorem}

\begin{proof}
As before, this is something standard, explained in chapter 4.
\end{proof}

With the above ingredients in hand, and getting back now to von Neumann algebras, and to our predual questions raised before, we first have the following result:

\index{predual}
\index{trace class operator}

\begin{theorem}
The linear space $B(H)_*\subset B(H)^*$ consisting of the linear forms $f:B(H)\to\mathbb C$ which are weakly continuous is given by
$$B(H)_*=\left\{T\to Tr(ST)\Big| S\in B_1(H)\right\}$$
and we have the following duality formula 
$$B(H)=(B(H)_*)^*$$
as a duality in the usual Banach space sense.
\end{theorem}

\begin{proof}
There are several things to be proved, the idea being as follows:

\medskip

(1) First of all, any linear form of type $T\to Tr(ST)$, with $S$ being trace class, is weakly continuous. Thus, if we denote by $B(H)_\circ$ the subspace of $B(H)$ in the statement, consisting of such linear forms, we have an inclusion as follows:
$$B(H)_\circ\subset B(H)_*$$

(2) In order to prove now the reverse inclusion, consider an arbitrary weakly continuous linear form $f\in B(H)_*$. We can then find vectors $(x_i)$ and $(y_i)$ such that:
$$f(T)=\sum_i<Tx_i,y_i>$$

Let us consider now the following operators, going by definition from the Hilbert space $l^2(\mathbb N)$ to our Hilbert space $H$, and which are both Hilbert-Schmidt:
$$Q:e_i\to x_i\quad,\quad 
R:e_i\to y_i$$

In terms of these operators, our linear form can be written as follows:
$$f(T)=Tr(R^*TQ)$$

On the other hand, by using the formula in Theorem 9.24 we obtain:
$$Tr(R^*TQ)=Tr(TQR^*)$$

Thus, with $S=QR^*$, which is trace class, we have the following formula:
$$f(T)=Tr(TS)$$

Thus, we have proved that we have an inclusion as follows:
$$B(H)_*\subset B(H)_\circ$$

(3) Summing up, from (1) and (2) we conclude that we have an equality as follows, which proves the first assertion in the statement:
$$B(H)_*=B(H)_\circ$$

(4) It remains to prove that $B(H)$ is indeed the dual of $B(H)_*$. For this purpose, we use the above identification, which ultimately identifies $B(H)_*$ with the space of trace class operators $B_1(H)$. So, assume that we have a linear form, as follows:
$$f:B_1(H)\to\mathbb C$$

It is then routine to show that $f$ must come from evaluation on a certain operator $T\in B(H)$, and this leads to the conclusion that $B(H)$ is indeed the dual of $B(H)_*$.
\end{proof}

More generally now, for the arbitrary von Neumann algebras $A\subset B(H)$, we have:

\begin{theorem}
Given a von Neumann algebra $A\subset B(H)$, if we set
$$A_*=\Big\{f:A\to\mathbb C,\ {\rm weakly\ continuous}\Big\}$$
regarded as a linear subspace, $A_*\subset A^*$, of the usual dual, given by:
$$A^*=\Big\{f:A\to\mathbb C,\ {\rm norm\ continuous}\Big\}$$
then we have the duality formula $A=(A_*)^*$, in the usual Banach space sense.
\end{theorem}

\begin{proof}
This can be proved in several steps, as follows:

\medskip

(1) First of all, we know from the above that the result holds for the von Neumann algebra $A=B(H)$ itself, in the sense that we have:
$$B(H)=(B(H)_*)^*$$

(2) The point now is that for any von Neumann subalgebra $A\subset B(H)$, or more generally for any weakly closed linear subspace $A\subset B(H)$, we have an equality as follows, coming as a consequence of the Hahn-Banach theorem:
$$A=A^{\perp\perp}$$

(3) Thus, modulo some standard algebra, and some standard identifications for quotient spaces and their duals, we are led to the conclusion in the statement.
\end{proof}

In fact, we have the following result, due to Sakai:

\index{predual}
\index{Sakai theorem}

\begin{theorem}
The von Neumann algebras are exactly the $C^*$-algebras which have a predual, in the above sense.
\end{theorem}

\begin{proof}
This is a variation of the above, which caps the above series of results, and closes any further discussions, and for details here, we refer to Sakai's book \cite{sak}.
\end{proof}

There are many other things that can be said, of purely abstract nature, on the von Neumann algebras. We will be back to this, from time to time, in what follows.

\section*{9e. Exercises} 

Things have been quite tricky in this chapter, with a number of detours, and by avoiding some difficulties, and as a first exercise, which is quite difficult, we have:

\begin{exercise}
Look up and learn von Neumann's reduction theory, stating that given a von Neumann algebra $A\subset B(H)$, if we write its center as 
$$Z(A)=L^\infty(X)$$
then we have a decomposition as follows, with the fibers $A_x$ having trivial center,
$$A=\int_XA_x\,dx$$
and then write down a brief account of what you learned.
\end{exercise} 

This is something very fundamental and instructive, and is actually fundamental to the point that we should have normally talked about this in the present chapter 9. However, in practice, this is not really doable, things very technical, and although foundational, this remains something quite advanced. We will be actually back to this later in this book, but only under the assumption that the algebra has a trace, $tr:A\to\mathbb C$, which simplifies a number of things, and rather with explanations, instead of a complete proof.

\begin{exercise}
Given two projections $P,Q\in B(H)$, with $H$ being infinite dimensional, find an elementary proof for the fact that we have, for any $x\in H$,
$$(PQ)^nx\to (P\wedge Q)x$$
but the operators $(PQ)^n$ do not necessarily converge in norm.
\end{exercise}

This is something that we proved in the above, but the problem now is that of finding an elementary proof of this fact, by using whatever tricks of your choice.

\begin{exercise}
Given a commutative von Neumann algebra, written as
$$A=L^\infty(X)$$
with $X$ being a measured space, write, by using the Gelfand theorem,
$$A=C(\widehat{X})$$
with $\widehat{X}$ being a compact space, and understand the correspondence $X\to\widehat{X}$. 
\end{exercise} 

As a bonus exercise here, try understanding as well what happens for the non-unital commutative $C^*$-algebras $A\subset B(H)$, and their weak closures $A''\subset B(H)$.

\chapter{Finite factors}

\section*{10a. Finite factors}

In this chapter we go for the real thing, namely the study of the ${\rm II}_1$ factors, following Murray and von Neumann \cite{mv1}, \cite{mv2}, \cite{mv3}, \cite{vn1}, \cite{vn2}, which is the basis for everything more advanced, in relation with operator algebras. We will only present the very basic theory of the ${\rm II}_1$ factors, with the idea in mind of using them later for doing subfactors a la Jones. We will mainly follow the simplified approach from Jones' book \cite{jo6}, with sometimes a look into Blackadar \cite{bla}, both books that we recommend for more.

\bigskip

Let us first talk about general factors. There are several possible ways of introducing them, and dividing them into several classes, for further study. In what concerns us, we will use a rather intuitive approach. The general idea, which is quite natural, is that among the von Neumann algebras $A\subset B(H)$, of particular interest are the ``free'' ones, having trivial center, $Z(A)=\mathbb C$. These algebras are called factors:

\begin{definition}
A factor is a von Neumann algebra $A\subset B(H)$ whose center
$$Z(A)=A\cap A'$$
which is a commutative von Neumann algebra, reduces to the scalars, $Z(A)=\mathbb C$.
\end{definition}

This notion is in fact something that we already met in the above, in the context of various comments or exercises, and time now to clarify all this. The idea is that there are two main motivations for the study of factors, with each of them being more than enough, as to serve as a strong motivation. First, at the intuitive level, we have:

\begin{principle}[Freeness]
The following happen:
\begin{enumerate}
\item The condition $Z(A)=\mathbb C$ defining the factors is, obviously, opposite to the condition $Z(A)=A$ defining the commutative von Neumann algebras. 

\item Therefore, the factors are the von Neumann algebras which are ``free'', meaning as far as possible from the commutative ones.

\item Equivalently, with $A=L^\infty(X)$, the quantum spaces $X$ coming from factors are those which are ``free'', meaning as far as possible from the classical spaces.
\end{enumerate}
\end{principle}

So, this was for our first principle, which is something reasonable, intuitive, and self-explanatory, and which can surely serve as a strong motivation for the study of factors. In fact, all that has being said above comes straight from the structure theorem for the commutative von Neumann algebras, $A=L^\infty(X)$, with $X$ being a measured space, that we know since chapter 5, and the above principle is just a corollary of that theorem.

\bigskip

At a more advanced level, another motivation for the study of factors, which among others justifies the name ``factors'' for them, comes from the reduction theory of von Neumann \cite{vn3}, which is something non-trivial, that can be summarized as follows:

\begin{principle}[Reduction theory]
Given a von Neumann algebra $A\subset B(H)$, if we write its center $Z(A)\subset A$, which is a commutative von Neumann algebra, as 
$$Z(A)=L^\infty(X)$$
with $X$ being a measured space, then the whole algebra decomposes as
$$A=\int_XA_x\,dx$$
with the fibers $A_x$ being factors, that is, satisfying $Z(A_x)=\mathbb C$.
\end{principle}

As a first comment, we have already seen an instance of such decomposition results in chapter 5, when talking about finite dimensional algebras. Indeed, such algebras decompose, in agreement with the above, as direct sums of matrix algebras, as follows:
$$A=\bigoplus_xM_{n_x}(\mathbb C)$$

In general, however, things are more complicated than this, and technically speaking, and as opposed to Principle 10.2, which was more of a triviality, Principle 10.3 is a tough theorem, due to von Neumann \cite{vn3}. More on this later, in chapter 11 below.

\bigskip

This was for the story, and let us close this philosophical discussion with:

\begin{conclusion}
Regardless of the approach and technical level, be that beginner or advanced, the von Neumann factors are the algebras that matter.
\end{conclusion}

Getting to work now, there are many things that can be said about factors. In order to get started, as a direct continuation of the work from chapter 9, for the general von Neumann algebras, let us first study their projections. We will see that many interesting things happen here, with everything coming from the following technical result:

\begin{proposition}
Given two projections $p,q\neq0$ in a factor $A$, we have
$$puq\neq0$$
for a certain unitary $u\in A$.
\end{proposition}

\begin{proof}
Assume by contradiction $puq=0$, for any unitary $u\in A$. This gives:
$$u^*puq=0$$

By using this for all the unitaries $u\in A$, we obtain the following formula:
$$\left(\bigvee_{u\in U_A}u^*pu\right)q=0$$

On the other hand, from $p\neq0$ we obtain, by factoriality of $A$:
$$\bigvee_{u\in U_A}u^*pu=1$$

Thus, our previous formula is in contradiction with $q\neq0$, as desired.
\end{proof}

Getteing back now to the order on projections from chapter 9, and to the whole von Neumann projection philosophy, in the case of factors things simplify, as follows:

\index{projections in factors}
\index{lattice of projections}

\begin{theorem}
Given two projections $p,q\in A$ in a factor, we have
$$p\preceq q\quad {\rm or}\quad q\preceq p$$
and so $\preceq$ is a total order on the equivalence classes of projections $p\in A$.
\end{theorem}

\begin{proof}
This basically follows from Proposition 10.5, and from the Zorn lemma, by using some standard functional analysis arguments. To be more precise:

\medskip

(1) Consider indeed the following set of partial isometries:
$$S=\left\{u\Big|uu^*\leq p,u^*u\leq q\right\}$$

We can then order this set $S$ by saying that we have $u\leq v$ when $u^*u\leq v^*v$, and when $u=v$ holds on the initial domain $u^*uH$ of $u$. With this convention made, the Zorn lemma applies, and provides us with a maximal element $u\in S$.

\medskip

(2) In the case where this maximal element $u\in S$ satisfies $uu^*=p$ or $u^*u=q$, we are led to one of the conditions $p\preceq q$ or $q\preceq p$ in the statement, and we are done.

\medskip

(3) So, assume that we are in the case left, $uu^*\neq p$ and $u^*u\neq q$. By Proposition 10.5 we obtain a unitary $v\neq0$ satisfying the following conditions:
$$vv^*\leq p-uu^*$$
$$v^*v\leq q-u^*u$$

But these conditions show that the element $u+v\in S$ is strictly bigger than $u\in S$, which is a contradiction, and we are done.
\end{proof}

Moving ahead now, as explained time and again throughout this book, for a variety of reasons, which can be elementary or advanced, and also mathematical or physical, we are mainly interested in the case where our algebras have traces:
$$tr:A\to\mathbb C$$

And in relation with the factors, by leaving aside the rather trivial case of the matrix algebras $A=M_N(\mathbb C)$, we are led in this way to the following key notion:

\begin{definition}
A ${\rm II}_1$ factor is a von Neumann algebra $A\subset B(H)$ which:
\begin{enumerate}
\item Is infinite dimensional, $\dim A=\infty$.

\item Has trivial center, $Z(A)=\mathbb C$.

\item Has a trace $tr:A\to\mathbb C$.
\end{enumerate}
\end{definition}

Here the order of the axioms is a bit random, with any of the possible $3!=6$ choices making sense, and corresponding to a slightly different vision on what the ${\rm II}_1$ factors truly are. With the above order, with (1) we are making it clear, right from the beginning, that we are not here for revolutionizing linear algebra. Then with (2) we adhere to Definition 10.1, and to what was said next about it, on freeness and reduction. And finally with (3) we adhere to the above principle, that von Neumann algebras must have traces.

\bigskip

More technically now, and leaving aside anything subjective, the above definition is motivated by the heavy classification work of Murray, von Neumann and Connes \cite{co1}, \cite{co2}, \cite{mv1}, \cite{mv2}, \cite{mv3}, \cite{vn1}, \cite{vn2}, \cite{vn3}, whose conclusion is more or less that everything in von Neumann algebras reduces, via some quite complicated procedures, we should mention that, to the study of the ${\rm II}_1$ factors. With the mantra here being as follows:

\begin{fact}
The ${\rm II}_1$ factors are the building blocks of the whole von Neumann algebra theory.
\end{fact}

To be more precise, this statement, that we will get to understand later, is something widely agreed upon, at least among operator algebra experts who are familiar with von Neumann algebras, and with this agreement being something great. What remains controversial, however, is how to start playing with these Lego bricks that we have:

\bigskip

(1) A first option is that of adding the matrix algebras $M_N(\mathbb C)$, not to be forgotten, and then stacking together such Lego bricks. According to the von Neumann reduction theory, this leads to the von Neumann algebras having traces, $tr:A\to\mathbb C$.

\bigskip

(2) A second option, perhaps even more playful, is that of taking crossed products of such Lego bricks by their automorphisms scaling the trace, or performing more general constructions inspired by advanced ergodic theory. This leads to general factors.

\bigskip

(3) And the third option is that of being a bad kid, or perhaps some kind of nerd, engineer in the becoming, and picking such a Lego brick, or a handful of them, and breaking them, see what's inside. Good option too, and more on this later.

\bigskip

Getting to work now, in practice, and forgetting about reduction theory, which raises the possibility of decomposing any tracial von Neumann algebra into factors, in order to obtain explicit examples of ${\rm II}_!$ factors, it is not even clear that such beasts exist. Fortunately the group von Neumann algebras are there, and we have the following result, which provides us with some examples of ${\rm II}_1$ factors, to start with:

\index{group algebra}
\index{center of algebra}
\index{conjugacy classes}
\index{ICC property}
\index{factoriality}

\begin{theorem}
The center of a group von Neumann algebra $L(\Gamma)$ is
$$Z(L(\Gamma))=\left\{\sum_g\lambda_gg\Big|\lambda_{gh}=\lambda_{hg}\right\}''$$
and if $\Gamma\neq\{1\}$ has infinite conjugacy classes, in the sense that 
$$\Big|\{ghg^{-1}|g\in G\}\Big|=\infty\quad,\quad\forall h\neq1$$
with this being called ICC property, the algebra $L(\Gamma)$ is a ${\rm II}_1$ factor.
\end{theorem}

\begin{proof}
There are two assertions here, the idea being as follows:

\medskip

(1) Consider a linear combination of group elements, which is in the weak closure of $\mathbb C[\Gamma]$, and so defines an element of the group von Neumann algebra $L(\Gamma)$:
$$a=\sum_g\lambda_gg$$

By linearity, this element $a\in L(\Gamma)$ belongs to the center of $L(\Gamma)$ precisely when it commutes with all the group elements $h\in\Gamma$, and this gives:
\begin{eqnarray*}
a\in Z(A)
&\iff&ah=ha\\
&\iff&\sum_g\lambda_ggh=\sum_g\lambda_ghg\\
&\iff&\sum_k\lambda_{kh^{-1}}k=\sum_k\lambda_{h^{-1}k}k\\
&\iff&\lambda_{kh^{-1}}=\lambda_{h^{-1}k}
\end{eqnarray*}

Thus, we obtain the formula for $Z(L(\Gamma))$ in the statement. 

\medskip

(2) We have to examine the 3 conditions defining the ${\rm II}_1$ factors. We already know from chapter 7 that the group algebra $L(G)$ has a trace, given by:
$$tr(g)=\delta_{g,1}$$

Regarding now the center, the condition $\lambda_{gh}=\lambda_{hg}$ that we found is equivalent to the fact that $g\to\lambda_g$ is constant on the conjugacy classes, and we obtain:
$$Z(L(\Gamma))=\mathbb C\iff \Gamma={\rm ICC}$$

Finally, assuming that this ICC condition is satisfied, with $\Gamma\neq\{1\}$, then our group $\Gamma$ is infinite, and so the algebra $L(\Gamma)$ is infinite dimensional, as desired.
\end{proof}

In order to look now for more examples of ${\rm II}_1$ factors, an idea would be that of attempting to decompose into factors the group von Neumann algebras $L(\Gamma)$, but this is something difficult, and in fact we won't really exit the group world in this way. Difficult as well is to investigate the factoriality of the von Neumann algebras of discrete quantum groups $L(\Gamma)$, because the basic computations from the proof of Theorem 10.9 won't extend to this setting, where the group elements $g\in\Gamma$ become corepresentations $g\in M_N(L(\Gamma))$. Despite years of efforts, it is presently not known at all what the ``quantum ICC'' condition should mean, and the problem comes from this. But more on this later.

\bigskip

In short, we have to stop here the construction of examples, and Theorem 10.9 will be what we have, at least for the moment. With this being actually not a big issue, the group factors $L(\Gamma)$ being known to be quite close to the generic ${\rm II}_1$ factors.

\section*{10b. Basic results}

Getting away now from the above difficulties, let us go back to the abstract ${\rm II}_1$ factors, as axiomatized in Definition 10.7. In order to investigate them, the idea will be that from chapter 9, namely looking at the projections, and their equivalence classes. 

\bigskip

In the case of the ${\rm II}_1$ factors, as a first interesting remark, the presence of the trace trivializes the proof of the main result that we have about projections, as follows:

\index{order of projections}
\index{lattice of projections}
\index{projections in factors}

\begin{theorem}
Given two projections $p,q\in A$ in a ${\rm II}_1$ factor we have, trivially
$$p\preceq q\quad {\rm or}\quad q\preceq p$$
and so $\preceq$ is a total order on the equivalence classes of projections $p\in A$.
\end{theorem}

\begin{proof}
This is something that we already know, from Theorem 10.6, and which actually holds for any factor, with the non-trivial part being the following implication:
$$p\preceq q,\ q\preceq p\implies p\simeq q$$

But this implication is clear in the present ${\rm II}_1$ factor setting, by using the trace.
\end{proof}

The above theorem and its proof, which are remarkable, are the first in a series of mysteries, in what concerns the special case of the ${\rm II}_1$ factors. More such mysteries to follow. In order to study now the trace of the ${\rm II}_1$ factors, we will need:

\index{closed ideal}

\begin{proposition}
Given a weakly closed left ideal $I\subset A$ in a von Neumann algebra, there exists a unique projection $p\in A$ such that:
$$I=Ap$$
Moreover, if $I\subset A$ is assumed to be a two-sided ideal, then $p\in Z(A)$.
\end{proposition}

\begin{proof}
We have several things to be proved, the idea being as follows:

\medskip

(1) Given an ideal $I\subset A$ as in the statement, consider the following intersection:
$$I\cap I^*\subset A$$

This is a weakly closed non-unital $*$-subalgebra of $A$, so if we denote by $p\in A$ its largest projection, or unit, then we have an inclusion $Ap\subset I$.

\medskip

(2) Conversely now, let us pick $x\in I$. By polar decomposition we can write $x=u|x|$, and we have the following implications, which prove the reverse inclusion $I\subset Ap$:
\begin{eqnarray*}
x\in I
&\implies&|x|=u^*x\in I\\
&\implies&|x|\in I\cap I^*\\
&\implies&|x|p=|x|\\
&\implies&x=u|x|=u|x|p\in Ap
\end{eqnarray*}

(3) The uniqueness assertion is clear from the comparison theorem for projections.

\medskip

(4) Regarding now the last assertion, assume that $I\subset A$ is a two-sided weakly closed ideal. Then for any unitary $u\in A$ we have:
\begin{eqnarray*}
I=uIu^*
&\implies&uIu^*=Ap\\
&\implies&I=Aupu^*
\end{eqnarray*}

Thus by uniqueness we obtain $upu^*=p$, and so $p\in Z(A)$, as desired.
\end{proof}

As a first main result now regarding the ${\rm II}_1$ factors, following the paper of Murray and von Neumann \cite{mv3}, which by the way is a must-read, we have:

\begin{theorem}
Given a ${\rm II}_1$ factor $A$, any weakly continuous positive trace
$$tr:A\to\mathbb C$$
is automatically faithful.
\end{theorem}

\begin{proof}
Consider the null space of the trace, which is by definition:
$$I=\left\{x\in A\Big|tr(x^*x)=0\right\}$$

We have the following inequality, which shows that $I$ is a left ideal:
$$x^*a^*ax\leq||a||^2x^*x$$

Now by using the trace condition $tr(ab)=tr(ba)$, we conclude that $I$ is a two-sided ideal. Also, the Cauchy-Schwarz inequality gives:
$$tr(x^*x)=0\iff tr(xy)=0,\forall y\in A$$

We conclude from this that $I$ is an intersection of kernels of weakly closed functionals, which are weakly closed, and so it is weakly closed. Thus the last assertion in Proposition 10.11 applies, and produces a projection $p\in Z(A)$ such that:
$$I=Ap$$

Now since $A$ was assumed to be a factor, we have $Z(A)=\mathbb C$. Thus $p=0$, and so the null ideal of the trace is $I=\{0\}$, and so our trace $tr$ is faithful, as desired.
\end{proof}

Our goal now will be that of proving that the trace on a ${\rm II}_1$ factor is unique, and takes on projections any value in $[0,1]$. Let us start with a technical result, as follows:

\begin{proposition}
Given a ${\rm II}_1$ factor $A$, the traces of the projections 
$$tr(p)\in[0,1]$$
can take arbitrarily small values.
\end{proposition}

\begin{proof}
Consider the set formed by all values of the trace on the projections:
$$S=\left\{tr(p)\Big|p^2=p=p^*\in A\right\}$$

We want to prove that the following number equals $0$:
$$c=\inf(S-\{0\})$$

In order to do so, assume by contradiction $c>0$, pick $\varepsilon>0$ small, and pick a projection $p\in A$ such that the following condition is satisfied:
$$tr(p)<c+\varepsilon$$

Since we are in a ${\rm II}_1$ factor, this projection $p\in A$ cannot be minimal, and so we can find another projection $q\in A$ satisfying $q<p$. Now observe that we have:
\begin{eqnarray*}
tr(p-q)
&=&tr(p)-tr(q)\\
&\leq&tr(p)-c\\
&\leq&\varepsilon
\end{eqnarray*}

Thus with $\varepsilon<c$ we obtain a contradiction, and so $c=0$, as desired.
\end{proof}

In order to prove our next main result, we will need as well:

\begin{proposition}
Given a ${\rm II}_1$ factor $A$ on a Hilbert space $H$ and a projection $p\in A$, the von Neumann algebra $pAp$ is a ${\rm II}_1$ factor on the Hilbert space $pH$.
\end{proposition}

\begin{proof}
We have to prove that the von Neumann algebra $pAp$ has a trace, and is infinite dimensional, and these two properties can be proved as follows:

\medskip

(1) In what regards the trace, we know that the trace $tr:A\to\mathbb C$ restricts to a trace $tr:pAp\to\mathbb C$, which must be nonzero, as desired.

\medskip

(2) In what regards the infinite dimensionality, this follows from the fact that a minimal projection in $pAp$ would be minimal in $A$, which is impossible.
\end{proof}

Still following the fundamental paper of Murray and von Neumann \cite{mv3}, we can now formulate a second main result regarding the ${\rm II}_1$ factors, as follows:

\index{traces of projections}
\index{projections in factors}
\index{continuous dimension}

\begin{theorem}
Given a ${\rm II}_1$ factor $A$, the traces of projections 
$$tr(p)\in[0,1]$$
can take any values in $[0,1]$.
\end{theorem}

\begin{proof}
Given a number $c\in[0,1]$, consider the following set:
$$S=\left\{p^2=p=p^*\in A\Big|tr(p)\leq c\right\}$$

This set satisfies the assumptions of the Zorn lemma, and so by this lemma we can find a maximal element $p\in S$. Assume by contradiction that we have:
$$tr(p)<c$$

The point now is that by using Proposition 10.13 and Proposition 10.14, we can slightly enlarge the trace of $p$, and we obtain a contradiction, as desired.
\end{proof}

As a third and last main result regarding the ${\rm II}_1$ factors, also from \cite{mv3}, we have:

\index{uniqueness of trace}

\begin{theorem}
The trace of a ${\rm II}_1$ factor
$$tr:A\to\mathbb C$$
is unique.
\end{theorem}

\begin{proof}
This can be proved in many ways, a standard one being that of proving that any two traces agree on the projections, as a consequence of the above results:

\medskip

(1) Assume indeed that we have a second trace $tr':A\to\mathbb C$. Since $A$ is generated by its projections, it is enough to show that we have $tr=tr'$ on projections.

\medskip

(2) As a first observation, since traces on matrix algebras are unique, we obtain that we have $tr=tr'$ on the projections $p\in A$ having rational trace, $tr(p)\in\mathbb Q$.

\medskip

(3) So, let us pick $p\in A$ having non-rational trace, $tr(p)\notin\mathbb Q$, and prove that we have $tr(p)=tr'(p)$. The idea will be that of using the result for the projections having rational traces, applied to an infinite direct sum of projections, converging to $p$.

\medskip

(4) To be more precise, assume that we have constructed our sequence $p_i\to p$ up to order $n\in\mathbb N$, and let us try to construct $p_{n+1}$. The idea is to use the following algebra:
$$A_n=(p-p_n)A(p-p_n)$$

(5) Indeed this algebra is a ${\rm II}_1$ factor, and we can choose inside it a projection $p_{n+1}$ satisfying $p_n\leq p_{n+1}\leq p$, such that $tr=tr'$ on it, and such that:
$$tr(p-p_{n+1})\leq\frac{1}{2}\cdot tr(p-p_n)$$

(6) According to our choices for these projections $p_n$, we have:
$$p=\bigvee_{n=1}^\infty p_n$$

Thus when evaluating $tr,tr'$ on $p$ we obtain the same result, as desired.
\end{proof}

In what regards illustrations for all this, as examples of ${\rm II}_1$ factors we have so far the group von Neumann algebras $L(\Gamma)$, with $\Gamma$ being an ICC group. In certain cases, it is possible to say more about all the above, and in particular about the projections, for instance with quite explicit procedures for constructing projections $p\in L(\Gamma)$ having an arbitrary prescribed trace $x\in[0,1]$. We will be back to this later, when discussing more in detail the group von Neumann algebras $L(\Gamma)$, and their generalizations.

\bigskip

Back to theory, we have seen that the ${\rm II}_1$ factors are very interesting objects, naturally lying above the matrix algebras $M_N(\mathbb C)$, which are type I factors. From this perspective, a ${\rm II}_1$ factor $A\subset B(H)$ is not really in need of the ambient Hilbert space $H$, and the question of ``representing'' it appears. We will discuss this question, in two steps:

\medskip

\begin{enumerate}
\item A first question is that of understanding the possible embeddings $A\subset B(H)$, with $H$ being a Hilbert space. The main result here will be the construction of a numeric invariant $\dim_AH$, called coupling constant.

\medskip

\item A second question is that of understanding the possible embeddings $A\subset B$, with $B$ being another ${\rm II}_1$ factor. By using the coupling constant for both $A,B$ we will construct a numeric invariant $[B:A]$, called index.
\end{enumerate}

\medskip

We will discuss now (1), and leave (2) for later, towards the end of this chapter. In order to get started, let us formulate the following definition:

\index{GNS construction}
\index{standard form}

\begin{definition}
Given a von Neumann algebra $A$ with a trace $tr:A\to\mathbb C$, the emdedding
$$A\subset B(L^2(A))$$
obtained by GNS construction is called standard form of $A$.
\end{definition}

Here we use the GNS construction, explained in chapter 7. As the name indicates, the standard representation is something ``standard'', to be compared with any other representation $A\subset B(H)$, in order to understand this latter representation.

\bigskip

As already seen in chapter 7, the GNS construction has a number of unique features, that can be exploited. In the present setting, the main result is as follows:

\index{J map}

\begin{theorem}
In the context of the standard representation we have
$$A'=JAJ$$
with $J:L^2(A)\to L^2(A)$ being the antilinear map given by $T\to T^*$.
\end{theorem} 

\begin{proof}
Observe first that any $T\in A$ can be regarded as a vector $T\in L^2(A)$, to which we can associate, in an antilinear way, the vector $T^*\in L^2(A)$. Thus we have indeed an antilinear map $J$ as in the statement. In terms of the standard cyclic and separating vector $\Omega$ for the GNS representation, the formula of this formula $J$ is:
$$J(x\Omega)=x^*\Omega$$

(1) Our first claim is that we have the following formula:
$$<J\xi,J\eta>=<\xi,\eta>$$

Indeed, with $\xi=x\Omega$ and $\eta=y\Omega$, we have the following computation:
\begin{eqnarray*}
<J\xi,J\eta>
&=&<yx^*\Omega,\Omega>\\
&=&tr(yx^*)\\
&=&<\xi,\eta>
\end{eqnarray*}

(2) Our second claim is that we have the following formula:
$$JxJ(y\Omega)=yx^*\Omega$$

Indeed, this follows from the following computation:
$$JxJ(y\Omega)
=J(xy^*\Omega)
=yx^*\Omega$$

(3) Our claim now is that we have an inclusion as follows:
$$JAJ\subset A'$$ 

Indeed, this follows from the formula obtained in (2).

\medskip

(4) In order to prove the reverse inclusion, our claim is that for $x\in A'$ we have:
$$Jx\Omega=x^*\Omega$$

Indeed, this follows from the following computation, valid for any $y\in A$:
\begin{eqnarray*}
<Jx\Omega,y\Omega>
&=&<Jy\Omega,x\Omega>\\
&=&<y^*\Omega,x\Omega>\\
&=&<\Omega,xy\Omega>\\
&=&<x^*\Omega,y\Omega>
\end{eqnarray*}

(5) Our claim now is that the following formula defines a trace on $A'$:
$$Tr(x)=<x\Omega,\Omega>$$

Indeed, for any two elements $x,y\in A'$ we have:
\begin{eqnarray*}
<xy\Omega,\Omega>
&=&<y\Omega,x^*\Omega>\\
&=&<y\Omega,Jx\Omega>\\
&=&<x\Omega,Jy\Omega>\\
&=&<x\Omega,y^*\Omega>\\
&=&<yx\Omega,\Omega>
\end{eqnarray*}

(6) We can now finish the proof. Indeed, by using the trace constructed in (5), we can apply our results obtained so far to $A'$, and we obtain $JA'J\subset A$, as desired.
\end{proof}

As a basic illustration for the above result, we have:

\begin{theorem}
The commutant of a von Neumann group algebra $L(\Gamma)$, which is obtained by definition by using the left regular representation, is the von Neumann group algebra $R(\Gamma)$, obtained by using the right regular representation.
\end{theorem}

\begin{proof}
We recall that the left and the right representations of a discrete group $\Gamma$ are given by the following formulae, by using the standard identification $\Gamma\subset l^2(\Gamma)$:
$$\lambda_g:h\to gh\quad,\quad\rho_g:h\to hg^{-1}$$

We have $Jg=g^{-1}$ for any group element $g\in\Gamma$, and by using this, we obtain:
\begin{eqnarray*}
J\lambda_gJh
&=&J\lambda_gh^{-1}\\
&=&Jgh^{-1}\\
&=&hg^{-1}\\
&=&\rho_gh
\end{eqnarray*}

Thus, the left and right representations are related by the following formula:
$$J\lambda_gJ=\rho_g$$

By using now Theorem 10.18 we can compute commutants, as follows:
$$L(\Gamma)'
=JL(\Gamma)J
=R(\Gamma)$$

Finally, we have $L(\Gamma)=R(\Gamma)'$ too, by taking the commutant.
\end{proof}

As another application of the standard representation, let us go back to the uniqueness of the trace, that we know from Theorem 10.16. There are several alternative proofs for this fact, which are all instructive. As a first such statement and proof, we have:

\begin{theorem}
Given a ${\rm II}_1$ factor $A$, and an element $a\in A$, we have the following Dixmier averaging property:
$$\overline{span\left\{uau^*\Big|u\in U_A\right\}}^{\,w}\cap\mathbb C1\neq\emptyset$$
In particular, the ${\rm II}_1$ factor trace $tr:A\to\mathbb C$ is unique.
\end{theorem}

\begin{proof}
We use the basic theory of the regular representation $A\subset L^2(A)$, with respect to the given trace $tr:A\to\mathbb C$, explained above. The proof goes as follows:

\medskip

(1) Given an element $a\in A$, consider the space in the statement, obtained as the weak closure of the space spanned by the spinned versions of $a$, namely:
$$K_a=\overline{span\left\{uau^*\Big|u\in U_A\right\}}^{\,w}$$

This linear space $K_a\subset A$ is by definition weakly closed, and it follows that the subset $K_a\Omega\subset L^2(A)$, where $\Omega\in L^2(A)$ is the canonical trace vector, is a weakly closed convex subset. In particular, we see that $K_a\Omega\subset L^2(A)$ is a norm closed convex subset. 

\medskip

(2) In view of this, we can consider the unique element $b\in K_a$ having the property that $b\Omega$ has a minimal norm. We have then the following formula, for any unitary $u\in U_A$, where $J:L^2(A)\to L^2(A)$ is the standard antilinear map, given by $T\to T^*$:
$$||uJuJb\Omega||=||b\Omega||$$

By uniqueness of $b$, it follows that for any unitary $u\in U_A$, we have:
$$uJuJb\Omega=b\Omega$$

But this shows that for any unitary $u\in U_A$, we have:
$$ubu^*=b$$

We conclude that we have $b\in\mathbb C1$, and this proves the first assertion.

\medskip

(3) Regarding now the second assertion, consider an arbitrary trace $tr:A\to\mathbb C$. By using $tr(uau^*)=tr(a)$, we conclude that this trace is constant on the following set:
$$K_a=\overline{span\left\{uau^*\Big|u\in U_A\right\}}^{\,w}$$

Now by using the first assertion, we conclude that we have the following formula:
$$\overline{span\left\{uau^*\Big|u\in U_A\right\}}^{\,w}\cap\mathbb C1=\big\{tr(a)1\big\}$$

Summarizing, we have obtained a purely algebraic formula for our trace $tr:A\to\mathbb C$, and it follows that this trace is indeed unique, as claimed.
\end{proof}

In relation with the above, let us mention that there is as well a third proof for the uniqueness of the trace, due to Yeadon, based on nothing or almost, meaning the definition of the ${\rm II}_1$ factors, along with some abstract functional analysis. For more on all this, basic theory of the ${\rm II}_1$ factors, we refer to the standard operator algebra books, with some good choices here being the books of Connes \cite{co3}, Jones \cite{jo6} and Blackadar \cite{bla}.

\bigskip

Before developing more general theory for the ${\rm II}_1$ factors, let us discuss the examples. We have so far only one class of examples, namely the group von Neumann algebras $L(\Gamma)$, which are ${\rm II}_1$ factors precisely when the discrete groups $\Gamma$ have the ICC property. This suggests doing several things, in order to have more examples, as follows:

\bigskip

(1) A first idea is that of looking at the algebras of discrete quantum groups, $A=L(\Gamma)$, or equivalently, $A=L^\infty(G)$, with $G=\widehat{\Gamma}$ being the compact quantum group dual to $\Gamma$. However, despite years of efforts, no one knows what ``quantum ICC'' should mean. 

\bigskip

(2) A more modest statement would be that if a compact quantum group $G\subset U_N^+$ appears as liberation of a classical group $G_{class}\subset U_N$, then the corresponding von Neumann algebra $A=L^\infty(G)$ should be a ${\rm II}_1$ factor. But this question is open, too.

\bigskip

(3) We can in fact conjecture that if a homogeneous space $X=G/H$, or a more general manifold $X$, appears as liberation of a homogeneous space $X_{class}=G_{class}/H_{class}$, or of a more general manifold $X_{class}$, then $A=L^\infty(X)$ should be a ${\rm II}_1$ factor.

\bigskip

(4) Along the same lines, but having this time von Neumann's reduction theory results in mind, we have the question of understanding how the various quantum group or quantum manifold algebras considered above decompose as sums of ${\rm II}_1$ factors.

\bigskip

Summarizing, we have many difficult questions here, with the Holy Grail being the reduction theory for the algebras of type $A=L^\infty(X)$, with $X$ being a quantum manifold. Fortunately, we have as well a series of alternative questions, also inspired by the group von Neumann algebras $L(\Gamma)$, and which are supposedly easier, as follows:

\bigskip

(5) A group von Neumann algebra $L(\Gamma)$ can be thought of as coming from the trivial action $\Gamma\curvearrowright\{.\}$, and the question is that of investigating von Neumann algebras associated to more general actions, $\Gamma\curvearrowright X$, by using various crossed product techniques.

\bigskip

(6) There are many natural examples of compact groups $G$ acting on von Neumann algebras $P$, and the question is that of understanding under which exact assumptions on the action $G\curvearrowright P$, the corresponding fixed point algebra $P^G$ is a factor.

\bigskip

(7) There are as well many examples of discrete groups $\Gamma$ acting on von Neumann algebras $R$, and the question is that of understanding under which exact assumptions on the action $G\curvearrowright R$, the corresponding crossed product algebra $R\rtimes\Gamma$ is a factor.

\bigskip

(8) Finally, the above questions are related to each other, and even more general questions come by looking at actions of compact quantum groups $G$, or discrete quantum groups $\Gamma$, on various quantum spaces $X$, or von Neumann algebras $P$ or $R$.

\bigskip

In order to get started, in connection with questions (1-2-3-4), let us first talk about free quantum groups, as a continuation of the material from chapters 7-8. The various combinatorial considerations there lead to the following conjecture:

\begin{conjecture}
Assuming that an easy quantum group $G\subset U_N^+$ is free, in the sense that it comes from a category of noncrossing partitions
$$D\subset NC$$
the associated von Neumann algebra $L^\infty(G)$ is a factor.
\end{conjecture}

This is something quite technical, motivated by the findings in \cite{bsp}, which state that the liberation operation $G\to G^+$ for the easy quantum groups corresponds to the operation $D\to D\cap NC$ at the level of the associated categories of partitions. 

\bigskip

As a more general, and also more elementary conjecture, we have:

\begin{conjecture}
Assuming that a closed subgroup $G\subset U_N^+$ satisfies
$$S_N^+\subset G\subset U_N^+$$
the associated von Neumann algebra $L^\infty(G)$ is a factor.
\end{conjecture}

In relation with this, let us recall that the free complexification of a Woronowicz algebra $(A,u)$ with $u\in M_N(A)$ is the Woronowicz algebra $(\tilde{A},\tilde{u})$ constructed as follows, where $z\in C(\mathbb T)$ is the standard generator, given by $x\to x$ for any $x\in\mathbb T$:
$$\tilde{A}=<\tilde{u}_{ij}>\subset C(\mathbb T)*A\quad,\quad 
\tilde{u}=zu\in M_N(\widetilde{A})$$

The point indeed with this notion is that, in the context of the liberation operation $G\to G^+$ discussed above, we usually have embeddings as follows:
$$G\subset\tilde{G}\subset G^+$$ 

Thus, we are led into a conjecture about free complexifications, as follows:

\index{free complexification}

\begin{conjecture}
Given a closed subgroup $G\subset U_N^+$, the von Neumann algebra
$$A=L^\infty(\tilde{G})$$
of $L^\infty$ functions on its free complexification $\tilde{G}\subset U_N^+$ is a factor.
\end{conjecture}

Many more things can be said here, and then about general quantum manifolds as well, with the global conjecture being that if such a manifold $X$ is free, in some suitable algebraic sense, then its associated von Neumann algebra $L^\infty(X)$ should be a factor.

\section*{10c. Type II factors}

Let us go back now to the general theory of the ${\rm II}_1$ factors, with the aim of talking about representations of such ${\rm II}_1$ factors, inside the category of the ${\rm II}_1$ factors, $A\subset B$. For this purpose we will need a key notion, called coupling constant. 

\bigskip

In order to discuss the construction of the coupling constant, we will need some further results on the type ${\rm II}$ factors, complementing those that we already have. The point indeed is that the class of ${\rm II}$ factors, to be axiomatized later, and with this being not something urgent, comprises, besides the ${\rm II}_1$ factors discussed above, the ${\rm II}_\infty$ factors as well:

\begin{definition}
A ${\rm II}_\infty$ factor is a von Neumann algebra of the form
$$B=A\otimes B(H)$$
with $A$ being a ${\rm II}_1$ factor, and with $H$ being an infinite dimensional Hilbert space.
\end{definition}

We should mention that there are several possible ways of defining the ${\rm II}_\infty$ factors, and the above definition is something rather intuitive, the point being that, once you learn the theory of the ${\rm II}_\infty$ factors, as we will do here, what you remember at the end of the day is what has been said above, $B=A\otimes B(H)$, with $A$ being a ${\rm II}_1$ factor. 

\bigskip

Getting started now, as a useful characterization of such factors, we have:

\begin{proposition}
For an infinite factor $B$, the following are equivalent:
\begin{enumerate}
\item There exists a projection $p\in B$ such that $pBp$ is a ${\rm II}_1$ factor.

\item $B$ is a ${\rm II}_\infty$ factor.
\end{enumerate}
\end{proposition}

\begin{proof}
This is something elementary, as follows:

\medskip

$(1)\implies(2)$ Assume indeed that $p\in B$ is a projection such that $pBp$ is a ${\rm II}_1$ factor. We choose a maximal family of pairwise orthogonal projections $\{p_i\}\subset B$ satisfying $p_i\simeq p$, for any $i$, and we consider the following projection, which satisfies $q\preceq p$:
$$q=1-\sum_ip_i$$

Since the indexing set for our set of projections $\{p_i\}$ must be infinite, we can use a strict embedding of this index set into itself, as to write a formula as follows:
\begin{eqnarray*}
1
&=&q+\sum_ip_i\\
&\preceq&p_0+\sum_{i\neq0}p_i\\
&\preceq&1
\end{eqnarray*}

Thus we have $\sum_ip_i\simeq1$, and we may further suppose that we have in fact:
$$\sum_ip_i=1$$

Thus the family $\{p_i\}$ can be used in order to construct a copy $B(H)\subset B$, with $H=l^2(\mathbb N)$, and we must have $B=A\otimes B(H)$, with $A$ being a ${\rm II}_1$ factor, as desired.

\medskip

$(2)\implies(1)$ This is clear, because when assuming $B=A\otimes B(H)$, as in Definition 10.24, we can take our projection $p\in B$ to be of the form $p=1\otimes q$, with $q\in B(H)$ being a rank 1 projection, and we have then $pBp=A$, which is a ${\rm II}_1$ factor, as desired.
\end{proof}

Getting back now to the original interpretation of the ${\rm II}_\infty$ factors, from Definition 10.24, the tensor product writing there $B=A\otimes B(H)$ suggests tensoring the trace of the ${\rm II}_1$ factor $A$ with the usual operator trace of $B(H)$. We are led in this way to:

\begin{definition}
Given a ${\rm II}_\infty$ factor $B$, written as $B=A\otimes B(H)$, with $A$ being a ${\rm II}_1$ factor and with $H$ being an infinite dimensional Hilbert space, we define a map
$$tr:B_+\to[0,\infty]\quad,\quad 
tr((x_{ij}))=\sum_itr(x_{ii})$$
where we have chosen a basis of $H$, as to have $H\simeq l^2(\mathbb N)$, and so $B(H)\subset M_\infty(\mathbb C)$.
\end{definition}

As an important observation, to start with, unlike in the ${\rm II}_1$ factor case, that of the factor $A$, or in the ${\rm I}_\infty$ factor case, that of the factor $B(H)$, it is not possible to suitably normalize the trace constructed above. This follows indeed from the results below.

\bigskip

On the positive side now, the trace that we constructed has all sorts of good properties, that we can use for various purposes, which can be summarized as follows:

\begin{proposition}
The ${\rm II}_\infty$ factor trace that we constructed above 
$$tr:B_+\to[0,\infty]$$
has the following properties:
\begin{enumerate}
\item $tr(x+y)=tr(x)+tr(y)$, and $tr(\lambda x)=\lambda tr(x)$ for $\lambda\geq0$.

\item If $x_i\nearrow x$ then $tr(x_i)\to tr(x)$.

\item $tr(xx^*)=tr(x^*x)$.

\item $tr(uxu^*)=tr(x)$ for any $u\in U_B$.
\end{enumerate}
\end{proposition}

\begin{proof}
All this is elementary, the idea being as follows:

\medskip

(1) This is clear from definitions.

\medskip

(2) This is again clear from definitions. 

\medskip

(3) This is something which is elementary as well.

\medskip

(4) This comes from (3), via the formula $uxu^*=u\sqrt{x}\cdot\sqrt{x}u^*$.
\end{proof}

As a main result now regarding the ${\rm II}_\infty$ factor trace, we have:

\begin{theorem}
The ${\rm II}_\infty$ factor trace $tr:B_+\to[0,\infty]$ constructed above, when restricted to the projections
$$tr:P(B)\to[0,\infty]$$
induces an isomorphism between the totally ordered set of equivalence classes of projections in $B$ and the interval $[0,\infty]$.
\end{theorem}

\begin{proof}
We have several things to be checked here, as follows:

\medskip

(1) Our first claim is that a projection $p\in B$ is finite precisely when $tr(p)<\infty$. 

\medskip

-- Indeed, in one sense, assume that we have $tr(p)<\infty$. If our projection $p$ was to be infinite, we would have a subprojection $q\leq p$ having the same trace as $p$, and so $r=p-q$ would be a projection of trace $0$, which is impossible. Thus $p$ is indeed finite.

\medskip

-- In the other sense now, assuming $tr(p)=\infty$, we have to prove that $p$ is infinite. For this purpose, let us pick a projection $q\leq p$ having finite trace. Then $r=p-q$ satisfies $tr(r)=\infty$, and so we can iterate the procedure, and we end up with an infinite sequence of pairwise orthogonal projections, which are all smaller than $p$. But this shows that $p$ dominates an infinite projection, and so that $p$ itself is infinite, as desired.

\medskip

(2) Our second claim is that if $p,q\in B$ are projections, with $p$ finite, then:
$$p\preceq q\iff tr(p)=tr(q)$$

But this follows exactly as in the ${\rm II}_1$ factor case, discussed above.

\medskip

(3) Our third and final claim, which will finish the proof, is that any infinite projection is equivalent to the identity. For this purpose, assume that $p\in B$ is infinite. By definition, this means that we can find a unitary $u\in B$ such that:
$$uu^*=p\quad,\quad 
u^*u\leq p\quad,\quad uu^*\neq p$$

But these conditions show that $(u^n)^iu^n$ is a strictly decreasing sequence of equivalent projections, and by using this sequence we conclude that we have $1\preceq p$, as desired.
\end{proof}

Moving ahead now, in order to further investigate the ${\rm II}_\infty$ factors, we will need:

\begin{theorem}
Given a ${\rm II}_1$ factor $A\subset B(H)$, there exists an isometry
$$u:H\to L^2(A)\otimes l^2(\mathbb N)$$
such that $ux=(x\otimes1)u$, for any $x\in A$.
\end{theorem}

\begin{proof}
We use a standard idea, that we used many times before, namely an amplification trick. Given a ${\rm II}_1$ factor $A\subset B(H)$, consider the following Hilbert space:
$$K=H\oplus L^2(A)\otimes l^2(\mathbb N)$$

Consider, as operators over this space $K$, the following projections:
$$p=id\oplus 0\quad,\quad 
q=0\oplus id$$

Both these projections $p,q$ belong then to $A'$, which is a type ${\rm II}_\infty$ factor. Now since $q\in A'$ is infinite, by Theorem 10.28 we can find a partial isometry $u\in A'$ such that:
$$u^*u=p\quad,\quad 
uu^*\leq q$$

Now let us represent this partial isometry $u\in B(K)$ as a $2\times2$ matrix, as follows:
$$u=\begin{pmatrix}a&b\\ c&d\end{pmatrix}$$

The above conditions $u^*u=p$ and $uu^*\leq q$ reformulate then as follows:
$$b^*b+d^*d=0\quad,\quad 
aa^*+bb^*=0$$

We conclude that our  partial isometry $u\in B(K)$ has the following special form:
$$u=\begin{pmatrix}0&0\\ c&0\end{pmatrix}$$

But the operator $c:H\to l^2(A)\otimes l^2(\mathbb N)$ that we found in this way must be an isometry, and from $u\in A'$ we obtain $ux=(x\otimes1)u$, for any $x\in A$, as desired.
\end{proof}

As a basic consequence of the above result, which is something good to know, and that we will use many times in what follows, we have:

\index{commutant of factor}

\begin{theorem}
The commutant of a ${\rm II}_1$ factor is a ${\rm II}_1$ factor, or a ${\rm II}_\infty$ factor.
\end{theorem}

\begin{proof}
This follows indeed from the explicit interpretation of the operator algebra embedding $A\subset B(H)$ of our ${\rm II}_1$ factor $A$, found in Theorem 10.29.
\end{proof}

Summarizing, we have an extension of the general theory of the ${\rm II}_1$ factors, developed before, to the general case of the type ${\rm II}$ factors, which comprises by definition the ${\rm II}_1$ factors and the ${\rm II}_\infty$ factors. All this is of course technically very useful.

\section*{10d. Coupling constant}

We are now in position of constructing the coupling constant. The idea here, following as usual the key paper of Murray and von Neumann \cite{mv3}, will be that given a representation of a ${\rm II}_1$ factor $A\subset B(H)$, we can try to understand how far is this representation from the standard form, where $H=L^2(A)$, from ``above'' or from ``below''. 

\bigskip

In order to discuss this, which is something quite technical, let us start with:

\begin{proposition}
Given a ${\rm II}_1$ factor $A\subset B(H)$, with its embedding into $B(H)$ being represented as above, in terms of an isometry
$$u:H\to L^2(A)\otimes l^2(\mathbb N)\quad,\quad 
ux=(x\otimes1)u$$
the following quantity does not depend on the choice of this isometry $u$:
$$C=tr(uu^*)$$
Moreover, for the standard form, where $H=L^2(A)$, this constant takes the value $1$.
\end{proposition}

\begin{proof}
Assume indeed that we have an isometry $u$ as in the statement, and that we have as well a second such isometry, of the same type, namely:
$$v:H\to L^2(A)\otimes l^2(\mathbb N)\quad,\quad 
vx=(x\otimes1)v$$

We have then $uu^*=uv^*vu^*$, and by using this, we obtain:
\begin{eqnarray*}
C_u
&=&tr(uu^*)\\
&=&tr(uv^*vu^*)\\
&=&tr(vu^*uv^*)\\
&=&tr(vv^*)\\
&=&C_v
\end{eqnarray*}

Thus, we are led to the conclusion in the statement. As for the last assertion, regarding the standard form, this is clear from definitions, because here we can take $u=1$.
\end{proof}

As a conclusion to all this, given a ${\rm II}_1$ factor $A\subset B(H)$, we know from  Theorem 10.29 that $H$ must appear as an ``inflated'' version of $L^2(A)$.  The corresponding inflation constant is a certain number, that we can call coupling constant, as follows:

\index{coupling constant}
\index{standard form}

\begin{definition}
Given a representation of a ${\rm II}_1$ factor $A\subset B(H)$, we can talk about the corresponding coupling constant, as being the number
$$\dim_AH\in(0,\infty]$$
constructed as follows, with $u:H\to L^2(A)\otimes l^2(\mathbb N)$ isometry satisfying $ux=(x\otimes1)u$:
$$\dim_AH=tr(uu^*)$$
For the standard form, where $H=L^2(A)$, this coupling constant takes the value $1$.
\end{definition}

This definition might seem a bit complicated, but things here are quite non-trivial, and there is no way of doing something substantially simpler. Alternatively, we can define the coupling constant via the following formula, after proving first that the number on the right is indeed independent of the choice on a nonzero vector $x\in H$:
$$\dim_AH=\frac{tr_A(P_{A'x})}{tr_{A'}(P_{Ax})}$$

This latter formula was in fact the original definition of the coupling constant, by Murray and von Neumann \cite{mv3}. However, technically speaking, things are slightly easier when using the approach in Definition 10.32. We will be back to this key formula of Murray and von Neumann, with full explanations, in a moment. 

\bigskip

Let us start our study of the coupling constant with some basic results, coming from definitions and from what we already have, as results, as follows:

\begin{proposition}
The coupling constant $\dim_AH\in(0,\infty]$ associated to a ${\rm II}_1$ factor representation $A\subset B(H)$ has the following properties:
\begin{enumerate}
\item For the standard form, $H=L^2(A)$, we have $\dim_AH=1$.

\item For the usual representation on $H=L^2(A)\otimes l^2(\mathbb N)$, we have $\dim_AH=\infty$.

\item We have $\dim_AH<\infty$ precisely when $A'$ is a ${\rm II}_1$ factor.

\item We have additivity, $\dim_A(\oplus_iH_i)=\sum_i\dim_AH_i$.

\item We have $\dim_A(L^2(A)p)=tr(p)$, for any projection $p\in A$.

\item The coupling constant can take any value in $(0,\infty]$.
\end{enumerate}
\end{proposition}

\begin{proof}
All these assertions are elementary, the idea being as follows:

\medskip

(1) This is something that we already know, coming from definitions.

\medskip

(2) This is something that comes from definitions too.

\medskip

(3) This comes from the general properties of the ${\rm II}_\infty$ factors, and their traces.

\medskip

(4) Again, this is clear from the definition of the coupling constant.

\medskip

(5) This follows by using $u(x)=x\otimes\xi$, with $\xi\in l^2(\mathbb N)$ being of norm 1.

\medskip

(6) This follows by starting with (5), and then making direct sums, as in (4).
\end{proof}

At a more advanced level now, in relation with projections and compressions, and getting towards the above-mentioned Murray-von Neumann approach, we have:

\begin{proposition}
We have the compression formula
$$\dim_{pAp}(pH)=\frac{\dim_AH}{tr_A(p)}$$
valid for any projection $p\in A$.
\end{proposition}

\begin{proof}
We can prove this result in two steps, as follows:

\medskip

(1) Assume that $H$ is as follows, with $q\in A$ being a projection satisfying $q\leq p$:
$$H=L^2(A)q$$

We can use the following unitary, intertwining the left and right actions of $pAp$:
$$L^2(pAp)\to pL^2(A)p\quad,\quad 
pxp\Omega\to p(x\Omega)p$$

Indeed, we obtain that the following algebras are unitarily equivalent:
$$pAp\subset B(pL^2(A)q)\quad,\quad 
pAp\subset B(L^2(pAp)q)$$

Thus, by using the formula (5) in Proposition 10.33 we obtain, as desired:
\begin{eqnarray*}
\dim_{pAp}(pH)
&=&tr_{pAp}(q)\\
&=&\frac{tr_A(q)}{tr_A(p)}\\
&=&\frac{\dim_AH}{tr_A(p)}
\end{eqnarray*}

(2) In the general case now, where $H$ is arbitrary, the result follows from what we proved above, and from the additivity property from Proposition 10.33 (4).
\end{proof}

With all these properties established, we can now recover, as a theorem, the original definition of the coupling constant, due to Murray and von Neumann, as follows:

\begin{theorem}
Given a ${\rm II}_1$ factor $A\subset B(H)$, with the commutant $A'\subset B(H)$ assumed to be finite, the corresponding coupling constant is finite, given by
$$\dim_AH=\frac{tr_A(P_{A'x})}{tr_{A'}(P_{Ax})}$$
with the number on the right being independent of the choice on a nonzero vector $x\in H$. In the case where $A'$ is infinite, the corresponding coupling constant is infinite.
\end{theorem}

\begin{proof}
There are several things to be proved here, the idea being as follows:

\medskip

(1) We know from Proposition 10.33 (3) that we have $\dim_AH<\infty$ precisely when the commutant $A'\subset B(H)$ is finite. Thus, we may assume that we are in this case.

\medskip

(2) Assuming so, we have the following formula, valid for any projection $p\in A'$, which follows from the basic properties of the coupling constant, established above:
$$\dim_{Ap}(pH)=tr_{A'}(p)\dim_AH$$

(3) Now with this formula in hand, the formula in the statement follows as well, once again by doing a number of standard amplification and compression manipulations.
\end{proof}

As an illustration for all this, given an inclusion of ICC groups $\Lambda\subset\Gamma$, whose group algebras are both ${\rm II}_1$ factors, we have the following formula:
$$\dim_{L(\Lambda)}L^2(\Gamma)=[\Gamma:\Lambda]$$

There are many other examples of explicit computations of the coupling constant, all leading into interesting mathematics. We will be back to this.

\bigskip

As a last topic for this chapter, given a ${\rm II}_1$ factor $A$, let us discuss now the representations of type $A\subset B$, with $B$ being another ${\rm II}_1$ factor. This is a quite natural notion, perhaps even more natural than the representations $A\subset B(H)$, because we have previously decided that the ${\rm II}_1$ factors $B$, and not the full operator algebras $B(H)$, are the correct infinite dimensional generalization of the usual matrix algebras $M_N(\mathbb C)$.

\bigskip

This was for the philosophy, and one can of course agree or not with this. Or at least agree or not at the present point of the presentation, because once we will get into the structure of the subfactors $A\subset B$, which is something amazing, there is no way back.

\bigskip

In practice now, given an inclusion of ${\rm II}_1$ factors $A\subset B$, a first question is that of defining its index, measuring how big is $B$ compared to $A$. The first thought here goes into defining the index of $A\subset B$ as being a purely algebraic quantity, as follows:
$$N=\dim_AB$$

However, this is non-trivial, due to the fact that we are in the ``continuous dimension'' setting, and so our algebraic intuition, where indices are always integers, will not help us much. We will be back to this question later, with a technical solution to it.

\bigskip

In order to solve our index problem, a much better approach is by using the ambient operator algebra $B(H)$, or rather the ambient Hilbert space $H$, as follows:

\index{index of subfactor}

\begin{theorem}
Given an inclusion of ${\rm II}_1$ factors $A\subset B$, the number
$$N=\frac{\dim_AH}{\dim_BH}$$
is independent of the ambient Hilbert space $H$, and is called index.
\end{theorem}

\begin{proof}
The fact that the index of the subfactor $A\subset B$, as defined by the above formula, is indeed independent of the ambient Hilbert space $H$, comes from the various basic properties of the coupling constant, established above.
\end{proof}

There are many examples of subfactors coming from groups, and every time we obtain the intuitive index. More suprisingly now, Jones proved in \cite{jo1} that the index, when small, is in fact ``quantized'', subject to the following unexpected restriction:
$$N\in\left\{4\cos^2\left(\frac{\pi}{n}\right)\Big|n\geq3\right\}\cup[4,\infty]$$

This is in fact part of a series of non-trivial results about the subfactors, due to Jones, and also Ocneanu, Popa, Wassermann and others, and involving as well the Temperley-Lieb algebra \cite{tli}, and many more. We will be back to this later, with the whole last part of the present book, chapters 13-16 below, being dedicated to subfactor theory.

\section*{10e. Exercises}

In relation with the general theory of factors, we have:

\begin{exercise}
Classify the von Neumann factors $A\subset B(H)$ by using as invariant the ordered semigroup formed by the equivalence classes of projections $p\in A$.
\end{exercise}

This is something quite tricky, and there are many things that can be done here, even in the context of the matrix algebras and type ${\rm II}$ factors. We will be back to this.

\begin{exercise}
Find a direct proof for the fact that the traces of projections in $L(\Gamma)$ can take any values in $[0,1]$, for an ICC group $\Gamma$ of your choice.
\end{exercise}

In other words, the question is that of picking a simple ICC group, such as $\Gamma=S_\infty$, and constructing projections in $L(\Gamma)$, whose traces converge to a given $c\in[0,1]$.

\begin{exercise}
Do something, statement and proof, even modest, in relation with the von Neumann algebras $L(\Gamma)=L^\infty(G)$ of the discrete quantum groups $\Gamma=\widehat{G}$. 
\end{exercise}

To be more precise, in what regards conjectures, there are many of them, and we discussed this in the above. The problem now is that of downgrading our dreams and expectations, and finding something which is doable. We will be back to this.

\begin{exercise}
Fully clarify the basic properties of the ${\rm II}_\infty$ factors, and the related construction of the coupling constant.
\end{exercise}

This is something that we already discussed in the above, but with a few details missing, and the problem now is that of clarifying all this. You can either go through the discussion which was made above, and come up with the missing details, or do something alternative, based on the various historical comments given above.

\begin{exercise}
Prove that we have the formula
$$\dim_{L(\Lambda)}L^2(\Gamma)=[\Gamma:\Lambda]$$
for any inclusion of ICC groups $\Lambda\subset\Gamma$
\end{exercise}

Normally this should not be difficult. We will be back to this.

\chapter{Advanced results}

\section*{11a. Reduction theory}

Welcome to advanced von Neumann algebra theory, and in the hope that we will survive, both you reader, and me author. Our main purpose here will be to discuss some key decomposition methods for the von Neumann algebras $A\subset B(H)$, in terms of von Neumann factors, $Z(A)=\mathbb C$, which altogether are called ``reduction theory''. 

\bigskip

The reduction theory, due to von Neumann himself, is something quite fundamental, to the point that it would have made sense to get into it right after the von Neumann algebra basics from chapter 5. However, the subject being quite technical, we have not done this after chapter 5, and nor in fact we will really do it here, with the presentation below being just a modest introduction to all this. So, this is the situation, and in the hope that the gods of mathematics, Bourbaki and others, will pardon us.

\bigskip

The story, first. Von Neumann started to work on operator algebras in the 1930s, and became increasingly convinced that these should be subject to a reduction theory theorem, with his interest in factors, which is obvious in his papers \cite{mv1}, \cite{mv2}, \cite{mv3}, \cite{vn1}, \cite{vn2}, basically coming from this. However, he was not able to come at that time, during his prime years of work, with a complete proof, and paper, on reduction theory. He only did that much later, after a break involving other things, like the Manhattan Project, game theory, computers and more, in his 1949 paper \cite{vn3}, written towards the end of his career. His main theorem in \cite{vn3}, which is quite easy to formulate, is as follows:

\begin{fact}[Reduction theory]
Given a von Neumann algebra $A\subset B(H)$, if we write its center $Z(A)\subset A$, which is a commutative von Neumann algebra, as 
$$Z(A)=L^\infty(X)$$
with $X$ being a measured space, then the whole algebra decomposes as
$$A=\int_XA_x\,dx$$
with the fibers $A_x$ being von Neumann algebra factors, $Z(A_x)=\mathbb C$.
\end{fact}

As a first comment, we have already seen an instance of such decomposition results in chapter 5, when talking about finite dimensional algebras. Indeed, such algebras decompose, in agreement with Fact 11.1, as direct sums of matrix algebras, as follows:
$$A=\bigoplus_xM_{n_x}(\mathbb C)$$

More generally, it is possible to axiomatize a certain class of ``type I algebras'', and then show that these algebras appear as direct integrals of matrix algebras:
$$A=\int_XM_{n_x}(\mathbb C)\,dx$$

Observe in particular that in the case where the decomposition is isotypic, $n_x=N$ for some $N\in\mathbb N$, we obtain the random matrix algebras studied in chapter 6:
$$A=M_N(L^\infty(X))$$

Beyond type I, however, things become quite complicated. Next in the hierarchy is the general ``finite case'', where the algebra is assumed to have a trace:
$$tr:A\to\mathbb C$$

Here the existence of the trace simplifies a bit things, although these still remain fairly  complicated, and actually adds to the final result, in the form of the supplementary formula, regarding its decomposition, the precise statement being as follows:

\begin{fact}[Reduction theory, finite case]
Given a von Neumann algebra $A\subset B(H)$ coming with a trace $tr:A\to\mathbb C$, if we write its center $Z(A)\subset A$ as 
$$Z(A)=L^\infty(X)$$
with $X$ being a measured space, then the whole algebra and its trace decompose as
$$A=\int_XA_x\,dx\quad,\quad tr=\int_Xtr_x\,dx$$
with the fibers $A_x$ being factors which are ``finite'', in the sense that they have traces, which in practice means that they can be usual matrix algebras, or ${\rm II}_1$ factors.
\end{fact}

As already mentioned, while some tricks are potentially available here, coming from the presence of the trace $tr:A\to\mathbb C$, this remains something complicated. As for the most general case, where the von Neumann algebra $A\subset B(H)$ is taken arbitrary, corresponding to Fact 11.1 in full generality, this is something even more complicated, with the only possible tools coming from advanced operator theory, and functional analysis. 

\bigskip

So, this is the situation, and what to do now. We cannot explain the above, because it is too complicated, but we cannot skip it either, because these are fundamentals. This situation has been known to generations of mathematicians, starting with von Neumann himself, who finished and published his reduction theory paper \cite{vn3} long after developing the basics of operator algebra theory, as mentioned above. The various books written afterwards, including Blackadar \cite{bla}, Connes \cite{co3}, Dixmier \cite{dix}, Jones \cite{jo6}, Kadison-Ringrose \cite{kri}, Sakai \cite{sak}, Str\u atil\u a-Zsido \cite{szs} and Takesaki \cite{tak} did not arrange things, being either evasive, or way too technical, not to say unreadable, on this subject. 

\bigskip

The present book won't be an exception to the rule. Our plan in what follows will be that of discussing a bit all this, reduction theory, notably with a study of examples:

\bigskip

(1) First we have the type I algebras, which are direct integrals of matrix algebras $M_{n_x}(\mathbb C)$, with the case $n_x=1$ corresponding to commutativity, the case $n_x\in\mathbb N$ corresponding to the ``type I finite case'', and with the general case being $n_x\in\mathbb N\cup\{\infty\}$. At the level of main examples, these come from finite groups and quantum groups.

\bigskip

(2) Then we have the type II algebras, where we can have both type I and type II factors in the decomposition. Of particular interest is the ``finite'' case, where the algebra is simply assumed to come with a trace, $tr:A\to\mathbb C$, and where the reduction theory result is Fact 11.2, with the factors being matrix algebras $M_N(\mathbb C)$, or ${\rm II}_1$ factors.

\bigskip

(3) Finally, we have the general type III case, with no assumption on the algebra $A\subset B(H)$, corresponding to Fact 11.1. Here the factors in the decomposition can be of type I, or of type II, or neither of type I or II, which are called by definition of type III. The interesting questions here regard the structure of the type III factors.

\bigskip

In order to get started, let us look at the commutative von Neumann algebras. Here we have the following result, that we basically know from chapter 5:

\begin{theorem}
The commutative von Neumann algebras are the algebras of type
$$A=L^\infty(X)$$
with $X$ being a measured space. Thus, we formally have for them the formula
$$A=\int_XA_x\,dx$$
with the fibers $A_x$ being trivial in this case, $A_x=\mathbb C$, for any $x\in X$.
\end{theorem}

\begin{proof}
We have several assertions to be proved, the idea being as follows:

\medskip

(1) In one sense, we must prove that given a measured space $X$, we can realize the commutative algebra $A=L^\infty(X)$ as a von Neumann algebra, on a certain Hilbert space $H$. But this is something that can be done via multiplicity operators, as follows: 
$$L^\infty(X)\subset B(L^2(X))$$

(2) In the other sense, given a commutative von Neumann algebra $A\subset B(H)$, we must construct a certain measured space $X$, and an identification $A=L^\infty(X)$. But this can be done by writing our von Neumann algebra as follows:
$$A=<T_i>$$

Indeed, no matter what particular family of generators $\{T_i\}$ we choose for our algebra $A$, these generators $T_i$ will be commuting normal operators. Thus the spectral theorem for such families of operators, from chapter 3, applies and gives the result.

\medskip

(3) In fact, by using the theory of projections from chapters 9-10, we can write our commutative von Neumann algebra $A\subset B(H)$ in singly generated form:
$$A=<T>$$

But this simplifies the situation, because the basic spectral theorem, for single normal operators, from chapter 3, applies to our generator $T$, and gives the result.

\medskip

(4) Finally, the last assertion, regarding the validity of the reduction theory result in this case, is something trivial, and of course without much practical interest.
\end{proof}

Moving forward, the above result is not the end of the story with the commutative von Neumann algebras, because we still have to understand how a given such algebra $A=L^\infty(X)$, or rather the weak topology isomorphism class of such an algebra, can be represented as an operator algebra, over the various Hilbert spaces $H$:
$$L^\infty(X)\subset B(H)$$

But this can be again solved by writing our algebra as $A=<T>$, and then applying the spectral theorem for normal operators, with the conclusion that the commutative von Neumann algebras are, up to spatial isomorphism, the algebras of the following form, with $X$ being a measured space, and with all this being up to a multiplicity:
$$L^\infty(X)\subset B(L^2(X))$$

With these results in hand, we are now in position of better understanding the idea behind von Neumann's reduction theory. Indeed, given an arbitrary von Neumann algebra $A\subset B(H)$, the idea is to consider its center, and write it as follows:
$$Z(A)=L^\infty(X)\subset B(H)$$

The point is then that everything will decompose over the measured space $X$, and in particular, the whole algebra $A$ itself will decompose as a direct integral of fibers:
$$A=\int_XA_x\,dx$$

As already mentioned, we will only partly explain this in what follows, and by insisting on examples. Also, we will do this slowly, following the type I, II, III hierarchy.

\section*{11b. Type I algebras}

In order to decompose our von Neumann algebras into factors, we must first make some upgrades to our terminology and notations regarding the factors, as follows:

\begin{definition}
The von Neumann algebras having trivial center, also called factors, can be divided into several types, as follows:
\begin{enumerate}
\item The matrix algebra $M_N(\mathbb C)$ is of type ${\rm I}_N$.

\item The operator algebra $B(H)$, with $H$ separable, is of type ${\rm I}_\infty$.

\item The factors which are infinite dimensional and have a trace are of type ${\rm II}_1$.

\item The tensor products $A\otimes B(H)$, with $A$ being a ${\rm II}_1$ factor, are of type ${\rm II}_\infty$.

\item As for the factors left, these are called of type ${\rm III}$.
\end{enumerate}
\end{definition}

It is possible to be more abstract here, but in practice, this is how these factors are best remembered. Now back to reduction theory, we will present it gradually, by following the general type I, II, III hierarchy for the von Neumann algebras, coming from the above classification of factors. Let us first discuss the type I case. Here as starting point we have the following result, which is something that we know well, from chapter 5:

\begin{theorem}
The finite dimensional von Neumann algebras $A\subset B(H)$ are exactly the direct sums of matrix algebras,
$$A=M_{n_1}(\mathbb C)\oplus\ldots\oplus M_{n_k}(\mathbb C)$$
with the summands coming by decomposing the unit into central minimal projections, $1=P_1+\ldots+P_k$. Thus, the general reduction theory formula, namely
$$A=\int_XA_x\,dx$$
holds for them, with the measured space $X$, coming via the formula $Z(A)=L^\infty(X)$, being in this case a finite space, $X=\{1,\ldots,k\}$, and with the fibers being matrix algebras. 
\end{theorem}

\begin{proof}
This is something that we know well from chapter 5. The center of $A$ is a finite dimensional commutative von Neumann algebra, of the following form:
$$Z(A)=\mathbb C^k$$

Now let $P_i$ be the Dirac mass at $i\in\{1,\ldots,k\}$. Then $P_i\in B(H)$ is an orthogonal projection, and these projections form a partition of unity. With $A_i=P_iAP_i$, it is elementary to check that we have a non-unital $*$-algebra decomposition, as follows:
$$A=A_1\oplus\ldots\oplus A_k$$

On the other hand, it follows from the minimality of each of the projections $P_i\in Z(A)$ that we have $A_i\simeq M_{n_i}(\mathbb C)$. Thus, we are led to the conclusion in the statement.
\end{proof}

It is possible to further build on the above result, in several directions, either by allowing the factors in the decomposition to be type ${\rm I}_\infty$ factors as well, that is, $A_x\simeq B(H)$, or by allowing the center to be an infinite measured space, $|X|=\infty$, or by allowing both. The first possible generalization is not very interesting. The second possible generalization, however, is something quite interesting, and we have here:

\index{reduction theory}
\index{type I algebra}
\index{type I factor}
\index{integral of algebras}

\begin{fact}[Reduction theory, type I finite case]
Given a von Neumann algebra $A\subset B(H)$ which is of discrete type, and has a trace $tr:A\to\mathbb C$, we can write
$$A=\int_XA_x\,dx$$
with $X$ coming via $Z(A)=L^\infty(X)$,  and the trace decomposes as well, as
$$tr=\int_Xtr_x\,dx$$
with the fibers $A_x$ being usual matrix algebras, $A_x=M_{n_x}(\mathbb C)$, with $n_x\in\mathbb N$.
\end{fact}

As a first observation, this statement generalizes both what we know about the commutative algebras, and the finite dimensional ones. However, having these two things jointly generalized is something quite technical, that we will not explain here in detail. The idea is of course first that of axiomatizing what ``discrete'' should mean in the above, say by looking at the finiteness properties of the projections $p\in A$, and then, once the statement properly formulated, to prove it by jointly generalizing what we know about the commutative algebras, and the finite dimensional ones. 

\bigskip

Moving ahead, let us lift now the assumption that the factors in the decomposition are of type ${\rm I}_N$, with $N<\infty$. We are led in this way to a general result, as follows:

\begin{fact}[Reduction theory, type I case]
Given a von Neumann algebra $A\subset B(H)$ which is of type I, in the sense that it is of a suitable discrete type, we can write
$$A=\int_XA_x\,dx$$
with $X$ coming via $Z(A)=L^\infty(X)$, and with the fibers $A_x$ being type ${\rm I}$ factors, meaning $A_x\simeq B(H_x)$, with each $H_x$ being either finite dimensional, or separable.
\end{fact}

As before with Fact 11.6, we will not attempt to explain this here. As a comment, however, this can only follow from Fact 11.6 applied to the ``finite'' part of the algebra, obtained by removing the infinite part, and after proving that this infinite part is something of type $L^\infty(Y)\otimes B(H)$, with $Y\subset X$, and with $H$ being separable.

\bigskip

All the above was quite abstract, and as something more concrete now, let us discuss the reduction theory for the group von Neumann algebras $L(\Gamma)$, in the finite case, $|\Gamma|<\infty$. For this purpose, it is convenient to change a bit our terminology and notations, making them more in tune with the quantum group formalism from chapter 7. First, we will denote our finite group $\Gamma$, which is at the same time discrete and compact, by $F$, and we will think of it as being the dual of a finite quantum group $G=\widehat{F}$. Also, since in the finite group case everything is automatically norm or weakly closed, we will use the more familiar notation $C^*(F)$ for the associated von Neumann algebra $L(F)$. With these conventions, we have the following result, which is standard: 

\begin{theorem}
Given a finite group $F$, the center of the associated von Neumann algebra is isomorphic to the algebra of central functions on $F$,
$$Z(C^*(F))\simeq C(F)_{central}$$
and the reduction theory applied to this algebra, which is a formula of type 
$$C^*(F)\simeq\bigoplus_{r\in X}M_{n_r}(\mathbb C)$$
appears by dualizing the Peter-Weyl decomposition of the usual function algebra
$$C(F)\simeq\bigoplus_{r\in Irr(F)}M_{\dim(r)}(\mathbb C)$$
via the standard identification between representations $r$ and their characters $\chi_r$.
\end{theorem}

\begin{proof}
In what concerns the first assertion, regarding the center, this is something that we already know, from chapter 10, coming from our study there of the general group algebras $L(\Gamma)$, with $\Gamma$ being a discrete group. To be more precise, in the case where $\Gamma=F$ is a finite group, the computation there gives the following formula for the center:
$$Z(C^*(F))=\left\{\sum_g\lambda_gg\Big|\lambda_{gh}=\lambda_{hg},\forall h\in F\right\}''$$

Now since on the right we have central functions on our group, $\lambda\in C(F)_{central}$, we obtain the isomorphism in the statement, namely:
$$Z(C^*(F))\simeq C(F)_{central}$$

Regarding now the second assertion, let us first recall that the Peter-Weyl theory applied to the finite group $F$ gives a direct sum decomposition as follows, which is technically an isomorphism of linear spaces, which is in addition a $*$-coalgebra isomorphism:
$$C(F)\simeq\bigoplus_{r\in Irr(F)}M_{\dim(r)}(\mathbb C)$$

Thus by dualizing, which is a standard functional analysis procedure, to be explained more in detail below, in a more general setting, we obtain a direct sum decomposition of the group algebra, as follows, which is this time a $*$-algebra isomorphism:
$$C^*(F)\simeq\bigoplus_{r\in Irr(F)}M_{\dim(r)}(\mathbb C)$$

Our claim now, which will finish the proof, is that this is exactly what comes out from von Neumann's reduction theory, applied to the von Neumann algebra $L(F)=C^*(F)$. Indeed, by using the standard identification between representations $r$ and their characters $\chi_r$, which are central functions on $F$, the center computation that we did above reads:
$$Z(C^*(F))\simeq L^\infty(Irr(F))$$

We conclude that von Neumann's reduction theory, applied to the von Neumann algebra $L(F)=C^*(F)$, gives a $*$-algebra isomorphism of the following type:
$$C^*(F)\simeq\bigoplus_{r\in Irr(F)}M_{n_r}(\mathbb C)$$

But a careful examination of the fibers appearing in this decomposition, based on their very definition, shows that these are precisely the above matrix blocks coming from Peter-Weyl. That is, we have $n_r=\dim(r)$ for any $r\in Irr(F)$, and we are done.
\end{proof}

Our next goal will be that of extending the above result to the finite quantum group setting. For this purpose, we will not really need the general compact quantum group formalism from chapter 7, and it is more convenient to start with:

\begin{definition}
A finite dimensional Hopf algebra is a finite dimensional $C^*$-algebra, with comultiplication, counit and antipode maps, satisfying the conditions
$$(\Delta\otimes id)\Delta=(id\otimes \Delta)\Delta$$
$$(\varepsilon\otimes id)\Delta=(id\otimes\varepsilon)\Delta=id$$
$$m(S\otimes id)\Delta=m(id\otimes S)\Delta=\varepsilon(.)1$$
along with the extra condition $S^2=id$. Given such an algebra we write 
$$A=C(G)=C^*(F)$$
and call $G,F$ finite quantum groups, dual to each other.
\end{definition}

In this definition everything is standard, except for the last axiom, $S^2=id$, which corresponds to the fact that, in the corresponding quantum group, we should have:
$$(g^{-1})^{-1}=g$$

It is possible to prove that this condition is automatic, in the present $C^*$-algebra setting. However, this is something non-trivial, and since all this is just an informative discussion, not needed later, we have opted for including $S^2=id$ in our axioms. 

\bigskip

We say that an algebra $A$ as above is cocommutative if $\Sigma\Delta=\Delta$, where $\Sigma(a\otimes b)=b\otimes a$ is the flip. With this convention made, we have the following result, which summarizes the basic theory of finite quantum groups, and justifies the terminology and axioms:

\begin{theorem}
The following happen:
\begin{enumerate}
\item If $G$ is a finite group then $C(G)$ is a commutative Hopf algebra, with
$$\Delta(\varphi)=(g,h)\to \varphi(gh)\quad,\quad 
\varepsilon(\varphi)=\varphi(1)\quad,\quad
S(\varphi)=g\to\varphi(g^{-1})$$
as structural maps. Any commutative Hopf algebra is of this form. 

\item If $F$ is a finite group then $C^*(F)$ is a cocommutative Hopf algebra, with
$$\Delta(g)=g\otimes g\quad,\quad 
\varepsilon(g)=1\quad,\quad
S(g)=g^{-1}$$
as structural maps. Any cocommutative Hopf algebra is of this form.

\item If $G,F$ are finite abelian groups, dual to each other via Pontrjagin duality,
$$C(G)=C^*(F)$$
as Hopf algebras, coming via a Fourier transform type operation.
\end{enumerate}
\end{theorem}

\begin{proof}
These results are all elementary, the idea being as follows:

\medskip

(1) The fact that $\Delta,\varepsilon,S$ satisfy the axioms is clear from definitions, and the converse follows from the Gelfand theorem, by working out the details, regarding $\Delta,\varepsilon,S$. 

\medskip

(2) Once again, the fact that $\Delta,\varepsilon,S$ satisfy the axioms is clear from definitions. For the converse, we use a trick. Let $A$ be an arbitrary finite dimensional Hopf algebra, as in Definition 11.9, and consider its comultiplication, counit, multiplication, unit and antipode maps. The transposes of these maps are then linear maps as follows:
$$\Delta^t:A^*\otimes A^*\to A^*$$
$$\varepsilon^t:\mathbb C\to A^*$$
$$m^t:A^*\to A^*\otimes A^*$$
$$u^t:A^*\to\mathbb C$$
$$S^t:A^*\to A^*$$

It is routine to check that these maps make $A^*$ into a Hopf algebra. Now assuming that  $A$ is cocommutative, it follows that $A^*$ is commutative, so by (1) we obtain $A^*=C(G)$ for a certain finite group $G$, which in turn gives $A=C^*(G)$, as desired.

\medskip

(3) This follows indeed from the discussion in the proof of (2), and from the general theory of Pontrjagin duality for finite abelian groups, explained in chapter 7.
\end{proof}

There are many other things that can be said about the finite dimensional Hopf algebras, and in what follows we will be particularly interested in the notion of corepresentation. These corepresentations can be introduced as follows:

\begin{definition}
A unitary corepresentation of a finite dimensional Hopf algebra $A$ is a unitary matrix $u\in M_n(A)$ satisfying the following conditions:
$$\Delta(u_{ij})=\sum_ku_{ik}\otimes u_{kj}\quad,\quad 
\varepsilon(u_{ij})=\delta_{ij}\quad,\quad 
S(u_{ij})=u_{ji}^*$$
We say that $u$ is irreducible, and we write $u\in Irr(A)$, when it has no nontrivial intertwiners, in the sense that $Tu=uT$ with $T\in M_n(\mathbb C)$ implies $T\in\mathbb C1$.
\end{definition}

Observe the similarity with the notions introduced in chapter 7, for the Woronowicz algebras. In fact, by using left regular representations we can see that any finite dimensional Hopf algebra in the sense of Definition 11.9 is a Woronowicz algebra in the sense of chapter 7. Thus, we can freely use here the results established in chapter 7, and in particular, we can use the Peter-Weyl type theory developed there.

\bigskip

In relation now with our von Neumann algebra questions, we have the following result, coming from that Peter-Weyl type theory, which generalizes Theorem 11.8:

\begin{theorem}
Given a finite quantum group $F$, the center of the associated von Neumann algebra is isomorphic to the algebra of central functions on $F$,
$$Z(C^*(F))\simeq C(F)_{central}$$
and the reduction theory applied to this algebra, which is a formula of type 
$$C^*(F)\simeq\bigoplus_{u\in X}M_{n_u}(\mathbb C)$$
appears by dualizing the Peter-Weyl decomposition of the usual function algebra
$$C(F)\simeq\bigoplus_{u\in Irr(F)}M_{\dim(u)}(\mathbb C)$$
via the standard identification between representations $u$ and their characters $\chi_u$.
\end{theorem}

\begin{proof}
The proof here is nearly identical to the proof of Theorem 11.8. To be more precise, with the more familiar notation $A=C^*(F)$, the proof goes as follows:

\medskip

(1) In what concerns the first assertion, regarding the center, we recall from Woronowicz \cite{wo1} that $A_{central}$ is by definition the subalgebra of $A$ appearing as follows:
$$A_{central}=\left\{a\in A\Big|\Delta a=a\right\}$$

But this shows, first by dualizing, and then by doing some computations similar to those that we did in chapter 10, when computing the centers of the usual group von Neumann algebras, that we have an isomorphism as in the statement, namely:
$$Z(A)\simeq(A^*)_{central}$$

(2) Regarding now the second assertion, we recall that the Peter-Weyl theory applied to Hopf algebra $A^*$ gives a direct sum decomposition as follows, which is technically an isomorphism of linear spaces, which is in addition a $*$-coalgebra isomorphism:
$$A^*\simeq\bigoplus_{u\in Irr(A^*)}M_{\dim(u)}(\mathbb C)$$

Thus by dualizing, we obtain a direct sum decomposition of the group algebra, as follows, which is this time a $*$-algebra isomorphism:
$$A\simeq\bigoplus_{u\in Irr(A^*)}M_{\dim(u)}(\mathbb C)$$

(3) Our claim now, which will finish the proof, is that this is exactly what comes out from von Neumann's reduction theory, applied to the algebra $A$. Indeed, by using the standard identification between corepresentations $u$ of $A^*$ and their characters $\chi_u$, which belong to the algebra $(A^*)_{central}$, the center computation that we did above reads:
$$Z(A)\simeq L^\infty(Irr(A^*))$$

We conclude that von Neumann's reduction theory, applied to the von Neumann algebra $A$, gives a $*$-algebra isomorphism of the following type:
$$A\simeq\bigoplus_{u\in Irr(A^*)}M_{n_u}(\mathbb C)$$

But a careful examination of the fibers shows that these are precisely the matrix blocks coming from Peter-Weyl. That is, $n_u=\dim(u)$ for any $u\in Irr(A^*)$, and we are done.
\end{proof}

All this is quite interesting, and it is possible to say more about it. However, when it comes to type I algebras, in general, the following comment is unavoidable:

\begin{comment}
The most interesting type I algebras are probably those having an isotypic decomposition, and so which can be written as follows:
$$A=M_N(L^\infty(X))$$
But these are precisely the random matrix algebras, that we investigated in great detail in chapter 6, right after introducing the von Neumann algebras. So, job done.
\end{comment}

Needless to say, this is something subjective. But, in any case, whether you agree or not with this, now you know more on the organization of the present book.

\section*{11c. Type II algebras}

Let us discuss now the type II case, where the truly interesting problems are. The central result here, that we already formulated in the beginning of this chapter, is:

\index{type II algebra}
\index{type II factor}
\index{integral of algebras}

\begin{fact}[Reduction theory, finite case]
Given a von Neumann algebra $A\subset B(H)$ coming with a trace $tr:A\to\mathbb C$, if we write its center $Z(A)\subset A$ as 
$$Z(A)=L^\infty(X)$$
with $X$ being a measured space, then the whole algebra and its trace decompose as
$$A=\int_XA_x\,dx\quad,\quad tr=\int_Xtr_x\,dx$$
with the fibers $A_x$ being either factors of type ${\rm I}_N$, with $N<\infty$, or of type ${\rm II}_1$.
\end{fact}

Regarding the proof, this is something quite technical, generalizing what we know, or rather what we don't, about the type I finite case, which is substantially easier. We refer here to Dixmier \cite{dix}, and with the comment that we will see soon examples of all this.

\bigskip

As before in the type I case, it is possible to add a bit of infinity in the above, and we have the following result, which is a bit more general, but more technical too:

\begin{fact}[Reduction theory, type II case]
Given a von Neumann algebra $A\subset B(H)$ which is of type ${\rm II}$, in a suitable sense, if we write its center $Z(A)\subset A$ as 
$$Z(A)=L^\infty(X)$$
with $X$ being a measured space, then the whole algebra decomposes as
$$A=\int_XA_x\,dx$$
with the fibers $A_x$ being von Neumann factors of type ${\rm I}$ or ${\rm II}$.
\end{fact}

As before with what happened in type I, the above results are particularly interesting in the case of the von Neumann algebras of the discrete groups, $A=L(\Gamma)$, and their generalizations. In order to discuss these questions, let us recall that the center of an arbitrary group von Neumann algebra $A=L(\Gamma)$ consists, up to some standard identifications, of the functions which are constant on the finite conjugacy classes. This suggests the following definition, which is something well-known in group theory:

\begin{definition}
A discrete group $F$ is said to have the FC property if all its conjugacy classes are finite. In other words, for any $g\in F$, we must have:
$$\left|\left\{hgh^{-1}\Big|h\in F\right\}\right|<\infty$$
If this finite conjugacy property is satisfied, we also say that $F$ is a FC group.
\end{definition}

As basic examples of FC groups, we have the finite groups, the abelian groups, and the products of such groups. Besides being stable under taking products, the class of FC groups is stable under a number of other basic operations, such as taking subgroups, or quotients. In connection now with our reduction theory questions, we have:

\begin{theorem}
Given a group $F$ having the FC property, the center of the associated von Neumann algebra is isomorphic to the algebra of central functions on $F$,
$$Z(L(F))\simeq C(F)_{central}$$
and the reduction theory applied to this algebra, which is a formula of type 
$$L(F)\simeq\int_{r\in X}A_r$$
appears in relation with the representation theory of $F$.
\end{theorem}

\begin{proof}
In what concerns the first assertion, regarding the center, this is something that we know  from chapter 10. Indeed, we have the following formula for the center:
$$Z(L(F))=\left\{\sum_g\lambda_gg\Big|\lambda_{gh}=\lambda_{hg},\forall h\in F\right\}''$$

Now since on the right we have central functions on our group, $\lambda\in C(F)_{central}$, we obtain the isomorphism in the statement, namely:
$$Z(L(F))\simeq C(F)_{central}$$

Regarding now the second assertion, this is something more tricky, as follows:

\medskip

(1) In the finite group case, we recall from Theorem 11.8 that, by using the standard identification between representations $r$ and their characters $\chi_r$, which are central functions on $F$, the center computation that we did above reads:
$$Z(L(F))\simeq L^\infty(Irr(F))$$

In order to discuss now the reduction theory for $L(F)$, we recall that the Peter-Weyl theory applied to $F$ gives a direct sum decomposition as follows, which is technically an isomorphism of linear spaces, which is in addition a $*$-coalgebra isomorphism:
$$L^\infty(F)\simeq\bigoplus_{r\in Irr(F)}M_{\dim(r)}(\mathbb C)$$

Thus by dualizing, we obtain a direct sum decomposition of the group von Neumann algebra as follows, which is this time a $*$-algebra isomorphism:
$$L(F)\simeq\bigoplus_{r\in Irr(F)}M_{\dim(r)}(\mathbb C)$$

But this is exactly what comes out from von Neumann's reduction theory, applied to the algebra $L(F)$, and so we are fully done with the finite group case.

\medskip

(2) As a second key particular case, let us discuss now the case where $F$ is abelian. In the simplest infinite group case, where our group is $F=\mathbb Z$, the group algebra is:
$$L(\mathbb Z)\simeq L^\infty(\mathbb T)$$

More generally, for the abelian groups $F=\mathbb Z^N$, which are those which are finitely generated and without torsion, we obtain the algebras of functions on various tori:
$$L(\mathbb Z^N)\simeq L^\infty(\mathbb T_N)$$

In general now, assuming that $F$ is finitely generated and abelian, here we know from Pontrjagin duality that we have an isomorphism as follows:
$$L(F)\simeq L^\infty(\widehat{F})$$

More explicitly now, let us write our finitely generated abelian group $F$ as a product of cyclic groups, possibly taken infinite, as follows:
$$F=\mathbb Z^N\times\left(\prod_i\mathbb Z_{n_i}\right)$$

The Pontrjagin dual of $F$ is then the following compact abelian group:
$$F=\mathbb T^N\times\left(\prod_i\mathbb Z_{n_i}\right)$$

Thus, things are very explicit here, and we are done with the abelian case too.

\medskip

(3) In the general case now, where our discrete group $F$ is only assumed to have the FC property, the reduction theory for the corresponding von Neumann algebra $L(F)$ appears somewhat as a mixture of what happens for the finite and for the abelian groups, discussed in (1) and (2) above. For more on all this, we refer to Dixmier \cite{dix}. 
\end{proof}

Regarding the corresponding problems for the discrete quantum groups, these are not solved yet. In fact, the knowledge here stops at a very basic level, with the analogue of the ICC property, leading to the factoriality of $L(\Gamma)$, not being known yet, and for more on all this, we refer to the discussion made in chapter 10.

\bigskip

Moving ahead from these difficulties, let us go back now to the usual group von Neumann algebras $L(\Gamma)$, and discuss what happens in general. Once again inspired by the basic computation that we have, namely that of the center of an arbitrary group algebra $L(\Gamma)$, let us formulate the following purely group-theoretical definition:

\begin{definition}
Given a discrete group $\Gamma$, its FC subgroup $F\subset\Gamma$ is the subgroup
$$F=\left\{g\in\Gamma\Big|\left|\left\{hgh^{-1}\Big|h\in\Gamma\right\}\right|<\infty\right\}$$
consisting of the elements in the finite conjugacy classes of $\Gamma$.
\end{definition}

Here the fact that $F$ is indeed a subgroup is clear from definitions, with the fact that $F$ is stable under multiplication coming from the following trivial observation:
$$h(gk)h^{-1}=hgh^{-1}\cdot hkh^{-1}$$

Observe that $\Gamma$ has the FC property, in the sense of Definition 11.16, precisely when the inclusion $F\subset\Gamma$ is an equality. As before with the FC groups, there are many known things about the FC subgroups $F\subset\Gamma$, and we refer here to the group theory literature. 

\bigskip

In connection now with our reduction theory questions, we have:

\begin{theorem}
Given a discrete group $\Gamma$, the center of the associated von Neumann algebra is isomorphic to the algebra of central functions on its FC subgroup $F\subset\Gamma$,
$$Z(L(\Gamma))\simeq C(F)_{central}$$
and the reduction theory applied to this algebra, which is a formula of type 
$$L(\Gamma)\simeq\int_{r\in X}A_r$$
appears in relation with the representation theory of $\Gamma$, and of its FC subgroup $F\subset\Gamma$.
\end{theorem}

\begin{proof}
In what concerns the first assertion, regarding the center, this is something that we know from chapter 10, coming from our study there of the general group algebras $L(\Gamma)$, with $\Gamma$ being a discrete group. To be more precise, we know from there that:
$$Z(L(\Gamma))=\left\{\sum_g\lambda_gg\Big|\lambda_{gh}=\lambda_{hg},\forall h\in F\right\}''$$

Now since on the right we have central functions on the FC subgroup, $\lambda\in C(F)_{central}$, we obtain the isomorphism in the statement, namely:
$$Z(L(\Gamma))\simeq C(F)_{central}$$

Regarding the second assertion, this is something more tricky, and we refer here to the relevant group theory and operator algebra literature, including Dixmier \cite{dix}.
\end{proof}

As a last topic for this section, let us briefly discuss the reduction theory in the general case, type III. In order to get started, we must discuss the type III factors, which are new to us. According to our various conventions above, these factors are defined as follows:

\index{type III factor}

\begin{definition}
A type ${\rm III}$ factor is a von Neumann algebra $A\subset B(H)$ which is a factor, $Z(A)=\mathbb C$, and satisfies one of the following equivalent conditions:
\begin{enumerate}
\item $A$ is not of type ${\rm I}$, or of type ${\rm II}$.

\item $A$ has no semifinite trace $tr:A\to\mathbb C$.

\item $A$ has no trace $tr:A\to\mathbb C$, and is not of type ${\rm I}_\infty$ or ${\rm II}_\infty$.
\end{enumerate}
\end{definition}

In order to investigate such factors, the general idea will be that of looking at the crossed products of type II factors, which can be lacking traces $tr:A\to\mathbb C$, and so which allow us to exit the type II world. In order to get started, however, we have:

\begin{theorem}
Any locally compact group $G$ has a left invariant Haar measure $\lambda$, and a right invariant Haar measure $\rho$, 
$$d\lambda(x)=d\lambda(yx)\quad,\quad d\rho(x)=d\lambda(xy)$$
which are unique up to multiplication by scalars. These two measures are absolutely continuous with respect to each other, and the Radon-Nikodym derivative 
$$m:G\to\mathbb R\quad,\quad m(x)=\frac{d\lambda(x)}{d\rho(x)}$$
well-defined up to multiplication by scalars, is called modulus of the group. The unimodular groups, for which $m=1$, include all compact groups, and all abelian groups.
\end{theorem}

\begin{proof}
There are many things here, with everything being very classical, and the proof, along with comments, examples and more theory, especially in what regards the unimodular groups, can be found in any good measure theory book.
\end{proof}

As it has become customary in this book, whenever talking about groups we must make some comments about quantum groups too. Things are quite interesting in connection with Theorem 11.21, because it is possible ``twist'' things in the compact case, as to have a notion of modulus there as well. We refer here to Woronowicz \cite{wo1} and related papers. In relation now with our factor questions, we have the following result:

\index{Tomita-Takesaki}
\index{KMS state}
\index{Connes classification}
\index{type III factor}

\begin{theorem}
The type ${\rm III}$ factors basically appear from the type ${\rm II}$ factors, via various crossed product constructions, and their generalizations.
\end{theorem}

\begin{proof}
This statement is obviously something quite informal, and we will certainly not attempt to explain the proof either. Here are however the main ideas, with the result itself being basically due to Connes \cite{co1}, along with some historical details:

\medskip

(1) First of all, Murray and von Neumann knew of course about such questions, but were quite evasive in their papers about type III, with the brief comment ``we don't know''. Whether they really worked or not on these questions, we'll never know.

\medskip

(2) Inspired by Theorem 11.21, it is possible to develop a whole machinery for the study of the non-tracial states $\varphi:A\to\mathbb C$, the main results here being the Kubo-Martin-Schwinger (KMS) condition, and the Tomita-Takesaki theory. See Takesaki \cite{tak}.

\medskip

(3) On the other hand, looking at type II factors and their crossed products by automorphisms, which are not necessarily of type II, leads to a lot of interesting theory as well, leading to large classes of type III factors, appearing from type II factors.

\medskip

(4) The above results are basically from the 50s and 60s, and Connes was able to put all this together, in the early 70s, via a series of quick, beautiful and surprising Comptes Rendus notes, eventually leading to his paper \cite{co1}, which is a must-read.
\end{proof}

In equivalent terms, and also by remaining a bit informal, we have:

\begin{theorem}
The von Neumann algebra factors can be classified as follows,
$${\rm I_N},{\rm I}_\infty$$
$${\rm II}_1,{\rm II}_\infty$$
$${\rm III}_0,{\rm III}_\lambda,{\rm III}_1$$
with the type ${\rm II}_1$ ones being the most important, basically producing the others too.
\end{theorem}

\begin{proof}
This follows by putting altogether what we have, results of Murray and von Neumann in type I and II, and then of Connes in type III. The last assertion is of course something quite informal, because the situation is not exactly as simple as that.
\end{proof}

Getting back now to our series of reduction theory results, we have:

\index{type III algebra}
\index{type III factor}
\index{reduction theory}

\begin{theorem}
Given an arbitrary von Neumann algebra $A\subset B(H)$, write its center as follows, with $X$ being a measured space:
$$Z(A)=L^\infty(X)$$ 
The whole algebra $A$ decomposes then over this measured space $X$, as a direct sum of fibers, taken in an appropriate sense,
$$A=\int_XA_x\,dx$$
with the fibers $A_x$ being von Neumann factors, which can be of type ${\rm I},{\rm II},{\rm III}$.
\end{theorem}

\begin{proof}
As before with other such results, this is something heavy, generalizing our previous knowledge in type I, and type II. The proof however is quite similar, basically using the same ideas. We refer here to the literature, for instance to Dixmier \cite{dix}.
\end{proof}

\section*{11d. Abelian subalgebras}

We would like to end this chapter with something more refreshing, in relation with the above, namely matrix models, and abelian subalgebras. These two topics are actually related, and in a quite subtle way, and we will provide an introduction to this.

\bigskip

Let us first discuss the matrix models, for the quantum groups. One interesting method for the study of the closed subgroups $G\subset U_N^+$ consists in modelling the coordinates $u_{ij}\in C(G)$ by concrete variables $U_{ij}\in B$. Indeed, assuming that the model is faithful in some suitable sense, and that the target algebra $B$ is something quite familiar, all questions about $G$ would correspond in this way to routine questions inside $B$. 

\bigskip

Regarding now the choice of the target algebra $B$, some very familiar and convenient algebras are the random matrix ones, $B=M_K(C(T))$, with $K\in\mathbb N$, and $T$ being a compact space. We are led in this way to the following definition:

\index{matrix model}
\index{random matrix model}

\begin{definition}
A matrix model for $G\subset U_N^+$ is a morphism of $C^*$-algebras
$$\pi:C(G)\to M_K(C(T))$$
where $T$ is a compact space, and $K\geq1$ is an integer.
\end{definition}

There are many examples of such models, and will discuss them later on. For the moment, let us develop some general theory. The question to be solved is that of understanding the suitable faithfulness assumptions needed on $\pi$, as for the model to ``remind'' the quantum group. The simplest situation is when $\pi$ is faithful in the usual sense. Let us introduce the following notion, which is related to faithfulness:

\index{stationary model}

\begin{definition}
A matrix model $\pi:C(G)\to M_K(C(T))$ is called stationary when
$$\int_G=\left(tr\otimes\int_T\right)\pi$$
where $\int_T$ is the integration with respect to a given probability measure on $T$.
\end{definition}

Here the term ``stationary'' comes from a functional analytic interpretation of all this, with a certain Ces\`aro limit being needed to be stationary, and this will be explained later. Yet another explanation comes from a certain relation with the lattice models, but this relation is rather something folklore, not axiomatized yet. We will be back to this. 

\bigskip

As a first result now, which is something which is not exactly trivial, and whose proof requires some functional analysis, the stationarity property implies the faithfulness:

\index{coamenability}
\index{type I algebra}
\index{stationarity}

\begin{theorem}
Assuming that $G\subset U_N^+$ has a stationary model,
$$\pi:C(G)\to M_K(C(T))\quad,\quad 
\int_G=\left(tr\otimes\int_T\right)\pi$$
it follows that $G$ is coamenable, and that the model is faithful, coming as:
$$C(G)\subset L^\infty(G)\subset M_K(L^\infty(T))$$
Moreover, we can have such models only when the algebra $L^\infty(G)$ is of type ${\rm I}$.
\end{theorem}

\begin{proof}
We use the basic theory of compact and discrete quantum groups, developed in chapter 7. Assume that we have a stationary model, as in the statement. By performing the GNS construction with respect to $\int_G$, we obtain a factorization as follows, which commutes with the respective canonical integration functionals:
$$\pi:C(G)\to C(G)_{red}\subset M_K(C(T))$$

Thus, in what regards the coamenability question, we can assume that $\pi$ is faithful. With this assumption made, observe that we have embeddings as follows:
$$C^\infty(G)\subset C(G)\subset M_K(C(T))$$

The point now is that the GNS construction gives a better embedding, as follows:
$$L^\infty(G)\subset M_K(L^\infty(T))$$

Now since the von Neumann algebra on the right is of type I, so must be its subalgebra $A=L^\infty(G)$. This means that, when writing the center of this latter algebra as  $Z(A)=L^\infty(X)$, the whole algebra decomposes over $X$, as an integral of type I factors:
$$L^\infty(G)=\int_XM_{K_x}(\mathbb C)\,dx$$

In particular, we can see from this that $C^\infty(G)\subset L^\infty(G)$ has a unique $C^*$-norm, and so $G$ is coamenable. Finally, the other assertions follow as well from the above, because our factorization of $\pi$ consists of the identity, and of an inclusion.
\end{proof}

More generally now, we can talk about matrix models for the algebraic submanifolds $X\subset S^{N-1}_{\mathbb C,+}$, in the obvious way, and we have the following result:

\index{stationary model}
\index{faithful model}

\begin{theorem}
Given a matrix model $\pi:C(X)\to M_K(C(T))$, with both $X,T$ being assumed to have integration functionals, the following are equivalent:
\begin{enumerate}
\item $\pi$ is stationary, in the sense that $\int_X=(tr\otimes\smallint_T)\pi$.

\item $\pi$ produces an inclusion $\pi':C_{red}(X)\subset M_K(X(T))$.

\item $\pi$ produces an inclusion $\pi'':L^\infty(X)\subset M_K(L^\infty(T))$.
\end{enumerate}
Moreover, in the quantum group case, these conditions imply that $\pi$ is faithful.
\end{theorem}

\begin{proof}
Consider the following diagram, with all the solid arrows being by definition the canonical maps between the various algebras concerned:
$$
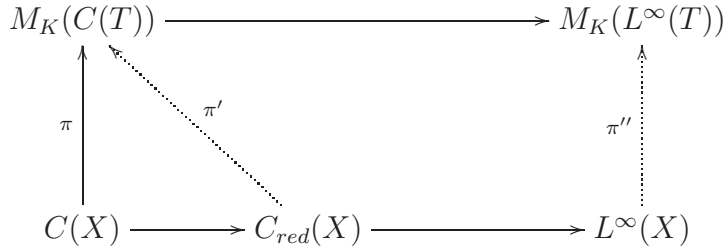
\xymatrix@R=60pt@C=30pt{
M_K(C(T))\ar[rrr]&&&M_K(L^\infty(T))\\
C(X)\ar[u]^\pi\ar[r]&C_{red}(X)\ar[rr]\ar@.[ul]_{\pi'}&&L^\infty(X)\ar@.[u]^{\pi''}
}$$

With this picture in hand, the equivalences $(1)\iff(2)\iff(3)$ between the above conditions (1,2,3) are all clear, coming from the basic properties of the GNS construction. As for the last assertion, this is something that we know from Theorem 11.27.
\end{proof}

Moving ahead now, our claim is that our modelling philosophy, with type ${\rm I}$ algebras as target, and more specifically with random matrix algebras as target, can perfectly apply, at least in the quantum group case, to the type ${\rm II}$ algebras as well. 

\bigskip

We have indeed the following result, which is something quite subtle:

\index{Hopf image}
\index{inner faithfulness}

\begin{theorem}
Given a matrix model $\pi:C(G)\to M_K(C(T))$, with $T$ being a probability space, there exists a smallest subgroup $G'\subset G$ producing a factorization
$$\pi:C(G)\to C(G')\to M_K(C(T))$$
with the intermediate algebra $C(G')$ being called Hopf image of $\pi$. When $\pi$ is inner faithful, in the sense that we have $G=G'$, we have the formula 
$$\int_G=\lim_{k\to\infty}\sum_{r=1}^k\varphi^{*r}$$
where $\varphi=(tr\otimes\smallint_T)\pi$ is the matrix model trace, and where $\phi*\psi=(\phi\otimes\psi)\Delta$. Also, the model $\pi$ is stationary precisely when this latter convergence is stationary.
\end{theorem}

\begin{proof}
All this is well-known, the idea being as follows:

\medskip

(1) The construction of the Hopf image can be done by dividing the algebra $C(G)$ by a suitable ideal, but for our purposes here it is more convenient to go via an alternative proof. Let us denote by $u=(u_{ij})$ the fundamental corepresentation of $G$, and consider the following vector spaces, taken in a formal sense, where $U_{ij}=\pi(u_{ij})$:
$$C_{kl}=Hom(U^{\otimes k},U^{\otimes l})$$

Since the morphisms increase the intertwining spaces, when defined either in a representation theory sense, or just formally, we have inclusions as follows:
$$Hom(u^{\otimes k},u^{\otimes l})\subset Hom(U^{\otimes k},U^{\otimes l})$$

More generally, we have such inclusions when replacing $(G,u)$ with any pair producing a factorization of $\pi$. Thus, by Woronowicz's Tannakian duality \cite{wo2}, the Hopf image must be given by the fact that the intertwining spaces must be the biggest, subject to the above inclusions. But since $u$ is biunitary, so is $U$, and it follows that the above spaces $C_{kl}$ form a Tannakian category, so have a quantum group $(G',v)$ given by:
$$Hom(v^{\otimes k},v^{\otimes l})=Hom(U^{\otimes k},U^{\otimes l})$$

By the above discussion, $C(G')$ follows to be the Hopf image of $\pi$, as claimed.

\medskip

(2) The formula for $\int_G$ follows by adapting Woronowicz's construction of the Haar integration functional, from \cite{wo1}, to the matrix model situation. If we denote by $\int_G'$ the limit in the statement, we must prove that this limit converges, and that we have:
$$\int_G'=\int_G$$

It is enough to check this on the coefficients of corepresentations, and if we let $w=u^{\otimes k}$ be one of the Peter-Weyl corepresentations, we must prove that we have:
$$\left(id\otimes\int_G'\right)w=\left(id\otimes\int_G\right)w$$

We know from chapter 7 that the matrix on the right is the orthogonal projection onto $Fix(w)$. Regarding now the matrix on the left, this is the orthogonal projection onto the $1$-eigenspace of $(id\otimes\varphi\pi)w$. Now observe that, if we set $W_{ij}=\pi(w_{ij})$, we have:
$$(id\otimes\varphi\pi)w=(id\otimes\varphi)W$$

Thus, exactly as in chapter 7, we conclude that the $1$-eigenspace that we are interested in equals $Fix(W)$. But, according to the proof of (1) above, we have:
$$Fix(W)=Fix(w)$$

Thus, we have proved that we have $\int_G'=\int_G$, as desired.
\end{proof}

The above result, with contributions by many people, and we refer to \cite{ba3} for the story, is quite important, for many reasons, mainly coming from the following fact:

\begin{fact}
There is no known restriction on the quantum groups having a model
$$\pi:C(G)\to M_K(C(T))$$
which is inner faithful, in the above sense.
\end{fact}

Which is obviously something interesting, conjecturally making Theorem 11.29 a clever way of passing from type ${\rm II}$ to type ${\rm I}$. There are also connections here with the Connes embedding problem, and with all sorts of questions from algebra, geometry, analysis and probability, coming from both mathematics and physics. We will be back to this.

\bigskip

In the general quantum algebraic manifold setting now, talking about inner faithfulness is in general not possible, unless our manifold $X\subset S^{N-1}_{\mathbb C,+}$ has some extra special structure, as for instance being an affine homogeneous space, and we refer here to \cite{ba3}.

\bigskip

Changing topics now, let us go back to the arbitrary von Neumann algebras, and explore some further perspectives opened by the various results that we know. Given a von Neumann algebra $A\subset B(H)$, looking at the center $Z(A)=A\cap A'$ is not the only possible way of getting to commutative, or abelian subalgebras, and we have as well:

\begin{definition}
Given a von Neumann algebra $A\subset B(H)$, an abelian subalgebra $B\subset A$ which is maximal, in the sense that there is no bigger abelian algebra
$$B\subset B'\subset A$$
is called maximal abelian subalgebra (MASA).
\end{definition}

It is possible to say many interesting things about the MASA, and skipping some details here, if we want to further build on this notion, we are led to:

\begin{definition}
Given a von Neumann algebra $A$ coming with a trace $tr:A\to\mathbb C$,
assume that we have a pair of maximal abelian subalgebras
$$B,C\subset A$$
satisfying the following orthogonality condition, with respect to the trace:
$$(B\ominus\mathbb C1)\perp(C\ominus\mathbb C1)$$
We say then that $B,C$ are orthogonal maximal abelian subalgebras.
\end{definition}

Here the scalar product is by definition $<b,c>=tr(bc^*)$, and by taking into account the multiples of the identity, the orthogonality condition reformulates as follows:
$$tr(bc)=tr(b)tr(c)$$

As explained by Popa in \cite{po1}, the interest in Definition 11.32 comes from the fact that a pair of orthogonal MASA brings some sort of 2D orientation inside the von Neumann algebra $A$, or at least inside the subalgebra $<B,C>\subset A$ generated by the MASA. There is also an obvious link with the notion of noncommutative independence discussed in chapter 8. But more on all this later, in chapter 15 below, when doing subfactors.

\bigskip

As a ``toy example'', we can try and see what happens for the simplest factor that we know, namely the matrix algebra $M_N(\mathbb C)$, endowed with its usual matrix trace. And in this context, we have the following surprising result of Popa \cite{po1}:

\index{Hadamard matrix}

\begin{theorem}
Up to a conjugation by a unitary, the pairs of orthogonal maximal abelian subalgebras in the simplest factor, namely $M_N(\mathbb C)$, are as follows,
$$A=\Delta\quad,\quad
B=H\Delta H^*$$
with $\Delta\subset M_N(\mathbb C)$ being the diagonal matrices, and with $H\in M_N(\mathbb C)$ being Hadamard, in the sense that $|H_{ij}|=1$ for any $i,j$, and the rows of $H$ are pairwise orthogonal.
\end{theorem}

\begin{proof}
Any maximal abelian subalgebra in $M_N(\mathbb C)$ being conjugated to $\Delta$, we can assume, up to conjugation by a unitary, that we have, with $U\in U_N$:
$$A=\Delta\quad,\quad 
B=U\Delta U^*$$  

Now observe that given two diagonal matrices $D,E\in\Delta$, we have:
\begin{eqnarray*}
tr(D\cdot UEU^*)
&=&\frac{1}{N}\sum_i(DUEU^*)_{ii}\\
&=&\frac{1}{N}\sum_{ij}D_{ii}U_{ij}E_{jj}\bar{U}_{ij}\\
&=&\frac{1}{N}\sum_{ij}D_{ii}E_{jj}|U_{ij}|^2
\end{eqnarray*}

Thus, the orthogonality condition $A\perp B$ reformulates as follows:
$$\frac{1}{N}\sum_{ij}D_{ii}E_{jj}|U_{ij}|^2=\frac{1}{N^2}\sum_{ij}D_{ii}E_{jj}$$

Thus the rescaled matrix $H=\sqrt{N}U$ must satisfy the following condition:
$$|H_{ij}|=1$$

Thus, we are led to the conclusion in the statement.
\end{proof}

The Hadamard matrices appearing in Theorem 11.33 are well-known objects, appearing in several branches of combinatorics, and quantum physics. The basic examples of such  matrices are the Fourier matrices of abelian groups, constructed as follows:

\index{Fourier matrix}

\begin{theorem}
Given a finite abelian group $G$, with dual group $\widehat{G}=\{\chi:G\to\mathbb T\}$, consider the Fourier coupling $\mathcal F_G:G\times\widehat{G}\to\mathbb T$, given by $(i,\chi)\to\chi(i)$.
\begin{enumerate}
\item Via the standard isomorphism $G\simeq\widehat{G}$, this Fourier coupling can be regarded as a square matrix, $F_G\in M_G(\mathbb T)$, which is a complex Hadamard matrix.

\item For the cyclic group $G=\mathbb Z_N$ we obtain in this way, via the standard identification $\mathbb Z_N=\{1,\ldots,N\}$, the standard Fourier matrix, $F_N=(w^{ij})$ with $w=e^{2\pi i/N}$.

\item In general, when using a decomposition $G=\mathbb Z_{N_1}\times\ldots\times\mathbb Z_{N_k}$, the corresponding Fourier matrix is given by $F_G=F_{N_1}\otimes\ldots\otimes F_{N_k}$.
\end{enumerate}
\end{theorem}

\begin{proof}
This follows indeed from some basic facts from group theory:

\medskip

(1) With the identification $G\simeq\widehat{G}$ made our matrix is given by $(F_G)_{i\chi}=\chi(i)$, and the scalar products between the rows are computed as follows:
$$<R_i,R_j>
=\sum_\chi\chi(i)\overline{\chi(j)}
=\sum_\chi\chi(i-j)
=|G|\cdot\delta_{ij}$$

Thus, we obtain indeed a complex Hadamard matrix.

\medskip

(2) This follows from the well-known and elementary fact that, via the identifications $\mathbb Z_N=\widehat{\mathbb Z_N}=\{1,\ldots,N\}$, the Fourier coupling here is as follows, with $w=e^{2\pi i/N}$:
$$(i,j)\to w^{ij}$$

(3) We use here the following well-known formula, for the duals of products: 
$$\widehat{H\times K}=\widehat{H}\times\widehat{K}$$

At the level of the corresponding Fourier couplings, we obtain from this:
$$F_{H\times K}=F_H\otimes F_K$$

Now by decomposing $G$ into cyclic groups, as in the statement, and by using (2) for the cyclic components, we obtain the formula in the statement.
\end{proof}

Summarizing, we have some interesting connections with finite group theory, and with the associated Fourier matrices. However, there are as well many exotic examples of Hadamard matrices, nor necessarily coming from finite groups, as in Theorem 11.34, and all this is quite of interest for us, in connection with Theorem 11.33.

\bigskip

We will be back to this later, with more results on the subject, in chapters 13-16, when talking about subfactors. Among others, we will see there that the combinatorics of the MASA associated to an Hadamard matrix comes from a certain quantum permutation groups, appearing as in Theorem 11.29, via a matrix model. More on this soon.

\section*{11e. Exercises} 

Things have been quite technical in this chapter, which was more of a survey than something else, and as a unique exercise on all this, we have:

\begin{exercise}
Learn some more basic von Neumann algebra theory, from the papers of von Neumann and Murray-von Neumann, then Tomita-Takesaki and Connes, and write down a brief account of what you learned.
\end{exercise}

In what follows we will avoid ourselves this type of exercise, basically by getting back to the material in chapter 10, and building on that, following Jones.

\chapter{Hyperfiniteness}

\section*{12a. The factor R}

Welcome to advanced operator algebra theory, again. What we saw in the previous chapter was in fact just half of the story, and the other half, regarding hyperfiniteness, still remains to be told. The idea indeed is that there has been a considerable amount of work on hyperfiniteness, comparable in size and difficulty with the general classification work for the factors, based on reduction theory, and we will discuss this here.

\bigskip

In practice, all this will be quite independent from what we did in chapter 11. What we have to do is to go back to the functional analysis methods for general von Neumann algebras developed in chapter 9, and to the theory of factors, and notably of the type ${\rm II}_1$ factors developed in chapter 10, with the aim of further building on this. Following old, classical work of Murray-von Neumann \cite{mv3}, our main object of study will be the central example of a ${\rm II}_1$ factor, namely the ``smallest'' one, the hyperfinite ${\rm II}_1$ factor $R$.

\bigskip

Once this factor $R$ introduced, and its basic theory understood, we will go on a more advanced discussion, including more theory of $R$, following Connes \cite{co2}, then a discussion of various quantum group aspects, as a continuation of what has been said in chapter 10, and finally a discussion of the connections with the material in chapter 11.

\bigskip

Needless to say, this chapter will be a bit like the previous one, more of a survey. Also, let us mention that afterwards, in chapters 13-16 below, we will go back to a more normal pace, with a standard introduction to the Jones theory of inclusions of ${\rm II}_1$ factors, with full details. The notion of hyperfiniteness and the factor $R$ will of course show up there, every now and then, but usually at the end of each chapter, and most of the time using actually only its basic theory, and not most of the advanced material below.

\bigskip

In order to get started now, let us formulate the following definition:

\begin{definition}
A von Neumann algebra $A\subset B(H)$ is called hyperfinite when it appears as the weak closure of an increasing limit of finite dimensional algebras:
$$A=\overline{\bigcup_iA_i}^{\,w}$$
When $A$ is a ${\rm II}_1$ factor, we call it hyperfinite ${\rm II}_1$ factor, and we denote it by $R$.
\end{definition}

As a first observation, there are many hyperfinite von Neumann algebras, for instance because any finite dimensional von Neumann algebra $A=\oplus_iM_{n_i}(\mathbb C)$ is such an algebra, as one can see simply by taking $A_i=A$ for any $i$, in the above definition. 

\bigskip

Also, given a measured space $X$, by using a dense sequence of points inside it, we can write $X=\bigcup_iX_i$ with $X_i\subset X$ being an increasing sequence of finite subspaces, and at the level of the corresponding algebras of functions this gives a decomposition as follows, which shows that the algebra $A=L^\infty(X)$ is hyperfinite, in the above sense:
$$L^\infty(X)=\overline{\bigcup_iL^\infty(X_i)}^{\,w}$$

The interesting point, however, is that when trying to construct ${\rm II}_1$ factors which are hyperfinite, all the possible constructions lead in fact to the same factor, denoted $R$. This is an old theorem of Murray and von Neumann \cite{mv3}, that we will explain now.

\bigskip

In order to get started, we will need a number of technical ingredients. Generally speaking, out main tool will be the expectation $E_i:A\to A_i$ from a hyperfinite von Neumann algebra $A$ onto its finite dimensional subalgebras $A_i\subset A$, so talking about such conditional expectations will be our first task. Let us start with:

\index{conditional expectation}

\begin{proposition}
Given an inclusion of finite von Neumann algebras $A\subset B$, there is a unique linear map
$$E:B\to A$$
which is positive, unital, trace-preserving and satisfies the following condition:
$$E(b_1ab_2)=b_1E(a)b_2$$
This map is called conditional expectation from $B$ onto $A$.
\end{proposition}

\begin{proof}
We make use of the standard representation of the finite von Neumann algebra $B$, with respect to its trace $tr:B\to\mathbb C$, as constructed in chapter 10:
$$B\subset L^2(B)$$

If we denote by $\Omega$ the cyclic and separating vector of $L^2(B)$, we have an identification of vector spaces $A\Omega=L^2(A)$. Consider now the following orthogonal projection:
$$e:L^2(B)\to L^2(A)$$

It follows from definitions that we have an inclusion $e(B\Omega)\subset A\Omega$, and so our projection $e$ induces by restriction a certain linear map, as follows:
$$E:B\to A$$

This linear map $E$ and the orthogonal projection $e$ are then related by:
$$exe=E(x)e$$

But this shows that the linear map $E$ satisfies the various conditions in the statement, namely positivity, unitality, trace preservation and bimodule property. As for the uniqueness assertion, this follows by using the same argument, applied backwards, the idea being that a map $E$ as in the statement must come from the projection $e$.
\end{proof}

Following Jones \cite{jo1}, who was a heavy user of such expectations, we will be often interested in what follows in the orthogonal projection $e:L^2(B)\to L^2(A)$ producing the expectation $E:B\to A$, rather than in $E$ itself. So, let us formulate:

\index{Jones projection}

\begin{definition}
Associated to any inclusion of finite von Neumann algebras $A\subset B$, as above, is the orthogonal projection
$$e:L^2(B)\to L^2(A)$$
producing the conditional expectation $E:B\to A$ via the following formula:
$$exe=E(x)e$$
This projection is called Jones projection for the inclusion $A\subset B$.
\end{definition}

We will heavily use Jones projections in chapters 13-16 below, in the context where both the algebras $A,B$ are ${\rm II}_1$ factors, when systematically studying the inclusions of such ${\rm II}_1$ factors $A\subset B$, called subfactors. In connection with our present hyperfiniteness questions, the idea, already mentioned above, will be that of using the conditional expectation $E_i:A\to A_i$ from a hyperfinite von Neumann algebra $A$ onto its finite dimensional subalgebras $A_i\subset A$, as well as its Jones projection versions $e_i:L^2(A)\to L^2(A_i)$. Let us start with a technical approximation result, as follows:

\begin{proposition}
Assume that a von Neumann algebra $A\subset B(H)$ appears as an increasing limit of von Neumann subalgebras
$$A=\overline{\bigcup_iA_i}^{\,w}$$
and denote by $E_i:A\to A_i$ the corresponding conditional expectations.
\begin{enumerate}
\item We have $||E_i(x)-x||\to0$, for any $x\in A$.

\item If $x_i\in A_i$ is a bounded sequence, satisfying $x_i=E_i(x_{i+1})$ for any $i$, then this sequence has a norm limit $x\in A$, satisfying $x_i=E_i(x)$ for any $i$.
\end{enumerate}
\end{proposition}

\begin{proof}
Both the assertions are elementary, as follows:

\medskip

(1) In terms of the Jones projections $e_i:L^2(A)\to L^2(A_i)$ associated to the expectations $E_i:A\to A_i$, the fact that the algebra $A$ appears as the increasing union of its subalgebras $A_i$ translates into the fact that the $e_i$ are increasing, and converging to $1$:
$$e_i\nearrow1$$

But this gives $||E_i(x)-x||\to0$, for any $x\in A$, as desired.

\medskip

(2) Let $\{x_i\}\subset A$ be a sequence as in the statement. Since this sequence was assumed to be bounded, we can pick a weak limit $x\in A$ for it, and we have then, for any $i$:
$$E_i(x)=x_i$$

Now by (1) we obtain from this $||x-x_n||\to0$, which gives the result.
\end{proof}

We have now all the needed ingredients for formulating a first key result, in connection with the hyperfinite ${\rm II}_1$ factors, due to Murray-von Neumann \cite{mv3}, as follows:

\begin{proposition}
Given an increasing union on matrix algebras, the following construction produces a hyperfinite ${\rm II}_1$ factor
$$R=\overline{\bigcup_{n_i}M_{n_i}(\mathbb C)}^{\,w}$$
called Murray-von Neumann hyperfinite factor.
\end{proposition}

\begin{proof}
This basically follows from the above, in two steps, as follows:

\medskip

(1) The von Neumann algebra $R$ constructed in the statement is hyperfinite by definition, with the remark here that the trace on it $tr:R\to\mathbb C$ comes as the increasing union of the traces on the matrix components $tr:M_{n_i}(\mathbb C)\to\mathbb C$, and with all the details here being elementary to check, by using the usual standard form technology.

\medskip

(2) Thus, it remains to prove that $R$ is a factor. For this purpose, pick an element belonging to its center, $x\in Z(R)$, and consider its expectation on $A_i=M_{n_i}(\mathbb C)$:
$$x_i=E_i(x)$$

We have then $x_i\in Z(A_i)$, and since the matrix algebra $A_i=M_{n_i}(\mathbb C)$ is a factor, we deduce from this that this expected value $x_i\in A_i$ is given by:
$$x_i=tr(x_i)1=tr(x)1$$

On the other hand, Proposition 12.4 applies, and shows that we have:
$$||x_i-x||=||E_i(x)-x||\to0$$

Thus our element is a scalar, $x=tr(x)1$, and so $R$ is a factor, as desired.
\end{proof}

Next, we have the following substantial improvement of the above result, also due to Murray-von Neumann \cite{mv3}, which will be our final saying on the subject:

\index{R}
\index{hyperfinite factor}
\index{Murray-von Neumann factor}
\index{limit of matrix algebras}

\begin{theorem}
There is a unique hyperfinite ${\rm II}_1$ factor, called Murray-von Neumann hyperfinite factor $R$, which appears as an increasing union on matrix algebras,
$$R=\overline{\bigcup_{n_i}M_{n_i}(\mathbb C)}^{\,w}$$
with the isomorphism class of this union not depending on the exact sizes of the matrix algebras involved, nor on the particular inclusions between them.
\end{theorem}

\begin{proof}
We already know from Proposition 12.5 that the union in the statement is a hyperfinite ${\rm II}_1$ factor, for any choice of the matrix algebras involved, and of the inclusions between them. Thus, in order to prove the result, it all comes down in proving the uniqueness of the hyperfinite ${\rm II}_1$ factor. But this can be proved as follows:

\medskip

(1) Given a ${\rm II}_1$ factor $A$, a von Neumann subalgebra $B\subset A$, and a subset $S\subset A$, let us write $S\subset_\varepsilon B$ when the following condition is satisfied, with $||x||_2=\sqrt{tr(x^*x)}$:
$$\forall x\in S,\exists y\in B, ||x-y||_2\leq\varepsilon$$

With this convention made, given a ${\rm II}_1$ factor $A$, the fact that this factor is hyperfinite in the sense of Definition 12.1 tells us that for any finite subset $S\subset A$, and any $\varepsilon>0$, we can find a finite dimensional von Neumann subalgebra $B\subset A$ such that:
$$S\subset_\varepsilon B$$

(2) With this observation made, assume that we are given a hyperfinite ${\rm II}_1$ factor $A$. Let us pick a dense sequence $\{x_k\}\subset A$, and let us set:
$$S_k=\{x_1,\ldots,x_k\}$$

By choosing $\varepsilon=1/k$ in the above, we can find, for any $k\in\mathbb N$, a finite dimensional von Neumann subalgebra $B_k\subset A$ such that the following condition is satisfied:
$$S_k\subset_{1/k}B_k$$

(3) Our first claim is that, by suitably choosing our subalgebra $B_k\subset A$, we can always assume that this is a matrix algebra, of the following special type:
$$B_k=M_{2^{n_k}}(\mathbb C)$$

But this is something which is quite routine, which can be proved by starting with a finite dimensional subalgebra $B_k\subset A$ as above, and then perturbing its set of minimal projections $\{e_i\}$ into a set of projections $\{e_i'\}$ which are close in norm, and have as traces multiples of $2^n$, with $n>>0$. Indeed, the algebra $B_k'\subset A$ having these new projections $\{e_i'\}$ as minimal projections will be then arbitrarily close to the algebra $B_k$, and so will still contain the subset $S_k$ in the above approximate sense, and due to our trace condition, will be contained in a subalgebra of type $B_k''\simeq M_{2^{n_k}}(\mathbb C)$, as desired.

\medskip

(4) Our next claim, whose proof is similar, by using standard perturbation arguments for the corresponding sets of minimal projections, is that in the above the sequence of subalgebras $\{B_k\}$ can be chosen increasing. Thus, up to a rescaling of everything, we can assume that our sequence of subalgebras $\{B_k\}$ is as follows:
$$B_k=M_{2^k}(\mathbb C)$$

(5) But this finishes the proof. Indeed, according to the above, we have managed to write our arbitrary hyperfinite ${\rm II}_1$ factor $A$ as a weak limit of the following type:
$$A=\overline{\bigcup_kM_{2^k}(\mathbb C)}^{\,w}$$

Thus we have uniqueness indeed, and our result is proved.
\end{proof}

The above result is something quite fundamental, and adds to a series of similar results, or rather philosophical conclusions, which are quite surprising, as follows:

\bigskip

(1) We have seen early on in this book that, up to isomorphism, there is only one Hilbert to be studied, namely the infinite dimensional separable Hilbert space, which can be taken to be, according to knowledge and taste, either $H=L^2(\mathbb R)$, or $H=l^2(\mathbb N)$.

\bigskip

(2) Regarding now the study of the operator algebras $A\subset B(H)$ over this unique Hilbert space, another somewhat surprising conclusion, from chapter 6, is that we won't miss much by assuming that $A=M_N(L^\infty(X))$ is a random matrix algebra.

\bigskip

(3) And now, guess what, what we just found is that when trying to get beyond random matrices, and what can be done with them, we are led to yet another unique von Neumann algebra, namely the above Murray-von Neumann hyperfinite ${\rm II}_1$ factor $R$.

\bigskip

(4) And for things to be complete, we will see later that when getting beyond type ${\rm II}_1$, things won't change, because the other types of hyperfinite factors, not necessarily of type ${\rm II}_1$, can be all shown to ultimately come from $R$, via various constructions. 

\bigskip

All this is certainly quite interesting, philosophically speaking. All in all, always the same conclusion, no need to go far to get to interesting algebras and questions: these interesting algebras and questions are just there, the most obvious ones.

\bigskip

Now back to more concrete things, one question is about how to best think of $R$, with Theorem 12.6 as stated not providing us with an answer. To be more precise, we would like to know what is the ``best model'' for $R$, that is, what exact matrix algebras should we use in practice, and with which inclusions between them. And here, a look at the proof of Theorem 12.6 suggests that the ``best writing'' of $R$ is as follows:
$$R=\overline{\bigcup_kM_{2^k}(\mathbb C)}^{\,w}$$

And we can in fact do even better, by observing that the inclusions between matrix algebras of size $2^k$ appear via tensor products, and formulating things as follows:

\begin{proposition}
The hyperfinite ${\rm II}_1$ factor $R$ appears as
$$R=\overline{\bigotimes_{r\in\mathbb N}M_2(\mathbb C)}^{\,w}$$
with the infinite tensor product being defined as an inductive limit, in the obvious way.
\end{proposition}

\begin{proof}
This follows from the above discussion, and with the remark that there is a binary choice there, of left/right type, to be made when constructing the inductive limit. And we prefer here not to make any choice, and leave things like this, because the best choice here always depends on the precise applications that you have in mind.  
\end{proof}

Along the same lines, we can ask as well for precise group algebra models for the hyperfinite ${\rm II}_1$ factor, $R=L(\Gamma)$, and the canonical choice here is as follows:

\begin{proposition}
The hyperfinite ${\rm II}_1$ factor $R$ appears as
$$R=L(S_\infty)$$
with $S_\infty=\bigcup_{r\in\mathbb N}S_r$ being the infinite symmetric group.
\end{proposition}

\begin{proof}
Consider indeed the infinite symmetric group $S_\infty$, which is by definition the group of permutations of $\{1,2,3,\ldots\}$ having finite support. Since such an infinite permutation with finite support must appear by extending a certain finite permutation $\sigma\in S_r$, with fixed points outside $\{1,\ldots,r\}$, we have then, as stated:
$$S_\infty=\bigcup_{r\in\mathbb N}S_r$$

But this shows that the von Neumann algebra $L(S_\infty)$ is hyperfinite. On the other hand $S_\infty$ has the ICC property, and so $L(S_\infty)$ is a ${\rm II}_1$ factor. Thus, $L(S_\infty)=R$.
\end{proof}

There are of course some more things that can be said here, because other groups of the same type as $S_\infty$, namely appearing as increasing limits of finite subgroups, and having the ICC property, will produce as well the hyperfinite factor, $L(\Gamma)=R$, and so there is some group theory to be done here, in order to fully understand such groups. However, we prefer to defer the discussion for later, after learning about amenability, which will lead to a substantial update of our theory, making such things obsolete.

\bigskip

As an interesting consequence of all this, however, let us formulate:

\begin{proposition}
Given two groups $\Gamma,\Gamma'$, each having the ICC property, and each appearing as an increasing union of finite subgroups, we have
$$L(\Gamma)\simeq L(\Gamma')$$
while the corresponding group algebras might not be isomorphic, $\mathbb C[\Gamma]\neq\mathbb C[\Gamma']$.
\end{proposition}

\begin{proof}
Here the first assertion follows from the above discusssion, the von Neumann algebra in question being the hyperfinite ${\rm II}_1$ factor $R$. As for the last assertion, there are countless counterexamples here, all coming from basic group theory.
\end{proof}

The point with the above result is that the isomorphisms of type $L(\Gamma)\simeq L(\Gamma')$ are in general impossible to prove with bare hands. Thus, we can see here the power of the Murray-von Neumann results in \cite{mv3}. And we can also see the magic of the weak topology, which by some kind of miracle, makes everyone equal in the end.

\section*{12b. Amenability}

The hyperfinite ${\rm II}_1$ factor $R$, which is a quite fascinating object, was heavily investigated by Murray-von Neumann \cite{mv3}, and then by Connes \cite{co2}. There are many things that can be said about it, which all interesting, but are usually quite technical as well. 

\bigskip

As a central result here, in what regards advanced hyperfiniteness theory, we have the following theorem of Connes \cite{co2}, whose proof is something remarkably heavy, and which is arguably the deepest result in operator algebra related functional analysis:

\begin{theorem}
For a finite von Neumann algebra $A$, the following are equivalent:
\begin{enumerate}
\item $A$ is hyperfinite in the usual sense, namely it appears as the weak closure of an increasing limit of finite dimensional algebras:
$$A=\overline{\bigcup_iA_i}^{\,w}$$

\item $A$ amenable, in the sense that the standard inclusion $A\subset B(H)$, with $H=L^2(A)$, admits a conditional expectation $E:B(H)\to A$.
\end{enumerate}
\end{theorem}

\begin{proof}
This result, due to Connes \cite{co2}, is something fairly heavy, that only a handful of people have really managed to understand, the idea being as follows:

\medskip

$(1)\implies(2)$ Assuming that the algebra $A$ is hyperfinite, let us write it as the weak closure of an increasing limit of finite dimensional subalgebras:
$$A=\overline{\bigcup_iA_i}^{\,w}$$

Consider the inclusion $A\subset B(H)$, with $H=L^2(A)$. In order to construct an expectation $E:B(H)\to A$, let us pick an ultrafilter $\omega$ on $\mathbb N$. Given $T\in B(H)$, we can define the following quantity, with $\mu_i$ being the Haar measure on the unitary group $U(A_i)$:
$$\psi(T)=\lim_{i\to\omega}\int_{U(A_i)}UTU^*\,d\mu_i(U)$$

With this construction made, by using now the standard involution $J:H\to H$, given by the formula $T\to T^*$, we can further define a map as follows:
$$E:B(H)\to A\quad,\quad 
E(T)=J\psi(T)J$$

But this is the expectation that we are looking for, with its left and right invariance properties coming from the left and right invariance of each Haar measure $\mu_i$.

\medskip

$(2)\implies(1)$ This is something heavy, using lots of advanced functional analysis, and for details here, we refer to Connes' original paper \cite{co2}.
\end{proof}

We should mention that Connes' results in \cite{co2}, besides proving the above implication $(2)\implies(1)$,  provide also a considerable extension of Theorem 2.10, with a number of further equivalent formulations of the notion of amenability, which are a bit more technical, but all good to know. The story here, still a bit simplified, is as follows:

\begin{fact}[Connes]
For a finite von Neumann algebra $A$, the following conditions are in fact equivalent:
\begin{enumerate}
\item $A$ is hyperfinite, in the sense that it appears as the weak closure of an increasing limit of finite dimensional algebras:
$$A=\overline{\bigcup_iA_i}^{\,w}$$

\item $A$ amenable, in the sense that the standard inclusion $A\subset B(H)$, with $H=L^2(A)$, admits a conditional expectation:
$$E:B(H)\to A$$

\item There exist unit vectors $\xi_n\in L^2(A)\otimes L^2(A)$ such that, for any $x\in A$:
$$||x\xi_n-\xi_nx||_2\to0\quad,\quad <x\xi_n,\xi_n>\to tr(x)$$

\item For any $x_1,\ldots,x_k\in A$ and $y_1,\ldots,y_k\in A$ we have:
$$\left|tr\left(\sum_ix_iy_i\right)\right|\leq\left|\left|\sum_ix_i\otimes y_i^{opp}\right|\right|_{min}$$
\end{enumerate}
\end{fact}

Again, this is something technical and advanced, that we won't get into, in this book. Let us mention however that the idea with all this is as follows:

\medskip

$(1)\implies(2)$ is elementary, as explained above.

\medskip

$(2)\implies(3)$ can be proved by using an inequality due to Powers-St\o rmer.

\medskip

$(3)\implies(4)$ is something quite technical, but doable as well. 

\medskip

$(4)\implies(2)$ is again something technical, but doable as well. 

\medskip

$(2)\implies(1)$ is, as before in Theorem 12.10, the difficult implication.

\medskip

Regarding the difficult implication, $(2)\implies(1)$, the difficulty here comes of course from the fact that, no matter what beautiful abstract functional analysis things you know about $A$, at some point you will have to get to work, and construct that finite dimensional subalgebras $A_i\subset A$, and it is not even clear where to start from. For a solution to this problem, and for more, we refer to Connes's article \cite{co2}, and also to his book \cite{co3}.

\bigskip

Getting back now to more everyday mathematics, the above results as stated remain something quite abstract, and advanced, and understanding their concrete implications will be our next task. In the case of the ${\rm II}_1$ factors, we have the following result:

\begin{theorem}
For a ${\rm II}_1$ factor $R$, the following are equivalent:
\begin{enumerate}
\item $R$ amenable, in the sense that we have an expectation, as follows:
$$E:B(L^2(R))\to R$$

\item $R$ is the Murray-von Neumann hyperfinite ${\rm II}_1$ factor.
\end{enumerate}
\end{theorem}

\begin{proof}
This follows indeed from Theorem 12.10, when coupled with the Murray-von Neumann uniqueness result for the hyperfinite ${\rm II}_1$ factor, from Theorem 12.6.
\end{proof}

As another application, getting back now to the general case, that of the finite von Neumann algebras, from Theorem 12.10 as stated, a first question is about how all this applies to the group von Neumann algebras, and more generally to the quantum group von Neumann algebras $L(\Gamma)$. In order to discuss this, let us start with the case of the usual discrete groups $\Gamma$. We will need the following result, which is standard:

\begin{theorem}
For a discrete group $\Gamma$, the following two conditions are equivalent, and if they are satisfied, we say that $\Gamma$ is amenable:
\begin{enumerate}
\item $\Gamma$ admits an invariant mean $m:l^\infty(\Gamma)\to\mathbb C$.

\item The projection map $C^*(\Gamma)\to C^*_{red}(\Gamma)$ is an isomorphism.
\end{enumerate}
Moreover, the class of amenable groups contains all the finite groups, all the abelian groups, and is stable under taking subgroups, quotients and products.
\end{theorem}

\begin{proof}
This is something very standard, the idea being as follows:

\medskip

(1) The equivalence $(1)\iff(2)$ is standard, with the amenability conditions (1,2) being in fact part of a much longer list of amenability conditions, including well-known criteria of F\o lner, Kesten and others. We will be back to this, with details, in a moment, directly in a more general setting, that of the discrete quantum groups.

\medskip

(2) As for the last assertion, regarding the finite groups, the abelian groups, and then the stability under taking subgroups, quotients and products, this is something elementary, which follows by using either of the above definitions of the amenability.
\end{proof}

Getting back now to operator algebras, we can complement Theorem 12.10 with:

\begin{theorem}
For a group von Neumann algebra $A=L(\Gamma)$, the following conditions are equivalent:
\begin{enumerate}
\item $A$ is hyperfinite.

\item $A$ amenable.

\item $\Gamma$ is amenable.
\end{enumerate}
\end{theorem}

\begin{proof}
The group von Neumann algebras $A=L(\Gamma)$ being by definition finite, Theorem 12.10 above applies, and gives the equivalence $(1)\iff(2)$. Thus, it remains to prove that we have $(2)\iff(3)$, and we can prove this as follows:

\medskip

$(2)\implies(3)$ This is something clear, because if we assume that $A=L(\Gamma)$ is amenable, we have by definition a conditional expectation $E:B(L^2(A))\to A$, and the restriction of this conditional expectation is the desired invariant mean $m:l^\infty(\Gamma)\to\mathbb C$.

\medskip

$(3)\implies(2)$ Assume that we are given a discrete amenable group $\Gamma$. In view of Theorem 12.13, this means that $\Gamma$ has an invariant mean, as follows:
$$m:l^\infty(\Gamma)\to\mathbb C$$

Consider now the Hilbert space $H=l^2(\Gamma)$, and for any operator $T\in B(H)$ consider the following map, which is a bounded sesquilinear form:
$$\varphi_T:H\times H\to\mathbb C$$
$$(\xi,\eta)\to m\left[\gamma\to<\rho_\gamma T\rho_\gamma^*\xi,\eta>\right]$$

By using the Riesz representation theorem, we conclude that there exists a certain operator $E(T)\in B(H)$, such that the following holds, for any two vectors $\xi,\eta$:
$$\varphi_T(\xi,\eta)=<E(T)\xi,\eta>$$ 

Summarizing, to any operator $T\in B(H)$ we have associated another operator, denoted $E(T)\in B(H)$, such that the following formula holds, for any two vectors $\xi,\eta$:
$$<E(T)\xi,\eta>=m\left[\gamma\to<\rho_\gamma T\rho_\gamma^*\xi,\eta>\right]$$ 

In order to prove now that this linear map $E$ is the desired expectation, observe that for any group element $g\in\Gamma$, and any two vectors $\xi,\eta\in H$, we have:
\begin{eqnarray*}
<\rho_gE(T)\rho_g^*\xi,\eta>
&=&<E(T)\rho_g^*\xi,\rho_g^*\eta>\\
&=&m\left[\gamma\to<\rho_\gamma T\rho_\gamma^*\rho_g^*\xi,\rho_g^*\eta>\right]\\
&=&m\left[\gamma\to<\rho_{g\gamma}T\rho_{g\gamma}^*\xi,\eta>\right]\\
&=&m\left[\gamma\to<\rho_\gamma T\rho_\gamma^*\xi,\eta>\right]\\
&=&<E(T)\xi,\eta>
\end{eqnarray*}

Since this is valid for any $\xi,\eta\in H$, we conclude that we have, for any $g\in\Gamma$:
$$\rho_gE(T)\rho_g^*=E(T)$$

But this shows that the element $E(T)\in B(H)$ is in the commutant of the right regular representation of $\Gamma$, and so belongs to the left regular group algebra of $\Gamma$:
$$E(T)\in L(\Gamma)$$

Summarizing, we have constructed a certain linear map $E:B(H)\to L(\Gamma)$. Now by using the above explicit formula of it, in terms of $m:l^\infty(\Gamma)\to\mathbb C$, which was assumed to be an invariant mean, we conclude that $E$ is indeed an expectation, as desired.
\end{proof}

As a very concrete application of all this technology, in relation now with the discrete group algebras which are ${\rm II}_1$ factors, the results that we have lead to:

\begin{theorem}
For a discrete group $\Gamma$, the following conditions are equivalent:
\begin{enumerate}
\item $\Gamma$ is amenable, and has the ICC property.

\item $A=L(\Gamma)$ is the hyperfinite ${\rm II}_1$ factor $R$.
\end{enumerate}
\end{theorem}

\begin{proof}
This follows indeed from Theorem 12.14, coupled with the standard fact, that we know well from chapter 10, that a group algebra $A=L(\Gamma)$ is a factor, and so a ${\rm II}_1$ factor, precisely when the group $\Gamma$ has the ICC property.
\end{proof}

As a comment here, this result, coming from Connes' Theorem 12.10, is far better than what we knew to come from Murray-von Neumann's Theorem 12.6, and with the statement itself being something elementary, not involving any kind of advanced functional analysis, such as the notion of amenability for von Neumann algebras. In fact, Murray-von Neumann knew about this statement, but their hunt for a proof proved to be unsuccessful, with the only possible proof being the one above, via advanced functional analysis. 

\bigskip

Summarizing, and to put things in context, Murray-von Neumann did great work with their papers \cite{mv1}, \cite{mv2}, \cite{mv3}, \cite{vn1}, \cite{vn2}, but were stuck with 3 questions, namely reduction theory, type ${\rm III}$ factors, and solutions of $L(\Gamma)=R$. And these questions were solved later by von Neumann himself \cite{vn3}, then Connes \cite{co1}, and Connes again \cite{co2}. 

\bigskip

Beautiful times these must have been, for mathematics and for mankind, and job for us, future generations, at least to write a complete von Neumann algebra book, clearly explaining all this fundamental material. And with many older people giving up, for various technical reasons, this will be most likely a job for you, young reader.

\section*{12c. Quantum groups}

Back now to work, we would like to discuss all sorts of questions, for the most open, or at least difficult, in relation with groups and quantum groups, taken finite, discrete or compact, and with more general quantum manifolds and quantum spaces, in connection with the Murray-von Neumann factor $R$, amenability and hyperfiniteness. As a first such question, in relation with the considerations from chapter 10, we would like to understand which discrete quantum groups $\Gamma$ produce group algebras as follows:
$$L(\Gamma)\simeq R$$

In terms of the compact quantum group duals $G=\widehat{\Gamma}$, the problem is that of understanding which compact quantum groups $G$ produce group algebras as follows:
$$L^\infty(G)\simeq R$$

In order to discuss this, we must first talk about amenability. We have here the following result, basically due to Woronowicz \cite{wo1}, and coming from the Peter-Weyl theory, extending to the discrete quantum groups the standard theory for discrete groups:

\index{amenability}

\begin{theorem}
Let $(A,u)$ with $u\in M_N(A)$ be a Woronowicz algebra, as axiomatized before. Let $A_{full}$ be the enveloping $C^*$-algebra of $\mathcal A=<u_{ij}>$, and let $A_{red}$ be the quotient of $A$ by the null ideal of the Haar integration. The following are then equivalent:
\begin{enumerate}
\item The Haar functional of $A_{full}$ is faithful.

\item The projection map $A_{full}\to A_{red}$ is an isomorphism.

\item The counit map $\varepsilon:A\to\mathbb C$ factorizes through $A_{red}$.

\item We have $N\in\sigma(Re(\chi_u))$, the spectrum being taken inside $A_{red}$.

\item $||ax_k-\varepsilon(a)x_k||\to0$ for any $a\in\mathcal A$, for certain norm $1$ vectors $x_k\in L^2(A)$.
\end{enumerate}
If this is the case, we say that the underlying discrete quantum group $\Gamma$ is amenable.
\end{theorem}

\begin{proof}
Before starting, we should mention that amenability and the present result are a bit like the spectral theorem, in the sense that knowing that the result formally holds does not help much, and in practice, one needs to remember the proof as well. For this reason, we will work out explicitly all the possible implications between (1-5), whenever possible, adding to the global formal proof, which will be linear, as follows:
$$(1)\implies(2)\implies(3)\implies(4)\implies(5)\implies(1)$$

In order to prove these implications, and the other ones too, the general idea is that this is is well-known in the group dual case, $A=C^*(\Gamma)$, with $\Gamma$ being a usual discrete group, and in general, the result follows by adapting the group dual case proof. 

\medskip

$(1)\iff(2)$ This follows from the fact that the GNS construction for the algebra $A_{full}$ with respect to the Haar functional produces the algebra $A_{red}$.

\medskip

$(2)\implies(3)$ This is trivial, because we have quotient maps $A_{full}\to A\to A_{red}$, and so our assumption $A_{full}=A_{red}$ implies that we have $A=A_{red}$. 

\medskip 

$(3)\implies(2)$ Assume indeed that we have a counit map, as follows:
$$\varepsilon:A_{red}\to\mathbb C$$

In order to prove $A_{full}=A_{red}$, we can use the right regular corepresentation. Indeed, we can define such a corepresentation by the following formula:
$$W(a\otimes x)=\Delta(a)(1\otimes x)$$

This corepresentation is unitary, so we can define a morphism as follows: 
$$\Delta':A_{red}\to A_{red}\otimes A_{full}\quad,\quad 
a\to W(a\otimes1)W^*$$

Now by composing with $\varepsilon\otimes id$, we obtain a morphism as follows:
$$(\varepsilon\otimes id)\Delta':A_{red}\to A_{full}\quad,\quad 
u_{ij}\to u_{ij}$$

Thus, we have our inverse for the canonical projection $A_{full}\to A_{red}$, as desired.

\medskip

$(3)\implies(4)$ This implication is clear, because we have:
\begin{eqnarray*}
\varepsilon(Re(\chi_u))
&=&\frac{1}{2}\left(\sum_{i=1}^N\varepsilon(u_{ii})+\sum_{i=1}^N\varepsilon(u_{ii}^*)\right)\\
&=&\frac{1}{2}(N+N)\\
&=&N
\end{eqnarray*}

Thus the element $N-Re(\chi_u)$ is not invertible in $A_{red}$, as claimed.

\medskip

$(4)\implies(3)$ In terms of the corepresentation $v=u+\bar{u}$, whose dimension is $2N$ and whose character is $2Re(\chi_u)$, our assumption $N\in\sigma(Re(\chi_u))$ reads:
$$\dim v\in\sigma(\chi_v)$$

By functional calculus the same must hold for $w=v+1$, and then once again by functional calculus, the same must hold for any tensor power of $w$:
$$w_k=w^{\otimes k}$$ 

Now choose for each $k\in\mathbb N$ a state $\varepsilon_k\in A_{red}^*$ having the following property:
$$\varepsilon_k(w_k)=\dim w_k$$

By Peter-Weyl we must have $\varepsilon_k(r)=\dim r$ for any $r\leq w_k$, and since any irreducible corepresentation appears in this way, the sequence $\varepsilon_k$ converges to a counit map: 
$$\varepsilon:A_{red}\to\mathbb C$$

$(4)\implies(5)$ Consider the following elements of $A_{red}$, which are positive:
$$a_i=1-Re(u_{ii})$$
 
Our assumption $N\in\sigma(Re(\chi_u))$ tells us that $a=\sum a_i$ is not invertible, and so there exists a sequence $x_k$ of norm one vectors in $L^2(A)$ such that: 
$$<ax_k,x_k>\to 0$$

Since the summands $<a_ix_k,x_k>$ are all positive, we must have, for any $i$:
$$<a_ix_k,x_k>\to0$$

We can go back to the variables $u_{ii}$ by using the following general formula:
$$||vx-x||^2=||vx||^2 +2<(1-Re(v))x,x>-1$$

Indeed, with $v=u_{ii}$ and $x=x_k$ the middle term on the right goes to 0, and so the whole term on the right becomes asymptotically negative, and so we must have:
$$||u_{ii}x_k-x_k||\to0$$

Now let $M_n(A_{red})$ act on $\mathbb C^n\otimes L^2(A)$. Since $u$ is unitary we have:
$$\sum_i||u_{ij}x_k||^2
=||u(e_j\otimes x_k)||
=1$$

From $||u_{ii}x_k||\to1$ we obtain $||u_{ij}x_k||\to0$ for $i\neq j$. Thus we have, for any $i,j$:
$$||u_{ij}x_k-\delta_{ij}x_k||\to0$$

Now by remembering that we have $\varepsilon(u_{ij})=\delta_{ij}$, this formula reads:
$$||u_{ij}x_k-\varepsilon(u_{ij})x_k||\to0$$

By linearity, multiplicativity and continuity, we must have, for any $a\in\mathcal  A$, as desired:
$$||ax_k-\varepsilon(a)x_k||\to0$$

$(5)\implies(1)$ This is something well-known, which follows via some standard functional analysis arguments, exactly as in the usual group case.

\medskip

$(1)\implies(5)$ Once again this is something well-known, which follows via some standard functional analysis arguments, exactly as in the usual group case.
\end{proof}

Before getting further, with advanced amenability and hyperfiniteness questions, and as a first application of the above, we can now advance on a problem that we left open before, in chapter 7, when talking about cocommutative Woronowicz algebras. Indeed, we can now state and prove the following result, which clarifies the situation:

\index{cocommutative algebra}

\begin{proposition}
The cocommutative Woronowicz algebras are the intermediate quotients of the following type, with $\Gamma=<g_1,\ldots,g_N>$ being a discrete group,
$$C^*(\Gamma)\to C^*_\pi(\Gamma)\to C^*_{red}(\Gamma)$$
and with $\pi$ being a unitary representation of $\Gamma$, subject to weak containment conditions of type $\pi\otimes\pi\subset\pi$ and $1\subset\pi$, which guarantee the existence of $\Delta,\varepsilon$.
\end{proposition}

\begin{proof}
We use the various findings from Theorem 12.16, following Woronowicz, the idea being to proceed in several steps, as follows:

\medskip

(1) Theorem 12.16 and standard functional analysis arguments show that the cocommutative Woronowicz algebras should appear as intermediate quotients, as follows:
$$C^*(\Gamma)\to A\to C^*_{red}(\Gamma)$$

(2) The existence of $\Delta:A\to A\otimes A$ requires our intermediate quotient to appear as follows, with $\pi$ being a unitary representation of $\Gamma$, satisfying the condition $\pi\otimes\pi\subset\pi$, taken in a weak containment sense, and with the tensor product $\otimes$ being taken here to be compatible with our usual maximal tensor product $\otimes$ for the $C^*$-algebras:
$$C^*(\Gamma)\to C^*_\pi(\Gamma)\to C^*_{red}(\Gamma)$$

(3) With this condition imposed, the existence of the antipode $S:A\to A^{opp}$ is then automatic, coming from the group antirepresentation $g\to g^{-1}$. 

\medskip

(4) The existence of the counit $\varepsilon:A\to\mathbb C$, however, is something non-trivial, related to amenability, and leading to a condition of type $1\subset\pi$, as in the statement.
\end{proof}  

Let us focus now on the Kesten amenability criterion, from Theorem 12.16 (4), which brings connections with interesting mathematics and physics, and which in practice will be our main amenability criterion. In order to discuss this, we will need:

\index{main character}

\begin{proposition}
Given a Woronowicz algebra $(A,u)$, with $u\in M_N(A)$, the moments of the main character $\chi=\sum_iu_{ii}$ are given by:
$$\int_G\chi^k=\dim\left(Fix(u^{\otimes k})\right)$$
In the case $u\sim\bar{u}$ the law of $\chi$ is a usual probability measure, supported on $[-N,N]$.
\end{proposition}

\begin{proof}
The first assertion follows from the Peter-Weyl theory, which tells us that we have the following formula, valid for any corepresentation $v\in M_n(A)$:
$$\int_G\chi_v=\dim(Fix(v))$$

Indeed, with $v=u^{\otimes k}$ we obtain the result. As for the second assertion, if we assume $u\sim\bar{u}$, then we have $\chi=\chi^*$, and so $law(\chi)$ is a real probability measure, supported by the spectrum of $\chi$. But, since the matrix $u\in M_N(A)$ is unitary, we have:
$$uu^*=1
\implies||u_{ij}||\leq 1,\forall i,j
\implies||\chi||\leq N$$

Thus the spectrum of the character satisfies $\sigma(\chi)\subset [-N,N]$, as desired. 
\end{proof}

In relation now with the notion of amenability, we have:

\index{Kesten amenability}
\index{amenability}

\begin{theorem}
A Woronowicz algebra $(A,u)$, with $u\in M_N(A)$, is amenable when
$$N\in supp\Big(law(Re(\chi))\Big)$$
and the support on the right depends only on $law(\chi)$.
\end{theorem}

\begin{proof}
There are two assertions here, the proof being as follows:

\medskip

(1) According to the Kesten amenability criterion, from Theorem 12.16 (4), the algebra $A$ is amenable when the following condition is satisfied:
$$N\in\sigma(Re(\chi))$$

Now since $Re(\chi)$ is self-adjoint, we know from spectral theory that the support of its spectral measure $law(Re(\chi))$ is precisely its spectrum $\sigma(Re(\chi))$, as desired:
$$supp(law(Re(\chi)))=\sigma(Re(\chi))$$

(2) Regarding the second assertion, once again the variable $Re(\chi)$ being self-adjoint, its law depends only on the moments $\int_GRe(\chi)^p$, with $p\in\mathbb N$. But, we have:
$$\int_GRe(\chi)^p
=\int_G\left(\frac{\chi+\chi^*}{2}\right)^p
=\frac{1}{2^p}\sum_{|k|=p}\int_G\chi^k$$

Thus $law(Re(\chi))$ depends only on $law(\chi)$, and this gives the result.
\end{proof}

Let us work out now in detail the group dual case. Here we obtain a very interesting measure, called Kesten measure of the group, as follows:

\index{Cayley graph}

\begin{proposition}
In the case $A=C^*(\Gamma)$ and $u=diag(g_1,\ldots,g_N)$, and with the normalization $1\in u=\bar{u}$ made, we have the formula
$$\int_{\widehat{\Gamma}}\chi^p=\#\left\{i_1,\ldots,i_p\Big|g_{i_1}\ldots g_{i_p}=1\right\}$$
counting the loops based at $1$, having length $p$, on the corresponding Cayley graph.
\end{proposition}

\begin{proof}
Consider indeed a discrete group $\Gamma=<g_1,\ldots,g_N>$. The main character of $A=C^*(\Gamma)$, with fundamental corepresentation $u=diag(g_1,\ldots,g_N)$, is then:
$$\chi=g_1+\ldots+g_N$$

Given a colored integer $k=e_1\ldots e_p$, the corresponding moment is given by:
$$\int_{\widehat{\Gamma}}\chi^k
=\int_{\widehat{\Gamma}}(g_1+\ldots+g_N)^k
=\#\left\{i_1,\ldots,i_p\Big|g_{i_1}^{e_1}\ldots g_{i_p}^{e_p}=1\right\}$$

In the self-adjoint case now, $u\sim\bar{u}$, as in the statement, we are only interested in the moments with respect to usual integers, $p\in\mathbb N$, and the above formula becomes:
$$\int_{\widehat{\Gamma}}\chi^p=\#\left\{i_1,\ldots,i_p\Big|g_{i_1}\ldots g_{i_p}=1\right\}$$

Assume now that we have in addition $1\in u$, so that the condition $1\in u=\bar{u}$ in the statement is satisfied. At the level of the generating set $S=\{g_1,\ldots,g_N\}$ this means:
$$1\in S=S^{-1}$$

Thus the corresponding Cayley graph is well-defined, with the elements of $\Gamma$ as vertices, and with the edges $g-h$ appearing when the following condition is satisfied:
$$gh^{-1}\in S$$

A loop on this graph based at 1, having length $p$, is then a sequence as follows:
$$(1)-(g_{i_1})-(g_{i_1}g_{i_2})-\ldots-(g_{i_1}\ldots g_{i_{p-1}})-(g_{i_1}\ldots g_{i_p}=1)$$

Thus the moments of $\chi$ count indeed such loops, as claimed.
\end{proof}

In order to generalize the above result to arbitrary Woronowicz algebras, we can use the discrete quantum group philosophy. The fundamental result here is as follows:

\index{Cayley graph}

\begin{theorem}
Let $(A,u)$ be a Woronowicz algebra, and assume, by enlarging if necessary $u$, that we have $1\in u=\bar{u}$. The following formula
$$d(v,w)=\min\left\{k\in\mathbb N\Big|1\subset\bar{v}\otimes w\otimes u^{\otimes k}\right\}$$
defines then a distance on $Irr(A)$, which coincides with the geodesic distance on the associated Cayley graph. In the group dual case we obtain the usual distance.
\end{theorem}

\begin{proof}
The fact that the lengths are finite follows from Woronowicz's analogue of Peter-Weyl theory, and the other verifications are as follows:

\medskip

(1) The symmetry axiom is clear.

\medskip

(2) The triangle inequality is elementary to establish as well. 

\medskip

(3) The Cayley graph assertion is something elementary as well.

\medskip

(4) Finally, in the group dual case, where our Woronowicz algebra is of the form $A=C^*(\Gamma)$, with $\Gamma=<S>$ being a finitely generated discrete group, our normalization condition $1\in u=\bar{u}$ means that the generating set must satisfy:
$$1\in S=S^{-1}$$

But this is precisely the normalization condition for the discrete groups, and the fact that we obtain the same metric space is clear.
\end{proof}

Summarizing, we have a good understanding of what a discrete quantum group is. We can now formulate a generalization of Proposition 12.20, as follows:

\begin{theorem}
Let $(A,u)$ be a Woronowicz algebra, with the normalization assumption $1\in u=\bar{u}$ made. The moments of the main character, 
$$\int_G\chi^p=\dim\left(Fix(u^{\otimes p})\right)$$
count then the loops based at $1$, having lenght $p$, on the corresponding Cayley graph.
\end{theorem}

\begin{proof}
Here the formula of the moments, with $p\in\mathbb N$, is the one coming from Proposition 12.18, and the Cayley graph interpretation comes from Theorem 12.21.
\end{proof}

As an application of this, we can introduce the notion of growth, as follows:

\index{growth}

\begin{definition}
Given a closed subgroup $G\subset U_N^+$, with $1\in u=\bar{u}$, consider the series  whose coefficients are the ball volumes on the corresponding Cayley graph,
$$f(z)=\sum_kb_kz^k\quad,\quad 
b_k=\sum_{l(v)\leq k}\dim(v)^2$$
and call it growth series of the discrete quantum group $\widehat{G}$. In the group dual case, $G=\widehat{\Gamma}$, we obtain in this way the usual growth series of $\Gamma$. 
\end{definition}

There are many things that can be said about the growth, and we will be back to this. As a first such result, in relation with the notion of amenability, we have:

\index{polynomial growth}
\index{amenability}

\begin{theorem}
Polynomial growth implies amenability.
\end{theorem}

\begin{proof}
We recall from Theorem 12.21 that the Cayley graph of $\widehat{G}$ has by definition the elements of $Irr(G)$ as vertices, and the distance is as follows:
$$d(v,w)=\min\left\{k\in\mathbb N\Big|1\subset\bar{v}\otimes w\otimes u^{\otimes k}\right\}$$

By taking $w=1$ and by using Frobenius reciprocity, the lenghts are given by:
$$l(v)=\min\left\{k\in\mathbb N\Big|v\subset u^{\otimes k}\right\}$$

By Peter-Weyl we have then a decomposition as follows, where $B_k$ is the ball of radius $k$, and where $m_k(v)\in\mathbb N$ are certain multiplicities:
$$u^{\otimes k}=\sum_{v\in B_k}m_k(v)\cdot v$$

By using now Cauchy-Schwarz, we obtain the following inequality:
\begin{eqnarray*}
m_{2k}(1)b_k
&=&\sum_{v\in B_k}m_k(v)^2\sum_{v\in B_k}\dim(v)^2\\
&\geq&\left(\sum_{v\in B_k}m_k(v)\dim(v)\right)^2\\
&=&N^{2k}
\end{eqnarray*}

But shows that if $b_k$ has polynomial growth, then the following happens:
$$\limsup_{k\to\infty}\, m_{2k}(1)^{1/2k}\geq N$$

Thus, the Kesten type criterion applies, and gives the result.
\end{proof}

There are many other things that can be said, as a continuation of the above, notably with explicit computations of growth exponents for all the discrete quantum groups that we know, and with some further generalities too, of functional analytic nature, and in relation with Lie theory and its generalizations too. For more on all this, you can check my quantum group textbook \cite{ba3}, and the quantum group literature cited there.

\bigskip

To summarize now, we have a decent understanding of what a discrete quantum group is, and also of what amenability means, in the discrete quantum group setting. However, all this does not exactly solve the von Neumann algebra questions, and we have:

\begin{question}
Which discrete quantum groups $\Gamma$ have the property $L(\Gamma)\simeq R$? Equivalently, which compact quantum groups $G$ have the property $L^\infty(G)\simeq R$?
\end{question}

Here the equivalence between the above two questions comes from the fact that, with $\Gamma=\widehat{G}$, we have $L(\Gamma)=L^\infty(G)$. As for the questions themselves, normally the hyperfiniteness part can be dealt with as in the classical group case, by using the amenability theory developed above, and the problem is with the ICC property, guaranteeing factoriality, with no one presently knowing what this ``quantum ICC'' property is.

\bigskip

As a funny comment here, the equation $L(\Gamma)\simeq R$ is precisely the one Murray and von Neumann were stuck with, in the classical group case, some 90 years ago. Some sort of Connes is needed, coming and solving this problem, with new ideas.

\bigskip

Finally, let us mention that in connection with amenability and hyperfiniteness, we have as well a series of further questions, in relation with the actions of quantum groups. To be more precise, the problems that we would like to solve are as follows:

\bigskip

(1) We would like to understand, given a compact group or quantum group acting on a von Neumann algebra, $G\curvearrowright P$, when the fixed point algebra $P^G$ is a factor. 

\bigskip

(2) More generally, we would like to understand under which assumptions on $G\curvearrowright P$ the fixed point algebra $(B\otimes P)^G$ is a factor, for any finite dimensional algebra $B$.

\bigskip

(3) In fact, we would like to understand when the fixed point algebra $P^G$, or more generally all the fixed point algebras $(B\otimes P)^G$, are the hyperfinite ${\rm II}_1$ factor $R$.

\bigskip

These questions are all of interest in subfactor theory, the idea being that a quite standard construction of subfactors is $(B_0\otimes P)^G\subset(B_1\otimes P)^G$, coming from a von Neumann algebra $P$, an inclusion of finite dimensional algebras $B_0\subset B_1$, and a compact quantum group $G$ acting on everything, provided that the fixed point algebras involved are indeed factors. And then, once such a subfactor constructed and studied, the main problem is that of understanding if this subfactor can be taken to be hyperfinite.

\bigskip

These are quite technical questions, to be discussed in chapters 13-16 below, when doing subfactors. Let us mention however, coming a bit in advance, that we have:

\begin{fact}
Assuming that $\Gamma=\widehat{G}$ has an outer action on the hyperfinite ${\rm II}_1$ factor
$$\Gamma\curvearrowright R$$
we can set $P=R\rtimes\Gamma$, and the answer to the above questions is yes.
\end{fact}

Which brings us into the very interesting question on whether we have such outer actions $\Gamma\curvearrowright R$, with the status of the subject being as follows:

\bigskip

(1) All this goes back to work in the 80s of Ocneanu, and Wassermann too, with Ocneanu eventually conjecturing that any discrete group $\Gamma$, and more generally any discrete quantum group $\Gamma$, should have such an action. This question is still open.

\bigskip

(2) In practice, the result is known in the finite case, $|\Gamma|<\infty$, and more generally in the case where $C^*(\Gamma)$ has an inner faithful matrix model, in the sense of chapter 11, with this being worked out in \cite{ba1} and its follow-ups, and then by Vaes in \cite{vae}.

\bigskip

(3) And there has been quite some work on this, since then. For the status of the question, and relations with other questions, such as the Connes embedding problem, Voiculescu microstates and more, we refer to Brannan-Chirvasitu-Freslon \cite{bcf}.

\bigskip

Summarizing, many things going on here, with the philosophy being somehow that, once we want our factors or subfactors to be hyperfinite, isomorphic to $R$, we are all of the sudden into all sorts of interesting questions, in relation with advanced mathematics and physics. But more on this later, in chapters 13-16 below, when doing subfactors.

\section*{12d. Hyperfinite factors}

Back to general theory, there are many other things that can be said, in relation with hyperfiniteness. We first have a reduction theory result, as follows:

\begin{theorem}
Any tracial hyperfinite von Neumann algebra appears as
$$A=\int_XA_x\,dx$$
with the factors $A_x$ being either usual matrix algebras, or the factor $R$.
\end{theorem}

\begin{proof}
This follows indeed by combining the von Neumann reduction theory from \cite{vn3} with the theory of $R$ of Murray-von Neumann \cite{mv3} and Connes \cite{co2}.
\end{proof}

More generally, we have the following result, this time in arbitrary type:

\index{hyperfinite algebra}
\index{reduction theory}

\begin{theorem}
Given a hyperfinite von Neumann algebra $A\subset B(H)$, write its center as follows, with $X$ being a measured space:
$$Z(A)=L^\infty(X)$$ 
The whole algebra $A$ decomposes then over this measured space $X$, as follows,
$$A=\int_XA_x\,dx$$
with the fibers $A_x$ being hyperfinite von Neumann factors, which can be of type ${\rm I},{\rm II},{\rm III}$.
\end{theorem}

\begin{proof}
This is again something heavy, combining the general reduction theory results of von Neumann with the work of Connes in the hyperfinite case.
\end{proof}

Which brings us into the question of classifying all hyperfinite factors. The result here, due to Connes \cite{co2}, with a key contribution by Haagerup \cite{ha1}, is as follows:

\index{hyperfinite factor}
\index{R}

\begin{theorem}
The hyperfinite factors are as follows, with $1$ factor in each class
$${\rm I_N},{\rm I}_\infty$$
$${\rm II}_1,{\rm II}_\infty$$
$${\rm III}_0,{\rm III}_\lambda,{\rm III}_1$$
and with the type ${\rm II}_1$ one $R$ being the most important, basically producing the others too.
\end{theorem}

\begin{proof}
This is again heavy, based on early work of Murray-von Neumann in type II \cite{mv3}, then on heavy work by Connes in type II and III \cite{co1}, \cite{co2}, basically finishing the classification, and with a final contribution by Haagerup in type ${\rm III}_1$ \cite{ha1}.
\end{proof}

Getting back now to the ${\rm II}_1$ factors, and beyond hyperfiniteness, where things are understood, with $R$ being the only example, there is a whole classification program here, by Popa and others, going on. Let us mention that a main open problem is that of deciding whether the free group factors are isomorphic or not:

\begin{question}
Are the von Neumann algebras of free groups isomorphic,
$$L(F_N)\simeq L(F_M)$$
for $N\neq M$, or not?
\end{question}

This question can be of course asked in crossed product form, in the spirit of the various crossed product results evoked above, and of advanced ergodic theory in general, with the space in question, producing the crossed product, being the point:
$$\{.\}\rtimes F_N\simeq^?\{.\}\rtimes F_M$$

This formulation, used by Popa, has the advantage of putting the above problem into a more conceptual framework, with lots of advanced machinery available around. However, it is not clear whether this formulation simplifies or not the original problem.

\bigskip

There are as well a number of alternative approaches to this question, and notably the Voiculescu one, using free probability, which is particularly conceptual and beautiful, the idea being that of recapturing the number $N\in\mathbb N$ from the knowledge of the von Neumann algebra $L(F_N)$, via an entropy-type invariant:
$$L(F_N)\rightsquigarrow N$$

This latter program, while not solving the original problem, due to technical difficulties, is however very successful, in the sense that it has led to a lot of interesting results and computations, in relation with a lot of mathematics and physics.

\bigskip

Is the free group factor problem something belonging to logic, as the difficult problems in functional analysis usually end up being? No one really knows the answer here.

\bigskip

Interestingly, the question is difficult to the point where the conjectural answer, yes or no, is not known. And even worse, excluding the many people who have spent considerable time on the matter, years or more, working on yes or no, most people familiar with the question don't even really know what to wish for, yes or no, as an answer.

\bigskip

In what concerns us, we have been quite close in this book to the ideas of Voiculescu, but, as a surprise, these very ideas of Voiculescu lead us into wishing for a yes answer to the above question, which is opposite to his no wish, and work using free entropy. Indeed, to put things in context, let us formulate the question in the following way:

\begin{question}
Is there a factor $F$, standing as a free counterpart for $R$?
\end{question}

And wouldn't you wish for a yes answer to this question, with $F$ being of course all the free group factors $L(F_N)$ combined, and probably many more, coming from all sorts of free quantum groups, free homogeneous spaces, or other free manifolds. It would be good to know in free geometry that what we get by default is this factor $F$.

\bigskip

As a last comment here, later on, when doing subfactors, we will see that the particular factor $F=L(F_\infty)$ quite does the job there, in subfactors, being more of less the only ``free factor'' that is needed, for that theory. But this does not really solve Question 12.31 in the context of subfactor theory because, ironically, the main questions there, including the ``free'' ones, rather concern the subfactors of the good old hyperfinite factor $R$.

\section*{12e. Exercises} 

Things have been very technical in this chapter, which was more of a survey than something else, and as a unique exercise on all this, we have:

\begin{exercise}
Learn some more basic hyperfinite von Neumann algebra theory, from the papers of von Neumann and Murray-von Neumann, then Connes, and then Haagerup and others, and write down a brief account of what you learned.
\end{exercise}

In what follows we will avoid ourselves this type of exercise, basically by getting back to the material in chapter 10, and building on that, following Jones.

\part{Subfactor theory}

\ \vskip50mm

\begin{center}
{\em Maria, you've gotta see her

Go insane and out of your mind

Latina, Ave Maria

A million and one candle lights}
\end{center}

\chapter{Subfactor theory}

\section*{13a. The Jones tower}

In this last part of the present book we discuss the basics of Jones' subfactor theory \cite{jo1}, \cite{jo2}, \cite{jo3}, \cite{jo4}, \cite{jo5}. The idea is that subfactors are quite subtle objects, generalizing various algebraic and combinatorial constructions from chapters 5-8, and coming from the functional analysis and operator theory considerations from chapters 9-12. Their study will bring us into a lot of advanced mathematics, mixing algebra, geometry, analysis and probability, and with everything being of modern physics flavor, often in relation with considerations from advanced statistical mechanics, and quantum mechanics.

\bigskip

We recall that a ${\rm II}_1$ factor is a von Neumann algebra $A\subset B(H)$ which has trivial center, $Z(A)=\mathbb C$, is infinite dimensional, and has a trace $tr:A\to\mathbb C$. For a number of reasons, ranging from simple and intuitive to fairly advanced, explained in chapters 9-12, such algebras are the core at the whole von Neumann algebra theory.

\bigskip

The world of ${\rm II}_1$ factors is a bit similar to the world of the usual matrix algebras $M_N(\mathbb C)$, which are actually called type ${\rm I}$ factors, in the sense that it is ``self-sufficient'', with no need to go further than that. In particular, a nice representation theory for such ${\rm II}_1$ factors can be obtained by staying inside the class of ${\rm II}_1$ factors, and we have the following definition to start with, which will keep us busy for the rest of this book:

\index{subfactor}

\begin{definition}
A subfactor is an inclusion of ${\rm II}_1$ factors $A\subset B$.
\end{definition}

We will see later some examples of such inclusions, along with motivations for their study. In order to get started now, the first thing to be done with such an inclusion is that of defining its index, as a quantity of the following type:
$$[B:A]=\dim_AB$$

Since both $A,B$ are infinite dimensional algebras, this is not exactly obvious. In addition, in view of our previous experience with the ${\rm II}_1$ factors, and notably with their ``continuous dimension'' features, we can only expect the index to range as follows:
$$[B:A]\in[1,\infty]$$

In order to discuss this, let us recall from chapter 10 that given a representation of a ${\rm II}_1$ factor $A\subset B(H)$, we can construct a number as follows, called coupling constant, which for the standard form, where $H=L^2(A)$, takes the value $1$, and which in general mesures how far is $A\subset B(H)$ from the standard form:
$$\dim_AH\in(0,\infty]$$

Getting now to the subfactors, in the sense of Definition 13.1, we have the following construction, that we know as well from chapter 10:

\index{index of subfactor}

\begin{theorem}
Given a subfactor $A\subset B$, the number
$$N=\frac{\dim_AH}{\dim_BH}\in[1,\infty]$$
is independent of the ambient Hilbert space $H$, and is called index.
\end{theorem}

\begin{proof}
This is something that we know from chapter 10, the idea being that the independence of the index from the choice of the ambient Hilbert space $H$ comes from the various basic properties of the coupling constant.
\end{proof}

There are many examples of subfactors, and we will discuss this gradually, in what follows. Following Jones \cite{jo1}, the most basic examples of subfactors are as follows:

\index{Jones subfactor}

\begin{proposition}
Assuming that $G$ is a compact group, acting on a ${\rm II}_1$ factor $P$ in a minimal way, in the sense that we have
$$(P^G)'\cap P=\mathbb C$$
and that $H\subset G$ is a closed subgroup of finite index, we have a subfactor
$$P^G\subset P^H$$
having index $N=[G:H]$, called Jones subfactor.
\end{proposition}

\begin{proof}
This is something standard, the idea being that the factoriality of $P^G,P^H$ comes from the minimality of the action, and that the index formula is clear. We will be back with full details about this in a moment, directly in a more general setting.
\end{proof}

In order to study the subfactors, let us start with the following standard result:

\index{conditional expectation}

\begin{proposition}
Given a subfactor $A\subset B$, there is a unique linear map
$$E:B\to A$$
which is positive, unital, trace-preserving and satisfies the following condition:
$$E(b_1ab_2)=b_1E(a)b_2$$
This map is called conditional expectation from $B$ onto $A$.
\end{proposition}

\begin{proof}
We make use of the standard representation of the ${\rm II}_1$ factor $B$, with respect to its unique trace $tr:B\to\mathbb C$, as constructed in chapter 10:
$$B\subset L^2(B)$$

If we denote by $\Omega$ the standard cyclic and separating vector of $L^2(B)$, we have an identification $A\Omega=L^2(A)$. Consider now the following orthogonal projection:
$$e:L^2(B)\to L^2(A)$$

It follows from definitions that we have an inclusion as follows:
$$e(B\Omega)\subset A\Omega$$ 

Thus $e$ induces by restriction a certain linear map $E:B\to A$. This linear map $E$ and the orthogonal projection $e$ are then related by:
$$exe=E(x)e$$

But this shows that the linear map $E$ satisfies the various conditions in the statement, namely positivity, unitality, trace preservation and bimodule property. As for the uniqueness assertion, this follows by using the same argument, applied backwards, the idea being that a map $E$ as in the statement must come from the projection $e$.
\end{proof}

Following Jones \cite{jo1}, we will be interested in what follows in the orthogonal projection $e:L^2(B)\to L^2(A)$ producing the expectation $E:B\to A$, rather than in $E$ itself:

\index{Jones projection}

\begin{definition}
Associated to any subfactor $A\subset B$ is the orthogonal projection
$$e:L^2(B)\to L^2(A)$$
producing the conditional expectation $E:B\to A$ via the following formula:
$$exe=E(x)e$$
This projection is called Jones projection for the subfactor $A\subset B$.
\end{definition}

Quite remarkably, the subfactor $A\subset B$, as well as its commutant, can be recovered from the knowledge of this projection, in the following way:

\begin{proposition}
Given a subfactor $A\subset B$, with Jones projection $e$, we have
$$A=B\cap\{e\}'$$
$$A'=(B'\cap\{e\})''$$
as equalities of von Neumann algebras, acting on the space $L^2(B)$.
\end{proposition}

\begin{proof}
These formulae basically follow from $exe=E(x)e$, as follows:

\medskip

(1) Let us first prove that we have $A\subset B\cap\{e\}'$. Given $x\in A$, we have:
$$xe=E(x)e=exe$$
$$x^*e=E(x^*)e=ex^*e$$

Thus, we obtain, as desired, that $x$ commutes with $e$:
$$ex
=(x^*e)^*
=(ex^*e)^*
=exe
=xe$$

(2) Let us prove now that $B\cap\{e\}'\subset A$. Assuming $ex=xe$, we have:
$$E(x)e=exe=xe^2=xe$$

We conclude from this that we have the following equality:
$$(E(x)-x)\Omega=(E(x)-x)e\Omega=0$$

Now since $\Omega$ is separating for $B$ we have, as desired:
$$x=E(x)\in A$$

(3) In order to prove now $A'=<B',e>$, observe that we have:
$$A
=B\cap\{e\}'
=B''\cap\{e\}'
=(B'\cap\{e\})'$$

Now by taking the commutant, we obtain $A'=(B'\cap\{e\})''$, as desired.
\end{proof}

Still following Jones \cite{jo1}, we are now ready to formulate a key definition:

\index{Jones projection}
\index{basic construction}

\begin{definition}
Associated to any subfactor $A\subset B$ is the basic construction
$$A\subset_eB\subset C$$
with $C=<B,e>$ being the algebra generated by $B$ and by the Jones projection
$$e:L^2(B)\to L^2(A)$$
acting on the Hilbert space $L^2(B)$.
\end{definition}

The idea in what follows will be that $B\subset C$ appears as a kind of ``reflection'' of $A\subset B$, and also that the basic construction can be iterated, with all this leading to nontrivial results. Let us start by further studying the basic construction:

\begin{theorem}
Given a subfactor $A\subset B$ having finite index, 
$$[B:A]<\infty$$
the basic construction $A\subset_eB\subset C$ has the following properties:
\begin{enumerate}
\item $C=JA'J$.

\item $C=\overline{B+Beb}$.

\item $C$ is a ${\rm II}_1$ factor.

\item $[C:B]=[B:A]$.

\item $eCe=Ae$.

\item $tr(e)=[B:A]^{-1}$.

\item $tr(xe)=tr(x)[B:A]^{-1}$, for any $x\in B$.
\end{enumerate}
\end{theorem}

\begin{proof}
All this is standard, the idea being as follows:

\medskip

(1) We have $JB'J=B$ and $JeJ=e$, which gives:
\begin{eqnarray*}
JA'J
&=&J<B',e>J\\
&=&<JB'J,JeJ>\\
&=&<B,e>\\
&=&C
\end{eqnarray*}

(2) This follows from the fact that the vector space $B+BeB$ is closed under multiplication, and from the fact that we have $exe=E(x)e$.

\medskip

(3) This follows from the fact, that we know from chapter 10, that our finite index assumption $[B:A]<\infty$ is equivalent to the fact that $A'$ is a factor. But this is in turn  equivalent to the fact that $C=JA'J$ is a factor, as desired.

\medskip

(4) We have indeed the folowing computation:
\begin{eqnarray*}
[C:B]
&=&\frac{\dim_BL^2(B)}{\dim_CL^2(B)}\\
&=&\frac{1}{\dim_CL^2(B)}\\
&=&\frac{1}{\dim_{JA'J}L^2(B)}\\
&=&\frac{1}{\dim_{A'}L^2(B)}\\
&=&\dim_AL^2(B)\\
&=&[B:A]
\end{eqnarray*}

(5) This follows indeed from (2) and from the formula $exe=E(x)e$.

\medskip

(6) We have the following computation:
\begin{eqnarray*}
1
&=&\dim_AL^2(A)\\
&=&\dim_A(eL^2(B))\\
&=&tr_{A'}(e)\dim_A(L^2(B))\\
&=&tr_{A'}(a)[B:A]
\end{eqnarray*}

Now since $C=JA'J$ and $JeJ=e$, we obtain from this, as desired:
\begin{eqnarray*}
tr(e)
&=&tr_{JA'J}(JeJ)\\
&=&tr_{A'}(e)\\
&=&[B:A]^{-1}
\end{eqnarray*}

(7) We already know from (6) that the formula in the statement holds for $x=1$. In order to discuss the general case, observe first that for $x,y\in A$ we have:
$$tr(xye)=tr(yex)=tr(yxe)$$

Thus the linear map $x\to tr(xe)$ is a trace on $A$, and by uniqueness of the trace on $A$, we must have, for a certain constant $c>0$:
$$tr(xe)=c\cdot tr(x)$$

Now by using (6) we obtain $c=[B:A]^{-1}$, so we have proved the formula in the statement for $x\in A$. The passage to the general case $x\in B$ can be done as follows:
\begin{eqnarray*}
tr(xe)
&=&tr(exe)\\
&=&tr(E(x)e)\\
&=&tr(E(x))c\\
&=&tr(x)c
\end{eqnarray*}

Thus, we have proved the formula in the statement, in general.
\end{proof}

The above result is quite interesting, so let us perform now twice the basic construction, and see what we get. The result here, which is more technical, is as follows:

\begin{proposition}
Associated to $A\subset B$ is the double basic construction
$$A\subset_eB\subset_fC\subset D$$
with $e,f$ being the following orthogonal projections,
$$e:L^2(B)\to L^2(A)$$
$$f:L^2(C)\to L^2(B)$$
having the following properties:
$$fef=[B:A]^{-1}f$$
$$efe=[B:A]^{-1}e$$
\end{proposition}

\begin{proof}
We have two formulae to be proved, the idea being as follows:

\medskip

(1) The first formula is clear, because we have:
\begin{eqnarray*}
fef
&=&E(e)f\\
&=&tr(e)f\\
&=&[B:A]^{-1}f
\end{eqnarray*}

(2) Regarding now the second formula, it is enough to check it on the dense subset $(B+BeB)\Omega$. Thus, we must show that for any $x,y,z\in B$, we have:
$$efe(x+yez)\Omega=[B:A]^{-1}e(x+yez)\Omega$$

For this purpose, we will prove that we have, for any $x,y,z\in B$:
$$efex\Omega=[B:A]^{-1}ex\Omega$$
$$efeyez\Omega=[B:A]^{-1}eyez\Omega$$

The first formula can be established as follows:
\begin{eqnarray*}
efex\Omega
&=&efexf\Omega\\
&=&eE(ex)f\Omega\\
&=&eE(e)xf\Omega\\
&=&[B:A]^{-1}exf\Omega\\
&=&[B:A]^{-1}ex\Omega
\end{eqnarray*}

The second formula can be established as follows:
\begin{eqnarray*}
efeyez\Omega
&=&efeyezf\Omega\\
&=&eE(eyez)f\Omega\\
&=&eE(eye)zf\Omega\\
&=&eE(E(y)e)zf\Omega\\
&=&eE(y)E(e)zf\Omega\\
&=&[B:A]^{-1}eE(y)zf\Omega\\
&=&[B:A]^{-1}eyezf\Omega\\
&=&[B:A]^{-1}eyez\Omega
\end{eqnarray*}

Thus, we are led to the conclusion in the statement.
\end{proof}

We can in fact perform the basic construction by recurrence, and we obtain:

\index{basic construction}
\index{Jones tower}
\index{Jones projection}

\begin{theorem}
Associated to any subfactor $A_0\subset A_1$ is the Jones tower
$$A_0\subset_{e_1}A_1\subset_{e_2}A_2\subset_{e_3}A_3\subset\ldots\ldots$$
with the Jones projections having the following properties:
\begin{enumerate}
\item $e_i^2=e_i=e_i^*$.

\item $e_ie_j=e_je_i$ for $|i-j|\geq2$.

\item $e_ie_{i\pm1}e_i=[B:A]^{-1}e_i$.

\item $tr(we_{n+1})=[B:A]^{-1}tr(w)$, for any word $w\in<e_1,\ldots,e_n>$.
\end{enumerate}
\end{theorem}

\begin{proof}
This follows from Theorem 13.8 and Proposition 13.9, because the triple basic construction does not need in fact any further study. See Jones \cite{jo1}.
\end{proof}

\section*{13b. Temperley-Lieb}

The relations found in Theorem 13.10 are in fact well-known, from the standard theory of the Temperley-Lieb algebra. This algebra, discovered by Temperley and Lieb in the context of statistical mechanics \cite{tli}, has a very simple definition, as follows:

\index{Temperley-Lieb}
\index{noncrossing pairings}

\begin{definition}
The Temperley-Lieb algebra of index $N\in[1,\infty)$ is defined as
$$TL_N(k)=span(NC_2(k,k))$$
with product given by vertical concatenation, with the rule
$$\bigcirc=N$$
for the closed circles that might appear when concatenating.
\end{definition}

In other words, the algebra $TL_N(k)$, depending on parameters $k\in\mathbb N$ and $N\in[1,\infty)$, is the formal linear span of the pairings $\pi\in NC_2(k,k)$. The product operation is obtained by linearity, for the pairings which span $TL_N(k)$ this being the usual vertical concatenation, with the conventions that things go ``from top to bottom'', and that each circle that might appear when concatenating is replaced by a scalar factor, equal to $N$.

\bigskip

In order to make the connection with subfactors, let us start with:

\begin{proposition}
The Temperley-Lieb algebra $TL_N(k)$ is generated by the diagrams
$$\varepsilon_1={\ }^\cup_\cap\quad,\quad 
\varepsilon_2=|\!{\ }^\cup_\cap\quad,\quad 
\varepsilon_3=||\!{\ }^\cup_\cap\quad,\quad 
\ldots$$
which are all multiples of projections, in the sense that their rescaled versions
$$e_i=N^{-1}\varepsilon_i$$
satisfy the abstract projection relations $e_i^2=e_i=e_i^*$.
\end{proposition}

\begin{proof}
We have two assertions here, the idea being as follows:

\medskip

(1) The fact that the algebra $TL_N(k)$ is indeed generated by the sequence of diagrams $\varepsilon_1,\varepsilon_2,\varepsilon_3,\ldots$ follows by drawing pictures, and more specifically by graphically decomposing each basis element $\pi\in NC_2(k,k)$ as a product of such elements $\varepsilon_i$.

\medskip

(2) Regarding now the projection assertion, when composing $\varepsilon_i$ with itself we obtain $\varepsilon_i$ itself, times a circle. Thus, according to our multiplication conventions, we have:
$$\varepsilon_i^2=N\varepsilon_i$$

Also, when turning upside-down $\varepsilon_i$, we obtain $\varepsilon_i$ itself. Thus, according to our involution convention for the Temperley-Lieb algebra, we have:
$$\varepsilon_i^*=\varepsilon_i$$

We conclude that the rescalings $e_i=N^{-1}\varepsilon_i$ satisfy $e_i^2=e_i=e_i^*$, as desired.
\end{proof}

As a second result now, making the link with Theorem 13.10, we have:

\index{Jones projection}

\begin{proposition}
The standard generators $e_i=N^{-1}\varepsilon_i$ of the Temperley-Lieb algebra $TL_N(k)$ have the following properties, where $tr$ is the trace obtained by closing:
\begin{enumerate}
\item $e_ie_j=e_je_i$ for $|i-j|\geq2$.

\item $e_ie_{i\pm1}e_i=[B:A]^{-1}e_i$.

\item $tr(we_{n+1})=[B:A]^{-1}tr(w)$, for any word $w\in<e_1,\ldots,e_n>$.
\end{enumerate}
\end{proposition}

\begin{proof}
This follows indeed by doing some elementary computations with diagrams, in the spirit of those performed in the proof of Proposition 13.12. Indeed:

\medskip

(1) This is clear from the definition of the diagrams $\varepsilon_i$.

\medskip

(2) This is clear as well from the definition of the diagrams $\varepsilon_i$.

\medskip

(3) This is something which is clear too, from the definition of $\varepsilon_{n+1}$.
\end{proof}

With the above results in hand, we can now reformulate our main finding about subfactors, namely Theorem 13.10, into something more conceptual, as follows:

\index{basic construction}
\index{Temperley-Lieb}

\begin{theorem}
Given a finite index subfactor $A_0\subset A_1$, with Jones tower
$$A_0\subset_{e_1}A_1\subset_{e_2}A_2\subset_{e_3}A_3\subset\ldots\ldots$$
the rescaled sequence of projections $e_1,e_2,e_3,\ldots\in B(H)$ produces a representation 
$$TL_N\subset B(H)$$
of the Temperley-Lieb algebra of index $N=[A_1:A_0]$.
\end{theorem}

\begin{proof}
The idea here is that Theorem 13.10, coming from the study of the basic construction, tells us that the rescaled sequence of projections $e_1,e_2,e_3,\ldots\in B(H)$ behaves algebrically exactly as the sequence of diagrams $\varepsilon_1,\varepsilon_2,\varepsilon_3,\ldots\in TL_N$ given by:
$$\varepsilon_1={\ }^\cup_\cap\quad,\quad 
\varepsilon_2=|\!{\ }^\cup_\cap\quad,\quad 
\varepsilon_3=||\!{\ }^\cup_\cap\quad,\quad 
\ldots$$

But these diagrams generate $TL_N$, and so we have an embedding $TL_N\subset B(H)$, where $H$ is the Hilbert space where our subfactor $A_0\subset A_1$ lives, as claimed.
\end{proof}

Before going further, with some examples, more theory, and consequences of Theorem 13.14, let us make the following key observation, also from Jones \cite{jo1}:

\index{relative commutant}
\index{higher commutant}

\begin{theorem} 
Given a finite index subfactor $A_0\subset A_1$, the graded algebra 
$$P=(P_k)$$
formed by the sequence of higher relative commutants
$$P_k=A_0'\cap A_k$$
contains the copy of the Temperley-Lieb algebra constructed above:
$$TL_N\subset P$$
This graded algebra $P=(P_k)$ is called ``planar algebra'' of the subfactor.
\end{theorem}

\begin{proof}
As a first observation, since the Jones projection $e_1:A_1\to A_0$ commutes with $A_0$, as was previously established in the above, we have:
$$e_1\in P_2'$$

By translation we obtain from this that we have, for any $k\in\mathbb N$:
$$e_1,\ldots,e_{k-1}\in P_k$$

Thus we have indeed an inclusion of graded algebras $TL_N\subset P$, as claimed.
\end{proof}

The point with the above result, which explains among others the terminology at the end, is that, in the context of Theorem 13.14, the planar algebra structure of $TL_N$, obtained by composing diagrams, extends into an abstract planar algebra structure of $P$. See \cite{jo3}. We will discuss all this, with full details, in the next chapter.

\section*{13c. Basic examples}

Let us discuss now some basic examples of subfactors, with concrete illustrations for all the above notions, constructions, and general theory. These examples will all come from group actions $G\curvearrowright P$, which are assumed to be minimal, in the sense that:
$$(P^G)'\cap P=\mathbb C$$ 

We will not provide proofs for the next few results to follow, the idea being that these constructions can be unified, and that we would like to keep the proofs for the unifications only. As a starting point, we have the following result, that we already know:

\index{Jones subfactor}

\begin{proposition}
Assuming that $G$ is a compact group, acting minimally on a ${\rm II}_1$ factor $P$, and that $H\subset G$ is a subgroup of finite index, we have a subfactor
$$P^G\subset P^H$$
having index $N=[G:H]$, called Jones subfactor.
\end{proposition}

\begin{proof}
This is something that we know, the idea being that the factoriality of $P^G,P^H$ comes from the minimality of the action, and that the index formula is clear.
\end{proof}

Along the same lines, we have the following result:

\index{Ocneanu subfactor}
\index{crossed product}

\begin{proposition}
Assuming that $G$ is a finite group, acting minimally on a ${\rm II}_1$ factor $P$, we have a subfactor as follows,
$$P\subset P\rtimes G$$
having index $N=|G|$, called Ocneanu subfactor.
\end{proposition}

\begin{proof}
This is standard as well, the idea being that the factoriality of $P\rtimes G$ comes from the minimality of the action, and that the index formula is clear.
\end{proof}

We have as well a third result of the same type, as follows:

\index{Wassermann subfactor}
\index{projective representation}

\begin{proposition}
Assuming that $G$ is a compact group, acting minimally on a ${\rm II}_1$ factor $P$, and that $G\to PU_n$ is a projective representation, we have a subfactor
$$P^G\subset (M_n(\mathbb C) \otimes P)^G$$
having index $N=n^2$, called Wassermann subfactor.
\end{proposition}

\begin{proof}
As before, the idea is that the factoriality of $P^G,(M_n(\mathbb C)\otimes P)^G$ comes from the minimality of the action, and the index formula is clear.
\end{proof}

The above subfactors look quite related, and indeed they are, due to:

\begin{theorem}
The Jones, Ocneanu and Wassermann subfactors are all of the same nature, and can be written as follows,
$$\left( P^G\subset P^H\right)\,\simeq\, \left( ({\mathbb C}\otimes P)^G\subset (l^\infty(G/H)\otimes P)^G\right)$$
$$\left( P\subset P\rtimes G\right)\,\simeq\,  \left( (l^\infty (G)\otimes P)^G\subset ({\mathcal L} (l^2(G))\otimes P)^G\right)$$
$$\left( P^G\subset (M_n(\mathbb C) \otimes P)^G\right)\,\simeq\, \left( ({\mathbb C}\otimes P)^G\subset (M_n(\mathbb C)\otimes P)^G\right)$$
with standard identifications for the various tensor products and fixed point algebras.
\end{theorem}

\begin{proof}
This is something very standard, modulo all kinds of standard identifications. We will explain all this more in detail later, after unifying these subfactors.
\end{proof}

In order to unify now the above constructions of subfactors, the idea is quite clear. Given a compact group $G$, acting minimally on a ${\rm II}_1$ factor $P$, and an inclusion of finite dimensional algebras $B_0\subset B_1$, endowed as well with an action of $G$, we would like to construct a kind of generalized Wassermann subfactor, as follows:
$$(B_0\otimes P)^G\subset (B_1\otimes P)^G$$

In order to do this, we must talk first about the finite dimensional algebras $B$, and about inclusions of such algebras $B_0\subset B_1$. Let us start with the following definition:

\begin{definition}
Associated to any finite dimensional algebra $B$ is its canonical trace, obtained by composing the left regular representation with the trace of $\mathcal L(B)$:
$$tr:B\subset\mathcal L(B)\to\mathbb C$$
We say that an inclusion of finite dimensional algebras $B_0\subset B_1$ is Markov if it commmutes with the canonical traces of $B_0,B_1$.
\end{definition}

In what regards the first notion, that of the canonical trace, this is something that we know well, from chapter 5. Indeed, as explained there, we can formally write $B=C(X)$, with $X$ being a finite quantum space, and the canonical trace $tr:B\to\mathbb C$ is then precisely the integration with respect to the ``counting measure'' on $X$.

\bigskip

In what regards the second notion, that of a Markov inclusion, this is something very natural too, and as a first example here, any inclusion of type $\mathbb C\subset B$ is Markov. In general, if we write $B_0=C(X_0)$ and $B_1=C(X_1)$, then the inclusion $B_0\subset B_1$ must come from a certain fibration $X_1\to X_0$, and the inclusion $B_0\subset B_1$ is Markov precisely when the fibration $X_1\to X_0$ commutes with the respective counting measures.

\bigskip

We will be back to Markov inclusions and their various properties on several occasions, in what follows. For our next purposes here, we just need the following result:

\begin{proposition}
Given a Markov inclusion of finite dimensional algebras $B_0\subset B_1$ we can perform to it the basic construction, as to obtain a Jones tower
$$B_0\subset_{e_1}B_1\subset_{e_2}B_2\subset_{e_3}B_3\subset\ldots\ldots$$
exactly as we did in the above for the inclusions of ${\rm II}_1$ factors.
\end{proposition}

\begin{proof}
This is something quite routine, by following the computations in the above, from the case of the ${\rm II}_1$ factors, and with everything extending well. It is of course possible to do something more general here, unifying the constructions for the inclusions of ${\rm II}_1$ factors $A_0\subset A_1$, and for the inclusions of Markov inclusions of finite dimensional algebras $B_0\subset B_1$, but we will not need this degree of generality, in what follows.
\end{proof}

With these ingredients in hand, getting back now to the Jones, Ocneanu and Wassermann subfactors, from Theorem 13.19, the point is that these constructions can be unified, and then further studied, the final result on the subject being as follows:

\index{fixed point subfactor}
\index{Markov inclusion}
\index{minimal action}
\index{centrally ergodic action}
\index{Wassermann subfactor}

\begin{theorem}
Let $G$ be a compact group, and $G\to Aut(P)$ be a minimal action on a ${\rm II}_1$ factor. Consider a Markov inclusion of finite dimensional algebras
$$B_0\subset B_1$$
and let $G\to Aut(B_1)$ be an action which leaves invariant $B_0$, and which is such that its restrictions to the centers of $B_0$ and $B_1$ are ergodic. We have then a subfactor
$$(B_0\otimes P)^G\subset (B_1\otimes P)^G$$
of index $N=[B_1:B_0]$, called generalized Wassermann subfactor, whose Jones tower is 
$$(B_1\otimes P)^G\subset(B_2\otimes P)^G\subset(B_3\otimes P)^G\subset\ldots$$
where $\{ B_i\}_{i\geq 1}$ are the algebras in the Jones tower for $B_0\subset B_1$, with the canonical actions of $G$ coming from the action $G\to Aut(B_1)$, and whose planar algebra is given by:
$$P_k=(B_0'\cap B_k)^G$$
These subfactors generalize the Jones, Ocneanu and Wassermann subfactors.
\end{theorem}

\begin{proof}
There are several things to be proved, the idea being as follows:

\medskip

(1) As before on various occasions, the idea is that the factoriality of the algebras $(B_i\otimes P)^G$ comes from the minimality of the action $G\to Aut(P)$, and that the index formula is clear as well, from the definition of the coupling constant and of the index. 

\medskip

(2) Regarding the Jones tower assertion, the precise thing to be checked here is that if $A\subset B\subset C$ is a basic construction, then so is the following sequence of inclusions:
$$(A\otimes P)^G\subset(B\otimes P)^G\subset(C\otimes P)^G$$

But this is something standard, which follows by verifying the basic construction conditions. We will be back to this in a moment, directly in a more general setting.

\medskip

(3) Next, regarding the planar algebra assertion, we have to prove here that for any indices $i\leq j$, we have the following equality between subalgebras of $B_j\otimes P$:
$$((B_i\otimes P)^G)'\cap(B_j\otimes P)^G=(B_i'\cap B_j^G)\otimes 1$$

But this is something which is routine too, following Wassermann \cite{was}, and we will be back to this in a moment, with full details, directly in a more general setting.

\medskip

(4) Finally, the last assertion, regarding the main examples of such subfactors, which are those of Jones, Ocneanu, Wassermann, follows from Theorem 13.19.
\end{proof}

In addition to the Jones, Ocneanu and Wassermann subfactors, discussed and unified in the above, we have the Popa subfactors, which are constructed as follows:

\index{Popa subfactor}
\index{diagonal subfactor}
\index{outer action}

\begin{proposition}
Given a discrete group $\Gamma=<g_1,\ldots,g_n>$, acting faithfully via outer automorphisms on a ${\rm II}_1$ factor $Q$, we have the following ``diagonal'' subfactor
$$\left\{ \begin{pmatrix}
g_1(q)\\
&\ddots\\
&& g_n(q)
\end{pmatrix} \Big| q\in Q\right\} \subset M_n(Q)$$
having index $N=n^2$, called Popa subfactor.
\end{proposition}

\begin{proof}
This is something standard, a bit as for the Jones, Ocneanu and Wassermann subfactors, with the result basically coming from the work of Popa, who was the main user of such subfactors. We will come in a moment with a more general result in this direction, involving discrete quantum groups, along with a complete proof.
\end{proof}

In order to unify now Theorem 13.22 and Proposition 13.23, observe that the diagonal subfactors can be written in the following way, by using a group dual:
$$(Q\rtimes\Gamma)^{\widehat{\Gamma}}\subset(M_n(\mathbb C)\otimes (Q\rtimes\Gamma))^{\widehat{\Gamma}}$$

Here the group dual $\widehat{\Gamma}$ acts on $P=Q\rtimes\Gamma$ via the dual of the action $\Gamma\subset Aut (Q)$, and on $M_n(\mathbb C)$ via the adjoint action of the following representation: 
$$\oplus g_i :\widehat{\Gamma}\to {\mathbb C}^n$$

Summarizing, we are led into quantum groups. Our plan in what follows will be that of discussing the quantum extension of Theorem 13.22, covering the diagonal subfactors as well, and this time with full details, and with examples and illustrations as well. 

\bigskip

We follow \cite{ba1}, where this extension of the Wassermann construction \cite{was} was developed. Let us start our discussion with some basic theory. We first have:

\index{coaction}

\begin{definition}
A coaction of a Woronowicz algebra $A$ on a finite von Neumann algebra $P$ is an injective morphism $\Phi:P\to P\otimes A''$ satisfying the following conditions:
\begin{enumerate}
\item Coassociativity: $(\Phi\otimes id)\Phi=(id\otimes\Delta)\Phi$.

\item Trace equivariance: $(tr\otimes id)\Phi=tr(.)1$.

\item Smoothness: $\overline{\mathcal P}^{\,w}=P$, where $\mathcal P=\Phi^{-1}(P\otimes_{alg}\mathcal A)$.
\end{enumerate}
\end{definition}

The above conditions come from what happens in the commutative case, $A=C(G)$, where they correspond to the usual associativity, trace equivariance and smoothness of the corresponding action $G\curvearrowright P$. Along the same lines, we have as well:

\index{minimal coaction}

\begin{definition}
A coaction $\Phi:P\to P\otimes A''$ as above is called:
\begin{enumerate}
\item Ergodic, if the algebra $P^\Phi=\left\{p\in P\big|\Phi(p)=p\otimes1\right\}$ reduces to $\mathbb C$.

\item Faithful, if the span of $\left\{(f\otimes id)\Phi(P)\big|f\in P_*\right\}$ is dense in $A''$.

\item Minimal, if it is faithful, and satisfies $(P^\Phi)'\cap P=\mathbb C$.
\end{enumerate} 
\end{definition}

Observe that the minimality of the action implies in particular that the fixed point algebra $P^\Phi$ is a factor. Thus, we are getting here to the case that we are interested in, actions producing factors, via their fixed point algebras. More on this later.

\bigskip

In order to prove our subfactor results, we need of some general theory regarding the minimal actions. Following Wassermann \cite{was}, let us start with the following definition:

\index{semidual coaction}

\begin{definition}
Let $\Phi:P\to P\otimes A''$ be a coaction. An eigenmatrix for a corepresentation $u\in B(H)\otimes A$ is an element $M\in B(H)\otimes P$ satisfying:
$$(id\otimes\Phi)M=M_{12}u_{13}$$
A coaction is called semidual if each corepresentation has a unitary eigenmatrix. 
\end{definition}

As a basic example here, the canonical coaction $\Delta:A\to A\otimes A$ is semidual. We will prove in what follows, following the work of Wassermann in the usual compact group case, that the minimal coactions of Woronowicz algebras are semidual. We first have:

\begin{proposition} 
If $\Phi:P\to P\otimes A''$ is a minimal coaction and $u\in Irr(A)$ is a corepresentation, then $u$ has a unitary eigenmatrix precisely when $P^u\neq\{ 0\}$.
\end{proposition}

\begin{proof}
Given $u\in M_n(A)$, consider the following unitary corepresentation:
$$u^+=(n\otimes 1)\oplus u=
\begin{pmatrix}1&0\\ 0&u\end{pmatrix}
\in M_2(M_n(\mathbb C)\otimes\mathcal A)
=M_2(\mathbb C)\otimes M_n(\mathbb C)\otimes\mathcal A$$

It is then routine to check, exactly as in \cite{was}, with the computation being explained in \cite{ba1}, that if the following algebra is a factor, then $u$ has a unitary eigenmatrix:
$$X_u=(M_2(\mathbb C)\otimes M_n(\mathbb C)\otimes P)^{\pi_{u^+}}$$

So, let us prove that $X_u$ is a factor. For this purpose, let $x\in Z(X_u)$. We have then $1\otimes1\otimes P^\Phi\subset X_u$, and from the irreducibility of the inclusion $P^\pi\subset P$ we obtain that:
$$x\in M_2(\mathbb C)\otimes M_n(\mathbb C)\otimes 1$$

On the other hand, we have the following formula:
\begin{eqnarray*}
X_u\cap M_2(\mathbb C)\otimes M_n(\mathbb C)\otimes 1
&=&(M_2(\mathbb C)\otimes M_n(\mathbb C))^{i_{u^+}}\otimes 1\\
&=&End(u^+)\otimes 1
\end{eqnarray*}

Since our corepresentation $u$ was chosen to be irreducible, it follows that $x$ must be of the following form, with $y\in M_n(\mathbb C)$, and with $\lambda\in\mathbb C$:
$$x=\begin{pmatrix}y&0\\0&\lambda I\end{pmatrix}\otimes 1$$

Now let us pick a nonzero element $p\in P^u$, and write:
$$\Phi(p)=\sum_{ij}p_{ij}\otimes u_{ij}$$

Then $\Phi(p_{ij})=\sum_kp_{kj}\otimes u_{ki}$ for any $i,j$, and so each column of $(p_{ij})_{ij}$ is a $u$-eigenvector. Choose such a nonzero column $l$ and let $m^i$ be the matrix having the $i$-th row equal to $l$, and being zero elsewhere. Then $m_i$ is a $u$-eigenmatrix for any $i$, and this implies that:
$$\begin{pmatrix}0&m^i\\0&0\end{pmatrix}\in X_u$$

The commutation relation of this matrix with $x$ is as follows:
$$\begin{pmatrix}y&0\\0&\lambda I\end{pmatrix}
\begin{pmatrix}0&m^i\\0&0\end{pmatrix}= 
\begin{pmatrix}0&m^i\\0&0\end{pmatrix}
\begin{pmatrix}y&0\\0&\lambda I\end{pmatrix}$$

But this gives $(y-\lambda I)m^i=0$. Now by definition of $m^i$, this shows that the $i$-th column of $y-\lambda I$ is zero. Thus $y-\lambda I=0$, and so $x=\lambda 1$, as desired. 
\end{proof}

We can now prove a main result about minimal coactions, as follows:

\index{minimal coaction}
\index{semidual coaction}

\begin{theorem}
The minimal coactions are semidual.
\end{theorem}

\begin{proof}
Let $K$ be the set of finite dimensional unitary corepresentations of $A$ which have unitary eigenmatrices. Then, according to the above, the following happen:

\medskip

(1) $K$ is stable under taking tensor products. Indeed, if $M,N$ are unitary eigenmatrices for $u,w$, then $M_{13}N_{23}$ is a unitary eigenmatrix for $u\otimes w$.

\medskip

(2) $K$ is stable under taking sums. Indeed, if $M_i$ are unitary eigenmatrices for $u_i$, then $diag(M_i)$ is a unitary eigenmatrix for $\oplus u_i$.

\medskip

(3) $K$ is stable under substractions. Indeed, if $M$ is an eigenmatrix for $U=\oplus_{i=1}^nu_i$, then the first $\dim(u_1)$ columns of $M$ are formed by elements of $P^{u_1}$, the next $\dim(u_2)$ columns of $M$ are formed by elements of $P^{u_2}$, and so on. Now if $M$ is unitary, it is in particular invertible, so all $P^{u_i}$ are different from $\{0\}$, and we may conclude that we can indeed substract corepresentations from $U$, by using Proposition 13.27.

\medskip

(4) $K$ is stable under complex conjugation. Indeed, by the above results we may restrict attention to irreducible corepresentations. Now if $u\in Irr(A)$ has a nonzero eigenmatrix $M$ then $\overline{M}$ is an eigenmatrix for $\overline{u}$. By Proposition 13.27 we obtain from this that $P^{\overline{u}}\neq\{0\}$, and we may conclude by using again Proposition 13.27.

\medskip

With this in hand, by using Peter-Weyl, we obtain the result. See \cite{ba1}.
\end{proof}

Let us construct now the fixed point subfactors. We first have:

\begin{proposition}
Consider a Woronowicz algebra $A=(A,\Delta,S)$, and denote by $A_\sigma$ the Woronowicz algebra $(A,\sigma\Delta ,S)$, where $\sigma$ is the flip. Given coactions
$$\beta:B\to B\otimes A$$
$$\pi:P\to P\otimes A_\sigma$$
with $B$ being finite dimensional, the following linear map, while not being multiplicative in general, is coassociative with respect to the comultiplication $\sigma\Delta$ of $A_\sigma$,
$$\beta\odot\pi:B\otimes P\to B\otimes P\otimes A_\sigma$$
$$b\otimes p\to \pi (p)_{23}((id\otimes S)\beta(b))_{13}$$
and its fixed point space, which is by definition the following linear space,
$$(B\otimes P)^{\beta\odot\pi}=\left\{x\in B\otimes P\Big|(\beta\odot\pi )x=x\otimes 1\right\}$$
is then a von Neumann subalgebra of $B\otimes P$. 
\end{proposition}

\begin{proof}
This is something standard, which follows from a straightforward algebraic verification, explained in \cite{ba1}. As mentioned in the statement, to be noted is that the tensor product coaction $\beta\odot\pi$ is not multiplicative in general. See \cite{ba1}.
\end{proof}

Our first task is to investigate the factoriality of such algebras, and we have here:

\begin{theorem}
If $\beta:B\to B\otimes A$ is a coaction and $\pi:P\to P\otimes A_\sigma$ is a minimal coaction, then the following conditions are equivalent:
\begin{enumerate}
\item The von Neumann algebra $(B\otimes P)^{\beta\odot\pi}$ is a factor.

\item The coaction $\beta$ is centrally ergodic, $Z(B)\cap B^\beta=\mathbb C$.
\end{enumerate}
\end{theorem}

\begin{proof}
This is something standard, from \cite{ba1}, the idea being as follows:

\medskip

(1) Our first claim, which is something whose proof is a routine verification, explained in \cite{ba1}, based on the semiduality of the minimal coaction $\pi$, that we know from Theorem 13.28, is that the following diagram is a non-degenerate commuting square:
$$\begin{matrix}
P&\subset&B\otimes P\\ 
\cup &\ &\cup \\
P^\pi&\subset&(B\otimes P)^{\beta\odot\pi}
\end{matrix}$$

(2) In order to prove now the result, it is enough to check the following equality, between von Neumann subalgebras of the algebra $B\otimes P$:
$$Z((B\otimes P)^{\beta\odot\pi})=(Z(B)\cap B^\beta)\otimes 1$$

So, let $x$ be in the algebra on the left. Then $x$ commutes with $1\otimes P^\pi$, so it has to be of the form $b\otimes 1$. Thus $x$ commutes with $1\otimes P$. But $x$ commutes with $(B\otimes P)^{\beta\odot\pi}$, and from the non-degeneracy of the above square, $x$ commutes with $B\otimes P$, and in particular with $B\otimes 1$. Thus we have $b\in Z(B)\cap B^\beta$. As for the other inclusion, this is obvious. 
\end{proof}

In view of the above result, we can talk about subfactors of type $(B_0\otimes P)^G\subset(B_1\otimes P)^G$. In order to investigate such subfactors, we will need the following technical result:

\begin{proposition}
Consider two commuting squares, as follows:
$$\begin{matrix}
F&\subset&E&\subset&D\\ 
\cup&&\cup&&\cup\\
A&\subset&B&\subset&C\\
\end{matrix}$$
\begin{enumerate}
\item If the square on the left and the big square are non-degenerate, then so is the square on the right.

\item If both squares are non-degenerate, $F\subset E\subset D$ is a basic construction, and the Jones projection $e\in D$ for this basic construction belongs to $C$, then the square on the right is the basic construction for the square on the left.
\end{enumerate}
\end{proposition}

\begin{proof}
We have several things to be proved, the idea being as follows:

\medskip

(1) This assertion is clear from the following computation:
$$D
={\overline{sp}^{\,w}\,}CF
={\overline{sp}^{\,w}\,}CBF
={\overline{sp}^{\,w}\,} CE$$

(2) Let $\Psi :D\to C$ be the expectation. By non-degeneracy, we have that: 
$$E={\overline{sp}^w\,} FB={\overline{sp}^w\,} BF$$

We also have $D={\overline{sp}^w\,} EeE$ by the basic construction, so we get that:
\begin{eqnarray*}
C
&=&\Psi(D)\\
&=&\Psi({\overline{sp}^{\,w}\,}EeE)\\
&=&\Psi({\overline{sp}^{\,w}\,}BFeFB)\\
&=&\Psi({\overline{sp}^{\,w}\,}BeFB)\\
&=&{\overline{sp}^{\,w}\,}Be\Psi(F)B\\
&=&{\overline{sp}^{\,w}\,}BeAB\\
&=&{\overline{sp}^{\,w}\,}BeB
\end{eqnarray*}

Thus the algebra $C$ is generated by $B$ and $e$, and this gives the result. 
\end{proof}

Next in line, we have the following key technical result:

\begin{proposition}
If $\beta:B\to B\otimes A$ is a coaction then
$$\begin{matrix}
A&\subset&B\otimes A\\ 
\cup&&\uparrow\beta\\ 
\mathbb C&\subset&B\\
\end{matrix}$$
is a non-degenerate commuting square.
\end{proposition}

\begin{proof}
From the $\beta$-equivariance of the trace we get that the inclusion on the left commutes with the traces, so that the above is a commuting diagram of finite von Neumann algebras. From the formula of the expectation $E_{\beta}=(id\otimes\int_A)\beta$ we get that this diagram is a commuting square. Choose now an orthonormal basis $\{b_i\}$ of $B$, write $\beta:b_i\to\sum_jb_j\otimes u_{ji}$, and consider the corresponding unitary corepresentation:
$$u_\beta=\sum e_{ij}\otimes u_{ij}$$

Then for any $k$ and any $a\in A$ we have the following computation:
\begin{eqnarray*}
\sum_i\beta(b_i)(1\otimes u_{ki}^*a)
&=&\sum_{ij}b_j\otimes u_{ji}u_{ki}^*a\\
&=&\sum_{ij}b_j\otimes \delta_{jk}a\\
&=&b_k\otimes a
\end{eqnarray*}

Thus our commuting square is non-degenerate, as claimed. 
\end{proof}

Getting now to the generalized Wassermann subfactors, we first have:

\begin{proposition}
Given a Markov inclusion of finite dimensional algebras $B_0\subset B_1$, construct its Jones tower, and denote it as follows:
$$B_0\subset B_1\subset_{e_1}B_2=<B_1,e_1> \subset_{e_2}B_3=<B_2,e_2> \subset_{e_3}\ldots$$
If $\beta_1:B_1\to B_1\otimes A$ is a coaction/anticoaction leaving $B_0$ invariant then there exists a unique sequence $\{\beta_i\}_{i\geq 0}$ of coactions/anticoactions 
$$\beta_i:B_i\to B_i\otimes A$$
such that each $\beta_i$ extends $\beta_{i-1}$ and leaves invariant the Jones projection $e_{i-1}$.
\end{proposition}

\begin{proof}
By taking opposite inclusions we see that the assertion for anticoactions is equivalent to the one for coactions, that we will prove now. The uniqueness is clear from $B_i=<B_{i-1},e_{i-1}>$. For the existence, we can apply Proposition 13.32 to:
$$\begin{matrix}
A&\subset&B_0\otimes A&\subset&B_1\otimes A\\ 
\cup&&\uparrow\beta_0&&\uparrow\beta_1\\ 
\mathbb C&\subset&B_0&\subset&B_1
\end{matrix}$$

Indeed, we get in this way that the square on the right is a non-degenerate. Now by performing basic constructions to it, we get a sequence as follows:
$$\begin{matrix}
B_0\otimes A&\subset&B_1\otimes A&\subset&B_2\otimes A&\subset&B_3\otimes A&\subset&\ldots\\ 
\uparrow\beta_0&&\uparrow\beta_1&&\uparrow\beta_2&&\uparrow\beta_3\\ 
B_0&\subset&B_1&\subset&B_2&\subset&B_3&\subset&\ldots
\end{matrix}$$

It is easy to see from definitions that the $\beta_i$ are coactions, that they extend each other, and that they leave invariant the Jones projections. But this gives the result.
\end{proof}

With the above technical results in hand, we can now formulate our main theorem regarding the fixed point subfactors, of the most possible general type, as follows:

\index{fixed point subfactor}
\index{Wassermann subfactor}

\begin{theorem}
Let $G$ be a compact quantum group, and $G\to Aut(P)$ be a minimal action on a ${\rm II}_1$ factor. Consider a Markov inclusion of finite dimensional algebras
$$B_0\subset B_1$$
and let $G\to Aut(B_1)$ be an action which leaves invariant $B_0$ and which is such that its restrictions to the centers of $B_0$ and $B_1$ are ergodic. We have then a subfactor
$$(B_0\otimes P)^G\subset (B_1\otimes P)^G$$
of index $N=[B_1:B_0]$, called generalized Wassermann subfactor, whose Jones tower is 
$$(B_1\otimes P)^G\subset(B_2\otimes P)^G\subset(B_3\otimes P)^G\subset\ldots$$
where $\{ B_i\}_{i\geq 1}$ are the algebras in the Jones tower for $B_0\subset B_1$, with the canonical actions of $G$ coming from the action $G\to Aut(B_1)$, and whose planar algebra is given by:
$$P_k=(B_0'\cap B_k)^G$$
These subfactors generalize the Jones, Ocneanu, Wassermann and Popa subfactors.
\end{theorem}

\begin{proof}
We have several things to be proved, the idea being as follows:

\medskip

(1) The first part of the statement, regarding the factoriality, the index and the Jones tower assertions, is something that follows exactly as in the classical group case. 

\medskip

(2) In order to prove now the planar algebra assertion, consider the following diagram, with $i<j$ being arbitrary integers:
$$\begin{matrix}
P&\subset&B_i\otimes P&\subset&B_j\otimes P\\ 
\cup&&\cup&&\cup\\ 
P^\pi&\subset&(B_i\otimes P)^{\beta_i\otimes\pi}&\subset&(B_j\otimes P)^{\beta_j\otimes\pi}
\end{matrix}$$

We know from Proposition 13.32 that the big square and the square on the left are both non-degenerate commuting squares. Thus Proposition 13.31 applies, and shows that the square on the right is a non-degenerate commuting square. 

\medskip

(3) Consider now the following sequence of non-degenerate commuting squares:
$$\begin{matrix}
B_0\otimes P&\subset&B_1\otimes P&\subset&B_2\otimes P&\subset&\ldots\\
\cup&&\cup&&\cup&&\\ 
(B_0\otimes P)^{\beta_0\otimes\pi}&\subset&(B_1\otimes P)^{\beta_1\otimes\pi}&\subset&(B_2\otimes P)^{\beta_2\otimes\pi}&\subset&\ldots
\end{matrix}$$

Since the Jones projections live in the lower line, Proposition 13.32 applies and shows that this is a sequence of basic constructions for non-degenerate commuting squares. In particular the lower line is a sequence of basic constructions, as desired. 

\medskip

(4) Finally, we already know from Theorem 13.22 that our construction generalizes the Jones, Ocneanu and Wassermann subfactors. As for the Popa subfactors, the result here follows from the discussion made after Proposition 13.23.
\end{proof}

\section*{13d. The index theorem}

Let us go back now to the arbitrary subfactors, with Theorem 13.14 being our main result. As an interesting consequence of the above results, somehow contradicting the ``continuous geometry'' philosophy that has being going on so far, in relation with the ${\rm II}_1$ factors, we have the following surprising result, also from Jones' original paper \cite{jo1}:

\index{index theorem}

\begin{theorem}
The index of subfactors $A\subset B$ is ``quantized'' in the $[1,4]$ range,
$$N\in\left\{4\cos^2\left(\frac{\pi}{n}\right)\Big|n\geq3\right\}\cup[4,\infty]$$
with the obstruction coming from the existence of the representation $TL_N\subset B(H)$.
\end{theorem}

\begin{proof}
This comes from the basic construction, and more specifically from the combinatorics of the Jones projections $e_1,e_2,e_3,\ldots$, the idea being as folows:

\medskip

(1) In order to best comment on what happens, when iterating the basic construction, let us record the first few values of the numbers in the statement:
$$4\cos^2\left(\frac{\pi}{3}\right)=1\quad,\quad 
4\cos^2\left(\frac{\pi}{4}\right)=2$$
$$4\cos^2\left(\frac{\pi}{5}\right)=\frac{3+\sqrt{5}}{2}\quad,\quad 
4\cos^2\left(\frac{\pi}{6}\right)=3$$
$$\ldots$$

(2) When performing a basic construction, we obtain, by trace manipulations on $e_1$:
$$N\notin(1,2)$$

With a double basic construction, we obtain, by trace manipulations on $<e_1,e_2>$:
$$N\notin\left(2,\frac{3+\sqrt{5}}{2}\right)$$

With a triple basic construction, we obtain, by trace manipulations on $<e_1,e_2,e_3>$:
$$N\notin\left(\frac{3+\sqrt{5}}{2},3\right)$$

Thus, we are led to the conclusion in the statement, by a kind of recurrence, involving a certain family of orthogonal polynomials.

\medskip

(3) In practice now, the most elegant way of proving the result is by using the fundamental fact, explained in Theorem 13.14, that that sequence of Jones projections $e_1,e_2,e_3,\ldots\subset B(H)$ generate a copy of the Temperley-Lieb algebra of index $N$:
$$TL_N\subset B(H)$$

With this result in hand, we must prove that such a representation cannot exist in index $N<4$, unless we are in the following special situation:
$$N=4\cos^2\left(\frac{\pi}{n}\right)$$

But this can be proved by using some suitable trace and positivity manipulations on $TL_N$, as in (2) above. For full details here, we refer to \cite{ghj}, \cite{jo1}, \cite{jsu}.
\end{proof}

The above result raises the question of understanding if there are further restrictions on the index of subfactors $A\subset B$, in the range found there, namely:
$$N\in\left\{4\cos^2\left(\frac{\pi}{n}\right)\Big|n\geq3\right\}\cup[4,\infty]$$

This question is quite tricky, because it depends on the ambient factor $B\subset B(H)$, and also on the irreducibility assumption on the subfactor, namely $A'\cap B=\mathbb C$, which is something quite natural, and can be added to the problem.

\bigskip

All this is quite technical, to be discussed later on, when doing more advanced subfactor theory. In the simplest formulation of the question, the answer is generally ``no'', as shown by the following result, also from Jones' original paper \cite{jo1}:

\index{R}
\index{hyperfinite subfactor}
\index{irreducible subfactor}

\begin{theorem}
Consider the Murray-von Neumann hyperfinite ${\rm II}_1$ factor $R$. Its subfactors $R_0\subset R$ are then as follows:
\begin{enumerate}
\item They exist for all admissible index values, $N\in\left\{4\cos^2\left(\frac{\pi}{n}\right)|n\geq3\right\}\cup[4,\infty]$.

\item In index $N\leq4$, they can be realized as irreducible subfactors, $R_0'\cap R=\mathbb C$. 

\item In index $N>4$, they can be realized as arbitrary subfactors.
\end{enumerate}
\end{theorem}

\begin{proof}
This is something quite tricky, worked out in Jones' original paper \cite{jo1}, and requiring some advanced algebra methods, the idea being as follows:

\medskip

(1) This basically follows by taking a copy of the Temperley-Lieb algebra $TL_N$, and then building a subfactor out of it, first by constructing a certain inclusion of inductive limits of finite dimensional algebras, $\mathcal A\subset\mathcal B$, and then by taking the weak closure, which produces copies of the Murray-von Neumann hyperfinite ${\rm II}_1$ factor, $A\simeq B\simeq R$.

\medskip

(2) This follows by examining and fine-tuning the construction in (1), which can be performed as to have control over the relative commutant.

\medskip

(3) This follows as well from (1), and with the simplest proof here being in fact quite simple, based on a projection trick. 
\end{proof}

As another application now, which is more theoretical, let us go back to the question of defining the index of a subfactor in a purely algebraic manner, which was open since chapter 10. The answer here, due to Pimsner and Popa \cite{ppo}, is as follows:

\index{Pimsner-Popa basis}

\begin{theorem}
Any finite index subfactor $A\subset B$ has an algebraic orthonormal basis, called Pimsner-Popa basis, which is constructed as follows:
\begin{enumerate}
\item In integer index, $N\in\mathbb N$, this is a usual basis, of type $\{b_1,\ldots,b_N\}$, whose length is exactly the index.

\item In non-integer index, $N\notin\mathbb N$, this is something of type $\{b_1,\ldots,b_n,c\}$, having length $n+1$, with $n=[N]$, and with $N-n\in(0,1)$ being related to $c$. 
\end{enumerate}
\end{theorem}

\begin{proof}
This is something quite technical, which follows from the basic theory of the basic construction. We refer here to the paper of Pimsner and Popa \cite{ppo}.
\end{proof}

\section*{13e. Exercises} 

There has been a lot of exciting theory in this chapter, leading us from functional analysis to concrete combinatorics, and as an exercise on all this, we have:

\begin{exercise}
Clarify all the details for the Jones index theorem, stating that
$$N\in\left\{4\cos^2\left(\frac{\pi}{n}\right)\Big|n\geq3\right\}\cup[4,\infty]$$
with the obstruction coming from the existence of the representation $TL_N\subset B(H)$.
\end{exercise}

This is something that we already discussed in the above, but with a few details missing. Time to have this understood, along the above lines.

\chapter{Planar algebras}

\section*{14a. Planar algebras}

We have seen the foundations of subfactor theory, and the main examples of subfactors, all having integer index. Following Jones' paper \cite{jo3}, in this chapter we go into the core of the theory, with the idea in mind of axiomatizing the combinatorics of a given subfactor $A_0\subset A_1$, via an object similar to the tensor categories for the quantum groups.

\bigskip

So, our starting point will be an arbitrary subfactor $A_0\subset A_1$, assumed to have finite index, $[A_1:A_0]<\infty$. Let us first review first what can be said about it, by using the Jones basic construction. We recall from chapter 13 that we have the following result:

\begin{theorem}
Given an inclusion of ${\rm II}_1$ factors $A_0\subset A_1$, with Jones tower
$$A_0\subset_{e_1}A_1\subset_{e_2}A_2\subset_{e_3}A_3\subset\ldots\ldots$$
the sequence of projections $e_1,e_2,e_3,\ldots\in B(H)$ produces a representation 
$$TL_N\subset B(H)$$
of the Temperley-Lieb algebra of index $N=[A_1:A_0]$. Moreover, we have
$$TL_N\subset P$$
where $P=(P_k)$ is the graded algebra formed by the commutants $P_k=A_0'\cap A_k$.
\end{theorem}

\begin{proof}
There are two statements here, both due to Jones \cite{jo1}, that we know from chapter 13 above, the idea for this, in short, being as follows:

\medskip

(1) A detailed study of the basic construction, performed in chapter 13, shows that the rescaled sequence of Jones projections $e_1,e_2,e_3,\ldots\in B(H)$ behaves algebrically exactly as the sequence of standard generators $\varepsilon_1,\varepsilon_2,\varepsilon_3,\ldots\in TL_N$. Thus we have an embedding $TL_N\subset B(H)$, where $H$ is the Hilbert space where our subfactor $A_0\subset A_1$ lives.

\medskip

(2) Once again by carefully looking at the Jones basic construction, the more precise conclusion is that the image of the representation $TL_N\subset B(H)$ constructed above lives indeed in the graded algebra $P=(P_k)$ formed by the commutants $P_k=A_0'\cap A_k$.
\end{proof}

Quite remarkably, the planar algebra structure of $TL_N$, taken in an intuitive sense, that of composing planar diagrams, extends to a planar algebra structure of $P$. In order to discuss this key result, let us start with the axioms for planar algebras. Following Jones' paper \cite{jo3}, we have the following definition:

\index{planar tangle}
\index{planar algebra}

\begin{definition}
The planar algebras are defined as follows:
\begin{enumerate}
\item We consider rectangles in the plane, with the sides parallel to the coordinate axes, and taken up to planar isotopy, and we call such rectangles boxes.

\item A labeled box is a box with $2k$ marked points on its boundary, $k$ on its upper side, and $k$ on its lower side, for some integer $k\in\mathbb N$.

\item A tangle is labeled box, containing a number of labelled boxes, with all marked points, on the big and small boxes, being connected by noncrossing strings.

\item A planar algebra is a sequence of finite dimensional vector spaces $P=(P_k)$, together with linear maps $P_{k_1}\otimes\ldots\otimes P_{k_r}\to P_k$, one for each tangle, such that the gluing of tangles corresponds to the composition of linear maps.
\end{enumerate}
\end{definition}

In this definition we are using rectangles, but everything being up to isotopy, we could have used instead circles with marked points, as in \cite{jo3}, which is the same thing. Our choice for using rectangles comes from the main examples that we have in mind, to be discussed below, where the planar algebra structure is best viewed by using rectangles.

\bigskip

Let us also mention that Definition 14.2 is something quite simplified, based on \cite{jo3}. As explained in \cite{jo3}, in order for subfactors to produce planar algebras and vice versa, there are quite a number of supplementary axioms that must be added, and in view of this, it is perhaps better to start with something stronger than Definition 14.2, as basic axioms. However, as before with rectangles vs circles, our axiomatic choices here are mainly motivated by the concrete examples that we have in mind.

\bigskip

As a basic example of a planar algebra, we have the Temperley-Lieb algebra:

\begin{theorem}
The Temperley-Lieb algebra $TL_N$, viewed as sequence of finite dimensional vector spaces
$$TL_N=(TL_N(k))_{k\in\mathbb N}$$
is a planar algebra in the above sense, with the corresponding linear maps associated to the planar tangles
$$TL_N(k_1)\otimes\ldots\otimes TL_N(k_r)\to TL_N(k)$$
appearing by putting the various $TL_N(k_i)$ diagrams into the small boxes of the given tangle, which produces a $TL_N(k)$ diagram.
\end{theorem}

\begin{proof}
This is something trivial, which follows from definitions:

\medskip

(1) Assume indeed that we are given a planar tangle $\pi$ in the sense of Definition 14.2, consisting of a box having $2k$ marked points on its boundary, and containing $r$ small boxes, having respectively $2k_1,\ldots,2k_r$ marked points on their boundaries, and then a total of $k+\Sigma k_i$ noncrossing strings, connecting the various $2k+\Sigma 2k_i$ marked points.

\medskip

(2) We want to associate to this planar tangle $\pi$ a linear map as follows:
$$T_\pi:TL_N(k_1)\otimes\ldots\otimes TL_N(k_r)\to TL_N(k)$$

For this purpose, by linearity, it is enough to construct elements as follows, for any choice of Temperley-Lieb diagrams $\sigma_i\in TL_N(k_i)$, with $i=1,\ldots,r$:
$$T_\pi(\sigma_1\otimes\ldots\otimes\sigma_r)\in TL_N(k)$$

(3) But constructing such an element is obvious, simply by putting the various diagrams $\sigma_i\in TL_N(k_i)$ into the small boxes the given tangle $\pi$. Indeed, this procedure produces a certain diagram in $TL_N(k)$, that we can call $T_\pi(\sigma_1\otimes\ldots\otimes\sigma_r)$, as above.

\medskip

(4) Finally, we still have to check that everything is well-defined up to planar isotopy, and that the gluing of tangles corresponds to the composition of linear maps. But both these checks are trivial, coming from the definition of $TL_N$, and we are done.
\end{proof}

As a conclusion, $P=TL_N$ is indeed a planar algebra, and of somewhat ``trivial'' type, with the triviality coming from the fact that, in this case, the elements of $P$ are planar diagrams themselves, and so the planar structure appears trivially. The Temperley-Lieb planar algebra $TL_N$ is however an important planar algebra, because it is the ``smallest'' one, appearing inside the planar algebra of any subfactor. But more on this later, when talking about the relation between planar algebras and subfactors.

\bigskip

Moving ahead, here is our second basic example of a planar algebra, due to Bisch-Jones \cite{bjo}, which is still ``trivial'' in the above sense, with the elements of the planar algebra being planar diagrams themselves, but which appears in a more complicated way:

\index{Fuss-Catalan algebra}

\begin{theorem}
The Fuss-Catalan algebra $FC_{N,M}$, which appears by coloring the Temperley-Lieb diagrams with black/white colors, clockwise, as follows 
$$\circ\bullet\bullet\circ\circ\bullet\bullet\circ\ldots\ldots\ldots\circ\bullet\bullet\circ$$
and keeping those diagrams whose strings connect either $\circ-\circ$ or $\bullet-\bullet$, is a planar algebra, with again the corresponding linear maps associated to the planar tangles
$$FC_{N,M}(k_1)\otimes\ldots\otimes FC_{N,M}(k_r)\to FC_{N,M}(k)$$
appearing by putting the various $FC_{N,M}(k_i)$ diagrams into the small boxes of the given tangle, which produces a $FC_{N,M}(k)$ diagram.
\end{theorem}

\begin{proof}
The proof here is nearly identical to the proof of Theorem 14.3, with the only change appearing at the level of the colors. To be more precise:

\medskip

(1) Forgetting about upper and lower sequences of points, which must be joined by strings, a Temperley-Lieb diagram can be thought of as being just a collection of strings, say black strings, which compose in the obvious way, with the rule that the value of the circle, which is now a black circle, is $N$. And it is this obvious composition rule that gives the planar algebra structure, as explained in the proof of Theorem 14.3. 

\medskip

(2) Similarly, forgetting about sequences of points, a Fuss-Catalan diagram can be thought of as being a collection of strings, which come now in two colors, black and white. These Fuss-Catalan diagrams compose in the obvious way, with the rule that the value of the black circle is $N$, and the value of the white circle is $M$. And it is this obvious composition rule that gives the planar algebra structure, as before for $TL_N$.
\end{proof}

The same comments as those for $TL_N$ apply. On one hand, $FC_{N,M}$ is by definition a ``trivial'' planar algebra, with the triviality coming from the fact that its elements are planar diagrams themselves. On the other hand, $FC_{N,M}$ is an important planar algebra, because it can be shown to appear inside the planar algebra of any subfactor $A\subset B$, assuming that an intermediate subfactor $A\subset C\subset B$ is present. But more on this later, when talking about the relation between planar algebras and subfactors.

\bigskip

Getting back now to generalities, and to Definition 14.2, that of a general planar algebra, we have so far two illustrations for it, which, while both important, are both ``trivial'', with the planar structure simply coming from the fact that, in both these illustrations, the elements of the planar algebra are planar diagrams themselves.

\bigskip

In general, the planar algebras are more complicated than this, and we will see some further examples in a moment. However, the idea is very simple, namely ``the elements of a planar algebra are not necessarily diagrams, but they behave like diagrams".

\bigskip

Let us begin with the construction of the tensor planar algebra $\mathcal T_N$, which is the third most important planar algebra, in our series of examples. This algebra is as follows:

\index{tensor planar algebra}
\index{Kronecker symbol}

\begin{definition}
The tensor planar algebra $\mathcal T_N$ is the sequence of vector spaces 
$$P_k=M_N(\mathbb C)^{\otimes k}$$
with the multilinear maps associated to the various $k$-tangles
$$T_\pi:P_{k_1}\otimes\ldots\otimes P_{k_r}\to P_k$$
being given by the following formula, in multi-index notation, 
$$T_\pi(e_{i_1}\otimes\ldots\otimes e_{i_r})=\sum_j\delta_\pi(i_1,\ldots,i_r:j)e_j$$
with the Kronecker symbols $\delta_\pi$ being $1$ if the indices fit, and being $0$ otherwise.
\end{definition}

In other words, we are using here a construction which is very similar to the construction $\pi\to T_\pi$ from easy quantum groups. We put the indices of the basic tensors on the marked points of the small boxes, in the obvious way, and the coefficients of the output tensor are then given by Kronecker symbols, exactly as in the easy case.

\bigskip

The fact that we have indeed a planar algebra, in the sense that the gluing of tangles corresponds to the composition of linear maps, as required by Definition 14.2, is something elementary, in the same spirit as the verification of the functoriality properties of the correspondence $\pi\to T_\pi$, discussed in chapter 8, and we refer here to Jones \cite{jo3}. 

\bigskip

Let us discuss now a second planar algebra of the same type, which is important as well for various reasons, namely the spin planar algebra $\mathcal S_N$. This planar algebra appears somewhat as a ``square root'' of the tensor planar algebra $\mathcal T_N$, and its construction is quite similar, but by using this time the algebra $\mathbb C^N$ instead of the algebra $M_N(\mathbb C)$.

\bigskip

There is one subtlety, however, coming from the fact that the general planar algebra formalism, from Definition 14.2 above, requires the tensors to have even length. Note that this was automatic for $\mathcal T_N$, where the tensors of $M_N(\mathbb C)$ have length 2.  In the case of the spin planar algebra $\mathcal S_N$, we want the vector spaces to be:
$$P_k=(\mathbb C^N)^{\otimes k}$$

Thus, we must double the indices of the tensors, in the following way:

\begin{definition}
We write the standard basis of $(\mathbb C^N)^{\otimes k}$ in $2\times k$ matrix form,
$$e_{i_1\ldots i_k}=
\begin{pmatrix}i_1 & i_1 &i_2&i_2&i_3&\ldots&\ldots\\
i_k&i_k&i_{k-1}&\ldots&\ldots&\ldots&\ldots 
\end{pmatrix}$$
by duplicating the indices, and then writing them clockwise, starting from top left.
\end{definition}

Now with this convention in hand for the tensors, we can formulate the construction of the spin planar algebra $\mathcal S_N$, also from Jones \cite{jo3}, as follows:

\index{spin planar algebra}
\index{Kronecker symbol}

\begin{definition}
The spin planar algebra $\mathcal S_N$ is the sequence of vector spaces 
$$P_k=(\mathbb C^N)^{\otimes k}$$
written as above, with the multilinear maps associated to the various $k$-tangles
$$T_\pi:P_{k_1}\otimes\ldots\otimes P_{k_r}\to P_k$$
being given by the following formula, in multi-index notation, 
$$T_\pi(e_{i_1}\otimes\ldots\otimes e_{i_r})=\sum_j\delta_\pi(i_1,\ldots,i_r:j)e_j$$
with the Kronecker symbols $\delta_\pi$ being $1$ if the indices fit, and being $0$ otherwise.
\end{definition}

In other words, we are using exactly the same construction as for the tensor planar algebra $\mathcal T_N$, which was itself very related to the easy quantum group formalism, but with $M_N(\mathbb C)$ replaced by $\mathbb C^N$, with the indices doubled, as in Definition 14.6. As before with the tensor planar algebra $\mathcal T_N$, the fact that the spin planar algebra $\mathcal S_N$ is indeed a planar algebra is something rather trivial, coming from definitions. 

\bigskip

Observe however that, unlike our previous planar algebras $TL_N$ and $FC_{N,M}$, which were ``trivial'' planar algebras, their elements being planar diagrams themselves, the planar algebras $\mathcal T_N$ and $\mathcal S_N$ are not trivial, their elements being not exactly planar diagrams. Let us also mention that the planar algebras $\mathcal T_N$ and $\mathcal S_N$ are important for a number of reasons, in the context of the fixed point subfactors, to be discussed later on.

\bigskip

Getting back now to the planar algebra structure of $\mathcal T_N$ and $\mathcal S_N$, which is something quite fundamental, worth being well understood, let us have here some more discussion. Generally speaking, the planar calculus for tensors is quite simple, and does not really require diagrams. Indeed, it suffices to imagine that the way various indices appear, travel around and dissapear is by following some obvious strings connecting them. Here are some illustrations for this principle, for the spin planar algebra $\mathcal S_N$:

\begin{example}
Identity, multiplication, inclusion.
\end{example}

The identity $1_k$ is the $(k,k)$-tangle having vertical strings only. The solutions of $\delta_{1_k}(x,y)=1$ being the pairs of the form $(x,x)$, this tangle $1_k$ acts by the identity:
$$1_k\begin{pmatrix}j_1 & \ldots & j_k\\ i_1 & \ldots & i_k\end{pmatrix}=\begin{pmatrix}j_1 & \ldots & j_k\\ i_1 & \ldots & i_k\end{pmatrix}$$

\medskip

The multiplication $M_k$ is the $(k,k,k)$-tangle having 2 input boxes, one on top of the other, and vertical strings only. It acts in the following way:
$$M_k\left( 
\begin{pmatrix}j_1 & \ldots & j_k\\ i_1 & \ldots & i_k\end{pmatrix}
\otimes\begin{pmatrix}l_1 & \ldots & l_k\\ m_1 & \ldots & m_k\end{pmatrix}
\right)=
\delta_{j_1m_1}\ldots \delta_{j_km_k}
\begin{pmatrix}l_1 & \ldots & l_k\\ i_1 & \ldots & i_k\end{pmatrix}$$

\medskip

The inclusion $I_k$ is the $(k,k+1)$-tangle which looks like $1_k$, but has one more vertical string, at right of the input box. Given $x$, the solutions of $\delta_{I_k}(x,y)=1$ are the elements $y$ obtained from $x$ by adding to the right a vector of the form $(^l_l)$, and so:
$${I_k}\begin{pmatrix}j_1 & \ldots & j_k\\ i_1 & \ldots & i_k\end{pmatrix}=
\sum_l\begin{pmatrix}j_1 & \ldots & j_k& l\\ i_1 & \ldots & i_k& l\end{pmatrix}$$

\medskip

Observe that $I_k$ is an inclusion of algebras, and that the various $I_k$ are compatible with each other. The union of the algebras $\mathcal S_N(k)$ is a graded algebra, denoted $\mathcal S_N$.

\bigskip

Along the same lines, some other important tangles are as follows:

\begin{example}
Expectation, Jones projection, trace.
\end{example}

The expectation $U_k$ is the $(k+1,k)$-tangle which looks like $1_k$, but has one more string, connecting the extra 2 input points, both at right of the input box:
$$U_k
\begin{pmatrix}j_1 & \ldots &j_k& j_{k+1}\\ i_1 & \ldots &i_k& i_{k+1}\end{pmatrix}=
\delta_{i_{k+1}j_{k+1}}
\begin{pmatrix}j_1 & \ldots & j_k\\ i_1 & \ldots & i_k\end{pmatrix}$$

Observe that $U_k$ is a bimodule morphism with respect to $I_k$.

\medskip

The Jones projection $E_k$ is a $(0,k+2)$-tangle, having no input box. There are $k$ vertical strings joining the first $k$ upper points to the first $k$ lower points, counting from left to right. The remaining upper 2 points are connected by a semicircle, and the remaining lower 2 points are also connected by a semicircle. We have the following formula:
$$E_k(1)=\sum_{ijl}\begin{pmatrix}i_1 & \ldots &i_k&j&j\\ i_1 & \ldots &i_k&l&l\end{pmatrix}$$

The elements $e_k=N^{-1}E_k(1)$ are projections, and define a representation of the infinite Temperley-Lieb algebra of index $N$ inside the inductive limit algebra $\mathcal S_N$.

\medskip

The trace $T_k$ is the $(k,0)$ tangle which ``closes the diagram'', by connecting upper points with lower points with noncrossing strings at right, in the obvious way:
$$T_k
\begin{pmatrix}j_1 & \ldots &j_k\\ i_1 & \ldots &i_k\end{pmatrix}=
\delta_{i_1j_1}\ldots\delta_{i_kj_k}$$

This tangle implements a trace on the planar algebra, and the expectations $U_k$ constructed above are then the conditional expectations with respect to this trace.

\bigskip

Finally, again along the same lines, we have the following two key tangles:

\index{rotation}

\begin{example}
Rotation, shift.
\end{example}

The rotation $R_k$ is the $(k,k)$-tangle which looks like $1_k$, but the first 2 input points are connected to the last 2 output points, and the same happens at right:
$$R_k=\begin{matrix}
\hskip 0.3mm\Cap \ |\ |\ |\ |\hskip -0.5mm |\cr
|\hskip -0.5mm |\hskip 10.3mm |\hskip -0.5mm |\cr
\hskip -0.3mm|\hskip -0.5mm |\ |\ |\ |\ \hskip -0.1mm\Cup
\end{matrix}$$

The action of $R_k$ on the standard basis is by rotation of the indices, as follows:
$$R_k(e_{i_1\ldots i_k})=e_{i_2i_3\ldots i_ki_1}$$

Thus $R_k$ acts by an order $k$ linear automorphism of $\mathcal S_N(k)$, also called rotation.

\medskip

As for the shift $S_k$, this is the $(k,k+2)$-tangle which looks like $1_k$, but has two more vertical strings, at left of the input box. This tangle acts as follows:
$$S_k\begin{pmatrix}j_1 & \ldots & j_k\\ i_1 & \ldots & i_k\end{pmatrix}=
\sum_{lm}\begin{pmatrix}l&m&j_1&\ldots&j_k\\l&m&i_1&\ldots&i_k\end{pmatrix}$$

Observe that $S_k$ is an inclusion of algebras, which is different from $I_{k+1}I_k$.

\bigskip

There are many other interesting examples of $k$-tangles, but in view of our present purposes, we can actually stop here, due to the following key fact, which basically reduces everything to the study of the above particular tangles, and that we will use many times in what follows, for the various planar algebra results that we will prove:

\begin{theorem}
The following tangles generate the set of all tangles, via gluing:
\begin{enumerate}
\item Multiplications, inclusions.

\item Expectations, Jones projections.

\item Rotations or shifts.
\end{enumerate}
\end{theorem}

\begin{proof}
As a first observation, the tangles in the statement are exactly those in the above examples, with the identity and trace tangles removed, due to the fact that these tangles won't bring anything new. Also, the statement itself consists in fact of 2 statements, depending on whether rotations and shifts are chosen in (3), with this being something technical, coming from the fact that we will need in what follows both these 2 statements. As for the proof, this is something elementary, obtained by ``chopping'' the various planar tangles into small pieces, belonging to the above list. See Jones \cite{jo3}.
\end{proof}

Finally, in order for things to be complete, we must talk as well about the $*$-structure. Once again this is constructed as in the easy quantum group calculus, as follows:
$$\begin{pmatrix}j_1 & \ldots & j_k\\ i_1 & \ldots & i_k\end{pmatrix}^*
=\begin{pmatrix}i_1 & \ldots & i_k\\ j_1 & \ldots & j_k\end{pmatrix}$$

Summarizing, the sequence of vector spaces $\mathcal S_N(k)=(\mathbb C^N)^{\otimes k}$ has a planar $*$-algebra structure, called spin planar algebra of index $N=|X|$. See Jones \cite{jo3}.

\bigskip

As a conclusion to all this, we have so far an abstract definition for the planar algebras, then two very basic examples, namely $TL_N$ and $FC_{N,M}$, where the elements of the planar algebra are actual diagrams, composing as the diagrams do, by gluing, and then two examples which are slightly more complicated, namely $\mathcal T_N$ and $\mathcal S_N$, where the planar algebra elements are tensors, composing according to the usual rules for the tensors.

\section*{14b. Higher commutants}

In relation now with subfactors, the result, which extends Theorem 14.1, and which was found by Jones in \cite{jo3}, almost 20 years after \cite{jo1}, is as follows:

\index{higher commutant}
\index{planar algebra}

\begin{theorem} 
Given a subfactor $A_0\subset A_1$, the collection $P=(P_k)$ of linear spaces 
$$P_k=A_0'\cap A_k$$
has a planar algebra structure, extending the planar algebra structure of $TL_N$.
\end{theorem}

\begin{proof}
We know from Theorem 14.1 that we have an inclusion as follows, coming from the basic construction, and with $TL_N$ itself being a planar algebra:
$$TL_N\subset P$$

Thus, the whole point is that of proving that the planar algebra structure of $TL_N$, obtained by composing diagrams, extends into a planar algebra structure of $P$. But this can be done via a long algebraic study, basically focusing on the basic tangles from Theorem 14.11, the idea here being as follows:

\medskip

(1) The multiplications and inclusions are the usual multiplications of the algebras $P_k=A_0'\cap A_k$, and the canonical inclusions $P_k\subset P_{k+1}$ between them.

\medskip

(2) The expectations and Jones projections are the usual expectations and Jones projections for the algebras $P_k=A_0'\cap A_k$, that we know from chapter 13.

\medskip

(3) As for rotations and shifts, things here are more tricky, the idea being that the algebras $P_k=A_0'\cap A_k$ have indeed some natural rotation and shift operations.

\medskip

In short, modulo some work needed for rotations and shifts, we know how the basic tangles act. Then, in order to make all the tangles act, we can invoke Theorem 14.11, along with a ``bubbling'' procedure in order to effectively construct the action, and to prove its uniqueness. And this ``bubbling'' procedure, which is something quite routine, but long and technical, taking about 10-20 pages, is explained in Jones' paper \cite{jo3}.
\end{proof}

So long for Jones' main result in \cite{jo3}. What has been said above is of course very far from a proof, and for this we refer of course to Jones' paper, but at least we have now an idea on what the result in \cite{jo3} really says. Regarding the reading of \cite{jo3}, which is a must-do thing if you want to fully understand subfactors, a few pieces of advice:

\bigskip

(1) Examples, examples, and examples again. The notion of planar algebra is something extremely general, somehow the idea being that planar algebras are to quantum groups what quantum groups are to groups, and 0 chances or almost to understand what Jones is doing in \cite{jo3}, without spending some substantial time on examples.

\bigskip

(2) And I'm saying this with knowledge of the matter, because back in 1999 when \cite{jo3} came out, I was postdoc at Berkeley with Jones and Voiculescu, and I saw quite a few  young people severely struggling with \cite{jo3}. For me things were easy because I was already familiar with examples coming from groups, and quantum groups.

\bigskip

(3) So, that would be a first way of getting introduced to the subject, via groups and quantum groups, and benefitting from what we already know from chapters 7-8, this is what we will do here, work out some examples coming from groups and quantum groups, as an introduction to Jones' paper \cite{jo3}, that you can read afterwards.

\bigskip

(4) But this is not the only way. As mentioned above, the subtlety comes from rotations and shifts, and understanding how these rotations and shifts work, directly in the subfactor context, in the spirit of what we did in chapter 13, is something that you can try too. Good references here are the texts of Ocneanu \cite{oc1}, \cite{oc2}.

\bigskip

Long story short, we are now into subtle mathematics, that takes some time to be understood. Back to work now, as a first illustration for Theorem 14.12, we have:

\index{Temperley-Lieb}
\index{Fuss-Catalan algebra}

\begin{theorem}
We have the following universality results:
\begin{enumerate}
\item The Temperley-Lieb algebra $TL_N$ appears inside the planar algebra of any subfactor $A\subset B$ having index $N$.

\item The Fuss-Catalan algebra $FC_{N,M}$ appears inside the planar algebra of any subfactor $A\subset B$, in the presence of an intermediate subfactor $A\subset C\subset B$.
\end{enumerate}
\end{theorem}

\begin{proof}
Here the first assertion is something that we already know, from Theorem 14.1, and the second assertion is something quite standard as well, by carefully working out the basic construction for $A\subset B$, in the presence of an intermediate subfactor $A\subset C\subset B$. For details here, we refer to the paper of Bisch and Jones \cite{bjo}.
\end{proof}

It is possible to prove as well that the tensor planar algebra $\mathcal T_N$ and the spin planar algebra $\mathcal S_N$ have similar universality properties, but this time being the biggest possible instead of the smallest possible, in the framework of some suitable fixed point subfactors. We will discuss all this in a moment, in the general context of fixed point subfactors.

\bigskip

All the above results raise the question on whether any planar algebra produces a subfactor. The answer here is yes, but with many subtleties, as follows:

\index{R}

\begin{theorem}
We have the following results:
\begin{enumerate}
\item Any planar algebra with positivity produces a subfactor.

\item In particular, we have $TL$ and $FC$ type subfactors.

\item In the amenable case, and with $A_1=R$, the correspondence is bijective.

\item In general, we must take $A_1=L(F_\infty)$, and we do not have bijectivity.

\item The axiomatization of $P$, in the case $A_1=R$, is not known.
\end{enumerate}
\end{theorem}

\begin{proof}
All this is quite heavy, mainly coming from the work of Popa in the 90s, using heavy functional analysis and operator theory \cite{po2}, \cite{po3}, \cite{po4}, completed with other papers like \cite{gjs}, \cite{jo3}, \cite{psh}, which are not any simpler either. In fact, understanding all this, with proofs, is a considerable investment, comparable to that of understanding the heavy papers of von Neumann and Connes \cite{co1}, \cite{co2}, \cite{mv3}, that we are stumbling upon all the time, in chapters 9-12. So, in the hope that you read that papers of von Neumann and Connes, in this way, reading Popa will look like a routine task.

\medskip

As an introduction to all this, following Ocneanu \cite{oc1}, \cite{oc2}, who first came upon such ideas, in the mid 80s, let us first talk about the finite depth case. The higher relative commutants $P_k=A_0\cap A_k'$ form an increasing sequence of algebras, as follows:
$$P_0\subset P_1\subset P_2\subset\ldots$$

The point now is that at each step, we have a copy of the basic construction which appears, in the sense that $P_{k+1}$ consists of a copy of the basic construction for $P_{k-1}\subset P_k$, colloquially called ``old stuff'', and of more things, called ``new stuff''. In case there is no new stuff inside $P_{k+1}$, there is no new stuff either inside $P_{k+2},P_{k+3},\ldots\,$, and the subfactor or planar algebra is called of ``finite depth''. And there are many examples here, such as the Ocneanu subfactors, the general idea being that finite depth means that the underlying ``generalized quantum group'', whatever that beast might be, is finite.

\medskip

The problem is now, given a planar algebra which has finite depth, how to construct a subfactor out of it. Due to the finite depth assumption, our data is simply:
$$P_0\subset P_1\subset\ldots\subset P_k$$

That is, our data is just a finite dimensional graded algebra $P_k$, and we are here into usual algebra, be that of quite complicated type. And Ocneanu's solution \cite{oc1}, \cite{oc2} was that of building out of this data, via various algebraic procedures, some further finite dimensional algebras, then taking inductive limits and closing under the weak topology, as to end up with a subfactor of type $A_0\subset A_1$, with $A_0\simeq A_1\simeq R$.

\medskip

This was for the general idea, in the finite depth case, but in practice, the above-mentioned ``various algebraic procedures'' are something quite complicated, involving a certain technical notion of ``commuting square'', which is something specialized, that we will discuss in chapter 15 below, and with the whole thing, complete theorem coming with complete proof, being something done by Popa, some time after Ocneanu, in \cite{po2}. 

\medskip

With this understood, and getting back to our theorem, all the items (1-5) there are extensions of this construction of Ocneanu and Popa, the idea being as follows:

\medskip

(1) As already mentioned in the comments after Definition 14.2, our planar algebra axioms here are something quite simplified, based on \cite{jo3}. However, by getting back to Theorem 14.12, and carefully looking at the planar algebras there, appearing from subfactors, the conclusion is that these subfactor planar algebras satisfy a number of supplementary ``positivity'' conditions, basically coming from the positivity of the ${\rm II}_1$ factor trace. And the point now is that, with these positivity conditions axiomatized, we reach to something which is equivalent to Popa's axiomatization of the lattice of higher relative commutants $A_i'\cap A_j$ of the finite index subfactors \cite{po4}, obtained in the 90s via heavy functional analysis. For the story here, and details, we refer to Jones \cite{jo3}.

\medskip

(2) We have been a bit quick in the above, and before anything, let us mention that our 4 main examples of planar algebras, namely $TL_N$ and $FC_{N,M}$, and then $\mathcal T_N$ and $\mathcal S_N$ too, do satisfy the positivity requirements needed in (1). Thus, there are subfactors associated to all of them. In practice now, the existence of the $TL_N$ subfactors, also known as ``$A_\infty$ subfactors'', is something which was known for some time, since some early work of Popa on the subject. As for the existence of the $FC_{N,M}$ subfactors, this can be shown by using the intermediate subfactor picture, $A\subset C\subset B$, by composing two $A_\infty$ subfactors of suitable indices, $A\subset C$ and $C\subset B$. For the story here, we refer to \cite{bjo}, \cite{jo3}.

\medskip

(3) This is something fairly heavy, as it is always the case with operator algebra results about hyperfiniteness and amenability, due to Popa \cite{po2}, \cite{po3}.

\medskip

(4) This is something more fashionable and recent, obtained by further building on the above-mentioned old constructions of Popa, and we refer here to \cite{psh}, \cite{gjs}.

\medskip

(5) This is the big open question in subfactors. The story here goes back to Jones' original  paper \cite{jo1}, which contains at the end the question, due to Connes, of finding the possible values of the index for the irreducible subfactors of $R$. This question, which certainly looks much easier than (5) in the statement, is in fact still open, now 40 years after its formulation, and with on one having any valuable idea of dealing with it.
\end{proof}

We refer to the original papers of Popa, and then to more recent papers by Jones, Popa and their collaborators for details on the above, which is quite heavy material.

\section*{14c. Fixed points}

We discuss now the connection of all the above with the main examples of subfactors. We recall from chapter 13 that the main examples of subfactors are all of integer index, and appear as fixed point subfactors, according to the following result:

\begin{theorem}
Let $G$ be a compact quantum group, and $G\to Aut(P)$ be a minimal action on a ${\rm II}_1$ factor. Consider a Markov inclusion of finite dimensional algebras
$$B_0\subset B_1$$
and let $G\to Aut(B_1)$ be an action which leaves invariant $B_0$ and which is such that its restrictions to the centers of $B_0$ and $B_1$ are ergodic. We have then a subfactor
$$(B_0\otimes P)^G\subset (B_1\otimes P)^G$$
of index $N=[B_1:B_0]$, called generalized Wassermann subfactor, whose Jones tower is 
$$(B_1\otimes P)^G\subset(B_2\otimes P)^G\subset(B_3\otimes P)^G\subset\ldots$$
where $\{ B_i\}_{i\geq 1}$ are the algebras in the Jones tower for $B_0\subset B_1$, with the canonical actions of $G$ coming from the action $G\to Aut(B_1)$, and whose planar algebra is given by:
$$P_k=(B_0'\cap B_k)^G$$
These subfactors generalize the Jones, Ocneanu, Wassermann and Popa subfactors.
\end{theorem}

\begin{proof}
This is something that we know well from chapter 13, whose proof basically comes by generalizing, several times, the results of Wassermann in \cite{was}.
\end{proof}

In view of the above result, what we have to do in relation with such subfactors is to further interpret the last formula there, that of the planar algebra, namely:
$$P_k=(B_0'\cap B_k)^G$$

To be more precise, we will show here that, under suitable assumptions on the original inclusion $B_0\subset B_1$, we can associate a certain combinatorial planar algebra $P(B_0\subset B_1)$ to this inclusion, and then the planar algebra associated to the fixed point subfactor itself appears as a fixed point subalgebra of this planar algebra, as follows:
$$P=P(B_0\subset B_1)^G$$

This is something quite technical, and we will do this in two steps. First we will explain, following Jones' paper \cite{jo4}, how to associate a planar algebra $P(B_0\subset B_1)$ to an inclusion of algebras $B_0\subset B_1$. And then we will explain, following \cite{ba1} and subsequent papers, and notably \cite{ba2}, how to prove the above formula $P=P(B_0\subset B_1)^G$.

\bigskip

Getting started now, the idea will be that $P(B_0\subset B_1)$ appears as a joint generalization of the spin and tensor planar algebras, discussed above, which appear as follows:
$$\mathcal S_N=P(\mathbb C\subset\mathbb C^N)$$
$$\mathcal T_N=P(\mathbb C\subset M_N(\mathbb C))$$

Thus, our first task will be that of getting back to the Markov inclusions $B_0\subset B_1$, from chapter 13, and further discuss the combinatorics of their basic construction, with planar algebra ideas in mind. As in chapter 13, it is most convenient to denote such inclusions by $A\subset B$, at least at a first stage of their study. Following the book of Goodman, de la Harpe and Jones \cite{ghj}, which is the standard reference for such things, we first have:

\begin{definition}
Associated to an inclusion $A\subset B$ of finite dimensional algebras are the following objects:
\begin{enumerate}
\item The column vector $(a_i)\in\mathbb N^s$ given by $A=\oplus_{i=1}^sM_{a_i}(\mathbb C)$.

\item The column vector $(b_j)\in\mathbb N^t$ given by $B=\oplus_{j=1}^tM_{b_j}(\mathbb C)$.

\item The inclusion matrix $(m_{ij})\in M_{s\times t}(\mathbb N)$, satisfying $m^ta=b$.
\end{enumerate}
\end{definition}

To be more precise here, in what regards the inclusion matrix, each minimal idempotent in $M_{a_i}(\mathbb C)\subset A$ splits, when regarded as an element of $B$, as a sum of minimal idempotents of $B$, and $m_{ij}\in\mathbb N$ is the number of such idempotents from $M_{b_j}(\mathbb C)$. We have the following result, bringing traces into picture:

\begin{proposition}
For an inclusion $A\subset B$, the following are equivalent:
\begin{enumerate}
\item $A\subset B$ commutes with the canonical traces.

\item We have $mb=ra$, where $r=||b||^2/||a||^2$.
\end{enumerate}
\end{proposition}

\begin{proof}
The weight vectors of the canonical traces of $A,B$ are given by:
$$\tau_i=\frac{a_i^2}{||a||^2}\quad,\quad 
\tau_j=\frac{b_j^2}{||b||^2}$$

We can plug these values into the following standard compatibility formula:
$$\frac{\tau_i}{a_i}=\sum_jm_{ij}\cdot\frac{\tau_j}{b_j}$$

We obtain in this way the condition in the statement.
\end{proof}

We will need as well the following basic facts, also from \cite{ghj}:

\index{Bratteli diagram}
\index{bipartite graph}

\begin{definition}
Associated to an inclusion of finite dimensional algebras $A\subset B$, with inclusion matrix $m\in M_{s\times t}(\mathbb N)$, are:
\begin{enumerate}
\item The Bratteli diagram: this is the bipartite graph $\Gamma$ having as vertices the sets $\{1,\ldots,s\}$ and $\{1,\ldots,t\}$, the number of edges between $i,j$ being $m_{ij}$.

\item The basic construction: this is the inclusion $B\subset A_1$ obtained from $A\subset B$ by reflecting the Bratteli diagram. 

\item The Jones tower: this is the tower of algebras $A\subset B\subset A_1\subset B_1\subset\ldots$ obtained by iterating the basic construction.
\end{enumerate}
\end{definition}

We know that for a Markov inclusion $A\subset B$ we have $m^ta=b$ and $mb=ra$, and so $mm^ta=ra$, which gives an eigenvector for the square matrix $mm^t\in M_s(\mathbb N)$. When this latter matrix has positive entries, by Perron-Frobenius we obtain:
$$||mm^t||=r$$

This equality holds in fact without assumptions on $m$, and we have: 

\index{Markov inclusion}

\begin{theorem}
Let $A\subset B$ be Markov, with inclusion matrix $m\in M_{s\times t}(\mathbb N)$.
\begin{enumerate}
\item $r=\dim(B)/\dim(A)$ is an integer.

\item $||m||=||m^t||=\sqrt{r}$.

\item $||\ldots mm^tmm^t\ldots||=r^{k/2}$, for any product of lenght $k$.
\end{enumerate}
\end{theorem}

\begin{proof}
Consider the vectors $a,b$, as in Definition 14.16. We know from definitions and from Proposition 14.17 that we have:
$$b=m^ta\quad,\quad 
mb=ra\quad,\quad 
r=||b||^2/||a||^2$$

(1) If we construct as above the Jones tower for $A\subset B$, we have, for any $k$:
$$\frac{\dim B_k}{\dim A_k}
=\frac{\dim A_k}{\dim B_{k-1}}
=r$$

On the other hand, we have as well the following well-known formula:
$$\lim_{k\to\infty}(\dim A_k)^{1/2k}
=\lim_{k\to\infty}(\dim B_k)^{1/2k}
=||mm^t||$$

By combining these two formulae we obtain the following formula:
$$||mm^t||=r$$

But from $r\in\mathbb Q$ and $(mm^t)^ka=r^ka$ for any $k\in\mathbb N$, we get $r\in\mathbb N$, and we are done.

\medskip

(2) This follows from the above equality $||mm^t||=r$, and from the following standard equalities, for any real rectangular matrix $r$:
$$||m||^2
=||m^t||^2
=||mm^t||$$

(3) Let $n$ be the length $k$ word in the statement. First, by applying the norm and by using the formula $||m||=||m^t||=\sqrt{r}$, we obtain the following inequality:
$$||n||\leq r^{k/2}$$

For the converse inequality, assume first that $k$ is even. Then $n$ has either $a$ or $b$ as eigenvector, depending on whether $n$ begins with a $m$ or with a $m^t$, in both cases with eigenvalue $r^{k/2}$, and this gives the desired inequality, namely:
$$||n||\geq r^{k/2}$$

Assume now that $k$ is odd, and let $\circ\in\{1,t\}$ be such that $n'=m^\circ n$ is alternating. Since $n'$ has even length, we already know that we have:
$$||n'||=r^{(k+1)/2}$$

On the other hand, we have as well the following estimate:
$$||n'||
\leq||m^\circ||\cdot ||n||
=\sqrt{r}||n||$$

But this gives the reverse inequality $||n||\geq r^{k/2}$, as desired.
\end{proof}

The point now is that for a Markov inclusion, the basic construction and the Jones tower have a particularly simple form. Let us first work out the basic construction:

\begin{proposition}
The basic construction for a Markov inclusion $i:A\subset B$ of index $r\in\mathbb N$ is the inclusion $j:B\subset A_1$ obtained as follows:
\begin{enumerate}
\item $A_1=M_r(\mathbb C)\otimes A$, as an algebra.

\item $j:B\subset A_1$ is given by $mb=ra$.

\item $ji:A\subset A_1$ is given by $(mm^t)a=ra$.
\end{enumerate}
\end{proposition}

\begin{proof}
With notations from the above, the weight vector of the algebra $A_1$ appearing from the basic construction is $mb=ra$, and this gives the result.
\end{proof}

We fix a Markov inclusion $i:A\subset B$. We have the following result:

\begin{proposition}
The Jones tower associated to a Markov inclusion $i:A\subset B$, denoted as follows, with alternating letters,
$$A\subset B\subset A_1\subset B_1\subset\ldots$$
is given by the following formulae:
\begin{enumerate}
\item $A_k=M_r(\mathbb C)^{\otimes k}\otimes A$.

\item $B_k=M_r(\mathbb C)^{\otimes k}\otimes B$.

\item $A_k\subset B_k$ is $id_k\otimes i$.

\item $B_k\subset A_{k+1}$ is $id_k\otimes j$.
\end{enumerate}
\end{proposition}

\begin{proof}
This follows from Proposition 14.20, with the remark that if $i:A\subset B$ is Markov, then so is its basic construction $j:B\subset A_1$.
\end{proof}

Regarding now the relative commutants for this tower, we have here:

\begin{proposition}
The relative commutants for the Jones tower 
$$A\subset B\subset A_1\subset B_1\subset\ldots$$
associated to a Markov inclusion $A\subset B$ are given by:
\begin{enumerate}
\item $A_s'\cap A_{s+k}=M_r(\mathbb C)^{\otimes k}\otimes (A'\cap A)$.

\item $A_s'\cap B_{s+k}=M_r(\mathbb C)^{\otimes k}\otimes (A'\cap B)$.

\item $B_s'\cap A_{s+k}=M_r(\mathbb C)^{\otimes k}\otimes (B'\cap A)$.

\item $B_s'\cap B_{s+k}=M_r(\mathbb C)^{\otimes k}\otimes (B'\cap B)$.
\end{enumerate}
\end{proposition}

\begin{proof}
The above assertions are all elementary, as follows:

\medskip

(1,2) These assertions both follow from Proposition 14.21. 

\medskip

(3) In order to prove the formula in the statement, observe first that we have:
\begin{eqnarray*}
B'\cap A_1
&=&(B'\cap B_1)\cap A_1\\
&=&(M_r(\mathbb C)\otimes Z(B))\cap (M_r(\mathbb C)\otimes A)\\
&=&M_r(\mathbb C)\otimes (B'\cap A)
\end{eqnarray*}

But this proves the assertion at $s=0,k=1$, and the general case follows from it.

\medskip

(4) This is again clear, once again coming from Proposition 14.21.
\end{proof}

In order to further refine all this, let us formulate the following key definition:

\begin{definition}
We say that a Markov inclusion $A\subset B$ is abelian if $[A,B]=0$, with the commutant being computed inside $B$.
\end{definition}

In other words, we are asking for the commutation relation $ab=ba$, for any $a\in A$ and $b\in B$. Note that this is the same as asking that $B$ is an $A$-algebra, $A\subset Z(B)$. As basic examples, observe that all inclusions with $A=\mathbb C$ or with $B=\mathbb C^n$ are abelian. The point with this notion is that it leads to the following simple statement:

\begin{proposition}
With $\tilde{B}_k=M_r(\mathbb C)^{\otimes k}\otimes Z(B)$, the relative commutants for the Jones tower $A\subset B\subset A_1\subset B_1\subset\ldots$ of an abelian inclusion are given by:
\begin{enumerate}
\item $A_s'\cap A_{s+k}=A_k$.

\item $A_s'\cap B_{s+k}=B_k$.

\item $B_s'\cap A_{s+k}=A_k$.

\item $B_s'\cap B_{s+k}=\tilde{B}_k$.
\end{enumerate}
\end{proposition}

\begin{proof}
This follows from the fact that for an abelian inclusion we have:
$$Z(A)=A\quad,\quad 
A'\cap B=B\quad,\quad
B'\cap A=A$$

Thus, we are led to the conclusion in the statement.
\end{proof}

Getting back now to the fixed point subfactors, from Theorem 14.15, we can improve the planar algebra computation there, in the abelian case, as follows:

\begin{theorem}
The commutants for the tower $N\subset M\subset N_1\subset M_1\subset\ldots$ associated to an abelian fixed point subfactor $(A\otimes P)^G\subset(B\otimes P)^G$ are:
\begin{enumerate}
\item $N_s'\cap N_{s+k}=A_k^G$.

\item $N_s'\cap M_{s+k}=B_k^G$.

\item $M_s'\cap N_{s+k}=A_k^G$.

\item $M_s'\cap M_{s+k}=\tilde{B}_k^G$.
\end{enumerate}
\end{theorem}

\begin{proof}
This follows indeed by combining the planar algebra computation from Theorem 14.15 with the result about abelian inclusions from Proposition 14.24.
\end{proof}

In order to further advance now, the idea will be that of associating to the original inclusion $B_0\subset B_1$ a certain combinatorial planar algebra $P(B_0\subset B_1)$, as for the planar algebra associated to the fixed point subfactor itself to appear as follows:
$$P=P(B_0\subset B_1)^G$$

As already mentioned, the idea will be that $P(B_0\subset B_1)$ appears as a joint generalization of the spin and tensor planar algebras, which appear as follows:
$$\mathcal S_N=P(\mathbb C\subset\mathbb C^N)$$
$$\mathcal T_N=P(\mathbb C\subset M_N(\mathbb C))$$

In practice now, we will need for all this the notion of planar algebra of a bipartite graph, generalizing the algebras $\mathcal S_N,\mathcal T_N$, constructed by Jones in \cite{jo4}. So, let $\Gamma$ be a bipartite graph, with vertex set $\Gamma_a\cup\Gamma_b$. It is useful to think of $\Gamma$ as being the Bratteli diagram of an inclusion $A\subset B$, in the sense of Definition 14.16. Our first task is to define the graded vector space $P$. Since the elements of $P$ will be subject to a planar calculus, it is convenient to introduce them ``in boxes'', as follows:

\begin{definition}
Associated to $\Gamma$ is the abstract vector space $P_k$ spanned by the $2k$-loops based at points of $\Gamma_a$. The basis elements of $P_k$ will be denoted 
$$x=\begin{pmatrix}
e_1&e_2&\ldots&e_k\\
e_{2k}&e_{2k-1}&\ldots &e_{k+1}
\end{pmatrix}$$
where $e_1,e_2,\ldots,e_{2k}$ is the sequence of edges of the corresponding $2k$-loop.
\end{definition}

Consider now the adjacency matrix of $\Gamma$, which is of the following type:
$$M=\begin{pmatrix}0&m\\ m^t&0\end{pmatrix}$$

We pick an $M$-eigenvalue $\gamma\neq 0$, and then a $\gamma$-eigenvector, as follows:
$$\eta:\Gamma_a\cup\Gamma_b\to\mathbb C-\{0\}$$

With this data in hand, we have the following construction, due to Jones \cite{jo4}:

\begin{definition}
Associated to any tangle is the multilinear map 
$$T(x_1\otimes\ldots\otimes x_r)=\gamma^c\sum_x\delta(x_1,\ldots ,x_r,x)\prod_m\mu(e_m)^{\pm 1}x$$
where the objects on the right are as follows:
\begin{enumerate}
\item The sum is over the basis of $P_k$, and $c$ is the number of circles of $T$.

\item $\delta=1$ if all strings of $T$ join pairs of identical edges, and $\delta=0$ if not.

\item The product is over all local maxima and minima of the strings of $T$. 

\item $e_m$ is the edge of $\Gamma$ labelling the string passing through $m$ (when $\delta=1$).

\item $\mu(e)=\sqrt{\eta(e_f)/\eta(e_i)}$, where $e_i,e_f$ are the initial and final vertex of $e$.

\item The $\pm$ sign is $+$ for a local maximum, and $-$ for a local minimum.
\end{enumerate}
\end{definition}

This looks quite similar to the calculus for the tensor and spin planar algebras. Let us work out the precise formula of the action, for 6 carefully chosen tangles:

\medskip

(1) Let us look first at the identity $1_k$. This tangle acts by the identity:
$$1_k\begin{pmatrix}f_1&\ldots&f_k\\ e_1&\ldots&e_k\end{pmatrix}=
\begin{pmatrix}f_1&\ldots&f_k\\ e_1&\ldots&e_k\end{pmatrix}$$

(2) The multiplication tangle $M_k$ acts as follows:
$$M_k\left( 
\begin{pmatrix}f_1&\ldots&f_k\\ e_1&\ldots&e_k\end{pmatrix}
\otimes\begin{pmatrix}h_1&\ldots&h_k\\ g_1&\ldots&g_k\end{pmatrix}
\right)=
\delta_{f_1g_1}\ldots\delta_{f_kg_k}
\begin{pmatrix}h_1&\ldots&h_k\\ e_1&\ldots& e_k\end{pmatrix}$$

(3) Regarding now the inclusion $I_k$, the formula here is:
$$I_k\begin{pmatrix}f_1&\ldots&f_k\\ e_1&\ldots&e_k\end{pmatrix}=
\sum_g\begin{pmatrix}f_1&\ldots&f_k&g\\ e_1&\ldots&e_k&g\end{pmatrix}$$

(4) The expectation tangle $U_k$ acts with a spin factor, as follows:
$$U_k\begin{pmatrix}f_1&\ldots&f_k&h\\ e_1&\ldots&e_k&g\end{pmatrix}=
\delta_{gh}\mu(g)^2\begin{pmatrix}f_1&\ldots&f_k\\ e_1&\ldots&e_k\end{pmatrix}$$

(5) For the Jones projection $E_k$, the formula is as follows:
$$E_k(1)=\sum_{egh}\mu(g)\mu(h)\begin{pmatrix}e_1&\ldots &e_k&h&h\\ e_1&\ldots&e_k&g&g\end{pmatrix}$$

(6) As for the shift $J_k$, its action is given by:
$$J_k\begin{pmatrix}f_1&\ldots&f_k\\ e_1&\ldots&e_k\end{pmatrix}=
\sum_{gh}\begin{pmatrix}g&h&f_1&\ldots&f_k\\ g&h&e_1&\ldots&e_k\end{pmatrix}$$

Summarizing, we have here formulae which are quite similar to those for the tensor and spin planar algebras. We have the following result, from Jones' paper \cite{jo4}:

\begin{theorem}
The graded linear space $P=(P_k)$, together with the action of the planar tangles given above, is a planar algebra. 
\end{theorem}

\begin{proof}
This is something which is quite routine, starting from the above study of the main planar algebra tangles, which can be proved by using Theorem 14.11. Also, let us mention that  all this generalizes the previous constructions of the spin and tensor planar algebras $\mathcal S_N,\mathcal T_N$, which appear respectively from the Bratteli diagrams of the inclusions $\mathbb C\subset\mathbb C^N$ and $\mathbb C\subset M_N(\mathbb C)$. For full details on all this, we refer to Jones \cite{jo4}.
\end{proof}

Let us go back now to the Markov inclusions $A\subset B$, as before. We have here the following result, regarding such inclusions, also from Jones' paper \cite{jo4}:

\index{bipartite graph algebra}

\begin{theorem}
The planar algebra associated to the graph of $A\subset B$, with eigenvalue $\gamma=\sqrt{r}$ and eigenvector $\eta(i)=a_i/\sqrt{\dim A}$, $\eta(j)=b_j/\sqrt{\dim B}$, is as follows:
\begin{enumerate}
\item The graded algebra structure is given by $P_{2k}=A'\cap A_k$, $P_{2k+1}=A'\cap B_k$.

\item The elements $e_k$ are the Jones projections for $A\subset B\subset A_1\subset B_1\subset\ldots$

\item The expectation and shift are given by the above formulae.
\end{enumerate}
\end{theorem}

\begin{proof}
As a first observation, $\eta$ is indeed a $\gamma$-eigenvector for the adjacency matrix of the graph. Indeed, we have the following formulae:
$$m^ta=b\quad,\quad 
mb=ra\quad,\quad 
\sqrt{r}=||b||/||a||$$

By using these formulae, we have the following computation:
\begin{eqnarray*}
\begin{pmatrix}0&m\\ m^t&0\end{pmatrix}\begin{pmatrix}a/||a||\\ b/||b||\end{pmatrix}
&=&\begin{pmatrix}\gamma^2a/||b||\\ b/||a||\end{pmatrix}\\
&=&\gamma\begin{pmatrix}\gamma a/||b||\\ b/\gamma||a||\end{pmatrix}\\
&=&\gamma\begin{pmatrix}a/||a||\\ b/||b||\end{pmatrix}
\end{eqnarray*}

Since the algebra $A$ was supposed abelian, the Jones tower algebras $A_k,B_k$ are simply the span of the $4k$-paths, respectively $4k+2$-paths on $\Gamma$, starting at points of $\Gamma_a$. With this description in hand, when taking commutants with $A$ we have to just have to restrict attention from paths to loops, and we obtain the above spaces $P_{2k},P_{2k+1}$. See \cite{jo4}.
\end{proof}

In the particular case of the inclusions satisfying $[A,B]=0$, we have:

\begin{proposition}
The ``bipartite graph'' planar algebra $P(A\subset B)$ associated to an abelian inclusion $A\subset B$ can be described as follows:
\begin{enumerate}
\item As a graded algebra, this is the Jones tower $A\subset B\subset A_1\subset B_1\subset\ldots$

\item The Jones projections and expectations are the usual ones for this tower.

\item The shifts correspond to the canonical identifications $A_1'\cap P_{k+2}=P_k$.
\end{enumerate}
\end{proposition}

\begin{proof}
The first assertion is a reformulation of Theorem 14.28 in the abelian case, by using the identifications $A'\cap A_k=A_k$ and $A'\cap B_k=B_k$ from Proposition 14.24. The assertion on Jones projections follows as well from Theorem 14.28, and the assertion on expectations follows from the fact that their composition is the usual trace. Regarding now the third assertion, let us recall first from Proposition 14.24 that we have indeed identifications $A_1'\cap A_{k+1}=A_k$ and $A_1'\cap B_{k+1}=B_k$. By using the path model for these algebras, as in the proof of Theorem 14.28, we obtain the result.
\end{proof}

In order to formulate now our main result, regarding the subfactors associated to the compact quantum groups $G$, we will need a few abstract notions. Let us start with:

\begin{definition}
Let $P_1,P_2$ be two finite dimensional algebras, coming with coactions $\alpha_i:P_i\to P_i\otimes L^\infty(G)$, and let $T:P_1\to P_2$ be a linear map.
\begin{enumerate}
\item We say that $T$ is $G$-equivariant if $(T\otimes id)\alpha_1=\alpha_2T$.

\item We say that $T$ is weakly $G$-equivariant if $T(P_1^G)\subset P_2^G$.
\end{enumerate}
\end{definition}

Consider now a planar algebra $P=(P_k)$. The annular category over $P$ is the collection of maps $T:P_k\to P_l$ coming from the ``annular'' tangles, having at most one input box. These maps form sets $Hom(k,l)$, and these sets form a category. We have:

\begin{definition}
A coaction of $L^\infty(G)$ on $P$ is a graded algebra coaction 
$$\alpha:P\to P\otimes L^\infty(G)$$
such that the annular tangles are weakly $G$-equivariant.
\end{definition}

This is something a bit technical, coming out of the known examples that we have. In fact, as we will show below, the examples are basically those coming from actions of quantum groups on Markov inclusions $A\subset B$, under the assumption $[A,B]=0$. For the moment, at the generality level of Definition 14.31, we have:

\begin{proposition}
If $G$ acts on a planar algebra $P$, then $P^G$ is a planar algebra.
\end{proposition}

\begin{proof}
The weak equivariance condition tells us that the annular category is contained in the suboperad $\mathcal P'\subset\mathcal P$ consisting of tangles which leave invariant $P^G$. On the other hand the multiplicativity of $\alpha$ gives $M_k\in\mathcal P'$, for any $k$. Now since $\mathcal P$ is generated by multiplications and annular tangles, we get $\mathcal P'=\mathcal P$, and we are done.
\end{proof}

Let us go back now to the abelian inclusions. We have the following key result: 

\begin{proposition}
If $G$ acts on an abelian inclusion $A\subset B$, the canonical extension of this coaction to the Jones tower is a coaction of $G$ on the planar algebra $P(A\subset B)$.
\end{proposition}

\begin{proof}
We know from the above that, as a graded algebra, $P=P(A\subset B)$ coincides with the Jones tower for our inclusion, denoted as follows:
$$A\subset B\subset A_1\subset B_1\subset\ldots$$

Thus the coaction in the statement can be regarded as a graded coaction, as follows:
$$\alpha:P\to P\otimes L^\infty(G)$$

In order to finish, we have to prove that the annular tangles are weakly equivariant, as in Definition 14.31, and this can be done as follows:

\medskip

(1) First, since the annular category is generated by $I_k,E_k,U_k,J_k$, we just have to prove that these 4 particular tangles are weakly equivariant. Now since $I_k,E_k,U_k$ are plainly equivariant, by construction of the coaction of $G$ on the Jones tower, it remains to prove that the shift $J_k$ is weakly equivariant.

\medskip

(2) We know that the image of the fixed point subfactor shift $J_k'$ is formed by the $G$-invariant elements of the relative commutant $A_1'\cap P_{k+2}=P_k$. Now since this commutant is the image of the planar shift $J_k$, we have $Im(J_k)=Im(J_k')$, and this gives the result.
\end{proof}

With the above result in hand, we can now prove:

\begin{proposition}
Assume that $G$ acts on an abelian inclusion $A\subset B$. Then the graded vector space of fixed points $P(A\subset B)^G$ is a planar subalgebra of $P(A\subset B)$.
\end{proposition}

\begin{proof}
This follows indeed from Proposition 14.33 and Proposition 14.34.
\end{proof}

We are now in position of stating and proving a main result, from \cite{ba2}:

\begin{theorem}
In the abelian case, the planar algebra of the fixed point subfactor
$$(A\otimes P)^G\subset (B\otimes P)^G$$
is the fixed point algebra $P(A\subset B)^G$ of the bipartite graph algebra $P(A\subset B)$.
\end{theorem}

\begin{proof}
This basically follows from what we have, as follows:

\medskip

(1) Let $P=P(A\subset B)$, and let $Q$ be the planar algebra of the fixed point subfactor. We know that we have an equality of graded algebras $Q=P^G$. Thus, it remains to prove that the planar algebra structure on $Q$ coming from the fixed point subfactor agrees with the planar algebra structure of $P$, coming from Proposition 14.30.

\medskip

(2) Since $\mathcal P$ is generated by the annular category $\mathcal A$ and by the multiplication tangles $M_k$, we just have to check that the annular tangles agree on $P,Q$. Moreover, since $\mathcal A$ is generated by $I_k,E_k,U_k,J_k$, we just have to check that these tangles agree on $P,Q$.

\medskip

(3) We know that $Q\subset P$ is an inclusion of graded algebras, that all the Jones projections  for $P$ are contained in $Q$, and that the conditional expectations agree. Thus the tangles $I_k,E_k,U_k$ agree on $P,Q$, and the only verification left is that for the shift $J_k$.

\medskip

(4) Now by using either the axioms of Popa in \cite{po4}, or the construction of Jones in \cite{jo4}, it is enough to show that the image of the subfactor shift $J_k'$ coincides with that of the planar shift $J_k$. But this follows as in the proof of Proposition 14.34.
\end{proof}

\section*{14d. Tannakian results}

We discuss here some converses to the above results, which are rather specialized results, of Tannakian nature. We will first prove that any quantum permutation group $G\subset S_N^+$ produces a planar subalgebra of $\mathcal S_N$. In order to do so, we first have:

\index{quantum permutation}
\index{coaction}
\index{spin planar algebra}

\begin{theorem}
Given a quantum permutation group $G\subset S_N^+$, consider the associated coaction map on $C(X)$, where $X=\{1,\ldots,N\}$,
$$\Phi:C(X)\to C(X)\otimes C(G)\quad,\quad 
e_j\to\sum_je_j\otimes u_{ji}$$
and then consider the tensor powers of this coaction, which are the following linear maps:
$$\Phi^k:C(X^k)\to C(X^k)\otimes C(G)\quad,\quad 
e_{i_1\ldots i_k}\to\sum_{j_1\ldots j_k}e_{j_1\ldots j_k}\otimes u_{j_1i_1}\ldots u_{j_ki_k}$$
The fixed point spaces of these latter coactions are then given by the formula
$$P_k=Fix(u^{\otimes k})$$
and form a planar subalgebra of the spin planar algebra $\mathcal S_N$.
\end{theorem}

\begin{proof}
This can be done in several steps, as follows:

\medskip

(1) Since the map $\Phi$ is a coaction, its tensor powers $\Phi^k$ are coactions too, and at the level of the fixed point algebras we have the following formula, which is standard:
$$P_k=Fix(u^{\otimes k})$$

(2) In order to prove now the planar algebra assertion, we use the presentation result for the spin planar algebras established before, involving the multiplications, inclusions, expectations, Jones projections and rotations.

\medskip

(3) Consider the rotation $R_k$. Rotating, then applying $\Phi^k$, and rotating backwards by $R_k^{-1}$ is the same as applying $\Phi^k$, then rotating each $k$-fold product of coefficients of $\Phi$. 

\medskip

(4) Thus the elements obtained by rotating, then applying $\Phi^k$, or by applying $\Phi^k$, then rotating, differ by a sum of Dirac masses tensored with commutators in $A=C(G)$:
$$\Phi^kR_k(x)-(R_k\otimes id)\Phi^k(x)\in C(X^k)\otimes [A,A]$$

(5) Now let $\int_A$ be the Haar functional of $A$, and consider the conditional expectation onto the fixed point algebra $P_k$, which is given by the following formula:
$$\phi_k=\left(id\otimes\int_A\right)\Phi^k$$

The square of the antipode being the identity, the Haar integration $\int_A$ is a trace, so it vanishes on commutators. Thus $R_k$ commutes with $\phi_k$:
$$\phi_kR_k=R_k\phi_k$$

(6) The commutation relation $\phi_kT=T\phi_l$ holds in fact for any $(l,k)$-tangle $T$. These tangles are called annular, and the proof is by verification on generators of the annular category. In particular we obtain, for any annular tangle $T$:
$$\phi_kT\phi_l=T\phi_l$$

(7) We conclude from this that the annular category is contained in the suboperad $\mathcal P'\subset\mathcal P$ of the planar operad consisting of tangles $T$ satisfying the following condition, where $\phi =(\phi_k)$, and where $i(.)$ is the number of input boxes:
$$\phi T\phi^{\otimes i(T)}=T\phi^{\otimes i(T)}$$

On the other hand the multiplicativity of $\Phi^k$ gives $M_k\in\mathcal P'$. Since $\mathcal P$ is generated by multiplications and annular tangles, it follows that we have:
$$\mathcal P'=P$$

(8) Thus for any tangle $T$ the corresponding multilinear map between spaces $P_k(X)$ restricts to a multilinear map between spaces $P_k$. In other words, the action of the planar operad $\mathcal P$ restricts to $P$, and makes it a subalgebra of $\mathcal S_N$, as claimed.
\end{proof}

As a second result now, completing our study, we have:

\index{quantum permutation}
\index{Tannakian duality}
\index{annular category}

\begin{theorem}
Given a subalgebra $Q\subset\mathcal S_N$, there is a unique quantum group 
$$G\subset S_N^+$$
whose associated planar algebra is $Q$.
\end{theorem}

\begin{proof}
The idea is that this will follow by applying Tannakian duality to the annular category over $Q$. Let $n,m$ be positive integers. To any element $T_{n+m}\in Q_{n+m}$ we can associate a linear map $L_{nm}(T_{n+m}):P_n(X)\to P_m(X)$ in the following way:
$$L_{nm}\left(\begin{matrix}|\ |\ |\cr T_{n+m}\cr |\ |\ |\end{matrix}\right):
\left(\begin{matrix}|\cr a_n\cr |\end{matrix}\right)
\to \left(\begin{matrix}
\hskip 1.5mm |\hskip 3.0mm |\hskip 3.0mm \cap\cr
\ \ T_{n+m}\hskip 0.0mm  |\cr
\hskip 1.9mm |\hskip 1.2mm |\hskip 3.2mm |\hskip2.2mm |\cr
a_n|\hskip 3.2mm |\hskip 2.2mm |\cr
\hskip 2.1mm\cup \hskip3.5mm |\hskip 2.2mm |
\end{matrix}\right)$$

That is, we consider the planar $(n,n+m,m)$-tangle having an small input $n$-box, a big input $n+m$-box and an output $m$-box, with strings as on the picture of the right. This defines a certain multilinear map, as follows:
$$P_n(X)\otimes P_{n+m}(X)\to P_m(X)$$

Now let us put the element $T_{n+m}$ in the big input box. We obtain in this way a certain linear map $P_n(X)\to P_m(X)$, that we call $L_{nm}$. Now let us set:
$$Q_{nm}=\left\{ L_{nm}(T_{n+m}):P_n(X)\to P_m(X)\Big| T_{n+m}\in Q_{n+m}\right\}$$

These spaces form a Tannakian category, and so by \cite{wo2} we obtain a Woronowicz algebra $(A,u)$, such that the following equalities hold, for any $m,n$:
$$Hom(u^{\otimes m},u^{\otimes n})=Q_{mn}$$

We prove that $u$ is a magic unitary. We have $Hom(1,u^{\otimes 2})=Q_{02}=Q_2$, so the unit of $Q_2$ must be a fixed vector of $u^{\otimes 2}$. But $u^{\otimes 2}$ acts on the unit of $Q_2$ as follows:
$$u^{\otimes 2}(1)=\sum_{kl}\begin{pmatrix}k&k\\ l&l\end{pmatrix}\otimes (uu^t)_{kl}$$

From $u^{\otimes 2}(1)=1\otimes 1$ ve get that $uu^t$ is the identity matrix, and together with the unitarity of $u$, this gives $u^t=u^*=u^{-1}$. Consider now the Jones projection $E_1\in Q_3$. The linear map $M=L_{21}(E_1)$ is the multiplication $\delta_i\otimes\delta_j\to\delta_{ij}\delta_i$, and we have:
$$(M\otimes id)u^{\otimes 2}\left(\begin{pmatrix}i&i\\ j&j\end{pmatrix}\otimes 1\right)
=\sum_{k}\begin{pmatrix}k\\ k\end{pmatrix}\delta_k\otimes u_{ki}u_{kj}$$
$$u(M\otimes id)\left(\begin{pmatrix}i&i\\ j&j\end{pmatrix}\otimes 1\right)
=\sum_k\begin{pmatrix}k\\ k\end{pmatrix}\delta_k\otimes\delta_{ij}u_{ki}$$

Thus $u_{ki}u_{kj}=\delta_{ij}u_{ki}$ for any $i,j,k$, and we deduce from this that $u$ is a magic unitary. Now if $P$ is the planar algebra associated to $u$, we have $Hom(1,v^{\otimes n})=P_n=Q_n$, as desired. As for the uniqueness, this is clear from the Peter-Weyl theory from \cite{wo1}.
\end{proof}

The above results, following old papers from the early 00s, subsequent to \cite{ba1}, regarding the subgroups $G\subset S_N^+$, have several generalizations, to the subgroups $G\subset O_N^+$ and $G\subset U_N^+$, as well as subfactor versions, going beyond the purely combinatorial level. For the modern story, we refer here to Tarrago-Wahl \cite{twa} and related papers.

\section*{14e. Exercises} 

Things have been quite complicated in this chapter, and as a main exercise on all this, focusing on topics which were beyond our scope here, we have:

\begin{exercise}
Look up the theorem stating that any planar algebra produces a subfactor, and write down a brief account of what you learned.
\end{exercise}

As already mentioned in the above, there are several theorems here, which are all non-trivial. And there is a big open question too, concerning hyperfiniteness.

\chapter{Commuting squares}

\section*{15a. Commuting squares}

In this chapter and in the next one we discuss a number of more specialized aspects of subfactor theory, making the link with several advanced topics, such as quantum groups, noncommutative geometry, free probability, and more. We will mainly insist on the connections with quantum groups, and with the material from chapters 7-8.

\bigskip

A first question, to be discussed in the present chapter, is the explicit construction of subfactors by using some suitable combinatorial data, encoded in a structure called ``commuting square''. Let us start with the following definition:

\index{commuting square}
\index{conditional expectation}

\begin{definition}
A commuting square in the sense of subfactor theory is a commuting diagram of finite dimensional algebras with traces, as follows,
$$\xymatrix@R=40pt@C40pt{
C_{01}\ar[r]&C_{11}\\
C_{00}\ar[u]\ar[r]&C_{10}\ar[u]}$$
having the property that the conditional expectations $C_{11}\to C_{01}$ and $C_{11}\to C_{10}$ commute, and their product is the conditional expectation $C_{11}\to C_{00}$.
\end{definition}

This notion is in fact something that we already talked about, in chapter 14, when discussing the classification of the finite depth subfactors, following the work of Ocneanu \cite{oc1}, \cite{oc2} and Popa \cite{po2}, \cite{po3}. To be more precise, it is possible to prove that any finite depth subfactor of $R$ appears from a commuting square, and vice versa. And as a well-known consequence of this, the subfactors of $R$ having index $<4$, which are all of finite depth, can be shown to be classified by ADE diagrams. But more on this later.

\bigskip

Getting back now to Definition 15.1 as it is, something quite simple, not obviously subfactor related, the idea is that there are many examples of such commuting squares, always coming from subtle combinatorial data. As an illustration for this principle, we have for instance commuting squares associated to the complex Hadamard matrices, that we met in chapter 11, in the maximal abelian subalgebra (MASA) context. In order to discuss this, let us recall from there that, following Popa \cite{po1}, we have:

\begin{theorem}
Up to a conjugation by a unitary, the pairs of orthogonal MASA in the simplest factor, namely $M_N(\mathbb C)$, are as follows,
$$A=\Delta\quad,\quad
B=H\Delta H^*$$
with $\Delta\subset M_N(\mathbb C)$ being the diagonal matrices, and with $H\in M_N(\mathbb C)$ being Hadamard, in the sense that $|H_{ij}|=1$ for any $i,j$, and the rows of $H$ are pairwise orthogonal.
\end{theorem}

\begin{proof}
Any maximal abelian subalgebra in $M_N(\mathbb C)$ being conjugated to $\Delta$, we can assume, up to conjugation by a unitary, that we have, with $U\in U_N$:
$$A=\Delta\quad,\quad 
B=U\Delta U^*$$ 

But a straightforward computation, explained in chapter 11, shows that the orthogonality condition reformulates as $|U_{ij}|=1/\sqrt{N}$, which gives the result.
\end{proof}

As explained in chapter 11, while being something quite trivial, this result remains a statement which is fundamental, surprising, and very interesting, making the link between the general theory of von Neumann algebras, usually associated to rather lugubrious functional analysis computations, and the complex Hadamard matrices, which are a totally opposite beast, belonging to a wild area of linear algebra and combinatorics. As an illustration here, check the following matrix out, with $w=e^{2\pi i/N}$:
$$F_N=\begin{pmatrix}
1&1&1&\ldots&1\\
1&w&w^2&\ldots&w^{N-1}\\
1&w^2&w^4&\ldots&w^{2(N-1)}\\
\vdots&\vdots&\vdots&&\vdots\\
1&w^{N-1}&w^{2(N-1)}&\ldots&w^{(N-1)^2}
\end{pmatrix}$$

This matrix, which is obviously a very beautiful one, hope you agree with me, is called Fourier matrix, and is the most basic example of a complex Hadamard matrix. As explained in chapter 11, this is the matrix of the Fourier transform over the cyclic group $\mathbb Z_N$, and by taking tensor products of such matrices, we obtain the matrices of the Fourier transforms over arbitrary finite abelian groups $G=\mathbb Z_{N_1}\times\ldots\times\mathbb Z_{N_k}$:
$$F_G=F_{N_1}\otimes\ldots\otimes F_{N_k}$$

But the story does not stop here, with basic discrete Fourier analysis. The complex Hadamard matrices, which can be thought of as being ``generalized Fourier matrices'', can be far wilder than that. And among others, above everything, we have:

\begin{conjecture}[Hadamard Conjecture]
There is a real Hadamard matrix
$$H\in M_N(\pm1)$$
for any $N\in 4\mathbb N$.
\end{conjecture}

Here the condition at the end comes from the fact that, assuming $N\geq3$, the orthogonality conditions between the first 3 rows give $N\in 4\mathbb N$. Observe that the Fourier matrices solve this conjecture only at values $N=2^k$, by tensoring $F_2\in M_2(\pm1)$ with itself. For anything else, $N=12,20,24,28,36,40,44,48,52,\ldots\,$, all sorts of clever constructions are needed, whose complexity grows with $N$, and with open questions at $N>666$.

\bigskip

And the conjecture is more than 100 years old, seemingly undoable. Which puts us in a quite delicate situation with our general von Neumann algebra philosophy:

\bigskip

(1) Generally speaking, classical mathematics looks simpler than quantum mathematics, because you start learning one in high school, and the other one in graduate school. And exactly the same goes with classical mechanics vs quantum mechanics.

\bigskip

(2) At a more advanced level, however, classical mathematics turns to be something extremely complicated, wild and unpredictable, with all sorts of notorious no-go areas, such as the Riemann Hypothesis, the Jacobian Conjecture, and so on.

\bigskip

(3) Also at the more advanced level, quantum mathematics, like von Neumann algebras, while certainly difficult, looks plainly doable. Open problems always end up being solved, and you can always dismiss the few no-go areas as being ``uninteresting''.

\bigskip

(4) And so, we have here evidence that quantum mathematics, while being something complicated of course, is probably simpler than classical mathematics. Again, things difficult, but peaceful horizons, with no black holes like the Riemann Hypothesis.

\bigskip

(5) Which agrees with what happens in physics too, where advanced classical mechanics is the hell on Earth, as opposed to quantum mechanics, where the landscape is rather relaxed, with beautiful results promised to everyone willing to give a serious try.

\bigskip

And so, what to do with these Hadamard matrices, which come via Theorem 15.2 to perturb our philosophy. All of the sudden, our von Neumann algebra theory, or even foundations, have a hole in them. Job for us to find a way of dealing with these beasts in a conceptual way, and then either solving Conjecture 15.3, or dismissing it as being ``uninteresting''. In what regards the first task, subfactors come to the rescue, via:

\begin{theorem}
Given an Hadamard matrix $H\in M_N(\mathbb C)$, the diagram formed by the associated pair of orthogonal maximal abelian subalgebras of $M_N(\mathbb C)$,
$$\xymatrix@R=35pt@C35pt{
\Delta\ar[r]&M_N(\mathbb C)\\
\mathbb C\ar[u]\ar[r]&H\Delta H^*\ar[u] }$$ 
where $\Delta\subset M_N(\mathbb C)$ are the diagonal matrices, is a commuting square.
\end{theorem}

\begin{proof}
The expectation $E_\Delta:M_N(\mathbb C)\to\Delta$ is the operation $M\to M_\Delta$ which consists in keeping the diagonal, and erasing the rest. Consider now the other expectation:
$$E_{H\Delta H^*}:M_N(\mathbb C)\to H\Delta H^*$$

It is better to identify this with the following expectation, with $U=H/\sqrt{N}$: 
$$E_{U\Delta U^*}:M_N(\mathbb C)\to U\Delta U^*$$

This must be given by a formula of type $M\to UX_\Delta U^*$, with $X$ satisfying:
$$<M,UDU^*>=<UX_\Delta U^*,UDU^*>\quad,\quad\forall D\in\Delta$$

The scalar products being given by $<a,b>=tr(ab^*)$, this condition reads:
$$tr(MUD^*U^*)=tr(X_\Delta D^*)\quad,\quad\forall D\in\Delta$$

Thus $X=U^*MU$, and the formulae of our two expectations are as follows:
\begin{eqnarray*}
E_\Delta(M)&=&M_\Delta\\
E_{U\Delta U^*}(M)&=&U(U^*MU)_\Delta U^*
\end{eqnarray*}

With these formulae in hand, we have the following computation:
\begin{eqnarray*}
(E_\Delta E_{U\Delta U^*}M)_{ij}
&=&\delta_{ij}(U(U^*MU)_\Delta U^*)_{ii}\\
&=&\delta_{ij}\sum_kU_{ik}(U^*MU)_{kk}\bar{U}_{ik}\\
&=&\delta_{ij}\sum_k\frac{1}{N}\cdot(U^*MU)_{kk}\\
&=&\delta_{ij}tr(U^*MU)\\
&=&\delta_{ij}tr(M)\\
&=&(E_\mathbb CM)_{ij}
\end{eqnarray*}

As for the other composition, the computation here is similar, as follows:
\begin{eqnarray*}
(E_{U\Delta U^*}E_\Delta M)_{ij}
&=&(U(U^*M_\Delta U)_\Delta U^*)_{ij}\\
&=&\sum_kU_{ik}(U^*M_\Delta U)_{kk}\bar{U}_{jk}\\
&=&\sum_{kl}U_{ik}\bar{U}_{lk}M_{ll}U_{lk}\bar{U}_{jk}\\
&=&\frac{1}{N}\sum_{kl}U_{ik}M_{ll}\bar{U}_{jk}\\
&=&\delta_{ij}tr(M)\\
&=&(E_\mathbb CM)_{ij}
\end{eqnarray*}

Thus, we have indeed a commuting square, as claimed.
\end{proof}

To summarize our discussion so far, we had a big scare coming from Popa's Theorem 15.2, but Theorem 15.4, also due to Popa \cite{po1}, puts our von Neumann algebra theory back on tracks. We are doing things which are certainly difficult, but somehow ``trivial'', meaning never undoable in the long run, and that feared Hadamard matrices are simply particular cases of commuting squares. And so, further studying commuting squares will tell us what's interesting and what's not, regarding these matrices, and so on.

\bigskip

Getting back now to Definition 15.1 as it is, there are many other explicit examples of commuting squares, all coming from subtle combinatorial data, and more on this later. So, leaving aside now examples, let us explain the connection with subfactors. For this purpose, consider an arbitrary commuting square, as in Definition 15.1:
$$\xymatrix@R=40pt@C40pt{
C_{01}\ar[r]&C_{11}\\
C_{00}\ar[u]\ar[r]&C_{10}\ar[u]}$$

The point is that, under some suitable extra mild assumptions, any such square $C$ produces a subfactor of the hyperfinite ${\rm II}_1$ factor $R$. Indeed, by performing the basic construction, in finite dimensions, we obtain a whole array, as follows:
$$\xymatrix@R=35pt@C35pt{
A_0&A_1&A_2&\\
C_{02}\ar[r]\ar@.[u]&C_{12}\ar[r]\ar@.[u]&C_{22}\ar@.[r]\ar@.[u]&B_2\\
C_{01}\ar[r]\ar[u]&C_{11}\ar[r]\ar[u]&C_{21}\ar@.[r]\ar[u]&B_1\\
C_{00}\ar[u]\ar[r]&C_{10}\ar[u]\ar[r]&C_{20}\ar[u]\ar@.[r]&B_0}$$

To be more precise, by performing the basic construction in both possible directions, namely to the right and upwards, we obtain a whole array of finite dimensional algebras with traces, that we can denote $(C_{ij})_{i,j\geq0}$, as above. Once this done, we can further consider the von Neumann algebras obtained in the limit, via GNS construction, on each vertical and horizontal line, and denote them $A_i,B_j$, as above. 

\bigskip

With this convention, we have the following result, due to Ocneanu \cite{oc1}, \cite{oc2}:

\index{Ocneanu compactness}
\index{R}

\begin{theorem}
In the context of the above diagram, the limiting von Neumann algebras $A_i,B_j$ are all isomorphic to the hyperfinite ${\rm II}_1$ factor $R$, and:
\begin{enumerate}
\item $A_0\subset A_1$ is a subfactor, and $\{A_i\}$ is the Jones tower for it.

\item The corresponding planar algebra is given by $A_0'\cap A_k=C_{01}'\cap C_{k0}$.

\item A similar result holds for the ``horizontal'' subfactor $B_0\subset B_1$.
\end{enumerate}
\end{theorem}

\begin{proof}
This is something very standard, with the factoriality of the limiting von Neumann algebras $A_i,B_j$ coming as a consequence of the general commutant computation in (2), which is independent from it, with the hyperfiniteness of the same $A_i,B_j$ algebras being clear by definition, and with the idea for the rest being as follows:

\medskip

(1) This is somewhat clear from definitions, or rather from a quick verification of the basic construction axioms, as formulated in chapter 13, because the tower of algebras $\{A_i\}$ appears by definition as the $j\to\infty$ limit of the towers of algebras $\{C_{ij}\}$, which are all Jones towers. Thus the limiting tower $\{A_i\}$ is also a Jones tower.

\medskip

(2) This is the non-trivial result, called Ocneanu compactness theorem, and whose proof is by doing some linear algebra. To be more precise, in one sense the result is clear, because by definition of the algebras $\{A_i\}$, we have inclusions as follows:
$$A_0'\cap A_k\supset C_{01}'\cap C_{k0}$$

In the other sense things are more tricky, mixing standard linear algebra with some functional analysis too, and we refer here to Ocneanu's lecture notes \cite{oc1}, \cite{oc2}.

\medskip

(3) This follows from (1,2), by transposing the whole diagram. Indeed, given a commuting square as in Definition 15.1, its transpose is a commuting square as well:
$$\xymatrix@R=40pt@C40pt{
C_{10}\ar[r]&C_{11}\\
C_{00}\ar[u]\ar[r]&C_{01}\ar[u]}$$

Thus we can apply (1,2) above to this commuting square, and we obtain in this way Jones tower and planar algebra results for the ``horizontal'' subfactor $B_0\subset B_1$.
\end{proof}

In relation with the examples of commuting squares that we have so far, namely those coming from the Hadamard matrices, from Theorem 15.4, we can upgrade what we have so far into something more conceptual, due to Jones \cite{jo3}, as follows:

\index{Hadamard matrix}
\index{commuting square}

\begin{theorem}
Given a complex Hadamard matrix $H\in M_N(\mathbb C)$, the diagram formed by the associated pair of orthogonal maximal abelian subalgebras, namely
$$\xymatrix@R=35pt@C35pt{
\Delta\ar[r]&M_N(\mathbb C)\\
\mathbb C\ar[u]\ar[r]&H\Delta H^*\ar[u]}$$ 
is a commuting square in the sense of subfactor theory, and the associated planar algebra $P=(P_k)$ is given by the following formula, in terms of $H$ itself,
$$T\in P_k\iff T^\circ G^2=G^{k+2}T^\circ$$
where the objects on the right are constructed as follows:
\begin{enumerate}
\item $T^\circ=id\otimes T\otimes id$.

\item $G_{ia}^{jb}=\sum_kH_{ik}\bar{H}_{jk}\bar{H}_{ak}H_{bk}$.

\item $G^k_{i_1\ldots i_k,j_1\ldots j_k}=G_{i_ki_{k-1}}^{j_kj_{k-1}}\ldots G_{i_2i_1}^{j_2j_1}$.
\end{enumerate}
\end{theorem}

\begin{proof}
We have several assertions here, the idea being as follows:

\medskip

(1) The fact that we have indeed a commuting square is something quite elementary, that we already know, from Theorem 15.4.

\medskip

(2) The computation of the associated planar algebra, directly in terms of $H$, is something which is definitely possible, thanks to the formula in Theorem 15.5 (2). 

\medskip

(3) As for the precise formula of the planar algebra, which emerges by doing the computation, we will be back to it, with full details, later on.

\medskip

(4) The point indeed is that we want to first develop some better methods in dealing with the Hadamard matrices, and leave the computation of $P$ for later.
\end{proof}

Summarizing, we have so far an interesting combinatorial notion, that of a commuting square, and a method of producing subfactors and planar algebras out of it. We will further explore all the possibilities that this opens up, in what follows:

\bigskip

(1) In the remainder of this chapter we will keep working on the Hadamard matrix problem, following \cite{ba1} and subsequent papers. This might look of course a bit like mania, focusing just like that on a single class of commuting squares, but we are strongly motivated by all that has being said after Theorem 15.2 and Conjecture 15.3, with this being a matter of life and death to us. And don't worry, we will learn in this way useful techniques, that will apply to other commuting squares too. And also, following Jones \cite{jo2}, \cite{jo3} and others, all this is potentially related to some interesting physics too.

\bigskip

(2) And in chapter 16 below we will go back to general commuting squares, and to their more traditional usage, for classification problems for small index subfactors. 

\section*{15b. Matrix models}

Our objective now is to clarify the planar algebra computation for the commuting squares coming from Hadamard matrices, from Theorem 15.6. Our claim is that all this is related, and in a beautiful way, to the quantum permutation groups that we met in chapters 7-8, and at the end of chapter 14 as well. In order to discuss this, and to present as well some generalizations, we will need some preliminaries on the quantum permutation groups, and their matrix models. Let us recall from chapter 11 that we have:

\index{matrix model}

\begin{definition}
A matrix model for a Woronowicz algebra $A=C(G)$ is a morphism of $C^*$-algebras of the following type,
$$\pi:C(G)\to M_K(C(T))$$
with $T$ being a compact space, and $K\geq1$ being an integer.
\end{definition}

As explained in chapter 11, assuming that $\pi$ is faithful leads to the conclusion that $C(G)$ must be a type I algebra, and so that $G$ must be coamenable, and with this being something quite restrictive, excluding for instance all the free quantum groups.

\bigskip

The solution to this problem comes from a weaker notion of faithfulness, called ``inner faithfulness'', which still allows to recover the combinatorics of $G$ from the combinatorics of the model, but does not potentially exclude any quantum group. The theory here, briefly explained in chapter 11 too, starts with the following definition:

\index{Hopf image}
\index{inner faithfulness}

\begin{definition}
Let $\pi:C(G)\to M_K(C(T))$ be a matrix model. 
\begin{enumerate}
\item The Hopf image of $\pi$ is the smallest quotient Hopf $C^*$-algebra $C(G)\to C(H)$ producing a factorization of type $\pi:C(G)\to C(H)\to M_K(C(T))$.

\item When the inclusion $H\subset G$ is an isomorphism, i.e. when there is no non-trivial factorization as above, we say that $\pi$ is inner faithful.
\end{enumerate}
\end{definition}

As explained in \cite{bbi}, the existence and uniqueness of the Hopf image come by dividing $C(G)$ by a suitable ideal, although we will come in a moment with an explicit Tannakian construction as well, also from \cite{bbi}. As a basic illustration for these notions, we have two main examples, which are somehow dual to each other, as follows:

\bigskip

(1) In the case where $G=\widehat{\Gamma}$ is a group dual, $\pi$ must come from a group representation $\rho:\Gamma\to C(T,U_K)$. We conclude that in this case, the minimal factorization constructed in Definition 15.8 is simply the one obtained by taking the image:
$$\rho:\Gamma\to\Lambda\subset C(T,U_K)$$

Thus $\pi$ is inner faithful when our group satisfies $\Gamma\subset C(T,U_K)$. And we can see here that $\pi$, while not being faithful, clearly reminds all of $\Gamma$, and so of $G=\widehat{\Gamma}$ too.

\bigskip

(2) As a second illustration, given a compact group $G$, and elements $g_1,\ldots,g_K\in G$, we have a representation $\pi:C(G)\to\mathbb C^K$, given by $f\to(f(g_1),\ldots,f(g_K))$. The minimal factorization of $\pi$ is then via $C(H)$, with $H\subset G$ being the following subgroup:
$$H=\overline{<g_1,\ldots,g_K>}$$

Thus $\pi$ is inner faithful precisely when $G=\overline{<g_1,\ldots,g_K>}$. Again, we can see here that $\pi$, while not being faithful, clearly reminds all of $G$, and so of $\Gamma=\widehat{G}$ too.

\bigskip

Summarizing, our notion of inner faithfulness does the job, reminding the quantum groups $G$ and $\Gamma=\widehat{G}$, and not excluding anything on functional analysis grounds. Which brings us into the question of recapturing the algebraic and analytic properties of $G$ and $\Gamma=\widehat{G}$ out the combinatorics of the model. Regarding algebra, we have here:

\index{Tannakian duality}

\begin{theorem}
Assuming $G\subset U_N^+$, with fundamental corepresentation $u=(u_{ij})$, the Hopf image of $\pi:C(G)\to M_K(C(T))$ comes from the following Tannakian category,
$$C_{kl}=Hom(U^{\otimes k},U^{\otimes l})$$
where $U_{ij}=\pi(u_{ij})$, and where the spaces on the right are taken in a formal sense.
\end{theorem}

\begin{proof}
This is something that we know from chapter 11, but we will recall the proof here. Since the morphisms increase the intertwining spaces, when defined either in a representation theory sense, or just formally, we have inclusions as follows:
$$Hom(u^{\otimes k},u^{\otimes l})\subset Hom(U^{\otimes k},U^{\otimes l})$$

More generally, we have such inclusions when replacing $(G,u)$ with any pair producing a factorization of $\pi$. Thus, by Tannakian duality, the Hopf image must be given by the fact that the intertwining spaces must be the biggest, subject to the above inclusions. On the other hand, since $u$ is biunitary, so is $U$, and it follows that the spaces on the right form a Tannakian category. Thus, we have a quantum group $(H,v)$ given by:
$$Hom(v^{\otimes k},v^{\otimes l})=Hom(U^{\otimes k},U^{\otimes l})$$

By the above discussion, $C(H)$ follows to be the Hopf image of $\pi$, as claimed.
\end{proof}

In what regards now analysis, the result here is as follows:

\begin{theorem}
Given an inner faithful model $\pi:C(G)\to M_K(C(T))$, we have
$$\int_G=\lim_{k\to\infty}\frac{1}{k}\sum_{r=1}^k\int_G^r$$
where $\int_G^r=(\varphi\circ\pi)^{*r}$, with $\varphi=tr\otimes\int_T$ being the random matrix trace.
\end{theorem}

\begin{proof}
Again, this is something that we know from chapter 11. If we denote by $\int_G'$ the limit in the statement, we must prove that this limit converges, and that we have:
$$\int_G'=\int_G$$

It is enough to check this on the coefficients of corepresentations, and if we let $v=u^{\otimes k}$ be one of the Peter-Weyl corepresentations, we must prove that we have:
$$\left(id\otimes\int_G'\right)v=\left(id\otimes\int_G\right)v$$

We know from chapter 7 that the matrix on the right is the orthogonal projection onto $Fix(v)$. Regarding now the matrix on the left, this is the orthogonal projection onto the $1$-eigenspace of $(id\otimes\varphi\pi)v$. Now observe that, if we set $V_{ij}=\pi(v_{ij})$, we have:
$$(id\otimes\varphi\pi)v=(id\otimes\varphi)V$$

Thus, as in chapter 7, we conclude that the $1$-eigenspace that we are interested in equals $Fix(V)$. But, according to Theorem 15.9, we have:
$$Fix(V)=Fix(v)$$

Thus, we have proved that we have $\int_G'=\int_G$, as desired.
\end{proof}

\section*{15c. Hadamard models}

With this theory in hand, let us go back now to our von Neumann algebra and subfactor questions. In relation with the complex Hadamard matrices, the connection with the quantum permutations is immediate, coming from the following observation:

\index{Hadamard matrix}

\begin{proposition}
If $H\in M_N(\mathbb C)$ is Hadamard, the rank one projections 
$$P_{ij}=Proj\left(\frac{H_i}{H_j}\right)$$
where $H_1,\ldots,H_N\in\mathbb T^N$ are the rows of $H$, form a magic unitary.
\end{proposition}

\begin{proof}
This is clear, the verification for the rows being as follows:
\begin{eqnarray*}
\left<\frac{H_i}{H_j},\frac{H_i}{H_k}\right>
&=&\sum_l\frac{H_{il}}{H_{jl}}\cdot\frac{H_{kl}}{H_{il}}\\
&=&\sum_l\frac{H_{kl}}{H_{jl}}\\
&=&N\delta_{jk}
\end{eqnarray*}

As for the verification for the columns, this is similar, as follows:
\begin{eqnarray*}
\left<\frac{H_i}{H_j},\frac{H_k}{H_j}\right>
&=&\sum_l\frac{H_{il}}{H_{jl}}\cdot\frac{H_{jl}}{H_{kl}}\\
&=&\sum_l\frac{H_{il}}{H_{kl}}\\
&=&N\delta_{ik}
\end{eqnarray*}

Thus, we have indeed a magic unitary, as claimed. 
\end{proof}

We are led in this way into the following notion:

\begin{definition}
To any Hadamard matrix $H\in M_N(\mathbb C)$ we associate the quantum permutation group $G\subset S_N^+$ given by the following Hopf image factorization,
$$\xymatrix{C(S_N^+)\ar[rr]^{\pi}\ar[rd]&&M_N(\mathbb C)\\&C(G)\ar[ur]&}$$
where $\pi(u_{ij})=Proj(H_i/H_j)$, with $H_1,\ldots,H_N\in\mathbb T^N$ being the rows of $H$.
\end{definition}

Our claim now is that this construction $H\to G$ is something really useful, with the quantum group $G$ encoding the combinatorics of $H$. To be more precise, the idea will be that ``$H$ can be thought of as being a kind of Fourier matrix for $G$''. As an illustration for this principle, we first have the following result:

\index{Fourier matrix}

\begin{theorem}
The construction $H\to G$ has the following properties:
\begin{enumerate}
\item For a Fourier matrix $H=F_G$ we obtain the group $G$ itself, acting on itself.

\item For $H\not\in\{F_G\}$, the quantum group $G$ is not classical, nor a group dual.

\item For a tensor product $H=H'\otimes H''$ we obtain a product, $G=G'\times G''$.
\end{enumerate}
\end{theorem}

\begin{proof}
All this material is standard, and elementary, as follows:

\medskip

(1) Let us first discuss the cyclic group case, $H=F_N$. Here the rows of $H$ are given by $H_i=\rho^i$, where $\rho=(1,w,w^2,\ldots,w^{N-1})$. Thus, we have the following formula:
$$\frac{H_i}{H_j}=\rho^{i-j}$$

It follows that the corresponding rank 1 projections $P_{ij}=Proj(H_i/H_j)$ form a circulant matrix, all whose entries commute. Since the entries commute, the corresponding quantum group must satisfy $G\subset S_N$. Now by taking into account the circulant property of $P=(P_{ij})$ as well, we are led to the conclusion that we have $G=\mathbb Z_N$, as claimed.

\medskip

In the general case now, where $H=F_G$, with $G$ being an arbitrary finite abelian group, the result can be proved either by extending the above proof, of by decomposing $G=\mathbb Z_{N_1}\times\ldots\times\mathbb Z_{N_k}$ and using (3) below, whose proof is independent from the rest.

\medskip

(2) This is something more tricky, needing some general study of the representations whose Hopf images are commutative, or cocommutative. For details here, along with a number of supplementary facts on the construction $H\to G$, we refer to \cite{bbi}.

\medskip

(3) Assume that we have a tensor product $H=H'\otimes H''$, and let $G,G',G''$ be the associated quantum permutation groups. We have then a diagram as follows:
$$\xymatrix@R=45pt@C25pt{
C(S_{N'}^+)\otimes C(S_{N''}^+)\ar[r]&C(G')\otimes C(G'')\ar[r]&M_{N'}(\mathbb C)\otimes M_{N''}(\mathbb C)\ar[d]\\
C(S_N^+)\ar[u]\ar[r]&C(G)\ar[r]&M_N(\mathbb C)
}$$

Here all the maps are the canonical ones, with those on the left and on the right coming from $N=N'N''$. At the level of standard generators, the diagram is as follows:
$$\xymatrix@R=45pt@C65pt{
u_{ij}'\otimes u_{ab}''\ar[r]&w_{ij}'\otimes w_{ab}''\ar[r]&P_{ij}'\otimes P_{ab}''\ar[d]\\
u_{ia,jb}\ar[u]\ar[r]&w_{ia,jb}\ar[r]&P_{ia,jb}
}$$

Now observe that this diagram commutes. We conclude that the representation associated to $H$ factorizes indeed through $C(G')\otimes C(G'')$, and this gives the result.
\end{proof}

In order to discuss now the relation with the commuting squares and the subfactors, we can use Theorem 15.9, and we are led to the following result:

\begin{theorem}
The Tannakian category of the quantum group $G\subset S_N^+$ associated to a complex Hadamard matrix $H\in M_N(\mathbb C)$ is given by
$$T\in Hom(u^{\otimes k},u^{\otimes l})\iff T^\circ G^{k+2}=G^{l+2}T^\circ$$
where the objects on the right are constructed as follows:
\begin{enumerate}
\item $T^\circ=id\otimes T\otimes id$.

\item $G_{ia}^{jb}=\sum_kH_{ik}\bar{H}_{jk}\bar{H}_{ak}H_{bk}$.

\item $G^k_{i_1\ldots i_k,j_1\ldots j_k}=G_{i_ki_{k-1}}^{j_kj_{k-1}}\ldots G_{i_2i_1}^{j_2j_1}$.
\end{enumerate}
\end{theorem}

\begin{proof}
According to Theorem 15.9, and with the notations there, we have the following formula for the Tannakian category that we are interested in:
$$Hom(u^{\otimes k},u^{\otimes l})=Hom(U^{\otimes k},U^{\otimes l})$$

The vector space on the right, that we will compute now, consists by definition of the complex $N^l\times N^k$ matrices $T$ satisfying the following relation:
$$TU^{\otimes k}=U^{\otimes l}T$$ 

If we denote this equality by $L=R$, the left term $L$ is given by:
\begin{eqnarray*}
L_{ij}
&=&(TU^{\otimes k})_{ij}\\
&=&\sum_aT_{ia}U^{\otimes k}_{aj}\\
&=&\sum_aT_{ia}U_{a_1j_1}\ldots U_{a_kj_k}
\end{eqnarray*}

As for the right term $R$, this is given by a similar formula, as follows:
\begin{eqnarray*}
R_{ij}
&=&(U^{\otimes l}T)_{ij}\\
&=&\sum_bU^{\otimes l}_{ib}T_{bj}\\
&=&\sum_bU_{i_1b_1}\ldots U_{i_lb_l}T_{bj}
\end{eqnarray*}

Consider now the vectors $\xi_{ij}=H_i/H_j$. Since these vectors span the ambient Hilbert space, the equality $L=R$ is equivalent to the following equality:
$$<L_{ij}\xi_{pq},\xi_{rs}>=<R_{ij}\xi_{pq},\xi_{rs}>$$

We use now the following well-known formula, expressing a product of rank one projections $P_1,\ldots,P_k$ in terms of the corresponding image vectors $\xi_1,\ldots,\xi_k$:
$$<P_1\ldots P_kx,y>=<x,\xi_k><\xi_k,\xi_{k-1}>\ldots\ldots<\xi_2,\xi_1><\xi_1,y>$$

This gives the following formula for the left term $L$:
\begin{eqnarray*}
<L_{ij}\xi_{pq},\xi_{rs}>
&=&\sum_aT_{ia}<P_{a_1j_1}\ldots P_{a_kj_k}\xi_{pq},\xi_{rs}>\\
&=&\sum_aT_{ia}<\xi_{pq},\xi_{a_kj_k}>\ldots<\xi_{a_1j_1},\xi_{rs}>\\
&=&\sum_aT_{ia}G_{pa_k}^{qj_k}G_{a_ka_{k-1}}^{j_kj_{k-1}}\ldots G_{a_2a_1}^{j_2j_1}G_{a_1r}^{j_1s}\\
&=&\sum_aT_{ia}G^{k+2}_{rap,sjq}\\
&=&(T^\circ G^{k+2})_{rip,sjq}
\end{eqnarray*}

As for the right term $R$, this is given by the following formula:
\begin{eqnarray*}
<R_{ij}\xi_{pq},\xi_{rs}>
&=&\sum_b<P_{i_1b_1}\ldots P_{i_lb_l}\xi_{pq},\xi_{rs}>T_{bj}\\
&=&\sum_b<\xi_{pq},\xi_{i_lb_l}>\ldots<\xi_{i_1b_1},\xi_{rs}>T_{bj}\\
&=&\sum_bG_{pi_l}^{qb_l}G_{i_li_{l-1}}^{b_lb_{l-1}}\ldots G_{i_2i_1}^{b_2b_1}G_{i_1r}^{b_1s}T_{bj}\\
&=&\sum_bG^{l+2}_{rip,sbq}T_{bj}\\
&=&(G^{l+2}T^\circ)_{rip,sjq}
\end{eqnarray*}

Thus, we obtain the formula in the statement.
\end{proof}

The point now is that, with $k=0$, we obtain in this way precisely the planar algebra spaces $P_l$ computed by Jones in \cite{jo3}, for the corresponding commuting square, described in Theorem 15.6. Thus, we are led in this way to the following result: 

\index{planar algebra}

\begin{theorem}
Let $H\in M_N(\mathbb C)$ be a complex Hadamard matrix.
\begin{enumerate}
\item The planar algebra associated to $H$ is given by the formula
$$P_k=Fix(u^{\otimes k})$$
where $G\subset S_N^+$ is the associated quantum permutation group.

\item The Poincar\'e series $\sum_k\dim(P_k)z^k$ equals the Stieltjes transform 
$$f(z)=\int_G\frac{1}{1-z\chi}$$ 
of the law of the main character $\chi=\sum_iu_{ii}$.
\end{enumerate}
\end{theorem}

\begin{proof}
This follows as indicated above, by putting together what we have:

\medskip

(1) As already mentioned above, this simply follows by comparing Theorem 15.14 with the subfactor computation in \cite{jo3}, discussed in Theorem 15.6.

\medskip

(2) This follows from (1) and from the Peter-Weyl theory, with the statement itself being a nice and concrete application of our main result, (1) above.
\end{proof}

Summarizing, in connection with the commuting square problematics from the beginning of this chapter, the conclusion is that for the simplest such commuting squares, namely those coming from Hadamard matrices, the combinatorics ultimately comes from quantum permutation groups. This is something nice, and exploring improvements and generalizations of this will be our main purpose, in the remainder of this chapter.

\section*{15d. Fixed points}

We know that the planar algebra associated to an Hadamard matrix $H\in M_N(\mathbb C)$ appears in fact as the planar algebra associated to a certain related quantum permutation group $G\subset S_N^+$. In view of the various results from chapters 13-14, this suggests that the subfactor itself associated to $H$ should appear as a fixed point subfactor associated to $G$. We will prove here that this is indeed the case. To be more precise, following \cite{ba1} and subsequent papers, regarding the subfactor itself, the result here is as follows:

\index{fixed point subfactor}

\begin{theorem}
The subfactor associated to $H\in M_N(\mathbb C)$ is of the form
$$A^G\subset(\mathbb C^N\otimes A)^G$$
with $A=R\rtimes\widehat{G}$, where $G\subset S_N^+$ is the associated quantum permutation group.
\end{theorem}

\begin{proof}
This is something more technical, the idea being that the basic construction procedure for the commuting squares, explained before Theorem 15.5, can be performed in an ``equivariant setting'', for commuting squares having components as follows:
$$D\otimes_GE=(D\otimes(E\rtimes\widehat{G}))^G$$

To be more precise, starting with a commuting square formed by such algebras, we obtain by basic construction a whole array of commuting squares as follows, with $\{D_i\},\{E_i\}$ being by definition Jones towers, and with $D_\infty,E_\infty$ being their inductive limits:
$$
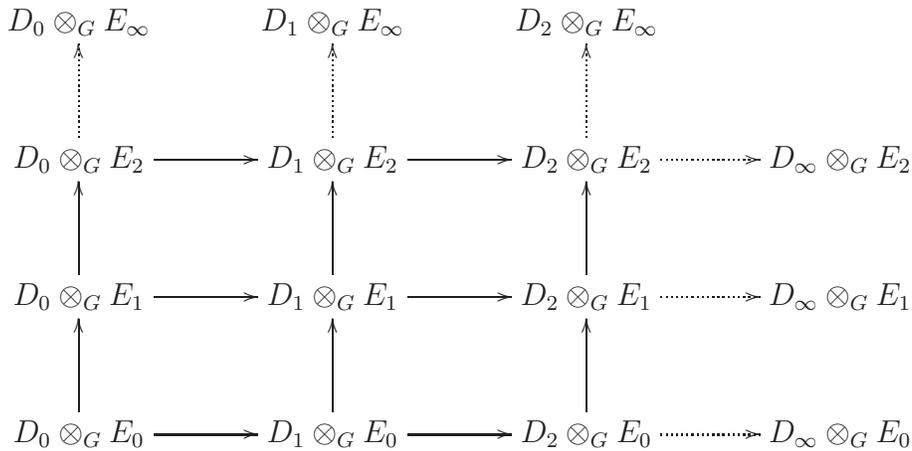
\xymatrix@R=35pt@C35pt{
D_0\otimes_GE_\infty&D_1\otimes_GE_\infty&D_2\otimes_GE_\infty\\
D_0\otimes_GE_2\ar@.[u]\ar[r]&D_1\otimes_GE_2\ar@.[u]\ar[r]&D_2\otimes_GE_2\ar@.[u]\ar@.[r]&D_\infty\otimes_GE_2\\
D_0\otimes_GE_1\ar[u]\ar[r]&D_1\otimes_GE_1\ar[u]\ar[r]&D_2\otimes_GE_1\ar[u]\ar@.[r]&D_\infty\otimes_GE_1\\
D_0\otimes_GE_0\ar[u]\ar[r]&D_1\otimes_GE_0\ar[u]\ar[r]&D_2\otimes_GE_0\ar[u]\ar@.[r]&D_\infty\otimes_GE_0}$$

The point now is that this quantum group picture works in fact for any commuting square having $\mathbb C$ in the lower left corner. In the Hadamard matrix case, that we are interested in here, the corresponding commuting square is as follows:
$$\xymatrix@R=35pt@C35pt{
\mathbb C\otimes_G\mathbb C^N\ar[r]&\mathbb C^N\otimes_G\mathbb C^N\\
\mathbb C\otimes_G\mathbb C\ar[u]\ar[r]&\mathbb C^N\otimes_G\mathbb C\ar[u] }$$ 

Thus, the subfactor obtained by vertical basic construction appears as follows:
$$\mathbb C\otimes_GE_\infty\subset\mathbb C^N\otimes_GE_\infty$$

But this gives the conclusion in the statement, with the ${\rm II}_1$ factor appearing there being by definition $A=E_\infty\rtimes\widehat{G}$, and with the remark that we have $E_\infty\simeq R$.
\end{proof}

All the above was of course quite brief, but we will discuss now all this with more details, directly in a more general setting, covering the Hadamard matrix situation. To be more precise, our claim is that the above fixed point subfactor techniques apply, more generally, to the commuting squares having $\mathbb C$ in the lower left corner:
$$\xymatrix@R=40pt@C40pt{
E\ar[r]&X\\
\mathbb C\ar[u]\ar[r]&D\ar[u]}$$

In order to discuss this, let us go back to the fixed point subfactors, from chapter 13. In what concerns the fixed point factors, we know from there that we have:

\begin{theorem}
Consider a Woronowicz algebra $A=(A,\Delta,S)$, and denote by $A_\sigma$ the Woronowicz algebra $(A,\sigma\Delta ,S)$, where $\sigma$ is the flip. Given coactions
$$\beta:B\to B\otimes A$$
$$\pi:P\to P\otimes A_\sigma$$
with $B$ being finite dimensional, the following linear map, while not being multiplicative in general, is coassociative with respect to the comultiplication $\sigma\Delta$ of $A_\sigma$,
$$\beta\odot\pi:B\otimes P\to B\otimes P\otimes A_\sigma$$
$$b\otimes p\to \pi (p)_{23}((id\otimes S)\beta(b))_{13}$$
and its fixed point space, which is by definition the following linear space,
$$(B\otimes P)^{\beta\odot\pi}=\left\{x\in B\otimes P\Big|(\beta\odot\pi )x=x\otimes 1\right\}$$
is then a von Neumann subalgebra of $B\otimes P$. Moreover, such algebras can be used in order to construct the generalized Wassermann subfactors, $(B_0\otimes P)^G\subset(B_1\otimes P)^G$.
\end{theorem}

\begin{proof}
This is something that we know from chapter 13, and for details, and comments in relation with the non-multiplicativity of $\beta\odot\pi$, we refer to the material there.
\end{proof}

Let $\int_A:A\to\mathbb C$ be the Haar functional, let $l^2(A)$ be its $l^2$-space and let $\widehat{A}\subset B(l^2(A))$ be the dual algebra. If $\alpha:E\to E\otimes\widehat A$ is a coaction of $\widehat{A}$ on a finite von Neumann algebra $E$, the crossed product $E\rtimes_\alpha\widehat{A}$ is the von Neumann subalgebra of $E\otimes B(l^2(A))$ generated by $\alpha(E)$ and by $1\otimes A$. There exists a unique coaction $\widehat{\alpha}$ of $A$ on $E\rtimes_\alpha\widehat{A}$ such that $(E\rtimes_\alpha\widehat{A})^{\widehat{\alpha}}=\alpha(E)$, and such that the copy $1\otimes A$ of $A$ is equivariant. With these conventions, again following \cite{ba1} and subsequent papers, we have the following result:
  
\begin{proposition}
Let $A$ be a Woronowicz algebra. If $\beta:D\to D\otimes A$ is a coaction on a finite dimensional finite von Neumann algebra and $\alpha:E\to E\otimes\widehat{A}_\sigma$ is a coaction on a finite von Neumann algebra then we have the equality
$$(D\otimes(E\rtimes_\alpha\widehat{A}_\sigma))^{\beta\odot\widehat{\alpha}}=\overline{sp}^{\,w}\Big\{\beta(D)_{13}\cdot\alpha(E)_{23}\Big\}$$
as linear subspaces of $D\otimes E\otimes B(l^2(A_\sigma))$. Moreover, the following diagram
$$\begin{matrix}
\alpha(E)_{23}&\subset&(D\otimes(E\rtimes_{\alpha}\widehat{A}_\sigma 
))^{\beta\odot\widehat{\alpha}}\\
\cup&&\cup\\ 
\mathbb C&\subset&\beta(D)_{13}
\end{matrix}$$
is a non-degenerate commuting square of finite von Neumann algebras.
\end{proposition}

\begin{proof}
By definition of the crossed product $E\rtimes_\alpha\widehat{A}_\sigma$, we have the following equalities between subalgebras of $D\otimes E\otimes B(l^2(A_\sigma))$:
\begin{eqnarray*}
D\otimes(E\rtimes_\alpha\widehat{A}_\sigma)
&=&D\otimes(\overline{sp}^{\,w}\{\alpha(E)\cdot (1\otimes A_\sigma)\})\\
&=&\overline{sp}^{\,w}\{(D\otimes A_\sigma )_{13}\cdot\alpha(E)_{23}\}
\end{eqnarray*}

On the other hand, since the coactions on the finite dimensional algebras
are automatically non-degenerate, we have as well the following equality:
$$D\otimes A_\sigma=sp\{(1\otimes A_\sigma)\cdot\beta(D)\}$$

Thus, we have the following equality of algebras:
$$D\otimes(E\rtimes_\alpha A_\sigma)=\overline{sp}^{\,w}\{(1\otimes
1\otimes A_\sigma)\cdot\beta(D)_{13}\cdot\alpha(E)_{23}\}$$

Let us compute now the restriction of the map $\beta\odot\widehat{\alpha}$ to the algebra $1\otimes1\otimes A_\sigma$, to the algebra $\beta(D)_{13}$, and to the algebra $\alpha (E)_{23}$. This can be done as follows:

\medskip

(1) The restriction of $\beta\odot\widehat{\alpha}$ to the algebra $1\otimes1\otimes A_\sigma$ is $1\otimes1\otimes\sigma\Delta$. In particular the map $\beta\odot\widehat{\alpha}$ has no fixed points in this algebra $1\otimes1\otimes A_\sigma$.

\medskip

(2) The algebra $\alpha(E)_{23}$ is by definition fixed by $\beta\odot\widehat{\alpha}$.

\medskip

(3) We prove now that the algebra $\beta(D)_{13}$ is also fixed by $\beta\odot\widehat{\alpha}$. For this purpose, let $\{u_{ij}\}$ be an orthonormal basis of $l^2(A_\sigma)$ consisting of coefficients of irreducible corepresentations of $A_\sigma$. Since we have $\beta(D)\subset D\otimes_{alg}A_\sigma$, for any $b\in D$ we can use the notation $\beta(b)=\sum_{uij}b^u_{ij}\otimes u_{ij}$. From the coassociativity of $\beta$ we obtain:
$$\sum_{uij}\beta(b_{ij}^u)\otimes u_{ij}=\sum_{uijk}b_{ij}^u\otimes u_{kj}\otimes u_{ik}$$

Thus we have $\beta (b_{ik}^u)=\sum_j b_{ij}^u\otimes u_{kj}$ for any
$u,i,k$, and so:
\begin{eqnarray*}
(id\otimes S)\beta(b^u_{ij})
&=&(id\otimes S)\left(\sum_s b_{is}^u\otimes u_{js}\right)\\
&=&\sum_sb_{is}^u\otimes u_{sj}^*
\end{eqnarray*}

Also, we have $\widehat{\alpha}(1\otimes u_{ij})=\sum_k1\otimes u_{ik}\otimes u_{kj}$, and we obtain from this that we have:
\begin{eqnarray*}
(\beta\odot\widehat{\alpha})(\beta (b)_{13})
&=&\sum_{uij}\left(\sum_k1\otimes1\otimes u_{ik}\otimes u_{kj}\right)\left(\sum_sb_{is}^u\otimes1\otimes1\otimes u_{sj}^*\right)\\
&=&\sum_{uijks}b_{is}^u\otimes 1\otimes u_{ik}\otimes u_{kj}u_{sj}^*
\end{eqnarray*}

By summing over $j$ the last term is replaced by $(uu^*)_{ks}=\delta_{k,s}1$. Thus we obtain, as desired, that our algebra consists indeed of fixed points:
\begin{eqnarray*}
(\beta\odot\widehat{\alpha})(\beta(b)_{13})
&=&\sum_{uik}b_{ik}^u\otimes 1\otimes u_{ik}\otimes 1\\
&=&(\beta(b)_{13})\otimes 1
\end{eqnarray*}

In order to finish now, observe that (1,2,3) above show that $(D\otimes(E\rtimes_\alpha\widehat{A}_\sigma))^{\beta\odot\widehat{\alpha}}$, which is the fixed point algebra of $\overline{sp}^{\,w}\{ (1\otimes
1\otimes A_\sigma)\cdot\beta(D)_{13}\cdot\alpha(E)_{23}\}$ under the coaction $\beta\odot\widehat{\alpha}$, is equal to $\overline{sp}^{\,w}\{ \beta  (D)_{13}\cdot\alpha(E)_{23}\}$. This finishes the proof of the first assertion, and proves as well the non-degeneracy of the diagram in the statement. 

\medskip

Finally, observe that the diagram in the statement is the dual of the square on the left in the following diagram, where $P=E\rtimes_{\alpha}\widehat{A}_\sigma$ and $\pi=\widehat{\alpha}$:
$$\begin{matrix}
D&\subset&(D\otimes P)^{\beta\odot\pi}&\subset&D\otimes P\\ 
\cup&&\cup&&\cup\\ 
\mathbb C&\subset&P^\pi&\subset&P
\end{matrix}$$

Since $\pi$ is dual, the square on the right is a non-degenerate commuting square. We also know that the rectangle is a non-degenerate commuting square. Thus if we denote by $E_X:D\otimes P\to D\otimes P$ the conditional expectation onto $X$, for any $X$, then for any $b\in D$ we have $E_{P^\pi}(b)=E_P(b)=E_\mathbb C(b)$, and this proves the commuting square condition.
\end{proof}

Let us denote now by ${\mathcal Alg}$ the category having as objects the finite dimensional $C^*$-algebras and having as arrows the inclusions of $C^*$-algebras which preserve the canonical traces. The above result suggests the following abstract definition:

\begin{definition}
Given objects $(D,\beta)\in A-{\mathcal Alg}$ and $(E,\alpha)\in\widehat{A}_\sigma -{\mathcal Alg}$, we let
$$D\Box_AE=(D\otimes(E\rtimes_\alpha\widehat{A}_\sigma ))^{\beta\odot\widehat{\alpha}}$$
be the object in ${\mathcal Alg}$, constructed as in Proposition 15.18 above.
\end{definition}

If $(D',\beta')\subset(D,\beta )$ is an arrow in $A-{\mathcal Alg}$ and $(E',\alpha')\subset(E,\alpha )$ is an arrow in $\widehat{A}_\sigma-{\mathcal Alg}$, then we have a canonical embedding, as follows:
$$D'\Box_AE'\subset D\Box_AE$$

Now since both $D'\Box_AE'$ and $D\Box_AE$ are endowed with their canonical traces, this inclusion is Markov. Thus, we have constructed a bifunctor, as follows:
$$\Box_A:A-{\mathcal Alg}\times\widehat{A}_\sigma-{\mathcal Alg}\to\mathcal Alg$$

With this convention, we have the following result:

\begin{theorem}
For any two arrows $D_0\subset D_1$ in $A-{\mathcal Alg}$ and $E_0\subset E_1$ in $\widehat{A}_\sigma -{\mathcal Alg}$,
$$\begin{matrix}
D_0\Box_AE_1&\subset&D_1\Box_AE_1\\ 
\cup&&\cup\\ 
D_0\Box_AE_0&\subset&D_1\Box_AE_0
\end{matrix}$$
is a non-degenerate commuting square of finite dimensional von Neumann algebras.
\end{theorem}

\begin{proof}
This can be proved in several steps, as follows:

\medskip

\underline{Step I}. In the simplest case, $D_0=E_0=\mathbb C$, this follows from the above.

\medskip

\underline{Step II}. We prove now the result in the general $E_0=\mathbb C$ case. Indeed, let $E=E_1$, and consider the following diagram:
$$\begin{matrix}
E&\subset&D_0\Box_AE&\subset&D_1\Box_AE\\ 
\cup&&\cup&&\cup\\ 
\mathbb C&\subset&D_0&\subset&D_1
\end{matrix}$$

By the result of Step I, both the square on the left and the rectangle are non-degenerate commuting squares. We want to prove that the square on the right is a non-degenerate commuting square. But the non-degeneracy condition follows from:
$$D_1\Box_AE=sp\{E\cdot D_1\}\subset sp\{D_0\Box_AE\cdot D_1\}$$

Now let $x\in D_0\Box_AE$ and write $x=\sum_ib_ia_i$ with $b_i\in D_0$ and $a_i\in E$. Then:
\begin{eqnarray*}
E_{D_1}(x)
&=&\sum_ib_iE_{D_1}(a_i)\\
&=&\sum_ib_iE_{\mathbb C}(a_i)\\
&=&\sum_ib_iE_{D_0}(a_i)\\
&=&E_{D_0}(x)
\end{eqnarray*}

But this proves the commuting square condition, and we are done.

\medskip

\underline{Step III}. A similar argument shows that the result holds in the case $D_0=\mathbb C$.

\medskip

\underline{Step IV}. General case. We will use many times the following diagram, in which all the rectangles and all the squares, except possibly for the square in the statement, are non-degenerate commuting squares, cf. the conclusions of Steps I, II, III:
$$\begin{matrix}
E_1&\subset&D_0\Box_AE_1&\subset&D_1\Box_AE_1\\
\cup&&\cup&&\cup\\
E_0&\subset&D_0\Box_AE_0&\subset&D_1\Box_AE_0\\
\cup&&\cup&&\cup\\
\mathbb C&\subset&D_0&\subset&D_1
\end{matrix}$$

The non-degeneracy condition follows from:
$$D_1\Box_AE_1=sp\{E_1\cdot D_1\}\subset sp\{D_0\Box_AE_1\cdot D_1\Box_AE_0\}$$

Now let $x\in D_0\Box_AE_1$ and write $x=\sum_ib_ia_i$ with $b_i\in D_0$ and $a_i\in E_1$. Then:
\begin{eqnarray*}
E_{D_1\Box_AE_0}(x)
&=&\sum_ib_iE_{D_1\Box_AE_0}(a_i)\\
&=&\sum_ib_iE_{E_0}(a_i)\\
&=&\sum_ib_iE_{D_0\Box_AE_0}(a_i)\\
&=&E_{D_0\Box_AE_0}(x)
\end{eqnarray*}

But this proves the commuting square condition, and we are done. 
\end{proof}

We show now that the bifunctor $\Box_A$ behaves well with respect to basic constructions. If $A$ is a Woronowicz algebra, a sequence of two arrows $D_0\subset D_1\subset D_2$ in $A-{\mathcal Alg}$ is called a basic construction if $D_0\subset D_1\subset D_2$ is a basic construction in ${\mathcal Alg}$ and if its Jones projection $e\in D_2$ is a fixed by the coaction $D_2\to D_2\otimes A$. An infinite sequence of basic constructions in $A-{\mathcal Alg}$ is called a Jones tower in $A-{\mathcal Alg}$. We have:

\begin{proposition}
If $D_0\subset D_1\subset D_2\subset D_3\subset\ldots$ is a Jones tower in $A-{\mathcal Alg}$ and $E_0\subset E_1\subset E_2\subset E_3\subset\ldots$ is a Jones tower in $\widehat{A}_\sigma-{\mathcal Alg}$ then
$$\begin{matrix}
\vdots &\ & \vdots &\ & \vdots &\ & \ \\ 
\cup &\ &\cup &\ &\cup\\
D_0\Box_AE_2&\subset&D_1\Box_AE_2&\subset&D_2\Box_AE_2&\subset&\cdots\\
\cup &\ &\cup &\ &\cup\\
D_0\Box_AE_1&\subset&D_1\Box_AE_1&\subset&D_2\Box_AE_1&\subset&\cdots\\
\cup &\ &\cup &\ &\cup\\ 
D_0\Box_AE_0&\subset&D_1\Box_AE_0&\subset&D_2\Box_AE_0&\subset&\cdots
\end{matrix}$$
is a lattice of basic constructions for non-degenerate commuting squares.
\end{proposition}

\begin{proof}
We prove only that the rows are Jones towers, the proof for the columns being similar. By restricting the attention to a pair of consecutive inclusions, it is enough to prove that if $D_0\subset D_1\subset D_2$ is a basic construction in $A-{\mathcal Alg}$ and $E$ is an object of $\widehat{A}_\sigma-{\mathcal Alg}$ then $D_0\Box_AE\subset D_1\Box_AE\subset D_2\Box_AE$ is a basic construction in ${\mathcal Alg}$. 

\medskip

For this purpose, we will use many times the following diagram, in which all squares and rectangles are non-degenerate commuting squares:
$$\begin{matrix}
E&\subset&D_0\Box_AE&\subset&D_1\Box_AE&\subset&D_2\Box_AE\\
\cup &\ &\cup &\ &\cup &\ &\cup\\
\mathbb C&\subset&D_0&\subset&D_1&\subset&D_2
\end{matrix}$$

We will use the abstract characterization of the basic construction, stating that $N\subset M\subset P$ is a basic construction, with Jones projection $e\in P$, precisely when: 

\medskip

(1) $P=sp\{ M\cdot e\cdot M\}$.

\medskip

(2) $[e,N]=0$.

\medskip

(3) $exe=E_N(x)e$ for any $x\in M$.

\medskip

(4) $tr(xe)=\lambda tr(x)$ for any $x\in M$, where $\lambda$ is the inverse of the index of $N\subset M$. 

\medskip

Let $e\in D_2$ be the Jones projection for the basic construction $D_0\subset D_1\subset D_2$. With $N=D_0$, $M=D_1$ and $P=D_2$ the verification of (1-4) is as follows:

\medskip

(1) This follows from the following computation:
\begin{eqnarray*}
D_2\Box_AE
&=&sp\{D_2\cdot E\}\\
&=&sp\{D_1\cdot e\cdot D_1\cdot E\}\\
&=&sp\{D_1\cdot e\cdot D_1\Box_AE\}
\end{eqnarray*}

(2) This follows from $D_0\Box_AE=sp\{D_0\cdot E\}$, from $[e,E]=0$ and from $[e,D_0]=0$.

\medskip

(3) Let $x\in D_1\Box_AE$, and write $x=\sum_ib_ia_i$ with $b_i\in D_1$ and $a_i\in E$. Then:
\begin{eqnarray*}
exe
&=&\sum_ieb_ia_ie\\
&=&\sum_ieb_iea_i\\
&=&\sum_iE_{D_0}(b_i)ea_i\\
&=&\sum_iE_{D_0}(b_i)a_ie
\end{eqnarray*}

On the other hand, we have as well the following computation:
\begin{eqnarray*}
E_{D_0\Box_AE}(x)e
&=&\sum_iE_{D_0\Box_AE}(b_ia_i)e\\
&=&\sum_iE_{D_0\Box_AE}(b_i)a_ie\\
&=&\sum_iE_{D_0}(b_i)a_ie
\end{eqnarray*}

(4) With the above notations, we have that:
\begin{eqnarray*}
E_{D_2}(xe)
&=&\sum_iE_{D_2}(b_ia_ie)\\
&=&\sum_ib_iE_{D_2}(a_i)e\\
&=&\sum_ib_iE_\mathbb C(a_i)e
\end{eqnarray*}

We also have $b_iE_\mathbb C(a_i)\in D_1$ for every $i$, and so:
\begin{eqnarray*}
tr_{D_2\Box_AE}(xe)
&=&tr_{D_2}(E_{D_2}(xe))\\
&=&\lambda\sum_itr_{D_1}(b_iE_\mathbb C(a_i))
\end{eqnarray*}

On the other hand, we have as well the following computation:
\begin{eqnarray*}
tr_{D_1\Box_AE}(x)
&=&tr_{D_1}(E_{D_1}(x))\\
&=&\sum_itr_{D_1}(b_iE_{D_1}(a_i))\\
&=&\sum_itr_{D_1}(b_iE_\mathbb C(a_i))
\end{eqnarray*}

Thus, the fourth condition for a basic construction is verified, as desired.
\end{proof}

With standard coaction conventions, from chapter 13, we have:

\begin{proposition}
Given a corepresentation and a representation, as follows,
$$v\in M_n(\mathbb C)\otimes A\quad,\quad \pi:A_\sigma\to M_k(\mathbb C)$$
consider, via some standard identifications, the associated objects 
$$(M_n(\mathbb C),\iota_v)\in A-{\mathcal Alg}\quad,\quad
(M_k(\mathbb C),\iota_{\check{\pi}})\in \widehat{A}_\sigma-{\mathcal Alg}$$
and form the corresponding algebra $M_n(\mathbb C)\Box_AM_k(\mathbb C)$. Then there exists an isomorphism
$$\begin{pmatrix}
M_k(\mathbb C)&\subset&M_n(\mathbb C)\Box_AM_k(\mathbb C)\\
\cup &\ &\cup\\
\mathbb C&\subset & M_n(\mathbb C)\end{pmatrix}
\simeq
\begin{pmatrix}
\mathbb C\otimes M_k(\mathbb C) &\subset & M_n(\mathbb C)\otimes M_k(\mathbb C)\\
\cup &\ &\cup\\
\mathbb C&\subset &u(M_n(\mathbb C)\otimes\mathbb C)u^*\end{pmatrix}$$
sending $z\to 1\otimes z$ for $z\in M_k(\mathbb C)$ and $y\to\iota_u (y)$ for $y\in M_n(\mathbb C)$, where $u=(id\otimes\pi )v$.
\end{proposition}

\begin{proof}
Consider the following $*$-morphism of algebras:
$$\Phi :M_n(\mathbb C)\otimes M_k(\mathbb C)\to M_n(\mathbb C)\otimes M_k(\mathbb C)\otimes B (l^2(A_\sigma ))$$
$$x\to ad(v_{13}\check{\pi}_{23}u_{12}^*)(x\otimes 1)$$

Since both the squares in the statement are non-degenerate commuting squares, all the assertions are consequences of the following formulae, that we will prove now:
$$\Phi (1\otimes z)=\iota_{\check{\pi}} (z)_{23}\quad,\quad 
\Phi (\iota_u (y))=\iota_v(y)_{13}$$

The second formula follows from the following computation:
$$\Phi (u(y\otimes 1)u^*)
=v_{13}(y\otimes 1\otimes 1)v_{13}^*
=(v(y\otimes 1)v^*)_{13}$$

For the first formula, what we have to prove is that:
$$v_{13}\check{\pi}_{23}u_{12}^*(1\otimes z\otimes 1)u_{12}\check{\pi}_{23}^*v_{13}^*=(\check{\pi}(z\otimes 1)\check{\pi}^*)_{23}$$

By moving the unitaries to the left and to the right we have to prove that:
$$\check{\pi}_{23}^*v_{13}\check{\pi}_{23}u_{12}^*\in (\mathbb C\otimes M_k(\mathbb C)\otimes\mathbb C)^\prime =M_n(\mathbb C)\otimes\mathbb C\otimes B(l^2(H_\sigma ))$$

Let us call this unitary $U$. Since $\check{\pi}=(\pi\otimes id )V^\prime$ we have:
$$U=(id\otimes\pi\otimes id)(V_{23}^*v_{13}V_{23}v_{12}^*)$$

The comultiplication of $H_\sigma$ is given by $\Delta(y)=V^*(1\otimes y)V$. On the other hand since $v^*$ is a corepresentation of $H_\sigma$, we have $(id\otimes\Delta )(v^*)=v^*_{12}v^*_{13}$. We get that:
\begin{eqnarray*}
V_{23}^*v_{13}V_{23}
&=&(V_{23}^*v_{13}^*V_{23})^*\\
&=&((id\otimes\Delta )(v^*))^*\\
&=&(v^*_{12}v^*_{13})^*\\
&=&v_{13}v_{12}
\end{eqnarray*}

Thus we have $U=v_{13}$, and we are done.
\end{proof}

We are now ready to formulate our main result, as follows:

\begin{theorem}
Any commuting square having $\mathbb C$ in the lower left corner,
$$\begin{matrix}
E&\subset & X\\
\cup &\ &\cup\\
\mathbb C &\subset&D\end{matrix}$$
must appear as follows, for a suitable Woronowicz algebra $A$, with actions on $D,E$,
$$\begin{matrix}
E&\subset&D\Box_AE\\
\cup &\ &\cup\\
\mathbb C&\subset&D\end{matrix}$$
and the vertical subfactor associated to it is isomorphic to
$$R\subset(D\otimes(R\rtimes_{\gamma_\pi}\widehat{A}_\sigma ))^{\beta_v\odot \widehat{\gamma_\pi}}$$
which is a fixed point subfactor, in the sense of chapter 13.
\end{theorem}

\begin{proof}
This is something quite technical, which basically follows by combining the above results, and for full details on this, we refer to \cite{ba1} and related papers.
\end{proof}

Summarizing, all the commuting squares having $\mathbb C$ in the lower left corner are described by quantum groups. This is of course something quite special, and we will study more general commuting squares, not coming from quantum groups, in the next chapter.

\section*{15e. Exercises} 

Things have been quite technical here, and as an exercise on this, we have:

\begin{exercise}
Given a commuting square of finite dimensional algebras
$$\xymatrix@R=40pt@C40pt{
C_{01}\ar[r]&C_{11}\\
C_{00}\ar[u]\ar[r]&C_{10}\ar[u]}$$
establish, with full details, the Ocneanu compactness formula
$$A_0'\cap A_k=C_{01}'\cap C_{k0}$$
for the associated vertical subfactor $A_0\subset A_1$.
\end{exercise}

This is something quite fundamental, that we discussed in the above, but with the details missing. Time to have this done, by working out the linear algebra.

\chapter{Spectral measures}

\section*{16a. Small index}

We have seen so far the foundations of Jones' subfactor theory, along with results regarding the most basic classes of such subfactors, namely those coming from compact groups, discrete group duals, and more generally compact quantum groups. These subfactors all have integer index, $N\in\mathbb N$, and appear as subfactors of the Murray-von Neumann hyperfinite ${\rm II}_1$ factor $R$, either by definition, or by theorem, or by conjecture.

\bigskip

This suggests looking into the classification of subfactors of integer index, or into the classification of the subfactors of $R$, or into the classification of the subfactors of $R$ having integer index. These are all good questions, that we will discuss here.

\bigskip

Before starting, however, and in order to have an idea on what we want to do, we should discuss the following question: should the index $N\in[1,\infty)$ be small, or big? This is something quite philosophical, and non-trivial, the situation being as follows:

\medskip

\begin{enumerate}
\item Mathematics and basic common sense suggest that subfactors should fall into two main classes, ``series'' and ``exceptional''. From this perspective, the series, corresponding to uniform values of the index, must be investigated first. 

\medskip

\item In practice now, passed a few simple cases, such as the FC or TL subfactors, we cannot hope for the index to take full uniform values. The more reasonable question here is that of looking at the case where $N\in\mathbb N$ is uniform.

\medskip

\item The problem now is that, in the lack of theory here, this basically brings us back to groups, group duals, and more generally compact quantum groups, whose combinatorics is notoriously simpler than that of the arbitrary subfactors.

\medskip

\item In short, naivity and pure mathematics tell us to investigate the ``big index'' case first, but with the remark however that we are missing something, and so that we must do in parallel some study in the ``small index'' case too. 
\end{enumerate}

\medskip

All this does not look very clear, and so after this discussion, we are basically still in the dark. So, should the answer come then from physics, and applications? 

\bigskip

Unfortunately, things here are quite complicated too, basically due to our current poor understanding of quantum mechanics, and of what precisely is to be done, in order to have things in physics moving. And in fact, things here are in fact split too, a bit in the same way as above, the situation being basically as follows:

\medskip

\begin{enumerate}
\item The very small index range, $N\in[1,4]$, is subject to the remarkable ``quantization'' result of Jones, stating that we should have $N=4\cos^2(\frac{\pi}{n})$, and has strong ties with a number of considerations in conformal field theory.

\medskip

\item In what concerns the other end, $N>>0$, this is in relation with statistical mechanics, once again following work of Jones on the subject, and with the index itself corresponding to physicists' famous ``big $N$'' variable.
\end{enumerate}

\medskip

In short, no hope for an answer here. At least with our current knowledge of the subject. Probably most illustrating here is the fact that the main experts, starting with Jones himself, have always being split, working on both small and big index.

\bigskip

Getting away now from these philosophical difficulties, and back to our present book, which is rather elementary and mathematical, in this final chapter we will survey the main structure and classification results available, both in small and big index. 

\bigskip

As already mentioned, we will focus on the subfactors of the Murray-von Neumann hyperfinite ${\rm II}_1$ factor $R$, by taking for granted the fact that these subfactors are the most ``important'', and related to physics. With the side remark, however, that this is actually subject to debate too, with many mathematicians opting for bigger factors like $L(F_\infty)$, and with some physicists joining them too. But let us not get into this here.

\bigskip

In order to get started now, in order to talk about classification, we need invariants for our subfactors. Which brings us into a third controversy, namely the choice between algebraic and analytic invariants. The situation here is as follows:

\index{principal graph}
\index{fusion algebra}
\index{Poincar\'e series}
\index{spectral measure}

\begin{definition}
Associated to any finite index subfactor $A\subset B$, having planar algebra $P=(P_k)$, are the following invariants:
\begin{enumerate}
\item Its principal graph $X$, which describes the inclusions $P_0\subset P_1\subset P_2\subset\ldots\,$, with the reflections coming from basic constructions removed.

\item Its fusion algebra $F$, which describes the fusion rules for the various types of bimodules that can appear, namely $A-A$, $A-B$, $B-A$, $B-B$.

\item Its Poincar\'e series $f$, which is the generating series of the graded components of the planar algebra, $f(z)=\sum_k\dim(P_k)z^k$.

\item Its spectral measure $\mu$, which is the probability measure having as moments the dimensions of the planar algebra components, $\int x^kd\mu(x)=\dim(P_k)$.
\end{enumerate}
\end{definition}

This definition is of course something a bit informal, and there is certainly some work to be done, in order to fully define all the above invariants $X,F,f,\mu$, and to work out the precise relation between them. We will be back to this later, but for the moment, let us keep in mind the fact that associated to a given subfactor $A\subset B$ are several combinatorial invariants, which are not exactly equivalent, but are definitely versions of the same thing, the ``combinatorics of the subfactor'', and which come in algebraic or analytic flavors.

\bigskip

So, what to use? As before, in relation with the previous controversies, the main experts, starting with Jones himself, have always being split themselves on this question, working with both algebraic and analytic invariants. Generally speaking, the algebraic invariants, which are (1) and (2) in the above list, tend to be more popular in small index, while the analytic invariants, (3) and (4), are definitely more popular in big index.

\bigskip

In order to get started now, let us first discuss the question of classifying the subfactors of the hyperfinite ${\rm II}_1$ factor $R$, up to isomorphism, having index $N\leq4$. 

\bigskip

This is something quite tricky, and the main idea here will be the fact, coming from the proof of the Jones index restriction theorem, explained in chapter 13 above, that the index $N\in(1,4]$ must be the squared norm of a certain graph:
$$N=||X||^2$$

Now with this observation in hand, the assumption $N\leq4$ forces $X$ to be one of the Coxeter-Dynkin graphs of type ADE, and then a lot of work, both of classification and exclusion, leads to an ADE classification for the subfactors of $R$ having index $N\leq 4$.

\bigskip

This was for the idea. More in detail now, let us begin by explaining in detail how our subfactor invariant here, which will be the principal graph $X$, is constructed.

\bigskip

Consider first an arbitrary finite index irreducible subfactor $A_0\subset A_1$, with associated planar algebra $P_k=A_0'\cap A_k$, and let us look at the following system of inclusions:
$$P_0\subset P_1\subset P_2\subset\ldots$$

By taking the Bratelli diagram of this system of inclusions, and then deleting the reflections coming from basic constructions, which automatically appear at each step, according to the various results from chapter 13, we obtain a certain graph $X$, called principal graph of $A_0\subset A_1$. The main properties of $X$ can be summarized as follows:

\index{amenable subfactor}

\begin{theorem}
The principal graph $X$ has the following properties:
\begin{enumerate}
\item The higher relative commutant $P_k=A_0'\cap A_k$ is isomorphic to the abstract vector space spanned by the $2k$-loops on $X$ based at the root.

\item In the amenable case, where $A_1=R$ and when the subfactor is ``amenable'', the index of $A_0\subset A_1$ is given by $N=||X||^2$.
\end{enumerate} 
\end{theorem}

\begin{proof}
This is something standard, the idea being as follows:

\medskip

(1) The statement here, which explains among others the relation between the principal graph $X$, and the other subfactor invariants, from Definition 16.1 above, comes from the definition of the principal graph, as a Bratelli diagram, with the reflections removed.

\medskip

(2) This is actually a quite subtle statement, but for our purposes here, we can take the equality $N=||X||^2$, which reminds a bit the Kesten amenability condition for discrete groups, as a definition for the amenability of the subfactor. With the remark that for the Popa diagonal subfactors what we have here is precisely the Kesten amenability condition for the underlying discrete group $\Gamma$, and that, more generally, for the arbitrary generalized Popa or Wassermann subfactors, what we have here is precisely the Kesten type amenability condition for the underlying discrete quantum group $\Gamma$.
\end{proof}

As a consequence of the above, in relation with classification questions, we have:

\index{ADE}
\index{Coxeter-Dynkin}

\begin{theorem}
The principal graph of a subfactor having index $N\leq 4$ must be one of the Coxeter-Dynkin graphs of type ADE.
\end{theorem}

\begin{proof}
This follows indeed from the formula $N=||X||^2$ from the above result, and from the considerations from the proof of the Jones index restriction theorem, explained in chapter 13 above. For full details on all this, we refer for instance to \cite{ghj}.
\end{proof}

More in detail now, the usual Coxeter-Dynkin graphs are as follows:
$$A_n=\bullet-\circ-\circ\cdots\circ-\circ-\circ\hskip20mm A_{\infty}=\bullet-\circ-\circ-\circ\cdots\hskip7mm$$
\vskip-7mm
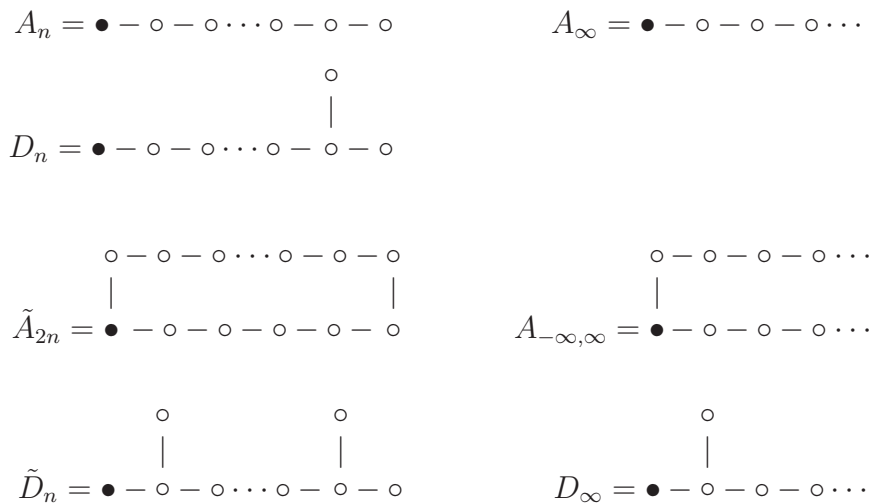
$$D_n=\bullet-\circ-\circ\dots\circ-
\begin{matrix}\ \circ\cr\ |\cr\ \circ \cr\ \cr\  \end{matrix}-\circ\hskip70mm$$
\vskip-7mm
$$\ \ \ \ \ \ \ \tilde{A}_{2n}=
\begin{matrix}
\circ&\!\!\!\!-\circ-\circ\cdots\circ-\circ-&\!\!\!\!\circ\cr
|&&\!\!\!\!|\cr
\bullet&\!\!\!\!-\circ-\circ-\circ-\circ-&\!\!\!\!\circ\cr\cr\cr\end{matrix}\hskip15mm A_{-\infty,\infty}=
\begin{matrix}
\circ&\!\!\!\!-\circ-\circ-\circ\cdots\cr
|&\cr
\bullet&\!\!\!\!-\circ-\circ-\circ\cdots\cr\cr\cr\end{matrix}
\hskip15mm$$
\vskip-9mm
$$\;\tilde{D}_n=\bullet-
\begin{matrix}\circ\cr|\cr\circ\cr\ \cr\ \end{matrix}-\circ\dots\circ-
\begin{matrix}\ \circ\cr\ |\cr\ \circ \cr\ \cr\  \end{matrix}-\circ \hskip20mm D_\infty=\bullet-
\begin{matrix}\circ\cr|\cr\circ\cr\ \cr\ \end{matrix}-\circ-\circ\cdots\hskip7mm$$
\vskip-7mm

Here the graphs $A_n$ with $n\geq 2$ and $D_n$ with $n\geq 3$ have by definition $n$ vertices each, $\tilde{A}_{2n}$ with $n\geq 1$ has $2n$ vertices, and $\tilde{D}_n$ with $n\geq 4$ has $n+1$ vertices. Thus, the first graph in each series is by definition as follows:
$$A_2=\bullet-\circ\hskip10mm 
D_3=\begin{matrix}\ \circ\cr\ |\cr\ \bullet \cr\ \cr\  \end{matrix}-\circ \hskip10mm
\tilde{A}_2=\begin{matrix}
\circ\cr
||\cr
\bullet\cr&\cr&\cr\end{matrix}\hskip10mm 
\tilde{D}_4=\bullet-\!\!\!\!\!\begin{matrix}
\circ\hskip5mm \circ\cr
\backslash\ \,\slash\cr
\circ\cr&\cr&\cr\end{matrix}\!\!\!\!\!\!\!\!\!\!-\circ$$
\vskip-7mm

There are also a number of exceptional Coxeter-Dynkin graphs. First we have:
$$E_6=\bullet-\circ-
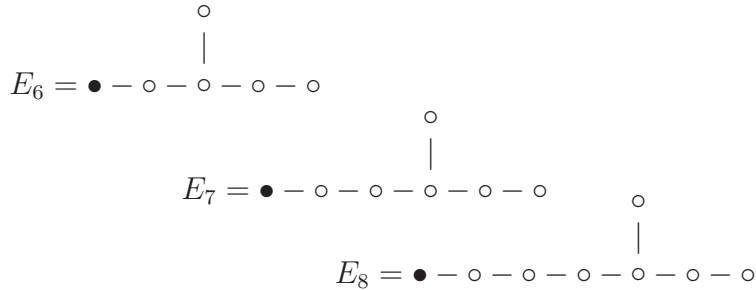
\begin{matrix}\circ\cr|\cr\circ\cr\ \cr\ \end{matrix}-
\circ-\circ\hskip71mm$$
\vskip-13mm
$$E_7=\bullet-\circ-\circ-
\begin{matrix}\circ\cr|\cr\circ\cr\ \cr\ \end{matrix}-
\circ-\circ\hskip18mm$$
\vskip-15mm
$$\hskip30mm E_8=\bullet-\circ-\circ-\circ-
\begin{matrix}\circ\cr|\cr\circ\cr\ \cr\ \end{matrix}-
\circ-\circ$$
\vskip-5mm

Also, we have as well index 4 versions of the above exceptional graphs, as follows:
$$\tilde{E}_6=\bullet-\circ-\begin{matrix}
\circ\cr|
\cr\circ\cr|&\cr\circ&\!\!\!\!-\ \circ\cr\ \cr\   \cr\ \cr\ \end{matrix}-\circ\hskip71mm$$
\vskip-22mm
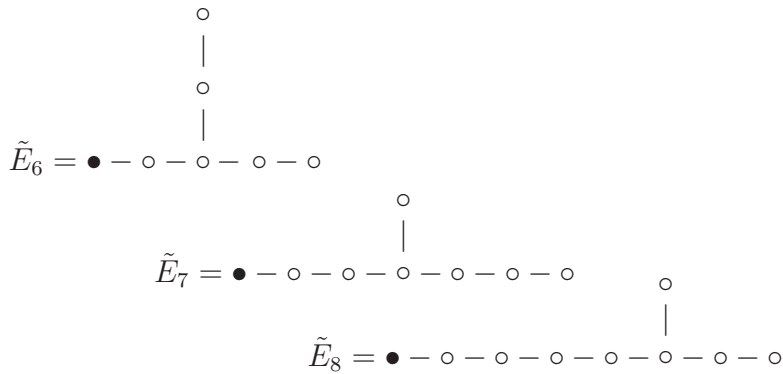
$$\tilde{E}_7=\bullet-\circ-\circ-
\begin{matrix}\circ\cr|\cr\circ\cr\ \cr\ \end{matrix}-
\circ-\circ-\circ\hskip18mm$$
\vskip-15mm
$$\hskip30mm \tilde{E}_8=\bullet-\circ-\circ-\circ-\circ-
\begin{matrix}\circ\cr|\cr\circ\cr\ \cr\ \end{matrix}-
\circ-\circ$$
\vskip-5mm

Getting back now to Theorem 16.3, with this list in hand, the story is not over here, because we still have to understand which of these graphs can really appear as principal graphs of subfactors. And, for those graphs which can appear, we must understand the structure and classification of the subfactors of $R$, having them as principal graphs.

\bigskip

In short, there is still a lot of work to be done, as a continuation of Theorem 16.3. The subfactors of index $\leq 4$ were intensively studied in the 80s and early 90s, and about 10 years after Jones' foundational paper \cite{jo1}, a complete classification result was found, with contributions by many authors. A simplified form of this result is as follows:

\index{ADE}
\index{principal graph}
\index{index theorem}

\begin{theorem}
The principal graphs of subfactors of index $\leq 4$ are:
\begin{enumerate}
\item Index $<4$ graphs: $A_n$, $D_{even}$, $E_6$, $E_8$. 

\item Index $4$ finite graphs: $\tilde{A}_{2n}$, $\tilde{D}_n$, $\tilde{E}_6$, $\tilde{E}_7$, $\tilde{E}_8$.

\item Index $4$ infinite graphs: $A_\infty$, $A_{-\infty,\infty}$, $D_\infty$.
\end{enumerate}
\end{theorem}

\begin{proof}
As already mentioned, this is something quite heavy, with contributions by many authors, and among the main papers to be read here, let us mention \cite{jo1}, \cite{oc1}, \cite{oc2}, \cite{po2}. Observe that the graphs $D_{odd}$ and $E_7$ don't appear in the above list. This is one of the subtle points of subfactor theory. For a discussion here, see \cite{eka}.
\end{proof}

There are many other things that can be said about the subfactors of index $N\leq4$, both at the theoretical level, of the finite depth and more generally of the amenable subfactors, and at the level of the ADE classification, which makes connections with other ADE classifications. We refer here to \cite{eka}, \cite{ghj}, \cite{oc1}, \cite{oc2}, \cite{po2}, \cite{po3}.

\bigskip

Regarding now the subfactors of index $N\in(4,5]$, and also of small index above 5, these can be classified, but this is a long and complicated story. Let us just record here the result in index 5, which is something quite easy to formulate, as follows:

\begin{theorem}
The principal graphs of the irreducible index $5$ subfactors are:
\begin{enumerate}
\item $A_\infty$, and a non-extremal perturbation of $A_\infty^{(1)}$.

\item The McKay graphs of $\mathbb Z_5,D_5,GA_1(5),A_5,S_5$.

\item The twists of the McKay graphs of $A_5,S_5$.
\end{enumerate}
\end{theorem}

\begin{proof}
This is a heavy result, and we refer to \cite{jms} for the whole story. The above formulation is the one from \cite{jms}, with the subgroup subfactors there replaced by fixed point subfactors, and with the cyclic groups denoted as usual by $\mathbb Z_N$.
\end{proof}

As a comment here, the above $N=5$ result was much harder to obtain than the classification in index $N=4$, obtained as a consequence of Theorem 16.4. However, at the level of the explicit construction of such subfactors, things are quite similar at $N=4$ and $N=5$, with the fixed point subfactors associated to quantum permutation groups $G\subset S_N^+$ providing most of the examples. We refer here to \cite{bbc} and related papers.

\bigskip

In index $N=6$ now, the subfactors cannot be classified, at least in general, due to several uncountable families, coming from groups, group duals, and more generally compact quantum groups. The exact assumption to be added is not known yet.

\bigskip

Summarizing, the current small index classification problem meets considerable difficulties in index $N=6$, and right below. In small index $N>6$ the situation is largely unexplored. We refer here to \cite{lmp} and the recent literature on the subject.

\section*{16b. Spectral measures}

Before getting into the case where the index is big, $N>>0$, let us comment on one of the key ingredients for the above classification results, at $N<6$. This is the Jones annular theory of subfactors, which is something very beautiful and useful, regarding the case where the index is arbitrary, $N\in[1,\infty)$. The main result is as follows:

\index{annular category}
\index{theta series}

\begin{theorem}
The theta series of a subfactor of index $N>4$, which is given by
$$\Theta (q)=q+\frac{1-q}{1+q}\, f\left( \frac{q}{(1+q)^2}\right)$$
with $f=\sum_k\dim(P_k)z^k$ being the Poincar\'e series, has positive coefficients.
\end{theorem}

\begin{proof}
This is something quite advanced, the idea being that $\Theta$ is the generating series of a certain series of multiplicities associated to the subfactor, and more specifically associated to the canonical inclusion $TL_N\subset P$. We refer here to Jones' paper \cite{jo5}.
\end{proof}

In relation to this, and to some questions from physics as well, coming from conformal field theory, an interesting question is that of computing the ``blowup'' of the spectral measure of the subfactor, via the Jones change of variables, namely:
$$z\to\frac{q}{1+q^2}$$

This question makes sense in any index, meaning both $N\in[1,4]$, where Theorem 16.6 does not apply, and $N\in(4,\infty)$, where Theorem 16.6 does apply. We will discuss in what follows both these questions, by starting with the small index one, $N\in[1,4]$.

\bigskip

Following \cite{bdb} and related papers, it is convenient to stay, at least for the beginning, at a purely elementary level, and associate such series to any rooted bipartite graph. Let us start with the following definition, which is something straightforward, inspired by the definition of the Poincar\'e series of a subfactor, and by Theorem 16.2:

\index{Poincar\'e series}

\begin{definition}
The Poincar\'e series of a rooted bipartite graph $X$ is
$$f(z)=\sum_{k=0}^\infty{\rm loop}_X(2k)z^k$$
where ${\rm loop}_X(2k)$ is the number of $2k$-loops based at the root.
\end{definition}

In the case where $X$ is the principal graph of a subfactor $A_0\subset A_1$, this series $f$ is the Poincar\'e series of the subfactor, in the usual sense:
$$f(z)=\sum_{k=0}^\infty\dim(A_0'\cap A_k)z^k$$

In general, the Poincar\'e series should be thought of as being a basic representation theory invariant of the underlying group-like object. For instance for the Wassermann type subfactor associated to a compact Lie group $G\subset U_N$, the Poincar\'e series is:
$$f(z)=\int_G\frac{1}{1-Tr(g)z}\,dg$$

Regarding now the theta series, this can introduced as a version of the Poincar\'e series, via the change of variables $z^{-1/2}=q^{1/2}+q^{-1/2}$, as follows:

\index{theta series}

\begin{definition}
The theta series of a rooted bipartite graph $X$ is
$$\Theta(q)=q+\frac{1-q}{1+q}f\left(\frac{q}{(1+q)^2}\right)$$
where $f$ is the Poincar\'e series.
\end{definition}

The theta series can be written as $\Theta(q)=\sum a_rq^r$, and it follows from the above formula, via some simple manipulations, that its coefficients are integers:
$$a_r\in\mathbb Z$$

In fact, we have the following explicit formula from Jones' paper \cite{jo5}, relating the coefficients of $\Theta(q)=\sum a_rq^r$ to those of the Poincar\'e series $f(z)=\sum c_kz^k$:
$$a_r=\sum_{k=0}^r(-1)^{r-k}\frac{2r}{r+k}\begin{pmatrix}r+k\cr r-k\end{pmatrix}c_k$$

In the case where $X$ is the principal graph of a subfactor $A_0\subset A_1$ of index $N>4$, it is known from \cite{jo5} that the numbers $a_r$ are certain multiplicities associated to the planar algebra inclusion $TL_N\subset P$, as explained in Theorem 16.6 and its proof. In particular, the coefficients of the theta series are in this case positive integers:
$$a_r\in\mathbb N$$

Before getting into computations, let us discuss as well the measure-theoretic versions of the above invariants. Once again, we start with an arbitrary rooted bipartite graph $X$. We can first introduce a measure $\mu$, whose Stieltjes transform is $f$, as follows:

\begin{definition}
The real measure $\mu$ of a rooted bipartite graph $X$ is given by
$$f(z)=\int_0^\infty\frac{1}{1-xz}\,d\mu(x)$$
where $f$ is the Poincar\'e series.
\end{definition}

In the case where $X$ is the principal graph of a subfactor $A_0\subset A_1$, we recover in this way the spectral measure of the subfactor, as introduced in Definition 16.1, with the remark however that the existence of such a measure $\mu$ was not discussed  there. In general, and so also in the particular subfactor case, clarifying the things here, the fact that $\mu$ as above exists indeed comes from the following simple fact:

\begin{proposition}
The real measure $\mu$ of a rooted bipartite graph $X$ is given by the following formula, where $L=MM^t$, with $M$ being the adjacency matrix of the graph,
$$\mu=law(L)$$
and with the probabilistic computation being with respect to the expectation 
$$A\to<A>$$
with $<A>$ being the $(*,*)$-entry of a matrix $A$, where $*$ is the root.
\end{proposition}

\begin{proof}
With the conventions in the statement, namely  $L=MM^t$, with $M$ being the adjacency matrix, and with $<A>$ being the $(*,*)$-entry of a matrix $A$, we have:
\begin{eqnarray*}
f(z)
&=&\sum_{k=0}^\infty{\rm loop}_X(2k)z^k\\
&=&\sum_{k=0}^\infty\left<L^k\right>z^k\\
&=&\left<\frac{1}{1-Lz}\right>
\end{eqnarray*}

But this shows that we have the formula $\mu=law(L)$, as desired.
\end{proof}

In the subfactor case some further interpretations are available as well. For instance in the case of the fixed point subfactors coming from of a compact group $G\subset U_N$, discussed after Definition 16.7 above, $\mu$ is the spectral measure of the main character:
$$\mu=law(\chi)$$

In relation now with the theta series, things are more tricky, in order to introduce its measure-theoretic version. Following \cite{bdb}, let us introduce the following notion:

\index{circular measure}

\begin{definition}
The circular measure $\varepsilon$ of a rooted bipartite graph $X$ is given by
$$d\varepsilon(q)=d\mu((q+q^{-1})^2)$$
where $\mu$ is the associated real measure.
\end{definition}

In other words, the circular measure $\varepsilon$ is by definition the pullback of the usual real measure $\mu$ via the following map, coming from the theory of the theta series in \cite{jo5}:
$$\mathbb R\cup\mathbb T\to\mathbb R_+$$
$$q\to (q+q^{-1})^2$$

As we will see, all this best works in index $N\in[1,4]$, with the circular measure $\varepsilon$ being here the best-looking invariant, among all subfactor invariants. In index $N>4$ things will turn to be quite complicated, but more on this later.

\bigskip

As a basic example for all this, assume that $\mu$ is a discrete measure, supported by $n$ positive numbers $x_1<\ldots<x_n$, with corresponding densities $p_1,\ldots,p_n$:
$$\mu=\sum_{i=1}^n p_i\delta_{x_i}$$

For each $i\in\{1,\ldots,n\}$ the equation $(q+q^{-1})^2=x_i$ has four solutions, that we can denote $q_i,q_i^{-1},-q_i,-q_i^{-1}$. With this notation, we have:
$$\varepsilon=\frac{1}{4}\sum_{i=1}^np_i\left(\delta_{q_i}+\delta_{q_i^{-1}}+\delta_{-q_i}+\delta_{-q_i^{-1}}\right)$$

In general, the basic properties of $\varepsilon$ can be summarized as follows:

\begin{proposition}
The circular measure has the following properties:
\begin{enumerate}
\item $\varepsilon$ has equal density at $q,q^{-1},-q,-q^{-1}$.

\item The odd moments of $\varepsilon$ are $0$.

\item The even moments of $\varepsilon$ are half-integers.

\item When $X$ has norm $\leq 2$, $\varepsilon$ is supported by the unit circle.

\item When $X$ is finite, $\varepsilon$ is discrete.

\item If $K$ is a solution of $L=(K+K^{-1})^2$, then $\varepsilon={\rm law}(K)$. 
\end{enumerate}
\end{proposition}

\begin{proof}
These results can be deduced from definitions, the idea being that (1-5) are trivial, and that (6) follows from the formula of $\mu$ from Proposition 16.10.
\end{proof}

In addition to the above, we have the following key formula, which gives the even moments of $\varepsilon$, and makes the connection with the Jones theta series:

\begin{theorem}
We have the Stieltjes transform type formula
$$2\int\frac{1}{1-qu^2}\,d\varepsilon(u)=1+T(q)(1-q)$$
where the $T$ series of a rooted bipartite graph $X$ is by definition given by
$$T(q)=\frac{\Theta(q)-q}{1-q}$$
with $\Theta$ being the associated theta series.
\end{theorem}

\begin{proof}
This follows by applying the change of variables $q\to (q+q^{-1})^2$ to the fact that $f$ is the Stieltjes transform of $\mu$. Indeed, we obtain in this way:
\begin{eqnarray*}
2\int\frac{1}{1-qu^2}\,d\varepsilon(u)
&=&1+\frac{1-q}{1+q}f\left(\frac{q}{(1+q)^2}\right)\\
&=&1+\Theta(q)-q\\
&=&1+T(q)(1-q)
\end{eqnarray*}

Thus, we are led to the conclusion in the statement.
\end{proof}

As a final theoretical result about all these invariants, which is this time something non-trivial, in the subfactor case, we have the following result, due to Jones \cite{jo5}:

\begin{theorem}
In the case where $X$ is the principal graph of an irreducible subfactor of index $>4$, the moments of $\varepsilon$ are positive numbers.
\end{theorem}

\begin{proof}
This follows indeed from the result in \cite{jo5} that the coefficients of $\Theta$ are positive numbers, as explained in Theorem 16.6, via the formula in Theorem 16.13.
\end{proof}

Summarizing, we have a whole menagery of subfactor, planar algebra and bipartite graph invariants, which come in several flavors, namely series and measures, and which can be linear or circular, and which all appear as versions of the Poincar\'e series.

\bigskip

Our claim now is that the circular measure $\varepsilon$ is the ``best'' invariant. As a first justification for this claim, let us compute $\varepsilon$ for the simplest possible graph in the index range $N\in[1,4]$, namely the graph $\tilde{A}_{2n}$. We obtain here something nice, as follows:

\begin{theorem}
The circular measure of the basic index $4$ graph, namely 
$$\begin{matrix}
&\circ&\!\!\!\!-\circ-\circ\cdots\circ-\circ-&\!\!\!\!\circ\cr
\tilde{A}_{2n}=&|&&\!\!\!\!|\cr
&\bullet&\!\!\!\!-\circ-\circ-\circ-\circ-&\!\!\!\!\circ\cr\cr\cr\end{matrix}$$
\vskip-7mm

\noindent is the uniform measure on the $2n$-roots of unity.
\end{theorem}

\begin{proof}
Let us identify the vertices of $X=\tilde{A}_{2n}$ with the group $\{w^k\}$ formed by the $2n$-th roots of unity in the complex plane, where $w=e^{\pi i/n}$. The adjacency matrix of $X$ acts then on the functions $f\in C(X)$ in the following way:
$$Mf(w^s)=f(w^{s-1})+f(w^{s+1})$$

But this shows that we have $M=K+K^{-1}$, where $K$ is given by:
$$Kf(w^s)=f(w^{s+1})$$

Thus we can use the last assertion in Proposition 16.12, and we get $\varepsilon={\rm law}(K)$, which is the uniform measure on the $2n$-roots of unity. See \cite{bdb} for details.
\end{proof}

In order to discuss all this more systematically, and for all the ADE graphs, the idea will be that of looking at the combinatorics of the roots of unity. Let us introduce:

\index{cyclotomic series}

\begin{definition}
The series of the form
$$\xi(n_1,\ldots,n_s:m_1,\ldots,m_t)=\frac{(1-q^{n_1})\ldots(1-q^{n_s})}{(1-q^{m_1})\ldots(1-q^{m_t})}$$
with $n_i,m_i\in\mathbb N$ are called cyclotomic.
\end{definition}

It is technically convenient to allow as well $1+q^n$ factors, to be designated by $n^+$ symbols in the above writing. For instance we have, by definition:
$$\xi(2^+:3)=\xi(4:2,3)$$

Also, it is convenient in what follows to use the following notations:
$$\xi'=\frac{\xi}{1-q}\quad,\quad \xi''=\frac{\xi}{1-q^2}$$

The Poincar\'e series of the ADE graphs are given by quite complicated formulae. However, the corresponding $T$ series are all cyclotomic, as follows:

\begin{theorem}
The $T$ series of the ADE graphs are as follows:
\begin{enumerate}
\item For $A_{n-1}$ we have $T=\xi(n-1:n)$.

\item For $D_{n+1}$ we have $T=\xi(n-1^+:n^+)$.

\item For $\tilde{A}_{2n}$ we have $T=\xi'(n^+:n)$.

\item For $\tilde{D}_{n+2}$ we have $T=\xi''(n+1^+:n)$.

\item For $E_6$ we have $T=\xi(8:3,6^+)$.

\item For $E_7$ we have $T=\xi(12:4,9^+)$.

\item For $E_8$ we have $T=\xi(5^+,9^+:15^+)$.

\item For $\tilde{E}_6$ we have $T=\xi(6^+:3,4)$.

\item For $\tilde{E}_7$ we have $T=\xi(9^+:4,6)$.

\item For $\tilde{E}_8$ we have $T=\xi(15^+:6,10)$.
\end{enumerate}
\end{theorem}

\begin{proof}
These formulae were obtained in \cite{bdb}, by counting loops, then by making the change of variables $z^{-1/2}=q^{1/2}+q^{-1/2}$, and factorizing the resulting series. An alternative proof for these formulae can be obtained by using planar algebra methods.
\end{proof}

Our purpose now will be that of converting the above technical results, regarding the $T$ series, into some final results, regarding the corresponding circular measures $\varepsilon$. For this purpose, we will use the conversion formula in Theorem 16.13.

\bigskip

In order to formulate our results, we will need some more theory. First, we have:

\index{cyclotomic measure}

\begin{definition}
A cyclotomic measure is a probability measure $\varepsilon$ on the unit circle, having the following properties:
\begin{enumerate}
\item  $\varepsilon$ is supported by the $2n$-roots of unity, for some $n\in\mathbb N$.

\item $\varepsilon$ has equal density at $q,q^{-1},-q,-q^{-1}$.
\end{enumerate}
\end{definition}

It follows from Theorem 16.17 that the circular measures of the finite ADE graphs are supported by certain roots of unity, hence are cyclotomic. We will be back to this.

\bigskip

At the general level now, let us introduce as well the following notion:

\begin{definition}
The $T$ series of a cyclotomic measure $\varepsilon$ is given by:
$$1+T(q)(1-q)=2\int\frac{1}{1-qu^2}\,d\varepsilon(u)$$
\end{definition}

Observe that this formula is nothing but the one in Theorem 16.13, written now in the other sense. In other words, if the cyclotomic measure $\varepsilon$ happens to be the circular measure of a rooted bipartite graph, then the $T$ series as defined above coincides with the $T$ series as defined before. This is useful for explicit computations.

\bigskip

We are now ready to discuss the circular measures of the various ADE graphs. The idea is that these measures are all cyclotomic, of level $\leq 3$, and can be expressed in terms of the basic polynomial densities of degree $\leq 6$, namely:
$$\alpha=Re(1-q^2)$$
$$\beta=Re(1-q^4)$$
$$\gamma=Re(1-q^6)$$

To be more precise, we have the following result, with $\alpha,\beta,\gamma$ being as above, with $d_n$ being the uniform measure on the $2n$-th roots of unity, and with $d_n'=2d_{2n}-d_n$ being the uniform measure on the odd $4n$-roots of unity:

\index{ADE}
\index{circular measure}

\begin{theorem}
The circular measures of the ADE graphs are given by:
\begin{enumerate}
\item $A_{n-1}\to\alpha_n$.

\item $\tilde{A}_{2n}\to d_n$.

\item $D_{n+1}\to\alpha_n'$.

\item $\tilde{D}_{n+2}\to (d_n+d_1')/2$.

\item $E_6\to\alpha_{12}+(d_{12}-d_6-d_4+d_3)/2$.

\item $E_7\to\beta_9'+(d_1'-d_3')/2$.

\item $E_8\to\alpha_{15}'+\gamma_{15}'-(d_5'+d_3')/2$.

\item $\tilde{E}_{n+3}\to (d_n+d_3+d_2-d_1)/2$.
\end{enumerate}
\end{theorem}

\begin{proof}
This follows from the $T$ series formulae in Theorem 16.17, via some routine manipulations, based on the general conversion formulae given above.
\end{proof}

It is possible to further build along the above lines, with a combinatorial refinement of the formulae in Theorem 16.20, making appear a certain connection with the Deligne work on the exceptional series of Lie groups, which is not understood yet.

\section*{16c. Measure blowup}

All the above, which was quite nice, was about index $N\in[1,4]$, where the Jones annular theory result from \cite{jo5} does not apply. In higher index now, $N\in(4,\infty)$, where the Jones result does apply, the precise correct ``blowup'' manipulation on the spectral measure is not known yet. The known results here are as follows:

\medskip

\begin{enumerate}
\item One one hand, there is as a computation for some basic Hadamard subfactors, with nice blowup, on a certain noncommutative manifold \cite{bbi}.

\medskip

\item On the other hand, there are many computations by Evans-Pugh, with quite technical blowup results, on some suitable real algebraic manifolds \cite{epu}. 
\end{enumerate}

\medskip

We will discuss in what follows (1), and to be more precise the computation of the spectral measure, and then the blowup problem, for the subfactors coming from the deformed Fourier matrices. Let us start with the following definition:

\index{deformed Fourier matrix}

\begin{definition}
Given two finite abelian groups $G,H$, we consider the corresponding deformed Fourier matrix, given by the formula
$$(F_G\otimes_Q F_H)_{ia,jb}=Q_{ib}(F_G)_{ij}(F_H)_{ab}$$
and we factorize the associated representation $\pi_Q$ of the algebra $C(S_{G\times H}^+)$,
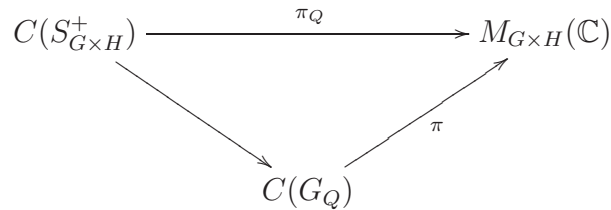
$$\xymatrix@R=40pt@C=40pt
{C(S_{G\times H}^+)\ar[rr]^{\pi_Q}\ar[rd]&&M_{G\times H}(\mathbb C)\\&C(G_Q)\ar[ur]_\pi&}$$
with $C(G_Q)$ being the Hopf image of this representation $\pi_Q$.
\end{definition}

Explicitely computing the above quantum permutation group $G_Q\subset S_{G\times H}^+$, as function of the parameter matrix $Q\in M_{G\times H}(\mathbb T)$, will be our main purpose, in what follows. In order to do so, we first have the following elementary result:

\begin{proposition}
We have a factorization as follows,
$$
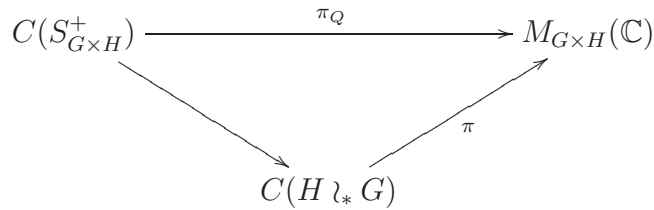
\xymatrix@R=40pt@C=40pt
{C(S_{G\times H}^+)\ar[rr]^{\pi_Q}\ar[rd]&&M_{G\times H}(\mathbb C)\\&C(H\wr_*G)\ar[ur]_\pi&}$$
given on the standard generators by the formulae
$$U_{ab}^{(i)}=\sum_jW_{ia,jb}\quad,\quad 
V_{ij}=\sum_aW_{ia,jb}$$
independently of $b$, where $W$ is the magic matrix producing $\pi_Q$.
\end{proposition}

\begin{proof}
With $K=F_G,L=F_H$ and $M=|G|,N=|H|$, the formula of the magic matrix $W\in M_{G\times H}(M_{G\times H}(\mathbb C))$ associated to $H=K\otimes_QL$ is as follows:
\begin{eqnarray*}
(W_{ia,jb})_{kc,ld}
&=&\frac{1}{MN}\cdot\frac{Q_{ic}Q_{jd}}{Q_{id}Q_{jc}}\cdot\frac{K_{ik}K_{jl}}{K_{il}K_{jk}}\cdot\frac{L_{ac}L_{bd}}{L_{ad}L_{bc}}\\
&=&\frac{1}{MN}\cdot\frac{Q_{ic}Q_{jd}}{Q_{id}Q_{jc}}\cdot K_{i-j,k-l}L_{a-b,c-d}
\end{eqnarray*}

Our claim now is that the representation $\pi_Q$ constructed in Definition 16.21 can be factorized in three steps, up to the factorization in the statement, as follows:
$$\xymatrix@R=60pt@C=50pt
{C(S_{G\times H}^+)\ar[rr]^{\pi_Q}\ar[d]&&M_{G\times H}(\mathbb C)\\
C(S_H^+\wr_*S_G^+)\ar[r]\ar@{.>}[rru]&C(S_H^+\wr_*G)\ar[r]\ar@{.>}[ur]&C(H\wr_*G)\ar@{.>}[u]}$$

Indeed, the construction of the map on the left is standard. Regarding the second factorization, this comes from the fact that since the elements $V_{ij}$ depend on $i-j$, they satisfy the defining relations for the quotient algebra $C(S_G^+)\to C(G)$. Finally, regarding the third factorization, observe that $W_{ia,jb}$ depends only on $i,j$ and on $a-b$. By summing over $j$ we obtain that the elements $U_{ab}^{(i)}$ depend only on $a-b$, and we are done.
\end{proof}

We have now all needed ingredients for refining Proposition 16.22, as follows:

\begin{proposition}
We have a factorization as follows,
$$\xymatrix@R=40pt@C=30pt
{C(S_{G\times H}^+)\ar[rr]^{\pi_Q}\ar[rd]&&M_{G\times H}(\mathbb C)\\&C^*(\Gamma_{G,H})\rtimes C(G)\ar[ur]_\rho&}$$
where the group on the bottom is given by:
$$\Gamma_{G,H}=H^{*G}\Big/\left<[c_1^{(i_1)}\ldots c_s^{(i_s)},d_1^{(j_1)}\ldots d_s^{(j_s)}]=1\Big|\sum_rc_r=\sum_rd_r=0\right>$$
\end{proposition}

\begin{proof}
Assume that we have a representation, as follows:
$$\pi:C^*(\Gamma)\rtimes C(G)\to M_L(\mathbb C)$$

Let $\Lambda$ be a $G$-stable normal subgroup of $\Gamma$, so that $G$ acts on $\Gamma/\Lambda$, and we can form the product $C^*(\Gamma/\Lambda)\rtimes C(G)$, and assume that $\pi$ is trivial on $\Lambda$. Then $\pi$ factorizes as:
$$\xymatrix@R=40pt@C=30pt
{C^*(\Gamma)\rtimes C(G)\ar[rr]^\pi\ar[rd]&&M_L(\mathbb C)\\&C^*(\Gamma/\Lambda)\rtimes C(G)\ar[ur]_\rho}$$

With $\Gamma=H^{*G}$, this gives the result.
\end{proof}

We have now all the needed ingredients for proving a main result, as follows:

\begin{theorem}
When $Q$ is generic, the minimal factorization for $\pi_Q$ is 
$$\xymatrix@R=45pt@C=40pt
{C(S_{G\times H}^+)\ar[rr]^{\pi_Q}\ar[rd]&&M_{G\times H}(\mathbb C)\\&C^*(\Gamma_{G,H})\rtimes C(G)\ar[ur]_\pi&}$$
where on the bottom
$$\Gamma_{G,H}\simeq\mathbb Z^{(|G|-1)(|H|-1)}\rtimes H$$
is the discrete group constructed above.
\end{theorem}

\begin{proof}
Consider the factorization in Proposition 16.23, which is as follows, where $L$ denotes the Hopf image of $\pi_Q$:
$$\theta : C^*(\Gamma_{G,H})\rtimes C(G)\to L$$

To be more precise, this morphism produces the following commutative diagram:
$$\xymatrix@R=40pt@C=40pt
{C(S_{G\times H}^+) \ar[rr]^{\pi_Q} \ar[dr]_{} \ar@/_/[ddr]_{}& & M_{G\times H}(\mathbb C) \\
& L \ar[ur]_{} & \\
& C^*(\Gamma_{G,H})\rtimes C(G) \ar@{-->}[u]_\theta \ar@/_/[uur]_{\pi}& 
}$$

The first observation is that the injectivity assumption on $C(G)$ holds by construction, and that for $f \in C(G)$, the matrix $\pi(f)$ is ``block scalar''. Now for $r \in \Gamma_{G,H}$ with $\theta(r\otimes 1)=\theta(1 \otimes f)$ for some $f \in C(G)$, we see, using the commutative diagram, that $\pi(r \otimes 1)$ is block scalar. Thus, modulo some standard algebra, we are done.
\end{proof}

Summarizing, we have computed the quantum permutation groups associated to the Di\c t\u a deformations of the tensor products of Fourier matrices, in the case where the deformation matrix $Q$ is generic. For some further computations, in the case where the deformation matrix $Q$ is no longer generic, we refer to the follow-ups of \cite{bbi}.

\bigskip

Let us compute now the Kesten measure $\mu=law(\chi)$, in the case where the deformation matrix is generic, as before. Our results here will be a combinatorial moment formula, a geometric interpretation of it, and an asymptotic result. We first have:

\begin{theorem}
We have the moment formula
$$\int\chi^p
=\frac{1}{|G|\cdot|H|}\#\left\{\begin{matrix}i_1,\ldots,i_p\in G\\ d_1,\ldots,d_p\in H\end{matrix}\Big|\begin{matrix}[(i_1,d_1),(i_2,d_2),\ldots,(i_p,d_p)]\ \ \ \ \\=[(i_1,d_p),(i_2,d_1),\ldots,(i_p,d_{p-1})]\end{matrix}\right\}$$
where the sets between square brackets are by definition sets with repetition.
\end{theorem}

\begin{proof}
According to the various formulae above, the factorization found in Theorem 16.24 is, at the level of standard generators, as follows:
$$\begin{matrix}
C(S_{G\times H}^+)&\to&C^*(\Gamma_{G,H})\otimes C(G)&\to&M_{G\times H}(\mathbb C)\\
u_{ia,jb}&\to&\frac{1}{|H|}\sum_cF_{b-a,c}c^{(i)}\otimes v_{ij}&\to&W_{ia,jb}
\end{matrix}$$

Thus, the main character of the quantum permutation group that we found in Theorem 16.24 is given by the following formula:
\begin{eqnarray*}
\chi
&=&\frac{1}{|H|}\sum_{iac}c^{(i)}\otimes v_{ii}\\
&=&\sum_{ic}c^{(i)}\otimes v_{ii}\\
&=&\left(\sum_{ic}c^{(i)}\right)\otimes\delta_1
\end{eqnarray*}

Now since the Haar functional of $C^*(\Gamma)\rtimes C(H)$ is the tensor product of the Haar functionals of $C^*(\Gamma),C(H)$, this gives the following formula, valid for any $p\geq1$:
$$\int\chi^p=\frac{1}{|G|}\int_{\widehat{\Gamma}_{G,H}}\left(\sum_{ic}c^{(i)}\right)^p$$

Consider the elements $S_i=\sum_cc^{(i)}$. With standard notations, we have:
$$S_i=\sum_c(b_{i0}-b_{ic},c)$$

Now observe that these elements multiply as follows:
$$S_{i_1}\ldots S_{i_p}=\sum_{c_1\ldots c_p}
\begin{pmatrix}
b_{i_10}-b_{i_1c_1}+b_{i_2c_1}-b_{i_2,c_1+c_2}&&\\
+b_{i_3,c_1+c_2}-b_{i_3,c_1+c_2+c_3}+\ldots\ldots&,&c_1+\ldots+c_p&\\
\ldots\ldots+b_{i_p,c_1+\ldots+c_{p-1}}-b_{i_p,c_1+\ldots+c_p}&&
\end{pmatrix}$$

In terms of the new indices $d_r=c_1+\ldots+c_r$, this formula becomes:
$$S_{i_1}\ldots S_{i_p}=\sum_{d_1\ldots d_p}
\begin{pmatrix}
b_{i_10}-b_{i_1d_1}+b_{i_2d_1}-b_{i_2d_2}&&\\
+b_{i_3d_2}-b_{i_3d_3}+\ldots\ldots&,&d_p&\\
\ldots\ldots+b_{i_pd_{p-1}}-b_{i_pd_p}&&
\end{pmatrix}$$

Now by integrating, we must have $d_p=0$ on one hand, and on the other hand:
$$[(i_1,0),(i_2,d_1),\ldots,(i_p,d_{p-1})]=[(i_1,d_1),(i_2,d_2),\ldots,(i_p,d_p)]$$

Equivalently, we must have $d_p=0$ on one hand, and on the other hand:
$$[(i_1,d_p),(i_2,d_1),\ldots,(i_p,d_{p-1})]=[(i_1,d_1),(i_2,d_2),\ldots,(i_p,d_p)]$$

Thus, by translation invariance with respect to $d_p$, we obtain:
$$\int_{\widehat{\Gamma}_{G,H}}S_{i_1}\ldots S_{i_p}
=\frac{1}{|H|}\#\left\{d_1,\ldots,d_p\in H\Big|\begin{matrix}[(i_1,d_1),(i_2,d_2),\ldots,(i_p,d_p)]\ \ \ \ \\=[(i_1,d_p),(i_2,d_1),\ldots,(i_p,d_{p-1})]\end{matrix}\right\}$$

It follows that we have the following moment formula:
$$\int_{\widehat{\Gamma}_{G,H}}\left(\sum_iS_i\right)^p
=\frac{1}{|H|}\#\left\{\begin{matrix}i_1,\ldots,i_p\in G\\ d_1,\ldots,d_p\in H\end{matrix}\Big|\begin{matrix}[(i_1,d_1),(i_2,d_2),\ldots,(i_p,d_p)]\ \ \ \ \\=[(i_1,d_p),(i_2,d_1),\ldots,(i_p,d_{p-1})]\end{matrix}\right\}$$

Now by dividing by $|G|$, we obtain the formula in the statement.
\end{proof} 

The formula in Theorem 16.25 can be further interpreted as follows:

\index{Gram matrix}

\begin{theorem}
With $M=|G|,N=|H|$ we have the formula
$$law(\chi)=\left(1-\frac{1}{N}\right)\delta_0+\frac{1}{N}law(A)$$
where the matrix
$$A\in C(\mathbb T^{MN},M_M(\mathbb C))$$
is given by $A(q)=$ Gram matrix of the rows of $q$.
\end{theorem}

\begin{proof}
According to Theorem 16.25, we have the following formula:
\begin{eqnarray*}
\int\chi^p
&=&\frac{1}{MN}\sum_{i_1\ldots i_p}\sum_{d_1\ldots d_p}\delta_{[i_1d_1,\ldots,i_pd_p],[i_1d_p,\ldots,i_pd_{p-1}]}\\
&=&\frac{1}{MN}\int_{\mathbb T^{MN}}\sum_{i_1\ldots i_p}\sum_{d_1\ldots d_p}\frac{q_{i_1d_1}\ldots q_{i_pd_p}}{q_{i_1d_p}\ldots q_{i_pd_{p-1}}}\,dq\\
&=&\frac{1}{MN}\int_{\mathbb T^{MN}}\sum_{i_1\ldots i_p}\left(\sum_{d_1}\frac{q_{i_1d_1}}{q_{i_2d_1}}\right)\left(\sum_{d_2}\frac{q_{i_2d_2}}{q_{i_3d_2}}\right)\ldots\left(\sum_{d_p}\frac{q_{i_pd_p}}{q_{i_1d_p}}\right)dq
\end{eqnarray*}

Consider now the Gram matrix in the statement, namely:
$$A(q)_{ij}=<R_i,R_j>$$

Here $R_1,\ldots,R_M$ are the rows of the following matrix:
$$q\in \mathbb T^{MN}\simeq M_{M\times N}(\mathbb T)$$

We have then the following computation:
\begin{eqnarray*}
\int\chi^p
&=&\frac{1}{MN}\int_{\mathbb T^{MN}}<R_{i_1},R_{i_2}><R_{i_2},R_{i_3}>\ldots<R_{i_p},R_{i_1}>\\
&=&\frac{1}{MN}\int_{\mathbb T^{MN}}A(q)_{i_1i_2}A(q)_{i_2i_3}\ldots A(q)_{i_pi_1}\\
&=&\frac{1}{MN}\int_{\mathbb T^{MN}}Tr(A(q)^p)dq\\
&=&\frac{1}{N}\int_{\mathbb T^{MN}}tr(A(q)^p)dq
\end{eqnarray*}

But this gives the formula in the statement, and we are done.
\end{proof}

In general, the moments of the Gram matrix $A$ are given by a quite complicated formula, and we cannot expect to have a refinement of Theorem 16.26, with $A$ replaced by a plain, non-matricial random variable, say over a compact abelian group. 

\bigskip

However, this kind of simplification does appear at $M=2$, and since this phenomenon is quite interesting, we will explain this now. We first have:

\begin{proposition}
For $F_2\otimes_QF_H$, with $Q\in M_{2\times N}(\mathbb T)$ generic, we have
$$N\int\left(\frac{\chi}{N}\right)^p=\int_{\mathbb T^N}\sum_{k\geq0}\binom{p}{2k}\left|\frac{a_1+\ldots+a_N}{N}\right|^{2k}da$$
where the integral on the right is with respect to the uniform measure on $\mathbb T^N$.
\end{proposition}

\begin{proof}
In order to prove the result, consider the following quantity, which appeared in the proof of Theorem 16.26:
$$\Phi(q)=\sum_{i_1\ldots i_p}\sum_{d_1\ldots d_p}\frac{q_{i_1d_1}\ldots q_{i_pd_p}}{q_{i_1d_p}\ldots q_{i_pd_{p-1}}}$$

We can ``half-dephase'' the matrix $q\in M_{2\times N}(\mathbb T)$ if we want to, as follows:
$$q=\begin{pmatrix}1&\ldots&1\\ a_1&\ldots&a_N\end{pmatrix}$$

Let us compute now the above quantity $\Phi(q)$, in terms of the numbers $a_1,\ldots,a_N$. Our claim is that we have the following formula:
$$\Phi(q)=2\sum_{k\geq0}N^{p-2k}\binom{p}{2k}\left|\sum_ia_i\right|^{2k}$$

Indeed, the idea is that:

\medskip

(1) The $2N^k$ contribution will come from $i=(1\ldots1)$ and $i=(2\ldots2)$.

\medskip

(2) Then we will have a $p(p-1)N^{k-2}|\sum_ia_i|^2$ contribution coming from indices of type $i=(2\ldots 21\ldots1)$, up to cyclic permutations.

\medskip

(3) Then we will have a $2\binom{p}{4}N^{p-4}|\sum_ia_i|^4$ contribution coming from indices of type $i=(2\ldots 21\ldots12\ldots21\ldots1)$.

\medskip

(4) And so on. 

\medskip

In practice now, in order to prove our claim, in order to find the $N^{p-2k}|\sum_ia_i|^{2k}$ contribution, we have to count the circular configurations consisting of $p$ numbers $1,2$, such that the $1$ values are arranged into $k$ non-empty intervals, and the $2$ values are arranged into $k$ non-empty intervals as well. Now by looking at the endpoints of these $2k$ intervals, we have $2\binom{p}{2k}$ choices, and this gives the above formula.

\medskip

Now by integrating, this gives the formula in the statement.
\end{proof}

Observe now that the integrals in Proposition 16.27 can be computed as follows:
\begin{eqnarray*}
\int_{\mathbb T^N}|a_1+\ldots+a_N|^{2k}da
&=&\int_{\mathbb T^N}\sum_{i_1\ldots i_k}\sum_{j_1\ldots j_k}\frac{a_{i_1}\ldots a_{i_k}}{a_{j_1}\ldots a_{j_k}}da\\
&=&\#\left\{i_1\ldots i_k,j_1\ldots j_k\Big|[i_1,\ldots,i_k]=[j_1,\ldots,j_k]\right\}\\
&=&\sum_{k=\sum r_i}\binom{k}{r_1,\ldots,r_N}^2
\end{eqnarray*}

We obtain in this way the following ``blowup'' result, for our measure:

\index{blowup}

\begin{proposition}
For $F_2\otimes_QF_H$, with $Q\in M_{2\times N}(\mathbb T)$ generic, we have
$$\mu=\left(1-\frac{1}{N}\right)\delta_0+\frac{1}{2N}\left(\Psi^+_*\varepsilon+\Psi^-_*\varepsilon\right)$$
where $\varepsilon$ is the uniform measure on $\mathbb T^N$, and where the blowup function is:
$$\Psi^\pm(a)=N\pm\left|\sum_ia_i\right|$$
\end{proposition}

\begin{proof}
We use the formula found in Proposition 16.27, along with the following standard identity, coming from the Taylor formula:
$$\sum_{k\geq0}\binom{p}{2k}x^{2k}=\frac{(1+x)^p+(1-x)^p}{2}$$

By using this identity, Proposition 16.27 reformulates as follows:
$$N\int\left(\frac{\chi}{N}\right)^p=\frac{1}{2}\int_{\mathbb T^N}\left(1+\left|\frac{\sum_ia_i}{N}\right|\right)^p+\left(1-\left|\frac{\sum_ia_i}{N}\right|\right)^p\,da$$

Now by multiplying by $N^{p-1}$, we obtain the following formula:
$$\int\chi^k=\frac{1}{2N}\int_{\mathbb T^N}\left(N+\left|\sum_ia_i\right|\right)^p+\left(N-\left|\sum_ia_i\right|\right)^p\,da$$

But this gives the formula in the statement, and we are done.
\end{proof}

We can further improve the above result, by reducing the maps $\Psi^\pm$ appearing there to a single one, and we are led to the following statement: 

\begin{theorem}
For $F_2\otimes_QF_H$, with $Q\in M_{2\times N}(\mathbb T)$ generic, we have
$$\mu=\left(1-\frac{1}{N}\right)\delta_0+\frac{1}{N}\Phi_*\varepsilon$$
where $\varepsilon$ is the uniform measure on $\mathbb Z_2\times\mathbb T^N$, and where the blowup map is:
$$\Phi(e,a)=N+e\left|\sum_ia_i\right|$$
\end{theorem}

\begin{proof}
This is clear indeed from Proposition 16.28.
\end{proof}

As already mentioned, the above results at $M=2$ are something quite special. In the general case, $M\in\mathbb N$, it is not clear how to construct a nice blowup of the measure.

\bigskip

Asymptotically, things are however quite simple. Let us go back indeed to the general case, where $M,N\in\mathbb N$ are both arbitrary. The problem that we would like to solve now is that of finding the good regime, of the following type, where the measure in Theorem 16.25 converges, after some suitable manipulations:
$$M=f(K)\quad,\quad 
N=g(K)\quad,\quad 
K\to\infty$$

In order to do so, we have to do some combinatorics. Let $NC(p)$ be the set of noncrossing partitions of $\{1,\ldots,p\}$, and for $\pi\in P(p)$ we denote by $|\pi|\in\{1,\ldots,p\}$ the number of blocks. With these conventions, we have the following result:

\begin{proposition}
With $M=\alpha K,N=\beta K$, $K\to\infty$ we have:
$$\frac{c_p}{K^{p-1}}\simeq\sum_{r=1}^p\#\left\{\pi\in NC(p)\Big||\pi|=r\right\}\alpha^{r-1}\beta^{p-r}$$
In particular, with $\alpha=\beta$ we have:
$$c_p\simeq\frac{1}{p+1}\binom{2p}{p}(\alpha K)^{p-1}$$
\end{proposition}

\begin{proof}
We use the combinatorial formula in Theorem 16.25. Our claim is that, with $\pi=\ker(i_1,\ldots,i_p)$, the corresponding contribution to $c_p$ is:
$$C_\pi\simeq
\begin{cases}
\alpha^{|\pi|-1}\beta^{p-|\pi|}K^{p-1}&{\rm if}\ \pi\in NC(p)\\
O(K^{p-2})&{\rm if}\ \pi\notin NC(p)
\end{cases}$$

As a first observation, the number of choices for a multi-index $(i_1,\ldots,i_p)\in X^p$ satisfying $\ker i=\pi$ is:
$$M(M-1)\ldots (M-|\pi|+1)\simeq M^{|\pi|}$$

Thus, we have the following estimate:
$$C_\pi\simeq M^{|\pi|-1}N^{-1}\#\left\{d_1,\ldots,d_p\in Y\Big|[d_\alpha|\alpha\in b]=[d_{\alpha-1}|\alpha\in b],\forall b\in\pi\right\}$$

Consider now the following partition:
$$\sigma=\ker d$$

The contribution of $\sigma$ to the above quantity $C_\pi$ is then given by:
$$\Delta(\pi,\sigma)N(N-1)\ldots(N-|\sigma|+1)\simeq\Delta(\pi,\sigma)N^{|\sigma|}$$

Here the quantities on the right are as follows:
$$\Delta(\pi,\sigma)=\begin{cases}
1&{\rm if}\ |b\cap c|=|(b-1)\cap c|,\forall b\in\pi,\forall c\in\sigma\\
0&{\rm otherwise}
\end{cases}$$

We use now the standard fact that for $\pi,\sigma\in P(p)$ satisfying $\Delta(\pi,\sigma)=1$ we have:
$$|\pi|+|\sigma|\leq p+1$$

In addition, the equality case is well-known to happen when $\pi,\sigma\in NC(p)$ are inverse to each other, via Kreweras complementation. This shows that for $\pi\notin NC(p)$ we have:
$$C_\pi=O(K^{p-2})$$

Also, this shows that for $\pi\in NC(p)$ we have:
\begin{eqnarray*}
C_\pi
&\simeq&M^{|\pi|-1}N^{-1}N^{p-|\pi|-1}\\
&=&\alpha^{|\pi|-1}\beta^{p-|\pi|}K^{p-1}
\end{eqnarray*}

Thus, we have obtained the result.
\end{proof}

We denote by $D$ the dilation operation, given by:
$$D_r(law(X))=law(rX)$$

With this convention, we have the following result:

\begin{theorem}
With $M=\alpha K,N=\beta K$, $K\to\infty$ we have:
$$\mu=\left(1-\frac{1}{\alpha\beta K^2}\right)\delta_0+\frac{1}{\alpha\beta K^2}D_{\frac{1}{\beta K}}(\pi_{\alpha/\beta})$$
In particular with $\alpha=\beta$ we have:
$$\mu=\left(1-\frac{1}{\alpha^2K^2}\right)\delta_0+\frac{1}{\alpha^2K^2}D_{\frac{1}{\alpha K}}(\pi_1)$$
\end{theorem}

\begin{proof}
At $\alpha=\beta$, this follows from Proposition 16.30. In general now, we have:
\begin{eqnarray*}
\frac{c_p}{K^{p-1}}
&\simeq&\sum_{\pi\in NC(p)}\alpha^{|\pi|-1}\beta^{p-|\pi|}\\
&=&\frac{\beta^p}{\alpha}\sum_{\pi\in NC(p)}\left(\frac{\alpha}{\beta}\right)^{|\pi|}\\
&=&\frac{\beta^p}{\alpha}\int x^pd\pi_{\alpha/\beta}(x)
\end{eqnarray*}

When $\alpha\geq\beta$, where $d\pi_{\alpha/\beta}(x)=\varphi_{\alpha/\beta}(x)dx$ is continuous, we obtain:
\begin{eqnarray*}
c_p
&=&\frac{1}{\alpha K}\int(\beta Kx)^p\varphi_{\alpha/\beta}(x)dx\\
&=&\frac{1}{\alpha\beta K^2}\int x^p\varphi_{\alpha/\beta}\left(\frac{x}{\beta K}\right)dx
\end{eqnarray*}

But this gives the formula in the statement. When $\alpha\leq\beta$ the computation is similar, with a Dirac mass as 0 dissapearing and reappearing, and gives the same result.
\end{proof}

We refer to \cite{bbi} and related papers for more on the above.

\section*{16d. Big index}

In big index now, the philosophy is that the index of subfactors $N\in[1,\infty)$ should be regarded as being the well-known $N$ variable from physics, which must be big: 
$$N\to\infty$$

More precisely, the idea is that the constructions involving groups, group duals, or more generally compact quantum groups, producing subfactors of integer index, $N\in\mathbb N$, can be used with ``uniform objects'' as input, and so produce an asymptotic theory. 

\bigskip

The problem however is how to axiomatize the uniformity notion which is needed, in order to have some control on the resulting planar algebra $P=(P_k)$. The answer here comes from the notion of easiness, that we already met in chapter 8, and its various technical extensions, which are in fact not currently unified, or even fully axiomatized.

\bigskip

The main technical questions here are the classification of the easy quantum groups on one hand, and the axiomatization of the super-quizzy quantum groups on the other hand. We also have the question of better understanding the relation between easiness, subfactors, planar algebras, noncommutative geometry and free probability, and we refer here to \cite{bco}, \cite{bcu}, \cite{bgo}, \cite{bpa}, \cite{csn}, \cite{dif}, \cite{lmp}, \cite{nsp}, \cite{twa}, \cite{wen}.

\bigskip

\index{R}

Summarizing, we have many interesting questions, both in small and big index. As a common ground here, both these questions happen inside the Murray-von Neumann factor $R$, although this is conjectural in big index, related to existence questions for outer actions and matrix models. Thus, as a good problem to finish with, which is from Jones' original subfactor paper \cite{jo1}, and is due to Connes, we have the question of axiomatizing the finite index subfactors of the Murray-von Neumann hyperfinite factor $R$. 

\bigskip

As already mentioned on several occasions, this longstanding question is in need of some new, brave functional analysis input, in relation with the notion of hyperfiniteness, which is probably of quite difficult type, beyond what the current experts can do.

\section*{16e. Exercises} 

Congratulations for having read this book, and no exercises here. But, as mentioned above, some good, difficult questions regarding $R$ are waiting for input from you.

\baselineskip=14pt

\printindex

\end{document}